\documentclass[oneside,english,reqno]{amsart}
\usepackage[T1]{fontenc}
\usepackage[utf8]{inputenc}
\usepackage{geometry}
\usepackage{color}
\usepackage{babel}
\usepackage{mathrsfs}
\usepackage{mathtools}
\usepackage{enumitem}
\usepackage{amstext}
\usepackage{amsthm}
\usepackage{amssymb}
\usepackage{stmaryrd}
\usepackage{wasysym}
\makeatletter
\numberwithin{equation}{section}
\numberwithin{figure}{section}
  \theoremstyle{plain}
  \newtheorem*{thm*}{\protect\theoremname}
\theoremstyle{plain}
\newtheorem{thm}{\protect\theoremname}[section]
  \theoremstyle{definition}
  \newtheorem{defn}[thm]{\protect\definitionname}
  \theoremstyle{remark}
  \newtheorem*{rem*}{\protect\remarkname}
  \theoremstyle{plain}
  \newtheorem{lem}[thm]{\protect\lemmaname}
  \theoremstyle{plain}
  \newtheorem{cor}[thm]{\protect\corollaryname}
  \theoremstyle{remark}
  \newtheorem{rem}[thm]{\protect\remarkname}
  \theoremstyle{plain}
  \newtheorem{prop}[thm]{\protect\propositionname}
  \theoremstyle{plain}
  \newtheorem{assumption}[thm]{\protect\assumptionname}

\usepackage[T1]{fontenc}
\usepackage{dsfont}

\usepackage{epigraph}

\let\emptyset\varnothing

\DeclareMathOperator*{\esssup}{ess\,sup}

\theoremstyle{remark}
\newtheorem{remx}[thm]{Remark}
\renewenvironment{rem}
  {\pushQED{\qed}\remx}
  {\popQED\endremx}

\theoremstyle{remark}
\newtheorem*{remxx}{Remark}
\renewenvironment{rem*}
  {\pushQED{\qed}\remxx}
  {\popQED\endremxx}

\theoremstyle{definition}
\newtheorem{assumex}[thm]{Assumption}
\renewenvironment{assumption}
  {\pushQED{\qed}\assumex}
  {\popQED\endassumex}

\theoremstyle{definition}
\newtheorem{defnx}[thm]{Definition}
\renewenvironment{defn}
  {\pushQED{\qed}\defnx}
  {\popQED\enddefnx}

\theoremstyle{plain}
\newtheorem{lemx}[thm]{Lemma}
\renewenvironment{lem}
  {\pushQED{\qed}\lemx}
  {\popQED\endlemx}

\theoremstyle{plain}
\newtheorem{corx}[thm]{Corollary}
\renewenvironment{cor}
  {\pushQED{\qed}\corx}
  {\popQED\endcorx}

\theoremstyle{plain}
\newtheorem{propx}[thm]{Proposition}
\renewenvironment{prop}
  {\pushQED{\qed}\propx}
  {\popQED\endpropx}

\theoremstyle{definition}
\newtheorem{thmx}[thm]{Theorem}
\renewenvironment{thm}
  {\pushQED{\qed}\thmx}
  {\popQED\endthmx}

\usepackage[ddmmyyyy,hhmmss]{datetime}
\newdateformat{daymonthyeardate}{%
  \THEDAY.\THEMONTH.\THEYEAR}
\usepackage{fancyhdr}
\makeatletter

\newsavebox{\@linebox}
\savebox{\@linebox}[3em][t]{\parbox[t]{3em}{%
  \@tempcnta\@ne\relax
  \loop{\underline{\scriptsize\the\@tempcnta}}\\
    \advance\@tempcnta by \@ne\ifnum\@tempcnta<55\repeat}}

\makeatother

\makeatother

  \providecommand{\assumptionname}{Assumption}
  \providecommand{\corollaryname}{Corollary}
  \providecommand{\definitionname}{Definition}
  \providecommand{\lemmaname}{Lemma}
  \providecommand{\propositionname}{Proposition}
  \providecommand{\remarkname}{Remark}
  \providecommand{\theoremname}{Theorem}
\providecommand{\theoremname}{Theorem}

\begin{document}

\title{Structured, Compactly Supported Banach Frame Decompositions\\
of Decomposition Spaces}

\date{}{}


\author{Felix Voigtlaender}

\maketitle
\global\long\def\vertiii#1{{\left\vert \kern-0.25ex  \left\vert \kern-0.25ex  \left\vert #1\right\vert \kern-0.25ex  \right\vert \kern-0.25ex  \right\vert }}
\global\long\def\essup{\esssup}
\global\long\def\with{\,\middle|\,}
\global\long\def\DecompSp#1#2#3#4{{\mathcal{D}\left({#1},L_{#4}^{#2},{#3}\right)}}
\global\long\def\FourierDecompSp#1#2#3#4{{\mathcal{D}_{\mathcal{F}}\left({#1},L_{#4}^{#2},{#3}\right)}}
\global\long\def\BAPUFourierDecompSp#1#2#3#4{{\mathcal{D}_{\mathcal{F},{#4}}\left({#1},L^{#2},{#3}\right)}}
\global\long\def\R{\mathbb{R}}
\global\long\def\Compl{\mathbb{C}}
\global\long\def\N{\mathbb{N}}
\global\long\def\Z{\mathbb{Z}}
\global\long\def\CalQ{\mathcal{Q}}
\global\long\def\CalP{\mathcal{P}}
\global\long\def\CalR{\mathcal{R}}
\global\long\def\CalS{\mathcal{S}}
\global\long\def\CalD{\mathcal{D}}
\global\long\def\BesovInhom#1#2#3{\mathcal{B}_{#3}^{#1,#2}}
\global\long\def\BesovHom#1#2#3{\dot{\mathcal{B}}_{#3}^{#1,#2}}
\global\long\def\dimension{d}
\global\long\def\ModSpace#1#2#3{M_{#3}^{#1,#2}}
\global\long\def\AlphaModSpace#1#2#3#4{M_{#3,#4}^{#1,#2}}
\global\long\def\ShearletSmoothness#1#2#3{S_{#3}^{#1,#2}}
\global\long\def\GL{\mathrm{GL}}
\global\long\def\CalO{\mathcal{O}}
\global\long\def\Coorbit{\mathrm{Co}}
\global\long\def\d{\operatorname{d}}
\global\long\def\TestFunctionSpace#1{C_{c}^{\infty}\left(#1\right)}
\global\long\def\DistributionSpace#1{\mathcal{D}'\left(#1\right)}
\global\long\def\SpaceTestFunctions#1{Z\left(#1\right)}
\global\long\def\SpaceReservoir#1{Z'\left(#1\right)}
\global\long\def\Schwartz{\mathcal{S}}
\global\long\def\Fourier{\mathcal{F}}
\global\long\def\supp{\operatorname{supp}}
\global\long\def\dist{\operatorname{dist}}
\global\long\def\Xhookrightarrow#1{\xhookrightarrow{#1}}
\global\long\def\Xmapsto#1{\xmapsto{#1}}
\global\long\def\Indicator{{\mathds{1}}}
\global\long\def\identity{\operatorname{id}}
\global\long\def\RaisedSup#1{\raisebox{0.15cm}{\mbox{\ensuremath{{\displaystyle \sup_{#1}}}}}}
\global\long\def\Raised#1{\raisebox{0.12cm}{\mbox{\ensuremath{#1}}}}
\global\long\def\osc#1{\operatorname{osc}_{#1}}
\global\long\def\mybullet{\bullet}

\begin{abstract}
We present a very general framework for the construction of structured,
possibly compactly supported Banach frames and atomic decompositions
for a given decomposition space. Here, a decomposition space $\DecompSp{\CalQ}p{\ell_{w}^{q}}{}$
is defined essentially like a classical Besov space, but the usual
dyadic covering is replaced by an (almost) arbitrary covering $\CalQ=\left(Q_{i}\right)_{i\in I}$.
Thus, if $\Phi=\left(\varphi_{i}\right)_{i\in I}$ is a suitable partition
of unity subordinate to $\CalQ$, then $\left\Vert g\right\Vert _{\DecompSp{\CalQ}p{\ell_{w}^{q}}{}}:=\left\Vert \left(\left\Vert \Fourier^{-1}\left(\varphi_{i}\widehat{g}\right)\right\Vert _{L^{p}}\right)_{i\in I}\right\Vert _{\ell_{w}^{q}}$.
Special cases include the class of Besov spaces and ($\alpha$-)modulation
spaces, as well as a large class of wavelet-type coorbit spaces and
so-called shearlet smoothness spaces.

Assuming that the covering $\CalQ$ is of the regular form $\CalQ=\left(T_{i}Q+b_{i}\right)_{i\in I}$,
with $T_{i}\in\GL\left(\R^{\dimension}\right),b_{i}\in\R^{\dimension}$,
we fix a \emph{prototype function} $\gamma\in L^{1}\left(\R^{\dimension}\right)$
and consider the structured generalized shift invariant system
\[
\Psi_{\delta}:=\left(L_{\delta\cdot T_{i}^{-T}k}\:\gamma^{\left[i\right]}\right)_{i\in I,k\in\Z^{\dimension}}\quad\text{ with }\quad\gamma^{\left[i\right]}:=\left|\det T_{i}\right|^{1/2}\cdot M_{b_{i}}\left(\gamma\circ T_{i}^{T}\right),
\]
where $L_{x}$ and $M_{\xi}$ denote translation and modulation, respectively.
The main contribution of the paper is to provide verifiable conditions
on the prototype $\gamma$ which ensure that $\Psi_{\delta}$ forms,
respectively, a Banach frame or an atomic decomposition for the space
$\DecompSp{\CalQ}p{\ell_{w}^{q}}{}$, for sufficiently small \emph{sampling
density} $\delta>0$. Crucially, while the decomposition space $\DecompSp{\CalQ}p{\ell_{w}^{q}}{}$
is defined using the \emph{bandlimited} family $\Phi$, the construction
presented here usually allows for the prototype $\gamma$ to be \emph{compactly
supported} in space. We emphasize that the theory presented here can
cover the whole range $p,q\in\left(0,\infty\right]$ and not only
the case $p,q\in\left[1,\infty\right]$ of Banach spaces.

An important feature of our theory is that in many cases, the system
$\Psi_{\delta}$ will \emph{simultaneously} form a Banach frame and
an atomic decomposition for $\DecompSp{\CalQ}p{\ell_{w}^{q}}{}$.
This implies that for frames of the form $\Psi_{\delta}$, \emph{analysis
sparsity is equivalent to synthesis sparsity}, i.e., the analysis
coefficients $\left(\left\langle f,\,L_{\delta\cdot T_{i}^{-T}k}\:\gamma^{\left[i\right]}\right\rangle \right)_{i\in I,k\in\Z^{\dimension}}$
lie in $\ell^{p}\left(I\times\Z^{\dimension}\right)$ if and only
if $f$ is an element of a certain decomposition space, if and only
if $f=\sum_{i\in I,k\in\Z^{\dimension}}\left[c_{k}^{\left(i\right)}\cdot L_{\delta\cdot T_{i}^{-T}k}\gamma^{\left[i\right]}\right]$
for some sequence $\left(\smash{c_{k}^{\left(i\right)}}\right)_{i\in I,k\in\Z^{\dimension}}\in\ell^{p}\left(I\times\Z^{\dimension}\right)$.
This is very convenient, since for many frame constructions—like shearlets—one
only knows that the \emph{analysis coefficients} for a class of ``nice''
signals are sparse. This, however, only entails synthesis sparsity
with respect to the \emph{dual} frame, about which often only limited
knowledge is available. Using the theory presented here, one can derive
synthesis sparsity with respect to the \emph{primal} frame, for which
one has an explicit formula and whose properties like smoothness and
time-frequency localization are well understood.

As a sample application, we show that the developed theory applies
to $\alpha$-modulation spaces and to (inhomogeneous) Besov spaces.
In a companion paper, we also show that the theory applies to shearlet
smoothness spaces.
\end{abstract}

\section{Introduction}

\label{sec:Introduction}In this section, we first motivate and describe
our approach for the construction of structured Banach frame decompositions
for decomposition spaces and compare our results to the known literature.
Then, we introduce a few standard and non-standard conventions and
standing assumptions which are used in the remainder of the paper.
Finally, we give a brief overview over the structure of the paper.

\subsection{Motivation and comparison to known results}

Given a Banach space $X$, a family $\Psi=\left(\psi_{i}\right)_{i\in I}$
in $X'$ is called a \textbf{Banach frame\cite{GroechenigDescribingFunctions}}
for $X$ if there is a \textbf{solid sequence space} $Y\leq\Compl^{I}$
such that
\begin{itemize}
\item the \textbf{analysis operator} $A_{\Psi}:X\to Y,x\mapsto\left(\left\langle x,\psi_{i}\right\rangle _{X,X'}\right)_{i\in I}$
is well-defined and bounded,
\item there is a bounded linear \textbf{reconstruction operator} $R:Y\to X$
satisfying $R\circ A_{\Psi}=\identity_{X}$.
\end{itemize}
In particular, this implies $\left\Vert x\right\Vert _{X}\asymp\left\Vert A_{\Psi}x\right\Vert _{Y}$
uniformly over $x\in X$. Here, a Banach space $Y\leq\Compl^{I}$
is called \textbf{solid} if for all sequences $x=\left(x_{i}\right)_{i\in I}$
and $y=\left(y_{i}\right)_{i\in I}$ with $y\in Y$ and $\left|x_{i}\right|\leq\left|y_{i}\right|$
for all $i\in I$, it follows that $x\in Y$ with $\left\Vert x\right\Vert _{Y}\leq\left\Vert y\right\Vert _{Y}$.

Dual to the notion of a Banach frame, a family $\Phi=\left(\varphi_{i}\right)_{i\in I}$
in $X$ is called an \textbf{atomic decomposition\cite{GroechenigDescribingFunctions}}
for $X$ if there is a solid sequence space $Z\leq\Compl^{I}$ such
that
\begin{itemize}
\item the \textbf{synthesis operator} $S_{\Phi}:Z\to X,\left(x_{i}\right)_{i\in I}\mapsto\sum_{i\in I}x_{i}\varphi_{i}$
is well-defined and bounded, where convergence of the series occurs
in a suitable (weak) sense,
\item there is a bounded linear \textbf{coefficient operator} $C:X\to Z$
satisfying $S_{\Phi}\circ C=\identity_{X}$.
\end{itemize}
In particular, this implies that every $x\in X$ can be written as
$x=\sum_{i\in I}c_{i}\varphi_{i}$ for a suitable sequence $c=\left(c_{i}\right)_{i\in I}=Cx$.

The existence of nice Banach frames and atomic decompositions for
a given (family of) Banach space(s) is extremely convenient, since
the study of many properties like existence of embeddings, boundedness
of operators and description of interpolation spaces, etc., of the
Banach spaces under consideration can be reduced to studying these
properties for the associated sequence spaces, which are often much
easier to understand.

For this reason, much effort has been spent to derive existence of
Banach frames and atomic decompositions for many well-known spaces
like Besov- and Sobolev spaces. The most well-known types of (Banach)
frames are probably the \textbf{wavelet characterization} of Besov
spaces (see e.g.\@ \cite[Theorem 1.64]{TriebelTheoryOfFunctionSpaces3}),
the closely related characterization of these spaces using the $\varphi$-transform
\cite{FrazierJawerthDecompositionOfBesovSpaces,FrazierJawerthDiscreteTransform},
as well as the existence of \textbf{Gabor frames} for modulation spaces\cite{GroechenigTimeFrequencyAnalysis}.

\subsubsection{Classical group-based coorbit theory}

By generalizing the similarities between the theories of wavelet-
and Gabor frames, Feichtinger and Gröchenig initiated the study of
so-called \textbf{coorbit spaces}\cite{FeichtingerCoorbit0,FeichtingerCoorbit1,FeichtingerCoorbit2,GroechenigDescribingFunctions},
which provide a systematic way of obtaining Banach frames and atomic
decompositions for certain Banach spaces. Precisely, one starts with
an irreducible, (square)-integrable representation $\pi:G\to\mathcal{U}\left(\mathcal{H}\right)$
of some locally compact Hausdorff (LCH) topological group $G$. This
representation induces for each $g\in\mathcal{H}$ an associated \textbf{voice
transform}
\[
V_{g}:\mathcal{H}\to C\left(G\right),f\mapsto V_{g}f\qquad\text{ where }\qquad\left(V_{g}f\right)\left(x\right)=\left\langle f,\,\pi\left(x\right)g\right\rangle _{\mathcal{H}}.
\]
For an \textbf{admissible vector} $\psi\in\mathcal{H}\setminus\left\{ 0\right\} $
(which means $V_{\psi}\psi\in L^{2}\left(G\right)$), it follows\cite{DufloMoore}
that $V_{\psi}:\mathcal{H}\to L^{2}\left(G\right)$ is (a scalar multiple
of) an isometry, so that $\left(\pi\left(x\right)\psi\right)_{x\in G}$
is a \textbf{tight continuous frame} for $\mathcal{H}$, since
\[
\left\Vert f\right\Vert _{\mathcal{H}}^{2}=C_{\psi}\cdot\int_{G}\left|\left(V_{\psi}f\right)\left(x\right)\right|^{2}\d\mu\left(x\right)\qquad\forall f\in\mathcal{H}.
\]
In particular, this identity implies $\mathcal{H}=\left\{ f\with V_{\psi}f\in L^{2}\left(G\right)\right\} $.
In generalization of this identity, coorbit theory shows that for
``good enough'' \emph{analyzing windows} $\psi$ and each suitable,
\emph{solid} function space $Y\leq L_{{\rm loc}}^{1}\left(G\right)$,
one can define the associated \textbf{coorbit space} as
\[
{\rm Co}\left(Y\right):=\left\{ f\in\CalR\with V_{\psi}f\in Y\right\} \qquad\text{ with norm }\qquad\left\Vert f\right\Vert _{{\rm Co}\left(Y\right)}=\left\Vert V_{\psi}f\right\Vert _{Y}.
\]
Here, $\CalR=\CalR_{Y}$ is a suitable \emph{reservoir}. Informally,
$\CalR$ corresponds to the set of (tempered) distributions; but due
to the generality in which coorbit spaces are defined, one has to
use a slightly different definition, intrinsic to the group $G$,
cf.\@ \cite[Section 4]{FeichtingerCoorbit1}.

The main statement of coorbit theory is that one can \emph{discretize}
the (continuous, tight) frame $\left(\pi\left(x\right)\psi\right)_{x\in G}$,
to obtain discrete Banach frames and atomic decompositions, simultaneously
for all spaces ${\rm Co}\left(Y\right)$, where $Y$ ranges over a
suitable set of solid function spaces on $Y$. More precisely, the
following are true:
\begin{itemize}
\item For each translation invariant, solid function space $Y\leq L_{{\rm loc}}^{1}\left(G\right)$,
there is a so-called \textbf{control weight} $w=w_{Y}:G\to\left(0,\infty\right)$,
cf.\@ \cite[equation (4.10)]{FeichtingerCoorbit1}. For the following
statements, we always assume that $w$ is a control weight for $Y$.
\item Associated to each control weight $w$, there is a class $\mathcal{B}_{w}\subset\mathcal{H}$
of \textbf{good (analyzing) vectors} such that for each two $\psi_{1},\psi_{2}\in\mathcal{B}_{w}\setminus\left\{ 0\right\} $,
the identity
\[
\left\{ f\in\mathcal{R}\with V_{\psi_{1}}f\in Y\right\} ={\rm Co}\left(Y\right)=\left\{ f\in\mathcal{R}\with V_{\psi_{2}}f\in Y\right\} 
\]
holds, i.e., one has a \emph{consistency statement}.
\item For each control weight $w$ and each $\psi\in\mathcal{B}_{w}\setminus\left\{ 0\right\} $,
there is a unit neighborhood $U=U\left(\psi,w\right)\subset G$, such
that for every $U$-dense and \textbf{relatively separated} family
$X=\left(x_{i}\right)_{i\in I}$ in $G$, the family $\left(\pi\left(x_{i}\right)\psi\right)_{i\in I}$
forms an \textbf{atomic decomposition} of ${\rm Co}\left(Y\right)$,
i.e., there is a solid, discrete sequence space $Y_{d}\left(X\right)\leq\Compl^{I}$
associated to $Y$ such that the following hold (see \cite[Theorem 6.1 and the associated remark]{FeichtingerCoorbit1}):

\begin{itemize}
\item the synthesis operator
\[
S:Y_{d}\left(X\right)\to{\rm Co}\left(Y\right),\left(\lambda_{i}\right)_{i\in I}\mapsto\sum_{i\in I}\lambda_{i}\cdot\pi\left(x_{i}\right)\psi
\]
is well-defined and bounded (with convergence in the weak-$\ast$-topology
of the reservoir $\mathcal{R}$),
\item there is a bounded linear operator $C:{\rm Co}\left(Y\right)\to Y_{d}\left(X\right)$
satisfying $S\circ C=\identity_{{\rm Co}\left(Y\right)}$.
\end{itemize}
\item For each control weight $w$ and each $\psi\in\mathcal{B}_{w}\setminus\left\{ 0\right\} $,
there is a unit neighborhood $V=V\left(\psi,w\right)\subset G$ such
that for every $V$-dense and relatively separated family $X=\left(x_{i}\right)_{i\in I}$
in $G$, the family $\left(\pi\left(x_{i}\right)\psi\right)_{i\in I}$
forms a \textbf{Banach frame} for ${\rm Co}\left(Y\right)$, i.e.,
with the same solid sequence space $Y_{d}\left(X\right)$ as above,
the following hold (see \cite[Theorem 5.3]{GroechenigDescribingFunctions}):

\begin{itemize}
\item the analysis operator
\[
A:{\rm Co}\left(Y\right)\to Y_{d}\left(X\right),f\mapsto\left(\left\langle f,\pi\left(x_{i}\right)\psi\right\rangle \right)_{i\in I}
\]
is well-defined and bounded,
\item there is a bounded linear operator $R:Y_{d}\left(X\right)\to{\rm Co}\left(Y\right)$
satisfying $R\circ A=\identity_{{\rm Co}\left(Y\right)}$.
\end{itemize}
\end{itemize}
Here, a family $X=\left(x_{i}\right)_{i\in I}$ is called \textbf{$V$-dense}
if $G=\bigcup_{i\in I}x_{i}V$ and \textbf{relatively separated} if
it is a finite union of separated sets, where a family $Z=\left(z_{j}\right)_{j\in J}$
is called \textbf{separated} if there is a unit neighborhood $W\subset G$
satisfying $z_{j}W\cap z_{\ell}W=\emptyset$ for $j\neq\ell$.

Among other examples, this \emph{group-based} coorbit theory can be
used to obtain Banach frames and atomic decompositions for modulation
spaces as well as for \emph{homogeneous} Besov spaces. There are also
several extensions, for example to the setting of Quasi-Banach spaces\cite{RauhutCoorbitQuasiBanach},
and to the setting of possibly reducible or non-integrable group representations\cite{CoorbitWithVoiceInFrechetSpace}.

The main limitation of this theory, however, is that many relevant
spaces like \emph{inhomogeneous} Besov spaces are \emph{not} covered
by it.

\subsubsection{Generalized coorbit theory}

To overcome this limitation, Fornasier, Rauhut and Ullrich\cite{GeneralizedCoorbit1,GeneralizedCoorbit2}
developed what is called \textbf{generalized coorbit theory}; see
also \cite{NickiErrataForRauhut} for some corrections and extensions
and \cite{QuasiBanachGeneralCoorbit} for a generalization to Quasi-Banach
spaces. For generalized coorbit theory, one starts from a Hilbert
space $\mathcal{H}$, for which one is given a \textbf{continuous
frame} $\Psi=\left(\psi_{x}\right)_{x\in X}$ which is indexed by
some locally compact measure space $X$, equipped with a Radon measure
$\mu$. Formally, this means that for each $f\in\mathcal{H}$, the
function $X\to\Compl,x\mapsto\left\langle f,\,\psi_{x}\right\rangle _{\mathcal{H}}$
is measurable and there are constants $0<A\leq B$ satisfying
\[
A\cdot\left\Vert f\right\Vert _{\mathcal{H}}^{2}\leq\int_{X}\left|\left\langle f,\,\psi_{x}\right\rangle _{\mathcal{H}}\right|^{2}\d\mu\left(x\right)\leq B\cdot\left\Vert f\right\Vert _{\mathcal{H}}^{2}\qquad\forall f\in\mathcal{H}.
\]
In this case, the \textbf{frame operator} $S:\mathcal{H}\to\mathcal{H},f\mapsto\int_{X}\left\langle f,\psi_{x}\right\rangle _{\mathcal{H}}\cdot\psi_{x}\d\mu\left(x\right)$
(with the integral understood in the weak sense) is well-defined,
self-adjoint and positive and thus invertible. Hence, one can form
the \textbf{canonical dual frame} $\tilde{\Psi}=\left(\smash{\tilde{\psi_{x}}}\right)_{x\in X}=\left(S^{-1}\psi_{x}\right)_{x\in X}$.
If the frame $\Psi$ is \textbf{tight}, one can choose $A=B$ and
the dual frame $\tilde{\Psi}$ is simply a scalar multiple of $\Psi$,
but in general, $\Psi$ and $\tilde{\Psi}$ might be very different.
For generalized coorbit theory to be applicable at all, one requires
the \textbf{cross-gramian kernel}
\[
R:X\times X\to\Compl,\left(x,y\right)\mapsto\left\langle \psi_{y},\,S^{-1}\psi_{x}\right\rangle =\left\langle \psi_{y},\,\tilde{\psi_{x}}\right\rangle 
\]
to have certain decay/mapping properties; precisely, one requires
$R\in\mathcal{A}_{m}$, where $m=m_{Y}$ is a given control weight
associated to the solid function space $Y\leq L_{{\rm loc}}^{1}\left(X\right)$
in which one is interested. Here, $\mathcal{A}_{m}$ is a suitable
algebra of kernels, cf.\@ \cite[Section 3]{GeneralizedCoorbit1}.

With these two frames $\Psi,\tilde{\Psi}$, there are \emph{two} associated
voice transforms, given by
\[
V_{\Psi}f\left(x\right):=\left\langle f,\,\psi_{x}\right\rangle _{\mathcal{H}}\qquad\text{ and }\qquad W_{\Psi}f\left(x\right):=\left\langle f,\,\smash{\tilde{\psi_{x}}}\right\rangle _{\mathcal{H}}=\left(V_{\Psi}\left[S^{-1}f\right]\right)\left(x\right),
\]
and then (cf.\@ \cite[equation (3.8) and Definition 3.1]{GeneralizedCoorbit1})
also \emph{two} reservoirs $\mathcal{R}_{1}:=\left(\mathcal{K}_{v}^{1}\right)^{\neg}$
and $\mathcal{R}_{2}:=\left(\mathcal{H}_{v}^{1}\right)^{\neg}$ and
\emph{two} coorbit spaces
\[
{\rm Co}\left(Y\right):=\left\{ f\in\mathcal{R}_{1}\with V_{\Psi}f\in Y\right\} \qquad\text{ and }\qquad\widetilde{{\rm Co}}\left(Y\right):=\left\{ f\in\mathcal{R}_{2}\with W_{\Psi}f\in Y\right\} .
\]

Then, if the frame $\Psi$ is ``good enough'' (the precise meaning
of which depends on the space $Y$), one can again discretize the
continuous frame $\Psi$ to obtain atomic decompositions and Banach
frames. However, one has to be a bit careful; assuming that the family
$\left(x_{i}\right)_{i\in I}$ is ``dense enough in $X$'' (cf.\@
\cite[Theorem 5.7]{GeneralizedCoorbit1} for the details), we have
the following:
\begin{itemize}
\item the family $\left(\psi_{x_{i}}\right)_{i\in I}$ is an atomic decomposition
of $\widetilde{{\rm Co}}\left(Y\right)$ with corresponding sequence
space $Y^{\natural}$,
\item the family $\left(\psi_{x_{i}}\right)_{i\in I}$ is a Banach frame
for ${\rm Co}\left(Y\right)$ with corresponding sequence space $Y^{\flat}$.
\end{itemize}
Thus, although generalized coorbit theory is immensely powerful and
general, its main limitation is that one essentially has to start
from a \emph{tight} continuous frame for a Hilbert space $\mathcal{H}$,
since in the non-tight case one faces several limitations:
\begin{itemize}
\item In most cases of continuous \emph{non-tight} frames, one knows very
little about the properties of the (canonical) dual frame $\tilde{\Psi}$,
which makes it hard to verify that the kernel $R\left(x,y\right)=\left\langle \psi_{y},\,\tilde{\psi_{x}}\right\rangle $
satisfies $R\in\mathcal{A}_{m}$.
\item As seen above, one is faced with two \emph{distinct} coorbit spaces
${\rm Co}\left(Y\right)$ and $\widetilde{{\rm Co}}\left(Y\right)$
and obtains a Banach frame for ${\rm Co}\left(Y\right)$ and an atomic
decomposition for $\widetilde{{\rm Co}}\left(Y\right)$. In many cases,
however, it is desired to \emph{simultaneously} have a Banach frame
and an atomic decomposition for \emph{one} common space.

We mention that \cite[Section 4]{GeneralizedCoorbit1} provides criteria
which ensure ${\rm Co}\left(Y\right)=\widetilde{{\rm Co}}\left(Y\right)$,
namely if $\Psi$ and $\tilde{\Psi}$ are $\mathcal{A}_{m}$-self-localized.
To show that this is true, however, one again needs to know a lot
about the dual frame $\tilde{\Psi}$, which in general one does not.
The most convenient way out (outlined in \cite[Theorem 4.7 and the comments afterward]{GeneralizedCoorbit1})
is to find a suitable \textbf{spectral} subalgebra $\mathcal{A}$
of $\mathcal{A}_{m}$ and then to show that the kernel $K\left(x,y\right)=\left\langle \psi_{y},\,\psi_{x}\right\rangle $
satisfies $K\in\mathcal{A}$. Once this is shown, \cite[Theorem 4.7]{GeneralizedCoorbit1}
yields ${\rm Co}\left(Y\right)=\widetilde{{\rm Co}}\left(Y\right)$
as well as $R\in\mathcal{A}\subset\mathcal{A}_{m}$, so that coorbit
theory is applicable. The main limitation of this approach is that
not too many spectral algebras of kernels are known.
\end{itemize}
In total, there are two desirable use cases of (generalized) coorbit
theory in which an actual application is \emph{difficult, or even
impossible:}
\begin{enumerate}
\item In the first case, one is given a (family of) Banach space(s) $B$
and wants to find Banach frames and atomic decompositions for $B$.
To achieve this via (generalized) coorbit theory, one has to find
a (preferably tight) \emph{continuous} frame $\Psi=\left(\psi_{x}\right)_{x\in X}$
and a (family of) solid Banach function space(s) $Y\leq L_{{\rm loc}}^{1}\left(X\right)$
such that $B={\rm Co}\left(Y\right)=\left\{ f\with V_{\Psi}f\in Y\right\} $.
Furthermore, one has to verify that $\Psi$ indeed satisfies all prerequisites
for the application of generalized coorbit theory. Finally, if $\Psi$
is non-tight, one has to verify ${\rm Co}\left(Y\right)=\widetilde{{\rm Co}}\left(Y\right)$,
for example by using the approach using spectral algebras which we
outlined above.
\item In the second case, which occurs e.g.\@ if one wants to study the
approximation theoretic properties of discrete, cone-adapted shearlet
frames\cite{CompactlySupportedShearletFrames}, one starts with a
\emph{discrete} frame $\Psi_{d}=\left(\psi_{i}\right)_{i\in I}$ (or
with a family of such discrete frames, e.g., parametrized by the sampling
density) for a Hilbert space $\mathcal{H}$ and one wants to understand
the space of those functions which are \emph{analysis-sparse} with
respect to this frame, e.g., the space
\[
B_{q}:=\left\{ f\in\mathcal{H}\with\left(\left\langle f,\,\psi_{i}\right\rangle \right)_{i\in I}\in\ell^{q}\left(I\right)\right\} \qquad\text{ for }\qquad q<2.
\]
An important property one might be interested in is \emph{whether
analysis sparsity is equivalent to synthesis sparsity}, i.e., whether
every $f\in B_{q}$ admits an expansion $f=\sum_{i\in I}c_{i}\psi_{i}$
for a sequence $c=\left(c_{i}\right)_{i\in I}\in\ell^{q}\left(I\right)$.

To derive such a statement using coorbit theory, one needs to find
a continuous (preferably tight) frame $\Psi=\left(\psi_{x}\right)_{x\in X}$
for $\mathcal{H}$ such that the discretization $\left(\psi_{x_{i}}\right)_{i\in I}$
of this frame (in the sense of generalized coorbit theory) is equal
to $\Psi_{d}$. Then, provided that ${\rm Co}\left(Y\right)=\widetilde{{\rm Co}}\left(Y\right)=B_{q}$,
coorbit theory will yield the desired statement.

The main problem here—as witnessed by the example of discrete cone-adapted
shearlets—is that it can often be very hard or even impossible to
find such a continuous frame $\Psi$, much less a tight one. There
are tight continuous shearlet frames, e.g.\@ those related to shearlet
coorbit spaces\cite{Dahlke_etal_sh_coorbit1,Dahlke_etal_sh_coorbit2,DahlkeShearletArbitraryDimension,DahlkeShearletCoorbitEmbeddingsInHigherDimensions},
but a discretization of these frames does \emph{not} yield discrete
\emph{cone-adapted} shearlet systems.
\end{enumerate}
As we will see now, our approach does \emph{not} require to have a
continuous frame which can then be discretized. Thus, in this aspect,
our approach improves upon (generalized) coorbit theory. As we will
see in the companion paper \cite{StructuredBanachFrames2}, we are
in particular able to handle discrete cone-adapted shearlet frames;
and for this case, our theory indeed shows that \emph{analysis sparsity
is equivalent to synthesis sparsity}.

\subsubsection{Our approach using decomposition spaces}

For our approach, we start with a \textbf{structured covering} $\CalQ$
of (an open subset $\CalO$ of) the frequency space $\R^{\dimension}$.
More precisely (see Subsection \ref{subsec:DecompSpaceDefinitionStandingAssumptions}
for the completely formal assumptions), we assume that
\begin{equation}
\CalQ=\left(Q_{i}\right)_{i\in I}=\left(T_{i}Q+b_{i}\right)_{i\in I}\label{eq:IntroductionStructuredCoveringAssumption}
\end{equation}
for a fixed open, precompact set $Q\subset\R^{\dimension}$ and certain
linear maps $T_{i}\in\GL\left(\R^{\dimension}\right)$ and translations
$b_{i}\in\R^{\dimension}$. Then, given a suitable partition of unity
$\Phi=\left(\varphi_{i}\right)_{i\in I}$ subordinate to $\CalQ$
and a suitable weight $w=\left(w_{i}\right)_{i\in I}$ on $I$, as
well as $p,q\in\left(0,\infty\right]$, one defines the decomposition
space (quasi)-norm of a distribution $g$ as
\[
\left\Vert g\right\Vert _{\DecompSp{\CalQ}p{\ell_{w}^{q}}{}}:=\left\Vert \left(\left\Vert \Fourier^{-1}\left(\varphi_{i}\cdot\widehat{g}\right)\right\Vert _{L^{p}}\right)_{i\in I}\right\Vert _{\ell_{w}^{q}}=\left\Vert \left(\left\Vert \left(\Fourier^{-1}\varphi_{i}\right)\ast g\right\Vert _{L^{p}}\right)_{i\in I}\right\Vert _{\ell_{w}^{q}}\,,
\]
while the \textbf{decomposition space} $\DecompSp{\CalQ}p{\ell_{w}^{q}}{}$
consists of all distributions for which this (quasi)-norm is finite.
For the exact interpretation of ``distribution'' in this context,
we refer to Subsection \ref{subsec:DecompSpaceDefinitionStandingAssumptions}.
In words, the decomposition space norm is computed by first decomposing
$g$ \emph{in frequency} according to the covering $\CalQ$ to obtain
the pieces $g_{i}=\Fourier^{-1}\left(\varphi_{i}\cdot\widehat{g}\right)$.
Each of these pieces is then measured in $L^{p}$ and the overall
norm is a certain $\ell_{w}^{q}$-norm over all of these contributions.
In most of the paper, we will even consider the weighted $L^{p}$-spaces
$L_{v}^{p}$ instead of $L^{p}$. But in this introduction, we will
mostly stick to the setting just described, for the sake of simplicity.

Our general aim is to show that one can obtain compactly supported
Banach frames and atomic decompositions $\Psi$ of a very special,
structured form for the decomposition space $\DecompSp{\CalQ}p{\ell_{w}^{q}}{}$.
In fact, it will turn out that the system $\Psi$ can be taken to
be a generalized shift invariant system generated by a \emph{single
prototype function} $\gamma$, similar to the way in which a prototype
function can generate Gabor, wavelet and shearlet systems.

To see exactly how such a system $\Psi$ might look like, let us write
$S_{i}\xi:=T_{i}\xi+b_{i}$ for $i\in I$. Note that if $\supp\widehat{\gamma}\subset Q$,
then $\supp\left[\widehat{\gamma}\circ S_{i}^{-1}\right]\subset Q_{i}$.
The same remains true in a weak sense if the strict inclusion $\supp\widehat{\gamma}\subset Q$
is replaced by requiring that $\widehat{\gamma}$ be \emph{essentially}
supported in $Q$, which can even hold if $\gamma$ is not band-limited.
Now, note $\Fourier\gamma^{\left(i\right)}=\widehat{\gamma}\circ S_{i}^{-1}$
for $\gamma^{\left(i\right)}:=\left|\det T_{i}\right|\cdot M_{b_{i}}\left[\gamma\circ T_{i}^{T}\right]$.
For consistency with the $L^{2}$-setting, we also consider 
\begin{equation}
\gamma^{\left[i\right]}:=\left|\det T_{i}\right|^{1/2}\cdot M_{b_{i}}\left[\gamma\circ T_{i}^{T}\right].\label{eq:IntroductionL2NormalizedFrameElement}
\end{equation}
In fact, we will even allow the generator $\gamma$ to vary with $i\in I$,
i.e., $\gamma^{\left[i\right]}=\left|\det T_{i}\right|^{1/2}\cdot M_{b_{i}}\left[\gamma_{i}\circ T_{i}^{T}\right]$.
An example where this is useful is an inhomogeneous wavelet system:
If the generator $\gamma$ is required to be \emph{independent} of
$i\in I$, the ``low-pass part'' of the wavelet system needs to
be obtained by a frequency shift (i.e., by a modulation) from the
mother wavelet $\gamma$. Indeed, since we consider only affine dilations
of $\widehat{\gamma}$ and since any \emph{linear} dilation of $\widehat{\gamma}$
will vanish at the origin, this is the only way in which one can cover
the origin of the frequency domain. In most cases, the exact shape
of the low-pass part is not important, so that taking a modulation
of the mother wavelet is acceptable. But in other cases, one might
desire more specific properties of the low-pass part; for example,
one could want it to be real-valued. In this case, the added flexibility
of allowing $\gamma$ to depend on $i\in I$ might be valuable. In
this introduction, however, we will only consider the case in which
$\gamma_{i}=\gamma$ is independent of $i\in I$, for the sake of
simplicity.

Now, since the family $\left(\widehat{\gamma^{\left(i\right)}}\right)_{i\in I}$
behaves similarly to the family $\left(\varphi_{i}\right)_{i\in I}$
(at least with respect to the (essential) frequency support), one
could be tempted to conjecture that
\begin{equation}
\left\Vert g\right\Vert _{\DecompSp{\CalQ}p{\ell_{w}^{q}}{}}\asymp\left\Vert \left(\left\Vert \smash{\gamma^{\left(i\right)}}\ast g\right\Vert _{L^{p}}\right)_{i\in I}\right\Vert _{\ell_{w}^{q}}.\label{eq:IntroductionSemiDiscreteBanachFrame}
\end{equation}
For the special case of $\alpha$-modulation spaces, this statement
was established (for (almost) arbitrary $\gamma\in\Schwartz\left(\R^{\dimension}\right)$)
in \cite{EmbeddingsOfAlphaModulationIntoSobolev}. Our first result
(cf.\@ Section \ref{sec:SemiDiscreteBanachFrames}) will be to show
that for $p\in\left[1,\infty\right]$, equation (\ref{eq:IntroductionSemiDiscreteBanachFrame})
is indeed valid under suitable assumptions on $\gamma$. Furthermore,
for $p\in\left(0,1\right)$, we have the slightly modified statement
\begin{equation}
\left\Vert g\right\Vert _{\DecompSp{\CalQ}p{\ell_{w}^{q}}{}}\asymp\left\Vert \left(\left\Vert \smash{\gamma^{\left(i\right)}}\ast g\right\Vert _{W_{T_{i}^{-T}\left[-1,1\right]^{\dimension}}\left(L^{p}\right)}\right)_{i\in I}\right\Vert _{\ell_{w}^{q}},\label{eq:IntroductionSemiDiscreteQuasiBanachFrame}
\end{equation}
where $W_{T_{i}^{-T}\left[-1,1\right]^{\dimension}}\left(L^{p}\right)$
is a so-called \textbf{Wiener amalgam space} (originally introduced
by Feichtinger\cite{FeichtingerWienerSpaces}). We refer to the results
(\ref{eq:IntroductionSemiDiscreteBanachFrame})-(\ref{eq:IntroductionSemiDiscreteQuasiBanachFrame})
as stating that the family $\left(\gamma^{\left(i\right)}\right)_{i\in I}$
forms a \textbf{\emph{semi-discrete}}\textbf{ Banach frame} for $\DecompSp{\CalQ}p{\ell_{w}^{q}}{}$.
The reason for this nomenclature is that the index set of the family
$\left(\left(\gamma^{\left(i\right)}\ast g\right)\left(x\right)\right)_{i\in I,x\in\R^{\dimension}}$
has the discrete part $I$, but also the continuous part $\R^{\dimension}$.

Our next results are concerned with a further discretization of this
semi-discrete Banach frame. Indeed, under more stringent assumptions
on $\gamma$, we show in Section \ref{sec:FullyDiscreteBanachFrames}
for $\delta>0$ sufficiently small that the \textbf{structured generalized
shift-invariant system}
\begin{equation}
\Psi_{\delta}:=\left(L_{\delta\cdot T_{i}^{-T}k}\:\widetilde{\gamma^{\left[i\right]}}\right)_{i\in I,k\in\Z^{\dimension}}\qquad\text{ with }\qquad\tilde{f}\left(x\right):=f\left(-x\right)\label{eq:IntroductionBanachFrameFamily}
\end{equation}
generates a Banach frame for $\DecompSp{\CalQ}p{\ell_{w}^{q}}{}$,
with the associated discrete sequence space
\begin{equation}
Y\!:=\ell_{\left(\left|\det T_{i}\right|^{\frac{1}{2}-\frac{1}{p}}\cdot w_{i}\!\right)_{\!i\in I}}^{q}\!\!\!\!\!\!\!\!\left(\!\left[\ell^{p}\!\left(\Z^{\dimension}\right)\right]_{i\in I}\right)\quad\text{ where }\quad\left\Vert \left(\smash{c_{k}^{\left(i\right)}}\right)_{i\in I,k\in\Z^{\dimension}}\right\Vert _{Y}=\left\Vert \!\left(\left|\det T_{i}\right|^{\frac{1}{2}-\frac{1}{p}}\!\cdot\!w_{i}\!\cdot\!\left\Vert \left(\smash{c_{k}^{\left(i\right)}}\right)_{k\in\Z^{\dimension}}\right\Vert _{\ell^{p}}\right)_{\!i\in I}\right\Vert _{\ell^{q}}.\label{eq:IntroductionSequenceSpaceDefinition}
\end{equation}
Since the system $\Psi_{\delta}$ is generated in a very structured
way—similar to the usual definition of Gabor, wavelet or shearlet
frames—from a \emph{single} prototype function, we call $\Psi_{\delta}$
a \textbf{structured Banach frame} for $\DecompSp{\CalQ}p{\ell_{w}^{q}}{}$.

Finally, we show in Section \ref{sec:AtomicDecompositions}—again
under slightly different assumptions on $\gamma$—that the family
$\left(L_{\delta\cdot T_{i}^{-T}k}\:\gamma^{\left[i\right]}\right)_{\!i\in I,k\in\Z^{\dimension}}$
forms an \textbf{atomic decomposition} for $\DecompSp{\CalQ}p{\ell_{w}^{q}}{}$,
with the same associated sequence space $Y$ as above. As above, we
call this family a \textbf{structured atomic decomposition}. Hence,
at least if $\gamma$ is symmetric and fulfills certain technical
conditions, the family $\Psi_{\delta}$ will \emph{simultaneously}
form a Banach frame, as well as an atomic decomposition for $\DecompSp{\CalQ}p{\ell_{w}^{q}}{}$;
in particular, this implies that \emph{analysis sparsity is equivalent
to synthesis sparsity} for $\Psi_{\delta}$.

We remark that the assumptions placed on the prototype function $\gamma$
are quite technical, even though we will achieve a significant simplification
of these conditions in Section \ref{sec:SimplifiedCriteria}. Indeed,
a slightly simplified version of our main theorem concerning Banach
frames reads as follows:
\begin{thm*}
(cf.\@ Corollary \ref{cor:BanachFrameSimplifiedCriteria} for the
precise statement)

Recall from equation (\ref{eq:IntroductionStructuredCoveringAssumption})
that $\CalQ=\left(T_{i}Q+b_{i}\right)_{i\in I}$. Assume that there
is an open set $P\subset\R^{\dimension}$ with $\overline{P}\subset Q$
and $\CalO=\bigcup_{i\in I}T_{i}P+b_{i}$. Let $p,q\in\left(0,\infty\right]$.
Then there are explicitly given $N\in\N$ and $\sigma,\tau>0$, depending
on $\dimension,p,q$, with the following property:

If $w=\left(w_{i}\right)_{i\in I}$ is a $\CalQ$-moderate\footnote{cf.\@ Section \ref{subsec:DecompSpaceDefinitionStandingAssumptions},
equation (\ref{eq:IntroductionModerateWeightDefinition}) for the
precise definition.} weight and if $\gamma\in L^{1}\left(\R^{\dimension}\right)$ satisfies
the following:

\begin{enumerate}
\item We have $\widehat{\gamma}\in C^{\infty}\left(\R^{\dimension}\right)$,
where all partial derivatives of $\widehat{\gamma}$ are polynomially
bounded.
\item We have $\gamma\in C^{1}\left(\R^{\dimension}\right)$ and $\partial_{\ell}\gamma\in L^{1}\left(\R^{\dimension}\right)\cap L^{\infty}\left(\R^{\dimension}\right)$
for all $\ell\in\left\{ 1,\dots,\dimension\right\} $.
\item We have $\widehat{\gamma}\left(\xi\right)\neq0$ for all $\xi\in\overline{Q}$.
\item We have
\begin{equation}
C_{1}:=\sup_{i\in I}\,\sum_{j\in I}M_{j,i}<\infty\qquad\text{ and }\qquad C_{2}:=\sup_{j\in I}\sum_{i\in I}M_{j,i}<\infty\label{eq:IntroductionSimplifiedCondition}
\end{equation}
with
\[
\qquad M_{j,i}:=\left(\frac{w_{j}}{w_{i}}\right)^{\tau}\cdot\left(1+\left\Vert T_{j}^{-1}T_{i}\right\Vert \right)^{\sigma}\cdot\max_{\left|\beta\right|\leq1}\left(\left|\det T_{i}\right|^{-1}\cdot\int_{Q_{i}}\max_{\left|\alpha\right|\leq N}\left|\left(\partial^{\alpha}\widehat{\partial^{\beta}\gamma}\right)\!\!\left(T_{j}^{-1}\left(\xi-b_{j}\right)\right)\right|\d\xi\right)^{\tau}.
\]
\end{enumerate}
Then, for $\delta>0$ sufficiently small, the family $\Psi_{\delta}$
from equation (\ref{eq:IntroductionBanachFrameFamily}) is a Banach
frame for $\DecompSp{\CalQ}p{\ell_{w}^{q}}{}$, with associated sequence
space $Y$ as given in equation (\ref{eq:IntroductionSequenceSpaceDefinition}).
\end{thm*}
The conditions which ensure that $\gamma$ generates an atomic decomposition
are similar, but slightly more complicated, cf.\@ Corollary \ref{cor:AtomicDecompositionSimplifiedCriteria}.
For the sake of brevity, we omit them in this introduction.

Of course, condition (\ref{eq:IntroductionSimplifiedCondition}) is
quite technical. The main reason for this is that hugely different
coverings $\CalQ$ are treated using the same theory. Thus, given
a specific covering (e.g.\@ the ones used to define Besov spaces
or $\alpha$-modulation spaces), the difficulty consists in reducing
the general, abstract criteria provided by the theory to readily verifiable
criteria involving only the smoothness, decay and Fourier decay of
$\gamma$. As we will see in Sections \ref{sec:CompactlySupportedAlphaModulationFrames}
and \ref{sec:BesovFrames}, this is indeed possible for Besov spaces
and $\alpha$-modulation spaces. In addition, in the companion paper
\cite{StructuredBanachFrames2} we will show that the theory also
applies to shearlet smoothness spaces. Furthermore, it turns out that
in each of these cases one can find \emph{compactly supported} prototype
functions $\gamma$ which fulfill the relevant criteria. Thus, although
the decomposition spaces $\DecompSp{\CalQ}p{\ell_{w}^{q}}{}$ are
defined using the \emph{bandlimited} partition of unity $\left(\varphi_{i}\right)_{i\in I}$,
it is usually possible to give alternative characterizations in terms
of \emph{compactly supported} functions.

\medskip{}

At a first glance, the difficulty pertaining to the technical conditions
described above seems be a major drawback of the theory presented
here in comparison to coorbit theory. But in fact, coorbit theory
faces the same problem: In essentially every example where coorbit
theory is applicable, one has a systematic way of assigning to each
``prototype'' $\psi$ a whole family $\Psi=\left(\psi_{x}\right)_{x\in X}$.
Then, one has to obtain a profound understanding of the mapping $\psi\mapsto\left(\psi_{x}\right)_{x\in X}$
in order to derive readily verifiable conditions on $\psi$ which
ensure that the family $\Psi$ is suitable for the application of
coorbit theory, in particular to ensure that $\Psi$ is a (Hilbert
space) frame and that the kernel $R\left(x,y\right)=\left\langle \psi_{y},\,S^{-1}\psi_{x}\right\rangle $
belongs to $\mathcal{A}_{m}$. As examples for the effort one still
has to put in to apply coorbit theory in specific situations, we mention
\cite{FuehrCoorbit1,FuehrCoorbit2,FuehrGeneralizedCalderonConditions,FuehrSimplifiedVanishingMomentCriteria,GeneralizedCoorbit2,Dahlke_etal_sh_coorbit1,Dahlke_etal_sh_coorbit2,DahlkeShearletArbitraryDimension,DahlkeToeplitzShearletTransform,FuehrVoigtlaenderCoorbitSpacesAsDecompositionSpaces,FuehrWaveletFramesAndAdmissibility,UllrichContinuousCharacterizationsOfBesovTriebelLizorkin}.
We emphasize that this effort should not be seen as a shortcoming
of the mentioned papers or of (generalized) coorbit theory, but rather
as showing that despite the tremendous simplifications coorbit theory
has to offer, one still has to put in work to apply it in concrete
situations. The same is true of the results in this paper.

In fact, there is an intimate connection between the decomposition
space setting considered here and the coorbit setting considered in
\cite{Dahlke_etal_sh_coorbit1,DahlkeShearletArbitraryDimension,FuehrCoorbit2,FuehrSimplifiedVanishingMomentCriteria,UllrichContinuousCharacterizationsOfBesovTriebelLizorkin}:
In all of these papers, the authors consider coorbit spaces of a semi-direct
product $\R^{\dimension}\rtimes H$ for suitable \textbf{dilation
groups} $H\leq\GL\left(\R^{\dimension}\right)$, where the associated
unitary representation $\pi:\R^{\dimension}\rtimes H\to\mathcal{U}\left(L^{2}\left(\R^{\dimension}\right)\right),\left(x,h\right)\mapsto L_{x}D_{h}$
is the \textbf{quasi-regular representation}, i.e., the natural action
of $\R^{\dimension}\rtimes H$ on $L^{2}\left(\R^{\dimension}\right)$
in terms of the translations $L_{x}$ and the dilations $D_{h}$ with
$D_{h}f=\left|\det h\right|^{-1/2}\cdot\left(f\circ h^{-1}\right)$.
The mentioned papers contain—typically somewhat technical and lengthy—sufficient
criteria which ensure that a given \emph{mother wavelet} can serve
as an atom in the coorbit scheme. These conditions heavily depend
on the considered dilation group $H$ and also on the weight $w:\R^{\dimension}\rtimes H\to\left(0,\infty\right)$
which is used for the weighted mixed Lebesgue space $L_{w}^{p,q}\left(\R^{\dimension}\rtimes H\right)$.
For a given mother wavelet $g$ satisfying these criteria, the theory
of coorbit spaces implies that each sufficiently densely sampled family
$\left(\pi\left(x_{j},h_{j}\right)g\right)_{j\in J}$ yields an atomic
decomposition, as well as a Banach frame for the coorbit space ${\rm Co}\left(L_{w}^{p,q}\left(\R^{\dimension}\rtimes H\right)\right)$.

But as shown in \cite{FuehrVoigtlaenderCoorbitSpacesAsDecompositionSpaces}
and in \cite[Section 4]{VoigtlaenderPhDThesis}, we have ${\rm Co}\left(L_{w}^{p,q}\left(\R^{\dimension}\rtimes H\right)\right)=\DecompSp{\CalQ_{H}}p{\ell_{\tilde{w}}^{q}}{}$
up to canonical identifications, at least if the weight $w=w\left(x,h\right)$
only depends on the second factor, i.e., if $w=w\left(h\right)$.
Here, the so-called \textbf{induced covering} $\CalQ_{H}=\left(h_{i}^{-T}Q\right)_{i\in I}$
of the \textbf{dual orbit} $\CalO=H^{T}\xi_{0}\subset\R^{\dimension}$
is determined by an arbitrary well-spread family $\left(h_{i}\right)_{i\in I}$
in $H$. Given this identification, one can then apply the theory
developed in this paper to derive conditions on the prototype $\gamma$
which ensure that the family 
\[
\left(L_{\delta\cdot h_{i}k}\:\gamma^{\left[i\right]}\right)_{i\in I,k\in\Z^{\dimension}}=\left(\left|\det h_{i}\right|^{-1/2}\cdot L_{\delta\cdot h_{i}k}\left[\gamma\circ h_{i}^{-1}\right]\right)_{i\in I,k\in\Z^{\dimension}}=\left(\pi\left(\delta h_{i}k,\,h_{i}\right)\gamma\right)_{i\in I,k\in\Z^{\dimension}}
\]
forms a Banach frame, or an atomic decomposition for the decomposition
space $\DecompSp{\CalQ_{H}}p{\ell_{\tilde{w}}^{q}}{}$ and thus also
for the coorbit space ${\rm Co}\left(L_{w}^{p,q}\left(\R^{\dimension}\rtimes H\right)\right)$.
Note that the family $\left[\left(\delta h_{i}k,\,h_{i}\right)\right]_{i\in I,k\in\Z^{\dimension}}$
is well-spread in $\R^{\dimension}\rtimes H$ since $\left(h_{i}\right)_{i\in I}$
is well-spread in $H$. Hence, the theory developed in this paper
yields Banach frames and atomic decompositions which are \emph{of
the same form} as those obtained via coorbit theory. As future work,
we plan a systematic comparison of the conditions imposed on the prototype
$\gamma$ by coorbit theory (as in \cite{Dahlke_etal_sh_coorbit1,DahlkeShearletArbitraryDimension,FuehrCoorbit2,FuehrSimplifiedVanishingMomentCriteria,UllrichContinuousCharacterizationsOfBesovTriebelLizorkin})
on the one hand and by the theory developed in this paper on the other
hand.

\medskip{}

In spite of the strong connection between coorbit theory and the theory
developed in this paper, they differ in some important aspects:

As a first difference, we observe that coorbit theory requires to
pass from the given continuous frame $\left(\psi_{x}\right)_{x\in X}$
to a sufficiently densely sampled version $\left(\psi_{x_{i}}\right)_{i\in I}$.
This will usually not only require a sufficiently dense sampling \emph{in
the space domain} (which corresponds to $\delta$ in equation (\ref{eq:IntroductionBanachFrameFamily})),
but also to a rather dense sampling \emph{in the frequency domain}.
In contrast, for our approach only the sampling density \emph{in space}
needs to be sufficiently high. The ``frequency sampling density''
is fixed a priori by choosing the covering $\CalQ=\left(T_{i}Q+b_{i}\right)_{i\in I}$.

Next, the main advantage of our approach in comparison to coorbit
theory is that one does \emph{not} need to start from a given \emph{continuous}
frame $\left(\psi_{x}\right)_{x\in X}$ which is then discretized.
In fact, one can even start from a given discrete frame which is of
the form (\ref{eq:IntroductionBanachFrameFamily}). As long as the
family $\CalQ=\left(T_{i}Q+b_{i}\right)_{i\in I}$ forms a suitable
covering, one can then consider the associated decomposition spaces
$\DecompSp{\CalQ}p{\ell_{w}^{q}}{}$ and use the theory presented
here to justify that the discrete frame one started with forms a Banach
frame and an atomic decomposition for $\DecompSp{\CalQ}p{\ell_{w}^{q}}{}$,
possibly after adjusting the sampling density.

Probably, this intuition is what originally lead Labate et al.\@
to the introduction of the \textbf{shearlet smoothness spaces}\cite{Labate_et_al_Shearlet},
although they did not have the machinery to rigorously prove that
the usual discrete, cone-adapted shearlet systems indeed yield Banach
frames and atomic decompositions for the shearlet smoothness spaces.
Using the theory developed here, we will see in the companion paper
\cite{StructuredBanachFrames2} that this is indeed the case. Furthermore,
we will employ our results to show that suitable discrete, cone-adapted
shearlet systems achieve an almost optimal approximation rate for
the class of \textbf{cartoon-like functions}. At a first glance, this
might appear to be a well-known statement, but a closer inspection
of the classical results about approximation of cartoon-like functions
by shearlets (see e.g.\@ \cite{ShearletsAndOptimallySparseApproximation,CompactlySupportedShearlets,CompactlySupportedShearletsAreOptimallySparse,OptimallySparseMultidimensionalRepresentationUsingShearlets})
reveals that these papers in fact only show that the $N$-term approximation
$f_{N}$ \emph{with respect to the dual frame} of the shearlet frame
satisfies the (almost optimal) rate $\left\Vert f-f_{N}\right\Vert _{L^{2}}\lesssim N^{-1}\cdot\left(\log N\right)^{\theta}$
for suitable $\theta>0$.

\subsubsection{Comparison to other constructions of Banach frame decompositions
of decomposition spaces}

One of the first general constructions of atomic decompositions for
decomposition spaces was given by Borup and Nielsen in \cite{BorupNielsenDecomposition}.
The main difference between their approach and ours is that our frame
elements $\gamma^{\left[i\right]}$ can be chosen to be compactly
supported, while Borup and Nielsen purely focus on bandlimited frame
elements.

There is also a more recent paper by Nielsen and Rasmussen \cite{CompactlySupportedFramesForDecompositionSpaces}
in which they construct compactly supported frames for certain decomposition
spaces. In comparison to that paper, our assumptions concerning the
covering $\CalQ$ are more general, while our conclusions are more
specific:
\begin{itemize}[leftmargin=0.5cm]
\item In \cite{CompactlySupportedFramesForDecompositionSpaces}, the authors
only consider coverings $\CalQ$ which are induced by considering
$\R^{\dimension}$ in a certain way as a space of homogeneous type:
More precisely, $\CalQ$ is assumed to satisfy $\CalQ=\left(Q_{k}\right)_{k\in\Z^{\dimension}}=\left(\mathcal{B}_{A}\left(\xi_{k},\varrho\cdot h\left(\xi_{k}\right)\right)\right)_{k\in\Z^{\dimension}}$,
where the balls $\mathcal{B}_{A}\left(\xi,r\right)=\left\{ \zeta\in\R^{\dimension}\with\left|\zeta-\xi\right|_{A}<r\right\} $
are defined using the quasi-metric $\left|\cdot\right|_{A}$, which
is induced in a certain way (cf.\@ \cite[Definition 2.1]{CompactlySupportedFramesForDecompositionSpaces})
by the one-parameter group of dilations $\left(\delta_{t}\right)_{t>0}$
where $\delta_{t}=\exp\left(A\cdot\ln t\right)$ for a fixed matrix
$A\in\R^{\dimension\times\dimension}$ with positive eigenvalues.
As shown between \cite[Lemma 2.6]{CompactlySupportedFramesForDecompositionSpaces}
and \cite[Definition 2.7]{CompactlySupportedFramesForDecompositionSpaces},
we have
\[
Q_{k}=\delta_{h\left(\xi_{k}\right)}\left[\mathcal{B}_{A}\left(0,\varrho\right)\right]+\xi_{k}\qquad\forall k\in\Z^{\dimension},
\]
so that all sets $Q_{k}$ of the covering $\CalQ$ are affine images
of a fixed set, where the linear parts of the affine maps are all
elements of the one-parameter family $\left(\delta_{t}\right)_{t>0}$.
Note with $\nu:={\rm trace}\,A>0$ that $\det\delta_{t}=t^{\nu}$
for all $t>0$, so that $\delta_{t}$ is uniquely determined by its
determinant. Since the covering used to define shearlet smoothness
spaces uses affine transformations for which many \emph{different}
linear parts have the \emph{same} determinant, this shows—or at least
very strongly indicates—that the covering used to define the shearlet
smoothness spaces does \emph{not} satisfy the assumptions imposed
in \cite{CompactlySupportedFramesForDecompositionSpaces}, while our
theory is able to handle these spaces.

Below, we will give another more rigorous argument which shows that
the theory developed in \cite{CompactlySupportedFramesForDecompositionSpaces}
does in fact \emph{neither} apply to the usual dyadic covering which
is used to define (inhomogeneous) Besov spaces, \emph{nor} to the
covering used to define shearlet smoothness spaces.
\item While each of the compactly supported Banach frames constructed in
\cite{CompactlySupportedFramesForDecompositionSpaces} is a union
of generalized shift invariant systems, it is \emph{not} true that
the frames are generated from a single prototype function in the same
structured way as in our paper. In contrast, the Banach frames constructed
in \cite{CompactlySupportedFramesForDecompositionSpaces} are of the
form
\[
\qquad\left(\psi_{k,n}\right)_{k,n\in\Z^{\dimension}}=\left(\left[h\left(\xi_{k}\right)\right]^{\nu/2}\cdot\tau_{k}\left(\delta_{h\left(\xi_{k}\right)}^{T}\bullet-\frac{\pi}{a}n\right)e^{i\left\langle \cdot,\xi_{k}\right\rangle }\right)_{k,n\in\Z^{\dimension}}\quad\text{ where }\quad\tau_{k}=\sum_{i=1}^{K}a_{i}^{\left(k\right)}g_{m}\left(\bullet+\smash{b_{i}^{\left(k\right)}}\right),
\]
for suitable $K,m\in\N$ and with $g_{m}=C_{g}m^{\nu}\cdot g\circ\delta_{m}^{T}$.
Hence, using notation as in eq.\@ (\ref{eq:IntroductionL2NormalizedFrameElement})
with the covering $\CalQ$ as defined above and with $b_{k}:=\xi_{k}$,
as well as $T_{k}:=\delta_{h\left(\xi_{k}\right)}$, we have $\left(\psi_{k,n}\right)_{k,n\in\Z^{\dimension}}=\left(L_{\frac{\pi}{a}T_{k}^{-T}n}\:\tau_{k}^{\left[k\right]}\right)_{k,n\in\Z^{\dimension}}$,
while the structured family $\Psi_{\delta}$ defined in equation (\ref{eq:IntroductionBanachFrameFamily})
satisfies $\Psi_{\delta}=\left(\!L_{\delta\cdot T_{k}^{-T}n}\:\widetilde{\psi^{\left[k\right]}}\right)_{\!k,n\in\Z^{\dimension}}$.
In other words, while the structured Banach frames constructed in
this paper arise from a \emph{single} prototype function by translations,
modulations and dilations, the frames constructed in \cite{CompactlySupportedFramesForDecompositionSpaces}
do \emph{not} satisfy this property.

In particular, if the covering $\CalQ$ is the usual dyadic covering
of $\R^{\dimension}$ used to define (inhomogeneous) Besov spaces,
then $\Psi_{\delta}$ will be an (inhomogeneous) wavelet frame, while
this is \emph{not} in general true of the frame constructed in \cite{CompactlySupportedFramesForDecompositionSpaces}.
Additionally, the results in \cite{CompactlySupportedFramesForDecompositionSpaces}
are not applicable in this setting, as we will see below.
\end{itemize}
This last defect—that the resulting Banach frame is not generated
from a \emph{single} prototype—is addressed in the follow-up paper
\cite{NielsenSinglyGeneratedFrames}. There, Morten Nielsen considers
the same general setting as described above. He then constructs a
\emph{bandlimited} Banach frame for the associated decomposition spaces
which is generated by a \emph{single} prototype function in the same
structured way as proposed in the present paper. Furthermore, Nielsen
then uses a distortion argument to show that one can also obtain a
structured Banach frame with a \emph{single, compactly supported}
generator. Hence, at a first glance, it might seem that all results
of the present paper are already contained in \cite{NielsenSinglyGeneratedFrames}.
This, however, is \emph{not} true for the following reasons:
\begin{itemize}[leftmargin=0.5cm]
\item As already observed above, the coverings considered in \cite{NielsenSinglyGeneratedFrames}
and \cite{CompactlySupportedFramesForDecompositionSpaces} are quite
restricted. They have to be of the form $\CalQ=\left(Q_{k}\right)_{k\in\Z^{\dimension}}=\left(\mathcal{B}_{A}\left(\xi_{k},\varrho\cdot h\left(\xi_{k}\right)\right)\right)_{k\in\Z^{\dimension}}$,
where the balls $\mathcal{B}_{A}\left(\xi,r\right)=\left\{ \zeta\in\R^{\dimension}\with\left|\zeta-\xi\right|_{A}<r\right\} $
are defined using the quasi-metric $\left|\cdot\right|_{A}$, which
is determined by a suitable matrix $A$.

In particular, the setting considered in \cite{NielsenSinglyGeneratedFrames}
does \emph{neither} include the case of homogeneous or inhomogeneous
Besov spaces, nor the case of shearlet smoothness spaces. To see this,
note that \cite[Proposition 3.6]{NielsenSinglyGeneratedFrames} does
not impose any vanishing moment conditions on the prototype $\gamma$
(which is called $g$ in the notation of \cite{NielsenSinglyGeneratedFrames}).
In fact, it is even required that $\widehat{\gamma}\left(0\right)\neq0$.
But it is folklore that the generator of an (inhomogeneous or homogeneous)
wavelet frame for $L^{2}\left(\R\right)$ has to satisfy certain vanishing
moment conditions; the proof for homogeneous wavelet frames is given
in \cite[Theorem 3.3.1]{DaubechiesTenLecturesOnWavelets}. A proof
of the corresponding statement for discrete cone-adapted shearlet
frames is given in Appendix \ref{sec:ShearletFrameVanishingMomentNecessity}.

In stark contrast, the theory developed in the present paper \emph{is}
able to handle Besov spaces (cf.\@ Section \ref{sec:BesovFrames}),
as well as shearlet smoothness spaces (cf.\@ the companion paper
\cite{StructuredBanachFrames2}).
\item Since a distortion argument is used to obtain a compactly supported
Banach frame from a bandlimited frame, the choice of the generator
$\gamma$ in \cite{NielsenSinglyGeneratedFrames} is quite restricted;
$\gamma$ has to be close enough to the generator of the bandlimited
frame.

In contrast, the assumptions imposed on $\gamma$ in the present paper
are quite mild. In most concrete cases (in particular for $\alpha$-modulation
spaces, Besov spaces and shearlet smoothness spaces), the conditions
reduce to suitable smoothness, decay and vanishing moment criteria,
in conjunction with a certain nonvanishing condition for the Fourier
transform $\widehat{\gamma}$.
\item In the present paper, we also consider the decomposition spaces $\DecompSp{\CalQ}p{\ell_{w}^{q}}v$
where a \emph{weighted} Lebesgue space $L_{v}^{p}\left(\R^{\dimension}\right)$
is used. In contrast, \cite{NielsenSinglyGeneratedFrames} only considers
the \emph{unweighted} case.
\end{itemize}
We remark however that \cite{NielsenSinglyGeneratedFrames} jointly
considers Triebel-Lizorkin type, as well as Besov type decomposition
spaces. In contrast, at least in its present state, the approach developed
in this paper only applies to the Besov type decomposition spaces.

Finally, we mention the recent paper \cite{CharacterizationOfSparseNonstationaryGaborExpansions}
in which Ottosen and Nielsen take the ``reverse'' of the usual approach:
Instead of starting with a given function space $X$ and then constructing
Banach frames or atomic decompositions for this space, the authors
start with a given \textbf{painless nonstationary Gabor frame} $\left(h_{i,k}\right)_{i,k\in\Z^{\dimension}}$
satisfying 
\[
h_{i,k}=L_{a_{i}\cdot k}\,h_{i}\quad\text{ and }\quad\supp\widehat{h_{i}}\subset\left[0,\,a_{i}^{-1}\right]^{\dimension}+b_{i}\quad\text{ for certain }\quad a_{i}>0\text{ and }b_{i}\in\R^{\dimension}.
\]
Under suitable assumptions on the (slightly enlarged) covering 
\[
\CalQ=\left(Q_{i}\right)_{i\in\Z^{\dimension}}\qquad\text{ with }\qquad Q_{i}=a_{i}^{-1}\cdot\left(-\delta,\,1+\delta\right)^{\dimension}+b_{i},
\]
Ottosen and Nielsen then show that the renormalized family $\left(h_{i,k}^{\left(p\right)}\right)_{i,k\in\Z^{\dimension}}$
defined by $h_{i,k}^{\left(p\right)}=a_{i}^{\frac{1}{p}-\frac{1}{2}}\cdot h_{i,k}$
forms a Banach frame for the decomposition space $\DecompSp{\CalQ}p{\ell_{\omega^{s}}^{q}}{}$,
where $\omega_{i}=1+\left\Vert \xi_{i}\right\Vert ^{2}$ for suitable
$\xi_{i}\in Q_{i}$. In addition, it is shown in \cite[Theorem 6.1]{CharacterizationOfSparseNonstationaryGaborExpansions}
for $p,q\in\left(0,\infty\right)$ that every $f\in\DecompSp{\CalQ}p{\ell_{\omega^{s}}^{q}}{}$
admits an expansion of the form
\begin{equation}
f=\sum_{i,k\in\Z^{\dimension}}\left\langle h,\,h_{i,k}\right\rangle \cdot\tilde{h}_{i,k},\label{eq:OttosenNielsenBanachFrameExpansion}
\end{equation}
where $\left(\smash{\tilde{h}_{i,k}}\right)_{i,k\in\Z^{\dimension}}$
is the canonical dual frame of the nonstationary Gabor frame $\left(h_{i,k}\right)_{i,k\in\Z^{\dimension}}$.

In summary, the paper \cite{CharacterizationOfSparseNonstationaryGaborExpansions}
starts with a given painless nonstationary Gabor frame and then shows
that the space of analysis-sparse signals w.r.t.\@ the frame coincides
with a suitably defined decomposition space. Note that the painless
nonstationary Gabor frames are always bandlimited. Using the theory
developed in this paper, it should be possible (perhaps with the cost
of changing the sampling density in comparison to the original frame)
to show similar results for nonstationary Gabor frames with \emph{compactly
supported} generators. Furthermore, while the results in \cite{CharacterizationOfSparseNonstationaryGaborExpansions}
only show that each $f\in\DecompSp{\CalQ}p{\ell_{w}^{q}}{}$ admits
a sparse expansion in terms of the \emph{dual frame} $\left(\smash{\tilde{h}_{i,k}}\right)_{i,k\in\Z^{\dimension}}$,
our results would yield a sparse expansion in terms of the frame itself,
so that analysis sparsity is equivalent to synthesis sparsity.

\subsection{Notation and conventions}

\label{subsec:Notation}We write $\N=\Z_{\geq1}$ for the set of \textbf{natural
numbers} and $\N_{0}=\Z_{\geq0}$ for the set of natural numbers including
$0$. For a matrix $A\in\Compl^{\dimension\times\dimension}$, we
denote by $A^{T}$ the transpose of $A$. The norm $\left\Vert A\right\Vert $
of $A$ is the usual \textbf{operator norm} of $A$, acting on $\R^{\dimension}$
equipped with the usual euclidean norm $\left|\cdot\right|=\left\Vert \cdot\right\Vert _{2}$.
The \textbf{open euclidean ball} of radius $r>0$ around $x\in\R^{\dimension}$
is denoted by $B_{r}\left(x\right)$. For a linear (bounded) operator
$T:X\to Y$ between (quasi)-normed spaces $X,Y$, we denote the \textbf{operator
norm} of $T$ by 
\[
\vertiii T:=\vertiii T_{X\to Y}:=\sup_{\left\Vert x\right\Vert _{X}\leq1}\left\Vert Tx\right\Vert _{Y}.
\]

For an arbitrary set $M$, we let $\left|M\right|\in\N_{0}\cup\left\{ \infty\right\} $
denote the number of elements of the set. For $n\in\N_{0}=\Z_{\geq0}$,
we write $\underline{n}:=\left\{ 1,\dots,n\right\} $; in particular,
$\underline{0}=\emptyset$. For the \textbf{closure} of a subset $M$
of some topological space, we write $\overline{M}$.

The $\dimension$-dimensional \textbf{Lebesgue measure} of a (measurable)
set $M\subset\R^{\dimension}$ is denoted by $\lambda\left(M\right)$
or by $\lambda_{\dimension}\left(M\right)$. Furthermore, for $M\subset\R^{\dimension}$,
we define the \textbf{indicator function} (or \textbf{characteristic
function}) $\Indicator_{M}$ of the set $M$ by
\[
\Indicator_{M}:\R^{\dimension}\to\left\{ 0,1\right\} ,x\mapsto\begin{cases}
1, & \text{if }x\in M,\\
0, & \text{otherwise}.
\end{cases}
\]
For two subsets $A,B\subset\R^{\dimension}$, we define the \textbf{Minkowski
sum} and the \textbf{Minkowski difference} of $A,B$ by
\[
A+B:=\left\{ a+b\with a\in A,\,b\in B\right\} \qquad\text{ and }\qquad A-B:=\left\{ a-b\with a\in A,\,b\in B\right\} .
\]
The Minkowski difference $A-B$ should be distinguished from the \textbf{set-theoretic
difference} $A\setminus B=\left\{ a\in A\with a\notin B\right\} $.

The \textbf{translation} and \textbf{modulation} of a function $f:\R^{\dimension}\to\Compl^{k}$
by $x\in\R^{\dimension}$ or $\xi\in\R^{\dimension}$ are, respectively,
denoted by 
\[
L_{x}f:\R^{\dimension}\to\Compl^{k},y\mapsto f\left(y-x\right),\qquad\text{ and }\qquad M_{\xi}f:\R^{\dimension}\to\Compl^{k},y\mapsto e^{2\pi i\left\langle \xi,y\right\rangle }f\left(y\right).
\]

For the \textbf{Fourier transform}, we use the convention $\widehat{f}\left(\xi\right):=\left(\Fourier f\right)\left(\xi\right):=\int_{\R^{\dimension}}f\left(x\right)\cdot e^{-2\pi i\left\langle x,\xi\right\rangle }\d x$
for $f\in L^{1}\left(\R^{\dimension}\right)$. It is well-known that
the Fourier transform extends to a unitary automorphism $\Fourier:L^{2}\left(\R^{\dimension}\right)\to L^{2}\left(\R^{\dimension}\right)$.
The inverse of this map is the continuous extension of the inverse
Fourier transform, given by $\left(\Fourier^{-1}f\right)\left(x\right)=\int_{\R^{\dimension}}f\left(\xi\right)e^{2\pi i\left\langle x,\xi\right\rangle }\d\xi$
for $f\in L^{1}\left(\R^{\dimension}\right)$. We will make frequent
use of the space $\Schwartz\left(\R^{\dimension}\right)$ of \textbf{Schwartz
functions} and its dual space $\Schwartz'\left(\R^{\dimension}\right)$,
the space of \textbf{tempered distributions}. For more details on
these spaces, we refer to \cite[Section 9]{FollandRA}; in particular,
we note that the Fourier transform restricts to a linear homeomorphism
$\Fourier:\Schwartz\left(\R^{\dimension}\right)\to\Schwartz\left(\R^{\dimension}\right)$;
by duality, we can thus define $\Fourier:\Schwartz'\left(\R^{\dimension}\right)\to\Schwartz'\left(\R^{\dimension}\right)$
by $\Fourier\varphi=\varphi\circ\Fourier$ for $\varphi\in\Schwartz'\left(\R^{\dimension}\right)$.

Given an open subset $U\subset\R^{\dimension}$, we let $\DistributionSpace U$
denote the space of \textbf{distributions} on $U$, i.e., the topological
dual space of $\TestFunctionSpace U$. For the precise definition
of the topology on $\TestFunctionSpace U$, we refer to \cite[Chapter 6]{RudinFA}.
We remark that the dual pairings $\left\langle \cdot,\cdot\right\rangle _{\CalD',\CalD}$
and $\left\langle \cdot,\cdot\right\rangle _{\Schwartz',\Schwartz}$
are always taken to be \emph{bilinear} instead of sesquilinear.

We write $v_{\dimension}:=\lambda_{\dimension}\left(B_{1}\left(0\right)\right)$
for the $\dimension$-dimensional Lebesgue measure of the euclidean
unit ball. An easy, but sometimes useful estimate is that $v_{\dimension}\leq2^{\dimension}$,
since $B_{1}\left(0\right)\subset\left[-1,1\right]^{\dimension}$.
Furthermore, we let $s_{\dimension}:=\mathcal{H}_{\dimension-1}\left(S^{\dimension-1}\right)$
denote the surface measure of the unit sphere. It is well-known that
$s_{\dimension}=\dimension\cdot v_{\dimension}\leq\dimension\cdot2^{\dimension}\leq2^{2\dimension}$,
since $\dimension\leq2^{\dimension}$. Finally, we have $B_{1}^{\left\Vert \cdot\right\Vert _{\infty}}\left(0\right)\subset B_{\sqrt{\dimension}}\left(0\right)$
and thus $2^{\dimension}=\lambda\left(B_{1}^{\left\Vert \cdot\right\Vert _{\infty}}\left(0\right)\right)\leq\lambda\left(B_{\sqrt{\dimension}}\left(0\right)\right)=v_{\dimension}\cdot\dimension^{\dimension/2}$,
which implies $v_{\dimension}\geq\left(2/\sqrt{\dimension}\right)^{\dimension}$.

The constant $s_{\dimension}$ be important for us due to the following:
For $p\in\left(0,\infty\right)$ and $N>\dimension/p$, we get using
polar coordinates that
\begin{align*}
\left\Vert \left(1+\left|\mybullet\right|\right)^{-N}\right\Vert _{L^{p}}^{p}=\int_{\R^{\dimension}}\left(1+\left|x\right|\right)^{-Np}\d x & =\int_{0}^{\infty}r^{\dimension-1}\int_{S^{\dimension-1}}\left(1+\left|r\xi\right|\right)^{-Np}\d\mathcal{H}_{\dimension-1}\left(\xi\right)\,\d r\\
 & =\mathcal{H}_{\dimension-1}\left(S^{\dimension-1}\right)\cdot\int_{0}^{\infty}r^{\dimension-1}\cdot\left(1+r\right)^{-Np}\d r\\
 & \leq s_{\dimension}\cdot\int_{0}^{\infty}\left(1+r\right)^{\dimension-Np-1}\d r\\
\left({\scriptstyle \text{since }\dimension-Np<0}\right) & =s_{\dimension}\cdot\frac{\left(1+r\right)^{\dimension-Np}}{\dimension-Np}\bigg|_{0}^{\infty}=\frac{s_{\dimension}}{Np-\dimension},
\end{align*}
and hence
\begin{equation}
\left\Vert \left(1-\left|\mybullet\right|\right)^{-N}\right\Vert _{L^{p}}\leq\left(\frac{1}{p}\cdot\frac{s_{\dimension}}{N-\frac{\dimension}{p}}\right)^{1/p}\qquad\forall N>\dimension/p,\label{eq:StandardDecayLpEstimate}
\end{equation}
which also remains valid (with the interpretation $x^{0}=1$ for arbitrary
$x\geq0$) for $p=\infty$.

\subsection{Definition of decomposition spaces and standing assumptions}

\label{subsec:DecompSpaceDefinitionStandingAssumptions}For the whole
paper, we fix a \textbf{semi-structured admissible covering} $\CalQ=\left(Q_{i}\right)_{i\in I}$
of an open subset $\CalO\subset\R^{\dimension}$. Precisely this means
that for each $i\in I$ there is a measurable subset $Q_{i}'\subset\R^{\dimension}$,
an invertible linear map $T_{i}\in\GL\left(\R^{\dimension}\right)$
and a translation $b_{i}\in\R^{\dimension}$ such that $Q_{i}=S_{i}Q_{i}'=T_{i}Q_{i}'+b_{i}$
for the affine transformation $S_{i}:\R^{\dimension}\to\R^{\dimension},\xi\mapsto T_{i}\xi+b_{i}$
and such that the following properties are fulfilled:
\begin{enumerate}
\item $\CalQ$ covers $\CalO$, i.e., $\CalO=\bigcup_{i\in I}Q_{i}$.
\item $\CalQ$ is \textbf{admissible}, i.e., we have $\left|i^{\ast}\right|\leq N_{\CalQ}<\infty$
for all $i\in I$, where
\begin{equation}
i^{\ast}:=\left\{ \ell\in I\with Q_{\ell}\cap Q_{i}\neq\emptyset\right\} .\label{eq:IndexClusterDefinition}
\end{equation}
\item There is some $R_{\CalQ}>0$ satisfying $Q_{i}'\subset\overline{B_{R_{\CalQ}}}\left(0\right)$
for all $i\in I$.
\item There is some $C_{\CalQ}>0$ satisfying $\left\Vert T_{i}^{-1}T_{\ell}\right\Vert \leq C_{\CalQ}$
for all $i\in I$ and all $\ell\in i^{\ast}$.
\end{enumerate}
The most common form of decomposition spaces uses a (quasi)-norm of
the form $\left\Vert \left(\left\Vert \Fourier^{-1}\left(\varphi_{i}\cdot\widehat{g}\right)\right\Vert _{L^{p}}\right)_{i\in I}\right\Vert _{\ell_{w}^{q}}$,
i.e., the frequency-localized pieces $g_{i}=\Fourier^{-1}\left(\varphi_{i}\cdot\widehat{g}\right)$
of $g$ are measured in $L^{p}\left(\R^{\dimension}\right)$. To achieve
even greater flexibility, we will allow weighted Lebesgue spaces of
the form $L_{v}^{p}\left(\R^{\dimension}\right)$ instead of $L^{p}\left(\R^{\dimension}\right)$.
Here, we write
\[
L_{v}^{p}\left(\smash{\R^{\dimension}}\right):=\left\{ f:\R^{\dimension}\to\Compl\with f\text{ measurable and }v\cdot f\in L^{p}\left(\smash{\R^{\dimension}}\right)\right\} ,
\]
equipped with the natural (quasi)-norm $\left\Vert f\right\Vert _{L_{v}^{p}}:=\left\Vert v\cdot f\right\Vert _{L^{p}}$.
In order to still obtain reasonable spaces and results, we assume
the following:
\begin{enumerate}[resume]
\item The weights $v,v_{0}:\R^{\dimension}\to\left(0,\infty\right)$ are
measurable and satisfy the following:

\begin{enumerate}
\item $v_{0}\geq1$ and\footnote{One can always assume $v_{0}\geq1$ without loss of generality, since
all properties of $v_{0}$ (including submultiplicativity) are also
fulfilled for $\widetilde{v_{0}}:=1+v_{0}\geq1$, where possibly $\Omega_{1}$
has to be enlarged, since $\widetilde{v_{0}}\left(x\right)\leq1+\Omega_{1}\cdot\left(1+\left|x\right|\right)^{K}\leq\left(1+\Omega_{1}\right)\cdot\left(1+\left|x\right|\right)^{K}$.
Likewise, by switching to $\tilde{v_{0}}\left(x\right):=v_{0}\left(x\right)+v_{0}\left(-x\right)$,
one can always assume $v_{0}$ to be symmetric.} $v_{0}$ is symmetric, i.e., $v_{0}\left(-x\right)=v_{0}\left(x\right)$
for all $x\in\R^{\dimension}$.
\item $v_{0}$ is submultiplicative, i.e., $v_{0}\left(x+y\right)\leq v_{0}\left(x\right)\cdot v_{0}\left(y\right)$
for all $x,y\in\R^{\dimension}$.
\item $v$ is $v_{0}$-moderate, i.e., $v\left(x+y\right)\leq v\left(x\right)\cdot v_{0}\left(y\right)$
for all $x,y\in\R^{\dimension}$.
\item There is some $K\geq0$ and some $\Omega_{1}\geq1$ satisfying $v_{0}\left(x\right)\leq\Omega_{1}\cdot\left(1+\left|x\right|\right)^{K}$
for all $x\in\R^{\dimension}$.
\item The constant $K$ from the previous step satisfies $K=0$ or there
is a constant $\Omega_{0}\geq1$ satisfying $\left\Vert T_{i}^{-1}\right\Vert \leq\Omega_{0}$
for all $i\in I$.
\end{enumerate}
\item There is a \textbf{$\CalQ$-$v_{0}$-BAPU} (bounded admissible partition
of unity) $\Phi=\left(\varphi_{i}\right)_{i\in I}$ for $\CalQ$,
which means that:

\begin{enumerate}
\item $\varphi_{i}\in\TestFunctionSpace{\CalO}$ for all $i\in I$ and furthermore
$\varphi_{i}\equiv0$ on $\CalO\setminus Q_{i}$.
\item $\sum_{i\in I}\varphi_{i}\equiv1$ on $\CalO$.
\item For each $p\in\left(0,\infty\right]$, the following expression (then
a constant) is finite:
\[
C_{\CalQ,\Phi,v_{0},p}:=\sup_{i\in I}\left[\left|\det T_{i}\right|^{\max\left\{ \frac{1}{p},1\right\} -1}\cdot\left\Vert \Fourier^{-1}\varphi_{i}\right\Vert _{L_{v_{0}}^{\min\left\{ 1,p\right\} }}\right].
\]
\end{enumerate}
\end{enumerate}
Clearly, if one chooses $K=0$ and $v=v_{0}\equiv1$, then one obtains
the usual decomposition spaces, as considered e.g.\@ in \cite{BorupNielsenDecomposition,VoigtlaenderPhDThesis,DecompositionEmbedding,DecompositionIntoSobolev}.
This will be the most common case. Note that in this case, we do \emph{not}
need to assume $\left\Vert T_{i}^{-1}\right\Vert \leq\Omega_{0}$
for all $i\in I$, i.e., the covering $\CalQ$ can be very general.

We observe for later use that the preceding assumptions imply
\begin{equation}
\left(1+\left|x\right|\right)^{K}\leq\Omega_{0}^{K}\cdot\left(1+\left|T_{i}^{T}x\right|\right)^{K}\qquad\forall x\in\R^{\dimension}.\label{eq:WeightLinearTransformationsConnection}
\end{equation}
Indeed, in case of $K=0$, this is trivial. In case of $K>0$, our
assumptions imply 
\[
\left|x\right|=\left|T_{i}^{-T}T_{i}^{T}x\right|\leq\left\Vert T_{i}^{-T}\right\Vert \cdot\left|T_{i}^{T}x\right|=\left\Vert T_{i}^{-1}\right\Vert \cdot\left|T_{i}^{T}x\right|\leq\Omega_{0}\cdot\left|T_{i}^{T}x\right|
\]
and hence $1+\left|x\right|\leq1+\Omega_{0}\cdot\left|T_{i}^{T}x\right|\leq\Omega_{0}\cdot\left(1+\left|T_{i}^{T}x\right|\right)$,
where the last step used that $\Omega_{0}\geq1$. This easily shows
that equation (\ref{eq:WeightLinearTransformationsConnection}) remains
valid also for $K>0$.

Finally, we observe for later use the convolution relation $L_{v_{0}}^{1}\left(\R^{\dimension}\right)\ast L_{v}^{p}\left(\R^{\dimension}\right)\hookrightarrow L_{v}^{p}\left(\R^{\dimension}\right)$
for $p\in\left[1,\infty\right]$. Indeed, we have
\begin{align*}
v\left(x\right)\cdot\left|\left(f\ast g\right)\left(x\right)\right| & \leq v\left(x\right)\cdot\int_{\R^{\dimension}}\left|f\left(y\right)\right|\cdot\left|g\left(x-y\right)\right|\d y\\
\left({\scriptstyle \text{since }v\left(x\right)=v\left(x-y+y\right)\leq v\left(x-y\right)\cdot v_{0}\left(y\right)}\right) & \leq\int_{\R^{\dimension}}\left|\left(v_{0}\cdot f\right)\left(y\right)\right|\cdot\left|\left(v\cdot g\right)\left(x-y\right)\right|\d y,
\end{align*}
so that Minkowski's inequality for integrals (cf.\@ \cite[Theorem (6.19)]{FollandRA}),
together with the isometric translation invariance of $L^{p}\left(\R^{\dimension}\right)$,
yields
\begin{align}
\left\Vert f\ast g\right\Vert _{L_{v}^{p}}=\left\Vert v\cdot\left(f\ast g\right)\right\Vert _{L^{p}} & \leq\left\Vert x\mapsto\int_{\R^{\dimension}}\left|\left(v_{0}\cdot f\right)\left(y\right)\right|\cdot\left|\left(v\cdot g\right)\left(x-y\right)\right|\d y\right\Vert _{L^{p}}\nonumber \\
 & \leq\int_{\R^{\dimension}}\left\Vert x\mapsto\left|\left(v_{0}\cdot f\right)\left(y\right)\right|\cdot\left|\left(v\cdot g\right)\left(x-y\right)\right|\right\Vert _{L^{p}}\d y\nonumber \\
 & =\int_{\R^{\dimension}}\left|\left(v_{0}\cdot f\right)\left(y\right)\right|\d y\cdot\left\Vert v\cdot g\right\Vert _{L^{p}}=\left\Vert f\right\Vert _{L_{v_{0}}^{1}}\cdot\left\Vert g\right\Vert _{L_{v}^{p}}<\infty.\label{eq:WeightedYoungInequality}
\end{align}
We will call this estimate the \textbf{weighted Young inequality}.
In particular, it shows that $\left(\left|f\right|\ast\left|g\right|\right)\left(x\right)<\infty$
for almost all $x\in\R^{\dimension}$.

\medskip{}

Given a $\CalQ$-$v_{0}$-BAPU $\Phi=\left(\varphi_{i}\right)_{i\in I}$,
we define the \textbf{clustered version} of $\Phi$ as $\Phi^{\ast}=\left(\varphi_{i}^{\ast}\right)_{i\in I}$,
where $\varphi_{i}^{\ast}:=\sum_{\ell\in i^{\ast}}\varphi_{\ell}$.
Because of $\sum_{i\in I}\varphi_{i}\equiv1$ on $\CalO\supset Q_{i}$
and since $\varphi_{\ell}\equiv0$ on $Q_{i}$ for all $\ell\in I\setminus i^{\ast}$,
it is not hard to see $\varphi_{i}^{\ast}\equiv1$ on $\overline{Q_{i}}$,
a property which we will use frequently. In particular, since $\varphi_{i}^{\ast}\in\TestFunctionSpace{\CalO}$
as a finite sum of elements of $\TestFunctionSpace{\CalO}$, we see
that $\overline{Q_{i}}\subset\CalO$ is compact.

\medskip{}

Next, we fix a \textbf{$\CalQ$-moderate weight} $w=\left(w_{i}\right)_{i\in I}$,
which means that $w_{i}\in\left(0,\infty\right)$ for each $i\in I$
and that there is a constant $C_{\CalQ,w}>0$ such that
\begin{equation}
w_{i}\leq C_{\CalQ,w}\cdot w_{\ell}\qquad\forall\:i\in I\text{ and }\ell\in i^{\ast}.\label{eq:IntroductionModerateWeightDefinition}
\end{equation}
Under these assumptions, it follows from \cite[Lemma 4.13]{DecompositionEmbedding}
that the \textbf{$\CalQ$-clustering map}
\begin{equation}
\Gamma_{\CalQ}:\ell_{w}^{q}\left(I\right)\to\ell_{w}^{q}\left(I\right),\left(c_{i}\right)_{i\in I}\mapsto\left(c_{i}^{\ast}\right)_{i\in I}\qquad\text{ with }\qquad c_{i}^{\ast}:=\sum_{\ell\in i^{\ast}}c_{\ell}\label{eq:QClusteringMapDefinition}
\end{equation}
is well-defined and bounded with 
\begin{equation}
\vertiii{\Gamma_{\CalQ}}\leq C_{\CalQ,w}\cdot N_{\CalQ}^{1+\frac{1}{q}}.\label{eq:WeightedSequenceSpaceClusteringMapNormEstimate}
\end{equation}
Here, the \textbf{weighted sequence space $\ell_{w}^{q}\left(I\right)$}
is given by
\[
\ell_{w}^{q}\left(I\right):=\left\{ c=\left(c_{i}\right)_{i\in I}\in\Compl^{I}\with\left\Vert c\right\Vert _{\ell_{w}^{q}}:=\left\Vert \left(w_{i}\cdot c_{i}\right)_{i\in I}\right\Vert _{\ell^{q}}<\infty\right\} ,
\]
for arbitrary $q\in\left(0,\infty\right]$.

\medskip{}

Given all of these assumptions, we define for $p,q\in\left(0,\infty\right]$
the \textbf{Fourier-side decomposition space} associated to $\CalQ$
and the parameters $p,q,v,w$ as
\[
\FourierDecompSp{\CalQ}p{\ell_{w}^{q}}v:=\left\{ f\in\CalD'\left(\CalO\right)\with\left\Vert f\right\Vert _{\FourierDecompSp{\CalQ}p{\ell_{w}^{q}}v}:=\left\Vert \left(\left\Vert \Fourier^{-1}\left(\varphi_{i}f\right)\right\Vert _{L_{v}^{p}}\right)_{i\in I}\right\Vert _{\ell_{w}^{q}}<\infty\right\} .
\]
Finally, we set $Z\left(\CalO\right):=\Fourier\left(\TestFunctionSpace{\CalO}\right)$,
equipped with the unique topology which makes the Fourier transform
$\Fourier:\TestFunctionSpace{\CalO}\to Z\left(\CalO\right)$ a topological
isomorphism. Then, with $Z'\left(\CalO\right)$ denoting the topological
dual space of $Z\left(\CalO\right)$, we define the (space-side) \textbf{decomposition
space} associated to $\CalQ$ and the parameters $p,q,v,w$ as
\[
\DecompSp{\CalQ}p{\ell_{w}^{q}}v:=\left\{ g\in Z'\left(\CalO\right)\with\left\Vert g\right\Vert _{\DecompSp{\CalQ}p{\ell_{w}^{q}}v}:=\left\Vert \widehat{g}\right\Vert _{\FourierDecompSp{\CalQ}p{\ell_{w}^{q}}v}<\infty\right\} ,
\]
where $\Fourier g:=\widehat{g}:=g\circ\Fourier\in\DistributionSpace{\CalO}$
for $g\in Z'\left(\CalO\right)$. It is not hard to see that the Fourier
transform $\Fourier:Z'\left(\CalO\right)\to\DistributionSpace{\CalO}$
restricts to an isometric isomorphism $\Fourier:\DecompSp{\CalQ}p{\ell_{w}^{q}}v\to\FourierDecompSp{\CalQ}p{\ell_{w}^{q}}v$
and that we have $Z'\left(\CalO\right)=\Fourier^{-1}\left(\CalD'\left(\CalO\right)\right)$.

For an explanation for the choice of the reservoirs $\CalD'\left(\CalO\right)$
and $Z'\left(\CalO\right)$, we refer to \cite[Remark 3.13]{DecompositionEmbedding}.
Finally, we mention that \cite[Section 8]{DecompositionEmbedding}
provides a convenient criterion which ensures that each $f\in\DecompSp{\CalQ}p{\ell_{w}^{q}}{}$
extends to a tempered distribution. In particular, if $v\gtrsim1$,
then clearly $\DecompSp{\CalQ}p{\ell_{w}^{q}}v\hookrightarrow\DecompSp{\CalQ}p{\ell_{w}^{q}}{}$.
Hence, if the previously mentioned criterion is fulfilled and if $\CalO=\R^{\dimension}$,
we have (up to trivial identifications) that
\[
\DecompSp{\CalQ}p{\ell_{w}^{q}}v=\left\{ g\in\Schwartz'\left(\smash{\R^{\dimension}}\right)\with\left\Vert g\right\Vert _{\DecompSp{\CalQ}p{\ell_{w}^{q}}v}=\left\Vert \left(\left\Vert \Fourier^{-1}\left(\varphi_{i}\widehat{g}\right)\right\Vert _{L_{v}^{p}}\right)_{i\in I}\right\Vert _{\ell_{w}^{q}}<\infty\right\} .
\]

We remark that the usual papers treating general decomposition spaces
(for general $p,q\in\left(0,\infty\right]$) do usually only consider
the case $v\equiv1$. Hence, it is not entirely clear that the spaces
defined here are indeed well-defined (Quasi)-Banach spaces for $v\not\equiv1$.
We will see below (cf.\@ Proposition \ref{prop:WeightedDecompositionSpaceWellDefined}
and Lemma \ref{lem:WeightedDecompositionSpaceComplete}) that this
is indeed the case.

\subsection{Structure of the paper}

The theory of decomposition spaces is highly dependent on convolutions,
since the very definition of the norm involves quantities of the form
\[
\left\Vert \Fourier^{-1}\left(\varphi_{i}\cdot\widehat{g}\right)\right\Vert _{L_{v}^{p}}=\left\Vert \left(\Fourier^{-1}\varphi_{i}\right)\ast g\right\Vert _{L_{v}^{p}}.
\]
For the range $p\in\left[1,\infty\right]$, Young's inequality $L^{1}\ast L^{p}\hookrightarrow L^{p}$
is usually sufficient to handle such convolutions. But in the range
$p\in\left(0,1\right)$, Young's inequality breaks down completely.
For the usual theory of decomposition spaces, one instead invokes
the convolution relation
\[
\left\Vert f\ast g\right\Vert _{L^{p}}\leq C_{p,\dimension}\cdot R^{\dimension\left(\frac{1}{p}-1\right)}\cdot\left\Vert f\right\Vert _{L^{p}}\cdot\left\Vert g\right\Vert _{L^{p}}\qquad\text{ assuming }\qquad\supp\widehat{f}\subset B_{R}\left(\xi_{1}\right)\text{ and }\supp\widehat{g}\subset B_{R}\left(\xi_{2}\right)
\]
for certain $\xi_{1},\xi_{2}\in\R^{\dimension}$. Note though that
this convolution relation only applies to \emph{band-limited} functions.
But since we are interested in characterizations of decomposition
spaces using (possibly) \emph{compactly supported} functions, this
is not of much use to us.

To overcome this problem, we will invoke the theory of the \textbf{Wiener
amalgam spaces} $W_{Q}\left(L^{\infty},L_{v}^{p}\right)$ which were
originally introduced by Feichtinger\cite{FeichtingerWienerSpaces}.
The main idea is to associate to a (measurable) function $f:\R^{\dimension}\to\Compl$
the local maximal function
\[
M_{Q}f:\R^{\dimension}\to\left[0,\infty\right],x\mapsto\left\Vert \Indicator_{x+Q}\cdot f\right\Vert _{L^{\infty}}
\]
and to define the Wiener amalgam (quasi)-norm of $f$ as $\left\Vert f\right\Vert _{W_{Q}\left(L^{\infty},L_{v}^{p}\right)}=\left\Vert M_{Q}f\right\Vert _{L_{v}^{p}}$.
Broadly speaking, functions in $W_{Q}\left(L^{\infty},L_{v}^{p}\right)$
are locally in $L^{\infty}$ and globally in $L_{v}^{p}$. For brevity,
we will simply write $W_{Q}\left(L_{v}^{p}\right):=W_{Q}\left(L^{\infty},L_{v}^{p}\right)$.
For these spaces, convolution relations are known, cf.\@ \cite{RauhutWienerAmalgam}
and \cite[Section 2.3]{VoigtlaenderPhDThesis}. For our purposes,
however, these results are not sufficient: They establish estimates
of the form
\[
\left\Vert f\ast g\right\Vert _{W_{Q}\left(L_{v}^{p}\right)}\leq C_{p,Q,v}\cdot\left\Vert f\right\Vert _{W_{Q}\left(L_{v_{0}}^{p}\right)}\cdot\left\Vert g\right\Vert _{W_{Q}\left(L_{v}^{p}\right)},
\]
where the constant $C_{p,Q,v}$ depends heavily—and \emph{in an unspecified
way}—on $Q$. But for our purposes, we will consider the spaces $W_{T_{i}^{-T}\left[-1,1\right]^{\dimension}}\left(L^{\infty},L_{v}^{p}\right)$
where $i\in I$ varies; see for example equation (\ref{eq:IntroductionSemiDiscreteQuasiBanachFrame}).
Then, we will need estimates of the form
\[
\left\Vert f\ast g\right\Vert _{W_{T_{i}^{-T}\left[-1,1\right]^{\dimension}}\left(L_{v}^{p}\right)}\leq C_{i,j,\ell,p,v}\cdot\left\Vert f\right\Vert _{W_{T_{j}^{-T}\left[-1,1\right]^{\dimension}}\left(L_{v_{0}}^{p}\right)}\cdot\left\Vert g\right\Vert _{W_{T_{\ell}^{-T}\left[-1,1\right]^{\dimension}}\left(L_{v}^{p}\right)},
\]
with precise control on the constant $C_{i,j,\ell,p,v}$. Hence, in
Section \ref{sec:QuasiBanachConvolutionWienerAmalgam}, we redevelop
parts of the theory of Wiener amalgam spaces, paying close attention
to the dependence of certain constants on the base-set $Q$.

Next, in Section \ref{sec:SemiDiscreteBanachFrames}, we derive assumptions
on the prototype function $\gamma$ which ensure that the norm equivalences
given in equations (\ref{eq:IntroductionSemiDiscreteBanachFrame})
and (\ref{eq:IntroductionSemiDiscreteQuasiBanachFrame}) are true.
More precisely, we will show that the map
\[
\DecompSp{\CalQ}p{\ell_{w}^{q}}v\to\ell_{w}^{q}\left(\left[V_{i}\right]_{i\in I}\right),g\mapsto\left(\gamma^{\left(i\right)}\ast g\right)_{i\in I}
\]
forms a Banach frame, where $V_{i}:=L_{v}^{p}\left(\R^{\dimension}\right)$
in case of $p\in\left[1,\infty\right]$ and $V_{i}:=W_{T_{i}^{-T}\left[-1,1\right]^{\dimension}}\left(L_{v}^{p}\right)$
in case of $p\in\left(0,1\right)$ and where finally
\[
\ell_{w}^{q}\left(\left[V_{i}\right]_{i\in I}\right)=\left\{ \left(g_{i}\right)_{i\in I}\with\left(\forall i\in I:g_{i}\in V_{i}\right)\text{ and }\left(\left\Vert g_{i}\right\Vert _{V_{i}}\right)_{i\in I}\in\ell_{w}^{q}\left(I\right)\right\} .
\]
Part of the problem is to explain how the convolution $\gamma^{\left(i\right)}\ast g$
can be interpreted, especially in case of $\CalO\subsetneq\R^{\dimension}$,
since then each element $g$ of the decomposition space $\DecompSp{\CalQ}p{\ell_{w}^{q}}v$
is the inverse Fourier transform of the distribution $\widehat{g}\in\DistributionSpace{\CalO}$,
so that it is not obvious how $\gamma^{\left(i\right)}\ast g$ can
be understood.

In Section \ref{sec:FullyDiscreteBanachFrames}, we further discretize
the Banach frame $\left(\gamma^{\left(i\right)}\right)_{i\in I}$
from above: Under slightly more strict assumptions on $\gamma$ than
before, we will be able to show that the family $\Psi_{\delta}=\left(L_{\delta\cdot T_{i}^{-T}k}\:\widetilde{\gamma^{\left[i\right]}}\right)_{i\in I,k\in\Z^{\dimension}}$
forms a \textbf{Banach frame} for $\DecompSp{\CalQ}p{\ell_{w}^{q}}v$,
once $\delta>0$ is chosen small enough. Our proof technique is similar
to that of coorbit theory: We use the partition of unity $\left(\varphi_{i}\right)_{i\in I}$
associated to the covering $\CalQ$ to obtain a kind of reproduction
formula, which we then discretize. The details, however, are quite
technical.

Next, in Section \ref{sec:AtomicDecompositions} we establish the
dual statement that the family $\left(L_{\delta\cdot T_{i}^{-T}k}\:\gamma^{\left[i\right]}\right)_{i\in I,k\in\Z^{\dimension}}$
forms an \textbf{atomic decomposition} for $\DecompSp{\CalQ}p{\ell_{w}^{q}}v$.
As above, this is based on a suitable discretization of a certain
reproduction formula.

Finally, since the varying assumptions placed on the prototype $\gamma$
are quite technical and hard to verify, Section \ref{sec:SimplifiedCriteria}
is devoted to a considerable simplification of these conditions. While
not exactly straightforward to verify, these conditions can be verified
in practice, where the degree of difficulty mainly depends on the
given covering $\CalQ=\left(T_{i}Q_{i}'+b_{i}\right)_{i\in I}$.

As a litmus test of our theory, we show in Section \ref{sec:CompactlySupportedAlphaModulationFrames}
that it can be used to obtain compactly supported structured Banach
frames and atomic decompositions for the $\alpha$-modulation spaces
$\AlphaModSpace pq{\alpha}s\left(\R^{\dimension}\right)$, even for
$p,q<1$, thereby extending the state of the art. Furthermore, in
Section \ref{sec:BesovFrames}, we show that our theory can be used
to establish that certain compactly supported wavelet systems generate
Banach frames and atomic decompositions for inhomogeneous Besov spaces.

We emphasize that we consider these two specific examples since they
can be handled with reasonably low effort, but still indicate that—\emph{and
how}—the general theory can be filled with life for concrete special
cases. The theory presented here certainly has more interesting and
more novel applications, in particular to the theory of shearlets.
But in order to keep the size of this paper somewhat manageable, we
postpone these applications to the companion paper \cite{StructuredBanachFrames2}.

\subsection*{Credit where credit is due}

\epigraph{``[...] virtually all of our techniques already exist in some antecedent form. Nevertheless their particular combination here leads to new conclusions and to sharpened versions of known results. Moreover, our presentation reveals a[...] structure underlying a diverse range of topics in harmonic analysis.''}{M.\@ Frazier and B.\@ Jawerth, \cite[Page 36]{FrazierJawerthDiscreteTransform}}

The results and proof techniques employed in this paper were heavily
inspired by a number of earlier results:

The first impulse for writing this paper was caused by my reading
of the paper \cite{EmbeddingsOfAlphaModulationIntoSobolev}. In that
paper, the author characterizes the existence of embeddings between
$\alpha$-modulation spaces and Sobolev spaces. As an intermediate
result, he also proves
\begin{equation}
\left\Vert g\right\Vert _{\AlphaModSpace pq{\alpha}s}\asymp\left\Vert \left(\left\Vert \smash{\gamma^{\left(i\right)}}\ast g\right\Vert _{L^{p}}\right)_{i\in\Z^{\dimension}}\right\Vert _{\ell_{\left(\left\langle k\right\rangle ^{s}\right)_{k\in\Z^{\dimension}}}^{q}}\label{eq:AlphaModulationNormCharacterization}
\end{equation}
for arbitrary $p,q\in\left(0,\infty\right]$, $\alpha\in\left[0,1\right)$
and $s\in\R$, as well as $g\in\AlphaModSpace pq{\alpha}s\left(\R^{\dimension}\right)$,
if the prototype function $\gamma\in\Schwartz\left(\R^{\dimension}\right)$
is chosen suitably. Here, the functions $\gamma^{\left(i\right)}$
for $i\in\Z^{\dimension}$ are formed from $\gamma$ as described
before equation (\ref{eq:IntroductionL2NormalizedFrameElement}),
where $\CalQ=\CalQ^{\left(\alpha\right)}=\left(\left\langle k\right\rangle ^{\frac{\alpha}{1-\alpha}}\cdot B_{R}\left(0\right)+\left\langle k\right\rangle ^{\frac{\alpha}{1-\alpha}}k\right)_{k\in\Z^{\dimension}}$
is the usual covering used to define $\alpha$-modulation spaces;
see also Section \ref{sec:CompactlySupportedAlphaModulationFrames}.
Note that—at least for $p\in\left[1,\infty\right]$—this result is
a special case of the results about semi-discrete Banach frames from
Section \ref{sec:SemiDiscreteBanachFrames}. Specifically, the paper
\cite{EmbeddingsOfAlphaModulationIntoSobolev} caused me to investigate
whether a norm characterization as in equation (\ref{eq:AlphaModulationNormCharacterization})
was also possible in the more general setting of (essentially) arbitrary
decomposition spaces and not only for $\alpha$-modulation spaces.
In particular, it caused me to consider the structured families of
the form $\left(\gamma^{\left(i\right)}\right)_{i\in I}$ with $\gamma^{\left(i\right)}=\left|\det T_{i}\right|\cdot M_{b_{i}}\left[\gamma\circ T_{i}^{T}\right]$,
where $\CalQ=\left(T_{i}Q_{i}'+b_{i}\right)_{i\in I}$.

Furthermore, an investigation of the proofs in \cite{EmbeddingsOfAlphaModulationIntoSobolev}
lead me to consider assumptions similar to those stated in Assumption
\ref{assu:MainAssumptions} below. Specifically, it caused me to impose
boundedness of the operator associated to the infinite matrix $\left(\left\Vert \Fourier^{-1}\left(\varphi_{i}\cdot\widehat{\gamma^{\left(j\right)}}\right)\right\Vert _{L^{1}}\right)_{j,i\in I}$.
In summary, at least for the case $p\in\left[1,\infty\right]$, the
results about semi-discrete Banach frames for decomposition spaces
in this paper (cf.\@ Section \ref{sec:SemiDiscreteBanachFrames})
can be seen as a slight generalization of the results in \cite{EmbeddingsOfAlphaModulationIntoSobolev}.

\medskip{}

For the case $p\in\left(0,1\right)$, however, I was not able to adapt
the techniques used in \cite{EmbeddingsOfAlphaModulationIntoSobolev}
to the general setting of decomposition spaces. In fact, for $p\in\left(0,1\right)$,
the results derived in \cite{EmbeddingsOfAlphaModulationIntoSobolev}
differ from those in Section \ref{sec:SemiDiscreteBanachFrames}:
While the characterization from \cite{EmbeddingsOfAlphaModulationIntoSobolev}
(cf.\@ equation (\ref{eq:AlphaModulationNormCharacterization}))
considers the usual $L^{p}$ norm of the convolutions $\gamma^{\left(i\right)}\ast f$,
in Section \ref{sec:SemiDiscreteBanachFrames} we show for $p\in\left(0,1\right)$
that
\[
\left\Vert g\right\Vert _{\DecompSp{\CalQ}p{\ell_{w}^{q}}{}}\asymp\left\Vert \left(\left\Vert \smash{\gamma^{\left(i\right)}}\ast g\right\Vert _{W_{T_{i}^{-T}\left[-1,1\right]^{\dimension}}\left(L^{p}\right)}\right)_{i\in I}\right\Vert _{\ell_{w}^{q}},
\]
i.e., we use \textbf{Wiener amalgam spaces} instead of the spaces
$L^{p}$ themselves.

Here, again, I was inspired heavily by earlier results: The main limitation
of the spaces $L^{p}\left(\R^{\dimension}\right)$ for $p\in\left(0,1\right)$
in the present setting is that there are no meaningful convolution
relations for them, partly since we do not even have $L^{p}\left(\R^{\dimension}\right)\hookrightarrow L_{{\rm loc}}^{1}\left(\R^{\dimension}\right)$.
Luckily, Holger Rauhut\cite{RauhutCoorbitQuasiBanach} had already
observed—while generalizing coorbit theory \cite{FeichtingerCoorbit0,FeichtingerCoorbit1,FeichtingerCoorbit2,GroechenigDescribingFunctions}
to the setting of Quasi-Banach spaces—that these limitations can be
avoided by considering the Wiener amalgam spaces $W_{Q}\left(L^{p}\right)$
instead of $L^{p}$ itself. Rauhut had also already developed associated
convolution relations\cite{RauhutWienerAmalgam} for these spaces.
Of course, all of this was based on the original invention of Wiener
Amalgam spaces which is due to Hans Feichtinger\cite{FeichtingerWienerInterpolation,FeichtingerWienerSpaces}.

All in all, given these earlier papers, it was natural to consider
Wiener amalgam spaces. The (as far as I know) novel idea was to consider
these Wiener amalgam spaces $W_{Q}\left(L^{p}\right)$ with a definite
choice of the base set $Q$, which was allowed to heavily vary with
$i\in I$. Further, it seems to be a new (or at least not completely
well-known) fact that suitably bandlimited $L^{p}$ functions automatically
belong to $W_{Q}\left(L^{p}\right)$, where this statement comes with
a precise estimate for the Wiener amalgam norm in terms of $Q$ and
the Fourier support of the function.

\medskip{}

At this point, I had managed to generalize the results about semi-discrete
Banach frames developed in \cite{EmbeddingsOfAlphaModulationIntoSobolev}
to the setting of general decomposition spaces. One of my main goals,
however, was a better understanding of the approximation theoretic
properties of discrete, cone-adapted shearlet systems. To achieve
this, a further discretization of these \emph{semi}-discrete Banach
frames was necessary. The inspiration for treating this additional
discretization step came from the theory of coorbit spaces as developed
by Feichtinger and Gröchenig\cite{FeichtingerCoorbit0,FeichtingerCoorbit1,FeichtingerCoorbit2,GroechenigDescribingFunctions}
and also (in more generalized form) by Rauhut, Fornasier and Ullrich\cite{RauhutCoorbitQuasiBanach,GeneralizedCoorbit1,GeneralizedCoorbit2}.
The underlying important idea of coorbit theory is to transfer the
study of certain function spaces via a suitable transform to the study
of certain Banach spaces which have a \textbf{reproducing property}.
Formally, one employs the so-called voice transform $V$ to establish
an isomorphism between the coorbit space ${\rm Co}\left(Y\right)$
and its image $Z:=V\left[{\rm Co}\left(Y\right)\right]$ under the
voice transform. The crucial property of the space $Z$ is that we
have the \textbf{reproducing formula}
\[
F=F\ast G\qquad\forall F\in Z
\]
for a suitable kernel $G$. In fact, in the setting of generalized
coorbit theory, the convolution with $G$ needs to be replaced by
a more general integral operator.

If the kernel $G$ is regular enough, the reproducing formula allows
to show that a sufficiently dense sampling of $F\in Z$ suffices to
reconstruct $F$ uniquely. Proving this is based on a (suitable) notion
of the \textbf{oscillation} of a function. This sampling result can
then be transferred to the coorbit space ${\rm Co}\left(Y\right)$
to obtain Banach frames and atomic decompositions. Similar techniques
are also used in \cite{GroechenigNonuniformSampling2}.

The new contribution was thus to derive a suitable reproducing formula
in the general setting of decomposition spaces, cf.\@ Lemma \ref{lem:SpecialProjection}.
Once this was established, existing ideas and techniques could be
used to obtain the desired discrete Banach frames and atomic decompositions.
We remark, however, that the established reproducing formula for decomposition
spaces is highly nontrivial.

\medskip{}

In total, the present paper would not have been possible without inspiration
from existing results, concepts and techniques (Wiener amalgam spaces
and their convolution relations, oscillation of a function, semi-discrete
Banach frames for $\alpha$-modulation spaces, etc.). The contribution
of the paper is that these results and techniques are combined and
refined to achieve novel and nontrivial results which—due to their
generality—apply in a wide variety of settings.

\subsection*{A comment on constants}

Instead of using only implied constants of the form $C=C\left(\dimension,p,\CalQ,...\right)$,
in this paper we try to provide explicit constants whenever possible.
In principle, this allows one e.g.\@ to determine an \emph{explicit}
$\delta_{0}>0$ such that the family $\Psi_{\delta}$ defined in equation
(\ref{eq:IntroductionBanachFrameFamily}) yields a Banach frame for
the decomposition space under consideration for $0<\delta\leq\delta_{0}$.
We make no effort, however, to produce the optimal (or even good)
constants. Occasionally, we even enlarge appearing constants just
to make the expressions for the constants in question more optically
pleasing (i.e., shorter). Due to these reasons, the resulting sampling
density $\delta_{0}$ will probably be of size $\delta_{0}\approx2^{-1000}$
or even smaller.

Thus, our leading philosophy is that \emph{an arbitrarily bad explicit
constant is still (much) better than an implicit constant which one
does not know at all}.

\section{Convolution in $L^{p},p\in\left(0,1\right)$ and Wiener Amalgam spaces}

\label{sec:QuasiBanachConvolutionWienerAmalgam}The well-known Young
inequality $\left\Vert f\ast g\right\Vert _{L^{p}}\leq\left\Vert f\right\Vert _{L^{1}}\cdot\left\Vert g\right\Vert _{L^{p}}$
fails for $p\in\left(0,1\right)$, cf.\@ \cite[Example 3.1]{DecompositionEmbedding}.
One can solve this in two ways: The first way is given in \cite[Proposition 1.5.1]{TriebelTheoryOfFunctionSpaces},
where it is shown that
\[
\left\Vert f\ast g\right\Vert _{L^{p}}\lesssim\left\Vert f\right\Vert _{L^{p}}\cdot\left\Vert g\right\Vert _{L^{p}}
\]
if $f$ and $g$ are both bandlimited. This theorem, however, has
two disadvantages:
\begin{itemize}
\item The restriction to bandlimited $f,g$ is rather severe; in particular
in our present setting, since we are interested in compactly supported
functions, which can never be bandlimited.
\item The implicit constant in the estimate above depends in a nontrivial
way on the frequency supports of $f,g$.
\end{itemize}
To overcome these limitations, we will develop an improved theory
of convolution for $L^{p},p\in\left(0,1\right)$ using the theory
of \textbf{Wiener amalgam spaces}. As a special case, we will recover
the estimate from above.

Before developing the theory, we remark that essentially everything
mentioned in this section is already known in one form or another.
In particular, Wiener amalgam spaces were originally invented by Feichtinger\cite{FeichtingerWienerInterpolation,FeichtingerWienerSpaces}
and later generalized to Quasi-Banach spaces by Rauhut\cite{RauhutCoorbitQuasiBanach}.
The use of these spaces—and of the oscillation of a function—for obtaining
Banach frames and atomic decompositions for certain spaces goes back
to the theory of coorbit spaces\cite{FeichtingerCoorbit0,FeichtingerCoorbit1,FeichtingerCoorbit2,GroechenigDescribingFunctions,GeneralizedCoorbit1,GeneralizedCoorbit2}
and was also exploited in \cite{GroechenigNonuniformSampling2}. Therefore,
no originality is claimed.

The usual treatments, however, mostly ignore or suppress the dependence
of the Wiener amalgam spaces on the chosen unit neighborhood (see
below for details), whereas this dependence is crucial for us. Hence,
we provide full proofs.

\subsection{Definition of Wiener amalgam spaces}

All of the theory of Wiener amalgam spaces is centered around the
notion of a certain maximal function:
\begin{defn}
\label{def:MaximalFunctionDefinition}(cf.\@ \cite{FeichtingerWienerSpaces},
\cite[Definition 2.2.2]{HeilPhDThesisWienerAmalgam}, \cite{RauhutWienerAmalgam}
and \cite[Definition 2.3.1]{VoigtlaenderPhDThesis}) Let $Q\subset\R^{\dimension}$
be a Borel measurable unit neighborhood and let $f:\R^{\dimension}\to\Compl^{k}$
be Borel measurable. We then define the \textbf{$Q$-maximal function}
of $f$ as
\[
M_{Q}f:\R^{\dimension}\to\left[0,\infty\right],x\mapsto\essup_{y\in x+Q}\left|f\left(y\right)\right|=\essup_{a\in Q}\left|f\left(x+a\right)\right|=\left\Vert L_{-x}f\right\Vert _{L^{\infty}\left(Q\right)}.
\]
For a given $p\in\left(0,\infty\right]$ and a (measurable) weight
$u:\R^{\dimension}\to\left(0,\infty\right)$, we define the \textbf{Wiener
amalgam space} with window $Q$, local component $L^{\infty}$ and
global component $L_{u}^{p}$ as
\[
W_{Q}^{k}\left(L_{u}^{p}\right):=W_{Q}^{k}\left(L^{\infty},L_{u}^{p}\right):=\left\{ f:\R^{\dimension}\to\Compl^{k}\with f\text{ measurable and }M_{Q}f\in L_{u}^{p}\left(\smash{\R^{\dimension}}\right)\right\} ,
\]
with the natural (quasi)-norm $\left\Vert f\right\Vert _{W_{Q}^{k}\left(L_{u}^{p}\right)}:=\left\Vert M_{Q}f\right\Vert _{L_{u}^{p}}$.
In the most common case $k=1$, we omit the exponent and write $W_{Q}\left(L_{u}^{p}\right)$
instead of $W_{Q}^{1}\left(L_{u}^{p}\right)$.
\end{defn}
\begin{rem*}

\begin{itemize}[leftmargin=0.4cm]
\item One can show for suitable weights $u$ (and we will do so in Lemma
\ref{lem:WienerAmalgamNormEquivalence}) that the space $W_{Q}\left(L_{u}^{p}\right)$
is independent of the choice of the \emph{bounded} measurable unit
neighborhood $Q\subset\R^{\dimension}$, with equivalent quasi-norms
for different choices. Hence, $Q$ is often suppressed in the literature
dealing with Wiener amalgam spaces. For us, however, the precise choice
of $Q$ will be crucial, since we will choose $Q_{i}=T_{i}^{-T}\left[-1,1\right]^{\dimension}$,
so that the sets $Q_{i}$ vary wildly with $i\in I$. Since the constants
appearing in the norm equivalences for different choices of $Q$ depend
heavily on the actual choices of $Q$, we will almost never use the
equivalence for different choices of $Q$, or only in very carefully
chosen ways.
\item Note that $M_{Q}f$ is always a Borel measurable function. Indeed,
since $L^{1}\left(\R^{\dimension}\right)$ is separable, there is
a countable dense family $\left(g_{n}\right)_{n\in\N}$ in $\Gamma:=\left\{ g\in L^{1}\left(\R^{\dimension}\right)\with g\geq0\text{ and }\left\Vert g\right\Vert _{L^{1}}\leq1\right\} $.
Then, we have for an arbitrary measurable function $f$ that
\begin{equation}
\left\Vert f\right\Vert _{L^{\infty}}=\sup_{n\in\N}\int_{\R^{\dimension}}g_{n}\left(x\right)\cdot\left|f\left(x\right)\right|\d x.\label{eq:LInftyNormUsingCountableFamily}
\end{equation}
For $f\in L^{\infty}\left(\R^{\dimension}\right)$, this follows from
the usual characterization of the $L^{\infty}$-norm by duality (cf.\@
\cite[Theorem 6.14]{FollandRA}). In case of $\left\Vert f\right\Vert _{L^{\infty}}=\infty$,
the same theorem shows that for each $M>0$, there is some $g\in\Gamma$
satisfying $\int_{\R^{\dimension}}g\left(x\right)\cdot\left|f\left(x\right)\right|\d x\geq M$.
But by density of the family $\left(g_{n}\right)_{n\in\N}$, there
is then a sequence $\left(n_{k}\right)_{k\in\N}$ such that $g_{n_{k}}\to g$
in $L^{1}\left(\R^{\dimension}\right)$. By switching to a subsequence,
we can also assume $g_{n_{k}}\to g$ almost everywhere. Now, Fatou's
Lemma yields
\[
M\leq\int_{\R^{\dimension}}g\left(x\right)\cdot\left|f\left(x\right)\right|\d x=\int_{\R^{\dimension}}\liminf_{k\to\infty}g_{n_{k}}\left(x\right)\cdot\left|f\left(x\right)\right|\d x\leq\liminf_{k\to\infty}\int_{\R^{\dimension}}g_{n_{k}}\left(x\right)\cdot\left|f\left(x\right)\right|\d x.
\]
Since $M>0$ was arbitrary, this easily yields $\sup_{n\in\N}\int_{\R^{\dimension}}g_{n}\left(x\right)\cdot\left|f\left(x\right)\right|\d x=\infty=\left\Vert f\right\Vert _{L^{\infty}}$.

Now, as a consequence of equation (\ref{eq:LInftyNormUsingCountableFamily}),
we get
\begin{align*}
\left(M_{Q}f\right)\left(x\right) & =\left\Vert \Indicator_{Q}\cdot L_{-x}f\right\Vert _{L^{\infty}}\\
 & =\sup_{n\in\N}\int_{\R^{\dimension}}g_{n}\left(y\right)\cdot\Indicator_{Q}\left(y\right)\cdot\left|\left(L_{-x}f\right)\left(y\right)\right|\d y\\
 & =\sup_{n\in\N}\int_{\R^{\dimension}}g_{n}\left(y\right)\cdot\Indicator_{Q}\left(y\right)\cdot\left|f\left(x+y\right)\right|\d y.
\end{align*}
But the function $\left(x,y\right)\mapsto g_{n}\left(y\right)\cdot\Indicator_{Q}\left(y\right)\cdot\left|f\left(x+y\right)\right|$
is Borel measurable, so that measurability of the integrated function
$x\mapsto\int_{\R^{\dimension}}g_{n}\left(y\right)\cdot\Indicator_{Q}\left(y\right)\cdot\left|f\left(x+y\right)\right|\d y$
follows from the Fubini-Tonelli theorem. Hence, $M_{Q}f$ is Borel
measurable.\qedhere
\end{itemize}
\end{rem*}
It is easy to see that $\left\Vert \cdot\right\Vert _{W_{Q}\left(L_{u}^{p}\right)}$
satisfies the (quasi)-triangle inequality, since $M_{Q}\left(f+g\right)\leq M_{Q}f+M_{Q}g$.
The remaining properties of a (quasi)-norm are also easy to check,
possibly with the exception of definiteness. But this is a consequence
of the following lemma:
\begin{lem}
\label{lem:MaximalFunctionDominatesF}For each Borel measurable unit
neighborhood $Q$ and each measurable $f:\R^{\dimension}\to\Compl^{k}$,
we have 
\[
\left|f\left(x\right)\right|\leq\left(M_{Q}f\right)\left(x\right)\qquad\text{ for almost all }x\in\R^{\dimension}.
\]
In particular, $\left\Vert f\right\Vert _{L_{u}^{p}}\leq\left\Vert f\right\Vert _{W_{Q}\left(L_{u}^{p}\right)}$
and hence $f=0$ almost everywhere if $\left\Vert f\right\Vert _{W_{Q}\left(L_{u}^{p}\right)}=0$.
\end{lem}
\begin{proof}
Since $Q$ is a unit neighbhorhood, there is $\varepsilon>0$ with
$B_{2\varepsilon}\left(0\right)\subset Q$. Since $\R^{\dimension}$
is second countable and since $\left(x+B_{\varepsilon}\left(0\right)\right)_{x\in\R^{\dimension}}$
is an open cover of $\R^{\dimension}$, there is a countable family
$\left(x_{n}\right)_{n\in\N}$ satisfying $\R^{\dimension}=\bigcup_{n\in\N}\left(x_{n}+B_{\varepsilon}\left(0\right)\right)$.

Now, for arbitrary $x\in\R^{\dimension}$ and $y\in x+B_{\varepsilon}\left(0\right)$,
we have $y+B_{\varepsilon}\left(0\right)\subset x+B_{2\varepsilon}\left(0\right)\subset x+Q$
and hence
\[
\left(M_{B_{\varepsilon}\left(0\right)}f\right)\left(y\right)=\left\Vert f\cdot\Indicator_{y+B_{\varepsilon}\left(0\right)}\right\Vert _{L^{\infty}}\leq\left\Vert f\cdot\Indicator_{x+Q}\right\Vert _{L^{\infty}}=M_{Q}f\left(x\right)\qquad\forall x\in\R^{\dimension}\,\forall y\in x+B_{\varepsilon}\left(0\right).
\]

Next, for each $n\in\N$, there is a null-set $N_{n}\subset x_{n}+B_{\varepsilon}\left(0\right)$
such that 
\[
\left|f\left(x\right)\right|\leq\left\Vert f\cdot\Indicator_{x_{n}+B_{\varepsilon}\left(0\right)}\right\Vert _{L^{\infty}}=\left(M_{B_{\varepsilon}\left(0\right)}f\right)\left(x_{n}\right)\qquad\text{ for all }\qquad x\in\left[x_{n}+B_{\varepsilon}\left(0\right)\right]\setminus N_{n}.
\]
But for each such $x$, there is some $\gamma\in B_{\varepsilon}\left(0\right)$
such that $x=x_{n}+\gamma$ and thus $x_{n}=x-\gamma\in x+B_{\varepsilon}\left(0\right)$,
so that the equation from above yields $\left|f\left(x\right)\right|\leq\left(M_{B_{\varepsilon}\left(0\right)}f\right)\left(x_{n}\right)\leq M_{Q}f\left(x\right)$.
Recall that this estimate holds for all $x\in\left[x_{n}+B_{\varepsilon}\left(0\right)\right]\setminus N_{n}$.
Since $N:=\bigcup_{n\in\N}N_{n}$ is a null-set and since $\R^{\dimension}=\bigcup_{n\in\N}\left(x_{n}+B_{\varepsilon}\left(0\right)\right)$,
this completes the proof.
\end{proof}
Although easy to prove, the following lemma is frequently helpful,
since it shows that Schwartz functions are contained in arbitrary
Wiener amalgam spaces.
\begin{lem}
\label{lem:SchwartzFunctionsAreWiener}For arbitrary $N\geq0$, we
have
\[
\left[M_{\left[-1,1\right]^{\dimension}}\left(1+\left|\mybullet\right|\right)^{-N}\right]\left(x\right)\leq\left(1+2\sqrt{\dimension}\right)^{N}\cdot\left(1+\left|x\right|\right)^{-N}\qquad\forall x\in\R^{\dimension}.
\]
In particular, if $p\in\left(0,\infty\right]$ is arbitrary and if
we set $\left\Vert f\right\Vert _{N}:=\sup_{x\in\R^{\dimension}}\left(1+\left|x\right|\right)^{N}\cdot\left|f\left(x\right)\right|$
for measurable $f:\R^{\dimension}\to\Compl$, then
\[
\left\Vert f\right\Vert _{W_{\left[-1,1\right]^{\dimension}}\left(L_{\left(1+\left|\cdot\right|\right)^{K}}^{p}\right)}\leq\left(1+2\sqrt{\dimension}\right)^{N}\cdot\left(\frac{1}{p}\frac{s_{\dimension}}{N-K-\frac{\dimension}{p}}\right)^{1/p}\cdot\left\Vert f\right\Vert _{N}\qquad\text{ as soon as }\qquad N>K+\frac{\dimension}{p}.
\]
Hence, $\Schwartz\left(\R^{\dimension}\right)\hookrightarrow W_{\left[-1,1\right]^{\dimension}}\left(L_{u}^{p}\right)$
for all $p\in\left(0,\infty\right]$ and $u\in\left\{ v,v_{0},\left(1+\left|\mybullet\right|\right)^{K}\right\} $.
\end{lem}
\begin{proof}
To prove the first claim, we distinguish two cases. For $\left|x\right|\leq2\sqrt{\dimension}$,
note that $\left(1+\left|x+a\right|\right)^{-N}\leq1$ for all $a\in\left[-1,1\right]^{\dimension}$,
so that we get 
\[
\left[M_{\left[-1,1\right]^{\dimension}}\left(1+\left|\mybullet\right|\right)^{-N}\right]\left(x\right)\leq1\leq\left(1+2\sqrt{\dimension}\right)^{N}\cdot\left(1+\left|x\right|\right)^{-N}.
\]
Otherwise, if $\left|x\right|\geq2\sqrt{\dimension}$, we have for
$a\in\left[-1,1\right]^{\dimension}$ that
\[
\left|x-a\right|\geq\left|x\right|-\left|a\right|\geq\left|x\right|-\sqrt{\dimension}=\frac{\left|x\right|}{2}+\frac{\left|x\right|}{2}-\sqrt{\dimension}\geq\frac{\left|x\right|}{2}
\]
and hence $\left(1+\left|x+a\right|\right)^{-N}\leq\left(1+\frac{\left|x\right|}{2}\right)^{-N}\leq\left(\frac{1}{2}\left(1+\left|x\right|\right)\right)^{-N}=2^{N}\left(1+\left|x\right|\right)^{-N}$.
Since $2^{N}\leq\left(1+2\sqrt{\dimension}\right)^{N}=:C$, we get
all in all that $\left(M_{\left[-1,1\right]^{\dimension}}\left(1+\left|\mybullet\right|\right)^{-N}\right)\left(x\right)\leq C\cdot\left(1+\left|x\right|\right)^{-N}$
for all $x\in\R^{\dimension}$, as claimed.

For the next claim, we can clearly assume $\left\Vert f\right\Vert _{N}<\infty$.
In this case, we have $\left|f\left(x\right)\right|\leq\left\Vert f\right\Vert _{N}\cdot\left(1+\left|x\right|\right)^{-N}$
for all $x\in\R^{\dimension}$ and hence 
\begin{align*}
\left(1+\left|x\right|\right)^{K}\cdot\left(M_{\left[-1,1\right]^{\dimension}}f\right)\left(x\right) & \leq\left\Vert f\right\Vert _{N}\cdot\left(1+\left|x\right|\right)^{K}\cdot\left(M_{\left[-1,1\right]^{\dimension}}\left(1+\left|\mybullet\right|\right)^{-N}\right)\left(x\right)\\
 & \leq C\cdot\left\Vert f\right\Vert _{N}\cdot\left(1+\left|x\right|\right)^{-\left(N-K\right)}\qquad\forall x\in\R^{\dimension}.
\end{align*}
This yields the claim, since equation (\ref{eq:StandardDecayLpEstimate})
shows $\left\Vert \left(1+\left|\mybullet\right|\right)^{-\left(N-K\right)}\right\Vert _{L^{p}}\leq\left(\frac{1}{p}\frac{s_{\dimension}}{N-K-\frac{\dimension}{p}}\right)^{1/p}<\infty$,
as long $N>K+\frac{\dimension}{p}$. The embedding $\Schwartz\left(\R^{\dimension}\right)\hookrightarrow W_{\left[-1,1\right]^{\dimension}}\left(L_{u}^{p}\right)$
for all $p\in\left(0,\infty\right]$ and $u\in\left\{ v,v_{0},\left(1+\left|\mybullet\right|\right)^{K}\right\} $
is now trivial, since the norm $\left\Vert \mybullet\right\Vert _{N}$
is continuous with respect to the topology on $\Schwartz\left(\R^{\dimension}\right)$
and since we have
\[
v\left(x\right)=v\left(0+x\right)\leq v\left(0\right)\cdot v_{0}\left(x\right)\leq\Omega_{1}v\left(0\right)\cdot\left(1+\left|x\right|\right)^{K}\qquad\forall x\in\R^{\dimension}.\qedhere
\]
\end{proof}
\begin{lem}
\label{lem:WienerTransformationFormula}For $k\in\N$, a measurable
$f:\R^{\dimension}\to\Compl^{k}$, $T\in\GL\left(\R^{\dimension}\right)$
and a measurable $Q\subset\R^{\dimension}$, we have
\[
M_{Q}\left(f\circ T\right)=\left(M_{TQ}f\right)\circ T
\]
and hence
\[
\left\Vert f\circ T\right\Vert _{W_{Q}\left(L^{p}\right)}=\left|\det T\right|^{-1/p}\cdot\left\Vert f\right\Vert _{W_{TQ}\left(L^{p}\right)}.\qedhere
\]
\end{lem}
\begin{proof}
Since $T$ and $T^{-1}$ map null-sets to null-sets, we have
\begin{align*}
\left[M_{Q}\left(f\circ T\right)\right]\left(x\right) & =\essup_{a\in Q}\left|\left(f\circ T\right)\left(x+a\right)\right|\\
 & =\essup_{b\in TQ}\left|f\left(b+Tx\right)\right|\\
 & =\left(M_{TQ}f\right)\left(Tx\right)
\end{align*}
for all $x\in\R^{\dimension}$. The final identity is a consequence
of the definitions and of $\left\Vert f\circ T\right\Vert _{L^{p}}=\left|\det T\right|^{-1/p}\cdot\left\Vert f\right\Vert _{L^{p}}$,
which follows easily from the change-of-variables formula.
\end{proof}
Next, we show that iterated applications of $M_{Q}$ can be estimated
using a single $M_{Q'}$.
\begin{lem}
\label{lem:IteratedMaximalFunction}Let $k\in\N$ and assume that
$f:\R^{\dimension}\to\Compl^{k}$ and $Q_{1},Q_{2}\subset\R^{\dimension}$
are measurable and that $Q_{1}+Q_{2}$ is also measurable. Then
\[
M_{Q_{1}}\left[M_{Q_{2}}f\right]\leq M_{Q_{1}+Q_{2}}f.
\]
In particular, for any measurable $u:\R^{\dimension}\to\left(0,\infty\right)$,
we have
\[
\left\Vert M_{Q_{2}}f\right\Vert _{W_{Q_{1}}\left(L_{u}^{p}\right)}\leq\left\Vert f\right\Vert _{W_{Q_{1}+Q_{2}}^{k}\left(L_{u}^{p}\right)}.\qedhere
\]
\end{lem}
\begin{proof}
For $a\in Q_{1}$, we have $a+Q_{2}\subset Q_{1}+Q_{2}$ and hence
\[
\left(M_{Q_{2}}f\right)\left(x+a\right)=\left\Vert f\cdot\Indicator_{x+a+Q_{2}}\right\Vert _{L^{\infty}}\leq\left\Vert f\cdot\Indicator_{x+Q_{1}+Q_{2}}\right\Vert _{L^{\infty}}=\left(M_{Q_{1}+Q_{2}}f\right)\left(x\right).
\]
Since this holds for all $a\in Q_{1}$, we get $\left(M_{Q_{1}}\left[M_{Q_{2}}f\right]\right)\left(x\right)=\essup_{a\in Q_{1}}\left(M_{Q_{2}}f\right)\left(x+a\right)\leq\left(M_{Q_{1}+Q_{2}}f\right)\left(x\right)$,
as claimed.
\end{proof}
The next three lemmas are important for us, since they imply $\left\Vert f\right\Vert _{W_{T_{j}^{-T}\left[-1,1\right]^{\dimension}}\left(L_{v}^{p}\right)}\!\leq\!C_{i,j,p,v}\cdot\left\Vert f\right\Vert _{W_{T_{i}^{-T}\left[-1,1\right]^{\dimension}}\left(L_{v}^{p}\right)}$,
where the constant $C_{i,j,p,v}$ is explicitly known, cf.\@ Corollary
\ref{cor:WienerLinearCubeNormEstimate}. We begin with an estimate
for the norm of the translation operators on $L_{v}^{p}\left(\R^{\dimension}\right)$.
\begin{lem}
\label{lem:WeightedLpTranslationNorm}For each $y\in\R^{\dimension}$,
the left-translation operator $L_{y}:L_{v}^{p}\left(\R^{\dimension}\right)\to L_{v}^{p}\left(\R^{\dimension}\right)$
is well-defined and bounded with
\[
\vertiii{L_{y}}\leq v_{0}\left(y\right)\leq\Omega_{1}\cdot\left(1+\left|y\right|\right)^{K}.\qedhere
\]
\end{lem}
\begin{rem*}
The only property of $v$ which is used in the proof is that $v\left(x+y\right)\leq v\left(x\right)v_{0}\left(y\right)$.
By submultiplicativity of $v_{0}$, the same estimate holds for $v_{0}$
instead of $v$, so that the claim of the lemma also holds for $v_{0}$
instead of $v$.
\end{rem*}
\begin{proof}
Let $f\in L_{v}^{p}\left(\R^{\dimension}\right)$ be arbitrary and
note
\begin{align*}
v\left(x\right)\cdot\left|\left(L_{y}f\right)\left(x\right)\right| & =v\left(\left(x-y\right)+y\right)\cdot\left|f\left(x-y\right)\right|\\
 & \leq v_{0}\left(y\right)\cdot\left|\left(v\cdot f\right)\left(x-y\right)\right|\\
 & =v_{0}\left(y\right)\cdot\left[L_{y}\left(v\cdot f\right)\right]\left(x\right).
\end{align*}
By solidity and translation invariance of $L^{p}\left(\R^{\dimension}\right)$,
this implies
\[
\left\Vert L_{y}f\right\Vert _{L_{v}^{p}}=\left\Vert v\cdot L_{y}f\right\Vert _{L^{p}}\leq v_{0}\left(y\right)\cdot\left\Vert L_{y}\left(v\cdot f\right)\right\Vert _{L^{p}}=v_{0}\left(y\right)\cdot\left\Vert v\cdot f\right\Vert _{L^{p}}=v_{0}\left(y\right)\cdot\left\Vert f\right\Vert _{L_{v}^{p}}<\infty.\qedhere
\]
\end{proof}
Now we can derive a first estimate which will allow us to switch from
one ``base set'' $Q$ to another one.
\begin{lem}
\label{lem:WienerAmalgamNormEquivalence}Let $Q_{1},Q_{2}\subset\R^{\dimension}$
and assume that there are $x_{1},\dots,x_{N}\in\R^{\dimension}$ such
that $Q_{1}\subset\bigcup_{i=1}^{N}\left(x_{i}+Q_{2}\right)$. Let
$p\in\left(0,\infty\right]$ and set $s:=\min\left\{ 1,p\right\} $.
Then we have
\[
\left\Vert f\right\Vert _{W_{Q_{1}}^{k}\left(L_{v}^{p}\right)}\leq\left(\sum_{i=1}^{N}\left[v_{0}\left(-x_{i}\right)\right]^{s}\right)^{1/s}\cdot\left\Vert f\right\Vert _{W_{Q_{2}}^{k}\left(L_{v}^{p}\right)}\leq\Omega_{1}\cdot\left(\sum_{i=1}^{N}\left(1+\left|x_{i}\right|\right)^{sK}\right)^{1/s}\cdot\left\Vert f\right\Vert _{W_{Q_{2}}^{k}\left(L_{v}^{p}\right)}.
\]
for all measurable $f:\R^{\dimension}\to\Compl^{k}$.
\end{lem}
\begin{rem*}

\begin{itemize}[leftmargin=0.4cm]
\item As for the previous lemma, the statement of the lemma remains true
for $v_{0}$ instead of $v$.
\item Note that if $Q_{1},Q_{2}\subset\R^{\dimension}$ are two (Borel measurable)
bounded unit-neighborhoods, compactness of $\overline{Q_{1}}$ yields
finitely many $x_{1},\dots,x_{N}\in\R^{\dimension}$ satisfying $Q_{1}\subset\overline{Q_{1}}\subset\bigcup_{i=1}^{N}\left(x_{i}+Q_{2}^{\circ}\right)\subset\bigcup_{i=1}^{N}\left(x_{i}+Q_{2}\right)$,
so that the preceding lemma yields $\left\Vert f\right\Vert _{W_{Q_{1}}^{k}\left(L_{v}^{p}\right)}\lesssim\left\Vert f\right\Vert _{W_{Q_{2}}^{k}\left(L_{v}^{p}\right)}$,
where the implied constant is independent of $f$. By symmetry, this
argument shows $W_{Q_{1}}^{k}\left(L_{v}^{p}\right)=W_{Q_{2}}^{k}\left(L_{v}^{p}\right)$,
with equivalent (quasi)-norms. But since the constants of the (quasi)-norm
equivalence depend heavily on $Q_{1},Q_{2}$, this statement is not
of too much value for us.\qedhere
\end{itemize}
\end{rem*}
\begin{proof}
We have for any measurable $f:\R^{\dimension}\to\Compl^{k}$ that
\begin{align*}
\left(M_{Q_{1}}f\right)\left(x\right)=\left\Vert f\cdot\Indicator_{x+Q_{1}}\right\Vert _{L^{\infty}} & \leq\left\Vert f\cdot\Indicator_{x+\bigcup_{i=1}^{N}\left(x_{i}+Q_{2}\right)}\right\Vert _{L^{\infty}}\\
 & \leq\sum_{i=1}^{N}\left\Vert f\cdot\Indicator_{x+x_{i}+Q_{2}}\right\Vert _{L^{\infty}}\\
 & =\sum_{i=1}^{N}\left(M_{Q_{2}}f\right)\left(x+x_{i}\right)\\
 & =\sum_{i=1}^{N}\left(L_{-x_{i}}\left[M_{Q_{2}}f\right]\right)\left(x\right).
\end{align*}
For $p\geq1$, we can thus use the triangle inequality for $L^{p}$
and the estimate for $\vertiii{L_{y}}$ from Lemma \ref{lem:WeightedLpTranslationNorm},
as well as solidity of $L^{p}$ to derive
\[
\left\Vert f\right\Vert _{W_{Q_{1}}^{k}\left(L_{v}^{p}\right)}=\left\Vert M_{Q_{1}}f\right\Vert _{L_{v}^{p}}\leq\sum_{i=1}^{N}\left\Vert L_{-x_{i}}\left(M_{Q_{2}}f\right)\right\Vert _{L_{v}^{p}}\leq\left[\sum_{i=1}^{N}v_{0}\left(-x_{i}\right)\right]\cdot\left\Vert M_{Q_{2}}f\right\Vert _{L_{v}^{p}}=\left[\sum_{i=1}^{N}v_{0}\left(-x_{i}\right)\right]\cdot\left\Vert f\right\Vert _{W_{Q_{2}}^{k}\left(L_{v}^{p}\right)}.
\]
Similarly, for $p\in\left(0,1\right)$, we use the $p$-triangle inequality
(i.e., $\left\Vert \sum_{i=1}^{n}f_{i}\right\Vert _{L^{p}}^{p}\leq\sum_{i=1}^{n}\left\Vert f_{i}\right\Vert _{L^{p}}^{p}$)
to derive
\[
\left\Vert f\right\Vert _{W_{Q_{1}}^{k}\left(L_{v}^{p}\right)}^{p}=\left\Vert M_{Q_{1}}f\right\Vert _{L_{v}^{p}}^{p}\leq\sum_{i=1}^{N}\left\Vert L_{-x_{i}}\left(M_{Q_{2}}f\right)\right\Vert _{L_{v}^{p}}^{p}\leq\left[\sum_{i=1}^{N}\left(v_{0}\left(-x_{i}\right)\right)^{p}\right]\cdot\left\Vert M_{Q_{2}}f\right\Vert _{L_{v}^{p}}^{p}=\left[\sum_{i=1}^{N}\left(v_{0}\left(-x_{i}\right)\right)^{p}\right]\cdot\left\Vert f\right\Vert _{W_{Q_{2}}^{k}\left(L_{v}^{p}\right)}^{p}.\qedhere
\]
\end{proof}
In view of the preceding lemma, our next result becomes relevant:
\begin{lem}
\label{lem:BallCoveringNumber}Let $\left\Vert \cdot\right\Vert $
be any norm on $\R^{\dimension}$ and let $R>0$. For any $r>0$ and
$N:=\left\lfloor \left(3+2r\right)^{\dimension}\right\rfloor $, there
are $x_{1},\dots,x_{N}\in B_{\left(1+r\right)R}^{\left\Vert \cdot\right\Vert }\left(0\right)$
satisfying
\[
B_{\left(1+r\right)R}^{\left\Vert \cdot\right\Vert }\left(0\right)\subset\bigcup_{i=1}^{N}\left[x_{i}+B_{R}^{\left\Vert \cdot\right\Vert }\left(0\right)\right],
\]
where $B_{s}^{\left\Vert \cdot\right\Vert }\left(0\right)=\left\{ x\in\R^{\dimension}\with\left\Vert x\right\Vert <s\right\} $.
\end{lem}
\begin{proof}
First of all, assume we are given $x_{1},\dots,x_{M}\in B_{\left(1+r\right)R}^{\left\Vert \cdot\right\Vert }\left(0\right)$
such that $\left(x_{i}+B_{R/2}^{\left\Vert \cdot\right\Vert }\left(0\right)\right)_{i\in\underline{M}}$
is pairwise disjoint. Because of $x_{i}\in B_{\left(1+r\right)R}^{\left\Vert \cdot\right\Vert }\left(0\right)$,
we have $x_{i}+B_{R/2}^{\left\Vert \cdot\right\Vert }\left(0\right)\subset B_{\left(\frac{3}{2}+r\right)R}^{\left\Vert \cdot\right\Vert }\left(0\right)$,
so that additivity and translation invariance of the Lebesgue-measure
yields
\begin{align*}
M\cdot\left(R/2\right)^{\dimension}\cdot\lambda\left(B_{1}^{\left\Vert \cdot\right\Vert }\left(0\right)\right) & =\sum_{i=1}^{M}\lambda\left(x_{i}+B_{R/2}^{\left\Vert \cdot\right\Vert }\left(0\right)\right)\\
 & =\lambda\left(\biguplus_{i=1}^{M}x_{i}+B_{R/2}^{\left\Vert \cdot\right\Vert }\left(0\right)\right)\\
 & \leq\lambda\left(B_{\left(\frac{3}{2}+r\right)R}^{\left\Vert \cdot\right\Vert }\left(0\right)\right)\\
 & =\left[\left(\frac{3}{2}+r\right)R\right]^{\dimension}\cdot\lambda\left(B_{1}^{\left\Vert \cdot\right\Vert }\left(0\right)\right)
\end{align*}
and thus $M\leq\left(3+2r\right)^{\dimension}$. Since $M\in\N$,
we even get $M\leq N$. In particular, there can be at most a finite
number of such $x_{i}$.

Now (e.g.\@ using Zorn's Lemma), we can find a \emph{maximal} family
$\left(x_{i}\right)_{i\in\underline{M}}$ in $B_{\left(1+r\right)R}^{\left\Vert \cdot\right\Vert }\left(0\right)$
such that the family of sets $\left(x_{i}+B_{R/2}^{\left\Vert \cdot\right\Vert }\left(0\right)\right)_{i\in\underline{M}}$
is pairwise disjoint. As seen above, $M\leq N$.

It remains to show $B_{\left(1+r\right)R}^{\left\Vert \cdot\right\Vert }\left(0\right)\subset\bigcup_{i=1}^{M}\left[x_{i}+B_{R}^{\left\Vert \cdot\right\Vert }\left(0\right)\right]=:\Gamma$.
Thus, let $x\in B_{\left(1+r\right)R}^{\left\Vert \cdot\right\Vert }\left(0\right)$
be arbitrary. In case of $x\in\left\{ x_{1},\dots,x_{M}\right\} $,
we clearly have $x\in\Gamma$. But for $x\notin\left\{ x_{1},\dots,x_{M}\right\} $,
we see by maximality of the family $\left(x_{i}\right)_{i\in\underline{M}}$
that there is some $i\in\underline{M}$ satisfying $\left(x+B_{R/2}^{\left\Vert \cdot\right\Vert }\left(0\right)\right)\cap\left(x_{i}+B_{R/2}^{\left\Vert \cdot\right\Vert }\left(0\right)\right)\neq\emptyset$.
But this easily yields $x\in x_{i}+B_{R}^{\left\Vert \cdot\right\Vert }\left(0\right)\subset\Gamma$,
as desired.
\end{proof}
As announced above, we can now derive a completely quantitative version
of the (quasi)-norm equivalence between $W_{T_{j}^{-T}\left[-R,R\right]^{\dimension}}\left(L_{v}^{p}\right)$
and $W_{T_{i}^{-T}\left[-L,L\right]^{\dimension}}\left(L_{v}^{p}\right)$.
\begin{cor}
\label{cor:WienerLinearCubeNormEstimate}Let $i,j\in I$ and $p\in\left(0,\infty\right]$,
let $R,L\in\left[1,\infty\right)$ and let $f:\R^{\dimension}\to\Compl^{k}$
be measurable. Then we have
\[
\left\Vert f\right\Vert _{W_{T_{j}^{-T}\left[-R,R\right]^{\dimension}}^{k}\left(L_{v}^{p}\right)}\leq\Omega_{0}^{K}\Omega_{1}\cdot\left(3\dimension\left(L+R\right)\right)^{K+\dimension\cdot\max\left\{ 1,\frac{1}{p}\right\} }\cdot\left(1+\left\Vert T_{j}^{-1}T_{i}\right\Vert \right)^{K+\dimension\cdot\max\left\{ 1,\frac{1}{p}\right\} }\cdot\left\Vert f\right\Vert _{W_{T_{i}^{-T}\left[-L,L\right]^{\dimension}}^{k}\left(L_{v}^{p}\right)}.\qedhere
\]
\end{cor}
\begin{rem*}
As for the preceding results, the statement of the corollary remains
valid if $v$ is replaced by $v_{0}$.

Finally, we explicitly state the two most important special cases
of the preceding corollary:

\begin{itemize}
\item We have $i=j$ and $R=2$, as well as $L=1$. In this case, the corollary
yields
\begin{equation}
\left\Vert f\right\Vert _{W_{T_{i}^{-T}\left[-2,2\right]^{\dimension}}^{k}\left(L_{v}^{p}\right)}\leq\Omega_{0}^{K}\Omega_{1}\cdot\left(18\dimension\right)^{K+\dimension\cdot\max\left\{ 1,\frac{1}{p}\right\} }\cdot\left\Vert f\right\Vert _{W_{T_{i}^{-T}\left[-1,1\right]^{\dimension}}^{k}\left(L_{v}^{p}\right)}.\label{eq:WienerLinearCubeEnlargement}
\end{equation}
\item We have $R=L=1$. In this case, the corollary yields
\begin{equation}
\left\Vert f\right\Vert _{W_{T_{j}^{-T}\left[-1,1\right]^{\dimension}}^{k}\left(L_{v}^{p}\right)}\leq\Omega_{0}^{K}\Omega_{1}\cdot\left(6\dimension\right)^{K+\dimension\cdot\max\left\{ 1,\frac{1}{p}\right\} }\cdot\left(1+\left\Vert T_{j}^{-1}T_{i}\right\Vert \right)^{K+\dimension\cdot\max\left\{ 1,\frac{1}{p}\right\} }\cdot\left\Vert f\right\Vert _{W_{T_{i}^{-T}\left[-1,1\right]^{\dimension}}^{k}\left(L_{v}^{p}\right)}.\qedhere\label{eq:WienerLinearCubeTransformationChange}
\end{equation}
\end{itemize}
\end{rem*}
\begin{proof}
For brevity, set $R':=\left\Vert \left(T_{i}^{-T}\right)^{-1}T_{j}^{-T}\right\Vert _{\ell^{\infty}\to\ell^{\infty}}\cdot R=\left\Vert T_{j}^{-1}T_{i}\right\Vert _{\ell^{1}\to\ell^{1}}\cdot R$
and note
\[
\left(T_{i}^{-T}\right)^{-1}T_{j}^{-T}\left[-R,R\right]^{\dimension}\subset\left[-R',R'\right]^{\dimension}=\overline{B_{R'}^{\left\Vert \cdot\right\Vert _{\infty}}}\left(0\right)\subset\overline{B_{\left(1+\frac{R'}{L}\right)L}^{\left\Vert \cdot\right\Vert _{\infty}}}\left(0\right).
\]
But Lemma \ref{lem:BallCoveringNumber} yields certain $x_{1},\dots,x_{N}\in\overline{B_{\left(1+\frac{R'}{L}\right)L}^{\left\Vert \cdot\right\Vert _{\infty}}}\left(0\right)$
satisfying
\[
\overline{B_{\left(1+\frac{R'}{L}\right)L}^{\left\Vert \cdot\right\Vert _{\infty}}}\left(0\right)\subset\bigcup_{i=1}^{N}\left(x_{i}+\overline{B_{L}^{\left\Vert \cdot\right\Vert _{\infty}}}\left(0\right)\right)=\bigcup_{i=1}^{N}\left(x_{i}+\left[-L,L\right]^{\dimension}\right),
\]
where 
\[
N\leq\left(3+2\frac{R'}{L}\right)^{\dimension}=\left(3+2\frac{R}{L}\left\Vert T_{j}^{-1}T_{i}\right\Vert _{\ell^{1}\to\ell^{1}}\right)^{\dimension}\leq3^{\dimension}\left(1+\frac{R}{L}\right)^{\dimension}\left(1+\left\Vert T_{j}^{-1}T_{i}\right\Vert _{\ell^{1}\to\ell^{1}}\right)^{\dimension}.
\]
Next, note $\left\Vert T_{j}^{-1}T_{i}\right\Vert _{\ell^{1}\to\ell^{1}}\leq\sqrt{\dimension}\cdot\left\Vert T_{j}^{-1}T_{i}\right\Vert $
and hence $N\leq\left(3\sqrt{\dimension}\right)^{\dimension}\left(1+\frac{R}{L}\right)^{\dimension}\left(1+\left\Vert T_{j}^{-1}T_{i}\right\Vert \right)^{\dimension}$.

By putting together what we derived above, we get
\[
T_{j}^{-T}\left[-R,R\right]^{\dimension}\subset\bigcup_{i=1}^{N}\left(T_{i}^{-T}x_{i}+T_{i}^{-T}\left[-L,L\right]^{\dimension}\right).
\]
Next, we set $s:=\min\left\{ 1,p\right\} $, note 
\begin{align*}
\left\Vert x_{i}\right\Vert _{\infty} & \leq\left(1+\frac{R'}{L}\right)L=\left(L+R'\right)\\
 & =\left(L+R\cdot\left\Vert T_{j}^{-1}T_{i}\right\Vert _{\ell^{1}\to\ell^{1}}\right)\\
 & \leq\left(L+R\right)\left(1+\left\Vert T_{j}^{-1}T_{i}\right\Vert _{\ell^{1}\to\ell^{1}}\right)\\
 & \leq\sqrt{\dimension}\cdot\left(L+R\right)\left(1+\left\Vert T_{j}^{-1}T_{i}\right\Vert \right)
\end{align*}
and thus $1+\left|x_{i}\right|\leq1+\dimension\cdot\left(L+R\right)\left(1+\left\Vert T_{j}^{-1}T_{i}\right\Vert \right)\leq\dimension\cdot\left(1+L+R\right)\left(1+\left\Vert T_{j}^{-1}T_{i}\right\Vert \right)$
and recall equation (\ref{eq:WeightLinearTransformationsConnection})
to derive
\begin{align*}
\left[\sum_{i=1}^{N}\left(1+\left|T_{i}^{-T}x_{i}\right|\right)^{sK}\right]^{1/s} & \leq\Omega_{0}^{K}\cdot\left[\sum_{i=1}^{N}\left(1+\left|x_{i}\right|\right)^{sK}\right]^{1/s}\\
 & \leq\left[\dimension\Omega_{0}\left(1+L+R\right)\left(1+\left\Vert T_{j}^{-1}T_{i}\right\Vert \right)\right]^{K}\cdot N^{1/s}\\
\left({\scriptstyle \text{since }R,L\geq1}\right) & \leq3^{K+\frac{\dimension}{s}}\dimension^{K+\frac{\dimension}{2s}}\left(L+R\right)^{K+\frac{\dimension}{s}}\cdot\Omega_{0}^{K}\cdot\left(1+\left\Vert T_{j}^{-1}T_{i}\right\Vert \right)^{K+\frac{\dimension}{s}}\\
 & \leq\Omega_{0}^{K}\cdot\left(3\dimension\left(L+R\right)\right)^{K+\frac{\dimension}{s}}\cdot\left(1+\left\Vert T_{j}^{-1}T_{i}\right\Vert \right)^{K+\frac{\dimension}{s}}.
\end{align*}

All in all, Lemma \ref{lem:WienerAmalgamNormEquivalence} implies
\[
\left\Vert f\right\Vert _{W_{T_{j}^{-T}\left[-R,R\right]^{\dimension}}^{k}\left(L_{v}^{p}\right)}\leq\Omega_{0}^{K}\Omega_{1}\cdot\left(3\dimension\left(L+R\right)\right)^{K+\frac{\dimension}{s}}\cdot\left(1+\left\Vert T_{j}^{-1}T_{i}\right\Vert \right)^{K+\frac{\dimension}{s}}\cdot\left\Vert f\right\Vert _{W_{T_{i}^{-T}\left[-L,L\right]^{\dimension}}^{k}\left(L_{v}^{p}\right)},
\]
as desired.
\end{proof}

\subsection{The oscillation of a function}

Later in the paper, we will need to discretize certain reproducing
formulas involving convolutions. As observed in \cite{FeichtingerCoorbit0,FeichtingerCoorbit1,FeichtingerCoorbit2,GroechenigDescribingFunctions,GroechenigNonuniformSampling2},
a central tool for these discretizations is the oscillation of a function
and certain properties of and estimates for it. The goal of this subsection
is to collect these properties and estimates.
\begin{defn}
\label{def:Oscillation}Let $f:\R^{\dimension}\to\Compl^{k}$ and
let $\emptyset\neq Q\subset\R^{d}$. We define the \textbf{$Q$-oscillation}
of $f$ by
\[
\osc Qf:\R^{\dimension}\to\left[0,\infty\right],x\mapsto\sup_{y,z\in x+Q}\left|f\left(y\right)-f\left(z\right)\right|=\sup_{a,b\in Q}\left|f\left(x+a\right)-f\left(x+b\right)\right|.\qedhere
\]
\end{defn}
\begin{rem*}
Note that if $f$ is continuous, then so is $x\mapsto\left|f\left(x+a\right)-f\left(x+b\right)\right|$,
so that $\osc Qf$ is lower semicontinuous and hence measurable.
\end{rem*}
As our first step, we investigate some elementary properties of the
oscillation, in particular the behaviour of the oscillation under
a linear change of variables and under convolution.
\begin{lem}
\label{lem:OscillationLinearChange}Let $f:\R^{\dimension}\to\Compl^{k}$,
$T\in\GL\left(\R^{\dimension}\right)$ and let $\emptyset\neq Q\subset\R^{\dimension}$.
Then
\[
\osc Q\left(f\circ T\right)=\left(\osc{TQ}f\right)\circ T.\qedhere
\]
\end{lem}
\begin{proof}
We have
\begin{align*}
\left[\osc Q\left(f\circ T\right)\right]\left(x\right) & =\sup_{a,b\in Q}\left|\left(f\circ T\right)\left(x+a\right)-\left(f\circ T\right)\left(x+b\right)\right|\\
 & =\sup_{\alpha,\beta\in TQ}\left|f\left(\alpha+Tx\right)-f\left(\beta+Tx\right)\right|\\
 & =\left(\osc{TQ}f\right)\left(Tx\right).\qedhere
\end{align*}
\end{proof}
\begin{lem}
\label{lem:OscillationConvolution}Let $f,g:\R^{\dimension}\to\Compl$
be measurable, let $\emptyset\neq Q\subset\R^{\dimension}$ and assume
that $\osc Qf$ is measurable and that $\left(\left|f\right|\ast\left|g\right|\right)\left(x\right)<\infty$
for every $x\in\R^{\dimension}$. Then
\[
\left[\osc Q\left(f\ast g\right)\right]\left(x\right)\leq\left[\left(\osc Qf\right)\ast\left|g\right|\right]\left(x\right)\qquad\forall x\in\R^{\dimension}.\qedhere
\]
\end{lem}
\begin{proof}
Let $x\in\R^{\dimension}$ and fix $a,b\in Q$. Then
\begin{align*}
\left|\left(f\ast g\right)\left(x+a\right)-\left(f\ast g\right)\left(x+b\right)\right| & =\left|\int_{\R^{\dimension}}f\left(\left(x+a\right)-y\right)g\left(y\right)\d y-\int_{\R^{\dimension}}f\left(\left(x+b\right)-y\right)g\left(y\right)\d y\right|\\
 & \leq\int_{\R^{\dimension}}\left|f\left(\left(x-y\right)+a\right)-f\left(\left(x-y\right)+b\right)\right|\cdot\left|g\left(y\right)\right|\d y\\
 & \leq\int_{\R^{\dimension}}\left(\osc Qf\right)\left(x-y\right)\cdot\left|g\left(y\right)\right|\d y\\
 & =\left[\left(\osc Qf\right)\ast\left|g\right|\right]\left(x\right),
\end{align*}
as claimed.
\end{proof}
Intuitively, it should be true that smooth functions have a small
oscillation if their derivative is small. The next two lemmas make
this precise:
\begin{lem}
\label{lem:OscillationEstimatedByWienerDerivative}Let $f\in C^{1}\left(\R^{\dimension};\Compl\right)$.
Then we have for every bounded, convex set $Q\subset\R^{\dimension}$
with nonempty interior that
\[
\osc Qf\leq{\rm diam}\left(Q\right)\cdot M_{Q}\left(\nabla f\right).\qedhere
\]
\end{lem}
\begin{proof}
For $x\in\R^{\dimension}$ and $a,b\in Q$, the fundamental theorem
of calculus yields
\begin{align*}
\left|f\left(x+b\right)-f\left(x+a\right)\right| & =\left|\int_{0}^{1}\frac{\d}{\d t}\bigg|_{t=s}f\left(x+a+t\left(b-a\right)\right)\d s\right|\\
 & =\left|\int_{0}^{1}\left\langle \nabla f\left(x+a+s\left(b-a\right)\right),\,b-a\right\rangle \d s\right|\\
 & \leq{\rm diam}\left(Q\right)\cdot\sup_{s\in\left[0,1\right]}\left|\nabla f\left(x+sb+\left(1-s\right)a\right)\right|\\
\left({\scriptstyle Q\text{ convex}}\right) & \leq{\rm diam}\left(Q\right)\cdot\sup_{c\in Q}\left|\nabla f\left(x+c\right)\right|\\
 & \overset{\left(\dagger\right)}{\leq}{\rm diam}\left(Q\right)\cdot\left[M_{Q}\left(\nabla f\right)\right]\left(x\right).
\end{align*}
Here, it only remains to justify the last step, i.e.\@ that $\left[M_{Q}\left(\nabla f\right)\right]\left(x\right)=\essup_{c\in Q}\left|\nabla f\left(x+c\right)\right|\overset{!}{=}\sup_{c\in Q}\left|\nabla f\left(x+c\right)\right|$.
Here, only ``$\geq$'' is nontrivial. But by continuity of $\nabla f$,
and since nonempty open sets have positive measure, it is not hard
to see
\[
\essup_{c\in Q}\left|\nabla f\left(x+c\right)\right|\geq\sup_{c\in Q^{\circ}}\left|\nabla f\left(x+c\right)\right|,
\]
so that it suffices (by continuity) to show that $\overline{Q^{\circ}}\supset Q$.
But for arbitrary $a\in Q^{\circ}$ and $b\in Q$, we have $B_{\varepsilon}\left(a\right)\subset Q$
for some $\varepsilon>0$. For $t\in\left(0,1\right)$, this implies
\[
ta+\left(1-t\right)b\in t\cdot B_{\varepsilon}\left(a\right)+\left(1-t\right)b\subset Q
\]
and hence $ta+\left(1-t\right)b\in Q^{\circ}$. Because of $ta+\left(1-t\right)b\xrightarrow[t\to0]{}b$,
we conclude $b\in\overline{Q^{\circ}}$, as desired.
\end{proof}
\begin{lem}
\label{lem:OscillationSchwartzFunction}Let $f\in C^{1}\left(\R^{\dimension};\Compl\right)$
and $N\geq0$ and set $C:=\left(3\sqrt{\dimension}\right)^{N+1}$.
Then
\begin{equation}
\left(\osc{\delta\cdot\left[-1,1\right]^{\dimension}}f\right)\left(x\right)\leq C\cdot\left\Vert \nabla f\right\Vert _{N}\cdot\delta\cdot\left(1+\left|x\right|\right)^{-N}\qquad\forall x\in\R^{\dimension}\:\forall\delta\in\left(0,1\right],\label{eq:OscillationPointwiseEstimate}
\end{equation}
where $\left\Vert \nabla f\right\Vert _{N}:=\sup_{x\in\R^{\dimension}}\left(1+\left|x\right|\right)^{N}\left|\nabla f\left(x\right)\right|$.
\end{lem}
\begin{proof}
By Cauchy-Schwarz, we have ${\rm diam}\left(\delta\cdot\left[-1,1\right]^{\dimension}\right)\leq2\sqrt{\dimension}\cdot\delta$.
We can clearly assume $\left\Vert \nabla f\right\Vert _{N}<\infty$.
Since we also have $\delta\left[-1,1\right]^{\dimension}\subset\left[-1,1\right]^{\dimension}$,
we get from Lemmas \ref{lem:OscillationEstimatedByWienerDerivative}
and \ref{lem:SchwartzFunctionsAreWiener} that
\begin{align*}
\left(\osc{\delta\left[-1,1\right]^{\dimension}}f\right)\left(x\right) & \leq2\sqrt{\dimension}\cdot\delta\cdot\left[M_{\delta\left[-1,1\right]^{\dimension}}\left(\nabla f\right)\right]\left(x\right)\\
 & \leq2\sqrt{\dimension}\cdot\delta\cdot\left[M_{\left[-1,1\right]^{\dimension}}\left(\nabla f\right)\right]\left(x\right)\\
 & \leq2\sqrt{\dimension}\cdot\delta\cdot\left\Vert \nabla f\right\Vert _{N}\cdot\left[M_{\left[-1,1\right]^{\dimension}}\left(1+\left|\mybullet\right|\right)^{-N}\right]\left(x\right)\\
\left({\scriptstyle \text{Lemma }\ref{lem:SchwartzFunctionsAreWiener}}\right) & \leq2\sqrt{\dimension}\cdot\left(1+2\sqrt{\dimension}\right)^{N}\cdot\delta\cdot\left\Vert \nabla f\right\Vert _{N}\cdot\left(1+\left|x\right|\right)^{-N}\\
\left({\scriptstyle \text{since }1\leq\sqrt{\dimension}}\right) & \leq\left(3\sqrt{\dimension}\right)^{N+1}\cdot\delta\cdot\left\Vert \nabla f\right\Vert _{N}\cdot\left(1+\left|x\right|\right)^{-N}
\end{align*}
for all $x\in\R^{\dimension}$.
\end{proof}

\subsection{Self-improving properties for bandlimited functions}

Our next aim is to show that bandlimited $L^{p}$-functions automatically
belong to $W_{Q}\left(L^{p}\right)$. More precisely, we will show
for $\supp\widehat{f}\subset T_{i}\left[-R,R\right]^{\dimension}+\xi_{0}$
that
\[
\left\Vert f\right\Vert _{W_{T_{i}^{-T}\left[-1,1\right]^{\dimension}}\left(L^{p}\right)}\leq C_{\dimension,p,R}\cdot\left\Vert f\right\Vert _{L^{p}},
\]
which we call a ``self-improving property'', since we can improve
a simple $L^{p}$ estimate to a Wiener-amalgam estimate, at least
for suitably bandlimited functions. In fact, we will develop a \emph{weighted
version} of the preceding estimate.

All of our results in this section are based on the following convolution
relation for bandlimited $L^{p}$-functions, which we take from \cite[Theorem 3.4]{DecompositionEmbedding}.
We remark that this \emph{pointwise} estimate already appears \emph{in
the proof} of \cite[Proposition 1.5.1]{TriebelTheoryOfFunctionSpaces},
but is not stated explicitly as a theorem.
\begin{thm}
\label{thm:PointwiseQuasiBanachBandlimitedConvolution}Let $Q,\Omega\subset\R^{\dimension}$
be compact and let $p\in\left(0,1\right]$. Furthermore, let $\psi\in L^{1}\left(\R^{\dimension}\right)$
with $\supp\psi\subset Q$ and such that $\Fourier^{-1}\psi\in L^{p}\left(\R^{\dimension}\right)$.

For each $f\in L^{p}\left(\R^{\dimension}\right)\cap\Schwartz'\left(\R^{\dimension}\right)$
with $\supp\widehat{f}\subset\Omega$, we have $\Fourier^{-1}\left(\psi\cdot\widehat{f}\right)=\left(\Fourier^{-1}\psi\right)\ast f\in L^{p}\left(\R^{\dimension}\right)$
with
\[
\left(\left|\Fourier^{-1}\psi\right|\ast\left|f\right|\right)\left(x\right)\leq\left[\lambda_{\dimension}\left(Q-\Omega\right)\right]^{\frac{1}{p}-1}\cdot\left[\int_{\R^{\dimension}}\left|\left(\Fourier^{-1}\psi\right)\left(y\right)\right|^{p}\cdot\left|f\left(x-y\right)\right|^{p}\d y\right]^{1/p}
\]
for all $x\in\R^{\dimension}$.
\end{thm}
In order to ``circumvent'' the assumption $f\in L^{p}\left(\R^{\dimension}\right)$,
we also need the following approximation result, a proof of which
can be found in \cite[Lemma 3.2]{DecompositionEmbedding}, or in \cite[Theorem 1.4.1]{TriebelTheoryOfFunctionSpaces}.
In fact, the proof given in \cite{DecompositionEmbedding} is based
on that in \cite{TriebelTheoryOfFunctionSpaces}.
\begin{lem}
\label{lem:BandlimitedPointwiseApproximation}Let $\Omega\subset\R^{\dimension}$
be compact and assume $f\in\Schwartz'\left(\R^{\dimension}\right)$
with $\supp\widehat{f}\subset\Omega$. Then $f$ is given by (integration
against) a smooth function $g\in C^{\infty}\left(\R^{\dimension}\right)$
with polynomially bounded derivatives of all orders.

Furthermore, there is a sequence of Schwartz functions $\left(g_{n}\right)_{n\in\N}$
with the following properties:

\begin{enumerate}
\item $\left|g_{n}\left(x\right)\right|\leq\left|g\left(x\right)\right|$
for all $x\in\R^{\dimension}$,
\item $g_{n}\left(x\right)\xrightarrow[n\to\infty]{}g\left(x\right)$ for
all $x\in\R^{\dimension}$,
\item $\supp\widehat{g_{n}}\subset B_{1/n}\left(\Omega\right)$, where $B_{1/n}\left(\Omega\right)$
is the $\frac{1}{n}$-neighborhood of $\Omega$, given by
\[
B_{1/n}\left(\Omega\right)=\left\{ \xi\in\R^{\dimension}\with\dist\left(\xi,\Omega\right)<\frac{1}{n}\right\} .
\]
\end{enumerate}
In the following, we will identify the bandlimited distribution $f$
with its ``smooth version'' $g$.
\end{lem}
Using the two preceding results, we can now establish our first ``self-improving
property''.
\begin{thm}
\label{thm:BandlimitedWienerAmalgamSelfImproving}For $p\in\left(0,\infty\right]$,
$\xi_{0}\in\R^{\dimension}$ and $f\in\Schwartz'\left(\R^{\dimension}\right)$
with $\supp\widehat{f}\subset T_{i}\left[-R,R\right]^{\dimension}+\xi_{0}$
for some $i\in I$, we have
\[
\left\Vert f\right\Vert _{W_{T_{i}^{-T}\left[-1,1\right]^{\dimension}}\left(L_{v}^{p}\right)}\leq C\cdot\left\Vert f\right\Vert _{L_{v}^{p}}
\]
with $s:=\min\left\{ 1,p\right\} $ and $N:=\left\lceil K+\frac{\dimension+1}{s}\right\rceil $,
as well as
\[
C=2^{4\left(1+\frac{\dimension}{s}\right)}s_{\dimension}^{\frac{1}{s}}\left(192\cdot\dimension^{3/2}\cdot N\right)^{N+1}\cdot\Omega_{0}^{K}\Omega_{1}\cdot\left(1+R\right)^{\frac{\dimension}{s}}.\qedhere
\]
\end{thm}
\begin{rem*}
The only specific property of $v$ which is used in the proof is that
$v\left(x+y\right)\leq v\left(x\right)v_{0}\left(y\right)$. By submultiplicativity
of $v_{0}$, the same remains true when $v$ is replaced by $v_{0}$.
Hence, we also have $\left\Vert f\right\Vert _{W_{T_{i}^{-T}\left[-1,1\right]^{\dimension}}\left(L_{v_{0}}^{p}\right)}\leq C\cdot\left\Vert f\right\Vert _{L_{v_{0}}^{p}}$
for $f\in\Schwartz'\left(\R^{\dimension}\right)$ with $\supp\widehat{f}\subset T_{i}\left[-R,R\right]^{\dimension}+\xi_{0}$.
\end{rem*}
\begin{proof}
We can clearly assume $\left\Vert f\right\Vert _{L_{v}^{p}}<\infty$,
since otherwise the claim is trivial. Using Lemma \ref{lem:BandlimitedPointwiseApproximation},
choose a sequence of Schwartz functions $\left(f_{n}\right)_{n\in\N}$
satisfying $\left|f_{n}\left(x\right)\right|\leq\left|f\left(x\right)\right|$,
as well as $f_{n}\left(x\right)\xrightarrow[n\to\infty]{}f\left(x\right)$
and furthermore $\supp\widehat{f_{n}}\subset B_{1/n}\left(T_{i}\left[-R,R\right]^{\dimension}+\xi_{0}\right)$
for all $n\in\N$. Note that $T_{i}\left(-\left(R+1\right),R+1\right)^{\dimension}+\xi_{0}$
is a neighborhood of the compact(!) set $T_{i}\left[-R,R\right]^{\dimension}+\xi_{0}$,
so that there is some $n_{0}\in\N$ satisfying $\supp\widehat{f_{n}}\subset T_{i}\left[-\left(R+1\right),R+1\right]^{\dimension}+\xi_{0}$
for all $n\geq n_{0}$. By dropping (or modifying) the first $n_{0}$
terms of the sequence $\left(f_{n}\right)_{n\in\N}$, we can assume
that this holds for all $n\in\N$.

Let $s=\min\left\{ 1,p\right\} $ and $N=\left\lceil K+\frac{\dimension+1}{s}\right\rceil \geq2$
as in the statement of the theorem. Using Lemma \ref{lem:SmoothCutOffFunctionConstants}
(with $R+1$ instead of $R$ and with $s=3$) and Corollary \ref{cor:CutoffInverseFourierEstimate},
we get a function $\psi\in\TestFunctionSpace{\R^{\dimension}}$ satisfying
$0\leq\psi\leq1$, as well as $\supp\psi\subset\left[-\left(R+4\right),R+4\right]^{\dimension}$
and $\psi\equiv1$ on $\left[-\left(R+1\right),R+1\right]^{\dimension}$
which also satisfies
\begin{align*}
\left|\left(\Fourier^{-1}\psi\right)\left(x\right)\right| & \leq2\pi\cdot2^{\dimension}\cdot\left(48\dimension\right)^{N+1}\left(N+2\right)!\cdot\left(R+4\right)^{\dimension}\cdot\left(1+\left|x\right|\right)^{-N}\\
 & =:C_{1}\cdot\left(1+\left|x\right|\right)^{-N}
\end{align*}
for all $x\in\R^{\dimension}$. Note because of $N\geq2$ that $\left(N+2\right)!=\prod_{\ell=1}^{N+2}\ell\leq\prod_{\ell=2}^{N+2}\left(N+2\right)=\left(N+2\right)^{N+1}\leq\left(2N\right)^{N+1}$
and hence
\[
C_{1}\leq2^{3\left(1+\dimension\right)}\left(96\dimension\cdot N\right)^{N+1}\cdot\left(1+R\right)^{\dimension}.
\]

\medskip{}

Next, define
\[
\varrho:=\Fourier^{-1}\left(L_{\xi_{0}}\left[\psi\circ T_{i}^{-1}\right]\right)=\left|\det T_{i}\right|\cdot M_{\xi_{0}}\left[\left(\Fourier^{-1}\psi\right)\circ T_{i}^{T}\right]
\]
and note $\widehat{\varrho}\equiv1$ on $T_{i}\left[-\left(R+1\right),R+1\right]^{\dimension}+\xi_{0}$,
so that $\widehat{f_{n}}=\widehat{\varrho}\cdot\widehat{f_{n}}$,
which implies $f_{n}=f_{n}\ast\varrho$. Next, we note for arbitrary
$x\in\R^{\dimension}$ and $y\in T_{i}^{-T}\left[-1,1\right]^{\dimension}$
because of
\begin{align*}
1+\left|T_{i}^{T}x\right| & \leq1+\left|T_{i}^{T}\left(x+y\right)\right|+\left|-T_{i}^{T}y\right|\\
 & \leq1+\sqrt{\dimension}+\left|T_{i}^{T}\left(x+y\right)\right|\\
 & \leq\left(1+\sqrt{\dimension}\right)\left(1+\left|T_{i}^{T}\left(x+y\right)\right|\right)
\end{align*}
that 
\begin{align}
\left|\varrho\left(x+y\right)\right| & =\left|\det T_{i}\right|\cdot\left|\left(\Fourier^{-1}\psi\right)\left(T_{i}^{T}\left(x+y\right)\right)\right|\nonumber \\
 & \leq C_{1}\cdot\left|\det T_{i}\right|\cdot\left(1+\left|T_{i}^{T}\left(x+y\right)\right|\right)^{-N}\nonumber \\
 & \leq\left(1+\sqrt{\dimension}\right)^{N}C_{1}\cdot\left|\det T_{i}\right|\cdot\left(1+\left|T_{i}^{T}x\right|\right)^{-N}.\label{eq:BandlimitedWienerSelfImprovementKernelEstimate}
\end{align}
Now, we distinguish the two cases $p\in\left[1,\infty\right]$ and
$p\in\left(0,1\right)$.

\textbf{Case 1}: We have $p\in\left[1,\infty\right]$ and hence $s=1$.
In this case, note because of
\begin{align}
v\left(x\right)=v\left(z+x-z\right) & \leq v\left(z\right)v_{0}\left(x-z\right)\nonumber \\
 & \leq\Omega_{1}\cdot v\left(z\right)\left(1+\left|x-z\right|\right)^{K}\nonumber \\
\left({\scriptstyle \text{eq. }\eqref{eq:WeightLinearTransformationsConnection}}\right) & \leq\Omega_{0}^{K}\Omega_{1}\cdot v\left(z\right)\cdot\left(1+\left|T_{i}^{T}\left(x-z\right)\right|\right)^{K}\label{eq:BandlimitedWienerSelfImprovementWeightProduct}
\end{align}
that
\begin{align*}
v\left(x\right)\cdot\left(M_{T_{i}^{-T}\left[-1,1\right]^{\dimension}}f_{n}\right)\left(x\right) & \leq v\left(x\right)\cdot\sup_{y\in T_{i}^{-T}\left[-1,1\right]^{\dimension}}\left|f_{n}\left(x+y\right)\right|\\
\left({\scriptstyle \text{since }f_{n}=f_{n}\ast\varrho}\right) & \leq v\left(x\right)\cdot\sup_{y\in T_{i}^{-T}\left[-1,1\right]^{\dimension}}\int_{\R^{\dimension}}\left|\varrho\left(\left(x-z\right)+y\right)\right|\cdot\left|f_{n}\left(z\right)\right|\d z\\
\left({\scriptstyle \text{eq. }\eqref{eq:BandlimitedWienerSelfImprovementKernelEstimate}\text{ with }x-z\text{ instead of }x}\right) & \leq\left(1+\sqrt{\dimension}\right)^{N}C_{1}\cdot\left|\det T_{i}\right|\cdot\int_{\R^{\dimension}}v\left(x\right)\cdot\left(1+\left|T_{i}^{T}\left(x-z\right)\right|\right)^{-N}\cdot\left|f_{n}\left(z\right)\right|\d z\\
\left({\scriptstyle \text{eq. }\eqref{eq:BandlimitedWienerSelfImprovementWeightProduct}}\right) & \leq\left(1+\sqrt{\dimension}\right)^{N}\Omega_{0}^{K}\Omega_{1}C_{1}\cdot\left|\det T_{i}\right|\cdot\int_{\R^{\dimension}}\left(1+\left|T_{i}^{T}\left(x-z\right)\right|\right)^{K-N}\cdot\left|\left(v\cdot f_{n}\right)\left(z\right)\right|\d z.
\end{align*}
Now, we simply use Young's inequality $\left\Vert f\ast g\right\Vert _{L^{p}}\leq\left\Vert f\right\Vert _{L^{1}}\cdot\left\Vert g\right\Vert _{L^{p}}$
to conclude
\begin{align*}
\left\Vert f_{n}\right\Vert _{W_{T_{i}^{-T}\left[-1,1\right]^{\dimension}}\left(L_{v}^{p}\right)} & =\left\Vert v\cdot M_{T_{i}^{-T}\left[-1,1\right]^{\dimension}}f_{n}\right\Vert _{L^{p}}\\
 & \leq\left(1+\sqrt{\dimension}\right)^{N}\Omega_{0}^{K}\Omega_{1}C_{1}\cdot\left|\det T_{i}\right|\cdot\left\Vert \left(1+\left|T_{i}^{T}\mybullet\right|\right)^{-\left(N-K\right)}\right\Vert _{L^{1}}\cdot\left\Vert v\cdot f_{n}\right\Vert _{L^{p}}\\
\left({\scriptstyle N-K\geq\dimension+1,\text{ eq. }\eqref{eq:StandardDecayLpEstimate}\text{ and }\left|f_{n}\right|\leq\left|f\right|}\right) & \leq\left(1+\sqrt{\dimension}\right)^{N}s_{\dimension}\Omega_{0}^{K}\Omega_{1}C_{1}\cdot\left\Vert f\right\Vert _{L_{v}^{p}}.
\end{align*}
But because of $f_{n}\to f$ pointwise, we get 
\[
\left(M_{T_{i}^{-T}\left[-1,1\right]^{\dimension}}f\right)\left(x\right)=\left\Vert \Indicator_{x+T_{i}^{-T}\left[-1,1\right]^{\dimension}}\cdot f\right\Vert _{L^{\infty}}\leq\liminf_{n\to\infty}\left\Vert \Indicator_{x+T_{i}^{-T}\left[-1,1\right]^{\dimension}}\cdot f_{n}\right\Vert _{L^{\infty}}=\liminf_{n\to\infty}\left(M_{T_{i}^{-T}\left[-1,1\right]^{\dimension}}f_{n}\right)\left(x\right),
\]
so that an application of Fatou's Lemma yields 
\[
\left\Vert f\right\Vert _{W_{T_{i}^{-T}\left[-1,1\right]^{\dimension}}\left(L_{v}^{p}\right)}\leq\liminf_{n\to\infty}\left\Vert f_{n}\right\Vert _{W_{T_{i}^{-T}\left[-1,1\right]^{\dimension}}\left(L_{v}^{p}\right)}\leq\left(1+\sqrt{\dimension}\right)^{N}s_{\dimension}\Omega_{0}^{K}\Omega_{1}C_{1}\cdot\left\Vert f\right\Vert _{L_{v}^{p}},
\]
as desired.

\medskip{}

\textbf{Case 2}: We have $p\in\left(0,1\right)$ and hence $s=p$.
In this case, we first note
\begin{align*}
\lambda_{\dimension}\left(\supp\widehat{f_{n}}-\supp\widehat{\varrho}\right) & \leq\lambda_{\dimension}\left(\left[T_{i}\left[-\left(R+1\right),R+1\right]^{\dimension}+\xi_{0}\right]-\left[T_{i}\left[-\left(R+4\right),R+4\right]^{\dimension}+\xi_{0}\right]\right)\\
 & \leq\lambda_{\dimension}\left(T_{i}\left[-\left(2R+5\right),2R+5\right]^{\dimension}\right)\\
 & =\left|\det T_{i}\right|\cdot\left(4R+10\right)^{\dimension}.
\end{align*}
For brevity, set $C_{2}:=\left(1+\sqrt{\dimension}\right)^{N}\left(4R+10\right)^{\dimension\left(\frac{1}{p}-1\right)}C_{1}$
and apply Theorem \ref{thm:PointwiseQuasiBanachBandlimitedConvolution}
to get
\begin{align*}
v\left(x\right)\cdot\left(M_{T_{i}^{-T}\left[-1,1\right]^{\dimension}}f_{n}\right)\left(x\right) & \leq v\left(x\right)\cdot\sup_{y\in T_{i}^{-T}\left[-1,1\right]^{\dimension}}\left|f_{n}\left(x+y\right)\right|\\
\left({\scriptstyle \text{since }f_{n}=f_{n}\ast\varrho}\right) & \leq v\left(x\right)\cdot\sup_{y\in T_{i}^{-T}\left[-1,1\right]^{\dimension}}\left(\left|f_{n}\right|\ast\left|\varrho\right|\right)\left(x+y\right)\\
\left({\scriptstyle \text{Theorem }\ref{thm:PointwiseQuasiBanachBandlimitedConvolution}}\right) & \leq\left[\left(4R\!+\!10\right)^{\dimension}\cdot\left|\det T_{i}\right|\right]^{\frac{1}{p}-1}\cdot v\left(x\right)\cdot\!\!\!\sup_{y\in T_{i}^{-T}\left[-1,1\right]^{\dimension}}\left(\int_{\R^{\dimension}}\left|\varrho\left(\left(x-z\right)+y\right)\right|^{p}\cdot\left|f_{n}\left(z\right)\right|^{p}\d z\right)^{1/p}\\
\left({\scriptstyle \text{eq. }\eqref{eq:BandlimitedWienerSelfImprovementKernelEstimate}}\right) & \leq C_{2}\cdot\left|\det T_{i}\right|^{1/p}\cdot v\left(x\right)\cdot\left(\int_{\R^{\dimension}}\left(1+\left|T_{i}^{T}\left(x-z\right)\right|\right)^{-Np}\cdot\left|f_{n}\left(z\right)\right|^{p}\d z\right)^{1/p}\\
\left({\scriptstyle \text{eq. }\eqref{eq:BandlimitedWienerSelfImprovementWeightProduct}}\right) & \leq C_{2}\Omega_{0}^{K}\Omega_{1}\cdot\left|\det T_{i}\right|^{1/p}\cdot\left(\int_{\R^{\dimension}}\left(1+\left|T_{i}^{T}\left(x-z\right)\right|\right)^{-p\left(N-K\right)}\cdot\left|\left(v\cdot f_{n}\right)\left(z\right)\right|^{p}\d z\right)^{1/p}.
\end{align*}
Finally, take the $L^{p}$-norm of the preceding estimate to conclude
\begin{align*}
\left\Vert f_{n}\right\Vert _{W_{T_{i}^{-T}\left[-1,1\right]^{\dimension}}\left(L_{v}^{p}\right)}^{p} & \leq\left(C_{2}\Omega_{0}^{K}\Omega_{1}\right)^{p}\cdot\left|\det T_{i}\right|\cdot\int_{\R^{\dimension}}\int_{\R^{\dimension}}\left(1+\left|T_{i}^{T}\left(x-z\right)\right|\right)^{-p\left(N-K\right)}\cdot\left|\left(v\cdot f_{n}\right)\left(z\right)\right|^{p}\d z\d x\\
\left({\scriptstyle \text{Fubini and }y=x-z}\right) & =\left(C_{2}\Omega_{0}^{K}\Omega_{1}\right)^{p}\cdot\left\Vert f\right\Vert _{L_{v}^{p}}^{p}\cdot\left|\det T_{i}\right|\cdot\int_{\R^{\dimension}}\left(1+\left|T_{i}^{T}y\right|\right)^{-p\left(N-K\right)}\d y\\
 & =\left(C_{2}\Omega_{0}^{K}\Omega_{1}\right)^{p}\cdot\left\Vert \left(1+\left|\mybullet\right|\right)^{-\left(N-K\right)}\right\Vert _{L^{p}}^{p}\cdot\left\Vert f\right\Vert _{L_{v}^{p}}^{p}\\
\left({\scriptstyle \text{eq. }\eqref{eq:StandardDecayLpEstimate}\text{ and }N-K-\frac{\dimension}{p}\geq\frac{1}{p}}\right) & \leq\left(C_{2}\Omega_{0}^{K}\Omega_{1}\right)^{p}\cdot s_{\dimension}\cdot\left\Vert f\right\Vert _{L_{v}^{p}}^{p}.
\end{align*}
The remainder of the proof is now as for $p\in\left[1,\infty\right]$,
but with a slightly different constant.
\end{proof}
Now, we establish our second ``self-improving property'', which
yields an estimate for the $L_{v}^{p}$-norm of the oscillation of
a band-limited function, only in terms of the $L_{v}^{p}$-norm of
the function.
\begin{thm}
\label{thm:BandlimitedOscillationSelfImproving}For each $p\in\left[0,\infty\right)$,
$i\in I$, $\xi_{0}\in\R^{\dimension}$, $R>0$ and $f\in\Schwartz'\left(\R^{\dimension}\right)$
with $\supp\widehat{f}\subset T_{i}\left[-R,R\right]^{\dimension}+\xi_{0}$,
we have
\[
\left\Vert \osc{\delta\cdot T_{i}^{-T}\left[-1,1\right]^{\dimension}}\left[M_{-\xi_{0}}f\right]\right\Vert _{W_{T_{i}^{-T}\left[-1,1\right]^{\dimension}}\left(L_{v}^{p}\right)}\leq C\cdot\delta\cdot\left\Vert f\right\Vert _{L_{v}^{p}}
\]
with
\[
C:=\Omega_{0}^{2K}\Omega_{1}^{2}\cdot\left(23040\cdot\dimension^{3/2}\cdot\left(K+1+\frac{\dimension+1}{\min\left\{ 1,p\right\} }\right)\right)^{K+2+\frac{\dimension+1}{\min\left\{ 1,p\right\} }}\cdot\left(1+R\right)^{1+\frac{\dimension}{\min\left\{ 1,p\right\} }}.\qedhere
\]
\end{thm}
\begin{rem*}
As usual, the claim remains valid when $v$ is replaced by $v_{0}$
throughout.
\end{rem*}
\begin{proof}
As usual, since $f$ is a \emph{bandlimited} tempered distribution,
it is actually given by integration against a smooth function with
polynomially bounded derivatives. Furthermore, $\supp\Fourier\left[M_{-\xi_{0}}f\right]=\supp L_{-\xi_{0}}\widehat{f}\subset T_{i}\left[-R,R\right]^{\dimension}$,
so that we can assume $\xi_{0}=0$ for the remainder of the proof.

The first part of the proof is now very similar to that of Theorem
\ref{thm:BandlimitedWienerAmalgamSelfImproving}: We can clearly assume
$\left\Vert f\right\Vert _{L_{v}^{p}}<\infty$, since otherwise the
claim is trivial. Using Lemma \ref{lem:BandlimitedPointwiseApproximation},
choose a sequence of Schwartz functions $\left(f_{n}\right)_{n\in\N}$
satisfying $\left|f_{n}\left(x\right)\right|\leq\left|f\left(x\right)\right|$,
as well as $f_{n}\left(x\right)\xrightarrow[n\to\infty]{}f\left(x\right)$
and furthermore $\supp\widehat{f_{n}}\subset B_{1/n}\left(T_{i}\left[-R,R\right]^{\dimension}\right)$
for all $n\in\N$. Note that $T_{i}\left(-\left(R+1\right),R+1\right)^{\dimension}$
is a neighborhood of the compact(!) set $T_{i}\left[-R,R\right]^{\dimension}$,
so that there is some $n_{0}\in\N$ satisfying $\supp\widehat{f_{n}}\subset T_{i}\left[-\left(R+1\right),R+1\right]^{\dimension}$
for all $n\geq n_{0}$. By dropping (or modifying) the first $n_{0}$
terms of the sequence $\left(f_{n}\right)_{n\in\N}$, we can assume
that this holds for all $n\in\N$.

Let $s:=\min\left\{ 1,p\right\} $ and $N:=\left\lceil K+\frac{\dimension+1}{s}\right\rceil \geq2$.
Using Lemma \ref{lem:SmoothCutOffFunctionConstants} (with $R+1$
instead of $R$ and with $s=3$) and Corollary \ref{cor:CutoffInverseFourierEstimate},
we get a function $\psi\in\TestFunctionSpace{\R^{\dimension}}$ satisfying
$0\leq\psi\leq1$, as well as $\supp\psi\subset\left[-\left(R+4\right),R+4\right]^{\dimension}$
and $\psi\equiv1$ on $\left[-\left(R+1\right),R+1\right]^{\dimension}$
which also satisfies
\begin{align}
\left|\left(\partial^{\alpha}\left[\Fourier^{-1}\psi\right]\right)\left(x\right)\right| & \leq2\pi\cdot2^{\dimension}\cdot\left(48\dimension\right)^{N+1}\left(N+2\right)!\cdot\left(R+5\right)^{\left|\alpha\right|}\left(R+4\right)^{\dimension}\cdot\left(1+\left|x\right|\right)^{-N}\nonumber \\
 & \leq2\pi\cdot2^{\dimension}\cdot\left(48\dimension\right)^{N+1}\left(N+2\right)!\cdot\left(R+5\right)^{\dimension+1}\cdot\left(1+\left|x\right|\right)^{-N}\nonumber \\
 & =:C_{1}\cdot\left(1+\left|x\right|\right)^{-N}\label{eq:SelfImprovingOscillationKernelDerivativeEstimate}
\end{align}
for all $x\in\R^{\dimension}$ and $\alpha\in\N_{0}^{\dimension}$
with $\left|\alpha\right|\leq1$. Using $\left(N+2\right)!=\prod_{\ell=1}^{N+2}\ell\leq\prod_{\ell=2}^{N+2}\left(N+2\right)=\left(N+2\right)^{N+1}\leq\left(2N\right)^{N+1}$,
we get
\begin{equation}
C_{1}\leq40\cdot10^{\dimension}\cdot\left(96\dimension\cdot N\right)^{N+1}\cdot\left(1+R\right)^{\dimension+1}.\label{eq:SelfImprovingOscillationC1Estimate}
\end{equation}

\medskip{}

Now, let $g_{n}:=f_{n}\circ T_{i}^{-T}$ for $n\in\N$. Using Lemmas
\ref{lem:OscillationLinearChange} and \ref{lem:OscillationEstimatedByWienerDerivative},
we see
\begin{align*}
\osc{\delta\cdot T_{i}^{-T}\left[-1,1\right]^{\dimension}}f_{n} & =\osc{\delta\cdot T_{i}^{-T}\left[-1,1\right]^{\dimension}}\left(g_{n}\circ T_{i}^{T}\right)\\
\left({\scriptstyle \text{Lemma }\ref{lem:OscillationLinearChange}}\right) & =\left(\osc{\delta\cdot\left[-1,1\right]^{\dimension}}g_{n}\right)\circ T_{i}^{T}\\
\left({\scriptstyle \text{Lemma }\ref{lem:OscillationEstimatedByWienerDerivative}}\right) & \leq2\sqrt{\dimension}\cdot\delta\cdot\left(M_{\delta\left[-1,1\right]^{\dimension}}\left[\nabla g_{n}\right]\right)\circ T_{i}^{T}.
\end{align*}
Based on this estimate, Lemmas \ref{lem:WienerTransformationFormula}
and \ref{lem:IteratedMaximalFunction} show
\begin{align}
M_{T_{i}^{-T}\left[-1,1\right]^{\dimension}}\left[\osc{\delta\cdot T_{i}^{-T}\left[-1,1\right]^{\dimension}}f_{n}\right] & \leq2\sqrt{\dimension}\cdot\delta\cdot M_{T_{i}^{-T}\left[-1,1\right]^{\dimension}}\left[\left(M_{\delta\left[-1,1\right]^{\dimension}}\left[\nabla g_{n}\right]\right)\circ T_{i}^{T}\right]\nonumber \\
\left({\scriptstyle \text{Lemma }\ref{lem:WienerTransformationFormula}}\right) & =2\sqrt{\dimension}\cdot\delta\cdot\left[M_{\left[-1,1\right]^{\dimension}}\left(M_{\delta\left[-1,1\right]^{\dimension}}\left[\nabla g_{n}\right]\right)\right]\circ T_{i}^{T}\nonumber \\
\left({\scriptstyle \text{Lemma }\ref{lem:IteratedMaximalFunction}\text{ and }\delta\leq1}\right) & \leq2\sqrt{\dimension}\cdot\delta\cdot\left(M_{\left[-2,2\right]^{\dimension}}\left[\nabla g_{n}\right]\right)\circ T_{i}^{T}\nonumber \\
\left({\scriptstyle \text{Lemma }\ref{lem:WienerTransformationFormula}}\right) & =2\sqrt{\dimension}\cdot\delta\cdot M_{T_{i}^{-T}\left[-2,2\right]^{\dimension}}\left[\left(\nabla g_{n}\right)\circ T_{i}^{T}\right].\label{eq:SelfImprovingOscillationPointwise}
\end{align}

Next, observe $\supp\widehat{g_{n}}=\supp\left[\left|\det T_{i}\right|\cdot\widehat{f_{n}}\circ T_{i}\right]=T_{i}^{-1}\supp\widehat{f_{n}}\subset\left[-\left(R+1\right),R+1\right]^{\dimension}$
for all $n\in\N$, so that we see $\widehat{g_{n}}=\widehat{g_{n}}\cdot\psi$.
Hence, $g_{n}=g_{n}\ast\Fourier^{-1}\psi$. Because of $g_{n},\Fourier^{-1}\psi\in\Schwartz\left(\R^{\dimension}\right)$,
this easily implies $\nabla g_{n}=g_{n}\ast\nabla\left(\Fourier^{-1}\psi\right)$.
But for arbitrary Schwartz functions $f,g$ and $T\in\GL\left(\R^{\dimension}\right)$,
we have
\begin{align*}
\left(f\ast g\right)\left(Tx\right) & =\int_{\R^{\dimension}}f\left(Tx-y\right)g\left(y\right)\d y\\
\left({\scriptstyle y=Tz}\right) & =\left|\det T\right|\cdot\int_{\R^{\dimension}}f\left(Tx-Tz\right)g\left(Tz\right)\d z\\
 & =\left|\det T\right|\cdot\left[\left(f\circ T\right)\ast\left(g\circ T\right)\right]\left(x\right),
\end{align*}
so that, if we understand the following equation componentwise,
\begin{align}
\left(\nabla g_{n}\right)\left(T_{i}^{T}x\right) & =\left|\det T_{i}\right|\cdot\left[\left(g_{n}\circ T_{i}^{T}\right)\ast\left(\left[\nabla\left(\Fourier^{-1}\psi\right)\right]\circ T_{i}^{T}\right)\right]\left(x\right)\nonumber \\
\left({\scriptstyle \text{with }\eta_{j}:=\partial_{j}\left(\Fourier^{-1}\psi\right)}\right) & =\left|\det T_{i}\right|\cdot\left[\left(f_{n}\ast\left[\eta_{j}\circ T_{i}^{T}\right]\right)\left(x\right)\right]_{j\in\underline{d}}.\label{eq:SelfImprovingOscillationGradientAsConvolution}
\end{align}

\medskip{}

Now, we divide the proof into the two cases $p\in\left[1,\infty\right]$
and $p\in\left(0,1\right)$. In the (easier) case $p\in\left[1,\infty\right]$,
we get
\begin{align*}
v\left(x\right)\cdot\left|\left(\partial_{j}g_{n}\right)\left(T_{i}^{T}x\right)\right| & =\left|\det T_{i}\right|\cdot v\left(x\right)\cdot\left|\left(f_{n}\ast\left[\eta_{j}\circ T_{i}^{T}\right]\right)\left(x\right)\right|\\
 & \leq\left|\det T_{i}\right|\cdot\int_{\R^{\dimension}}v\left(x\right)\cdot\left|f_{n}\left(y\right)\right|\cdot\left|\eta_{j}\left(T_{i}^{T}\left(x-y\right)\right)\right|\d y\\
\left({\scriptstyle \text{eq. }\eqref{eq:SelfImprovingOscillationKernelDerivativeEstimate}\text{ and }\eta_{j}=\partial_{j}\left[\Fourier^{-1}\psi\right]}\right) & \leq C_{1}\cdot\left|\det T_{i}\right|\cdot\int_{\R^{\dimension}}v\left(x\right)\cdot\left|f_{n}\left(y\right)\right|\cdot\left(1+\left|T_{i}^{T}\left(x-y\right)\right|\right)^{-N}\d y\\
\left({\scriptstyle \text{since }v\left(x\right)=v\left(x-y+y\right)\leq v\left(y\right)v_{0}\left(x-y\right)}\right) & \leq C_{1}\cdot\left|\det T_{i}\right|\cdot\int_{\R^{\dimension}}\left|\left(v\cdot f_{n}\right)\left(y\right)\right|\cdot v_{0}\left(x-y\right)\cdot\left(1+\left|T_{i}^{T}\left(x-y\right)\right|\right)^{-N}\d y\\
\left({\scriptstyle \text{assumption on }v_{0}\text{ and eq. }\eqref{eq:WeightLinearTransformationsConnection}}\right) & \leq\Omega_{0}^{K}\Omega_{1}C_{1}\cdot\left|\det T_{i}\right|\cdot\int_{\R^{\dimension}}\left|\left(v\cdot f_{n}\right)\left(y\right)\right|\cdot\left(1+\left|T_{i}^{T}\left(x-y\right)\right|\right)^{K-N}\d y
\end{align*}
for arbitrary $x\in\R^{\dimension}$.

But for $z\in T_{i}^{-T}\left[-2,2\right]^{\dimension}$, we have
\begin{align}
1+\left|T_{i}^{T}x\right| & \leq1+\left|T_{i}^{T}\left(x+z\right)\right|+\left|-T_{i}^{T}z\right|\nonumber \\
\left({\scriptstyle \left|T_{i}^{T}z\right|\leq2\sqrt{\dimension}\text{ since }T_{i}^{T}z\in\left[-2,2\right]^{\dimension}}\right) & \leq1+2\sqrt{\dimension}+\left|T_{i}^{T}\left(x+z\right)\right|\nonumber \\
 & \leq\left(1+2\sqrt{\dimension}\right)\left(1+\left|T_{i}^{T}\left(x+z\right)\right|\right)\label{eq:SelfImprovingOscillationChineseBracketShift}
\end{align}
and
\begin{align}
v\left(x\right) & =v\left(x+z-z\right)\leq v\left(x+z\right)\cdot v_{0}\left(-z\right)\nonumber \\
 & \leq\Omega_{1}\cdot v\left(x+z\right)\cdot\left(1+\left|-z\right|\right)^{K}\nonumber \\
\left({\scriptstyle \text{eq. }\eqref{eq:WeightLinearTransformationsConnection}}\right) & \leq\Omega_{0}^{K}\Omega_{1}\cdot v\left(x+z\right)\cdot\left(1+\left|T_{i}^{T}z\right|\right)^{K}\nonumber \\
 & \leq\Omega_{0}^{K}\Omega_{1}\left(1+2\sqrt{\dimension}\right)^{K}\cdot v\left(x+z\right).\label{eq:SelfImprovingOscillationSlightWeightShift}
\end{align}
By applying these two estimates and noting $K-N<0$, we get
\begin{align*}
v\left(x\right)\cdot\left|\left[\left(\partial_{j}g_{n}\right)\circ T_{i}^{T}\right]\left(x+z\right)\right| & \leq\Omega_{0}^{K}\Omega_{1}\cdot\left(1+2\sqrt{\dimension}\right)^{K}\cdot v\left(x+z\right)\cdot\left|\left[\left(\partial_{j}g_{n}\right)\circ T_{i}^{T}\right]\left(x+z\right)\right|\\
 & \leq\Omega_{0}^{2K}\Omega_{1}^{2}C_{1}\cdot\left(1+2\sqrt{\dimension}\right)^{K}\cdot\left|\det T_{i}\right|\cdot\int_{\R^{\dimension}}\left|\left(v\cdot f_{n}\right)\left(y\right)\right|\cdot\left(1+\left|T_{i}^{T}\left(x-y+z\right)\right|\right)^{K-N}\d y\\
\left({\scriptstyle \text{eq. }\eqref{eq:SelfImprovingOscillationChineseBracketShift}\text{ with }x-y\text{ instead of }x}\right) & \leq\Omega_{0}^{2K}\Omega_{1}^{2}C_{1}\cdot\left(1+2\sqrt{\dimension}\right)^{N}\cdot\left|\det T_{i}\right|\cdot\int_{\R^{\dimension}}\left|\left(v\cdot f_{n}\right)\left(y\right)\right|\cdot\left(1+\left|T_{i}^{T}\left(x-y\right)\right|\right)^{K-N}\d y.
\end{align*}
Noting that this holds for arbitrary $z\in T_{i}^{-T}\left[-2,2\right]^{\dimension}$
and by taking the $\ell^{2}$-norm over $j\in\underline{\dimension}$,
we conclude
\begin{align*}
v\left(x\right)\cdot\left[M_{T_{i}^{-T}\left[-2,2\right]^{\dimension}}\left(\left[\nabla g_{n}\right]\circ T_{i}^{T}\right)\right]\left(x\right) & \leq C_{2}\cdot\left|\det T_{i}\right|\cdot\int_{\R^{\dimension}}\left|\left(v\cdot f_{n}\right)\left(y\right)\right|\cdot\left(1+\left|T_{i}^{T}\left(x-y\right)\right|\right)^{K-N}\d y\\
 & =C_{2}\cdot\left|\det T_{i}\right|\cdot\left(\left|v\cdot f_{n}\right|\ast\left[\left(1+\left|\mybullet\right|\right)^{K-N}\circ T_{i}^{T}\right]\right)\left(x\right),
\end{align*}
where $C_{2}:=\Omega_{0}^{2K}\Omega_{1}^{2}\cdot C_{1}\sqrt{\dimension}\left(1+2\sqrt{\dimension}\right)^{N}$.

By taking the $L^{p}$-norm and using Young's theorem for convolutions,
we conclude
\begin{align*}
\left\Vert M_{T_{i}^{-T}\left[-2,2\right]^{\dimension}}\left(\left[\nabla g_{n}\right]\circ T_{i}^{T}\right)\right\Vert _{L_{v}^{p}} & \leq C_{2}\cdot\left|\det T_{i}\right|\cdot\left\Vert v\cdot f_{n}\right\Vert _{L^{p}}\cdot\left\Vert \left(1+\left|\mybullet\right|\right)^{K-N}\circ T_{i}^{T}\right\Vert _{L^{1}}\\
 & =C_{2}\cdot\left\Vert f_{n}\right\Vert _{L_{v}^{p}}\cdot\left\Vert \left(1+\left|\mybullet\right|\right)^{K-N}\right\Vert _{L^{1}}\\
\left({\scriptstyle \text{eq. }\eqref{eq:StandardDecayLpEstimate}\text{ and }K-N=K-\left\lceil K+\dimension+1\right\rceil \leq-\left(\dimension+1\right)}\right) & \leq s_{\dimension}C_{2}\cdot\left\Vert f_{n}\right\Vert _{L_{v}^{p}}\\
\left({\scriptstyle \text{since }\left|f_{n}\right|\leq\left|f\right|}\right) & \leq s_{\dimension}C_{2}\cdot\left\Vert f\right\Vert _{L_{v}^{p}}<\infty.
\end{align*}
In view of equation (\ref{eq:SelfImprovingOscillationPointwise}),
we have thus shown
\begin{align*}
\left\Vert \osc{\delta\cdot T_{i}^{-T}\left[-1,1\right]^{\dimension}}f_{n}\right\Vert _{W_{T_{i}^{-T}\left[-1,1\right]^{\dimension}}\left(L_{v}^{p}\right)} & \leq2\sqrt{\dimension}\cdot\delta\cdot\left\Vert M_{T_{i}^{-T}\left[-2,2\right]^{\dimension}}\left[\left(\nabla g_{n}\right)\circ T_{i}^{T}\right]\right\Vert _{L_{v}^{p}}\\
 & \leq2s_{\dimension}\sqrt{\dimension}C_{2}\cdot\delta\cdot\left\Vert f\right\Vert _{L_{v}^{p}}<\infty.
\end{align*}
Now, we note
\begin{align*}
2s_{\dimension}\sqrt{\dimension}C_{2} & =\Omega_{0}^{2K}\Omega_{1}^{2}\cdot C_{1}\cdot2s_{\dimension}\cdot\dimension\cdot\left(1+2\sqrt{\dimension}\right)^{N}\\
\left({\scriptstyle \text{eq. }\eqref{eq:SelfImprovingOscillationC1Estimate}}\right) & \leq\Omega_{0}^{2K}\Omega_{1}^{2}\cdot2s_{\dimension}\cdot\dimension\cdot\left(1+2\sqrt{\dimension}\right)^{N}\cdot40\cdot10^{\dimension}\cdot\left(96\dimension\cdot N\right)^{N+1}\cdot\left(1+R\right)^{\dimension+1}\\
\left({\scriptstyle \text{since }s_{\dimension}\leq4^{\dimension}\text{ and }\dimension\leq2^{\dimension}}\right) & \leq\Omega_{0}^{2K}\Omega_{1}^{2}\cdot2\cdot4^{\dimension}2^{\dimension}\left(1+2\sqrt{\dimension}\right)^{N}\cdot40\cdot10^{\dimension}\cdot\left(96\dimension\cdot N\right)^{N+1}\cdot\left(1+R\right)^{\dimension+1}\\
 & \leq\Omega_{0}^{2K}\Omega_{1}^{2}\cdot\left(1+2\sqrt{\dimension}\right)^{N}\cdot80^{\dimension+1}\cdot\left(96\dimension\cdot N\right)^{N+1}\cdot\left(1+R\right)^{\dimension+1}\\
 & \leq\Omega_{0}^{2K}\Omega_{1}^{2}\cdot\left(3\sqrt{\dimension}\right)^{N}\cdot80^{\dimension+1}\cdot\left(96\dimension\cdot N\right)^{N+1}\cdot\left(1+R\right)^{\dimension+1}\\
 & \leq\Omega_{0}^{2K}\Omega_{1}^{2}\cdot80^{\dimension+1}\cdot\left(288\cdot\dimension^{3/2}\cdot N\right)^{N+1}\cdot\left(1+R\right)^{\dimension+1}\\
\left({\scriptstyle \text{since }N\geq\dimension+1}\right) & \leq\Omega_{0}^{2K}\Omega_{1}^{2}\cdot\left(23040\cdot\dimension^{3/2}\cdot N\right)^{N+1}\cdot\left(1+R\right)^{\dimension+1}\\
\left({\scriptstyle \text{def. of }N}\right) & \leq\Omega_{0}^{2K}\Omega_{1}^{2}\cdot\left(23040\cdot\dimension^{3/2}\cdot\left(K+\dimension+2\right)\right)^{K+\dimension+3}\cdot\left(1+R\right)^{\dimension+1}.
\end{align*}

All that remains is to extend this estimate to $f$ instead of $f_{n}$.
But for $x\in\R^{\dimension}$ and arbitrary $y,z\in\delta\cdot T_{i}^{-T}\left[-1,1\right]^{\dimension}$,
we have
\[
\left|f\left(x+y\right)-f\left(x+z\right)\right|=\liminf_{n\to\infty}\left|f_{n}\left(x+y\right)-f_{n}\left(x+z\right)\right|\leq\liminf_{n\to\infty}\left(\osc{\delta\cdot T_{i}^{-T}\left[-1,1\right]^{\dimension}}f_{n}\right)\left(x\right)
\]
and thus $\left(\osc{\delta\cdot T_{i}^{-T}\left[-1,1\right]^{\dimension}}f\right)\left(x\right)\leq\liminf_{n\to\infty}\left(\osc{\delta\cdot T_{i}^{-T}\left[-1,1\right]^{\dimension}}f_{n}\right)\left(x\right)$
for all $x\in\R^{\dimension}$. A similar estimate holds for the maximal
function $M_{T_{i}^{-T}\left[-1,1\right]^{\dimension}}$ instead of
the oscillation. All in all, an application of Fatou's Lemma yields
\[
\left\Vert \osc{\delta\cdot T_{i}^{-T}\left[-1,1\right]^{\dimension}}f\right\Vert _{W_{T_{i}^{-T}\left[-1,1\right]^{\dimension}}\left(L_{v}^{p}\right)}\leq\liminf_{n\to\infty}\left\Vert \osc{\delta\cdot T_{i}^{-T}\left[-1,1\right]^{\dimension}}f_{n}\right\Vert _{W_{T_{i}^{-T}\left[-1,1\right]^{\dimension}}\left(L_{v}^{p}\right)}\leq2s_{\dimension}\sqrt{\dimension}C_{2}\cdot\delta\cdot\left\Vert f\right\Vert _{L_{v}^{p}}<\infty,
\]
as desired.

\medskip{}

Now, we consider the case $p\in\left(0,1\right)$. Here, we recall
$\eta_{j}=\partial_{j}\left(\Fourier^{-1}\psi\right)$ and observe
\begin{align*}
\supp\Fourier\left[\eta_{j}\circ T_{i}^{T}\right] & =\supp\left[\left|\det T_{i}\right|^{-1}\cdot\widehat{\eta_{j}}\circ T_{i}^{-1}\right]\\
 & =T_{i}\supp\widehat{\eta_{j}}\\
 & =T_{i}\supp\left(\xi\mapsto2\pi i\xi_{j}\cdot\Fourier\left[\Fourier^{-1}\psi\right]\left(\xi\right)\right)\\
 & \subset T_{i}\supp\psi\subset T_{i}\left[-\left(R+4\right),R+4\right]^{\dimension}
\end{align*}
and $\supp\widehat{f_{n}}\subset T_{i}\left[-\left(R+1\right),R+1\right]^{\dimension}$,
so that equation (\ref{eq:SelfImprovingOscillationGradientAsConvolution})
and Theorem \ref{thm:PointwiseQuasiBanachBandlimitedConvolution}
yield
\begin{align*}
 & \left|\left(\partial_{j}g_{n}\right)\left(T_{i}^{T}x\right)\right|^{p}\\
 & \leq\left|\det T_{i}\right|^{p}\cdot\left[\lambda_{\dimension}\left(T_{i}\left(\left[-\left(R\!+\!4\right),R\!+\!4\right]^{\dimension}-\left[-\left(R\!+\!1\right),R\!+\!1\right]^{\dimension}\right)\right)\right]^{1-p}\cdot\int_{\R^{\dimension}}\left|f_{n}\left(y\right)\right|^{p}\cdot\left|\left(\eta_{j}\circ T_{i}^{T}\right)\left(x-y\right)\right|^{p}\d y\\
\left({\scriptstyle \text{eq. }\eqref{eq:SelfImprovingOscillationKernelDerivativeEstimate}}\right) & \leq\left|\det T_{i}\right|\cdot C_{1}^{p}\left[2\cdot\left(2R+5\right)\right]^{\dimension\left(1-p\right)}\cdot\int_{\R^{\dimension}}\left|f_{n}\left(y\right)\right|^{p}\cdot\left(1+\left|T_{i}^{T}\left(x-y\right)\right|\right)^{-Np}\d y.
\end{align*}
Similar to the case $p\in\left[1,\infty\right]$, this implies
\begin{align*}
\left[v\left(x\right)\!\cdot\!\left|\left(\partial_{j}g_{n}\right)\left(T_{i}^{T}x\right)\right|\right]^{p} & \leq\left|\det T_{i}\right|\!\cdot C_{1}^{p}\left[2\cdot\left(2R+5\right)\right]^{\dimension\left(1-p\right)}\!\cdot\!\int_{\R^{\dimension}}\left|\left(v\!\cdot\!f_{n}\right)\left(y\right)\right|^{p}\!\cdot\!\left[v_{0}\left(x-y\right)\!\cdot\!\left(1\!+\!\left|T_{i}^{T}\!\left(x-y\right)\right|\right)^{-N}\right]^{p}\d y\\
\left({\scriptstyle \text{assump. on }v_{0}\text{ and eq. }\eqref{eq:WeightLinearTransformationsConnection}}\right) & \leq\left[\Omega_{0}^{K}\Omega_{1}\right]^{p}\cdot\left|\det T_{i}\right|\cdot C_{1}^{p}\left[4R+10\right]^{\dimension\left(1-p\right)}\!\cdot\!\int_{\R^{\dimension}}\left|\left(v\!\cdot\!f_{n}\right)\left(y\right)\right|^{p}\cdot\left[\left(1\!+\!\left|T_{i}^{T}\!\left(x-y\right)\right|\right)^{K-N}\right]^{p}\d y.
\end{align*}

As above, equations (\ref{eq:SelfImprovingOscillationChineseBracketShift})
and (\ref{eq:SelfImprovingOscillationSlightWeightShift}) show for
arbitrary $z\in T_{i}^{-T}\left[-2,2\right]^{\dimension}$ that
\begin{align*}
 & \left[v\left(x\right)\cdot\left|\left[\left(\partial_{j}g_{n}\right)\circ T_{i}^{T}\right]\left(x+z\right)\right|\right]^{p}\\
\left({\scriptstyle \text{eq. }\eqref{eq:SelfImprovingOscillationSlightWeightShift}}\right) & \leq\left[\Omega_{0}^{K}\Omega_{1}\left(1+2\sqrt{\dimension}\right)^{K}\right]^{p}\cdot\left|\left(v\cdot\left[\left(\partial_{j}g_{n}\right)\circ T_{i}^{T}\right]\right)\left(x+z\right)\right|^{p}\\
 & \leq\left[\Omega_{0}^{2K}\Omega_{1}^{2}\left(1+2\sqrt{\dimension}\right)^{K}\right]^{p}\!\!\cdot\!\left|\det T_{i}\right|\cdot C_{1}^{p}\left(10+4R\right)^{\dimension\left(1-p\right)}\!\cdot\!\int_{\R^{\dimension}}\left|\left(v\!\cdot\!f_{n}\right)\left(y\right)\right|^{p}\!\cdot\!\left[\left(1\!+\!\left|T_{i}^{T}\!\left(x+z-y\right)\right|\right)^{K-N}\right]^{p}\d y\\
\left({\scriptstyle \text{eq. }\eqref{eq:SelfImprovingOscillationChineseBracketShift}}\right) & \leq\left[\Omega_{0}^{2K}\Omega_{1}^{2}\left(1+2\sqrt{\dimension}\right)^{N}\right]^{p}\!\!\cdot\!\left|\det T_{i}\right|\cdot C_{1}^{p}\left(10+4R\right)^{\dimension\left(1-p\right)}\!\cdot\!\int_{\R^{\dimension}}\left|\left(v\!\cdot\!f_{n}\right)\left(y\right)\right|^{p}\!\cdot\!\left(1\!+\!\left|T_{i}^{T}\!\left(x-y\right)\right|\right)^{p\left(K-N\right)}\d y.
\end{align*}
Since this holds for arbitrary $z\in T_{i}^{-T}\left[-2,2\right]^{\dimension}$
and $j\in\underline{\dimension}$ and since
\[
\left|\left[\left(\nabla g_{n}\right)\circ T_{i}^{T}\right]\left(x+z\right)\right|\leq\sqrt{\dimension}\cdot\max_{j\in\underline{\dimension}}\left|\left[\left(\partial_{j}g_{n}\right)\circ T_{i}^{T}\right]\left(x+z\right)\right|,
\]
we conclude
\[
\left[v\!\cdot\!M_{T_{i}^{-T}\left[-2,2\right]^{\dimension}}\left(\left[\nabla g_{n}\right]\circ T_{i}^{T}\right)\right]^{p}\!\leq\!\left[\Omega_{0}^{2K}\Omega_{1}^{2}\sqrt{\dimension}\left(1\!+\!2\sqrt{\dimension}\right)^{\!N}\right]^{p}\!\cdot\left|\det T_{i}\right|\cdot C_{1}^{p}\left(10\!+\!4R\right)^{\dimension\left(1-p\right)}\cdot\left|v\!\cdot\!f_{n}\right|^{p}\ast\left[\left(1+\left|\mybullet\right|\right)^{p\left(K-N\right)}\circ T_{i}^{T}\right],
\]
so that Young's inequality $L^{1}\ast L^{1}\hookrightarrow L^{1}$
yields
\begin{align*}
\left\Vert \left[\nabla g_{n}\right]\!\circ\!T_{i}^{T}\right\Vert _{W_{T_{i}^{-T}\left[-2,2\right]^{\dimension}}\left(L_{v}^{p}\right)}^{p} & \leq\!\left[C_{1}\Omega_{0}^{2K}\Omega_{1}^{2}\sqrt{\dimension}\left(1\!+\!2\sqrt{\dimension}\right)^{\!\!N}\right]^{p}\!\!\left(10\!+\!4R\right)^{\dimension\left(1-p\right)}\!\cdot\!\left|\det T_{i}\right|\cdot\!\left\Vert \left|v\!\cdot\!f_{n}\right|^{p}\right\Vert _{L^{1}}\!\cdot\left\Vert \left(1\!+\!\left|\mybullet\right|\right)^{p\left(K-N\right)}\!\circ\!T_{i}^{T}\right\Vert _{L^{1}}\\
 & =\!\left[\Omega_{0}^{2K}\Omega_{1}^{2}\sqrt{\dimension}\left(1\!+\!2\sqrt{\dimension}\right)^{\!\!N}\right]^{p}\!\!\cdot C_{1}^{p}\left(10\!+\!4R\right)^{\dimension\left(1-p\right)}\cdot\left\Vert f_{n}\right\Vert _{L_{v}^{p}}^{p}\cdot\left\Vert \left(1+\left|\mybullet\right|\right)^{p\left(K-N\right)}\right\Vert _{L^{1}}\\
\left({\scriptstyle \text{eq. }\eqref{eq:StandardDecayLpEstimate}\text{ and }p\left(K-N\right)\leq-\left(\dimension+1\right)}\right) & \leq\!\left[\Omega_{0}^{2K}\Omega_{1}^{2}\sqrt{\dimension}\left(1\!+\!2\sqrt{\dimension}\right)^{\!\!N}\right]^{p}\!\!\cdot C_{1}^{p}\left(10\!+\!4R\right)^{\dimension\left(1-p\right)}\cdot\left\Vert f_{n}\right\Vert _{L_{v}^{p}}^{p}\cdot s_{\dimension}\\
\left({\scriptstyle \text{since }\left|f_{n}\right|\leq\left|f\right|\text{ and }s_{\dimension}\leq4^{\dimension}}\right) & \leq\!\left[\Omega_{0}^{2K}\Omega_{1}^{2}\sqrt{\dimension}\left(1\!+\!2\sqrt{\dimension}\right)^{\!\!N}4^{\dimension/p}C_{1}\left(10\!+\!4R\right)^{\dimension\left(\frac{1}{p}-1\right)}\right]^{p}\cdot\left\Vert f\right\Vert _{L_{v}^{p}}^{p}.
\end{align*}
In view of equation (\ref{eq:SelfImprovingOscillationPointwise}),
we arrive at
\begin{align*}
\left\Vert \osc{\delta\cdot T_{i}^{-T}\left[-1,1\right]^{\dimension}}f_{n}\right\Vert _{W_{T_{i}^{-T}\left[-1,1\right]^{\dimension}}\left(L_{v}^{p}\right)} & \leq2\sqrt{\dimension}\cdot\delta\cdot\left\Vert \left(\nabla g_{n}\right)\circ T_{i}^{T}\right\Vert _{W_{T_{i}^{-T}\left[-2,2\right]^{\dimension}}\left(L_{v}^{p}\right)}\\
 & \leq\delta\cdot2\dimension\cdot\Omega_{0}^{2K}\Omega_{1}^{2}\left(1+2\sqrt{\dimension}\right)^{N}4^{\dimension/p}C_{1}\left(10+4R\right)^{\dimension\left(\frac{1}{p}-1\right)}\cdot\left\Vert f\right\Vert _{L_{v}^{p}}.
\end{align*}
The remainder of the proof is now as for $p\in\left[1,\infty\right]$,
noting that
\begin{align*}
 & 2\dimension\cdot\Omega_{0}^{2K}\Omega_{1}^{2}\left(1+2\sqrt{\dimension}\right)^{N}4^{\dimension/p}\cdot C_{1}\left(10+4R\right)^{\dimension\left(\frac{1}{p}-1\right)}\\
\left({\scriptstyle \text{eq. }\eqref{eq:SelfImprovingOscillationC1Estimate}}\right) & \leq2\dimension\cdot\Omega_{0}^{2K}\Omega_{1}^{2}\left(1+2\sqrt{\dimension}\right)^{N}4^{\dimension/p}\cdot40\cdot10^{\dimension}\cdot\left(96\dimension\cdot N\right)^{N+1}\cdot10^{\dimension\left(\frac{1}{p}-1\right)}\left(1+R\right)^{1+\frac{\dimension}{p}}\\
\left({\scriptstyle \text{since }\dimension\leq2^{\dimension}\leq2^{\dimension/p}}\right) & \leq\Omega_{0}^{2K}\Omega_{1}^{2}\cdot\left(3\sqrt{\dimension}\right)^{N}\cdot\left(96\dimension\cdot N\right)^{N+1}\cdot80^{1+\frac{\dimension}{p}}\left(1+R\right)^{1+\frac{\dimension}{p}}\\
\left({\scriptstyle \text{since }N+1\geq\frac{\dimension+1}{p}+1\geq\frac{\dimension}{p}+1}\right) & \leq\Omega_{0}^{2K}\Omega_{1}^{2}\cdot\left(23040\cdot\dimension^{3/2}\cdot N\right)^{N+1}\cdot\left(1+R\right)^{1+\frac{\dimension}{p}}\\
\left({\scriptstyle \text{since }N\leq K+\frac{\dimension+1}{p}+1}\right) & \leq\Omega_{0}^{2K}\Omega_{1}^{2}\cdot\left(23040\cdot\dimension^{3/2}\cdot\left(K+1+\frac{\dimension+1}{p}\right)\right)^{K+2+\frac{\dimension+1}{p}}\cdot\left(1+R\right)^{1+\frac{\dimension}{p}}.\qedhere
\end{align*}
\end{proof}

\subsection{Convolution relation for Wiener amalgam spaces}

Finally, we come to the convolution relation for the Wiener Amalgam
spaces. The theorem stated here is a slight variation (and specialization)
of \cite[Theorem 2.3.24]{VoigtlaenderPhDThesis}, which originally
appeared in \cite{RauhutWienerAmalgam}.
\begin{thm}
\label{thm:WienerAmalgamConvolution}Let $Q_{1},Q_{2}\subset\R^{\dimension}$
be bounded, Borel measurable unit neighborhoods and assume that $Q_{1}-Q_{1}$
and $Q_{2}-Q_{1}$ are measurable. Let $p\in\left(0,\infty\right]$
and set $r:=\min\left\{ 1,p\right\} $.

Assume that there is a countable family $\left(x_{i}\right)_{i\in I}$
in $\R^{\dimension}$ satisfying $\R^{\dimension}=\bigcup_{i\in I}\left(x_{i}+Q_{1}\right)$
and
\[
N:=\sup_{x\in\R^{\dimension}}\left|\left\{ i\in I\with x\in x_{i}+Q_{1}\right\} \right|<\infty.
\]
Then we have for every $f\in W_{Q_{1}-Q_{1}}\left(L_{v_{0}}^{r}\right)$
and every $g\in W_{Q_{2}-Q_{1}}\left(L_{v}^{p}\right)$ that

\begin{itemize}
\item $f\in L_{v_{0}}^{1}\left(\R^{\dimension}\right)$ with $\left\Vert f\right\Vert _{L_{v_{0}}^{1}}\lesssim\left\Vert f\right\Vert _{W_{Q_{1}-Q_{1}}\left(L_{v_{0}}^{r}\right)}$,
where the implied constant only depends on $N,Q_{1},r,v_{0}$.
\item $g\in L_{v}^{\infty}\left(\R^{\dimension}\right)$ with $\left\Vert g\right\Vert _{L_{v}^{\infty}}\lesssim\left\Vert g\right\Vert _{W_{Q_{2}-Q_{1}}\left(L_{v}^{p}\right)}$,
where the implied constant only depends on $N,Q_{1},Q_{2},p,v_{0}$.
In particular, 
\[
W_{Q_{2}-Q_{1}}\left(L_{v}^{p}\right)\hookrightarrow L_{v}^{\infty}\left(\smash{\R^{\dimension}}\right)\hookrightarrow L_{\left(1+\left|\mybullet\right|\right)^{-K}}^{\infty}\left(\smash{\R^{\dimension}}\right)\hookrightarrow\Schwartz'\left(\smash{\R^{\dimension}}\right).
\]
\item The convolution $f\ast g:\R^{\dimension}\to\Compl$ is a well-defined
continuous function.
\item We have
\[
\left\Vert f\ast g\right\Vert _{W_{Q_{2}}\left(L_{v}^{p}\right)}\leq N^{\frac{1}{r}}\cdot\left[\sup_{x\in Q_{1}}v_{0}\left(x\right)\right]\cdot\left[\lambda_{\dimension}\left(Q_{1}\right)\right]^{1-\frac{1}{r}}\cdot\left\Vert f\right\Vert _{W_{Q_{1}-Q_{1}}\left(L_{v_{0}}^{r}\right)}\cdot\left\Vert g\right\Vert _{W_{Q_{2}-Q_{1}}\left(L_{v}^{p}\right)}.\qedhere
\]
\end{itemize}
\end{thm}
\begin{rem}
\label{rem:WeightedWienerAmalgamTemperedDistribution}Since one can
always choose a compact unit neighborhood $Q_{1}$ for which the assumptions
of the theorem are satisfied (choose e.g.\@ $Q_{1}=\left[-\frac{1}{2},\frac{1}{2}\right]^{\dimension}$
and $x_{i}:=i$ for $i\in I:=\Z^{\dimension}$), we see in view of
Lemma \ref{lem:WienerAmalgamNormEquivalence} that
\begin{equation}
W_{Q}\left(L_{v}^{p}\right)\hookrightarrow L_{v}^{\infty}\left(\smash{\R^{\dimension}}\right)\hookrightarrow L_{\left(1+\left|\mybullet\right|\right)^{-K}}^{\infty}\left(\smash{\R^{\dimension}}\right)\hookrightarrow\Schwartz'\left(\smash{\R^{\dimension}}\right)\label{eq:WeightedWienerAmalgamTemperedDistribution}
\end{equation}
holds for every bounded unit-neighborhood $Q\subset\R^{\dimension}$.
The same also holds with $v$ instead of $v_{0}$, since $v_{0}$
satisfies all properties that $v$ has.
\end{rem}
\begin{proof}
In the following, we will frequently use the discrete weights $v_{i}^{{\rm disc}}:=v\left(x_{i}\right)$
and $\left(v_{0}^{{\rm disc}}\right)_{i}:=v_{0}\left(x_{i}\right)$,
as well as the constant $C_{1}:=\sup_{x\in Q_{1}}v_{0}\left(x\right)\leq\Omega_{1}\cdot\sup_{x\in Q_{1}}\left(1+\left|x\right|\right)^{K}<\infty$,
which is finite since $Q_{1}$ is bounded. These quantities are important,
since we have for $y=x_{i}+q\in x_{i}+Q_{1}$ the estimates
\[
v\left(y\right)=v\left(x_{i}+q\right)\leq v\left(x_{i}\right)v_{0}\left(q\right)\leq C_{1}\cdot v_{i}^{{\rm disc}}\quad\text{ and similarly }\quad v_{0}\left(y\right)\leq C_{1}\cdot\left(v_{0}^{{\rm disc}}\right)_{i}.
\]
Likewise, by symmetry and submultiplicativity of $v_{0}$, we also
have
\[
\left(v_{0}^{{\rm disc}}\right)_{i}=v_{0}\left(x_{i}\right)=v_{0}\left(y-q\right)\leq v_{0}\left(y\right)\cdot v_{0}\left(-q\right)=v_{0}\left(y\right)\cdot v_{0}\left(q\right)\leq C_{1}\cdot v_{0}\left(y\right).
\]
Completely similar, we also get $v_{i}^{{\rm disc}}\leq C_{1}\cdot v\left(y\right)$
for all $y\in x_{i}+Q_{1}$.

\medskip{}

Now, we first show that we can write each $h\in W_{Q_{1}-Q_{1}}\left(L_{v_{0}}^{r}\right)$
as $h=\sum_{i\in I}h_{i}$, where\footnote{In this proof and the next, but not elsewhere in the paper, we write
$\supp f:=\left\{ x\in\R^{\dimension}\with f\left(x\right)\neq0\right\} $,
which is different from the usual meaning $\supp f:=\overline{\left\{ x\in\R^{\dimension}\with f\left(x\right)\neq0\right\} }$.} $\supp h_{i}\subset x_{i}+Q_{1}$ and where
\begin{equation}
\left\Vert \left(\left\Vert h_{i}\right\Vert _{L^{\infty}}\right)_{i\in I}\right\Vert _{\ell_{v_{0}^{{\rm disc}}}^{r}}\leq C_{1}\cdot\frac{N^{1/r}}{\left[\lambda_{\dimension}\left(Q_{1}\right)\right]^{1/r}}\cdot\left\Vert h\right\Vert _{W_{Q_{1}-Q_{1}}\left(L_{v_{0}}^{r}\right)}.\label{eq:WienerAmalgamDecomposition}
\end{equation}
Indeed, since $I$ is countable (and necessarily infinite, since $\R^{\dimension}=\bigcup_{i\in I}\left(x_{i}+Q_{1}\right)$,
with $Q_{1}$ bounded), we can assume $I=\N$. Then, define $h_{i}:=h\cdot\Indicator_{\left(x_{i}+Q_{1}\right)\setminus\bigcup_{j=1}^{i-1}\left(x_{j}+Q_{1}\right)}$.
Because of $\R^{\dimension}=\bigcup_{i\in I}\left(x_{i}+Q_{1}\right)$,
this easily yields $h=\sum_{i\in I}h_{i}$ and $\supp h_{i}\subset x_{i}+Q_{1}$
is trivial, so that we only need to verify estimate (\ref{eq:WienerAmalgamDecomposition}).

To this end, first note for $x\in x_{i}+Q_{1}$ that $x_{i}\in x-Q_{1}$
and hence $x_{i}+Q_{1}\subset x+Q_{1}-Q_{1}$, which yields
\[
\left\Vert h_{i}\right\Vert _{L^{\infty}}\leq\left\Vert h\cdot\Indicator_{x_{i}+Q_{1}}\right\Vert _{L^{\infty}}\leq\left\Vert h\cdot\Indicator_{x+Q_{1}-Q_{1}}\right\Vert _{L^{\infty}}=\left(M_{Q_{1}-Q_{1}}h\right)\left(x\right).
\]
Now, take the $r$-th power of this estimate, multiply both sides
with $v_{0}^{r}\left(x\right)\cdot\Indicator_{x_{i}+Q_{1}}\left(x\right)$
and sum over $i\in I$ to arrive at
\[
\sum_{i\in I}\left[\left\Vert h_{i}\right\Vert _{L^{\infty}}^{r}\cdot v_{0}^{r}\left(x\right)\cdot\Indicator_{x_{i}+Q_{1}}\left(x\right)\right]\leq\left[v_{0}\left(x\right)\cdot\left(M_{Q_{1}-Q_{1}}h\right)\left(x\right)\right]^{r}\cdot\sum_{i\in I}\Indicator_{x_{i}+Q_{1}}\left(x\right)\leq N\cdot\left[v_{0}\left(x\right)\cdot\left(M_{Q_{1}-Q_{1}}h\right)\left(x\right)\right]^{r}.
\]
As observed at the beginning of the proof, we have $\Indicator_{x_{i}+Q_{1}}\cdot\left(v_{0}^{{\rm disc}}\right)_{i}\leq C_{1}\cdot v_{0}\cdot\Indicator_{x_{i}+Q_{1}}$.
By combining this with the preceding estimate and integrating, we
get
\begin{align*}
\lambda_{\dimension}\left(Q_{1}\right)\cdot\left\Vert \left(\left\Vert h_{i}\right\Vert _{L^{\infty}}\right)_{i\in I}\right\Vert _{\ell_{v_{0}^{{\rm disc}}}^{r}}^{r} & =\sum_{i\in I}\left\Vert h_{i}\right\Vert _{L^{\infty}}^{r}\cdot\left(v_{0}^{{\rm disc}}\right)_{i}^{r}\cdot\lambda_{\dimension}\left(x_{i}+Q_{1}\right)\\
 & =\int_{\R^{\dimension}}\sum_{i\in I}\left(\left\Vert h_{i}\right\Vert _{L^{\infty}}\cdot\left(v_{0}^{{\rm disc}}\right)_{i}\cdot\Indicator_{x_{i}+Q_{1}}\left(x\right)\right)^{r}\d x\\
 & \leq C_{1}^{r}N\cdot\int_{\R^{\dimension}}\left[v_{0}\left(x\right)\cdot\left(M_{Q_{1}-Q_{1}}h\right)\left(x\right)\right]^{r}\d x\\
 & =C_{1}^{r}N\cdot\left\Vert h\right\Vert _{W_{Q_{1}-Q_{1}}\left(L_{v_{0}}^{r}\right)}^{r}<\infty.
\end{align*}
Rearranging shows that equation (\ref{eq:WienerAmalgamDecomposition})
is indeed satisfied.

\medskip{}

Now, since $\ell^{r}\left(I\right)\hookrightarrow\ell^{1}\left(I\right)$
and because of $\supp h_{i}\subset x_{i}+Q_{1}$, so that 
\[
v_{0}\cdot\left|h_{i}\right|\leq C_{1}\cdot\left(v_{0}^{{\rm disc}}\right)_{i}\cdot\left|h_{i}\right|\leq C_{1}\cdot\left(v_{0}^{{\rm disc}}\right)_{i}\cdot\left\Vert h_{i}\right\Vert _{L^{\infty}}\cdot\Indicator_{x_{i}+Q_{1}}\quad\text{ almost everywhere},
\]
we get
\begin{align*}
\left\Vert h\right\Vert _{L_{v_{0}}^{1}}\leq\sum_{i\in I}\left\Vert h_{i}\right\Vert _{L_{v_{0}}^{1}} & \leq\left[\sum_{i\in I}\left\Vert h_{i}\right\Vert _{L_{v_{0}}^{1}}^{r}\right]^{1/r}\\
 & \leq C_{1}\cdot\left[\sum_{i\in I}\left(v_{0}^{{\rm disc}}\right)_{i}^{r}\cdot\left\Vert h_{i}\right\Vert _{L^{\infty}}^{r}\cdot\left[\lambda\left(x_{i}+Q_{1}\right)\right]^{r}\right]^{1/r}\\
 & \leq C_{1}\cdot\lambda_{\dimension}\left(Q_{1}\right)\cdot\left\Vert \left(\left\Vert h_{i}\right\Vert _{L^{\infty}}\right)_{i\in I}\right\Vert _{\ell_{v_{0}^{{\rm disc}}}^{r}}\\
\left({\scriptstyle \text{eq. }\eqref{eq:WienerAmalgamDecomposition}}\right) & \leq C_{1}^{2}\cdot\left[\lambda_{\dimension}\left(Q_{1}\right)\right]^{1-\frac{1}{r}}\cdot N^{1/r}\cdot\left\Vert h\right\Vert _{W_{Q_{1}-Q_{1}}\left(L_{v_{0}}^{r}\right)}<\infty,
\end{align*}
which proves the first part of the theorem.

\medskip{}

Now, we want to prove the second part of the theorem. For $p=\infty$,
we have $\left\Vert g\right\Vert _{L_{v}^{\infty}}=\left\Vert g\right\Vert _{L_{v}^{p}}\leq\left\Vert g\right\Vert _{W_{Q_{2}-Q_{1}}\left(L_{v}^{p}\right)}$
by Lemma \ref{lem:MaximalFunctionDominatesF}, so that we can assume
$p\in\left(0,\infty\right)$.

Next, we define $g_{i}:=g\cdot\Indicator_{x_{i}+Q_{1}}$ for $i\in I$
and note for $x\in x_{i}+Q_{1}$ as above that $x_{i}+Q_{1}\subset x+Q_{1}-Q_{1}$,
so that
\[
\left\Vert g_{i}\right\Vert _{L^{\infty}}\leq\left\Vert g\cdot\Indicator_{x+Q_{1}-Q_{1}}\right\Vert _{L^{\infty}}=\left(M_{Q_{1}-Q_{1}}g\right)\left(x\right).
\]
Hence,
\begin{align*}
\frac{1}{C_{1}}\cdot v_{i}^{{\rm disc}}\cdot\left\Vert g_{i}\right\Vert _{L^{\infty}}\cdot\left[\lambda_{\dimension}\left(Q_{1}\right)\right]^{1/p} & =\frac{1}{C_{1}}\cdot\left[\int_{x_{i}+Q_{1}}\left(v_{i}^{{\rm disc}}\cdot\left\Vert g_{i}\right\Vert _{L^{\infty}}\right)^{p}\d x\right]^{1/p}\\
 & \leq\left[\int_{\R^{\dimension}}\left(v\left(x\right)\cdot\left\Vert g_{i}\right\Vert _{L^{\infty}}\cdot\Indicator_{x_{i}+Q_{1}}\left(x\right)\right)^{p}\d x\right]^{1/p}\\
 & \leq\left(\int_{\R^{\dimension}}\left[\left(v\cdot M_{Q_{1}-Q_{1}}g\right)\left(x\right)\right]^{p}\d x\right)^{1/p}\\
 & =\left\Vert g\right\Vert _{W_{Q_{1}-Q_{1}}\left(L_{v}^{p}\right)}\\
\left({\scriptstyle \text{Lemma }\ref{lem:WienerAmalgamNormEquivalence}}\right) & \leq C_{2}\cdot\left\Vert g\right\Vert _{W_{Q_{2}-Q_{1}}\left(L_{v}^{p}\right)}.
\end{align*}
Here, the last step used that $Q_{1}-Q_{1}$ and $Q_{2}-Q_{1}$ are
both measurable, bounded unit-neighborhoods, so that Lemma \ref{lem:WienerAmalgamNormEquivalence}
yields a constant $C_{2}=C_{2}\left(Q_{1},Q_{2},v_{0},p\right)>0$
satisfying $\left\Vert g\right\Vert _{W_{Q_{1}-Q_{1}}\left(L_{v}^{p}\right)}\leq C_{2}\cdot\left\Vert g\right\Vert _{W_{Q_{2}-Q_{1}}\left(L_{v}^{p}\right)}$.

But there is a null-set $N_{i}\subset x_{i}+Q_{1}$ satisfying $\left|g\left(x\right)\right|=\left|g_{i}\left(x\right)\right|\leq\left\Vert g_{i}\right\Vert _{L^{\infty}}$
for all $x\in\left(x_{i}+Q_{1}\right)\setminus N_{i}$. Hence,
\[
v\left(x\right)\cdot\left|g\left(x\right)\right|\leq C_{1}\cdot v_{i}^{{\rm disc}}\cdot\left\Vert g_{i}\right\Vert _{L^{\infty}}\leq\frac{C_{1}^{2}C_{2}}{\left[\lambda_{\dimension}\left(Q_{1}\right)\right]^{1/p}}\cdot\left\Vert g\right\Vert _{W_{Q_{2}-Q_{1}}\left(L_{v}^{p}\right)}
\]
for all $x\in\left(x_{i}+Q_{1}\right)\setminus N_{i}$. But since
$N:=\bigcup_{i\in I}N_{i}\subset\R^{\dimension}$ is a null-set and
since $\R^{\dimension}=\bigcup_{i\in I}\left(x_{i}+Q_{1}\right)$,
we get $\left\Vert g\right\Vert _{L_{v}^{\infty}}\leq\frac{C_{1}^{2}C_{2}}{\left[\lambda_{\dimension}\left(Q_{1}\right)\right]^{1/p}}\cdot\left\Vert g\right\Vert _{W_{Q_{2}-Q_{1}}\left(L_{v}^{p}\right)}$,
which proves the main part of the second part of the theorem for $p\in\left(0,\infty\right)$.

To establish the embedding $W_{Q_{2}-Q_{1}}\left(L_{v}^{p}\right)\hookrightarrow L_{v}^{\infty}\left(\R^{\dimension}\right)\hookrightarrow L_{\left(1+\left|\mybullet\right|\right)^{-K}}^{\infty}\left(\R^{\dimension}\right)\hookrightarrow\Schwartz'\left(\R^{\dimension}\right)$,
we first observe that $L_{\left(1+\left|\mybullet\right|\right)^{-K}}^{\infty}\left(\R^{\dimension}\right)\hookrightarrow\Schwartz'\left(\R^{\dimension}\right)$
is trivial. Furthermore,
\begin{equation}
v\left(0\right)=v\left(x+\left(-x\right)\right)\leq v\left(x\right)\cdot v_{0}\left(-x\right)\leq\Omega_{1}\left(1+\left|-x\right|\right)^{K}\cdot v\left(x\right)\qquad\forall x\in\R^{\dimension},\label{eq:WeightBoundedBelow}
\end{equation}
so that $v\left(x\right)\geq\frac{v\left(0\right)}{\Omega_{1}}\cdot\left(1+\left|x\right|\right)^{-K}$
and hence $W_{Q_{2}-Q_{1}}\left(L_{v}^{p}\right)\hookrightarrow L_{v}^{\infty}\left(\R^{\dimension}\right)\hookrightarrow L_{\left(1+\left|\mybullet\right|\right)^{-K}}^{\infty}\left(\R^{\dimension}\right)$,
as desired.

\medskip{}

Now, note for $f\in L_{v_{0}}^{1}\left(\R^{\dimension}\right)$ and
$g\in L_{v}^{\infty}\left(\R^{\dimension}\right)$ because of $v\left(x\right)=v\left(y+\left(x-y\right)\right)\leq v\left(y\right)\cdot v_{0}\left(x-y\right)$
that
\begin{align}
v\left(x\right)\cdot\int_{\R^{\dimension}}\left|f\left(x-y\right)\right|\cdot\left|g\left(y\right)\right|\d y & \leq\int_{\R^{\dimension}}\left|\left(v_{0}\cdot f\right)\left(x-y\right)\right|\cdot\left|\left(v\cdot g\right)\left(y\right)\right|\d y\nonumber \\
 & \leq\left\Vert g\right\Vert _{L_{v}^{\infty}}\cdot\left\Vert f\right\Vert _{L_{v_{0}}^{1}}<\infty\qquad\forall x\in\R^{\dimension}.\label{eq:WeightedLInftyConvolution}
\end{align}
Hence, $\left(f\ast g\right)\left(x\right)$ is well-defined for all
$x\in\R^{\dimension}$ and $\left\Vert f\ast g\right\Vert _{L_{\left(1+\left|\mybullet\right|\right)^{-K}}^{\infty}}\lesssim\left\Vert f\ast g\right\Vert _{L_{v}^{\infty}}\leq\left\Vert f\right\Vert _{L_{v_{0}}^{1}}\cdot\left\Vert g\right\Vert _{L_{v}^{\infty}}$.
Now, note that the subspace $C\left(\R^{\dimension}\right)\cap L_{\left(1+\left|\mybullet\right|\right)^{-K}}^{\infty}\left(\R^{\dimension}\right)$
of continuous functions in $L_{\left(1+\left|\mybullet\right|\right)^{-K}}^{\infty}\left(\R^{\dimension}\right)$
is a closed subspace of $L_{\left(1+\left|\mybullet\right|\right)^{-K}}^{\infty}\left(\R^{\dimension}\right)$.
Furthermore, $C_{c}\left(\R^{\dimension}\right)\subset L_{v_{0}}^{1}\left(\R^{\dimension}\right)$
is dense and $L_{v}^{\infty}\left(\R^{\dimension}\right)\hookrightarrow L_{\left(1+\left|\mybullet\right|\right)^{-K}}^{\infty}\left(\R^{\dimension}\right)\hookrightarrow L_{{\rm loc}}^{\infty}\left(\R^{\dimension}\right)$.
But for $f\in C_{c}\left(\R^{\dimension}\right)$ and $g\in L_{{\rm loc}}^{\infty}\left(\R^{\dimension}\right)$,
it is not hard to see that $f\ast g$ is continuous.

Altogether, the preceding properties show that $f\ast g\in C\left(\R^{\dimension}\right)\cap L_{\left(1+\left|\mybullet\right|\right)^{-K}}^{\infty}\left(\R^{\dimension}\right)$
is well-defined and continuous for all $f\in L_{v_{0}}^{1}\left(\R^{\dimension}\right)$
and $g\in L_{v}^{\infty}\left(\R^{\dimension}\right)$. But in the
setting of the theorem, we have $f\in W_{Q_{1}-Q_{1}}\left(L_{v_{0}}^{r}\right)\hookrightarrow L_{v_{0}}^{1}\left(\R^{\dimension}\right)$
and $g\in W_{Q_{2}-Q_{1}}\left(L_{v}^{p}\right)\hookrightarrow L_{v}^{\infty}\left(\R^{\dimension}\right)$,
so that the third part of the theorem is established.

\medskip{}

It remains to prove the last part of the theorem. To this end, recall
from equation (\ref{eq:WienerAmalgamDecomposition}) that we can write
$f=\sum_{i\in I}f_{i}$, where $\supp f_{i}\subset x_{i}+Q_{1}$ and
such that equation (\ref{eq:WienerAmalgamDecomposition}) is fulfilled,
with $f_{i}$ instead of $h_{i}$ and $f$ instead of $h$.

Next, we estimate $M_{Q_{2}}\left(f_{i}\ast g\right)$ for each $i\in I$
as follows: For $x\in\R^{\dimension}$ and $q\in Q_{2}$, we have
\begin{align*}
\left|\left(f_{i}\ast g\right)\left(x+q\right)\right|\leq\left(\left|f_{i}\right|\ast\left|g\right|\right)\left(x+q\right) & =\int_{\R^{\dimension}}\left|f_{i}\left(y\right)\right|\cdot\left|g\left(x+q-y\right)\right|\d y\\
 & \leq\left\Vert f_{i}\right\Vert _{L^{\infty}}\cdot\int_{\R^{\dimension}}\Indicator_{x_{i}+Q_{1}}\left(y\right)\cdot\left|g\left(x+q-y\right)\right|\d y\\
\left({\scriptstyle z=x+q-y}\right) & =\left\Vert f_{i}\right\Vert _{L^{\infty}}\cdot\int_{\R^{\dimension}}\Indicator_{x_{i}+Q_{1}}\left(x+q-z\right)\cdot\left|g\left(z\right)\right|\d z\\
\left({\scriptstyle x+q-z\in x_{i}+Q_{1}\text{ implies }z\in x-x_{i}+q-Q_{1}\subset x-x_{i}+Q_{2}-Q_{1}}\right) & \leq\left\Vert f_{i}\right\Vert _{L^{\infty}}\cdot\int_{\R^{\dimension}}\Indicator_{x_{i}+Q_{1}}\left(x+q-z\right)\d z\cdot\left\Vert g\cdot\Indicator_{x-x_{i}+Q_{2}-Q_{1}}\right\Vert _{L^{\infty}}\\
 & =\left\Vert f_{i}\right\Vert _{L^{\infty}}\cdot\lambda_{\dimension}\left(x+q-x_{i}-Q_{1}\right)\cdot\left(M_{Q_{2}-Q_{1}}g\right)\left(x-x_{i}\right)\\
 & =\lambda_{\dimension}\left(Q_{1}\right)\cdot\left\Vert f_{i}\right\Vert _{L^{\infty}}\cdot\left(L_{x_{i}}\left[M_{Q_{2}-Q_{1}}g\right]\right)\left(x\right).
\end{align*}
Since this holds for all $q\in Q_{2}$, we get
\[
\left[M_{Q_{2}}\left(\left|f_{i}\right|\ast\left|g\right|\right)\right]\left(x\right)\leq\lambda_{\dimension}\left(Q_{1}\right)\cdot\left\Vert f_{i}\right\Vert _{L^{\infty}}\cdot\left(L_{x_{i}}\left[M_{Q_{2}-Q_{1}}g\right]\right)\left(x\right)\qquad\forall x\in\R^{\dimension}.
\]
In view of Lemma \ref{lem:WeightedLpTranslationNorm} and by solidity
of $L_{v}^{p}\left(\R^{\dimension}\right)$, this implies
\begin{align*}
\left\Vert M_{Q_{2}}\left[\left|f_{i}\right|\ast\left|g\right|\right]\right\Vert _{L_{v}^{p}} & \leq\left\Vert f_{i}\right\Vert _{L^{\infty}}\cdot\lambda_{\dimension}\left(Q_{1}\right)\cdot\left\Vert L_{x_{i}}\left[M_{Q_{2}-Q_{1}}g\right]\right\Vert _{L_{v}^{p}}\\
 & \leq\left(v_{0}^{{\rm disc}}\right)_{i}\cdot\left\Vert f_{i}\right\Vert _{L^{\infty}}\cdot\lambda_{\dimension}\left(Q_{1}\right)\cdot\left\Vert g\right\Vert _{W_{Q_{2}-Q_{1}}\left(L_{v}^{p}\right)}.
\end{align*}

Next, it is not hard to see $M_{Q_{2}}\left(\sum_{i\in I}h_{i}\right)\leq\sum_{i\in I}M_{Q_{2}}h_{i}$,
so that we get because of
\[
\left|\left(f\ast g\right)\left(x\right)\right|\leq\left(\left|f\right|\ast\left|g\right|\right)\left(x\right)\leq\sum_{i\in I}\left(\left|f_{i}\right|\ast\left|g\right|\right)\left(x\right)
\]
that
\begin{align*}
\left\Vert M_{Q_{2}}\left(f\ast g\right)\right\Vert _{L_{v}^{p}}^{r}\leq\left\Vert M_{Q_{2}}\left[\left|f\right|\ast\left|g\right|\right]\right\Vert _{L_{v}^{p}}^{r} & \leq\left\Vert \sum_{i\in I}M_{Q_{2}}\left[\left|f_{i}\right|\ast\left|g\right|\right]\right\Vert _{L_{v}^{p}}^{r}\\
\left({\scriptstyle L_{v}^{p}\text{ satisfies the }r-\text{triangle inequality}}\right) & \leq\sum_{i\in I}\left\Vert M_{Q_{2}}\left[\left|f_{i}\right|\ast\left|g\right|\right]\right\Vert _{L_{v}^{p}}^{r}\\
 & \leq\left[\lambda_{\dimension}\left(Q_{1}\right)\cdot\left\Vert g\right\Vert _{W_{Q_{2}-Q_{1}}\left(L_{v}^{p}\right)}\right]^{r}\cdot\sum_{i\in I}\left(v_{0}^{{\rm disc}}\right)_{i}^{r}\cdot\left\Vert f_{i}\right\Vert _{L^{\infty}}^{r}\\
\left({\scriptstyle \text{eq. }\eqref{eq:WienerAmalgamDecomposition}}\right) & \leq\left[\lambda_{\dimension}\left(Q_{1}\right)\cdot\left\Vert g\right\Vert _{W_{Q_{2}-Q_{1}}\left(L_{v}^{p}\right)}\right]^{r}\cdot\left(C_{1}\cdot\frac{N^{1/r}}{\left[\lambda_{\dimension}\left(Q_{1}\right)\right]^{1/r}}\cdot\left\Vert f\right\Vert _{W_{Q_{1}-Q_{1}}\left(L_{v_{0}}^{r}\right)}\right)^{r},
\end{align*}
which finally yields
\[
\left\Vert f\ast g\right\Vert _{W_{Q_{2}}\left(L_{v}^{p}\right)}\leq N^{\frac{1}{r}}C_{1}\cdot\left[\lambda_{\dimension}\left(Q_{1}\right)\right]^{1-\frac{1}{r}}\cdot\left\Vert f\right\Vert _{W_{Q_{1}-Q_{1}}\left(L_{v_{0}}^{r}\right)}\cdot\left\Vert g\right\Vert _{W_{Q_{2}-Q_{1}}\left(L_{v}^{p}\right)}<\infty,
\]
as desired.
\end{proof}
With a very slight variant of the above proof, one can also show the
following modification of the theorem. For completeness, we provide
the proof, but with slightly less details than above.
\begin{prop}
\label{prop:AlternativeWienerAmalgamConvolution}Under the assumptions
of Theorem \ref{thm:WienerAmalgamConvolution}, if $p\in\left(0,1\right]$,
then
\[
\left\Vert f\ast g\right\Vert _{W_{Q_{2}}\left(L_{v}^{p}\right)}\leq N^{\frac{1}{p}}\cdot\left[\sup_{x\in Q_{1}}v_{0}\left(x\right)\right]\cdot\left[\lambda_{\dimension}\left(Q_{1}\right)\right]^{1-\frac{1}{p}}\cdot\left\Vert f\right\Vert _{W_{Q_{2}-Q_{1}}\left(L_{v_{0}}^{p}\right)}\cdot\left\Vert g\right\Vert _{W_{Q_{1}-Q_{1}}\left(L_{v}^{p}\right)}.\qedhere
\]
\end{prop}
\begin{proof}
As in the proof of Theorem \ref{thm:WienerAmalgamConvolution}, let
$C_{1}:=\sup_{x\in Q_{1}}v_{0}\left(x\right)$. Also as in that proof,
we can assume $I=\N$, so that we have $g=\sum_{i\in I}g_{i}$ with
$\supp g_{i}\subset x_{i}+Q_{1}$ for $g_{i}:=g\cdot\Indicator_{\left(x_{i}+Q_{1}\right)\setminus\bigcup_{j=1}^{i-1}\left(x_{j}+Q_{1}\right)}$.
Furthermore, for arbitrary $x\in x_{i}+Q_{1}$, we have $x_{i}+Q_{1}\subset x+Q_{1}-Q_{1}$
and thus
\[
\left\Vert g_{i}\right\Vert _{L^{\infty}}\leq\left\Vert g\cdot\Indicator_{x_{i}+Q_{1}}\right\Vert _{L^{\infty}}\leq\left\Vert g\cdot\Indicator_{x+Q_{1}-Q_{1}}\right\Vert _{L^{\infty}}=\left(M_{Q_{1}-Q_{1}}g\right)\left(x\right).
\]
Now, multiply both sides with $v\left(x\right)$, take the $p$th
power, multiply with $\Indicator_{x_{i}+Q_{1}}\left(x\right)$ and
sum over $i\in I$ to obtain
\begin{align*}
\sum_{i\in I}\left[v\left(x\right)\cdot\left\Vert g_{i}\right\Vert _{L^{\infty}}\right]^{p}\Indicator_{x_{i}+Q_{1}}\left(x\right) & \leq\sum_{i\in I}\left[v\left(x\right)\cdot\left(M_{Q_{1}-Q_{1}}g\right)\left(x\right)\right]^{p}\Indicator_{x_{i}+Q_{1}}\left(x\right)\\
 & \leq N\cdot\left[v\left(x\right)\cdot\left(M_{Q_{1}-Q_{1}}g\right)\left(x\right)\right]^{p}.
\end{align*}
But for $x\in x_{i}+Q_{1}$, i.e., $x=x_{i}+q$ with $q\in Q_{1}$,
we have
\[
v\left(x_{i}\right)=v\left(x-q\right)\leq v\left(x\right)\cdot v_{0}\left(-q\right)=v\left(x\right)\cdot v_{0}\left(q\right)\leq C_{1}\cdot v\left(x\right),
\]
so that we arrive at
\[
\sum_{i\in I}\left[v\left(x_{i}\right)\cdot\left\Vert g_{i}\right\Vert _{L^{\infty}}\right]^{p}\Indicator_{x_{i}+Q_{1}}\left(x\right)\leq C_{1}^{p}\cdot N\cdot\left[v\left(x\right)\cdot\left(M_{Q_{1}-Q_{1}}g\right)\left(x\right)\right]^{p}.
\]
Integrating this estimate over $x\in\R^{\dimension}$ finally yields
\begin{equation}
\lambda_{\dimension}\left(Q_{1}\right)\cdot\sum_{i\in I}\left[v\left(x_{i}\right)\cdot\left\Vert g_{i}\right\Vert _{L^{\infty}}\right]^{p}\leq C_{1}^{p}\cdot N\cdot\left\Vert g\right\Vert _{W_{Q_{1}-Q_{1}}\left(L_{v}^{p}\right)}^{p}<\infty.\label{eq:AlternativeWienerAmalgamConvolutionDecompositionStep}
\end{equation}

Now, let $x\in\R^{\dimension}$ and $q\in Q_{2}$ be arbitrary. Since
$\supp g_{i}\subset x_{i}+Q_{1}$, we have
\begin{align*}
\left(\left|f\right|\ast\left|g_{i}\right|\right)\left(x+q\right) & \leq\left\Vert g_{i}\right\Vert _{L^{\infty}}\cdot\int_{\R^{\dimension}}\Indicator_{x_{i}+Q_{1}}\left(y\right)\cdot\left|f\left(x+q-y\right)\right|\d y\\
\left({\scriptstyle z=x+q-y}\right) & =\left\Vert g_{i}\right\Vert _{L^{\infty}}\cdot\int_{\R^{\dimension}}\Indicator_{x_{i}+Q_{1}}\left(x+q-z\right)\cdot\left|f\left(z\right)\right|\d z\\
\left({\scriptstyle x+q-z\in x_{i}+Q_{1}\text{ implies }z\in x+q-x_{i}-Q_{1}\subset x-x_{i}+Q_{2}-Q_{1}}\right) & \leq\left\Vert g_{i}\right\Vert _{L^{\infty}}\cdot\int_{\R^{\dimension}}\Indicator_{x_{i}+Q_{1}}\left(x+q-z\right)\d z\cdot\left\Vert f\cdot\Indicator_{x-x_{i}+Q_{2}-Q_{1}}\right\Vert _{L^{\infty}}\\
 & \leq\left\Vert g_{i}\right\Vert _{L^{\infty}}\cdot\lambda_{\dimension}\left(x+q-x_{i}-Q_{1}\right)\cdot\left\Vert f\cdot\Indicator_{x-x_{i}+Q_{2}-Q_{1}}\right\Vert _{L^{\infty}}\\
 & =\lambda_{\dimension}\left(Q_{1}\right)\cdot\left\Vert g_{i}\right\Vert _{L^{\infty}}\cdot\left(M_{Q_{2}-Q_{1}}f\right)\left(x-x_{i}\right).
\end{align*}
Since this holds for arbitrary $q\in Q_{2}$, we have shown
\[
\left[M_{Q_{2}}\left(\left|f\right|\ast\left|g_{i}\right|\right)\right]\left(x\right)\leq\lambda_{\dimension}\left(Q_{1}\right)\cdot\left\Vert g_{i}\right\Vert _{L^{\infty}}\cdot\left(M_{Q_{2}-Q_{1}}f\right)\left(x-x_{i}\right).
\]
Hence,
\begin{align*}
v\left(x\right)\cdot\left[M_{Q_{2}}\left(\left|f\right|\ast\left|g_{i}\right|\right)\right]\left(x\right) & \leq\lambda_{\dimension}\left(Q_{1}\right)\cdot\left\Vert g_{i}\right\Vert _{L^{\infty}}\cdot v\left(x\right)\cdot\left(M_{Q_{2}-Q_{1}}f\right)\left(x-x_{i}\right)\\
\left({\scriptstyle \text{since }v\left(x\right)=v\left(x-x_{i}+x_{i}\right)\leq v_{0}\left(x-x_{i}\right)v\left(x_{i}\right)}\right) & \leq\lambda_{\dimension}\left(Q_{1}\right)\cdot v\left(x_{i}\right)\left\Vert g_{i}\right\Vert _{L^{\infty}}\cdot\left[v_{0}\cdot M_{Q_{2}-Q_{1}}f\right]\left(x-x_{i}\right).
\end{align*}
Taking the $L^{p}$ norm on both sides, and using the isometric translation
invariance of $L^{p}$, we conclude
\[
\left\Vert \left|f\right|\ast\left|g_{i}\right|\right\Vert _{W_{Q_{2}}\left(L_{v}^{p}\right)}\leq\lambda_{\dimension}\left(Q_{1}\right)\cdot v\left(x_{i}\right)\left\Vert g_{i}\right\Vert _{L^{\infty}}\cdot\left\Vert f\right\Vert _{W_{Q_{2}-Q_{1}}\left(L_{v_{0}}^{p}\right)}.
\]
Now, we finally combine the estimate $\left|\left(f\ast g\right)\left(x\right)\right|\leq\left(\left|f\right|\ast\left|g\right|\right)\left(x\right)\leq\sum_{i\in I}\left(\left|f\right|\ast\left|g_{i}\right|\right)\left(x\right)$
with solidity of $W_{Q_{2}}\left(L_{v}^{p}\right)$ and with the $p$-triangle
inequality for $W_{Q_{2}}\left(L_{v}^{p}\right)$ (which holds sine
$p\in\left(0,1\right]$) to deduce
\begin{align*}
\left\Vert f\ast g\right\Vert _{W_{Q_{2}}\left(L_{v}^{p}\right)}^{p} & \leq\left[\lambda_{\dimension}\left(Q_{1}\right)\right]^{p}\cdot\left\Vert f\right\Vert _{W_{Q_{2}-Q_{1}}\left(L_{v_{0}}^{p}\right)}^{p}\cdot\sum_{i\in I}\left[v\left(x_{i}\right)\left\Vert g_{i}\right\Vert _{L^{\infty}}\right]^{p}\\
\left({\scriptstyle \text{eq. }\eqref{eq:AlternativeWienerAmalgamConvolutionDecompositionStep}}\right) & \leq\left[\lambda_{\dimension}\left(Q_{1}\right)\right]^{p-1}\cdot C_{1}^{p}\cdot N\cdot\left\Vert f\right\Vert _{W_{Q_{2}-Q_{1}}\left(L_{v_{0}}^{p}\right)}^{p}\cdot\left\Vert g\right\Vert _{W_{Q_{1}-Q_{1}}\left(L_{v}^{p}\right)}^{p},
\end{align*}
which easily yields the claim.
\end{proof}
We now formulate an important special case of Theorem \ref{thm:WienerAmalgamConvolution}
as a corollary:
\begin{cor}
\label{cor:WienerAmalgamConvolutionSimplified}Let $i,j\in I$, $p\in\left(0,\infty\right]$,
 $f\in W_{T_{j}^{-T}\left[-1,1\right]^{\dimension}}\left(L_{v_{0}}^{r}\right)$
for $r:=\min\left\{ 1,p\right\} $ and $g\in W_{T_{j}^{-T}\left[-1,1\right]^{\dimension}}\left(L_{v}^{p}\right)$.
Then the convolution $f\ast g$ is pointwise defined and continuous
and we have
\[
\left\Vert f\ast g\right\Vert _{W_{T_{j}^{-T}\left[-1,1\right]^{\dimension}}\left(L_{v}^{p}\right)}\leq\Omega_{0}^{3K}\Omega_{1}^{3}C\cdot\left|\det T_{j}\right|^{\frac{1}{r}-1}\cdot\left\Vert f\right\Vert _{W_{T_{j}^{-T}\left[-1,1\right]^{\dimension}}\left(L_{v_{0}}^{r}\right)}\cdot\left\Vert g\right\Vert _{W_{T_{j}^{-T}\left[-1,1\right]^{\dimension}}\left(L_{v}^{p}\right)}
\]
for $C:=\dimension^{-\frac{\dimension}{2r}}\cdot\left(972\cdot\dimension^{5/2}\right)^{K+\frac{\dimension}{r}}$.
\end{cor}
\begin{proof}
We apply Theorem \ref{thm:WienerAmalgamConvolution} with $Q_{1}=Q_{2}=T_{j}^{-T}\left[-1,1\right]^{\dimension}$.
Note that we have
\[
\R=\bigcup_{k\in\Z}\left(2k+\left[-1,1\right]\right)\qquad\text{ and hence }\qquad\R^{\dimension}=\bigcup_{k\in\Z^{\dimension}}\left(2k+\left[-1,1\right]^{\dimension}\right).
\]
Furthermore, if $x\in\left(2k+\left[-1,1\right]^{\dimension}\right)\cap\left(2\ell+\left[-1,1\right]^{\dimension}\right)$,
we get $2k+\mu=x=2\ell+\nu$ for certain $\mu,\nu\in\left[-1,1\right]^{\dimension}$
and thus $\left\Vert k-\ell\right\Vert _{\infty}=\left\Vert \frac{\nu-\mu}{2}\right\Vert _{\infty}\leq1$.
Thus, we see (by fixing $k\in\Z^{\dimension}$ with $x\in2k+\left[-1,1\right]^{\dimension}$)
that $x\in2\ell+\left[-1,1\right]^{\dimension}$ can hold for at most
$3^{\dimension}$ values of $\ell$, namely for $\ell\in\prod_{j=1}^{\dimension}\left\{ k_{j}-1,k_{j},k_{j}+1\right\} $.
Since $T_{j}^{-T}:\R^{\dimension}\to\R^{\dimension}$ is bijective,
we see 
\begin{equation}
\R^{\dimension}=\bigcup_{k\in\Z^{\dimension}}\left(2T_{j}^{-T}k+T_{j}^{-T}\left[-1,1\right]^{\dimension}\right)\quad\text{ and }\quad N:=\sup_{x\in\R^{\dimension}}\left|\left\{ \ell\in\Z^{\dimension}\with x\in2T_{j}^{-T}\ell+T_{j}^{-T}\left[-1,1\right]^{\dimension}\right\} \right|\leq3^{\dimension}.\label{eq:LinearImageOfCubePartition}
\end{equation}

Furthermore, equation (\ref{eq:WeightLinearTransformationsConnection})
yields
\begin{align*}
\sup_{x\in Q_{1}}v_{0}\left(x\right)=\sup_{y\in\left[-1,1\right]^{\dimension}}v_{0}\left(T_{j}^{-T}y\right) & \leq\Omega_{1}\cdot\sup_{y\in\left[-1,1\right]^{\dimension}}\left(1+\left|T_{j}^{-T}y\right|\right)^{K}\\
\left({\scriptstyle \text{eq. }\eqref{eq:WeightLinearTransformationsConnection}}\right) & \leq\Omega_{0}^{K}\Omega_{1}\cdot\sup_{y\in\left[-1,1\right]^{\dimension}}\left(1+\left|y\right|\right)^{K}\\
 & \leq\left(2\sqrt{\dimension}\right)^{K}\Omega_{0}^{K}\Omega_{1}.
\end{align*}
All in all, Theorem \ref{thm:WienerAmalgamConvolution} shows
\begin{align*}
\left\Vert f\ast g\right\Vert _{W_{T_{j}^{-T}\left[-1,1\right]^{\dimension}}\left(L_{v}^{p}\right)} & =\left\Vert f\ast g\right\Vert _{W_{Q_{2}}\left(L_{v}^{p}\right)}\\
 & \leq3^{\frac{\dimension}{r}}\cdot\left(2\sqrt{\dimension}\right)^{K}\Omega_{0}^{K}\Omega_{1}\cdot\left[\lambda_{\dimension}\left(Q_{1}\right)\right]^{1-\frac{1}{r}}\cdot\left\Vert f\right\Vert _{W_{Q_{1}-Q_{1}}\left(L_{v_{0}}^{r}\right)}\cdot\left\Vert g\right\Vert _{W_{Q_{2}-Q_{1}}\left(L_{v}^{p}\right)}\\
\left({\scriptstyle Q_{2}-Q_{1}=Q_{1}-Q_{1}=T_{j}^{-T}\left[-2,2\right]^{\dimension}}\right) & \leq2^{K}3^{\frac{\dimension}{r}}\cdot\dimension^{\frac{K}{2}}\cdot2^{\dimension\left(1-\frac{1}{r}\right)}\cdot\Omega_{0}^{K}\Omega_{1}\!\cdot\!\left|\det T_{j}^{-T}\right|^{1-\frac{1}{r}}\!\cdot\left\Vert f\right\Vert _{W_{T_{j}^{-T}\left[-2,2\right]^{\dimension}}\left(L_{v_{0}}^{r}\right)}\cdot\left\Vert g\right\Vert _{W_{T_{j}^{-T}\left[-2,2\right]^{\dimension}}\left(L_{v}^{p}\right)}\\
\left({\scriptstyle \text{eq. }\eqref{eq:WienerLinearCubeEnlargement}}\right) & \leq\Omega_{0}^{3K}\Omega_{1}^{3}\cdot\dimension^{-\frac{\dimension}{2r}}\cdot\left(972\!\cdot\!\dimension^{\frac{5}{2}}\right)^{\!K+\frac{\dimension}{r}}\cdot\left|\det T_{j}\right|^{\frac{1}{r}-1}\!\cdot\left\Vert f\right\Vert _{W_{T_{j}^{-T}\left[-1,1\right]^{\dimension}}\left(L_{v_{0}}^{r}\right)}\cdot\left\Vert g\right\Vert _{W_{T_{j}^{-T}\left[-1,1\right]^{\dimension}}\left(L_{v}^{p}\right)}.\qedhere
\end{align*}
\end{proof}
Next, we establish a more quantitative—and weighted—version of the
convolution relation for (suitably) bandlimited functions given in
\cite[Corollary 3.14]{DecompositionEmbedding}, which is in turn a
specialized version of \cite[Proposition in §1.5.1]{TriebelTheoryOfFunctionSpaces}.

The following proposition uses the notation $Q_{i}^{n\ast}:=\bigcup_{\ell\in i^{n\ast}}Q_{\ell}$,
where $i^{1\ast}:=i^{\ast}$ and $i^{\left(n+1\right)\ast}:=\bigcup_{\ell\in i^{n\ast}}\ell^{\ast}$.
For the definition of $i^{\ast}$, cf.\@ equation (\ref{eq:IndexClusterDefinition}).
\begin{prop}
\label{prop:BandlimitedConvolution}Let $p\in\left(0,1\right]$ and
$n\in\N_{0}$. If $i\in I$ and

\begin{itemize}
\item if $\psi\in\TestFunctionSpace{\R^{\dimension}}$ with $\supp\psi\subset\overline{Q_{i}^{n\ast}}$
and
\item if $f\in\CalD'\left(\CalO\right)$ with $\supp f\subset\overline{Q_{i}^{n\ast}}$
and $\Fourier^{-1}f\in L_{v}^{p}\left(\R^{\dimension}\right)$,
\end{itemize}
then $\Fourier^{-1}\left(\psi f\right)=\left(\Fourier^{-1}\psi\right)\ast\left(\Fourier^{-1}f\right)\in L_{v}^{p}\left(\R^{\dimension}\right)$
with
\[
\left\Vert \Fourier^{-1}\left(\psi f\right)\right\Vert _{L_{v}^{p}}\leq\left[4R_{\CalQ}\cdot\left(3C_{\CalQ}\right)^{n}\right]^{\dimension\left(\frac{1}{p}-1\right)}\cdot\left|\det T_{i}\right|^{\frac{1}{p}-1}\cdot\left\Vert \Fourier^{-1}\psi\right\Vert _{L_{v_{0}}^{p}}\cdot\left\Vert \Fourier^{-1}f\right\Vert _{L_{v}^{p}}
\]
and
\[
\left\Vert \Fourier^{-1}\left(\psi f\right)\right\Vert _{W_{T_{i}^{-T}\left[-1,1\right]^{\dimension}}\left(L_{v}^{p}\right)}\leq C\cdot\left|\det T_{i}\right|^{\frac{1}{p}-1}\cdot\left\Vert \Fourier^{-1}\psi\right\Vert _{L_{v_{0}}^{p}}\cdot\left\Vert \Fourier^{-1}f\right\Vert _{L_{v}^{p}},
\]
where $C:=\Omega_{0}^{K}\Omega_{1}\cdot\left(2^{14}\!\cdot\!\dimension^{\frac{3}{2}}\!\cdot\!\left\lceil K\!+\!\frac{\dimension+1}{p}\right\rceil \right)^{\!\!K+\frac{\dimension+1}{p}+2}\left[1\!+\!4R_{\CalQ}\left(3C_{\CalQ}\right)^{n}\right]^{\dimension\left(\frac{2}{p}-1\right)}$.
\end{prop}
\begin{rem*}

\begin{itemize}[leftmargin=0.4cm]
\item Again, the only property of $v$ which we use is that $v$ is measurable
and $v\left(x+y\right)\leq v\left(x\right)v_{0}\left(y\right)$ for
all $x,y\in\R^{\dimension}$. Since this also holds for $v_{0}$ instead
of $v$, the claim also holds with $v$ replaced by $v_{0}$ everywhere.
\item Since $\overline{Q_{j}}\subset\CalO$ is compact for each $j\in I$,
the same is true of $\overline{Q_{i}^{n\ast}}\subset\CalO$. Hence,
the distribution $f\in\CalD'\left(\CalO\right)$ extends to a tempered
distribution $f\in\Schwartz'\left(\R^{\dimension}\right)$, so that
$\Fourier^{-1}f$ is well-defined and such that $\psi f\in\Schwartz'\left(\R^{\dimension}\right)$
is a tempered distribution with compact support, since $\psi\in\TestFunctionSpace{\R^{\dimension}}$.
Finally, it follows from \cite[Proposition 2.3.22(11)]{GrafakosClassicalFourierAnalysis}
that $\Fourier^{-1}\left(\psi f\right)=\Fourier^{-1}\psi\ast\Fourier^{-1}f$.\qedhere
\end{itemize}
\end{rem*}
\begin{proof}
First, we note that \cite[Lemma 2.7]{DecompositionEmbedding} yields
\[
Q_{j}\subset T_{i}\left[\overline{B_{\left(1+2C_{\CalQ}\right)^{n}R_{\CalQ}}}\left(0\right)\right]+b_{i}\qquad\forall j\in i^{n\ast}.
\]
Hence, setting $R:=\left(1+2C_{\CalQ}\right)^{n}R_{\CalQ}$, we have
\begin{equation}
\overline{Q_{i}^{n\ast}}\subset T_{i}\overline{B_{R}}\left(0\right)+b_{i}\subset T_{i}\left[-R,R\right]^{\dimension}+b_{i}=:\Omega.\label{eq:BandlimitedConvolutionStarredSetInclusion}
\end{equation}
Note that, once we have proved the first claimed estimate, the second
one is a consequence of Theorem \ref{thm:BandlimitedWienerAmalgamSelfImproving}
(and some simple estimates of the resulting constant, using $C_{\CalQ}\geq\left\Vert T_{i}^{-1}T_{i}\right\Vert =1$
and $s_{\dimension}\leq2^{2\dimension}$), since we have $\supp\Fourier\left[\Fourier^{-1}\left(\psi f\right)\right]\subset\supp\psi\subset\overline{Q_{i}^{n\ast}}\subset\Omega$.

As seen in the remark following the proposition, we have $\Fourier^{-1}f\in\Schwartz'\left(\R^{\dimension}\right)$
with $\supp\Fourier\left[\Fourier^{-1}f\right]\subset\overline{Q_{i}^{n\ast}}\subset\Omega$
and likewise $\Fourier^{-1}\psi\in\Schwartz\left(\R^{\dimension}\right)\subset\Schwartz'\left(\R^{\dimension}\right)$
with $\supp\Fourier\left[\Fourier^{-1}\psi\right]\subset\overline{Q_{i}^{n\ast}}\subset\Omega$.
In view of Theorems \ref{thm:BandlimitedWienerAmalgamSelfImproving}
and \ref{thm:WienerAmalgamConvolution}, we thus get $\Fourier^{-1}\psi\in W_{T_{i}^{-T}\left[-1,1\right]^{\dimension}}\left(L_{v_{0}}^{p}\right)\hookrightarrow L_{v_{0}}^{1}\left(\R^{\dimension}\right)$
(cf.\@ Lemmas \ref{lem:SchwartzFunctionsAreWiener} and \ref{lem:WienerAmalgamNormEquivalence})
and $\Fourier^{-1}f\in W_{T_{i}^{-T}\left[-1,1\right]^{\dimension}}\left(L_{v}^{p}\right)\hookrightarrow L_{v}^{\infty}\left(\R^{\dimension}\right)$,
so that $\Fourier^{-1}\psi\ast\Fourier^{-1}f$ is pointwise well-defined
by Corollary \ref{cor:WienerAmalgamConvolutionSimplified}.

Now, Lemma \ref{lem:BandlimitedPointwiseApproximation} ensures existence
of a sequence $\left(h_{n}\right)_{n\in\N}$ of Schwartz functions
satisfying $\left|h_{n}\left(x\right)\right|\leq\left|\left(\Fourier^{-1}f\right)\left(x\right)\right|$,
as well as $h_{n}\left(x\right)\xrightarrow[n\to\infty]{}\left(\Fourier^{-1}f\right)\left(x\right)$
for all $x\in\R^{\dimension}$ and finally $\supp\widehat{h_{n}}\subset B_{1/n}\left(\Omega\right)$
for all $n\in\N$. It is not hard to see $\Omega-B_{1/n}\left(\Omega\right)\subset B_{1/n}\left(\Omega-\Omega\right)$.
Furthermore, by compactness of $\Omega-\Omega$—and using continuity
of the Lebesgue measure from above, cf.\@ \cite[Theorem 1.8(d)]{FollandRA}—we
get 
\begin{align*}
\lambda_{\dimension}\left(B_{1/n}\left(\Omega-\Omega\right)\right)\xrightarrow[n\to\infty]{}\lambda_{\dimension}\left(\Omega-\Omega\right) & =\lambda_{\dimension}\left(\left[T_{i}\left[-R,R\right]^{\dimension}+b_{i}\right]-\left[T_{i}\left[-R,R\right]^{\dimension}+b_{i}\right]\right)\\
 & \leq\lambda_{\dimension}\left(T_{i}\left[-2R,2R\right]^{\dimension}\right)=\left(4R\right)^{\dimension}\cdot\left|\det T_{i}\right|.
\end{align*}

Now, since $h_{n},\Fourier^{-1}\psi\in\Schwartz\left(\R^{\dimension}\right)\subset L^{p}\left(\R^{\dimension}\right)$,
Theorem \ref{thm:PointwiseQuasiBanachBandlimitedConvolution} yields
\begin{align*}
v\left(x\right)\cdot\left(\left|\Fourier^{-1}\psi\right|\ast\left|h_{n}\right|\right)\left(x\right) & \leq\left[\lambda_{\dimension}\left(\supp\Fourier\left[\Fourier^{-1}\psi\right]-\supp\widehat{h_{n}}\right)\right]^{\frac{1}{p}-1}\cdot\left[\int_{\R^{\dimension}}\left[v\left(x\right)\right]^{p}\cdot\left|\Fourier^{-1}\psi\left(x-y\right)\right|^{p}\cdot\left|h_{n}\left(y\right)\right|^{p}\d y\right]^{1/p}\\
\left({\scriptstyle \text{since }v\left(x\right)\leq v_{0}\left(x-y\right)v\left(y\right)}\right) & \leq\left[\lambda_{\dimension}\left(B_{1/n}\left(\Omega-\Omega\right)\right)\right]^{\frac{1}{p}-1}\cdot\left[\int_{\R^{\dimension}}\left|\left(v_{0}\cdot\Fourier^{-1}\psi\right)\left(x-y\right)\right|^{p}\cdot\left|\left(v\cdot h_{n}\right)\left(y\right)\right|^{p}\d y\right]^{1/p}\\
\left({\scriptstyle \text{since }\left|h_{n}\right|\leq\left|\Fourier^{-1}f\right|}\right) & \leq\left[\lambda_{\dimension}\left(B_{1/n}\left(\Omega-\Omega\right)\right)\right]^{\frac{1}{p}-1}\cdot\left[\left(\left|v_{0}\cdot\Fourier^{-1}\psi\right|^{p}\ast\left|v\cdot\Fourier^{-1}f\right|^{p}\right)\left(x\right)\right]^{1/p}.
\end{align*}
Taking the limes inferior on both sides, we get
\[
\liminf_{n\to\infty}\left[v\left(x\right)\cdot\left(\left|\Fourier^{-1}\psi\right|\ast\left|h_{n}\right|\right)\left(x\right)\right]\leq\left[\left(4R\right)^{\dimension}\cdot\left|\det T_{i}\right|\right]^{\frac{1}{p}-1}\cdot\left[\left(\left|v_{0}\cdot\Fourier^{-1}\psi\right|^{p}\ast\left|v\cdot\Fourier^{-1}f\right|^{p}\right)\left(x\right)\right]^{1/p}.
\]

Next, since $h_{n}\to\Fourier^{-1}f$ pointwise, and since we saw
above that $\Fourier^{-1}\psi\ast\Fourier^{-1}f$ is pointwise well-defined,
we get by Fatou's Lemma that
\begin{align*}
\left|\left(\Fourier^{-1}\psi\ast\Fourier^{-1}f\right)\left(x\right)\right| & \leq\int_{\R^{\dimension}}\left|\left(\Fourier^{-1}\psi\right)\left(y\right)\right|\cdot\left|\left(\Fourier^{-1}f\right)\left(x-y\right)\right|\d y\\
 & =\int_{\R^{\dimension}}\liminf_{n\to\infty}\left[\left|\left(\Fourier^{-1}\psi\right)\left(y\right)\right|\cdot\left|h_{n}\left(x-y\right)\right|\right]\d y\\
 & \leq\liminf_{n\to\infty}\int_{\R^{\dimension}}\left|\left(\Fourier^{-1}\psi\right)\left(y\right)\right|\cdot\left|h_{n}\left(x-y\right)\right|\d y=\liminf_{n\to\infty}\left(\left|\Fourier^{-1}\psi\right|\ast\left|h_{n}\right|\right)\left(x\right)
\end{align*}
for all $x\in\R^{\dimension}$. Hence, we finally see
\begin{align*}
\left\Vert \Fourier^{-1}\left(\psi f\right)\right\Vert _{L_{v}^{p}} & =\left\Vert \Fourier^{-1}\psi\ast\Fourier^{-1}f\right\Vert _{L_{v}^{p}}\leq\left\Vert x\mapsto\liminf_{n\to\infty}v\left(x\right)\cdot\left(\left|\Fourier^{-1}\psi\right|\ast\left|h_{n}\right|\right)\left(x\right)\right\Vert _{L^{p}}\\
 & \leq\left[\left(4R\right)^{\dimension}\cdot\left|\det T_{i}\right|\right]^{\frac{1}{p}-1}\cdot\left\Vert x\mapsto\left[\left(\left|v_{0}\cdot\Fourier^{-1}\psi\right|^{p}\ast\left|v\cdot\Fourier^{-1}f\right|^{p}\right)\left(x\right)\right]^{1/p}\right\Vert _{L^{p}}\\
 & =\left[\left(4R\right)^{\dimension}\cdot\left|\det T_{i}\right|\right]^{\frac{1}{p}-1}\cdot\left\Vert \left|v_{0}\cdot\Fourier^{-1}\psi\right|^{p}\ast\left|v\cdot\Fourier^{-1}f\right|^{p}\right\Vert _{L^{1}}^{1/p}\\
\left({\scriptstyle \text{Young's inequality}}\right) & \leq\left[\left(4R\right)^{\dimension}\cdot\left|\det T_{i}\right|\right]^{\frac{1}{p}-1}\cdot\left\Vert \left|v_{0}\cdot\Fourier^{-1}\psi\right|^{p}\right\Vert _{L^{1}}^{1/p}\cdot\left\Vert \left|v\cdot\Fourier^{-1}f\right|^{p}\right\Vert _{L^{1}}^{1/p}\\
 & =\left[\left(4R\right)^{\dimension}\cdot\left|\det T_{i}\right|\right]^{\frac{1}{p}-1}\cdot\left\Vert \Fourier^{-1}\psi\right\Vert _{L_{v_{0}}^{p}}\cdot\left\Vert \Fourier^{-1}f\right\Vert _{L_{v}^{p}}<\infty.
\end{align*}
Since we have $C_{\CalQ}\geq\left\Vert T_{i}^{-1}T_{i}\right\Vert =1$,
we get $R=\left(1+2C_{\CalQ}\right)^{n}R_{\CalQ}\leq\left(3C_{\CalQ}\right)^{n}R_{\CalQ}$,
which easily yields the claim.
\end{proof}
As our last result in this section, we show—as a consequence of our
developed convolution relations—that the decomposition space $\DecompSp{\CalQ}p{\ell_{w}^{q}}v$
is well-defined, even if $v\not\equiv1$.
\begin{prop}
\label{prop:WeightedDecompositionSpaceWellDefined}Let $\Phi=\left(\varphi_{i}\right)_{i\in I}$
and $\Psi=\left(\psi_{i}\right)_{i\in I}$ be two $\CalQ$-$v_{0}$-BAPUs.
Then we have
\[
\left\Vert \left(\left\Vert \Fourier^{-1}\left(\varphi_{i}f\right)\right\Vert _{L_{v}^{p}}\right)_{i\in I}\right\Vert _{\ell_{w}^{q}}\asymp\left\Vert \left(\left\Vert \Fourier^{-1}\left(\psi_{i}f\right)\right\Vert _{L_{v}^{p}}\right)_{i\in I}\right\Vert _{\ell_{w}^{q}}
\]
uniformly over $f\in\DistributionSpace{\CalO}$. In particular, the
decomposition space $\DecompSp{\CalQ}p{\ell_{w}^{q}}v$ is independent
of the choice of the $\CalQ$-$v_{0}$-BAPU.
\end{prop}
\begin{proof}
By symmetry, it suffices to establish the estimate ``$\lesssim$''.
We can clearly assume $\left\Vert \left(\left\Vert \Fourier^{-1}\left(\psi_{i}f\right)\right\Vert _{L_{v}^{p}}\right)_{i\in I}\right\Vert _{\ell_{w}^{q}}<\infty$.
Since $L_{v}^{p}\left(\R^{\dimension}\right)$ is a quasi-normed space
and since we have the uniform estimate $\left|i^{\ast}\right|\leq N_{\CalQ}$
for all $i\in I$, we have
\[
d_{i}:=\left\Vert \Fourier^{-1}\left(\psi_{i}^{\ast}f\right)\right\Vert _{L_{v}^{p}}\leq C\cdot\sum_{\ell\in i^{\ast}}\left\Vert \Fourier^{-1}\left(\psi_{\ell}f\right)\right\Vert _{L_{v}^{p}}=C\cdot\left(\Gamma_{\CalQ}e\right)_{i},
\]
for a suitable constant $C=C\left(p,N_{\CalQ}\right)$, where $e=\left(e_{i}\right)_{i\in I}$
is defined by $e_{i}:=\left\Vert \Fourier^{-1}\left(\psi_{i}f\right)\right\Vert _{L_{v}^{p}}$
and where $\Gamma_{\CalQ}$ is the $\CalQ$-clustering map, as defined
in Section \ref{subsec:DecompSpaceDefinitionStandingAssumptions},
equation (\ref{eq:QClusteringMapDefinition}).

Now, as seen in Section \ref{subsec:DecompSpaceDefinitionStandingAssumptions},
we have $\psi_{i}^{\ast}\equiv1$ on $Q_{i}$ and thus $\varphi_{i}=\psi_{i}^{\ast}\varphi_{i}$
for all $i\in I$. Hence,
\[
c_{i}:=\left\Vert \Fourier^{-1}\left(\varphi_{i}f\right)\right\Vert _{L_{v}^{p}}=\left\Vert \Fourier^{-1}\left(\varphi_{i}\psi_{i}^{\ast}f\right)\right\Vert _{L_{v}^{p}}=\left\Vert \left[\Fourier^{-1}\varphi_{i}\right]\ast\Fourier^{-1}\left(\psi_{i}^{\ast}f\right)\right\Vert _{L_{v}^{p}}.
\]
In case of $p\in\left[1,\infty\right]$, we can now use the weighted
Young inequality (equation (\ref{eq:WeightedYoungInequality})) to
derive
\[
c_{i}\leq\left\Vert \Fourier^{-1}\varphi_{i}\right\Vert _{L_{v_{0}}^{1}}\cdot\left\Vert \Fourier^{-1}\left(\psi_{i}^{\ast}f\right)\right\Vert _{L_{v}^{p}}\leq C\cdot C_{\CalQ,\Phi,v_{0},p}\cdot\left(\Gamma_{\CalQ}e\right)_{i}.
\]
Otherwise, if $p\in\left(0,1\right)$, we use Proposition \ref{prop:BandlimitedConvolution}
(with $n=1$, since $\supp\varphi_{i}\subset\overline{Q_{i}^{\ast}}$
and $\supp\psi_{i}^{\ast}\subset\overline{Q_{i}^{\ast}}$) to derive
\begin{align*}
c_{i}=\left\Vert \Fourier^{-1}\left(\varphi_{i}\psi_{i}^{\ast}f\right)\right\Vert _{L_{v}^{p}} & \leq\left[12R_{\CalQ}C_{\CalQ}\right]^{\dimension\left(\frac{1}{p}-1\right)}\cdot\left|\det T_{i}\right|^{\frac{1}{p}-1}\cdot\left\Vert \Fourier^{-1}\varphi_{i}\right\Vert _{L_{v_{0}}^{p}}\cdot\left\Vert \Fourier^{-1}\left[\psi_{i}^{\ast}f\right]\right\Vert _{L_{v}^{p}}\\
 & \leq C\cdot\left[12R_{\CalQ}C_{\CalQ}\right]^{\dimension\left(\frac{1}{p}-1\right)}C_{\CalQ,\Phi,v_{0},p}\cdot\left(\Gamma_{\CalQ}e\right)_{i}.
\end{align*}

In summary, there is for arbitrary $p\in\left(0,\infty\right]$ a
constant $C'=C'\left(\CalQ,p,\Phi,v_{0}\right)>0$ satisfying $c_{i}\leq C'\cdot\left(\Gamma_{\CalQ}e\right)_{i}<\infty$
for all $i\in I$. By solidity of $\ell_{w}^{q}\left(I\right)$ and
by boundedness of $\Gamma_{\CalQ}$, this implies
\[
\left\Vert \left(\left\Vert \Fourier^{-1}\left(\varphi_{i}f\right)\right\Vert _{L_{v}^{p}}\right)_{i\in I}\right\Vert _{\ell_{w}^{q}}\leq C'\cdot\left\Vert \Gamma_{\CalQ}e\right\Vert _{\ell_{w}^{q}}\leq C'\cdot\vertiii{\smash{\Gamma_{\CalQ}}}\cdot\left\Vert e\right\Vert _{\ell_{w}^{q}}=C'\cdot\vertiii{\smash{\Gamma_{\CalQ}}}\cdot\left\Vert \left(\left\Vert \Fourier^{-1}\left(\psi_{i}f\right)\right\Vert _{L_{v}^{p}}\right)_{i\in I}\right\Vert _{\ell_{w}^{q}}.\qedhere
\]
\end{proof}

\section{Semi-discrete Banach Frames}

\label{sec:SemiDiscreteBanachFrames}
\begin{assumption}
\label{assu:MainAssumptions}In the remainder of the paper, we will
use the following assumptions and notations:

\begin{enumerate}
\item We are given a family $\Gamma=\left(\gamma_{i}\right)_{i\in I}$ of
functions $\gamma_{i}:\R^{\dimension}\to\Compl$ with the following
additional properties:

\begin{enumerate}
\item We have $\gamma_{i}\in L_{\left(1+\left|\mybullet\right|\right)^{K}}^{1}\left(\R^{d}\right)\hookrightarrow L_{v_{0}}^{1}\left(\R^{\dimension}\right)\hookrightarrow L^{1}\left(\R^{\dimension}\right)$
for all $i\in I$.
\item We have $\widehat{\gamma_{i}}\in C^{\infty}\left(\R^{d}\right)$ for
all $i\in I$, where all partial derivatives of $\widehat{\gamma_{i}}$
are polynomially bounded, i.e.,
\[
\qquad\qquad\qquad\left|\left(\partial^{\alpha}\widehat{\gamma_{i}}\right)\left(\xi\right)\right|\leq C_{\alpha,i}\cdot\left(1+\left|\xi\right|\right)^{N_{\alpha,i}}\qquad\forall\,\xi\in\R^{d}\,\forall\,\alpha\in\N_{0}^{d}\,\forall\,i\in I,\text{ for suitable }C_{\alpha,i}>0\text{ and }N_{\alpha,i}\in\N_{0}.
\]
\end{enumerate}
\item For $i\in I$, we define
\begin{equation}
\begin{split}\gamma^{\left(i\right)} & :=\Fourier^{-1}\left(\widehat{\gamma_{i}}\circ S_{i}^{-1}\right)\\
 & =\Fourier^{-1}\left(L_{b_{i}}\left(\widehat{\gamma_{i}}\circ T_{i}^{-1}\right)\right)\\
 & =M_{b_{i}}\left[\Fourier^{-1}\left(\widehat{\gamma_{i}}\circ T_{i}^{-1}\right)\right]\\
 & =\left|\det T_{i}\right|\cdot M_{b_{i}}\left[\gamma_{i}\circ T_{i}^{T}\right],
\end{split}
\label{eq:NonCompactFilterDefinition}
\end{equation}
as well as the $L^{2}$-normalized version
\begin{equation}
\gamma^{\left[i\right]}:=\left|\det T_{i}\right|^{1/2}\cdot M_{b_{i}}\left[\gamma_{i}\circ T_{i}^{T}\right]=\left|\det T_{i}\right|^{-1/2}\cdot\gamma^{\left(i\right)}.\label{eq:L2NormalizedFilterDefinition}
\end{equation}
\item For $i\in I$, we set
\[
V_{i}:=\begin{cases}
L_{v}^{p}\left(\R^{\dimension}\right), & \text{if }p\in\left[1,\infty\right],\\
W_{T_{i}^{-T}\left[-1,1\right]^{\dimension}}\left(L_{v}^{p}\right), & \text{if }p\in\left(0,1\right).
\end{cases}
\]
Furthermore, we will occasionally make use of the space
\[
V:=\ell_{w}^{q}\left(\left[V_{i}\right]_{i\in I}\right):=\left\{ \left(f_{i}\right)_{i\in I}\with\left(\forall i\in I:\,f_{i}\in V_{i}\right)\text{ and }\left(\left\Vert f_{i}\right\Vert _{V_{i}}\right)_{i\in I}\in\ell_{w}^{q}\left(I\right)\right\} ,
\]
equipped with the quasi-norm $\left\Vert \left(f_{i}\right)_{i\in I}\right\Vert _{\ell_{w}^{q}\left(\left[V_{i}\right]_{i\in I}\right)}:=\left\Vert \left(\left\Vert f_{i}\right\Vert _{V_{i}}\right)_{i\in I}\right\Vert _{\ell_{w}^{q}}$.
\item Finally, we set
\[
r:=\max\left\{ q,\frac{q}{p}\right\} =\begin{cases}
q, & \text{if }p\in\left[1,\infty\right],\\
\frac{q}{p}, & \text{if }p\in\left(0,1\right)
\end{cases}
\]
and
\[
A_{j,i}:=\begin{cases}
\left\Vert \Fourier^{-1}\left(\varphi_{i}\cdot\widehat{\gamma^{\left(j\right)}}\right)\right\Vert _{L_{v_{0}}^{1}}, & \text{if }p\in\left[1,\infty\right],\\
\left(1+\left\Vert T_{j}^{-1}T_{i}\right\Vert \right)^{\dimension}\cdot\left|\det T_{i}\right|^{1-p}\cdot\left\Vert \Fourier^{-1}\left(\varphi_{i}\cdot\widehat{\gamma^{\left(j\right)}}\right)\right\Vert _{L_{v_{0}}^{p}}^{p}, & \text{if }p\in\left(0,1\right)
\end{cases}
\]
for $i,j\in I$ and we assume that $\overrightarrow{A}$ is a bounded
operator $\overrightarrow{A}:\ell_{w^{\min\left\{ 1,p\right\} }}^{r}\left(I\right)\to\ell_{w^{\min\left\{ 1,p\right\} }}^{r}\left(I\right)$,
where
\[
\overrightarrow{A}\left(c_{i}\right)_{i\in I}:=\left(\sum_{i\in I}A_{j,i}\,c_{i}\right)_{j\in I}.\qedhere
\]
\end{enumerate}
\end{assumption}
\begin{rem}
\label{rem:MainAssumptionsRemark}

\begin{enumerate}
\item The most common case will be to have $\gamma_{i}=\gamma$ for all
$i\in I$, for a fixed prototype $\gamma$. The added flexibility
of allowing $\gamma_{i}$ to vary with $i\in I$ is only rarely needed.
In the cases where it is, we usually have $\gamma_{i}=\gamma_{n_{i}}$
with a given (finite) list of prototypes $\gamma_{1},\dots,\gamma_{N}$.
\item The assumption that $\partial^{\alpha}\widehat{\gamma_{i}}$ is polynomially
bounded for all $\alpha\in\N_{0}^{\dimension}$ is satisfied if $\gamma_{i}\in L^{1}\left(\R^{\dimension}\right)$
has compact support, say $\supp\gamma_{i}\subset\left[-R,R\right]^{\dimension}$
with $R\geq1$, since then differentiation under the integral yields
\begin{align*}
\left|\partial^{\alpha}\widehat{\gamma_{i}}\left(\xi\right)\right| & =\left|\int_{\R^{\dimension}}\gamma_{i}\left(x\right)\cdot\partial_{\xi}^{\alpha}e^{-2\pi i\left\langle x,\xi\right\rangle }\d x\right|\\
 & \leq\int_{\left[-R,R\right]^{\dimension}}\left|\gamma_{i}\left(x\right)\right|\cdot\left(2\pi\left|x\right|\right)^{\left|\alpha\right|}\d x\leq\left(2\pi R\right)^{\left|\alpha\right|}\cdot\left\Vert \gamma_{i}\right\Vert _{L^{1}}<\infty
\end{align*}
for all $\xi\in\R^{\dimension}$ and arbitrary $\alpha\in\N_{0}^{\dimension}$.
\item \label{enu:StructuredFamilyFourierTransformPolynomiallyBounded}Under
the above assumptions, the chain rule implies
\begin{align*}
\left|\left(\partial^{\alpha}\widehat{\gamma^{\left(i\right)}}\right)\left(x\right)\right| & =\left|\left(\partial^{\alpha}\left[\widehat{\gamma_{i}}\circ T_{i}^{-1}\right]\right)\left(x-b_{i}\right)\right|\\
 & \leq C^{\left(\alpha\right)}\cdot\left\Vert T_{i}^{-1}\right\Vert ^{\left|\alpha\right|}\cdot\max_{\left|\beta\right|\leq\left|\alpha\right|}\left|\left(\partial^{\beta}\widehat{\gamma_{i}}\right)\left(T_{i}^{-1}\left(x-b_{i}\right)\right)\right|\\
 & \leq\left(\max_{\left|\beta\right|\leq\left|\alpha\right|}C_{\beta,i}\right)\cdot C^{\left(\alpha\right)}\cdot\left\Vert T_{i}^{-1}\right\Vert ^{\left|\alpha\right|}\cdot\max_{\left|\beta\right|\leq\left|\alpha\right|}\left(1+\left|T_{i}^{-1}\left(x-b_{i}\right)\right|\right)^{N_{\beta,i}}\\
\left({\scriptstyle \text{with }N_{\alpha,i}':=\max_{\left|\beta\right|\leq\left|\alpha\right|}N_{\beta,i}}\right) & \leq C_{\alpha,i}'\cdot\left(1+\left|T_{i}^{-1}\left(x-b_{i}\right)\right|\right)^{N_{\alpha,i}'}\\
\left({\scriptstyle \text{for suitable }C_{\alpha,i}''>0}\right) & \leq C_{\alpha,i}''\cdot\left(1+\left|x\right|\right)^{N_{\alpha,i}'},
\end{align*}
where the last step used
\begin{align*}
1+\left|T_{i}^{-1}\left(x-b_{i}\right)\right| & \leq1+\left\Vert T_{i}^{-1}\right\Vert \left|x-b_{i}\right|\\
 & \leq1+\left\Vert T_{i}^{-1}\right\Vert \left|b_{i}\right|+\left\Vert T_{i}^{-1}\right\Vert \left|x\right|\\
 & \leq\left(1+\left\Vert T_{i}^{-1}\right\Vert \left|b_{i}\right|+\left\Vert T_{i}^{-1}\right\Vert \right)\cdot\left(1+\left|x\right|\right).
\end{align*}
Hence, all partial derivatives of each $\widehat{\gamma^{\left(i\right)}}$
are polynomially bounded.
\item Using $\gamma_{i}\in L_{\left(1+\left|\mybullet\right|\right)^{K}}^{1}\left(\R^{\dimension}\right)$,
we also get
\begin{align*}
\left\Vert \gamma^{\left(i\right)}\right\Vert _{L_{\left(1+\left|\mybullet\right|\right)^{K}}^{1}} & =\left|\det T_{i}\right|\cdot\left\Vert \left(1+\left|\mybullet\right|\right)^{K}\cdot\left(\gamma_{i}\circ T_{i}^{T}\right)\right\Vert _{L^{1}}\\
 & =\left\Vert \left(1+\left|T_{i}^{-T}\mybullet\right|\right)^{K}\cdot\gamma_{i}\right\Vert _{L^{1}}\\
\left({\scriptstyle \text{eq. }\eqref{eq:WeightLinearTransformationsConnection}}\right) & \leq\Omega_{0}^{K}\cdot\left\Vert \left(1+\left|\mybullet\right|\right)^{K}\cdot\gamma_{i}\right\Vert _{L^{1}}=\Omega_{0}^{K}\cdot\left\Vert \gamma_{i}\right\Vert _{L_{\left(1+\left|\mybullet\right|\right)^{K}}^{1}}<\infty
\end{align*}
and thus $\gamma^{\left(i\right)}\in L_{\left(1+\left|\mybullet\right|\right)^{K}}^{1}\left(\R^{\dimension}\right)\hookrightarrow L_{v_{0}}^{1}\left(\R^{\dimension}\right)\hookrightarrow L^{1}\left(\R^{\dimension}\right)$,
where the last embedding uses $v_{0}\geq1$.
\item Point (\ref{enu:StructuredFamilyFourierTransformPolynomiallyBounded})
from above shows $\widehat{\gamma^{\left(i\right)}}\cdot\widehat{f}\in\Schwartz'\left(\R^{d}\right)$
for arbitrary $f\in\Schwartz'\left(\R^{d}\right)$, so that $\gamma^{\left(i\right)}\ast f:=\Fourier^{-1}\left(\widehat{\gamma^{\left(i\right)}}\cdot\widehat{f}\right)$
is a well-defined tempered distribution. Of course, the same also
holds for $\gamma^{\left[i\right]}\ast f:=\Fourier^{-1}\left(\widehat{\gamma^{\left[i\right]}}\cdot\widehat{f}\right)$.
\item Since $\R^{\dimension}$ is $\sigma$-compact, it follows from \cite[Lemma 2.3.7]{VoigtlaenderPhDThesis}
(see also \cite[Theorem 2.6]{RauhutWienerAmalgam}) that for $p\in\left(0,1\right)$,
each of the spaces $V_{i}=W_{T_{i}^{-T}\left[-1,1\right]^{\dimension}}\left(L_{v}^{p}\right)$
is complete (and thus a Quasi-Banach space) for each $i\in I$. Furthermore,
\cite[Lemma 2.3.4]{VoigtlaenderPhDThesis} and \cite[Exercise 1.1.5(c)]{GrafakosClassicalFourierAnalysis}
show $\left\Vert f+g\right\Vert _{V_{i}}\leq2^{\frac{1}{p}-1}\cdot\left[\left\Vert f\right\Vert _{V_{i}}+\left\Vert g\right\Vert _{V_{i}}\right]$
for all $f,g\in V_{i}$. In case of $p\in\left[1,\infty\right]$,
it is clear that $V_{i}=L_{v}^{p}\left(\R^{\dimension}\right)$ is
a Banach space.\qedhere
\end{enumerate}
\end{rem}
Note that in the preceding remark, we \emph{defined} $\gamma^{\left(i\right)}\ast f:=\Fourier^{-1}\left(\widehat{\gamma^{\left(i\right)}}\cdot\widehat{f}\right)$.
We needed to do so, since the usual results about convolution in $\Schwartz'\left(\R^{\dimension}\right)$
only define $f\ast\varphi$ for $\varphi\in\Schwartz'\left(\R^{\dimension}\right)$
if $f\in\Schwartz\left(\R^{\dimension}\right)$ (cf.\@ \cite[Proposition (8.44)]{FollandRA})
or if $f$ is a distribution with compact support. (cf.\@ \cite[Chapter 8, Exercise 35]{FollandRA}).
But note that if we not only know $\varphi\in\Schwartz'\left(\R^{\dimension}\right)$,
but the stronger property $\varphi\in\left(L^{1}+L^{\infty}\right)\left(\R^{\dimension}\right)$
and if $f\in L^{1}\left(\R^{\dimension}\right)$, then $f\ast\varphi\in\left(L^{1}+L^{\infty}\right)\left(\R^{\dimension}\right)$
is already defined. Our next result shows that in this (and in a slightly
more general) case, the new definition is consistent.
\begin{lem}
\label{lem:SpecialConvolutionConsistent}Assume $\varphi\in L_{v}^{1}\left(\R^{\dimension}\right)+L_{v}^{\infty}\left(\R^{\dimension}\right)$
and assume that $f\in L_{v_{0}}^{1}\left(\R^{\dimension}\right)$
is such that $\widehat{f}\in C^{\infty}\left(\R^{\dimension}\right)$
and such that all partial derivatives of $\widehat{f}$ have at most
polynomial growth. Then $f\ast\varphi\in L_{v}^{1}\left(\R^{\dimension}\right)+L_{v}^{\infty}\left(\R^{\dimension}\right)\hookrightarrow\Schwartz'\left(\R^{\dimension}\right)$
and
\[
f\ast\varphi=\Fourier^{-1}\left[\widehat{f}\cdot\widehat{\varphi}\right].
\]
The assumption on $\varphi$ is in particular fulfilled if $\varphi\in V_{i}$
for some $i\in I$. More precisely,
\begin{equation}
V_{i}\hookrightarrow L_{v}^{1}\left(\smash{\R^{\dimension}}\right)+L_{v}^{\infty}\left(\smash{\R^{\dimension}}\right)\hookrightarrow\Schwartz'\left(\smash{\R^{\dimension}}\right).\qedhere\label{eq:ViAreTemperedDistributions}
\end{equation}
\end{lem}
\begin{rem}
\label{rem:WeightedSpacesYieldTemperedDistributions}We saw in the
proof of Theorem \ref{thm:WienerAmalgamConvolution} (cf.\@ equation
(\ref{eq:WeightBoundedBelow})) that $v\left(x\right)\gtrsim\left(1+\left|x\right|\right)^{-K}$.
Furthermore, $v_{0}\geq1$, and thus $L_{v}^{p}\left(\R^{\dimension}\right)\hookrightarrow\Schwartz'\left(\R^{\dimension}\right)$
and $L_{v_{0}}^{p}\left(\R^{\dimension}\right)\hookrightarrow L^{p}\left(\R^{\dimension}\right)\hookrightarrow\Schwartz'\left(\R^{\dimension}\right)$
for all $p\in\left[1,\infty\right]$. Hence, the expressions $\widehat{f}$
and $\widehat{\varphi}$ above are well-defined tempered distributions.

Since $\widehat{f}\in C^{\infty}\left(\R^{\dimension}\right)$ with
all derivatives of $\widehat{f}$ of at most polynomial growth, we
also see $\widehat{f}\cdot\widehat{\varphi}\in\Schwartz'\left(\R^{\dimension}\right)$
and thus also $\Fourier^{-1}\left[\widehat{f}\cdot\widehat{\varphi}\right]\in\Schwartz'\left(\R^{\dimension}\right)$.
\end{rem}
\begin{proof}
From the weighted Young inequality (equation (\ref{eq:WeightedYoungInequality})),
we know  $L_{v_{0}}^{1}\left(\R^{\dimension}\right)\ast L_{v}^{\infty}\left(\R^{\dimension}\right)\hookrightarrow L_{v}^{\infty}\left(\R^{\dimension}\right)$.
Likewise, the same inequality also yields $\left\Vert f\ast g\right\Vert _{L_{v}^{1}}\leq\left\Vert f\right\Vert _{L_{v_{0}}^{1}}\cdot\left\Vert g\right\Vert _{L_{v}^{1}}<\infty$
and thus in particular $\left(\left|f\right|\ast\left|g\right|\right)\left(x\right)<\infty$
for almost all $x\in\R^{\dimension}$ for $f\in L_{v_{0}}^{1}\left(\R^{\dimension}\right)$
and $g\in L_{v}^{1}\left(\R^{\dimension}\right)$. Hence, together
with Remark \ref{rem:WeightedSpacesYieldTemperedDistributions}, we
see indeed that $f\ast\varphi\in L_{v}^{1}\left(\R^{\dimension}\right)+L_{v}^{\infty}\left(\R^{\dimension}\right)\hookrightarrow\Schwartz'\left(\R^{\dimension}\right)$.

\medskip{}

Now, let $\psi\in\Schwartz\left(\R^{\dimension}\right)$ be arbitrary.
We have
\begin{align*}
\left\langle \Fourier^{-1}\left[\widehat{f}\cdot\widehat{\varphi}\right],\psi\right\rangle _{\Schwartz',\Schwartz} & =\left\langle \widehat{f}\cdot\widehat{\varphi},\,\Fourier^{-1}\psi\right\rangle _{\Schwartz',\Schwartz}\\
 & =\left\langle \widehat{\varphi},\,\widehat{f}\cdot\Fourier^{-1}\psi\right\rangle _{\Schwartz',\Schwartz}\\
 & =\left\langle \varphi,\,\Fourier\left[\widehat{f}\cdot\Fourier^{-1}\psi\right]\right\rangle _{\Schwartz',\Schwartz}\\
 & =\left\langle \varphi,\,\left(\Fourier^{-1}\left[\widehat{f}\cdot\Fourier^{-1}\psi\right]\right)\left(-\mybullet\right)\right\rangle _{\Schwartz',\Schwartz}.
\end{align*}
Recall that $v_{0}\geq1$, so that $f\in L_{v_{0}}^{1}\left(\R^{\dimension}\right)\hookrightarrow L^{1}\left(\R^{\dimension}\right)$.
Hence, $h:=f\ast\tilde{\psi}\in L^{1}\left(\R^{\dimension}\right)\ast L^{1}\left(\R^{\dimension}\right)\subset L^{1}\left(\R^{\dimension}\right)$,
where $\tilde{\psi}\left(x\right):=\psi\left(-x\right)$. Thus, the
convolution theorem yields $\widehat{h}=\widehat{f}\cdot\widehat{\tilde{\psi}}=\widehat{f}\cdot\tilde{\widehat{\psi}}=\widehat{f}\cdot\Fourier^{-1}\psi\in\Schwartz\left(\R^{\dimension}\right)\subset L^{1}\left(\R^{\dimension}\right)$,
since $\Fourier^{-1}\psi\in\Schwartz\left(\R^{\dimension}\right)$
and since all partial derivatives of $\widehat{f}$ are polynomially
bounded. By the Fourier inversion theorem, this implies $f\ast\tilde{\psi}=h=\Fourier^{-1}\widehat{h}=\Fourier^{-1}\left[\widehat{f}\cdot\Fourier^{-1}\psi\right]$,
so that we can continue the calculation from above as follows:
\begin{align*}
\left\langle \Fourier^{-1}\left[\widehat{f}\cdot\widehat{\varphi}\right],\psi\right\rangle _{\Schwartz',\Schwartz} & =\left\langle \varphi,\,\left(f\ast\tilde{\psi}\right)\left(-\mybullet\right)\right\rangle _{\Schwartz',\Schwartz}\\
 & =\int_{\R^{\dimension}}\varphi\left(x\right)\cdot\int_{\R^{\dimension}}f\left(-x-y\right)\cdot\tilde{\psi}\left(y\right)\d y\d x\\
\left({\scriptstyle z=-y}\right) & =\int_{\R^{\dimension}}\varphi\left(x\right)\cdot\int_{\R^{\dimension}}f\left(z-x\right)\cdot\psi\left(z\right)\d z\d x\\
\left({\scriptstyle \text{Fubini}}\right) & =\int_{\R^{\dimension}}\int_{\R^{\dimension}}\varphi\left(x\right)\cdot f\left(z-x\right)\d x\cdot\psi\left(z\right)\d z\\
 & =\left\langle f\ast\varphi,\,\psi\right\rangle _{\Schwartz',\Schwartz},
\end{align*}
which proves the claim.

All that remains is to justify the application of Fubini's theorem.
To this end, we can assume $\varphi\in L_{v}^{1}\left(\R^{\dimension}\right)$
or $\varphi\in L_{v}^{\infty}\left(\R^{\dimension}\right)$, since
then the general case follows by linearity. But for $\varphi\in L_{v}^{\infty}\left(\R^{\dimension}\right)$,
we have because of
\begin{align*}
v\left(0\right) & =v\left(x+\left(-x\right)\right)\\
 & \leq v\left(x\right)\cdot v_{0}\left(-x\right)\\
 & =v\left(x\right)\cdot v_{0}\left(z-x+\left(-z\right)\right)\\
 & \leq v\left(x\right)\cdot v_{0}\left(z-x\right)\cdot v_{0}\left(-z\right)\\
 & \leq\Omega_{1}\cdot v\left(x\right)\cdot v_{0}\left(z-x\right)\cdot\left(1+\left|z\right|\right)^{K}
\end{align*}
that
\begin{align*}
\int_{\R^{\dimension}}\int_{\R^{\dimension}}\left|\varphi\left(x\right)\cdot f\left(z-x\right)\cdot\psi\left(z\right)\right|\d x\d z & \leq\frac{\Omega_{1}}{v\left(0\right)}\cdot\int_{\R^{\dimension}}\int_{\R^{\dimension}}\left|\left(v\cdot\varphi\right)\left(x\right)\right|\cdot\left|\left(v_{0}\cdot f\right)\left(z-x\right)\right|\cdot\left(1+\left|z\right|\right)^{K}\left|\psi\left(z\right)\right|\d x\d z\\
 & \leq\frac{\Omega_{1}}{v\left(0\right)}\cdot\left\Vert \varphi\right\Vert _{L_{v}^{\infty}}\cdot\int_{\R^{\dimension}}\left(1+\left|z\right|\right)^{K}\left|\psi\left(z\right)\right|\int_{\R^{\dimension}}\left|\left(v_{0}\cdot f\right)\left(z-x\right)\right|\d x\d z\\
\left({\scriptstyle y=z-x}\right) & =\frac{\Omega_{1}}{v\left(0\right)}\cdot\left\Vert \varphi\right\Vert _{L_{v}^{\infty}}\left\Vert f\right\Vert _{L_{v_{0}}^{1}}\int_{\R^{\dimension}}\left(1+\left|z\right|\right)^{K}\left|\psi\left(z\right)\right|\d z<\infty.
\end{align*}
Furthermore, in case of $\varphi\in L_{v}^{1}\left(\R^{\dimension}\right)$,
we get with a similar estimate that
\begin{align*}
\int_{\R^{\dimension}}\int_{\R^{\dimension}}\left|\varphi\left(x\right)\cdot f\left(z-x\right)\cdot\psi\left(z\right)\right|\d x\d z & \leq\frac{\Omega_{1}}{v\left(0\right)}\cdot\int_{\R^{\dimension}}\int_{\R^{\dimension}}\left|\left(v\cdot\varphi\right)\left(x\right)\right|\cdot\left|\left(v_{0}\cdot f\right)\left(z-x\right)\right|\cdot\left(1+\left|z\right|\right)^{K}\left|\psi\left(z\right)\right|\d x\d z\\
 & \leq\frac{\Omega_{1}}{v\left(0\right)}\cdot\left[\sup_{z\in\R^{\dimension}}\left(1+\left|z\right|\right)^{K}\left|\psi\left(z\right)\right|\right]\cdot\int_{\R^{\dimension}}\left|\left(v\cdot\varphi\right)\left(x\right)\right|\int_{\R^{\dimension}}\left|\left(v_{0}\cdot f\right)\left(z-x\right)\right|\d z\d x\\
 & =\frac{\Omega_{1}}{v\left(0\right)}\cdot\left[\sup_{z\in\R^{\dimension}}\left(1+\left|z\right|\right)^{K}\left|\psi\left(z\right)\right|\right]\cdot\left\Vert \varphi\right\Vert _{L_{v}^{1}}\cdot\left\Vert f\right\Vert _{L_{v_{0}}^{1}}<\infty.
\end{align*}

For the proof of equation (\ref{eq:ViAreTemperedDistributions}),
note for $p\in\left[1,\infty\right]$ that $V_{i}=L_{v}^{p}\left(\R^{\dimension}\right)\hookrightarrow L_{v}^{1}\left(\R^{\dimension}\right)+L_{v}^{\infty}\left(\R^{\dimension}\right)$,
because of the well-known (cf.\@ \cite[Proposition (6.9)]{FollandRA})
embedding $L^{p}\left(\R^{\dimension}\right)\hookrightarrow L^{1}\left(\R^{\dimension}\right)+L^{\infty}\left(\R^{\dimension}\right)$.
But in case of $p\in\left(0,1\right)$, Theorem \ref{thm:WienerAmalgamConvolution}
and the ensuing remark yield $V_{i}=W_{T_{i}^{-T}\left[-1,1\right]^{\dimension}}\left(L_{v}^{p}\right)\hookrightarrow L_{v}^{\infty}\left(\R^{\dimension}\right)\hookrightarrow L_{v}^{1}\left(\R^{\dimension}\right)+L_{v}^{\infty}\left(\R^{\dimension}\right)$,
as desired.
\end{proof}
One of our aims in this section is to show under the conditions of
Assumption \ref{assu:MainAssumptions} (and certain additional assumptions,
cf.\@ Assumption \ref{assu:GammaCoversOrbit}) on $\CalQ,\Gamma=\left(\gamma_{i}\right)_{i\in I}$
and $p,q,v,w$ that we have
\begin{equation}
\left\Vert f\right\Vert _{\DecompSp{\CalQ}p{\ell_{w}^{q}}v}\asymp\left\Vert \left(\left\Vert \gamma^{\left(i\right)}\ast f\right\Vert _{V_{i}}\right)_{i\in I}\right\Vert _{\ell_{w}^{q}}=\left\Vert \left(\gamma^{\left(i\right)}\ast f\right)_{i\in I}\right\Vert _{V}\label{eq:DesiredCharacterization}
\end{equation}
for all $f\in\DecompSp{\CalQ}p{\ell_{w}^{q}}v$. Note though, that
it is not a priori clear how the convolution $\gamma^{\left(i\right)}\ast f$
can be interpreted, since we have $f\in\DecompSp{\CalQ}p{\ell_{w}^{q}}v\leq\Fourier^{-1}\left(\DistributionSpace{\CalO}\right)\nsubseteq\Schwartz'\left(\R^{\dimension}\right)$.
The purpose of the following result is to clarify how $\gamma^{\left(i\right)}\ast f$
can be interpreted and to establish the estimate ``$\gtrsim$''
in equation (\ref{eq:DesiredCharacterization}). We remark that the
theorem uses the notion of \textbf{normal convergence} of a series.
In our context, we say that a series $\sum_{i\in I}g_{i}$ converges
normally in $V_{j}$ if
\[
\begin{cases}
\sum_{i\in I}\left\Vert g_{i}\right\Vert _{V_{j}}<\infty, & \text{if }p\in\left[1,\infty\right],\\
\sum_{i\in I}\left\Vert g_{i}\right\Vert _{V_{j}}^{p}<\infty, & \text{if }p\in\left(0,1\right).
\end{cases}
\]

\begin{thm}
\label{thm:ConvolvingDecompositionSpaceWithGammaJ}If Assumption \ref{assu:MainAssumptions}
is fulfilled, the following hold:

For every $f\in\DecompSp{\CalQ}p{\ell_{w}^{q}}v$ and $j\in I$, the
distribution $\widehat{\gamma^{\left(j\right)}}\cdot\widehat{f}\in\DistributionSpace{\CalO}$
extends to a tempered distribution $f_{j}\in\Fourier\left(V_{j}\right)\subset\Schwartz'\left(\R^{\dimension}\right)$,
given by
\[
f_{j}:\Schwartz\left(\R^{\dimension}\right)\to\Compl,\phi\mapsto\sum_{i\in I}\left\langle \widehat{\gamma^{\left(j\right)}}\cdot\widehat{f},\,\varphi_{i}\phi\right\rangle _{\DistributionSpace{\CalO},\TestFunctionSpace{\CalO}}.
\]
Furthermore, the inverse Fourier transform $\Fourier^{-1}f_{j}\in V_{j}$
is given by
\begin{equation}
\left(\Fourier^{-1}f_{j}\right)\left(x\right)=\sum_{i\in I}\left[\Fourier^{-1}\left(\varphi_{i}\widehat{\gamma^{\left(j\right)}}\widehat{f}\right)\right]\left(x\right),\label{eq:SpecialConvolutionInterpretation}
\end{equation}
where the series converges normally in $V_{j}$ and absolutely almost
everywhere.

Finally, the linear map 
\[
{\rm Ana}_{\Gamma}:\DecompSp{\CalQ}p{\ell_{w}^{q}}v\to\ell_{w}^{q}\left(\left[V_{i}\right]_{i\in I}\right),f\mapsto\left(\Fourier^{-1}f_{j}\right)_{j\in I}
\]
is well-defined and bounded, with 
\[
\vertiii{{\rm Ana}_{\Gamma}}\leq C\cdot\vertiii{\smash{\Gamma_{\CalQ}}}\cdot\vertiii{\smash{\overrightarrow{A}}}_{\ell_{w^{\min\left\{ 1,p\right\} }}^{r}\left(I\right)\to\ell_{w^{\min\left\{ 1,p\right\} }}^{r}\left(I\right)}^{\max\left\{ 1,\frac{1}{p}\right\} },
\]
where
\[
C:=\begin{cases}
1, & \text{if }p\in\left[1,\infty\right],\\
N_{\CalQ}^{\frac{1}{p}-1}\cdot\left(12288\cdot\dimension^{3/2}\cdot\left\lceil K+\frac{\dimension+1}{p}\right\rceil \right)^{\left\lceil K+\frac{\dimension+1}{p}\right\rceil +1}\cdot\left(1+R_{\CalQ}\right)^{\dimension/p}\left(12R_{\CalQ}C_{\CalQ}\right)^{\dimension\left(\frac{1}{p}-1\right)}\cdot\Omega_{0}^{K}\Omega_{1}, & \text{if }p\in\left(0,1\right)
\end{cases}
\]
and where $\Gamma_{\CalQ}:\ell_{w}^{q}\left(I\right)\to\ell_{w}^{q}\left(I\right)$
is the \textbf{$\CalQ$-clustering map}, i.e., $\Gamma_{\CalQ}\left(c_{i}\right)_{i\in I}=\left(c_{i}^{\ast}\right)_{i\in I}$,
with $c_{i}^{\ast}:=\sum_{\ell\in i^{\ast}}c_{\ell}$.
\end{thm}
\begin{rem*}
In the following, we will use the notation $\gamma^{\left(j\right)}\ast f$
instead of $\Fourier^{-1}f_{j}$, so that we have 
\[
{\rm Ana}_{\Gamma}\,f=\left(\gamma^{\left(j\right)}\ast f\right)_{j\in I}.
\]
Likewise, because of $\gamma^{\left[j\right]}=\left|\det T_{j}\right|^{-1/2}\cdot\gamma^{\left(j\right)}$,
it is natural to define
\[
\gamma^{\left[j\right]}\ast f:=\left|\det T_{j}\right|^{-1/2}\cdot\gamma^{\left(j\right)}\ast f.
\]

This new notation $\gamma^{\left(j\right)}\ast f$ (and thus also
$\gamma^{\left[j\right]}\ast f$) is consistent in the following sense:
If we have $\CalO=\R^{\dimension}$ and $\DecompSp{\CalQ}p{\ell_{w}^{q}}v\hookrightarrow\Schwartz'\left(\R^{\dimension}\right)$
(i.e., if every $f\in\DecompSp{\CalQ}p{\ell_{w}^{q}}v\subset Z'\left(\R^{\dimension}\right)=\left[\Fourier\left(\TestFunctionSpace{\R^{\dimension}}\right)\right]'$
extends to a tempered distribution $f_{\Schwartz}$), then our new
definition of the convolution $\gamma^{\left(j\right)}\ast f:=\Fourier^{-1}f_{j}$
agrees with the usual interpretation of $\gamma^{\left(j\right)}\ast f:=\Fourier^{-1}\left(\widehat{\gamma^{\left(j\right)}}\cdot\widehat{f}\right)$
for $f\in\Schwartz'\left(\R^{\dimension}\right)\supset\DecompSp{\CalQ}p{\ell_{w}^{q}}v$,
as we will see now.

First note $\widehat{f_{\Schwartz}}|_{\TestFunctionSpace{\R^{\dimension}}}=\widehat{f}$,
where $\widehat{f}=f\circ\Fourier\in\DistributionSpace{\R^{\dimension}}$.
Thus, we have for arbitrary $\phi\in\Fourier\left(\TestFunctionSpace{\R^{\dimension}}\right)$
that
\begin{align*}
\left\langle \Fourier^{-1}f_{j},\,\phi\right\rangle _{\Schwartz',\Schwartz} & =\left\langle f_{j},\,\Fourier^{-1}\phi\right\rangle _{\Schwartz',\Schwartz}\\
 & =\sum_{i\in I}\left\langle \widehat{\gamma^{\left(j\right)}}\widehat{f},\,\varphi_{i}\cdot\Fourier^{-1}\phi\right\rangle _{\DistributionSpace{\R^{\dimension}},\TestFunctionSpace{\R^{\dimension}}}\\
\left({\scriptstyle \text{since }\Fourier^{-1}\phi\in\TestFunctionSpace{\smash{\R^{\dimension}}}\text{ and }\sum_{i\in I}\varphi_{i}\equiv1\text{ with a locally finite sum}}\right) & =\left\langle \widehat{\gamma^{\left(j\right)}}\cdot\widehat{f},\,\Fourier^{-1}\phi\right\rangle _{\DistributionSpace{\R^{\dimension}},\TestFunctionSpace{\R^{\dimension}}}\\
 & =\left\langle \widehat{f},\,\widehat{\gamma^{\left(j\right)}}\cdot\Fourier^{-1}\phi\right\rangle _{\DistributionSpace{\R^{\dimension}},\TestFunctionSpace{\R^{\dimension}}}\\
\left({\scriptstyle \widehat{\gamma^{\left(j\right)}}\cdot\Fourier^{-1}\phi\in\TestFunctionSpace{\smash{\R^{\dimension}}}\subset\Schwartz\left(\smash{\R^{\dimension}}\right),\text{ since }\widehat{\gamma^{\left(j\right)}}\in C^{\infty}\left(\smash{\R^{\dimension}}\right)\text{ and }\Fourier^{-1}\phi\in\TestFunctionSpace{\smash{\R^{\dimension}}}}\right) & =\left\langle \widehat{f_{\Schwartz}},\,\widehat{\gamma^{\left(j\right)}}\cdot\Fourier^{-1}\phi\right\rangle _{\Schwartz',\Schwartz}\\
\left({\scriptstyle \widehat{\gamma^{\left(j\right)}}\cdot\widehat{f_{\Schwartz}}\in\Schwartz'\left(\smash{\R^{\dimension}}\right),\text{ since }\widehat{f_{\Schwartz}}\in\Schwartz'\left(\smash{\R^{\dimension}}\right)\text{ and all derivatives of }\widehat{\gamma^{\left(j\right)}}\text{ pol. bounded}}\right) & =\left\langle \widehat{\gamma^{\left(j\right)}}\cdot\widehat{f_{\Schwartz}},\,\Fourier^{-1}\phi\right\rangle _{\Schwartz',\Schwartz}\\
 & =\left\langle \Fourier^{-1}\left[\widehat{\gamma^{\left(j\right)}}\cdot\widehat{f_{\Schwartz}}\right],\,\phi\right\rangle _{\Schwartz',\Schwartz}=\left\langle \gamma^{\left(j\right)}\ast f_{\Schwartz},\,\phi\right\rangle _{\Schwartz',\Schwartz}.
\end{align*}
Here, the last step uses the \emph{definition} $\gamma^{\left(j\right)}\ast f_{\Schwartz}:=\Fourier^{-1}\left[\widehat{\gamma^{\left(j\right)}}\cdot\widehat{f_{\Schwartz}}\right]$
from above. This definition coincides with the usual one if $\gamma_{i}\in\Schwartz\left(\R^{\dimension}\right)$
for all $i\in I$ (so that $\gamma^{\left(j\right)}\in\Schwartz\left(\R^{\dimension}\right)$)
or (by Lemma \ref{lem:SpecialConvolutionConsistent} and since $\gamma^{\left(j\right)}\in L_{v_{0}}^{1}\left(\R^{\dimension}\right)$
as seen in Remark \ref{rem:MainAssumptionsRemark}) if $f_{\Schwartz}\in\left(L_{v}^{1}+L_{v}^{\infty}\right)\left(\R^{\dimension}\right)$,
which is satisfied in many cases.

Now, since $\Fourier\left(\TestFunctionSpace{\R^{\dimension}}\right)$
is dense in $\Schwartz\left(\R^{\dimension}\right)$ (cf.\@ \cite[Proposition 9.9]{FollandRA})
and since we have $\Fourier^{-1}f_{j}\in V_{j}\subset\Schwartz'\left(\R^{\dimension}\right)$
and $\gamma^{\left(j\right)}\ast f_{\Schwartz}=\Fourier^{-1}\left(\widehat{\gamma^{\left(j\right)}}\cdot\widehat{f_{\Schwartz}}\right)\in\Schwartz'\left(\R^{\dimension}\right)$,
we conclude $\gamma^{\left(j\right)}\ast f_{\Schwartz}=\Fourier^{-1}f_{j}$,
as claimed.
\end{rem*}
\begin{proof}[Proof of Theorem \ref{thm:ConvolvingDecompositionSpaceWithGammaJ}]
Let $f\in\DecompSp{\CalQ}p{\ell_{w}^{q}}v$ be arbitrary and let
$c_{i}:=\left\Vert \Fourier^{-1}\left(\varphi_{i}^{\ast}\cdot\widehat{f}\right)\right\Vert _{L_{v}^{p}}$
for $i\in I$. Using the (quasi)-triangle inequality for $L^{p}\left(\R^{\dimension}\right)$
and the uniform estimate $\left|i^{\ast}\right|\leq N_{\CalQ}$, we
obtain a constant $C_{1}=C_{1}\left(p,\CalQ\right)>0$ satisfying
\[
c_{i}=\left\Vert \Fourier^{-1}\!\left(\varphi_{i}^{\ast}\cdot\widehat{f}\right)\right\Vert _{L_{v}^{p}}\leq C_{1}\cdot\sum_{\ell\in i^{\ast}}\left\Vert \Fourier^{-1}\!\left(\varphi_{\ell}\cdot\widehat{f}\right)\right\Vert _{L_{v}^{p}}=C_{1}\cdot\left(\Gamma_{\CalQ}d\right)_{i}\quad\text{for}\quad d=\left(d_{i}\right)_{i\in I},\text{ with }d_{i}:=\left\Vert \Fourier^{-1}\!\left(\varphi_{i}\cdot\widehat{f}\right)\right\Vert _{L_{v}^{p}}.
\]
In fact, as shown in \cite[Exercise 1.1.5(c)]{GrafakosClassicalFourierAnalysis},
we can choose
\[
C_{1}=\begin{cases}
1, & \text{if }p\in\left[1,\infty\right],\\
N_{\CalQ}^{\frac{1}{p}-1}, & \text{if }p\in\left(0,1\right).
\end{cases}
\]
Since $d\in\ell_{w}^{q}\left(I\right)$ with $\left\Vert d\right\Vert _{\ell_{w}^{q}}=\left\Vert f\right\Vert _{\DecompSp{\CalQ}p{\ell_{w}^{q}}v}$,
we get $c\in\ell_{w}^{q}\left(I\right)$ as well, and $\left\Vert c\right\Vert _{\ell_{w}^{q}}\leq C_{1}\cdot\vertiii{\Gamma_{\CalQ}}\cdot\left\Vert f\right\Vert _{\DecompSp{\CalQ}p{\ell_{w}^{q}}v}$.

Now, we distinguish the two cases $p\in\left[1,\infty\right]$ and
$p\in\left(0,1\right)$.

\medskip{}

In case of $p\in\left[1,\infty\right]$, we have $V_{j}=L_{v}^{p}\left(\R^{\dimension}\right)$.
Here, the weighted Young inequality (equation (\ref{eq:WeightedYoungInequality}))
yields
\begin{align*}
\left\Vert \Fourier^{-1}\left(\widehat{\gamma^{\left(j\right)}}\cdot\varphi_{i}\cdot\widehat{f}\right)\right\Vert _{L_{v}^{p}} & =\left\Vert \Fourier^{-1}\left(\widehat{\gamma^{\left(j\right)}}\cdot\varphi_{i}\cdot\varphi_{i}^{\ast}\cdot\widehat{f}\right)\right\Vert _{L_{v}^{p}}\\
 & =\left\Vert \Fourier^{-1}\left(\widehat{\gamma^{\left(j\right)}}\cdot\varphi_{i}\right)\ast\Fourier^{-1}\left(\varphi_{i}^{\ast}\cdot\widehat{f}\right)\right\Vert _{L_{v}^{p}}\\
 & \leq\left\Vert \Fourier^{-1}\left(\widehat{\gamma^{\left(j\right)}}\cdot\varphi_{i}\right)\right\Vert _{L_{v_{0}}^{1}}\cdot\left\Vert \Fourier^{-1}\left(\varphi_{i}^{\ast}\cdot\widehat{f}\right)\right\Vert _{L_{v}^{p}}\\
 & =A_{j,i}\cdot c_{i},
\end{align*}
with $A_{j,i}$ as in Assumption \ref{assu:MainAssumptions}. Hence,
we get
\begin{equation}
\sum_{i\in I}\left\Vert \Fourier^{-1}\left(\widehat{\gamma^{\left(j\right)}}\cdot\varphi_{i}\cdot\widehat{f}\right)\right\Vert _{L_{v}^{p}}\leq\sum_{i\in I}\left[A_{j,i}\cdot c_{i}\right]=\left(\overrightarrow{A}\cdot c\right)_{j}<\infty,\label{eq:GammaAnalysisBanachCaseMainEstimate}
\end{equation}
since we have $c\in\ell_{w}^{q}\left(I\right)$ and since Assumption
\ref{assu:MainAssumptions} includes (for $p\in\left[1,\infty\right]$)
the assumption that $\overrightarrow{A}:\ell_{w}^{q}\left(I\right)\to\ell_{w}^{q}\left(I\right)$
is well-defined and bounded. This implies that the function
\[
F_{j}:=\sum_{i\in I}\Fourier^{-1}\left(\widehat{\gamma^{\left(j\right)}}\cdot\varphi_{i}\cdot\widehat{f}\right)\in L_{v}^{p}\left(\smash{\R^{\dimension}}\right)=V_{j}
\]
is well-defined, with normal convergence in $V_{j}$ and with absolute
convergence a.e.\@ of the defining series and such that $\left\Vert F_{j}\right\Vert _{L_{v}^{p}}\leq\left(\overrightarrow{A}\cdot c\right)_{j}$
for all $j\in I$.

\medskip{}

Next, in case of $p\in\left(0,1\right)$, define $e_{i}:=c_{i}^{p}$
for $i\in I$ and note $e=\left(e_{i}\right)_{i\in I}\in\ell_{w^{p}}^{q/p}\left(I\right)=\ell_{w^{\min\left\{ 1,p\right\} }}^{r}\left(I\right)$,
with 
\[
\left\Vert e\right\Vert _{\ell_{w^{\min\left\{ 1,p\right\} }}^{r}\left(I\right)}=\left\Vert \left(w_{i}^{p}\cdot c_{i}^{p}\right)_{i\in I}\right\Vert _{\ell^{q/p}}=\left\Vert \left(w_{i}\cdot c_{i}\right)_{i\in I}\right\Vert _{\ell^{q}}^{p}=\left\Vert c\right\Vert _{\ell_{w}^{q}\left(I\right)}^{p}.
\]
Next, we note
\begin{align*}
\supp\left(\widehat{\gamma^{\left(j\right)}}\cdot\varphi_{i}\cdot\widehat{f}\right) & \subset\supp\varphi_{i}\subset\overline{Q_{i}}\subset T_{i}\overline{B_{R_{\CalQ}}\left(0\right)}+b_{i}\\
 & \subset T_{j}\left[T_{j}^{-1}T_{i}\overline{B_{R_{\CalQ}}\left(0\right)}\right]+b_{i}\\
 & \subset T_{j}\left[\left\Vert T_{j}^{-1}T_{i}\right\Vert \overline{B_{R_{\CalQ}}}\left(0\right)\right]+b_{i}\\
 & \subset T_{j}\left[-\left\Vert T_{j}^{-1}T_{i}\right\Vert R_{\CalQ},\,\left\Vert T_{j}^{-1}T_{i}\right\Vert R_{\CalQ}\right]^{\dimension}+b_{i},
\end{align*}
so that Theorem \ref{thm:BandlimitedWienerAmalgamSelfImproving} yields
for $C_{2}:=2^{4\left(1+\frac{\dimension}{p}\right)}s_{\dimension}^{\frac{1}{p}}\left(192\cdot\dimension^{\frac{3}{2}}\cdot\left\lceil K\!+\!\frac{\dimension+1}{p}\right\rceil \right)^{\left\lceil K+\frac{\dimension+1}{p}\right\rceil +1}\cdot\Omega_{0}^{K}\Omega_{1}$
and $C_{3}:=C_{2}\cdot\left(1+R_{\CalQ}\right)^{\dimension/p}$ that
\begin{align*}
\left\Vert \Fourier^{-1}\left(\widehat{\gamma^{\left(j\right)}}\cdot\varphi_{i}\cdot\widehat{f}\right)\right\Vert _{V_{j}} & =\left\Vert \Fourier^{-1}\left(\widehat{\gamma^{\left(j\right)}}\cdot\varphi_{i}\cdot\widehat{f}\right)\right\Vert _{W_{T_{j}^{-T}\left[-1,1\right]^{\dimension}}\left(L_{v}^{p}\right)}\\
 & \leq C_{2}\left(1+\left\Vert T_{j}^{-1}T_{i}\right\Vert R_{\CalQ}\right)^{\dimension/p}\cdot\left\Vert \Fourier^{-1}\left(\widehat{\gamma^{\left(j\right)}}\cdot\varphi_{i}\cdot\widehat{f}\right)\right\Vert _{L_{v}^{p}}\\
\left({\scriptstyle 1+ab\leq\left(1+a\right)\left(1+b\right)\text{ for }a,b\geq0}\right) & \leq C_{3}\cdot\left(1+\left\Vert T_{j}^{-1}T_{i}\right\Vert \right)^{\dimension/p}\cdot\left\Vert \Fourier^{-1}\left(\widehat{\gamma^{\left(j\right)}}\cdot\varphi_{i}\right)\ast\Fourier^{-1}\left(\varphi_{i}^{\ast}\cdot\widehat{f}\right)\right\Vert _{L_{v}^{p}}\\
\left({\scriptstyle \text{Prop. }\ref{prop:BandlimitedConvolution}\text{ with }n=1}\right) & \leq C_{3}\!\cdot\!\left(12R_{\CalQ}C_{\CalQ}\right)^{\dimension\left(\frac{1}{p}-1\right)}\!\cdot\!\left(1\!+\!\left\Vert T_{j}^{-1}T_{i}\right\Vert \right)^{\dimension/p}\left|\det T_{i}\right|^{\frac{1}{p}-1}\left\Vert \Fourier^{-1}\left(\varphi_{i}\,\widehat{\gamma^{\left(j\right)}}\right)\right\Vert _{L_{v_{0}}^{p}}\left\Vert \Fourier^{-1}\left(\varphi_{i}^{\ast}\,\widehat{f}\right)\right\Vert _{L_{v}^{p}}\\
 & \leq C_{3}\cdot\left(12R_{\CalQ}C_{\CalQ}\right)^{\dimension\left(\frac{1}{p}-1\right)}\cdot A_{j,i}^{1/p}\cdot c_{i}\\
 & =:C_{4}\cdot A_{j,i}^{1/p}\cdot c_{i}.
\end{align*}
Here, Proposition \ref{prop:BandlimitedConvolution} is applicable,
since $\varphi_{i}\in\TestFunctionSpace{\R^{\dimension}}$ and $\widehat{\gamma^{\left(j\right)}}\in C^{\infty}\left(\R^{\dimension}\right)$,
so that $\varphi_{i}\cdot\widehat{\gamma^{\left(j\right)}}\in\TestFunctionSpace{\R^{\dimension}}$
and since clearly $\supp\left[\varphi_{i}\widehat{\gamma^{\left(j\right)}}\right]\subset\overline{Q_{i}^{\ast}}$
and $\supp\left[\varphi_{i}^{\ast}\widehat{f}\right]\subset\overline{Q_{i}^{\ast}}$.

Consequently, we arrive at
\begin{equation}
\sum_{i\in I}\left\Vert \Fourier^{-1}\left(\widehat{\gamma^{\left(j\right)}}\cdot\varphi_{i}\cdot\widehat{f}\right)\right\Vert _{V_{j}}^{p}\leq C_{4}^{p}\cdot\sum_{i\in I}\left[A_{j,i}\cdot c_{i}^{p}\right]=C_{4}^{p}\cdot\left(\overrightarrow{A}\cdot e\right)_{j}<\infty,\label{eq:GammaAnalysisQuasiBanachCaseMainEstimate}
\end{equation}
since $\overrightarrow{A}:\ell_{w^{\min\left\{ 1,p\right\} }}^{r}\left(I\right)\to\ell_{w^{\min\left\{ 1,p\right\} }}^{r}\left(I\right)$
is well-defined and bounded and $e\in\ell_{w^{\min\left\{ 1,p\right\} }}^{r}\left(I\right)$.

Finally, we use the $p$-triangle inequality for $L^{p}\left(\R^{\dimension}\right)$
(yielding the $p$-triangle inequality for $V_{j}=W_{T_{j}^{-T}\left[-1,1\right]^{\dimension}}\left(L_{v}^{p}\right)$)
to conclude that $F_{j}:=\sum_{i\in I}\Fourier^{-1}\left(\widehat{\gamma^{\left(j\right)}}\cdot\varphi_{i}\cdot\widehat{f}\right)\in V_{j}$
is well-defined, with normal convergence in $V_{j}$ and a.e.\@ absolute
convergence of the defining series and with $\left\Vert F_{j}\right\Vert _{V_{j}}\leq C_{4}\cdot\left(\smash{\overrightarrow{A}}\cdot e\right)_{j}^{1/p}$.

\medskip{}

Our next goal is to show that the previous results imply that $f_{j}\in\Schwartz'\left(\R^{\dimension}\right)$
yields a well-defined tempered distribution. To this end, recall from
Lemma \ref{lem:SpecialConvolutionConsistent} that $V_{j}\hookrightarrow\Schwartz'\left(\R^{\dimension}\right)$
for all $p\in\left(0,\infty\right]$. Consequently, we get $F_{j}\in\Schwartz'\left(\R^{\dimension}\right)$
and $F_{j}=\sum_{i\in I}\Fourier^{-1}\left(\widehat{\gamma^{\left(j\right)}}\cdot\varphi_{i}\cdot\widehat{f}\right)$
with unconditional convergence in $V_{j}\hookrightarrow\Schwartz'\left(\R^{\dimension}\right)$,
which implies for $\phi\in\Schwartz\left(\R^{\dimension}\right)$
that
\begin{align*}
\left\langle \Fourier F_{j},\,\Fourier^{-1}\phi\right\rangle _{\Schwartz',\Schwartz}=\left\langle F_{j},\,\phi\right\rangle _{\Schwartz',\Schwartz} & =\sum_{i\in I}\left\langle \Fourier^{-1}\left(\widehat{\gamma^{\left(j\right)}}\cdot\varphi_{i}\cdot\widehat{f}\right),\,\phi\right\rangle _{\Schwartz',\Schwartz},\\
 & =\sum_{i\in I}\left\langle \widehat{\gamma^{\left(j\right)}}\cdot\widehat{f},\,\varphi_{i}\cdot\Fourier^{-1}\phi\right\rangle _{\DistributionSpace{\CalO},\TestFunctionSpace{\CalO}}\\
 & =\left\langle f_{j},\,\Fourier^{-1}\phi\right\rangle _{\Schwartz',\Schwartz},
\end{align*}
where the right-hand side is well-defined (with absolute convergence
of the series), since the left-hand side is. This shows that $f_{j}=\Fourier F_{j}\in\Fourier V_{j}\subset\Schwartz'\left(\R^{\dimension}\right)$
is a well-defined tempered distribution, as claimed. Finally, we have
$\Fourier^{-1}f_{j}=F_{j}=\sum_{i\in I}\Fourier^{-1}\left(\widehat{\gamma^{\left(j\right)}}\cdot\varphi_{i}\cdot\widehat{f}\right)$,
where the series converges normally in $V_{j}$ and absolutely a.e.,
as claimed.

\medskip{}

It remains to verify boundedness of ${\rm Ana}_{\Gamma}$. But for
$p\in\left[1,\infty\right]$, we have by solidity of $\ell_{w}^{q}\left(I\right)$
and by the triangle inequality for $L^{p}\left(\R^{\dimension}\right)$,
and since $C_{1}=1$ for $p\in\left[1,\infty\right]$, that
\begin{align*}
\left\Vert \left(\left\Vert \gamma^{\left(j\right)}\ast f\right\Vert _{V_{j}}\right)_{j\in I}\right\Vert _{\ell_{w}^{q}}=\left\Vert \left(\left\Vert F_{j}\right\Vert _{L_{v}^{p}}\right)_{j\in I}\right\Vert _{\ell_{w}^{q}} & \leq\left\Vert \left(\sum_{i\in I}\left\Vert \Fourier^{-1}\left(\widehat{\gamma^{\left(j\right)}}\cdot\varphi_{i}\cdot\widehat{f}\right)\right\Vert _{L_{v}^{p}}\right)_{j\in I}\right\Vert _{\ell_{w}^{q}}\\
\left({\scriptstyle \text{eq. }\eqref{eq:GammaAnalysisBanachCaseMainEstimate}}\right) & \leq\left\Vert \left[\left(\overrightarrow{A}\cdot c\right)_{j}\right]_{j\in I}\right\Vert _{\ell_{w}^{q}}\\
 & \leq\vertiii{\smash{\overrightarrow{A}}}\cdot\left\Vert c\right\Vert _{\ell_{w}^{q}}\\
 & \leq\vertiii{\smash{\Gamma_{\CalQ}}}\cdot\vertiii{\smash{\overrightarrow{A}}}\cdot\left\Vert f\right\Vert _{\DecompSp{\CalQ}p{\ell_{w}^{q}}v}<\infty,
\end{align*}
as desired.

Finally, in case of $p\in\left(0,1\right)$, the $p$-triangle inequality
for $W_{T_{j}^{-T}\left[-1,1\right]^{\dimension}}\left(L_{v}^{p}\right)$
yields
\begin{align*}
\left\Vert \gamma^{\left(j\right)}\ast f\right\Vert _{V_{j}}=\left\Vert F_{j}\right\Vert _{W_{T_{j}^{-T}\left[-1,1\right]^{\dimension}}\left(L_{v}^{p}\right)} & \leq\left[\sum_{i\in I}\left\Vert \Fourier^{-1}\left(\widehat{\gamma^{\left(j\right)}}\cdot\varphi_{i}\cdot\widehat{f}\right)\right\Vert _{W_{T_{j}^{-T}\left[-1,1\right]^{\dimension}}\left(L_{v}^{p}\right)}^{p}\right]^{1/p}\\
\left({\scriptstyle \text{eq. }\eqref{eq:GammaAnalysisQuasiBanachCaseMainEstimate}}\right) & \leq C_{4}\cdot\left(\overrightarrow{A}\cdot e\right)_{j}^{1/p}.
\end{align*}
By solidity of $\ell_{w}^{q}\left(I\right)$, this implies
\begin{align*}
\left\Vert \left(\left\Vert \gamma^{\left(j\right)}\ast f\right\Vert _{V_{j}}\right)_{j\in I}\right\Vert _{\ell_{w}^{q}} & \leq C_{4}\cdot\left\Vert \left(\overrightarrow{A}\cdot e\right)^{1/p}\right\Vert _{\ell_{w}^{q}}\\
 & =C_{4}\cdot\left\Vert \left(w^{p}\cdot\left[\overrightarrow{A}\cdot e\right]\right)^{1/p}\right\Vert _{\ell^{q}}\\
 & =C_{4}\cdot\left\Vert w^{\min\left\{ 1,p\right\} }\cdot\left[\overrightarrow{A}\cdot e\right]\right\Vert _{\ell^{q/p}}^{1/p}\\
 & =C_{4}\cdot\left\Vert \overrightarrow{A}\cdot e\right\Vert _{\ell_{w^{\min\left\{ 1,p\right\} }}^{r}}^{1/p}\\
 & \leq C_{4}\cdot\vertiii{\smash{\overrightarrow{A}}}^{1/p}\cdot\left\Vert e\right\Vert _{\ell_{w^{\min\left\{ 1,p\right\} }}^{r}}^{1/p}\\
 & =C_{4}\cdot\vertiii{\smash{\overrightarrow{A}}}^{1/p}\cdot\left\Vert c\right\Vert _{\ell_{w}^{q}}\\
 & \leq C_{1}C_{4}\cdot\vertiii{\smash{\Gamma_{\CalQ}}}\cdot\vertiii{\smash{\overrightarrow{A}}}^{1/p}\cdot\left\Vert f\right\Vert _{\DecompSp{\CalQ}p{\ell_{w}^{q}}v}<\infty,
\end{align*}
which completes the proof.
\end{proof}
Next, we establish the estimate ``$\lesssim$'' in equation (\ref{eq:DesiredCharacterization}),
under suitable assumptions on $\left(\gamma_{i}\right)_{i\in I}$.
Notice that up to now we have not excluded the case $\gamma_{i}\equiv0$
for all $i\in I$. But if equation (\ref{eq:DesiredCharacterization})
was true, we would need at least that the family of frequency supports
$\supp\widehat{\gamma^{\left(i\right)}}$, with $i\in I$, covers
all of $\CalO$. To ensure this, we introduce the following additional
assumption:
\begin{assumption}
\label{assu:GammaCoversOrbit}We assume that for each $i\in I$ there
is some function $\theta_{i}^{\natural}\in\TestFunctionSpace{\R^{\dimension}}$
such that the family $\theta=\left(\smash{\theta_{i}^{\natural}}\right)_{i\in I}$
satisfies the following properties:

\begin{enumerate}
\item We have $\theta_{i}^{\natural}\cdot\widehat{\gamma_{i}}\equiv1$ on
$Q_{i}'$ (and thus on $\overline{Q_{i}'}$) for all $i\in I$.
\item For each $p\in\left(0,\infty\right]$, the constant
\[
\Omega_{2}^{\left(p,K\right)}:=\Omega_{2}^{\left(p,K\right)}\left(\theta\right):=\begin{cases}
\sup_{i\in I}\left\Vert \Fourier^{-1}\theta_{i}^{\natural}\right\Vert _{W_{\left[-1,1\right]^{\dimension}}\left(L_{\left(1+\left|\mybullet\right|\right)^{K}}^{p}\right)}, & \text{if }p\in\left(0,1\right),\\
\sup_{i\in I}\left\Vert \Fourier^{-1}\theta_{i}^{\natural}\right\Vert _{L_{\left(1+\left|\mybullet\right|\right)^{K}}^{1}}, & \text{if }p\in\left[1,\infty\right]
\end{cases}
\]
is finite.
\end{enumerate}
We fix such a family $\theta=\left(\smash{\theta_{i}^{\natural}}\right)_{i\in I}$
and the constant $\Omega_{2}^{\left(p,K\right)}$ for the remainder
of the paper. Finally, we recall $S_{i}\xi=T_{i}\xi+b_{i}$ and define
\[
\theta_{i}:=\theta_{i}^{\natural}\circ S_{i}^{-1}\in\TestFunctionSpace{\smash{\R^{\dimension}}}\qquad\forall i\in I.\qedhere
\]
\end{assumption}
At least in the case where the set of prototypes $\left\{ \gamma_{i}\with i\in I\right\} $
is finite, the preceding assumption can be heavily simplified, as
we show now:
\begin{lem}
\label{lem:GammaCoversOrbitAssumptionSimplified}Assume that there
are $N$ functions $\gamma_{1}^{\left(0\right)},\dots,\gamma_{N}^{\left(0\right)}$
such that for each $i\in I$ we have $\gamma_{i}=\gamma_{n_{i}}^{\left(0\right)}$
for a suitable $n_{i}\in\underline{N}$. For $n\in\underline{N}$
let 
\[
Q^{\left(n\right)}:=\bigcup\left\{ Q_{i}'\with i\in I\text{ and }n_{i}=n\right\} .
\]
If there is some $c>0$ satisfying $\left|\left(\Fourier\smash{\gamma_{n}^{\left(0\right)}}\right)\left(\xi\right)\right|\geq c$
for all $\xi\in Q^{\left(n\right)}$, then the family $\left(\gamma_{i}\right)_{i\in I}$
satisfies Assumption \ref{assu:GammaCoversOrbit}.

In fact, for arbitrary $p_{0}\in\left(0,1\right]$ and $K^{\left(0\right)}\geq0$,
there is a constant $\Omega_{3}=\Omega_{3}\left(\CalQ,\gamma_{1}^{\left(0\right)},\dots,\gamma_{N}^{\left(0\right)},p_{0},K^{\left(0\right)},\dimension\right)>0$
satisfying
\[
\Omega_{2}^{\left(p,K\right)}\leq\Omega_{3}\qquad\forall p\geq p_{0}\text{ and }K\leq K^{\left(0\right)}.\qedhere
\]
\end{lem}
\begin{rem*}
If $\gamma_{i}=\gamma$ for all $i\in I$, then the above assumptions
reduce to $\left|\widehat{\gamma}\left(\xi\right)\right|\geq c>0$
for all $\xi\in Q:=\bigcup_{i\in I}Q_{i}'$.
\end{rem*}
\begin{proof}
Recall from Assumption \ref{assu:MainAssumptions} that we always
have $\widehat{\gamma_{i}}\in C^{\infty}\left(\R^{\dimension}\right)$.
Now, by possibly dropping some elements of the family $\gamma_{1}^{\left(0\right)},\dots,\gamma_{N}^{\left(0\right)}$,
we can assume that for each $n\in\underline{N}$, there is some $i\in I$
satisfying $n_{i}=n$ and thus $\gamma_{n}^{\left(0\right)}=\gamma_{i}$.
In particular, this implies $\widehat{\gamma_{n}^{\left(0\right)}}\in C^{\infty}\left(\R^{\dimension}\right)$
for all $n\in\underline{N}$.

By continuity of $\Fourier\gamma_{n}^{\left(0\right)}$, we get $\left|\left(\Fourier\smash{\gamma_{n}^{\left(0\right)}}\right)\left(\xi\right)\right|\geq c$
for all $\xi\in\overline{Q^{\left(n\right)}}$. Furthermore, recall
from Subsection \ref{subsec:DecompSpaceDefinitionStandingAssumptions}
that we have $Q_{i}'\subset\overline{B_{R_{\CalQ}}}\left(0\right)$
for all $i\in I$, so that each of the sets $Q^{\left(n\right)}$
is bounded. Hence, $\overline{Q^{\left(n\right)}}$ is compact. Again
by continuity of $\Fourier\gamma_{n}^{\left(0\right)}$, each of the
sets 
\[
U_{n}:=\left\{ \xi\in\R^{\dimension}\with\left|\left(\Fourier\smash{\gamma_{n}^{\left(0\right)}}\right)\left(\xi\right)\right|>\frac{c}{2}\right\} 
\]
is open with $\overline{Q^{\left(n\right)}}\subset U_{n}$. Thus,
the $C^{\infty}$-Urysohn-Lemma (cf.\@ \cite[Lemma 8.18]{FollandRA})
yields some $\eta_{n}\in\TestFunctionSpace{U_{n}}$ with $\eta_{n}|_{Q^{\left(n\right)}}\equiv1$.

Now, note that $\eta_{n}/\widehat{\gamma_{n}^{\left(0\right)}}\in\TestFunctionSpace{U_{n}}$
is well-defined, since $\widehat{\gamma_{n}^{\left(0\right)}}\neq0$
on $U_{n}$. Thus, the function
\[
\theta^{\left(n\right)}:\R^{\dimension}\to\Compl,\xi\mapsto\begin{cases}
\frac{\eta_{n}\left(\xi\right)}{\widehat{\gamma_{n}^{\left(0\right)}}\left(\xi\right)}, & \text{if }\xi\in U_{n},\\
0, & \text{if }\xi\notin U_{n}
\end{cases}
\]
is a smooth function $\theta^{\left(n\right)}\in\TestFunctionSpace{\R^{\dimension}}$
with $\supp\theta^{\left(n\right)}\subset U_{n}$ and with $\theta^{\left(n\right)}\cdot\widehat{\gamma_{n}^{\left(0\right)}}=\eta_{n}\equiv1$
on $Q^{\left(n\right)}$.

Now, define $\theta_{i}^{\natural}:=\theta^{\left(n_{i}\right)}\in\TestFunctionSpace{\R^{\dimension}}$
for $i\in I$. Then, for each $i\in I$, we have $\theta_{i}^{\natural}\cdot\widehat{\gamma_{i}}=\theta^{\left(n_{i}\right)}\cdot\widehat{\gamma_{n_{i}}^{\left(0\right)}}\equiv1$
on $Q^{\left(n_{i}\right)}\supset Q_{i}'$, cf.\@ the definition
of $Q^{\left(n\right)}$.

Finally, Lemma \ref{lem:SchwartzFunctionsAreWiener} (with $N=K^{\left(0\right)}+\frac{\dimension}{p_{0}}+1$)
yields for $p\geq p_{0}$ and $K\leq K^{\left(0\right)}$ the estimate
\begin{align*}
\left\Vert \Fourier^{-1}\theta_{i}^{\natural}\right\Vert _{L_{\left(1+\left|\cdot\right|\right)^{K}}^{p}}\leq\left\Vert \Fourier^{-1}\theta_{i}^{\natural}\right\Vert _{W_{\left[-1,1\right]^{\dimension}}\left(L_{\left(1+\left|\cdot\right|\right)^{K}}^{p}\right)} & \leq\left(1+2\sqrt{\dimension}\right)^{N}\cdot\left(\frac{1}{p}\frac{s_{\dimension}}{N-K-\frac{\dimension}{p}}\right)^{1/p}\cdot\left\Vert \Fourier^{-1}\theta_{i}^{\natural}\right\Vert _{N}\\
 & \leq\left(1+2\sqrt{\dimension}\right)^{N}\cdot\left(1+\frac{s_{\dimension}}{p_{0}}\right)^{1/p_{0}}\cdot\max_{n\in\underline{N}}\left\Vert \Fourier^{-1}\theta^{\left(n\right)}\right\Vert _{N}=:\Omega_{3}.
\end{align*}
Since $N$ only depends on $K^{\left(0\right)},\dimension,p_{0}$
and since $\theta^{\left(1\right)},\dots,\theta^{\left(N\right)}$
only depend on $\CalQ$ and on $\gamma_{1}^{\left(0\right)},\dots,\gamma_{N}^{\left(0\right)}$,
$\Omega_{3}$ is as claimed in the lemma. Note that each of the norms
$\left\Vert \Fourier^{-1}\theta^{\left(n\right)}\right\Vert _{N}$
is finite, since $\theta^{\left(n\right)}\in\TestFunctionSpace{\R^{\dimension}}$,
from which we get $\Fourier^{-1}\theta^{\left(n\right)}\in\Schwartz\left(\R^{\dimension}\right)$.
\end{proof}
Now, instead of just establishing equation (\ref{eq:DesiredCharacterization}),
we will actually show that the ``coefficient map'' ${\rm Ana}_{\Gamma}$
from Theorem \ref{thm:ConvolvingDecompositionSpaceWithGammaJ} yields
a semi-discrete Banach frame for $\DecompSp{\CalQ}p{\ell_{w}^{q}}v$.
By this we mean that there exists a bounded linear ``reconstruction''
map $R:V\to\DecompSp{\CalQ}p{\ell_{w}^{q}}v$ satisfying $R\circ{\rm Ana}_{\Gamma}=\identity_{\DecompSp{\CalQ}p{\ell_{w}^{q}}v}$.
For the construction of $R$, the following result will turn out to
be helpful:
\begin{lem}
\label{lem:LocalInverseConvolution}Assume that $\Gamma=\left(\gamma_{i}\right)_{i\in I}$
satisfies Assumption \ref{assu:GammaCoversOrbit} and let $\left(\theta_{i}\right)_{i\in I}$
be defined as in that assumption.

Then $\widehat{\gamma^{\left(i\right)}}\cdot\theta_{i}\equiv1$ on
$\overline{Q_{i}}$ for each $i\in I$ and each of the maps
\[
I_{i}:V_{i}\to V_{i},f\mapsto\left(\Fourier^{-1}\theta_{i}\right)\ast f
\]
is well-defined and bounded, with $\sup_{i\in I}\vertiii{I_{i}}\leq C<\infty$,
where
\[
C:=\begin{cases}
\Omega_{0}^{4K}\Omega_{1}^{4}\Omega_{2}^{\left(p,K\right)}\cdot\dimension^{-\frac{\dimension}{2p}}\cdot\left(972\cdot\dimension^{5/2}\right)^{K+\frac{\dimension}{p}}, & \text{if }p\in\left(0,1\right),\\
\Omega_{0}^{K}\Omega_{1}\Omega_{2}^{\left(p,K\right)}, & \text{if }p\in\left[1,\infty\right].
\end{cases}
\]

Hence, the map
\[
m_{\theta}:=\bigotimes_{i\in I}I_{i}:V\to V,\left(f_{i}\right)_{i\in I}\mapsto\left(\left(\Fourier^{-1}\theta_{i}\right)\ast f_{i}\right)_{i\in I}
\]
is well-defined and bounded as well, with $\vertiii{m_{\theta}}\leq C$.
\end{lem}
\begin{proof}
First, observe 
\[
\widehat{\gamma^{\left(i\right)}}\cdot\theta_{i}=\left(\widehat{\gamma_{i}}\circ S_{i}^{-1}\right)\cdot\left(\theta_{i}^{\natural}\circ S_{i}^{-1}\right)=\underbrace{\left(\widehat{\gamma_{i}}\cdot\theta_{i}^{\natural}\right)}_{\equiv1\text{ on }\overline{Q_{i}'}}\circ S_{i}^{-1}\equiv1\text{ on }S_{i}\overline{Q_{i}'}=\overline{Q_{i}},
\]
so that it remains to show that each of the maps $I_{i}$ is well-defined
and bounded, with the claimed estimate for the operator norm.

In case of $p\in\left[1,\infty\right]$, this is a consequence of
equation (\ref{eq:WeightedYoungInequality}), once we show that $\left\Vert \Fourier^{-1}\theta_{i}\right\Vert _{L_{v_{0}}^{1}}$
is uniformly bounded. But we simply have
\begin{equation}
\theta_{i}=\theta_{i}^{\natural}\circ S_{i}^{-1}=L_{b_{i}}\left(\theta_{i}^{\natural}\circ T_{i}^{-1}\right)\qquad\text{ and hence }\qquad\Fourier^{-1}\theta_{i}=\left|\det T_{i}\right|\cdot M_{b_{i}}\left[\left(\Fourier^{-1}\theta_{i}^{\natural}\right)\circ T_{i}^{T}\right],\label{eq:LocalInverseClosedForm}
\end{equation}
which implies
\begin{align*}
\left\Vert \Fourier^{-1}\theta_{i}\right\Vert _{L_{v_{0}}^{1}} & =\left|\det T_{i}\right|\cdot\left\Vert v_{0}\cdot\left[\left(\Fourier^{-1}\theta_{i}^{\natural}\right)\circ T_{i}^{T}\right]\right\Vert _{L^{1}}\\
 & =\left\Vert \left(v_{0}\circ T_{i}^{-T}\right)\cdot\left(\Fourier^{-1}\theta_{i}^{\natural}\right)\right\Vert _{L^{1}}\\
\left({\scriptstyle \text{assumption on }v_{0}}\right) & \leq\Omega_{1}\cdot\left\Vert x\mapsto\left(1+\left|T_{i}^{-T}x\right|\right)^{K}\cdot\left(\Fourier^{-1}\theta_{i}^{\natural}\right)\left(x\right)\right\Vert _{L^{1}}\\
\left({\scriptstyle \text{eq. }\eqref{eq:WeightLinearTransformationsConnection}}\right) & \leq\Omega_{0}^{K}\Omega_{1}\cdot\left\Vert \left(1+\left|\mybullet\right|\right)^{K}\cdot\Fourier^{-1}\theta_{i}^{\natural}\right\Vert _{L^{1}}\leq\Omega_{0}^{K}\Omega_{1}\cdot\Omega_{2}^{\left(p,K\right)}.
\end{align*}

Finally, for $p\in\left(0,1\right)$, we get from Corollary \ref{cor:WienerAmalgamConvolutionSimplified}
for $C_{1}:=\Omega_{0}^{3K}\Omega_{1}^{3}\cdot\dimension^{-\frac{\dimension}{2p}}\cdot\left(972\cdot\dimension^{5/2}\right)^{K+\frac{\dimension}{p}}$
that
\begin{align*}
\left\Vert \left(\Fourier^{-1}\theta_{i}\right)\ast f\right\Vert _{W_{T_{i}^{-T}\left[-1,1\right]^{\dimension}}\left(L_{v}^{p}\right)} & \leq C_{1}\cdot\left|\det T_{i}\right|^{\frac{1}{p}-1}\cdot\left\Vert \Fourier^{-1}\theta_{i}\right\Vert _{W_{T_{i}^{-T}\left[-1,1\right]^{\dimension}}\left(L_{v_{0}}^{p}\right)}\cdot\left\Vert f\right\Vert _{W_{T_{i}^{-T}\left[-1,1\right]^{\dimension}}\left(L_{v}^{p}\right)}\\
\left({\scriptstyle \text{eq. }\eqref{eq:LocalInverseClosedForm}\text{ and }\left\Vert M_{b}f\right\Vert _{W_{Q}\left(L_{v_{0}}^{p}\right)}=\left\Vert f\right\Vert _{W_{Q}\left(L_{v_{0}}^{p}\right)}}\right) & =C_{1}\cdot\left|\det T_{i}\right|^{\frac{1}{p}-1}\cdot\left|\det T_{i}\right|\cdot\left\Vert \left(\Fourier^{-1}\theta_{i}^{\natural}\right)\circ T_{i}^{T}\right\Vert _{W_{T_{i}^{-T}\left[-1,1\right]^{\dimension}}\left(L_{v_{0}}^{p}\right)}\cdot\left\Vert f\right\Vert _{V_{i}}\\
\left({\scriptstyle \text{Lemma }\ref{lem:WienerTransformationFormula}}\right) & =C_{1}\cdot\left|\det T_{i}\right|^{\frac{1}{p}}\cdot\left\Vert \left(M_{\left[-1,1\right]^{\dimension}}\left[\Fourier^{-1}\theta_{i}^{\natural}\right]\right)\circ T_{i}^{T}\right\Vert _{L_{v_{0}}^{p}}\cdot\left\Vert f\right\Vert _{V_{i}}\\
 & =C_{1}\cdot\left\Vert \left(v_{0}\circ T_{i}^{-T}\right)\cdot M_{\left[-1,1\right]^{\dimension}}\left[\Fourier^{-1}\theta_{i}^{\natural}\right]\right\Vert _{L^{p}}\cdot\left\Vert f\right\Vert _{V_{i}}\\
\left({\scriptstyle \text{asusmption on }v_{0}\text{ and eq. }\eqref{eq:WeightLinearTransformationsConnection}}\right) & \leq\Omega_{0}^{K}\Omega_{1}C_{1}\cdot\left\Vert \left(1+\left|\mybullet\right|\right)^{K}\cdot M_{\left[-1,1\right]^{\dimension}}\left[\Fourier^{-1}\theta_{i}^{\natural}\right]\right\Vert _{L^{p}}\cdot\left\Vert f\right\Vert _{V_{i}}\\
 & \leq\Omega_{0}^{K}\Omega_{1}\Omega_{2}^{\left(p,K\right)}\cdot C_{1}\cdot\left\Vert f\right\Vert _{V_{i}}.\qedhere
\end{align*}
\end{proof}
Our final ingredient for the construction of the ``reconstruction
map'' $R:V\to\DecompSp{\CalQ}p{\ell_{w}^{q}}v$ is the following
lemma.
\begin{lem}
\label{lem:DecompositionSynthesis}The map
\[
{\rm Synth}_{\CalD}:V\to\DecompSp{\CalQ}p{\ell_{w}^{q}}v,\left(f_{i}\right)_{i\in I}\mapsto\sum_{i\in I}\left[\left(\Fourier^{-1}\varphi_{i}\right)\ast f_{i}\right]\overset{\text{Lemma }\ref{lem:SpecialConvolutionConsistent}}{=}\sum_{i\in I}\left[\Fourier^{-1}\left(\varphi_{i}\cdot\widehat{f_{i}}\right)\right]
\]
is well-defined and bounded with unconditional convergence of the
series in $Z'\left(\CalO\right)$ and with $\vertiii{{\rm Synth}_{\CalD}}\leq\vertiii{\smash{\Gamma_{\CalQ}}}\cdot C$,
where
\[
C=\begin{cases}
\frac{\left(768/\sqrt{\dimension}\right)^{\dimension/p}}{59049\cdot12^{\dimension}\cdot\dimension^{5}}\cdot\!\left(\!186624\cdot\!\dimension^{4}\!\cdot\!\left\lceil K\!+\!\frac{\dimension+1}{p}\right\rceil \!\right)^{\!1+\left\lceil K+\frac{\dimension+1}{p}\right\rceil }\!\cdot\!\left(1\!+\!R_{\CalQ}C_{\CalQ}\right)^{\dimension\left(\frac{2}{p}-1\right)}\!\cdot\!\Omega_{0}^{4K}\Omega_{1}^{4}\cdot N_{\CalQ}^{\frac{1}{p}-1}C_{\CalQ,\Phi,v_{0},p}^{2}, & \text{if }p\in\left(0,1\right),\\
C_{\CalQ,\Phi,v_{0},p}^{2}, & \text{if }p\in\left[1,\infty\right]
\end{cases}
\]
and where $\Gamma_{\CalQ}:\ell_{w}^{q}\left(I\right)\to\ell_{w}^{q}\left(I\right),c\mapsto c^{\ast}$
denotes the $\CalQ$-clustering map, i.e., $c_{i}^{\ast}=\sum_{\ell\in i^{\ast}}c_{\ell}$.
\end{lem}
\begin{proof}
Recall that the Fourier transform $\Fourier:Z'\left(\CalO\right)\to\DistributionSpace{\CalO}$
is an isomorphism that restricts to an isometric isomorphism $\mathcal{F}:\DecompSp{\CalQ}p{\ell_{w}^{q}}v\to\FourierDecompSp{\CalQ}p{\ell_{w}^{q}}v$.
Hence, it suffices to show that the map
\[
\Theta:=\mathcal{F}\circ{\rm Synth}_{\CalD}:V\to\FourierDecompSp{\CalQ}p{\ell_{w}^{q}}v,\left(f_{i}\right)_{i\in I}\mapsto\sum_{i\in I}\varphi_{i}\widehat{f_{i}}
\]
is well-defined and bounded, with unconditional convergence of the
series in $\CalD'\left(\CalO\right)$. Since the $\left(\varphi_{i}\right)_{i\in I}$
form a \emph{locally finite} partition of unity on $\CalO$, the series
\emph{does} converge unconditionally in $\CalD'\left(\CalO\right)$,
given that each term is a well-defined element of $\CalD'\left(\CalO\right)$.
But this is an easy consequence of the embedding $V_{i}\hookrightarrow\Schwartz'\left(\R^{\dimension}\right)$,
which holds for each $i\in I$, according to Lemma \ref{lem:SpecialConvolutionConsistent}.

Now, define
\[
c_{i}:=\left\Vert f_{i}\right\Vert _{V_{i}}=\begin{cases}
\left\Vert f_{i}\right\Vert _{L_{v}^{p}}, & \text{if }p\in\left[1,\infty\right],\\
\left\Vert f_{i}\right\Vert _{W_{T_{i}^{-T}\left[-1,1\right]^{\dimension}}\left(L_{v}^{p}\right)}, & \text{if }p\in\left(0,1\right).
\end{cases}
\]
By definition of $V$, we have $c:=\left(c_{i}\right)_{i\in I}\in\ell_{w}^{q}\left(I\right)$
and $\left\Vert \left(f_{i}\right)_{i\in I}\right\Vert _{V}=\left\Vert c\right\Vert _{\ell_{w}^{q}}$.
Furthermore, since the $\CalQ$-clustering map $\Gamma_{\CalQ}$ is
bounded, it suffices to show $\left\Vert \Fourier^{-1}\left(\varphi_{j}\cdot\Theta\left(f_{i}\right)_{i\in I}\right)\right\Vert _{L_{v}^{p}}\leq C_{1}\cdot c_{j}^{\ast}$
for all $j\in I$ and a suitable constant $C_{1}>0$, since this implies
\[
\left\Vert \Theta\left(f_{i}\right)_{i\in I}\right\Vert _{\FourierDecompSp{\CalQ}p{\ell_{w}^{q}}v}=\left\Vert \left(\left\Vert \Fourier^{-1}\!\left(\varphi_{j}\cdot\Theta\left(f_{i}\right)_{i\in I}\right)\right\Vert _{L_{v}^{p}}\right)_{j\in I}\right\Vert _{\ell_{w}^{q}}\!\leq C_{1}\cdot\left\Vert c^{\ast}\right\Vert _{\ell_{w}^{q}}\leq C_{1}\vertiii{\smash{\Gamma_{\CalQ}}}\cdot\left\Vert c\right\Vert _{\ell_{w}^{q}}=C_{1}\vertiii{\smash{\Gamma_{\CalQ}}}\cdot\left\Vert \left(f_{i}\right)_{i\in I}\right\Vert _{V}.
\]
To show $\left\Vert \Fourier^{-1}\left(\varphi_{j}\cdot\Gamma\left(f_{i}\right)_{i\in I}\right)\right\Vert _{L_{v}^{p}}\leq C_{1}\cdot c_{j}^{\ast}$,
we distinguish two cases regarding $p$:

Let us start with the case $p\in\left[1,\infty\right]$. Since $\varphi_{j}\varphi_{\ell}\equiv0$
unless $\ell\in j^{\ast}$, we have
\begin{align*}
\left\Vert \Fourier^{-1}\left[\varphi_{j}\cdot\Theta\left(f_{\ell}\right)_{\ell\in I}\right]\right\Vert _{L_{v}^{p}} & =\left\Vert \Fourier^{-1}\left[\varphi_{j}\cdot\sum_{\ell\in I}\varphi_{\ell}\widehat{f_{\ell}}\right]\right\Vert _{L_{v}^{p}}\\
 & =\left\Vert \sum_{\ell\in j^{\ast}}\Fourier^{-1}\left[\varphi_{j}\varphi_{\ell}\widehat{f_{\ell}}\right]\right\Vert _{L_{v}^{p}}\\
 & \leq\sum_{\ell\in j^{\ast}}\left\Vert \left(\Fourier^{-1}\varphi_{j}\right)\ast\left(\Fourier^{-1}\varphi_{\ell}\right)\ast f_{\ell}\right\Vert _{L_{v}^{p}}\\
\left({\scriptstyle \text{eq. }\eqref{eq:WeightedYoungInequality}}\right) & \leq\sum_{\ell\in j^{\ast}}\left\Vert \Fourier^{-1}\varphi_{j}\right\Vert _{L_{v_{0}}^{1}}\left\Vert \Fourier^{-1}\varphi_{\ell}\right\Vert _{L_{v_{0}}^{1}}\cdot\left\Vert f_{\ell}\right\Vert _{L_{v}^{p}}\\
 & \leq C_{\CalQ,\Phi,v_{0},p}^{2}\cdot\sum_{\ell\in j^{\ast}}\left\Vert f_{\ell}\right\Vert _{L_{v}^{p}}=C_{\CalQ,\Phi,v_{0},p}^{2}\cdot c_{j}^{\ast},
\end{align*}
so that we can choose $C_{1}:=C_{\CalQ,\Phi,v_{0},p}^{2}$.

Now, we consider the case $p\in\left(0,1\right)$. Here, we can first
argue as before:
\begin{align*}
\left\Vert \Fourier^{-1}\left[\varphi_{j}\cdot\Theta\left(f_{\ell}\right)_{\ell\in I}\right]\right\Vert _{L_{v}^{p}} & =\left\Vert \sum_{\ell\in j^{\ast}}\Fourier^{-1}\left[\varphi_{j}\varphi_{\ell}\widehat{f_{\ell}}\right]\right\Vert _{L_{v}^{p}}\\
\left({\scriptstyle \text{quasi-triangle inequality and }\left|j^{\ast}\right|\leq N_{\CalQ}}\right) & \leq N_{\CalQ}^{\frac{1}{p}-1}\cdot\sum_{\ell\in j^{\ast}}\left\Vert \Fourier^{-1}\left(\varphi_{j}\varphi_{\ell}\right)\ast f_{\ell}\right\Vert _{L_{v}^{p}}\\
\left({\scriptstyle \text{Lemma }\ref{lem:MaximalFunctionDominatesF}}\right) & \leq N_{\CalQ}^{\frac{1}{p}-1}\cdot\sum_{\ell\in j^{\ast}}\left\Vert \Fourier^{-1}\left(\varphi_{j}\varphi_{\ell}\right)\ast f_{\ell}\right\Vert _{W_{T_{\ell}^{-T}\left[-1,1\right]^{\dimension}}\left(L_{v}^{p}\right)}.
\end{align*}
Then, we set $C_{2}:=\Omega_{0}^{3K}\Omega_{1}^{3}\cdot\dimension^{-\frac{\dimension}{2p}}\cdot\left(972\cdot\dimension^{5/2}\right)^{K+\frac{\dimension}{p}}$
and use Corollary \ref{cor:WienerAmalgamConvolutionSimplified} to
estimate each summand as follows:
\begin{align*}
 & \left\Vert \Fourier^{-1}\left(\varphi_{j}\varphi_{\ell}\right)\ast f_{\ell}\right\Vert _{W_{T_{\ell}^{-T}\left[-1,1\right]^{\dimension}}\left(L_{v}^{p}\right)}\\
\left({\scriptstyle \text{Corollary }\ref{cor:WienerAmalgamConvolutionSimplified}}\right) & \leq C_{2}\cdot\left|\det T_{\ell}\right|^{\frac{1}{p}-1}\cdot\left\Vert \Fourier^{-1}\left(\varphi_{j}\varphi_{\ell}\right)\right\Vert _{W_{T_{\ell}^{-T}\left[-1,1\right]^{\dimension}}\left(L_{v_{0}}^{p}\right)}\,\left\Vert f_{\ell}\right\Vert _{W_{T_{\ell}^{-T}\left[-1,1\right]^{\dimension}}\left(L_{v}^{p}\right)}\\
 & =C_{2}\cdot\left|\det T_{\ell}\right|^{\frac{1}{p}-1}\cdot\left\Vert \Fourier^{-1}\left(\varphi_{j}\varphi_{\ell}\right)\right\Vert _{W_{T_{\ell}^{-T}\left[-1,1\right]^{\dimension}}\left(L_{v_{0}}^{p}\right)}\cdot c_{\ell}.
\end{align*}
Now, note 
\[
\supp\left(\varphi_{j}\varphi_{\ell}\right)\subset\overline{Q_{\ell}}=T_{\ell}\overline{Q_{\ell}'}+b_{\ell}\subset T_{\ell}\left[-R_{\CalQ},R_{\CalQ}\right]^{\dimension}+b_{\ell},
\]
so that Theorem \ref{thm:BandlimitedWienerAmalgamSelfImproving} (with
$v_{0}$ instead of $v$) shows for 
\[
C_{3}:=2^{4\left(1+\frac{\dimension}{p}\right)}s_{\dimension}^{\frac{1}{p}}\left(192\cdot\dimension^{3/2}\cdot\left\lceil K+\frac{\dimension+1}{p}\right\rceil \right)^{\left\lceil K+\frac{\dimension+1}{p}\right\rceil +1}\cdot\Omega_{0}^{K}\Omega_{1}\cdot\left(1+R_{\CalQ}\right)^{\frac{\dimension}{p}}
\]
that
\begin{align*}
\left\Vert \Fourier^{-1}\left(\varphi_{j}\varphi_{\ell}\right)\right\Vert _{W_{T_{\ell}^{-T}\left[-1,1\right]^{\dimension}}\left(L_{v_{0}}^{p}\right)} & \leq C_{3}\cdot\left\Vert \Fourier^{-1}\left(\varphi_{j}\varphi_{\ell}\right)\right\Vert _{L_{v_{0}}^{p}}\\
\left({\scriptstyle \text{Proposition }\ref{prop:BandlimitedConvolution}\text{ and }\supp\varphi_{j}\subset\overline{Q_{j}}\subset\overline{Q_{j}^{\ast}}\text{ and }\supp\varphi_{\ell}\subset\overline{Q_{\ell}}\subset\overline{Q_{j}^{\ast}}}\right) & \leq C_{3}\left(12R_{\CalQ}C_{\CalQ}\right)^{\dimension\left(\frac{1}{p}-1\right)}\!\cdot\!\left|\det T_{j}\right|^{\frac{1}{p}-1}\!\cdot\!\left\Vert \Fourier^{-1}\varphi_{j}\right\Vert _{L_{v_{0}}^{p}}\!\cdot\!\left\Vert \Fourier^{-1}\varphi_{\ell}\right\Vert _{L_{v_{0}}^{p}}\\
 & \leq C_{4}\cdot\left\Vert \Fourier^{-1}\varphi_{\ell}\right\Vert _{L_{v_{0}}^{p}}
\end{align*}
for $C_{4}:=C_{3}\cdot\left(12R_{\CalQ}C_{\CalQ}\right)^{\dimension\left(\frac{1}{p}-1\right)}\cdot C_{\CalQ,\Phi,v_{0},p}$.
In total, we conclude
\begin{align*}
\left\Vert \Fourier^{-1}\left(\varphi_{j}\varphi_{\ell}\right)\ast f_{\ell}\right\Vert _{W_{T_{\ell}^{-T}\left[-1,1\right]^{\dimension}}\left(L_{v}^{p}\right)} & \leq C_{2}\cdot\left|\det T_{\ell}\right|^{\frac{1}{p}-1}\cdot\left\Vert \Fourier^{-1}\left(\varphi_{j}\varphi_{\ell}\right)\right\Vert _{W_{T_{\ell}^{-T}\left[-1,1\right]^{\dimension}}\left(L_{v_{0}}^{p}\right)}\cdot c_{\ell}\\
 & \leq C_{2}C_{4}\cdot\left|\det T_{\ell}\right|^{\frac{1}{p}-1}\cdot\left\Vert \Fourier^{-1}\varphi_{\ell}\right\Vert _{L_{v_{0}}^{p}}\cdot c_{\ell}\\
 & \leq C_{2}C_{4}C_{\CalQ,\Phi,v_{0},p}\cdot c_{\ell}
\end{align*}
and hence
\begin{align*}
\left\Vert \Fourier^{-1}\left[\varphi_{j}\cdot\Theta\left(f_{\ell}\right)_{\ell\in I}\right]\right\Vert _{L_{v}^{p}} & \leq N_{\CalQ}^{\frac{1}{p}-1}\cdot\sum_{\ell\in j^{\ast}}\left\Vert \Fourier^{-1}\left(\varphi_{j}\varphi_{\ell}\right)\ast f_{\ell}\right\Vert _{W_{T_{\ell}^{-T}\left[-1,1\right]^{\dimension}}\left(L_{v}^{p}\right)}\\
 & \leq N_{\CalQ}^{\frac{1}{p}-1}C_{2}C_{4}C_{\CalQ,\Phi,v_{0},p}\cdot\sum_{\ell\in j^{\ast}}c_{\ell}\\
 & =N_{\CalQ}^{\frac{1}{p}-1}C_{2}C_{4}C_{\CalQ,\Phi,v_{0},p}\cdot c_{j}^{\ast},
\end{align*}
so that the desired estimate from the start of the proof holds with
$C_{1}:=N_{\CalQ}^{\frac{1}{p}-1}C_{2}C_{4}C_{\CalQ,\Phi,v_{0},p}$.
Now, we finally set $N:=\left\lceil K+\frac{\dimension+1}{p}\right\rceil $
and observe because of $N\geq K+\frac{\dimension}{p}+1\geq\frac{\dimension}{p}+1$
and $s_{\dimension}\leq4^{\dimension}$, as well as $C_{\CalQ}\geq\left\Vert T_{i}^{-1}T_{i}\right\Vert \geq1$
that
\begin{align*}
C_{1} & =\dimension^{-\frac{\dimension}{2p}}\!\cdot\!\left(972\cdot\dimension^{\frac{5}{2}}\right)^{\!K+\frac{\dimension}{p}}2^{4\left(1+\frac{\dimension}{p}\right)}s_{\dimension}^{\frac{1}{p}}\left(192\!\cdot\!\dimension^{\frac{3}{2}}\cdot N\right)^{\!N+1}\!\cdot\left(1\!+\!R_{\CalQ}\right)^{\frac{\dimension}{p}}\left(12R_{\CalQ}C_{\CalQ}\right)^{\dimension\left(\frac{1}{p}-1\right)}\cdot\Omega_{0}^{4K}\Omega_{1}^{4}\cdot N_{\CalQ}^{\frac{1}{p}-1}C_{\CalQ,\Phi,v_{0},p}^{2}\\
 & \leq2^{4}\!\cdot\!\left(2^{6}/\sqrt{\dimension}\right)^{\dimension/p}\left(972\cdot\dimension^{\frac{5}{2}}\right)^{K+\frac{\dimension}{p}}\left(\!192\!\cdot\!\dimension^{\frac{3}{2}}\cdot N\right)^{\!N+1}\!\cdot12^{\dimension\left(\frac{1}{p}-1\right)}\left(1\!+\!R_{\CalQ}C_{\CalQ}\right)^{\dimension\left(\frac{2}{p}-1\right)}\!\cdot\Omega_{0}^{4K}\Omega_{1}^{4}\cdot N_{\CalQ}^{\frac{1}{p}-1}C_{\CalQ,\Phi,v_{0},p}^{2}\\
 & \leq\frac{\left(768/\sqrt{\dimension}\right)^{\dimension/p}}{59049\cdot12^{\dimension}\cdot\dimension^{5}}\cdot\left(186624\cdot\dimension^{4}\cdot N\right)^{N+1}\cdot\left(1+R_{\CalQ}C_{\CalQ}\right)^{\dimension\left(\frac{2}{p}-1\right)}\cdot\Omega_{0}^{4K}\Omega_{1}^{4}\cdot N_{\CalQ}^{\frac{1}{p}-1}C_{\CalQ,\Phi,v_{0},p}^{2}.\qedhere
\end{align*}
\end{proof}
Now, we can finally show existence of the reconstruction map $R$
and thus also derive the estimate ``$\lesssim$'' in equation (\ref{eq:DesiredCharacterization}).
\begin{thm}
\label{thm:SemiDiscreteBanachFrame}Assume that the family $\Gamma=\left(\gamma_{i}\right)_{i\in I}$
satisfies Assumptions \ref{assu:MainAssumptions} and \ref{assu:GammaCoversOrbit}.

Then, with $m_{\theta}$ as in Lemma \ref{lem:LocalInverseConvolution},
with ${\rm Synth}_{\CalD}$ as in Lemma \ref{lem:DecompositionSynthesis}
and with ${\rm Ana}_{\Gamma}$ as in Theorem \ref{thm:ConvolvingDecompositionSpaceWithGammaJ},
the map
\[
R:={\rm Synth}_{\CalD}\circ m_{\theta}:V\to\DecompSp{\CalQ}p{\ell_{w}^{q}}v
\]
is well-defined and bounded with $\vertiii R\leq\vertiii{\Gamma_{\CalQ}}C_{\CalQ,\Phi,v_{0},p}^{2}\cdot C$
for 
\[
C\!:=\!\!\begin{cases}
\frac{\left(768/\dimension\right)^{\dimension/p}}{2^{35}\cdot12^{\dimension}\cdot\dimension^{10}}\cdot\!\left(\!2^{28}\!\cdot\!\dimension^{\frac{13}{2}}\!\cdot\!\left\lceil \!K\!+\!\frac{\dimension+1}{p}\right\rceil \!\right)^{\!1+\left\lceil K+\frac{\dimension+1}{p}\right\rceil }\!\!\cdot\!\left(1\!+\!R_{\CalQ}C_{\CalQ}\right)^{\dimension\left(\frac{2}{p}-1\right)}\!\cdot\Omega_{0}^{8K}\Omega_{1}^{8}\Omega_{2}^{\left(p,K\right)}\cdot N_{\CalQ}^{\frac{1}{p}-1}, & \text{if }p<1,\\
\Omega_{0}^{K}\Omega_{1}\Omega_{2}^{\left(p,K\right)}, & \text{if }p\geq1,
\end{cases}
\]
where as usual $\Gamma_{\CalQ}:\ell_{w}^{q}\left(I\right)\to\ell_{w}^{q}\left(I\right)$
denotes the $\CalQ$-clustering map.

Furthermore, $R$ satisfies
\begin{equation}
R\circ{\rm Ana}_{\Gamma}=\identity_{\DecompSp{\CalQ}p{\ell_{w}^{q}}v}.\label{eq:SemiDiscreteBanachFrame}
\end{equation}
In particular, equation (\ref{eq:DesiredCharacterization}) is satisfied,
i.e., we have
\[
\left\Vert f\right\Vert _{\DecompSp{\CalQ}p{\ell_{w}^{q}}v}\asymp\left\Vert \left(\smash{\gamma^{\left(i\right)}}\ast f\right)_{i\in I}\right\Vert _{V}\qquad\forall f\in\DecompSp{\CalQ}p{\ell_{w}^{q}}v.\qedhere
\]
\end{thm}
\begin{proof}
Boundedness of ${\rm Synth}_{\CalD}$ and $m_{\theta}$ and thus
of $R$ is a consequence of Lemmas \ref{lem:DecompositionSynthesis}
and \ref{lem:LocalInverseConvolution}, respectively, so that it suffices
to prove equation (\ref{eq:SemiDiscreteBanachFrame}).

To see this, we again use the isomorphism $\Fourier:\DecompSp{\CalQ}p{\ell_{w}^{q}}v\to\FourierDecompSp{\CalQ}p{\ell_{w}^{q}}v$.
Recall from Lemma \ref{lem:DecompositionSynthesis} that we have
\[
\left(\Fourier\circ{\rm Synth}_{\CalD}\right)\left(f_{i}\right)_{i\in I}=\sum_{i\in I}\left(\varphi_{i}\cdot\widehat{f_{i}}\right)\qquad\text{ for }\qquad\left(f_{i}\right)_{i\in I}\in V=\ell_{w}^{q}\left(\left[V_{i}\right]_{i\in I}\right),
\]
where it is used that $f_{i}\in V_{i}\hookrightarrow\Schwartz'\left(\R^{\dimension}\right)$
for all $i\in I$.

Hence, for $f\in\DecompSp{\CalQ}p{\ell_{w}^{q}}v$, we have
\begin{align*}
\left(\Fourier\circ R\circ{\rm Ana}_{\Gamma}\right)f & =\sum_{i\in I}\left[\varphi_{i}\cdot\Fourier\left[\left(m_{\theta}\circ{\rm Ana}_{\Gamma}\right)f\right]_{i}\right]\\
 & =\sum_{i\in I}\left[\varphi_{i}\cdot\Fourier\left[\left(\Fourier^{-1}\theta_{i}\right)\ast\left(\gamma^{\left(i\right)}\ast f\right)\right]\right]\\
\left({\scriptstyle \text{Lemma }\ref{lem:SpecialConvolutionConsistent}}\right) & =\sum_{i\in I}\left[\varphi_{i}\cdot\theta_{i}\cdot\widehat{\gamma^{\left(i\right)}\ast f}\right]\\
\left({\scriptstyle \text{Special Def. of }\gamma^{\left(i\right)}\ast f=\Fourier^{-1}f_{i}\text{, cf. Theorem }\ref{thm:ConvolvingDecompositionSpaceWithGammaJ}}\right) & =\sum_{i\in I}\left[\varphi_{i}\cdot\theta_{i}\cdot f_{i}\right],
\end{align*}
where
\[
f_{i}:\Schwartz\left(\R^{\dimension}\right)\to\Compl,\phi\mapsto\sum_{\ell\in I}\left\langle \widehat{\gamma^{\left(i\right)}}\cdot\widehat{f},\,\varphi_{\ell}\phi\right\rangle _{\DistributionSpace{\CalO},\TestFunctionSpace{\CalO}}.
\]

Thus, we have for arbitrary $\phi\in\TestFunctionSpace{\CalO}$ that
\begin{align*}
\left\langle \left(\Fourier\circ R\circ{\rm Ana}_{\Gamma}\right)f,\,\phi\right\rangle _{\DistributionSpace{\CalO},\TestFunctionSpace{\CalO}} & =\sum_{i\in I}\left\langle f_{i},\,\varphi_{i}\cdot\theta_{i}\cdot\phi\right\rangle _{\Schwartz',\Schwartz}\\
 & =\sum_{i\in I}\sum_{\ell\in I}\left\langle \widehat{\gamma^{\left(i\right)}}\cdot\widehat{f},\,\varphi_{\ell}\varphi_{i}\cdot\theta_{i}\cdot\phi\right\rangle _{\DistributionSpace{\CalO},\TestFunctionSpace{\CalO}}\\
 & =\sum_{i\in I}\sum_{\ell\in I}\left\langle \widehat{f},\,\varphi_{\ell}\varphi_{i}\cdot\widehat{\gamma^{\left(i\right)}}\cdot\theta_{i}\cdot\phi\right\rangle _{\DistributionSpace{\CalO},\TestFunctionSpace{\CalO}}\\
\left({\scriptstyle \widehat{\gamma^{\left(i\right)}}\cdot\theta_{i}\equiv1\text{ on }\overline{Q_{i}}\supset\supp\varphi_{i},\text{ cf. Lemma }\ref{lem:LocalInverseConvolution}}\right) & =\sum_{i\in I}\sum_{\ell\in I}\left\langle \widehat{f},\,\varphi_{\ell}\varphi_{i}\cdot\phi\right\rangle _{\DistributionSpace{\CalO},\TestFunctionSpace{\CalO}}\\
\left({\scriptstyle \phi\in\TestFunctionSpace{\CalO}\text{ and }\left(\varphi_{i}\right)_{i\in I}\text{ loc. finite part. of unity on }\CalO}\right) & =\sum_{i\in I}\left\langle \widehat{f},\,\varphi_{i}\cdot\phi\right\rangle _{\DistributionSpace{\CalO},\TestFunctionSpace{\CalO}}\\
\left({\scriptstyle \text{as above}}\right) & =\left\langle \widehat{f},\,\phi\right\rangle _{\DistributionSpace{\CalO},\TestFunctionSpace{\CalO}}.
\end{align*}
Hence, we have shown $\Fourier\circ R\circ{\rm Ana}_{\Gamma}=\Fourier$
on $\DecompSp{\CalQ}p{\ell_{w}^{q}}v$. Since $\Fourier:\DecompSp{\CalQ}p{\ell_{w}^{q}}v\to\FourierDecompSp{\CalQ}p{\ell_{w}^{q}}v$
is an isomorphism, this implies $R\circ{\rm Ana}_{\Gamma}=\identity_{\DecompSp{\CalQ}p{\ell_{w}^{q}}v}$,
as desired.
\end{proof}

\section{Fully Discrete Banach Frames}

\label{sec:FullyDiscreteBanachFrames}In the previous section, we
obtained \emph{semi-discrete} Banach frames for $\DecompSp{\CalQ}p{\ell_{w}^{q}}v$;
in particular, we showed
\[
\left\Vert f\right\Vert _{\DecompSp{\CalQ}p{\ell_{w}^{q}}v}\asymp\left\Vert \left(\left\Vert \smash{\gamma^{\left(i\right)}}\ast f\right\Vert _{V_{i}}\right)_{i\in I}\right\Vert _{\ell_{w}^{q}}.
\]
We call such a frame \textbf{semi-discrete}, because while the index
set $I$ is discrete, the convolutions $\gamma^{\left(i\right)}\ast f$
are treated as genuine functions, which are defined on the continuous
(non-discrete) index set $\R^{\dimension}$.

In this section, our aim is a further discretization of these frames,
so that we will in the end obtain a (quasi)-norm equivalence of the
form
\begin{equation}
\left\Vert f\right\Vert _{\DecompSp{\CalQ}p{\ell_{w}^{q}}v}\asymp\left\Vert \left(\left\Vert \left[\left(\gamma^{\left[j\right]}\ast f\right)\!\!\left(\delta\cdot T_{j}^{-T}k\right)\right]_{k\in\Z^{\dimension}}\right\Vert _{C_{j}^{\left(\delta\right)}}\right)_{\!\!j\in I}\right\Vert _{\ell_{u_{q}}^{q}}\qquad\text{ for a suitable weight }u_{q}\text{ on }I,\label{eq:DesiredBanachFrameNormEquivalence}
\end{equation}
where for each $\delta\in\left(0,1\right]$ and $j\in I$, the \textbf{coefficient
space} $C_{j}^{\left(\delta\right)}$ is given by
\begin{equation}
C_{j}^{\left(\delta\right)}:=\ell_{v^{\left(j,\delta\right)}}^{p}\left(\smash{\Z^{\dimension}}\right)\quad\text{ with }\quad v_{k}^{\left(j,\delta\right)}=v\left(\delta\cdot T_{j}^{-T}k\right)\quad\text{ for }k\in\Z^{\dimension}.\label{eq:CoefficientSpaceDefinition}
\end{equation}
For simplicity, the reader should keep in mind the important special
case $v\equiv1$, for which $C_{j}^{\left(\delta\right)}=\ell^{p}\left(\Z^{\dimension}\right)$,
independently of $j,\delta$.

To ensure that equation (\ref{eq:DesiredBanachFrameNormEquivalence})
holds, we will introduce suitable assumptions on $\left(\gamma_{i}\right)_{i\in I}$
and $\delta>0$. In particular, for the formula above to make sense,
we also have to establish (at least) continuity of $\gamma^{\left(j\right)}\ast f$
(and thus of $\gamma^{\left[j\right]}\ast f$), so that the pointwise
evaluations $\left(\gamma^{\left[j\right]}\ast f\right)\left(\delta\cdot T_{j}^{-T}k\right)$
are meaningful; note that up to now, we only know (for $p\in\left[1,\infty\right]$)
that $\gamma^{\left[j\right]}\ast f\in L_{v}^{p}\left(\R^{\dimension}\right)$.

To ensure this continuity, we introduce a new set of additional assumptions
and notations. In these assumptions, the $L_{v_{0}}^{p}$ (quasi)-norm
of certain vector valued functions $g:\R^{\dimension}\to\Compl^{k}$
is considered. This has to be understood as $\left\Vert g\right\Vert _{L_{v_{0}}^{p}}:=\left\Vert \,\left|g\right|\,\right\Vert _{L_{v_{0}}^{p}}$,
where as usual $\left|g\right|\left(x\right):=\left|g\left(x\right)\right|=\left\Vert g\left(x\right)\right\Vert _{2}$
denotes the euclidean norm of $g\left(x\right)$. Furthermore, for
such a function $g=\left(g_{1},\dots,g_{k}\right)$, expressions as
the (inverse) Fourier transform $\Fourier^{-1}g:=\left(\Fourier^{-1}g_{1},\dots,\Fourier^{-1}g_{k}\right)$
are always understood in a coordinatewise sense.
\begin{assumption}
\label{assu:DiscreteBanachFrameAssumptions}In addition to Assumption
\ref{assu:MainAssumptions}, we assume the following:

\begin{enumerate}
\item We have $\gamma_{i}\in C^{1}\left(\R^{\dimension}\right)$ for all
$i\in I$ and the gradient $\phi_{i}:=\nabla\gamma_{i}$ satisfies
the following:

\begin{enumerate}
\item $\phi_{i}$ is bounded for each $i\in I$,
\item We have $\phi_{i}\in L_{v_{0}}^{1}\left(\R^{\dimension};\Compl^{\dimension}\right)\hookrightarrow L^{1}\left(\R^{\dimension};\Compl^{\dimension}\right)$
for all $i\in I$,
\item We have $\widehat{\phi_{i}}\in C^{\infty}\left(\R^{\dimension};\Compl^{\dimension}\right)$
for all $i\in I$.
\end{enumerate}
\item For $j\in I$, we define
\[
\phi^{\left(j\right)}:=\Fourier^{-1}\left(\widehat{\phi_{j}}\circ S_{j}^{-1}\right)=\left|\det T_{j}\right|\cdot M_{b_{j}}\left[\phi_{j}\circ T_{j}^{T}\right],
\]
so that $\phi^{\left(j\right)}$ is to $\phi_{j}$ as $\gamma^{\left(j\right)}$
is to $\gamma_{j}$.
\item For $j,i\in I$, set
\[
B_{j,i}:=\begin{cases}
\left(1+\left\Vert T_{j}^{-1}T_{i}\right\Vert \right)^{K+\dimension}\cdot\left\Vert \Fourier^{-1}\left(\varphi_{i}\cdot\widehat{\phi^{\left(j\right)}}\right)\right\Vert _{L_{v_{0}}^{1}}, & \text{if }p\in\left[1,\infty\right],\\
\left(1+\left\Vert T_{j}^{-1}T_{i}\right\Vert \right)^{pK+\dimension}\cdot\left|\det T_{i}\right|^{1-p}\cdot\left\Vert \Fourier^{-1}\left(\varphi_{i}\cdot\widehat{\phi^{\left(j\right)}}\right)\right\Vert _{L_{v_{0}}^{p}}^{p}, & \text{if }p\in\left(0,1\right).
\end{cases}
\]
\item With 
\[
r:=\max\left\{ q,\frac{q}{p}\right\} =\begin{cases}
q, & \text{if }p\in\left[1,\infty\right],\\
\frac{q}{p}, & \text{if }p\in\left(0,1\right)
\end{cases}
\]
as in Assumption \ref{assu:MainAssumptions}, we assume that the operator
$\overrightarrow{B}$ induced by $\left(B_{j,i}\right)_{j,i\in I}$,
i.e.\@
\[
\overrightarrow{B}\left(c_{i}\right)_{i\in I}:=\left(\sum_{i\in I}B_{j,i}\,c_{i}\right)_{j\in I}
\]
defines a well-defined, bounded operator $\overrightarrow{B}:\ell_{w^{\min\left\{ 1,p\right\} }}^{r}\left(I\right)\to\ell_{w^{\min\left\{ 1,p\right\} }}^{r}\left(I\right)$.
\item For $j\in I$ and $\delta\in\left(0,1\right]$, we let the $j$-th
\textbf{coefficient space} $C_{j}^{\left(\delta\right)}$ be defined
as in equation (\ref{eq:CoefficientSpaceDefinition}) and set
\[
W_{j}\!:=\!\left\{ f:\R^{\dimension}\to\Compl\with f\text{ continuous and }\left\Vert f\right\Vert _{W_{j}}\!<\!\infty\right\} ,
\]
where
\[
\left\Vert f\right\Vert _{W_{j}}:=\left\Vert f\right\Vert _{V_{j}}+\sup_{0<\delta\leq1}\!\frac{1}{\delta}\left\Vert \osc{\delta\cdot T_{j}^{-T}\left[-1,1\right]^{\dimension}}\left[M_{-b_{j}}f\right]\right\Vert _{V_{j}}.
\]
\item Finally, we define
\[
\ell_{w}^{q}\left(\smash{\left[W_{j}\right]_{j\in I}}\right):=\left\{ \left(f_{j}\right)_{j\in I}\with\left(\forall j\in I:\:f_{j}\in W_{j}\right)\text{ and }\left(\smash{\left\Vert f_{j}\right\Vert _{W_{j}}}\right)_{j\in I}\in\ell_{w}^{q}\left(I\right)\right\} ,
\]
equipped with the natural (quasi)-norm $\left\Vert \smash{\left(f_{j}\right)_{j\in I}}\right\Vert _{\ell_{w}^{q}\left(\left[W_{j}\right]_{j\in I}\right)}:=\left\Vert \left(\smash{\left\Vert f_{j}\right\Vert _{W_{j}}}\vphantom{\sum}\right)_{j\in I}\right\Vert _{\ell_{w}^{q}}$.\qedhere
\end{enumerate}
\end{assumption}
\begin{rem*}
Note that $\varphi_{i}\cdot\widehat{\phi^{\left(j\right)}}\in\TestFunctionSpace{\R^{\dimension}}$
for all $i,j\in I$, since $\varphi_{i}\in\TestFunctionSpace{\R^{\dimension}}$
and since $\widehat{\phi^{\left(j\right)}}=\widehat{\phi_{j}}\circ S_{j}^{-1}$
is smooth, because $\widehat{\phi_{j}}$ is. Hence, $\Fourier^{-1}\left(\varphi_{i}\cdot\widehat{\phi^{\left(j\right)}}\right)\in\Schwartz\left(\R^{\dimension}\right)$,
so that $B_{j,i}<\infty$ for all $i,j\in I$, cf.\@ Lemma \ref{lem:SchwartzFunctionsAreWiener}.

Although we again stated the assumption in the general case where
the prototype $\gamma_{i}$ may depend on $i\in I$, the reader should
keep in mind the most important special case where $\gamma_{i}=\gamma$
is independent of $i\in I$.
\end{rem*}
Given these assumptions, we now want to show in particular that $\gamma^{\left[j\right]}\ast f$
is continuous for each $f\in\DecompSp{\CalQ}p{\ell_{w}^{q}}v$. The
following lemma makes an important step in that direction.
\begin{lem}
\label{lem:MainOscillationEstimate}Assume that $\Gamma=\left(\gamma_{i}\right)_{i\in I}$
satisfies Assumption \ref{assu:DiscreteBanachFrameAssumptions}.
Then the following hold:

For $p\in\left[1,\infty\right]$ and
\[
C:=\frac{2^{6\dimension}}{\sqrt{\dimension}}\cdot\left(1152\cdot\dimension^{5/2}\cdot\left\lceil K+\dimension+1\right\rceil \right)^{\left\lceil K\right\rceil +\dimension+2}\cdot\Omega_{0}^{2K}\Omega_{1}^{2}\cdot\left(1+R_{\CalQ}\right)^{\dimension},
\]
we have
\[
\left\Vert {\rm osc}_{\delta\cdot T_{j}^{-T}\left[-1,1\right]^{\dimension}}\left[M_{-b_{j}}\Fourier^{-1}\left(\widehat{\gamma^{\left(j\right)}}\cdot\varphi_{i}\widehat{f}\right)\right]\right\Vert _{L_{v}^{p}}\leq C\cdot\delta\cdot\left(1+\left\Vert T_{j}^{-1}T_{i}\right\Vert \right)^{K+\dimension}\cdot\left\Vert \Fourier^{-1}\left(\varphi_{i}\cdot\smash{\widehat{\phi^{\left(j\right)}}}\right)\right\Vert _{L_{v_{0}}^{1}}\cdot\left\Vert \Fourier^{-1}\left(\varphi_{i}^{\ast}\smash{\widehat{f}}\,\right)\right\Vert _{L_{v}^{p}}
\]
for all $0<\delta\leq1$, all $f\in Z'\left(\CalO\right)$ and all
$i,j\in I$.

Likewise, for $p\in\left(0,1\right)$ and
\[
C:=\frac{2^{16\frac{\dimension}{p}}\cdot\left(1+C_{\CalQ}R_{\CalQ}\right)^{\frac{2\dimension}{p}}}{370\cdot\dimension^{11/2}\cdot\dimension^{\dimension/2p}}\cdot\left(4032\cdot\dimension^{3}\cdot\left\lceil K+\frac{\dimension+1}{p}\right\rceil \right)^{2\left\lceil K+\frac{\dimension+1}{p}\right\rceil +2}\cdot\Omega_{0}^{5K}\Omega_{1}^{5},
\]
we have
\begin{align*}
 & \left\Vert {\rm osc}_{\delta\cdot T_{j}^{-T}\left[-1,1\right]^{\dimension}}\!\left[M_{-b_{j}}\Fourier^{-1}\!\left(\varphi_{i}\widehat{\gamma^{\left(j\right)}}\smash{\widehat{f}}\right)\right]\right\Vert _{W_{T_{j}^{-T}\left[-1,1\right]^{\dimension}}\left(L_{v}^{p}\right)}\\
 & \leq\vphantom{\sum^{j}}C\cdot\delta\cdot\left|\det T_{i}\right|^{\frac{1}{p}-1}\cdot\left(1+\left\Vert T_{j}^{-1}T_{i}\right\Vert \right)^{K+\frac{\dimension}{p}}\cdot\left\Vert \Fourier^{-1}\!\left(\varphi_{i}\smash{\widehat{\phi^{\left(j\right)}}}\right)\right\Vert _{L_{v_{0}}^{p}}\cdot\left\Vert \Fourier^{-1}\!\left(\varphi_{i}^{\ast}\smash{\widehat{f}}\,\right)\right\Vert _{L_{v}^{p}}
\end{align*}
for all $0<\delta\leq1$, all $f\in Z'\left(\CalO\right)$ and all
$i,j\in I$.
\end{lem}
\begin{proof}
First of all, note that $\widehat{f}\in\CalD'\left(\CalO\right)$
for $f\in Z'\left(\CalO\right)$. Because of $\varphi_{i}\in\TestFunctionSpace{\CalO}$,
this implies that $\varphi_{i}\cdot\widehat{f}\in\Schwartz'\left(\R^{\dimension}\right)$
is a well-defined tempered distribution with compact support, so that
the same holds for $\varphi_{i}\cdot\widehat{\gamma^{\left(j\right)}}\cdot\widehat{f}$,
since $\widehat{\gamma^{\left(j\right)}}\in C^{\infty}\left(\R^{\dimension}\right)$.
Hence, by the Paley-Wiener theorem, $\Fourier^{-1}\left(\widehat{\gamma^{\left(j\right)}}\cdot\varphi_{i}\cdot\widehat{f}\right)$
is a smooth (even analytic) function with polynomially bounded derivatives
of all orders. In particular, expressions like $\left({\rm osc}_{\delta\cdot T_{j}^{-T}\left[-1,1\right]^{\dimension}}\!\left[M_{-b_{j}}\Fourier^{-1}\!\left(\varphi_{i}\widehat{\gamma^{\left(j\right)}}\smash{\widehat{f}}\right)\right]\right)\left(x\right)$
are well-defined for every $x\in\R^{\dimension}$.

Let $f\in Z'\left(\CalO\right)$ be arbitrary. We can clearly assume
$\left\Vert \Fourier^{-1}\left(\varphi_{i}^{\ast}\cdot\smash{\widehat{f}}\,\right)\right\Vert _{L_{v}^{p}}<\infty$,
since otherwise, the claim is trivial. Now, note that $\widehat{\gamma^{\left(j\right)}}\cdot\varphi_{i}\in\TestFunctionSpace{\R^{\dimension}}\subset\Schwartz\left(\R^{\dimension}\right)$
and $\varphi_{i}^{\ast}\widehat{f}\in\Schwartz'\left(\R^{\dimension}\right)$,
as well as
\begin{align*}
M_{-b_{j}}\left[\Fourier^{-1}\left(\widehat{\gamma^{\left(j\right)}}\cdot\varphi_{i}\cdot\widehat{f}\right)\right] & =\Fourier^{-1}\left[L_{-b_{j}}\left(\widehat{\gamma^{\left(j\right)}}\varphi_{i}\cdot\varphi_{i}^{\ast}\widehat{f}\right)\right]\\
 & =\Fourier^{-1}\left[L_{-b_{j}}\left(\widehat{\gamma^{\left(j\right)}}\cdot\varphi_{i}\right)\right]\ast\Fourier^{-1}\left[L_{-b_{j}}\left(\varphi_{i}^{\ast}\smash{\widehat{f}}\,\right)\right].
\end{align*}
In particular, the convolution is pointwise well-defined, so that
Lemma \ref{lem:OscillationConvolution} shows
\begin{equation}
{\rm osc}_{\delta\cdot T_{j}^{-T}\left[-1,1\right]^{\dimension}}M_{-b_{j}}\left[\Fourier^{-1}\left(\widehat{\gamma^{\left(j\right)}}\cdot\varphi_{i}\smash{\widehat{f}}\,\right)\right]\leq\left(\osc{\delta\cdot T_{j}^{-T}\left[-1,1\right]^{\dimension}}\Fourier^{-1}\left[L_{-b_{j}}\left(\smash{\widehat{\gamma^{\left(j\right)}}}\cdot\varphi_{i}\right)\right]\right)\ast\left|\Fourier^{-1}\left[L_{-b_{j}}\left(\varphi_{i}^{\ast}\smash{\widehat{f}}\,\right)\right]\right|.\label{eq:MainOscillationEstimateOscillationToTheLeft}
\end{equation}

\medskip{}

Now, for $p\in\left(0,1\right)$, we want to apply Proposition \ref{prop:AlternativeWienerAmalgamConvolution}
with $Q_{1}=T_{i}^{-T}\left[-1,1\right]^{\dimension}$, $Q_{2}=T_{j}^{-T}\left[-1,1\right]^{\dimension}$
and 
\[
g=\left|\Fourier^{-1}\left[L_{-b_{j}}\left(\varphi_{i}^{\ast}\widehat{f}\right)\right]\right|,\qquad\text{ as well as }\qquad f=\osc{\delta\cdot T_{j}^{-T}\left[-1,1\right]^{\dimension}}\Fourier^{-1}\left[L_{-b_{j}}\left(\widehat{\gamma^{\left(j\right)}}\cdot\varphi_{i}\right)\right].
\]
To this end, first note just as in the proof of Corollary \ref{cor:WienerAmalgamConvolutionSimplified}
(cf.\@ equation (\ref{eq:LinearImageOfCubePartition})) that with
this choice of $Q_{1}$ and suitable choices of $\left(x_{i}\right)_{i\in\Z^{\dimension}}$,
the constant $N$ from Theorem \ref{thm:WienerAmalgamConvolution}
satisfies $N\leq3^{\dimension}$.

Next, we use the identities $Q_{1}-Q_{1}=T_{i}^{-T}\left[-2,2\right]^{\dimension}$
and $\left|\Fourier^{-1}\left[L_{b}h\right]\right|=\left|M_{b}\left[\Fourier^{-1}h\right]\right|=\left|\Fourier^{-1}h\right|$
and equation (\ref{eq:WienerLinearCubeEnlargement}), as well as
Theorem \ref{thm:BandlimitedWienerAmalgamSelfImproving} to get
\begin{align*}
\left\Vert \,\left|\Fourier^{-1}\left[L_{-b_{j}}\left(\varphi_{i}^{\ast}\smash{\widehat{f}}\,\right)\right]\right|\,\right\Vert _{W_{Q_{1}-Q_{1}}\left(L_{v}^{p}\right)} & =\left\Vert \Fourier^{-1}\left(\varphi_{i}^{\ast}\smash{\widehat{f}}\,\right)\right\Vert _{W_{T_{i}^{-T}\left[-2,2\right]^{\dimension}}\left(L_{v}^{p}\right)}\\
\left({\scriptstyle \text{eq. }\eqref{eq:WienerLinearCubeEnlargement}}\right) & \leq\Omega_{0}^{K}\Omega_{1}\cdot\left(18\dimension\right)^{K+\frac{\dimension}{p}}\cdot\left\Vert \Fourier^{-1}\left(\varphi_{i}^{\ast}\smash{\widehat{f}}\,\right)\right\Vert _{W_{T_{i}^{-T}\left[-1,1\right]^{\dimension}}\left(L_{v}^{p}\right)}\\
\left({\scriptstyle \text{Theorem }\ref{thm:BandlimitedWienerAmalgamSelfImproving}}\right) & \leq\Omega_{0}^{K}\Omega_{1}\cdot\left(18\dimension\right)^{K+\frac{\dimension}{p}}C_{1}\cdot\left\Vert \Fourier^{-1}\left(\varphi_{i}^{\ast}\smash{\widehat{f}}\,\right)\right\Vert _{L_{v}^{p}}
\end{align*}
for
\[
C_{1}:=2^{4\left(1+\frac{\dimension}{p}\right)}s_{\dimension}^{\frac{1}{p}}\left(192\cdot\dimension^{3/2}\cdot\left\lceil K+\frac{\dimension+1}{p}\right\rceil \right)^{\left\lceil K+\frac{\dimension+1}{p}\right\rceil +1}\cdot\Omega_{0}^{K}\Omega_{1}\cdot\left(1+\left(1+2C_{\CalQ}\right)R_{\CalQ}\right)^{\frac{\dimension}{p}},
\]
since \cite[Lemma 2.7]{DecompositionEmbedding} shows $\supp\left(\varphi_{i}^{\ast}\widehat{f}\right)\subset\overline{Q_{i}^{\ast}}\subset T_{i}\left[\overline{B_{R}}\left(0\right)\right]+b_{i}\subset T_{i}\left[-R,R\right]^{\dimension}+b_{i}$
for $R=\left(1+2C_{\CalQ}\right)R_{\CalQ}$.

All in all, we now set $C_{2}:=\Omega_{0}^{K}\Omega_{1}\cdot\left(18\dimension\right)^{K+\frac{\dimension}{p}}C_{1}$
and use equation (\ref{eq:MainOscillationEstimateOscillationToTheLeft}),
Proposition \ref{prop:AlternativeWienerAmalgamConvolution}, and
the identity $\widehat{\gamma^{\left(j\right)}}=L_{b_{j}}\left(\widehat{\gamma_{j}}\circ T_{j}^{-1}\right)$
to conclude because of
\begin{align*}
\sup_{x\in Q_{1}}v_{0}\left(x\right) & \leq\Omega_{1}\cdot\sup_{x\in T_{i}^{-T}\left[-1,1\right]^{\dimension}}\left(1+\left|x\right|\right)^{K}\\
\left({\scriptstyle \text{eq. }\eqref{eq:WeightLinearTransformationsConnection}}\right) & \leq\Omega_{0}^{K}\Omega_{1}\cdot\sup_{x\in T_{i}^{-T}\left[-1,1\right]^{\dimension}}\left(1+\left|T_{i}^{T}x\right|\right)^{K}\\
 & =\Omega_{0}^{K}\Omega_{1}\cdot\sup_{y\in\left[-1,1\right]^{\dimension}}\left(1+\left|y\right|\right)^{K}\\
 & \leq\Omega_{0}^{K}\Omega_{1}\left(1+\sqrt{\dimension}\right)^{K}
\end{align*}
that
\begin{align*}
 & \left\Vert {\rm osc}_{\delta\cdot T_{j}^{-T}\left[-1,1\right]^{\dimension}}\left[M_{-b_{j}}\Fourier^{-1}\left(\widehat{\gamma^{\left(j\right)}}\cdot\varphi_{i}\widehat{f}\right)\right]\right\Vert _{W_{T_{j}^{-T}\left[-1,1\right]^{\dimension}}\left(L_{v}^{p}\right)}\\
 & \leq3^{\frac{\dimension}{p}}\Omega_{0}^{K}\Omega_{1}\left(1+\sqrt{\dimension}\right)^{K}\cdot\left[\lambda_{\dimension}\left(Q_{1}\right)\right]^{1-\frac{1}{p}}\cdot\left\Vert \left|\Fourier^{-1}\left[L_{-b_{j}}\left(\varphi_{i}^{\ast}\smash{\widehat{f}}\,\right)\right]\right|\right\Vert _{W_{Q_{1}-Q_{1}}\left(L_{v}^{p}\right)}\\
 & \phantom{\leq}\cdot\left\Vert {\rm osc}_{\delta\cdot T_{j}^{-T}\left[-1,1\right]^{\dimension}}\left(\Fourier^{-1}\left[L_{-b_{j}}\left(\widehat{\gamma^{\left(j\right)}}\cdot\varphi_{i}\right)\right]\right)\right\Vert _{W_{T_{j}^{-T}\left[-1,1\right]^{\dimension}-T_{i}^{-T}\left[-1,1\right]^{\dimension}}\left(L_{v_{0}}^{p}\right)}\\
 & \leq3^{\frac{\dimension}{p}}\Omega_{0}^{K}\Omega_{1}\left(1+\sqrt{\dimension}\right)^{K}C_{2}\cdot2^{\dimension\left(1-\frac{1}{p}\right)}\left|\det T_{i}\right|^{\frac{1}{p}-1}\cdot\left\Vert \Fourier^{-1}\left(\varphi_{i}^{\ast}\smash{\widehat{f}}\,\right)\right\Vert _{L_{v}^{p}}\!\\
 & \phantom{=}\cdot\left\Vert {\rm osc}_{\delta\cdot T_{j}^{-T}\left[-1,1\right]^{\dimension}}\left(\Fourier^{-1}\left[L_{-b_{j}}\left(L_{b_{j}}\left(\widehat{\gamma_{j}}\circ T_{j}^{-1}\right)\cdot\varphi_{i}\right)\right]\right)\right\Vert _{W_{T_{j}^{-T}\left[-1,1\right]^{\dimension}-T_{i}^{-T}\left[-1,1\right]^{\dimension}}\left(L_{v_{0}}^{p}\right)}\\
 & \leq\left(2\sqrt{\dimension}\right)^{K}3^{\frac{\dimension}{p}}C_{2}\cdot\Omega_{0}^{K}\Omega_{1}\cdot\left|\det T_{i}\right|^{\frac{1}{p}-1}\cdot\left\Vert \Fourier^{-1}\left(\varphi_{i}^{\ast}\smash{\widehat{f}}\,\right)\right\Vert _{L_{v}^{p}}\\
 & \phantom{\leq}\cdot\left\Vert {\rm osc}_{\delta\cdot T_{j}^{-T}\left[-1,1\right]^{\dimension}}\left(\Fourier^{-1}\left[\left(\widehat{\gamma_{j}}\cdot\left[\left(L_{-b_{j}}\varphi_{i}\right)\circ T_{j}\right]\right)\circ T_{j}^{-1}\right]\right)\right\Vert _{W_{T_{j}^{-T}\left[-1,1\right]^{\dimension}-T_{i}^{-T}\left[-1,1\right]^{\dimension}}\left(L_{v_{0}}^{p}\right)}.
\end{align*}

Now, we recall $\phi_{j}=\nabla\gamma_{j}$ and estimate
\begin{align*}
 & \left\Vert {\rm osc}_{\delta\cdot T_{j}^{-T}\left[-1,1\right]^{\dimension}}\left(\Fourier^{-1}\left[\left(\widehat{\gamma_{j}}\cdot\left[\left(L_{-b_{j}}\varphi_{i}\right)\circ T_{j}\right]\right)\circ T_{j}^{-1}\right]\right)\right\Vert _{W_{T_{j}^{-T}\left[-1,1\right]^{\dimension}-T_{i}^{-T}\left[-1,1\right]^{\dimension}}\left(L_{v_{0}}^{p}\right)}\\
 & =\left|\det T_{j}\right|\cdot\left\Vert {\rm osc}_{\delta\cdot T_{j}^{-T}\left[-1,1\right]^{\dimension}}\left[\left(\Fourier^{-1}\left[\widehat{\gamma_{j}}\cdot\left(\left[L_{-b_{j}}\varphi_{i}\right]\circ T_{j}\right)\right]\right)\circ T_{j}^{T}\right]\right\Vert _{W_{T_{j}^{-T}\left[-1,1\right]^{\dimension}-T_{i}^{-T}\left[-1,1\right]^{\dimension}}\left(L_{v_{0}}^{p}\right)}\\
\left({\scriptstyle \text{Lem. }\ref{lem:OscillationLinearChange},\,\ref{lem:WienerTransformationFormula}}\right) & =\left|\det T_{j}\right|\cdot\left\Vert \left(M_{\left[-1,1\right]^{\dimension}-T_{j}^{T}T_{i}^{-T}\left[-1,1\right]^{\dimension}}\left[{\rm osc}_{\delta\left[-1,1\right]^{\dimension}}\left(\Fourier^{-1}\left[\widehat{\gamma_{j}}\cdot\left(\left[L_{-b_{j}}\varphi_{i}\right]\circ T_{j}\right)\right]\right)\right]\right)\circ T_{j}^{T}\right\Vert _{L_{v_{0}}^{p}}\\
 & =\left|\det T_{j}\right|^{1-\frac{1}{p}}\left\Vert \left(v_{0}\circ T_{j}^{-T}\right)\cdot M_{\left[-1,1\right]^{\dimension}-T_{j}^{T}T_{i}^{-T}\left[-1,1\right]^{\dimension}}\left[{\rm osc}_{\delta\left[-1,1\right]^{\dimension}}\left(\Fourier^{-1}\left[\widehat{\gamma_{j}}\cdot\left(\left[L_{-b_{j}}\varphi_{i}\right]\circ T_{j}\right)\right]\right)\right]\right\Vert _{L^{p}}\\
\left({\scriptstyle \text{Lem. }\ref{lem:OscillationEstimatedByWienerDerivative}}\right) & \leq2\sqrt{\dimension}\cdot\delta\cdot\left|\det T_{j}\right|^{1-\frac{1}{p}}\!\cdot\!\left\Vert \left(v_{0}\circ T_{j}^{-T}\right)\!\cdot M_{\left[-1,1\right]^{\dimension}-T_{j}^{T}T_{i}^{-T}\left[-1,1\right]^{\dimension}}\!\left[M_{\delta\left[-1,1\right]^{\dimension}}\nabla\left(\Fourier^{-1}\!\left[\widehat{\gamma_{j}}\cdot\!\left(\left[L_{-b_{j}}\varphi_{i}\right]\!\circ\!T_{j}\right)\right]\right)\right]\right\Vert _{L^{p}}\!\!\!.
\end{align*}
Since we have $\delta\leq1$, Lemma \ref{lem:IteratedMaximalFunction}
allows us to continue the estimate as follows:
\begin{align}
\dots & \leq2\sqrt{\dimension}\cdot\delta\cdot\left|\det T_{j}\right|^{1-\frac{1}{p}}\!\cdot\!\left\Vert \left(v_{0}\circ T_{j}^{-T}\right)\!\cdot M_{\left[-2,2\right]^{\dimension}-T_{j}^{T}T_{i}^{-T}\left[-1,1\right]^{\dimension}}\left[\nabla\left(\gamma_{j}\ast\Fourier^{-1}\left[\left(L_{-b_{j}}\varphi_{i}\right)\circ T_{j}\right]\right)\right]\right\Vert _{L^{p}}\nonumber \\
\left({\scriptstyle \nabla\left(f\ast g\right)=\left(\nabla f\right)\ast g}\right) & \overset{\left(\ast\right)}{=}2\sqrt{\dimension}\cdot\delta\cdot\left|\det T_{j}\right|^{1-\frac{1}{p}}\!\cdot\!\left\Vert \left(v_{0}\circ T_{j}^{-T}\right)\!\cdot M_{\left[-2,2\right]^{\dimension}-T_{j}^{T}T_{i}^{-T}\left[-1,1\right]^{\dimension}}\left(\phi_{j}\ast\Fourier^{-1}\left[\left(L_{-b_{j}}\varphi_{i}\right)\circ T_{j}\right]\right)\right\Vert _{L^{p}}\nonumber \\
 & =2\sqrt{\dimension}\cdot\delta\cdot\left|\det T_{j}\right|^{1-\frac{1}{p}}\!\cdot\!\left\Vert \left(v_{0}\circ T_{j}^{-T}\right)\!\cdot M_{\left[-2,2\right]^{\dimension}-T_{j}^{T}T_{i}^{-T}\left[-1,1\right]^{\dimension}}\!\left[\Fourier^{-1}\!\left(\left[\left(\widehat{\phi_{j}}\circ T_{j}^{-1}\right)\!\cdot\!\left(L_{-b_{j}}\varphi_{i}\right)\right]\circ T_{j}\right)\right]\right\Vert _{L^{p}}\nonumber \\
 & =2\sqrt{\dimension}\cdot\delta\cdot\left|\det T_{j}\right|^{-\frac{1}{p}}\!\cdot\!\left\Vert \left(v_{0}\circ T_{j}^{-T}\right)\!\cdot M_{\left[-2,2\right]^{\dimension}-T_{j}^{T}T_{i}^{-T}\left[-1,1\right]^{\dimension}}\!\left[\left(\Fourier^{-1}\!\left[\left(\widehat{\phi_{j}}\circ T_{j}^{-1}\right)\!\cdot\!\left(L_{-b_{j}}\varphi_{i}\right)\right]\right)\!\circ T_{j}^{-T}\right]\right\Vert _{L^{p}}\nonumber \\
\left({\scriptstyle \text{Lemma }\ref{lem:WienerTransformationFormula}}\right) & =2\sqrt{\dimension}\cdot\delta\cdot\left|\det T_{j}\right|^{-\frac{1}{p}}\!\cdot\!\left\Vert \left[v_{0}\cdot M_{T_{j}^{-T}\left[-2,2\right]^{\dimension}-T_{i}^{-T}\left[-1,1\right]^{\dimension}}\left(\Fourier^{-1}\left[\left(\widehat{\phi_{j}}\circ T_{j}^{-1}\right)\cdot\left(L_{-b_{j}}\varphi_{i}\right)\right]\right)\right]\circ T_{j}^{-T}\right\Vert _{L^{p}}\nonumber \\
 & =2\sqrt{\dimension}\cdot\delta\cdot\left\Vert \Fourier^{-1}\left[\left(\widehat{\phi_{j}}\circ T_{j}^{-1}\right)\cdot\left(L_{-b_{j}}\varphi_{i}\right)\right]\right\Vert _{W_{T_{j}^{-T}\left[-2,2\right]^{\dimension}-T_{i}^{-T}\left[-1,1\right]^{\dimension}}^{\dimension}\left(L_{v_{0}}^{p}\right)}\label{eq:MainOscillationEstimateQuasiBanachCalculationPart1}
\end{align}
Here, the step marked with $\left(\ast\right)$ is justified, since
$\Fourier^{-1}\left[\left(L_{-b_{j}}\varphi_{i}\right)\circ T_{j}\right]\in\Schwartz\left(\R^{\dimension}\right)$
and since $\gamma_{j}\in L^{1}\left(\R^{\dimension}\right)\cap C^{1}\left(\R^{\dimension}\right)$,
where $\phi_{j}=\nabla\gamma_{j}$ is bounded by Assumptions \ref{assu:MainAssumptions}
and \ref{assu:DiscreteBanachFrameAssumptions}.

Next, we observe
\begin{align*}
T_{j}^{-T}\left[-2,2\right]^{\dimension}-T_{i}^{-T}\left[-1,1\right]^{\dimension} & =T_{i}^{-T}\left(T_{i}^{T}T_{j}^{-T}\left[-2,2\right]^{\dimension}-\left[-1,1\right]^{\dimension}\right)\\
 & \subset T_{i}^{-T}\left(2\left\Vert \left(T_{j}^{-1}T_{i}\right)^{T}\right\Vert _{\ell^{\infty}\to\ell^{\infty}}\left[-1,1\right]^{\dimension}-\left[-1,1\right]^{\dimension}\right)\\
 & \subset T_{i}^{-T}\left[-\left(1+2\left\Vert T_{j}^{-1}T_{i}\right\Vert _{\ell^{1}\to\ell^{1}}\right),\,1+2\left\Vert T_{j}^{-1}T_{i}\right\Vert _{\ell^{1}\to\ell^{1}}\right]^{\dimension}.
\end{align*}
Consequently, if we set $R:=1+2\left\Vert T_{j}^{-1}T_{i}\right\Vert _{\ell^{1}\to\ell^{1}}$
for brevity, then Corollary \ref{cor:WienerLinearCubeNormEstimate}
(with $v_{0}$ instead of $v$, with $i=j$ and with $L=1$) yields
for arbitrary measurable $h:\R^{\dimension}\to\Compl^{k}$ the estimate
\begin{align*}
\left\Vert h\right\Vert _{W_{T_{j}^{-T}\left[-2,2\right]^{\dimension}-T_{i}^{-T}\left[-1,1\right]^{\dimension}}^{k}\left(L_{v_{0}}^{p}\right)} & \leq\left\Vert h\right\Vert _{W_{T_{i}^{-T}\left[-R,R\right]^{\dimension}}^{k}\left(L_{v_{0}}^{p}\right)}\\
 & \leq\Omega_{0}^{K}\Omega_{1}\cdot\left[3\dimension\left(1+1+2\left\Vert T_{j}^{-1}T_{i}\right\Vert _{\ell^{1}\to\ell^{1}}\right)\right]^{K+\frac{\dimension}{p}}\cdot\left(1+1\right)^{K+\frac{\dimension}{p}}\cdot\left\Vert h\right\Vert _{W_{T_{i}^{-T}\left[-1,1\right]^{\dimension}}^{k}\left(L_{v_{0}}^{p}\right)}\\
\left({\scriptstyle \text{since }\left\Vert A\right\Vert _{\ell^{1}\to\ell^{1}}\leq\sqrt{\dimension}\left\Vert A\right\Vert }\right) & \leq\Omega_{0}^{K}\Omega_{1}\cdot\left[12\cdot\dimension^{\frac{3}{2}}\left(1+\left\Vert T_{j}^{-1}T_{i}\right\Vert \right)\right]^{K+\frac{\dimension}{p}}\cdot\left\Vert h\right\Vert _{W_{T_{i}^{-T}\left[-1,1\right]^{\dimension}}^{k}\left(L_{v_{0}}^{p}\right)}.
\end{align*}

Now, we use this estimate and standard properties of the Fourier transform
to further estimate the right-hand side of equation (\ref{eq:MainOscillationEstimateQuasiBanachCalculationPart1})
as follows:
\begin{align}
\text{r.h.s.}\eqref{eq:MainOscillationEstimateQuasiBanachCalculationPart1} & =2\sqrt{\dimension}\cdot\delta\cdot\left\Vert \Fourier^{-1}\left(L_{-b_{j}}\left[\varphi_{i}\cdot L_{b_{j}}\left(\widehat{\phi_{j}}\circ T_{j}^{-1}\right)\right]\right)\right\Vert _{W_{T_{j}^{-T}\left[-2,2\right]^{\dimension}-T_{i}^{-T}\left[-1,1\right]^{\dimension}}^{\dimension}\left(L_{v_{0}}^{p}\right)}\nonumber \\
\left({\scriptstyle \widehat{\phi^{\left(j\right)}}=L_{b_{j}}\left(\widehat{\phi_{j}}\circ T_{j}^{-1}\right)}\right) & =2\sqrt{\dimension}\cdot\delta\cdot\left\Vert \Fourier^{-1}\left(L_{-b_{j}}\left[\varphi_{i}\cdot\widehat{\phi^{\left(j\right)}}\right]\right)\right\Vert _{W_{T_{j}^{-T}\left[-2,2\right]^{\dimension}-T_{i}^{-T}\left[-1,1\right]^{\dimension}}^{\dimension}\left(L_{v_{0}}^{p}\right)}\nonumber \\
\left({\scriptstyle \left|\Fourier^{-1}\left[L_{b}h\right]\right|=\left|\Fourier^{-1}h\right|}\right) & =2\sqrt{\dimension}\cdot\delta\cdot\left\Vert \Fourier^{-1}\left[\widehat{\phi^{\left(j\right)}}\cdot\varphi_{i}\right]\right\Vert _{W_{T_{j}^{-T}\left[-2,2\right]^{\dimension}-T_{i}^{-T}\left[-1,1\right]^{\dimension}}^{\dimension}\left(L_{v_{0}}^{p}\right)}\nonumber \\
 & \leq2\sqrt{\dimension}\cdot\Omega_{0}^{K}\Omega_{1}\cdot\left[12\cdot\dimension^{\frac{3}{2}}\left(1+\left\Vert T_{j}^{-1}T_{i}\right\Vert \right)\right]^{K+\frac{\dimension}{p}}\cdot\delta\cdot\left\Vert \Fourier^{-1}\left[\widehat{\phi^{\left(j\right)}}\cdot\varphi_{i}\right]\right\Vert _{W_{T_{i}^{-T}\left[-1,1\right]^{\dimension}}^{\dimension}\left(L_{v_{0}}^{p}\right)}.\label{eq:MainOscillationEstimateQuasiBanachCalculationPart2}
\end{align}
Recall that we are in the case $p\in\left(0,1\right)$. In particular,
we have $\left|y\right|\leq\left\Vert y\right\Vert _{\ell^{p}}$ for
each $y\in\R^{\dimension}$. For a vector-valued function $f:\R^{\dimension}\to\R^{k}$
and any (measurable) weight $u:\R^{\dimension}\to\left(0,\infty\right)$,
this implies 
\begin{align*}
\left\Vert f\right\Vert _{W_{Q}^{k}\left(L_{u}^{p}\right)}^{p} & =\int_{\R^{\dimension}}\left[u\left(x\right)\cdot\left|\left(M_{Q}f\right)\left(x\right)\right|\right]^{p}\d x\\
 & =\int_{\R^{\dimension}}\left[u\left(x\right)\right]^{p}\cdot\essup_{y\in x+Q}\left|f\left(y\right)\right|^{p}\d x\\
 & \leq\int_{\R^{\dimension}}\left[u\left(x\right)\right]^{p}\cdot\essup_{y\in x+Q}\left\Vert f\left(y\right)\right\Vert _{\ell^{p}}^{p}\d x\\
 & =\int_{\R^{\dimension}}\left[u\left(x\right)\right]^{p}\cdot\essup_{y\in x+Q}\sum_{\ell=1}^{k}\left|f_{\ell}\left(y\right)\right|^{p}\d x\\
 & \leq\sum_{\ell=1}^{k}\int_{\R^{\dimension}}\left[u\left(x\right)\right]^{p}\cdot\essup_{y\in x+Q}\left|f_{\ell}\left(y\right)\right|^{p}\d x\\
 & \leq k\cdot\max_{\ell\in\underline{k}}\left\Vert f_{\ell}\right\Vert _{W_{Q}\left(L_{u}^{p}\right)}^{p}.
\end{align*}
In other words, we have shown
\begin{equation}
\left\Vert f\right\Vert _{W_{Q}^{k}\left(L_{u}^{p}\right)}\leq k^{1/p}\cdot\max_{\ell\in\underline{k}}\left\Vert f_{\ell}\right\Vert _{W_{Q}\left(L_{u}^{p}\right)}.\label{eq:VectorValuedWienerEstimateQuasiBanach}
\end{equation}
Using this inequality (with $k=\dimension$), we can further estimate
the right-hand side of equation (\ref{eq:MainOscillationEstimateQuasiBanachCalculationPart2})
as follows:
\begin{align*}
\text{r.h.s.}\eqref{eq:MainOscillationEstimateQuasiBanachCalculationPart2} & \leq2\dimension^{\frac{1}{2}+\frac{1}{p}}\cdot\Omega_{0}^{K}\Omega_{1}\cdot\left[12\cdot\dimension^{\frac{3}{2}}\left(1+\left\Vert T_{j}^{-1}T_{i}\right\Vert \right)\right]^{K+\frac{\dimension}{p}}\cdot\delta\cdot\max_{\ell\in\underline{\dimension}}\left\Vert \Fourier^{-1}\left(\left[\vphantom{\phi^{\left(j\right)}}\smash{\widehat{\phi^{\left(j\right)}}}\right]_{\ell}\cdot\varphi_{i}\right)\right\Vert _{W_{T_{i}^{-T}\left[-1,1\right]^{\dimension}}\left(L_{v_{0}}^{p}\right)}\\
\left({\scriptstyle \text{Theorem }\ref{thm:BandlimitedWienerAmalgamSelfImproving}}\right) & \leq2\dimension^{\frac{1}{2}+\frac{1}{p}}\cdot\Omega_{0}^{K}\Omega_{1}\cdot\left[12\cdot\dimension^{\frac{3}{2}}\left(1+\left\Vert T_{j}^{-1}T_{i}\right\Vert \right)\right]^{K+\frac{\dimension}{p}}C_{3}\cdot\delta\cdot\max_{\ell\in\underline{\dimension}}\left\Vert \Fourier^{-1}\left(\left[\vphantom{\phi^{\left(j\right)}}\smash{\widehat{\phi^{\left(j\right)}}}\right]_{\ell}\cdot\varphi_{i}\right)\right\Vert _{L_{v_{0}}^{p}}\\
\left({\scriptstyle \text{since }\dimension\leq2^{\dimension}}\right) & \leq2\dimension^{\frac{1}{2}}2^{\frac{\dimension}{p}}\cdot\Omega_{0}^{K}\Omega_{1}\cdot\left[12\cdot\dimension^{\frac{3}{2}}\left(1+\left\Vert T_{j}^{-1}T_{i}\right\Vert \right)\right]^{K+\frac{\dimension}{p}}C_{3}\cdot\delta\cdot\left\Vert \Fourier^{-1}\left[\widehat{\phi^{\left(j\right)}}\cdot\varphi_{i}\right]\right\Vert _{L_{v_{0}}^{p}}\\
 & =:C_{4}\cdot\delta\cdot\left(1+\left\Vert T_{j}^{-1}T_{i}\right\Vert \right)^{K+\frac{\dimension}{p}}\cdot\left\Vert \Fourier^{-1}\left[\widehat{\phi^{\left(j\right)}}\cdot\varphi_{i}\right]\right\Vert _{L_{v_{0}}^{p}}.
\end{align*}
Here,
\[
C_{3}=2^{4\left(1+\frac{\dimension}{p}\right)}s_{\dimension}^{\frac{1}{p}}\left(192\cdot\dimension^{3/2}\cdot\left\lceil K+\frac{\dimension+1}{p}\right\rceil \right)^{\left\lceil K+\frac{\dimension+1}{p}\right\rceil +1}\cdot\Omega_{0}^{K}\Omega_{1}\cdot\left(1+R_{\CalQ}\right)^{\frac{\dimension}{p}},
\]
cf.\@ Theorem \ref{thm:BandlimitedWienerAmalgamSelfImproving}, since
we have for arbitrary $\ell\in\underline{\dimension}$ that
\[
\supp\left(\left[\vphantom{\phi^{\left(j\right)}}\smash{\widehat{\phi^{\left(j\right)}}}\right]_{\ell}\cdot\varphi_{i}\right)\subset\overline{Q_{i}}\subset T_{i}\left[\overline{B_{R_{\CalQ}}}\left(0\right)\right]+b_{i}\subset T_{i}\left[-R_{\CalQ},R_{\CalQ}\right]^{\dimension}+b_{i}.
\]

Putting everything together, we arrive at
\begin{align*}
 & \left\Vert {\rm osc}_{\delta\cdot T_{j}^{-T}\left[-1,1\right]^{\dimension}}\left[M_{-b_{j}}\Fourier^{-1}\left(\widehat{\gamma^{\left(j\right)}}\cdot\varphi_{i}\smash{\widehat{f}}\,\right)\right]\right\Vert _{W_{T_{j}^{-T}\left[-1,1\right]^{\dimension}}\left(L_{v}^{p}\right)}\\
 & \leq\left(2\sqrt{\dimension}\right)^{K}3^{\frac{\dimension}{p}}C_{2}\Omega_{0}^{K}\Omega_{1}\cdot\left|\det T_{i}\right|^{\frac{1}{p}-1}\cdot\left\Vert \Fourier^{-1}\!\left(\varphi_{i}^{\ast}\smash{\widehat{f}}\,\right)\right\Vert _{L_{v}^{p}}\\
 & \phantom{\leq}\cdot\left\Vert {\rm osc}_{\delta\cdot T_{j}^{-T}\left[-1,1\right]^{\dimension}}\left(\Fourier^{-1}\!\left[\left(\widehat{\gamma_{j}}\cdot\left[\left(L_{-b_{j}}\varphi_{i}\right)\circ T_{j}\right]\right)\!\circ\!T_{j}^{-1}\right]\right)\right\Vert _{W_{T_{j}^{-T}\left[-1,1\right]^{\dimension}-T_{i}^{-T}\left[-1,1\right]^{\dimension}}\left(L_{v_{0}}^{p}\right)}\\
 & \leq\left(2\sqrt{\dimension}\right)^{K}3^{\frac{\dimension}{p}}C_{2}C_{4}\cdot\Omega_{0}^{K}\Omega_{1}\cdot\left|\det T_{i}\right|^{\frac{1}{p}-1}\cdot\delta\cdot\left(1+\left\Vert T_{j}^{-1}T_{i}\right\Vert \right)^{K+\frac{\dimension}{p}}\cdot\left\Vert \Fourier^{-1}\left(\varphi_{i}^{\ast}\smash{\widehat{f}}\,\right)\right\Vert _{L_{v}^{p}}\cdot\left\Vert \Fourier^{-1}\left[\widehat{\phi^{\left(j\right)}}\cdot\varphi_{i}\right]\right\Vert _{L_{v_{0}}^{p}}.
\end{align*}
This establishes the claim for $p\in\left(0,1\right)$, since we have
$C_{\CalQ}\geq\left\Vert T_{i}^{-1}T_{i}\right\Vert \geq1$ and $s_{\dimension}\leq2^{2\dimension}$
and hence
\begin{align*}
 & \left(2\sqrt{\dimension}\right)^{K}3^{\frac{\dimension}{p}}C_{2}C_{4}\cdot\Omega_{0}^{K}\Omega_{1}\\
 & =C_{1}\cdot2^{5}\dimension^{1/2}\cdot\left(2\sqrt{\dimension}\right)^{K}96^{\frac{\dimension}{p}}\cdot\left(216\cdot\dimension^{\frac{5}{2}}\right)^{K+\frac{\dimension}{p}}\cdot s_{\dimension}^{\frac{1}{p}}\left(192\cdot\dimension^{3/2}\cdot\left\lceil K\!+\!\frac{\dimension+1}{p}\right\rceil \right)^{\left\lceil K+\frac{\dimension+1}{p}\right\rceil +1}\cdot\Omega_{0}^{4K}\Omega_{1}^{4}\cdot\left(1\!+\!R_{\CalQ}\right)^{\frac{\dimension}{p}}\\
 & \leq2^{9}\dimension^{\frac{1}{2}}\cdot2^{15\frac{\dimension}{p}}\left(2\sqrt{\dimension}\right)^{-\frac{\dimension}{p}}\!\cdot\left(432\cdot\dimension^{3}\right)^{K+\frac{\dimension}{p}}\!\cdot\!\left(192\cdot\dimension^{3/2}\cdot\left\lceil K\!+\!\frac{\dimension+1}{p}\right\rceil \right)^{2\left\lceil K+\frac{\dimension+1}{p}\right\rceil +2}\cdot\Omega_{0}^{5K}\Omega_{1}^{5}\cdot\left(1\!+\!R_{\CalQ}\right)^{\frac{\dimension}{p}}\left(1\!+\!3C_{\CalQ}R_{\CalQ}\right)^{\frac{\dimension}{p}}\\
 & \leq2^{9}\dimension^{1/2}\cdot2^{17\frac{\dimension}{p}}\left(2\sqrt{\dimension}\right)^{-\frac{\dimension}{p}}\cdot\left(21\cdot\dimension^{3/2}\right)^{-4}\cdot\left(4032\cdot\dimension^{3}\cdot\left\lceil K\!+\!\frac{\dimension+1}{p}\right\rceil \right)^{2\left\lceil K+\frac{\dimension+1}{p}\right\rceil +2}\cdot\Omega_{0}^{5K}\Omega_{1}^{5}\cdot\left(1\!+\!C_{\CalQ}R_{\CalQ}\right)^{\frac{2\dimension}{p}}\\
 & \leq\frac{2^{16\frac{\dimension}{p}}\cdot\left(1+C_{\CalQ}R_{\CalQ}\right)^{\frac{2\dimension}{p}}}{370\cdot\dimension^{11/2}\cdot\dimension^{\dimension/2p}}\cdot\left(4032\cdot\dimension^{3}\cdot\left\lceil K+\frac{\dimension+1}{p}\right\rceil \right)^{2\left\lceil K+\frac{\dimension+1}{p}\right\rceil +2}\cdot\Omega_{0}^{5K}\Omega_{1}^{5}.
\end{align*}

\medskip{}

For $p\in\left[1,\infty\right]$, the proof is simpler: We use the
weighted Young inequality (equation (\ref{eq:WeightedYoungInequality}))
and equation (\ref{eq:MainOscillationEstimateOscillationToTheLeft})
to derive
\begin{align*}
 & \left\Vert {\rm osc}_{\delta\cdot T_{j}^{-T}\left[-1,1\right]^{\dimension}}\left(M_{-b_{j}}\left[\Fourier^{-1}\left(\widehat{\gamma^{\left(j\right)}}\cdot\varphi_{i}\widehat{f}\right)\right]\right)\right\Vert _{L_{v}^{p}}\\
\left({\scriptstyle \text{eqs. }\eqref{eq:MainOscillationEstimateOscillationToTheLeft},\eqref{eq:WeightedYoungInequality}}\right) & \leq\left\Vert \osc{\delta\cdot T_{j}^{-T}\left[-1,1\right]^{\dimension}}\Fourier^{-1}\left[L_{-b_{j}}\left(\widehat{\gamma^{\left(j\right)}}\cdot\varphi_{i}\right)\right]\right\Vert _{L_{v_{0}}^{1}}\cdot\left\Vert \Fourier^{-1}\left[L_{-b_{j}}\left(\varphi_{i}^{\ast}\smash{\widehat{f}}\,\right)\right]\right\Vert _{L_{v}^{p}}\\
\left({\scriptstyle \widehat{\gamma^{\left(j\right)}}=L_{b_{j}}\left(\widehat{\gamma_{j}}\circ T_{j}^{-1}\right)}\right) & =\left\Vert \osc{\delta\cdot T_{j}^{-T}\left[-1,1\right]^{\dimension}}\Fourier^{-1}\left[\left(\widehat{\gamma_{j}}\circ T_{j}^{-1}\right)\cdot\left(L_{-b_{j}}\varphi_{i}\right)\right]\right\Vert _{L_{v_{0}}^{1}}\cdot\left\Vert \Fourier^{-1}\left[L_{-b_{j}}\left(\varphi_{i}^{\ast}\smash{\widehat{f}}\,\right)\right]\right\Vert _{L_{v}^{p}}\\
\left({\scriptstyle \left|\Fourier^{-1}\left[L_{b}h\right]\right|=\left|\Fourier^{-1}h\right|}\right) & =\left\Vert \osc{\delta\cdot T_{j}^{-T}\left[-1,1\right]^{\dimension}}\Fourier^{-1}\left(\left[\widehat{\gamma_{j}}\cdot\left(\left[L_{-b_{j}}\varphi_{i}\right]\circ T_{j}\right)\right]\circ T_{j}^{-1}\right)\right\Vert _{L_{v_{0}}^{1}}\cdot\left\Vert \Fourier^{-1}\left(\varphi_{i}^{\ast}\smash{\widehat{f}}\,\right)\right\Vert _{L_{v}^{p}}\\
 & =\left|\det T_{j}\right|\cdot\left\Vert \osc{\delta\cdot T_{j}^{-T}\left[-1,1\right]^{\dimension}}\left[\left(\Fourier^{-1}\left[\widehat{\gamma_{j}}\cdot\left(\left[L_{-b_{j}}\varphi_{i}\right]\circ T_{j}\right)\right]\right)\circ T_{j}^{T}\right]\right\Vert _{L_{v_{0}}^{1}}\cdot\left\Vert \Fourier^{-1}\left(\varphi_{i}^{\ast}\smash{\widehat{f}}\,\right)\right\Vert _{L_{v}^{p}}\\
\left({\scriptstyle \text{Lemma }\ref{lem:OscillationLinearChange}}\right) & =\left|\det T_{j}\right|\cdot\left\Vert \left(\osc{\delta\cdot\left[-1,1\right]^{\dimension}}\Fourier^{-1}\left[\widehat{\gamma_{j}}\cdot\left(\left[L_{-b_{j}}\varphi_{i}\right]\circ T_{j}\right)\right]\right)\circ T_{j}^{T}\right\Vert _{L_{v_{0}}^{1}}\cdot\left\Vert \Fourier^{-1}\left(\varphi_{i}^{\ast}\smash{\widehat{f}}\,\right)\right\Vert _{L_{v}^{p}}\\
 & =\left\Vert \left(v_{0}\circ T_{j}^{-T}\right)\cdot\osc{\delta\cdot\left[-1,1\right]^{\dimension}}\left[\gamma_{j}\ast\Fourier^{-1}\left(\left[L_{-b_{j}}\varphi_{i}\right]\circ T_{j}\right)\right]\right\Vert _{L^{1}}\cdot\left\Vert \Fourier^{-1}\left(\varphi_{i}^{\ast}\smash{\widehat{f}}\,\right)\right\Vert _{L_{v}^{p}}\\
\left({\scriptstyle \text{Lemma }\ref{lem:OscillationEstimatedByWienerDerivative}}\right) & \leq2\delta\sqrt{\dimension}\cdot\left\Vert \left(v_{0}\circ T_{j}^{-T}\right)\cdot M_{\delta\left[-1,1\right]^{\dimension}}\left(\nabla\left[\gamma_{j}\ast\Fourier^{-1}\left(\left[L_{-b_{j}}\varphi_{i}\right]\circ T_{j}\right)\right]\right)\right\Vert _{L^{1}}\cdot\left\Vert \Fourier^{-1}\left(\varphi_{i}^{\ast}\smash{\widehat{f}}\,\right)\right\Vert _{L_{v}^{p}}\\
\left({\scriptstyle \text{since }\delta\leq1}\right) & \leq2\delta\sqrt{\dimension}\cdot\left\Vert \left(v_{0}\circ T_{j}^{-T}\right)\cdot M_{\left[-1,1\right]^{\dimension}}\left(\nabla\left[\gamma_{j}\ast\Fourier^{-1}\left(\left[L_{-b_{j}}\varphi_{i}\right]\circ T_{j}\right)\right]\right)\right\Vert _{L^{1}}\cdot\left\Vert \Fourier^{-1}\left(\varphi_{i}^{\ast}\smash{\widehat{f}}\,\right)\right\Vert _{L_{v}^{p}}\\
\left({\scriptstyle \nabla\left(\gamma\ast h\right)=\left(\nabla\gamma\right)\ast h}\right) & =2\delta\sqrt{\dimension}\cdot\left\Vert \left(v_{0}\circ T_{j}^{-T}\right)\cdot M_{\left[-1,1\right]^{\dimension}}\left[\left(\nabla\gamma_{j}\right)\ast\left(\Fourier^{-1}\left[\left(L_{-b_{j}}\varphi_{i}\right)\circ T_{j}\right]\right)\right]\right\Vert _{L^{1}}\cdot\left\Vert \Fourier^{-1}\left(\varphi_{i}^{\ast}\smash{\widehat{f}}\,\right)\right\Vert _{L_{v}^{p}}.
\end{align*}
Here, the last step is justified just as for $p\in\left(0,1\right)$.
Now, we recall $\phi_{j}=\nabla\gamma_{j}$ and continue our estimate:
\begin{align*}
\dots & =2\delta\sqrt{\dimension}\cdot\left\Vert \left(v_{0}\circ T_{j}^{-T}\right)\cdot M_{\left[-1,1\right]^{\dimension}}\left[\Fourier^{-1}\left(\widehat{\phi_{j}}\cdot\left[\left(L_{-b_{j}}\varphi_{i}\right)\circ T_{j}\right]\right)\right]\right\Vert _{L^{1}}\cdot\left\Vert \Fourier^{-1}\left(\varphi_{i}^{\ast}\smash{\widehat{f}}\,\right)\right\Vert _{L_{v}^{p}}\\
 & =2\delta\sqrt{\dimension}\cdot\left\Vert \left(v_{0}\circ T_{j}^{-T}\right)\cdot M_{\left[-1,1\right]^{\dimension}}\left[\Fourier^{-1}\left(\left[\left(\widehat{\phi_{j}}\circ T_{j}^{-1}\right)\cdot\left(L_{-b_{j}}\varphi_{i}\right)\right]\circ T_{j}\right)\right]\right\Vert _{L^{1}}\cdot\left\Vert \Fourier^{-1}\left(\varphi_{i}^{\ast}\smash{\widehat{f}}\,\right)\right\Vert _{L_{v}^{p}}\\
 & =2\delta\sqrt{\dimension}\cdot\left|\det T_{j}\right|^{-1}\left\Vert \left(v_{0}\!\circ\!T_{j}^{-T}\right)\!\cdot\!M_{\left[-1,1\right]^{\dimension}}\!\left[\!\left(\Fourier^{-1}\!\left[L_{-b_{j}}\!\left[\varphi_{i}\cdot L_{b_{j}}\!\left(\widehat{\phi_{j}}\circ T_{j}^{-1}\right)\!\right]\right]\right)\!\circ\!T_{j}^{-T}\right]\!\right\Vert _{L^{1}}\left\Vert \Fourier^{-1}\left(\varphi_{i}^{\ast}\smash{\widehat{f}}\,\right)\right\Vert _{L_{v}^{p}}\\
\left({\scriptstyle \text{Lem. }\ref{lem:WienerTransformationFormula}}\right) & =2\delta\sqrt{\dimension}\cdot\left|\det T_{j}\right|^{-1}\left\Vert \left[v_{0}\!\cdot\!M_{T_{j}^{-T}\left[-1,1\right]^{\dimension}}\left(\Fourier^{-1}\!\left[L_{-b_{j}}\left[\varphi_{i}\cdot L_{b_{j}}\!\left(\widehat{\phi_{j}}\circ T_{j}^{-1}\right)\right]\right]\right)\right]\circ T_{j}^{-T}\right\Vert _{L^{1}}\left\Vert \Fourier^{-1}\left(\varphi_{i}^{\ast}\smash{\widehat{f}}\,\right)\right\Vert _{L_{v}^{p}}\\
 & =2\delta\sqrt{\dimension}\cdot\left\Vert \Fourier^{-1}\!\left(L_{-b_{j}}\left[\varphi_{i}\cdot\widehat{\phi^{\left(j\right)}}\right]\right)\right\Vert _{W_{T_{j}^{-T}\left[-1,1\right]^{\dimension}}^{\dimension}\left(L_{v_{0}}^{1}\right)}\left\Vert \Fourier^{-1}\left(\varphi_{i}^{\ast}\smash{\widehat{f}}\,\right)\right\Vert _{L_{v}^{p}}\\
 & =2\delta\sqrt{\dimension}\cdot\left\Vert \Fourier^{-1}\left[\varphi_{i}\cdot\widehat{\phi^{\left(j\right)}}\right]\right\Vert _{W_{T_{j}^{-T}\left[-1,1\right]^{\dimension}}^{\dimension}\left(L_{v_{0}}^{1}\right)}\cdot\left\Vert \Fourier^{-1}\left(\varphi_{i}^{\ast}\smash{\widehat{f}}\,\right)\right\Vert _{L_{v}^{p}}.
\end{align*}
Here, the last step used that $\left|\Fourier^{-1}\left[L_{b}h\right]\right|=\left|\Fourier^{-1}h\right|$.

Now, we need an analog of equation (\ref{eq:VectorValuedWienerEstimateQuasiBanach})
for the case $p\in\left[1,\infty\right]$. But for an arbitrary (measurable)
weight $u:\R^{\dimension}\to\left(0,\infty\right)$ and any $q\in\left[1,\infty\right]$,
the solidity of $W_{Q}\left(L_{u}^{q}\right)$ and the triangle inequality
for the associated norm yield for any measurable vector-valued function
$f=\left(f_{1},\dots,f_{k}\right):\R^{\dimension}\to\Compl^{k}$ that
\begin{equation}
\left\Vert f\right\Vert _{W_{Q}^{k}\left(L_{u}^{q}\right)}=\left\Vert \left|f\right|\right\Vert _{W_{Q}\left(L_{u}^{q}\right)}\leq\left\Vert \sum_{\ell=1}^{k}\left|f_{\ell}\right|\right\Vert _{W_{Q}\left(L_{u}^{q}\right)}\leq\sum_{\ell=1}^{k}\left\Vert f_{\ell}\right\Vert _{W_{Q}\left(L_{u}^{q}\right)}\leq k\cdot\max_{\ell\in\underline{k}}\left\Vert f_{\ell}\right\Vert _{W_{Q}\left(L_{u}^{q}\right)}.\label{eq:VectorValuedWienerEstimateBanach}
\end{equation}
We now use this estimate (with $k=\dimension$), as well as  equation
(\ref{eq:WienerLinearCubeTransformationChange}) and Theorem \ref{thm:BandlimitedWienerAmalgamSelfImproving}
(both with $v_{0}$ instead of $v$) to conclude
\begin{align*}
\left\Vert \Fourier^{-1}\left[\varphi_{i}\cdot\widehat{\phi^{\left(j\right)}}\right]\right\Vert _{W_{T_{j}^{-T}\cdot\left[-1,1\right]^{\dimension}}^{\dimension}\left(L_{v_{0}}^{1}\right)} & \leq\Omega_{0}^{K}\Omega_{1}\cdot\left(6\dimension\right)^{K+\dimension}\!\cdot\!\left(1\!+\!\left\Vert T_{j}^{-1}T_{i}\right\Vert \right)^{\!K+\dimension}\cdot\left\Vert \Fourier^{-1}\left[\varphi_{i}\cdot\widehat{\phi^{\left(j\right)}}\right]\right\Vert _{W_{T_{i}^{-T}\left[-1,1\right]^{\dimension}}^{\dimension}\left(L_{v_{0}}^{1}\right)}\\
 & \leq\dimension\cdot\Omega_{0}^{K}\Omega_{1}\cdot\left(6\dimension\right)^{K+\dimension}\!\cdot\!\left(1\!+\!\left\Vert T_{j}^{-1}T_{i}\right\Vert \right)^{\!K+\dimension}\cdot\max_{\ell\in\underline{\dimension}}\left\Vert \Fourier^{-1}\left[\varphi_{i}\cdot\left(\widehat{\phi^{\left(j\right)}}\right)_{\ell}\right]\right\Vert _{W_{T_{i}^{-T}\left[-1,1\right]^{\dimension}}\left(L_{v_{0}}^{1}\right)}\\
\left({\scriptstyle \text{Thm. }\ref{thm:BandlimitedWienerAmalgamSelfImproving}}\right) & \leq C_{5}\cdot\dimension\cdot\Omega_{0}^{K}\Omega_{1}\cdot\left(6\dimension\right)^{K+\dimension}\!\cdot\!\left(1\!+\!\left\Vert T_{j}^{-1}T_{i}\right\Vert \right)^{\!K+\dimension}\cdot\left\Vert \Fourier^{-1}\left[\varphi_{i}\cdot\widehat{\phi^{\left(j\right)}}\right]\right\Vert _{L_{v_{0}}^{1}}.
\end{align*}
Here, Theorem \ref{thm:BandlimitedWienerAmalgamSelfImproving} is
applicable, since we have $\supp\left(\varphi_{i}\cdot\left(\widehat{\phi^{\left(j\right)}}\right)_{\ell}\right)\subset\overline{Q_{i}}\subset T_{i}\left[-R_{\CalQ},R_{\CalQ}\right]^{\dimension}+b_{i}$.
Hence, that theorem justifies the last step in the estimate above,
with constant
\[
C_{5}:=2^{4\left(1+\dimension\right)}s_{\dimension}\left(192\cdot\dimension^{3/2}\cdot\left\lceil K+\dimension+1\right\rceil \right)^{\left\lceil K+\dimension+1\right\rceil +1}\cdot\Omega_{0}^{K}\Omega_{1}\cdot\left(1+R_{\CalQ}\right)^{\dimension}.
\]
It is not hard to see that this implies the claim for $p\in\left[1,\infty\right]$.
\end{proof}
Next, we show that the map ${\rm Ana}_{\Gamma}$ considered in Theorem
\ref{thm:ConvolvingDecompositionSpaceWithGammaJ} is not merely bounded
as a map into $\ell_{w}^{q}\left(\left[V_{j}\right]_{j\in I}\right)$,
but even as a map into the smaller space $\ell_{w}^{q}\left(\left[W_{j}\right]_{j\in I}\right)$.
In particular, this establishes continuity of $\gamma^{\left(j\right)}\ast f$
for every $j\in I$ and arbitrary $f\in\DecompSp{\CalQ}p{\ell_{w}^{q}}v$.
\begin{lem}
\label{lem:OscillationForFree}Let $p,q\in\left(0,\infty\right]$
and assume that $\Gamma=\left(\gamma_{i}\right)_{i\in I}$ fulfills
Assumption \ref{assu:DiscreteBanachFrameAssumptions}.

Then, the map
\[
{\rm Ana}_{{\rm osc}}:\DecompSp{\CalQ}p{\ell_{w}^{q}}v\to\ell_{w}^{q}\left(\left[W_{j}\right]_{j\in I}\right),f\mapsto\left(\gamma^{\left(j\right)}\ast f\right)_{j\in I}
\]
is well-defined and bounded, with
\[
\vertiii{{\rm Ana}_{{\rm osc}}}\leq C\cdot2^{\max\left\{ 0,\frac{1}{q}-1\right\} }\vertiii{\smash{\Gamma_{\CalQ}}}\cdot\left(\vertiii{\smash{\overrightarrow{A}}}^{\max\left\{ 1,\frac{1}{p}\right\} }+\vertiii{\smash{\overrightarrow{B}}}^{\max\left\{ 1,\frac{1}{p}\right\} }\right),
\]
where $\Gamma_{\CalQ}:\ell_{w}^{q}\left(I\right)\to\ell_{w}^{q}\left(I\right),c\mapsto c^{\ast}$
denotes the $\CalQ$-clustering map, i.e., $c_{i}^{\ast}=\sum_{\ell\in i^{\ast}}c_{\ell}$
and where
\[
C:=\begin{cases}
N_{\CalQ}^{\frac{1}{p}-1}\cdot\frac{2^{16\frac{\dimension}{p}}\cdot\left(1+C_{\CalQ}R_{\CalQ}\right)^{\frac{2\dimension}{p}}}{370\cdot\dimension^{11/2}\cdot\dimension^{\dimension/2p}}\cdot\left(4032\cdot\dimension^{3}\cdot\left\lceil K+\frac{\dimension+1}{p}\right\rceil \right)^{2\left\lceil K+\frac{\dimension+1}{p}\right\rceil +2}\cdot\Omega_{0}^{5K}\Omega_{1}^{5}, & \text{if }p\in\left(0,1\right),\\
\frac{2^{6\dimension}}{\sqrt{\dimension}}\cdot\left(1152\cdot\dimension^{5/2}\cdot\left\lceil K+\dimension+1\right\rceil \right)^{\left\lceil K\right\rceil +\dimension+2}\cdot\Omega_{0}^{2K}\Omega_{1}^{2}\cdot\left(1+R_{\CalQ}\right)^{\dimension}, & \text{if }p\in\left[1,\infty\right].
\end{cases}
\]

Furthermore, we have
\begin{equation}
\left(\gamma^{\left(j\right)}\ast f\right)\left(x\right)=\sum_{i\in I}\Fourier^{-1}\left(\widehat{\gamma^{\left(j\right)}}\cdot\varphi_{i}\cdot\widehat{f}\right)\left(x\right)\qquad\forall x\in\R^{\dimension}\qquad\forall f\in\DecompSp{\CalQ}p{\ell_{w}^{q}}v,\label{eq:SpecialConvolutionPointwiseDefinition}
\end{equation}
with \emph{locally} uniform convergence of the series.
\end{lem}
\begin{proof}
Recall from Theorem \ref{thm:ConvolvingDecompositionSpaceWithGammaJ}
and from the ensuing remark (which contains the definition of $\gamma^{\left(j\right)}\ast f$)
that
\begin{equation}
\left(\gamma^{\left(j\right)}\ast f\right)\left(x\right)=\sum_{i\in I}\Fourier^{-1}\left(\widehat{\gamma^{\left(j\right)}}\cdot\varphi_{i}\cdot\widehat{f}\right)\left(x\right)\qquad\forall f\in\DecompSp{\CalQ}p{\ell_{w}^{q}}v,\label{eq:OscillationForFreeSeriesRepresentation}
\end{equation}
where we already know that the series converges absolutely almost
everywhere. Next, note that each of the summands of the series above
is a smooth function; this follows from the Paley-Wiener theorem,
since $\widehat{\gamma^{\left(j\right)}}\cdot\varphi_{i}\cdot\widehat{f}$
is a (tempered) distribution with compact support. Thus, to prove
continuity of $\gamma^{\left(j\right)}\ast f$, it suffices to show
that the series actually converges \emph{locally} uniformly; by continuity
of the summands, for this it suffices to have convergence in $L_{\left(1+\left|\mybullet\right|\right)^{-K}}^{\infty}\left(\R^{\dimension}\right)$.
We will prove this convergence in $L_{\left(1+\left|\mybullet\right|\right)^{-K}}^{\infty}\left(\R^{\dimension}\right)$
simultaneously with the boundedness of ${\rm Ana}_{{\rm osc}}$.

\medskip{}

Let us first consider the case $p\in\left[1,\infty\right]$. Here,
we let $C_{1}>0$ be the constant provided by Lemma \ref{lem:MainOscillationEstimate}
(for $p\in\left[1,\infty\right]$), so that we get for arbitrary
$0<\delta\leq1$ the estimate
\begin{align}
 & \frac{1}{\delta}\sum_{i\in I}\left\Vert \osc{\delta\cdot T_{j}^{-T}\left[-1,1\right]^{\dimension}}\left(M_{-b_{j}}\left[\Fourier^{-1}\left(\widehat{\gamma^{\left(j\right)}}\cdot\varphi_{i}\cdot\widehat{f}\right)\right]\right)\right\Vert _{L_{v}^{p}}\nonumber \\
\left({\scriptstyle \text{Lemma }\ref{lem:MainOscillationEstimate}}\right) & \leq C_{1}\cdot\sum_{i\in I}\left[\left(1+\left\Vert T_{j}^{-1}T_{i}\right\Vert \right)^{K+\dimension}\cdot\left\Vert \Fourier^{-1}\left[\varphi_{i}\cdot\widehat{\phi^{\left(j\right)}}\right]\right\Vert _{L_{v_{0}}^{1}}\cdot\left\Vert \Fourier^{-1}\left(\varphi_{i}^{\ast}\smash{\widehat{f}}\,\right)\right\Vert _{L_{v}^{p}}\right]\nonumber \\
 & =C_{1}\cdot\sum_{i\in I}\left[B_{j,i}\cdot c_{i}\right]=C_{1}\cdot\left(\overrightarrow{B}c\right)_{j},\label{eq:OscillationSumEstimateBanachCase}
\end{align}
where we defined $c_{i}:=\left\Vert \Fourier^{-1}\left(\varphi_{i}^{\ast}\smash{\widehat{f}}\,\right)\right\Vert _{L_{v}^{p}}$
for all $i\in I$.

Setting $d_{i}:=\left\Vert \Fourier^{-1}\left(\varphi_{i}\smash{\widehat{f}}\,\right)\right\Vert _{L_{v}^{p}}$
for $i\in I$ and using the triangle inequality for $L_{v}^{p}$,
we get $c_{i}\leq\left(\Gamma_{\CalQ}\,d\right)_{i}$ for $i\in I$.
By solidity of $\ell_{w}^{q}\left(I\right)$, this allows us to conclude
\begin{align}
C_{1}\cdot\left\Vert \overrightarrow{B}c\right\Vert _{\ell_{w}^{q}} & \leq C_{1}\cdot\vertiii{\smash{\overrightarrow{B}}}\cdot\left\Vert c\right\Vert _{\ell_{w}^{q}}\nonumber \\
 & \leq C_{1}\cdot\vertiii{\smash{\overrightarrow{B}}}\cdot\left\Vert \Gamma_{\CalQ}\,d\right\Vert _{\ell_{w}^{q}}\nonumber \\
 & \leq C_{1}\cdot\vertiii{\smash{\Gamma_{\CalQ}}}\cdot\vertiii{\smash{\overrightarrow{B}}}\cdot\left\Vert d\right\Vert _{\ell_{w}^{q}}\nonumber \\
 & =C_{1}\cdot\vertiii{\smash{\Gamma_{\CalQ}}}\cdot\vertiii{\smash{\overrightarrow{B}}}\cdot\left\Vert f\right\Vert _{\DecompSp{\CalQ}p{\ell_{w}^{q}}v}<\infty.\label{eq:OscillationForFreeSequenceNormEstimateBanachCase}
\end{align}

In particular, we get $\left(\smash{\overrightarrow{B}}c\right)_{j}<\infty$
for all $j\in I$, so that the right-hand side of equation (\ref{eq:OscillationSumEstimateBanachCase})
is finite. We now use this estimate for $\delta=1$: For arbitrary
$x\in\R^{\dimension}$ and $a\in T_{j}^{-T}\left[-1,1\right]^{\dimension}$,
we have $x,\,x+a\in x+T_{j}^{-T}\left[-1,1\right]^{\dimension}$ and
hence
\begin{align*}
 & \left|\left[\Fourier^{-1}\left(\widehat{\gamma^{\left(j\right)}}\cdot\varphi_{i}\cdot\widehat{f}\right)\right]\left(x+a\right)\right|\\
 & \leq\left|\left(M_{-b_{j}}\left[\Fourier^{-1}\left(\widehat{\gamma^{\left(j\right)}}\cdot\varphi_{i}\cdot\widehat{f}\right)\right]\right)\left(x+a\right)-\left(M_{-b_{j}}\left[\Fourier^{-1}\left(\widehat{\gamma^{\left(j\right)}}\cdot\varphi_{i}\cdot\widehat{f}\right)\right]\right)\left(x\right)\right|+\left|\left(M_{-b_{j}}\left[\Fourier^{-1}\left(\widehat{\gamma^{\left(j\right)}}\cdot\varphi_{i}\cdot\widehat{f}\right)\right]\right)\left(x\right)\right|\\
 & \leq\left[\osc{T_{j}^{-T}\left[-1,1\right]^{\dimension}}\left(M_{-b_{j}}\left[\Fourier^{-1}\left(\widehat{\gamma^{\left(j\right)}}\cdot\varphi_{i}\cdot\widehat{f}\right)\right]\right)\right]\left(x\right)+\left|\left[\Fourier^{-1}\left(\widehat{\gamma^{\left(j\right)}}\cdot\varphi_{i}\cdot\widehat{f}\right)\right]\left(x\right)\right|,
\end{align*}
which yields
\[
M_{T_{j}^{-T}\left[-1,1\right]^{\dimension}}\left[\Fourier^{-1}\left(\widehat{\gamma^{\left(j\right)}}\cdot\varphi_{i}\cdot\widehat{f}\right)\right]\leq\osc{T_{j}^{-T}\left[-1,1\right]^{\dimension}}\left(M_{-b_{j}}\left[\Fourier^{-1}\left(\widehat{\gamma^{\left(j\right)}}\cdot\varphi_{i}\cdot\widehat{f}\right)\right]\right)+\left|\Fourier^{-1}\left(\widehat{\gamma^{\left(j\right)}}\cdot\varphi_{i}\cdot\widehat{f}\right)\right|.
\]
Using the triangle inequality for $L_{v}^{p}\left(\R^{\dimension}\right)$
and solidity of $L_{v}^{p}\left(\R^{\dimension}\right)$, this yields
\begin{align*}
 & \sum_{i\in I}\left\Vert \Fourier^{-1}\left(\widehat{\gamma^{\left(j\right)}}\cdot\varphi_{i}\cdot\widehat{f}\right)\right\Vert _{W_{T_{j}^{-T}\left[-1,1\right]^{\dimension}}\left(L_{v}^{p}\right)}\\
 & =\sum_{i\in I}\left\Vert M_{T_{j}^{-T}\left[-1,1\right]^{\dimension}}\left[\Fourier^{-1}\left(\widehat{\gamma^{\left(j\right)}}\cdot\varphi_{i}\cdot\widehat{f}\right)\right]\right\Vert _{L_{v}^{p}}\\
 & \leq\sum_{i\in I}\left(\left\Vert \Fourier^{-1}\left(\widehat{\gamma^{\left(j\right)}}\cdot\varphi_{i}\cdot\widehat{f}\right)\right\Vert _{L_{v}^{p}}+\left\Vert \osc{T_{j}^{-T}\left[-1,1\right]^{\dimension}}\left(M_{-b_{j}}\left[\Fourier^{-1}\left(\widehat{\gamma^{\left(j\right)}}\cdot\varphi_{i}\cdot\widehat{f}\right)\right]\right)\right\Vert _{L_{v}^{p}}\right)<\infty.
\end{align*}
Here, finiteness of the right-hand side follows from equation (\ref{eq:OscillationSumEstimateBanachCase})
(with $\delta=1$) and from Theorem \ref{thm:ConvolvingDecompositionSpaceWithGammaJ},
where we saw that the series $\sum_{i\in I}\Fourier^{-1}\left(\widehat{\gamma^{\left(j\right)}}\cdot\varphi_{i}\cdot\widehat{f}\right)$
converges normally in $L_{v}^{p}$.

But it follows from equation (\ref{eq:WeightedWienerAmalgamTemperedDistribution})
that $W_{T_{j}^{-T}\left[-1,1\right]^{\dimension}}\left(L_{v}^{p}\right)\hookrightarrow L_{\left(1+\left|\mybullet\right|\right)^{-K}}^{\infty}\left(\R^{\dimension}\right)$,
where the norm of the embedding might heavily depend on $j$. Setting
$\left\Vert h\right\Vert _{\ast}:=\sup_{x\in\R^{\dimension}}\left(1+\left|x\right|\right)^{-K}\left|h\left(x\right)\right|$,
this allows us to conclude by continuity that
\[
\sum_{i\in I}\left\Vert \Fourier^{-1}\left(\widehat{\gamma^{\left(j\right)}}\cdot\varphi_{i}\cdot\widehat{f}\right)\right\Vert _{\ast}=\sum_{i\in I}\left\Vert \Fourier^{-1}\left(\widehat{\gamma^{\left(j\right)}}\cdot\varphi_{i}\cdot\widehat{f}\right)\right\Vert _{L_{\left(1+\left|\mybullet\right|\right)^{-K}}^{\infty}}\lesssim_{j}\,\,\sum_{i\in I}\left\Vert \Fourier^{-1}\left(\widehat{\gamma^{\left(j\right)}}\cdot\varphi_{i}\cdot\widehat{f}\right)\right\Vert _{W_{T_{j}^{-T}\left[-1,1\right]^{\dimension}}\left(L_{v}^{p}\right)}<\infty,
\]
so that the series in equation (\ref{eq:OscillationForFreeSeriesRepresentation})
indeed converges \emph{locally} uniformly. Hence, $\gamma^{\left(j\right)}\ast f$
is continuous for every $j\in I$ and arbitrary $f\in\DecompSp{\CalQ}p{\ell_{w}^{q}}v$.

\medskip{}

Now, it is not hard to see $\osc Q\left(\sum_{i\in I}f_{i}\right)\leq\sum_{i\in I}\left(\osc Qf_{i}\right)$
for each pointwise convergent series $\sum_{i\in I}f_{i}$. Hence,
equation (\ref{eq:OscillationForFreeSeriesRepresentation}) and the
triangle inequality for $L_{v}^{p}\left(\R^{\dimension}\right)$ imply
\begin{align*}
\frac{1}{\delta}\left\Vert \osc{\delta\cdot T_{j}^{-T}\left[-1,1\right]^{\dimension}}\left[M_{-b_{j}}\left(\gamma^{\left(j\right)}\ast f\right)\right]\right\Vert _{L_{v}^{p}} & \leq\frac{1}{\delta}\sum_{i\in I}\left\Vert \osc{\delta\cdot T_{j}^{-T}\left[-1,1\right]^{\dimension}}\left(M_{-b_{j}}\left[\Fourier^{-1}\left(\widehat{\gamma^{\left(j\right)}}\cdot\varphi_{i}\cdot\widehat{f}\right)\right]\right)\right\Vert _{L_{v}^{p}}\\
\left({\scriptstyle \text{eq. }\eqref{eq:OscillationSumEstimateBanachCase}}\right) & \leq C_{1}\cdot\left(\overrightarrow{B}c\right)_{j}<\infty
\end{align*}
for all $j\in I$ and $\delta\in\left(0,1\right]$. By equation (\ref{eq:OscillationForFreeSequenceNormEstimateBanachCase})
and by solidity of $\ell_{w}^{q}\left(I\right)$, this yields
\[
\left\Vert \left(\sup_{0<\delta\leq1}\frac{1}{\delta}\left\Vert \osc{\delta\cdot T_{j}^{-T}\left[-1,1\right]^{\dimension}}\left[M_{-b_{j}}\left(\gamma^{\left(j\right)}\ast f\right)\right]\right\Vert _{L_{v}^{p}}\right)_{j\in I}\right\Vert _{\ell_{w}^{q}}\leq C_{1}\cdot\left\Vert \overrightarrow{B}c\right\Vert _{\ell_{w}^{q}}\leq C_{1}\cdot\vertiii{\smash{\Gamma_{\CalQ}}}\cdot\vertiii{\smash{\overrightarrow{B}}}\cdot\left\Vert f\right\Vert _{\DecompSp{\CalQ}p{\ell_{w}^{q}}v}<\infty.
\]
Finally, Theorem \ref{thm:ConvolvingDecompositionSpaceWithGammaJ}
shows
\[
\left\Vert \left(\left\Vert \gamma^{\left(j\right)}\ast f\right\Vert _{L_{v}^{p}}\right)_{j\in I}\right\Vert _{\ell_{w}^{q}}\leq\vertiii{\smash{\Gamma_{\CalQ}}}\cdot\vertiii{\smash{\overrightarrow{A}}}\cdot\left\Vert f\right\Vert _{\DecompSp{\CalQ}p{\ell_{w}^{q}}v}<\infty.
\]
It is not hard to see that this implies boundedness of ${\rm Ana}_{{\rm osc}}$,
with a bound for the operator norm as in the statement of the lemma,
since $2^{\max\left\{ 0,\frac{1}{q}-1\right\} }$ is a valid triangle
constant for $\ell_{w}^{q}\left(I\right)$ and since $C_{1}\geq1$.

\medskip{}

In case of $p\in\left(0,1\right)$, we first note that equation (\ref{eq:WeightedWienerAmalgamTemperedDistribution})
yields $V_{j}=W_{T_{j}^{-T}\left[-1,1\right]^{\dimension}}\left(L_{v}^{p}\right)\hookrightarrow L_{\left(1+\left|\mybullet\right|\right)^{-K}}^{\infty}\left(\R^{\dimension}\right)$,
where again the norm of the embedding might depend heavily on the
choice of $j\in I$. But as seen in Theorem \ref{thm:ConvolvingDecompositionSpaceWithGammaJ},
the series in equation (\ref{eq:OscillationForFreeSeriesRepresentation})
converges in $V_{j}$ and hence in $L_{\left(1+\left|\mybullet\right|\right)^{-K}}^{\infty}\left(\R^{\dimension}\right)$,
which yields \emph{locally} uniform convergence, since each summand
of the series is continuous. In particular, we get continuity of $\gamma^{\left(j\right)}\ast f$
for each $j\in I$.

The remainder of the argument is similar as that for $p\in\left[1,\infty\right]$.
Nevertheless, it needs to be modified slightly, since for $p\in\left(0,1\right)$,
$L_{v}^{p}\left(\R^{\dimension}\right)$ does not satisfy the triangle
inequality, but instead the so-called $p$-triangle inequality, i.e.,
$\left\Vert f+g\right\Vert _{L_{v}^{p}}^{p}\leq\left\Vert f\right\Vert _{L_{v}^{p}}^{p}+\left\Vert g\right\Vert _{L_{v}^{p}}^{p}$.
Precisely, using equation (\ref{eq:OscillationForFreeSeriesRepresentation})
and the estimates $\osc Q\left(\sum_{i\in I}f_{i}\right)\leq\sum_{i\in I}\left(\osc Qf_{i}\right)$
and $M_{Q}\left(\sum_{i\in I}f_{i}\right)\leq\sum_{i\in I}\left(M_{Q}f_{i}\right)$,
as well as the $p$-triangle inequality for $L_{v}^{p}\left(\R^{\dimension}\right)$,
we get for arbitrary $0<\delta\leq1$ that
\begin{align*}
 & \left(\frac{1}{\delta}\left\Vert \osc{\delta\cdot T_{j}^{-T}\left[-1,1\right]^{\dimension}}\left[M_{-b_{j}}\left(\gamma^{\left(j\right)}\ast f\right)\right]\right\Vert _{W_{T_{j}^{-T}\left[-1,1\right]^{\dimension}}\left(L_{v}^{p}\right)}\right)^{p}\\
 & \leq\frac{1}{\delta^{p}}\sum_{i\in I}\left\Vert \osc{\delta\cdot T_{j}^{-T}\left[-1,1\right]^{\dimension}}\left(M_{-b_{j}}\left[\Fourier^{-1}\left(\widehat{\gamma^{\left(j\right)}}\cdot\varphi_{i}\cdot\smash{\widehat{f}}\,\right)\right]\right)\right\Vert _{W_{T_{j}^{-T}\left[-1,1\right]^{\dimension}}\left(L_{v}^{p}\right)}^{p}\\
\left({\scriptstyle \text{Lemma }\ref{lem:MainOscillationEstimate}}\right) & \leq C_{2}^{p}\cdot\sum_{i\in I}\left[\left|\det T_{i}\right|^{1-p}\cdot\left(1+\left\Vert T_{j}^{-1}T_{i}\right\Vert \right)^{pK+\dimension}\cdot\left\Vert \Fourier^{-1}\left(\varphi_{i}^{\ast}\smash{\widehat{f}}\,\right)\right\Vert _{L_{v}^{p}}^{p}\cdot\left\Vert \Fourier^{-1}\left[\widehat{\phi^{\left(j\right)}}\cdot\varphi_{i}\right]\right\Vert _{L_{v_{0}}^{p}}^{p}\right]\\
\left({\scriptstyle \text{with }\theta_{i}:=c_{i}^{p}=\left\Vert \Fourier^{-1}\left(\varphi_{i}^{\ast}\smash{\widehat{f}}\,\right)\right\Vert _{L_{v}^{p}}^{p}}\right) & =C_{2}^{p}\cdot\sum_{i\in I}\left[B_{j,i}\theta_{i}\right]=C_{2}^{p}\cdot\left(\overrightarrow{B}\theta\right)_{j},
\end{align*}
where the constant $C_{2}>0$ is provided by Lemma \ref{lem:MainOscillationEstimate}.

We conclude using the solidity of $\ell_{w}^{q}\left(I\right)$ that
\begin{align*}
 & \left\Vert \left(\sup_{0<\delta\leq1}\frac{1}{\delta}\left\Vert \osc{\delta\cdot T_{j}^{-T}\left[-1,1\right]^{\dimension}}\left[M_{-b_{j}}\left(\gamma^{\left(j\right)}\ast f\right)\right]\right\Vert _{W_{T_{j}^{-T}\left[-1,1\right]^{\dimension}}\left(L_{v}^{p}\right)}\right)_{j\in I}\right\Vert _{\ell_{w}^{q}}\\
 & \leq C_{2}\cdot\left\Vert \left(\overrightarrow{B}\theta\right)^{1/p}\right\Vert _{\ell_{w}^{q}}=C_{2}\cdot\left\Vert \left(w^{p}\cdot\overrightarrow{B}\theta\right)^{1/p}\right\Vert _{\ell^{q}}\\
 & =C_{2}\cdot\left\Vert w^{p}\cdot\overrightarrow{B}\theta\right\Vert _{\ell^{q/p}}^{1/p}=C_{2}\cdot\left\Vert \overrightarrow{B}\theta\right\Vert _{\ell_{w^{\min\left\{ 1,p\right\} }}^{r}}^{1/p}\\
 & \leq C_{2}\cdot\left(\vertiii{\smash{\overrightarrow{B}}}\cdot\left\Vert \theta\right\Vert _{\ell_{w^{\min\left\{ 1,p\right\} }}^{r}}\right)^{1/p}\\
 & =C_{2}\cdot\vertiii{\smash{\overrightarrow{B}}}^{1/p}\cdot\left\Vert w^{p}\cdot\theta\right\Vert _{\ell^{q/p}}^{1/p}\\
 & =C_{2}\cdot\vertiii{\smash{\overrightarrow{B}}}^{1/p}\cdot\left\Vert w\cdot\theta^{1/p}\right\Vert _{\ell^{q}}\\
 & =C_{2}\cdot\vertiii{\smash{\overrightarrow{B}}}^{1/p}\cdot\left\Vert c\right\Vert _{\ell_{w}^{q}}.
\end{align*}

Finally, using the quasi-triangle inequality $\left\Vert \sum_{i=1}^{N}f_{i}\right\Vert _{L^{p}}\leq N^{\frac{1}{p}-1}\cdot\sum_{i=1}^{N}\left\Vert f_{i}\right\Vert _{L^{p}}$
(cf.\@ \cite[Exercise 1.1.5(c)]{GrafakosClassicalFourierAnalysis})
and the estimate $\left|i^{\ast}\right|\leq N_{\CalQ}$ for all $i\in I$,
we also get $c_{i}\leq N_{\CalQ}^{\frac{1}{p}-1}\cdot\left(\Gamma_{\CalQ}\,d\right)_{i}$
for all $i\in I$ and $d_{i}:=\left\Vert \Fourier^{-1}\left(\varphi_{i}\cdot\widehat{f}\right)\right\Vert _{L_{v}^{p}}$.
Hence,
\[
\left\Vert \left(\sup_{0<\delta\leq1}\frac{1}{\delta}\left\Vert \osc{\delta\cdot T_{j}^{-T}\left[-1,1\right]^{\dimension}}\left[M_{-b_{j}}\left(\gamma^{\left(j\right)}\ast f\right)\right]\right\Vert _{W_{T_{j}^{-T}\left[-1,1\right]^{\dimension}}\left(L_{v}^{p}\right)}\right)_{j\in I}\right\Vert _{\ell_{w}^{q}}\leq C_{2}\cdot N_{\CalQ}^{\frac{1}{p}-1}\cdot\vertiii{\smash{\Gamma_{\CalQ}}}\cdot\vertiii{\smash{\overrightarrow{B}}}^{1/p}\cdot\left\Vert d\right\Vert _{\ell_{w}^{q}}.
\]
Because of $\left\Vert d\right\Vert _{\ell_{w}^{q}}=\left\Vert f\right\Vert _{\DecompSp{\CalQ}p{\ell_{w}^{q}}v}$
and in combination with Theorem \ref{thm:ConvolvingDecompositionSpaceWithGammaJ},
we can now derive the claim using the same arguments as for $p\in\left[1,\infty\right]$.
Here, we use that $C_{\CalQ}\geq\left\Vert T_{i}^{-1}T_{i}\right\Vert =1$,
so that the constant $C_{3}>0$ provided by Theorem \ref{thm:ConvolvingDecompositionSpaceWithGammaJ}
(for $p\in\left(0,1\right)$) satisfies
\begin{align*}
C_{3} & =N_{\CalQ}^{\frac{1}{p}-1}\cdot\left(12288\cdot\dimension^{3/2}\cdot\left\lceil K+\frac{\dimension+1}{p}\right\rceil \right)^{\left\lceil K+\frac{\dimension+1}{p}\right\rceil +1}\cdot\left(1+R_{\CalQ}\right)^{\dimension/p}\left(12R_{\CalQ}C_{\CalQ}\right)^{\dimension\left(\frac{1}{p}-1\right)}\cdot\Omega_{0}^{K}\Omega_{1}\\
\left({\scriptstyle \text{since }\Omega_{0},\Omega_{1}\geq1}\right) & \leq N_{\CalQ}^{\frac{1}{p}-1}\cdot12^{\dimension\left(\frac{1}{p}-1\right)}\left(12288\cdot\dimension^{3/2}\cdot\left\lceil K+\frac{\dimension+1}{p}\right\rceil \right)^{\left\lceil K+\frac{\dimension+1}{p}\right\rceil +1}\cdot\left(1+C_{\CalQ}R_{\CalQ}\right)^{\dimension\left(\frac{2}{p}-1\right)}\cdot\Omega_{0}^{5K}\Omega_{1}^{5},\\
 & \leq N_{\CalQ}^{\frac{1}{p}-1}\cdot12^{\dimension\left(\frac{1}{p}-1\right)}\left(12288\cdot\dimension^{3/2}\cdot\left\lceil K+\frac{\dimension+1}{p}\right\rceil \right)^{\left\lceil K+\frac{\dimension+1}{p}\right\rceil +1}\cdot\left(1+C_{\CalQ}R_{\CalQ}\right)^{\frac{2\dimension}{p}}\cdot\Omega_{0}^{5K}\Omega_{1}^{5},
\end{align*}
so that
\begin{align*}
\frac{C_{3}}{N_{\CalQ}^{\frac{1}{p}-1}C_{2}} & \leq370\cdot\dimension^{11/2}\cdot\dimension^{\dimension/2p}\cdot\frac{12^{\dimension\left(\frac{1}{p}-1\right)}\left(12288\cdot\dimension^{3/2}\cdot\left\lceil K+\frac{\dimension+1}{p}\right\rceil \right)^{\left\lceil K+\frac{\dimension+1}{p}\right\rceil +1}}{2^{16\frac{\dimension}{p}}\cdot\left(4032\cdot\dimension^{3}\cdot\left\lceil K+\frac{\dimension+1}{p}\right\rceil \right)^{2\left\lceil K+\frac{\dimension+1}{p}\right\rceil +2}}\\
 & =\frac{370\cdot\dimension^{11/2}\cdot\dimension^{\dimension/2p}\cdot12^{\dimension\left(\frac{1}{p}-1\right)}}{2^{16\frac{\dimension}{p}}\left(4032\cdot\dimension^{3}\cdot\left\lceil K+\frac{\dimension+1}{p}\right\rceil \right)^{\left\lceil K+\frac{\dimension+1}{p}\right\rceil +1}}\cdot\left(\frac{12288\cdot\dimension^{3/2}}{4032\cdot\dimension^{3}}\right)^{\left\lceil K+\frac{\dimension+1}{p}\right\rceil +1}\\
 & \leq\frac{370\cdot\dimension^{11/2}\cdot\dimension^{\dimension/2p}\cdot12^{\frac{\dimension}{p}}}{2^{16\frac{\dimension}{p}}}\cdot\left(\frac{4}{4032\cdot\dimension^{9/2}}\right)^{\left\lceil K+\frac{\dimension+1}{p}\right\rceil +1}\\
\left({\scriptstyle \text{since }\left\lceil K+\frac{\dimension+1}{p}\right\rceil \geq\frac{\dimension+1}{p}\geq\frac{\dimension}{p}+1}\right) & \leq\frac{370\cdot\dimension^{11/2}\cdot\dimension^{\dimension/2p}}{2^{12\frac{\dimension}{p}}}\cdot\left(\frac{1}{1000\cdot\dimension^{9/2}}\right)^{\frac{\dimension}{p}+2}\\
 & \leq\frac{370}{1000000}\frac{\dimension^{11/2}\cdot\dimension^{\dimension/2p}}{2^{12\frac{\dimension}{p}}}\cdot\dimension^{-\frac{9}{2}\frac{\dimension}{p}}\dimension^{-9}=\frac{370}{1000000}\frac{1}{2^{12\frac{\dimension}{p}}}\cdot\dimension^{-4\frac{\dimension}{p}}\dimension^{-7/2}\leq1.\qedhere
\end{align*}
\end{proof}
In view of the preceding lemma, we know that (if $\Gamma=\left(\gamma_{i}\right)_{i\in I}$
fulfills Assumption \ref{assu:DiscreteBanachFrameAssumptions}) each
of the functions $\gamma^{\left(j\right)}\ast f$ and hence also $\gamma^{\left[j\right]}\ast f=\left|\det T_{j}\right|^{-1/2}\cdot\gamma^{\left(j\right)}\ast f$
is continuous, so that the coefficient mapping
\[
\DecompSp{\CalQ}p{\ell_{w}^{q}}v\ni f\mapsto\left[\left(\gamma^{\left[j\right]}\ast f\right)\!\left(\delta\cdot T_{j}^{-T}k\right)\right]_{j\in I,\,k\in\Z^{\dimension}}\in\Compl^{I\times\Z^{\dimension}}
\]
is well-defined. But eventually, we want to show that this map yields
a Banach frame for $\DecompSp{\CalQ}p{\ell_{w}^{q}}v$, so that we
also need to construct a suitable ``reconstruction mapping'', which
can recover $f$ from these coefficients. The following lemma is an
important ingredient for the construction of this reconstruction map.
\begin{lem}
\label{lem:DiscreteSynthesisOperatorIsAlmostIsometric}For $i\in I$
and $0<\delta\leq1$, let
\[
{\rm Synth}_{\delta,i}:\Compl^{\Z^{\dimension}}\to\left\{ f:\R^{\dimension}\to\Compl\with f\text{ measurable}\right\} ,\left(c_{k}\right)_{k\in\Z^{\dimension}}\mapsto M_{b_{i}}\left[\sum_{k\in\Z^{\dimension}}c_{k}\cdot e^{-2\pi i\left\langle b_{i},\delta\cdot T_{i}^{-T}k\right\rangle }\Indicator_{\delta\cdot T_{i}^{-T}\left(k+\left[0,1\right)^{\dimension}\right)}\right].
\]
Then ${\rm Synth}_{\delta,i}$ is well-defined and yields a bounded
operator ${\rm Synth}_{\delta,i}:C_{i}^{\left(\delta\right)}\to V_{i}$,
where the $i$-th coefficient space $C_{i}^{\left(\delta\right)}$
is defined as in equation (\ref{eq:CoefficientSpaceDefinition}).
More precisely, we have
\[
\frac{\delta^{\dimension/p}}{\Omega_{0}^{K}\Omega_{1}\cdot\left(1+\sqrt{\dimension}\right)^{K}}\cdot\left|\det T_{i}\right|^{-1/p}\cdot\left\Vert c\right\Vert _{C_{i}^{\left(\delta\right)}}\leq\left\Vert {\rm Synth}_{\delta,i}\,c\right\Vert _{V_{i}}\leq C_{\dimension,p,\delta,K}\cdot\left|\det T_{i}\right|^{-1/p}\cdot\left\Vert c\right\Vert _{C_{i}^{\left(\delta\right)}}\qquad\forall c\in\Compl^{\Z^{\dimension}},
\]
with
\[
C_{\dimension,p,\delta,K}=\begin{cases}
\left(1+\sqrt{\dimension}\right)^{K}\cdot\Omega_{0}^{K}\Omega_{1}\cdot\delta^{\dimension/p}, & \text{if }p\in\left[1,\infty\right],\\
4^{\dimension/p}\cdot\left(1+2\sqrt{\dimension}\right)^{K}\cdot\Omega_{0}^{K}\Omega_{1}, & \text{if }p\in\left(0,1\right).
\end{cases}\qedhere
\]
\end{lem}
\begin{proof}
First note that ${\rm Synth}_{\delta,i}$ is well-defined, since the
sets $\left(\delta\cdot T_{i}^{-T}\left(k+\left[0,1\right)^{\dimension}\right)\right)_{k\in\Z^{\dimension}}$
are pairwise disjoint. Also, we can ignore the modulation $M_{b_{i}}$
in the following, since $\left\Vert M_{b_{i}}f\right\Vert _{V_{i}}=\left\Vert f\right\Vert _{V_{i}}$,
because of $\left\Vert f\right\Vert _{V_{i}}=\left\Vert g\right\Vert _{V_{i}}$
for measurable $f,g$ satisfying $\left|f\right|=\left|g\right|$.
Furthermore, since we have for $x\in\delta\cdot T_{i}^{-T}\left(k+\left[0,1\right)^{\dimension}\right)$,
i.e., for $x=\delta\cdot T_{i}^{-T}k+\delta T_{i}^{-T}q$ with $q\in\left[0,1\right)^{\dimension}$
that
\begin{align*}
v_{k}^{\left(i,\delta\right)} & =v\left(\delta\cdot T_{i}^{-T}k\right)=v\left(x-\delta T_{i}^{-T}q\right)\leq v\left(x\right)\cdot v_{0}\left(-\delta T_{i}^{-T}q\right)\\
\left({\scriptstyle \text{assump. on }v_{0}}\right) & \leq\Omega_{1}\cdot\left(1+\left|\delta T_{i}^{-T}q\right|\right)^{K}\cdot v\left(x\right)\\
\left({\scriptstyle \text{eq. }\eqref{eq:WeightLinearTransformationsConnection}}\right) & \leq\Omega_{0}^{K}\Omega_{1}\cdot\left(1+\left|\delta q\right|\right)^{K}\cdot v\left(x\right)\\
\left({\scriptstyle \text{since }\delta\leq1}\right) & \leq\Omega_{0}^{K}\Omega_{1}\cdot\left(1+\sqrt{\dimension}\right)^{K}\cdot v\left(x\right),
\end{align*}
Lemma \ref{lem:MaximalFunctionDominatesF} implies
\begin{equation}
\begin{split}\left\Vert {\rm Synth}_{\delta,i}\left(c_{k}\right)_{k\in\Z^{\dimension}}\right\Vert _{V_{i}} & \geq\left\Vert {\rm Synth}_{\delta,i}\left(c_{k}\right)_{k\in\Z^{\dimension}}\right\Vert _{L_{v}^{p}}\\
\left({\scriptstyle \text{pairwise disjointness}}\right) & =\left[\sum_{k\in\Z^{\dimension}}\left|c_{k}\right|^{p}\int_{\delta\cdot T_{i}^{-T}\left(k+\left[0,1\right)^{\dimension}\right)}\left[v\left(x\right)\right]^{p}\d x\right]^{1/p}\\
 & \geq\frac{1}{\Omega_{0}^{K}\Omega_{1}\cdot\left(1+\sqrt{\dimension}\right)^{K}}\cdot\left[\sum_{k\in\Z^{\dimension}}\left|v_{k}^{\left(i,\delta\right)}\cdot c_{k}\right|^{p}\cdot\lambda_{\dimension}\left(\delta\cdot T_{i}^{-T}\left[k+\left[0,1\right)^{\dimension}\right]\right)\right]^{1/p}\\
 & =\frac{\delta^{\dimension/p}\cdot\left|\det T_{i}\right|^{-1/p}}{\Omega_{0}^{K}\Omega_{1}\cdot\left(1+\sqrt{\dimension}\right)^{K}}\cdot\left\Vert \left(c_{k}\right)_{k\in\Z^{\dimension}}\right\Vert _{C_{i}^{\left(\delta\right)}}.
\end{split}
\label{eq:DiscreteSynthesisOperatorLowerBound}
\end{equation}
This proves the lower bound.

\medskip{}

Now, we establish the reverse inequality for $p\in\left[1,\infty\right]$:
For $x=\delta\cdot T_{i}^{-T}k+\delta T_{i}^{-T}q\in\delta\cdot T_{i}^{-T}\left(k+\left[-L,L\right]^{\dimension}\right)$
with $L\geq1$, we have as above that
\begin{align}
v\left(x\right) & =v\left(\delta\cdot T_{i}^{-T}k+\delta\cdot T_{i}^{-T}q\right)\nonumber \\
 & \leq v\left(\delta\cdot T_{i}^{-T}k\right)\cdot v_{0}\left(\delta\cdot T_{i}^{-T}q\right)\nonumber \\
\left({\scriptstyle \text{assump. on }v_{0}\text{ and eq. }\eqref{eq:WeightLinearTransformationsConnection}}\right) & \leq\Omega_{0}^{K}\Omega_{1}\cdot v_{k}^{\left(i,\delta\right)}\cdot\left(1+\left|\delta\cdot q\right|\right)^{K}\nonumber \\
\left({\scriptstyle \text{since }\delta\leq1}\right) & \leq\left(1+L\sqrt{\dimension}\right)^{K}\cdot\Omega_{0}^{K}\Omega_{1}\cdot v_{k}^{\left(i,\delta\right)}.\label{eq:DiscreteSynthesisOperatorUpperWeightEstimate}
\end{align}
Furthermore, since $p\in\left[1,\infty\right]$, the first inequality
in estimate (\ref{eq:DiscreteSynthesisOperatorLowerBound}) from above
is actually an equality, so that
\begin{align*}
\left\Vert {\rm Synth}_{\delta,i}\left(c_{k}\right)_{k\in\Z^{\dimension}}\right\Vert _{V_{i}} & =\left[\sum_{k\in\Z^{\dimension}}\left|c_{k}\right|^{p}\int_{\delta\cdot T_{i}^{-T}\left(k+\left[0,1\right)^{\dimension}\right)}\left[v\left(x\right)\right]^{p}\d x\right]^{1/p}\\
 & \leq\left(1+\sqrt{\dimension}\right)^{K}\cdot\Omega_{0}^{K}\Omega_{1}\cdot\left[\sum_{k\in\Z^{\dimension}}\left|v_{k}^{\left(i,\delta\right)}c_{k}\right|^{p}\cdot\lambda_{\dimension}\left(\delta\cdot T_{i}^{-T}\left[k+\left[0,1\right)^{\dimension}\right]\right)\right]^{1/p}\\
 & =\left(1+\sqrt{\dimension}\right)^{K}\cdot\Omega_{0}^{K}\Omega_{1}\cdot\delta^{\dimension/p}\cdot\left|\det T_{i}\right|^{-1/p}\cdot\left\Vert \left(c_{k}\right)_{k\in\Z^{\dimension}}\right\Vert _{C_{i}^{\left(\delta\right)}}.
\end{align*}

\medskip{}

Finally, for $p\in\left(0,1\right)$, we use the estimate $M_{Q}\left(\sum_{i\in I}f_{i}\right)\leq\sum_{i\in I}M_{Q}f_{i}$
and the $p$-triangle inequality for $L_{v}^{p}\left(\R^{\dimension}\right)$
to deduce
\begin{align*}
\left\Vert {\rm Synth}_{\delta,i}\left(c_{k}\right)_{k\in\Z^{\dimension}}\right\Vert _{V_{i}}^{p} & \leq\left\Vert \sum_{k\in\Z^{\dimension}}M_{T_{i}^{-T}\left[-1,1\right]^{\dimension}}\left(c_{k}\cdot e^{-2\pi i\left\langle b_{i},\delta\cdot T_{i}^{-T}k\right\rangle }\Indicator_{\delta\cdot T_{i}^{-T}\left(k+\left[0,1\right)^{\dimension}\right)}\right)\right\Vert _{L_{v}^{p}}^{p}\\
 & \leq\sum_{k\in\Z^{\dimension}}\left[\left|c_{k}\right|^{p}\cdot\left\Vert M_{T_{i}^{-T}\left[-1,1\right]^{\dimension}}\Indicator_{\delta\cdot T_{i}^{-T}\left(k+\left[0,1\right)^{\dimension}\right)}\right\Vert _{L_{v}^{p}}^{p}\right].
\end{align*}
Next, observe
\[
\left(M_{T_{i}^{-T}\left[-1,1\right]^{\dimension}}\Indicator_{\delta\cdot T_{i}^{-T}\left(k+\left[0,1\right)^{\dimension}\right)}\right)\left(x\right)=\left\Vert \Indicator_{\delta\cdot T_{i}^{-T}\left(k+\left[0,1\right)^{\dimension}\right)}\cdot\Indicator_{x+T_{i}^{-T}\left[-1,1\right]^{\dimension}}\right\Vert _{L^{\infty}}\leq1.
\]
Furthermore, if the function inside the $\left\Vert \mybullet\right\Vert _{L^{\infty}}$
norm does not vanish identically, we have
\[
x\in\delta T_{i}^{-T}k+T_{i}^{-T}\left(\delta\left[0,1\right)^{\dimension}-\left[-1,1\right]^{\dimension}\right)\subset\delta T_{i}^{-T}k+T_{i}^{-T}\left[-2,2\right]^{\dimension}.
\]
Hence, equation (\ref{eq:DiscreteSynthesisOperatorUpperWeightEstimate})
yields $v\left(x\right)\leq\left(1+2\sqrt{\dimension}\right)^{K}\Omega_{0}^{K}\Omega_{1}\cdot v_{k}^{\left(i,\delta\right)}$
and thus
\begin{align*}
\left\Vert M_{T_{i}^{-T}\left[-1,1\right]^{\dimension}}\Indicator_{\delta\cdot T_{i}^{-T}\left(k+\left[0,1\right)^{\dimension}\right)}\right\Vert _{L_{v}^{p}}^{p} & \leq\left(\left(1+2\sqrt{\dimension}\right)^{K}\Omega_{0}^{K}\Omega_{1}\right)^{p}\cdot\left[v_{k}^{\left(i,\delta\right)}\right]^{p}\cdot\lambda_{\dimension}\left(\delta T_{i}^{-T}k+T_{i}^{-T}\left[-2,2\right]^{\dimension}\right)\\
 & =4^{\dimension}\cdot\left(\left(1+2\sqrt{\dimension}\right)^{K}\Omega_{0}^{K}\Omega_{1}\right)^{p}\cdot\left|\det T_{i}\right|^{-1}\cdot\left[v_{k}^{\left(i,\delta\right)}\right]^{p},
\end{align*}
so that we get
\[
\left\Vert {\rm Synth}_{\delta,i}\left(c_{k}\right)_{k\in\Z^{\dimension}}\right\Vert _{V_{i}}^{p}\leq4^{\dimension}\cdot\left(\left(1+2\sqrt{\dimension}\right)^{K}\Omega_{0}^{K}\Omega_{1}\right)^{p}\cdot\left|\det T_{i}\right|^{-1}\cdot\left\Vert \left(c_{k}\right)_{k\in\Z^{\dimension}}\right\Vert _{C_{i}^{\left(\delta\right)}}^{p},
\]
as claimed.
\end{proof}
Below, we will employ a Neumann series argument to construct the reconstruction
operator $R$. To this end, we need to know that the space $\ell_{w}^{q}\left(\left[V_{i}\right]_{i\in I}\right)$
is a Quasi-Banach space.
\begin{lem}
\label{lem:IteratedSequenceSpaceComplete}Let $\left(X_{i}\right)_{i\in I}$
be a sequence of Quasi-Banach spaces, each with $\left\Vert x+y\right\Vert _{X_{i}}\leq C_{i}\cdot\left(\left\Vert x\right\Vert _{X_{i}}+\left\Vert y\right\Vert _{X_{i}}\right)$
for all $x,y\in X_{i}$ and suitable $C_{i}\geq1$. Assume that $C:=\sup_{i\in i}C_{i}$
is finite and that each quasi-norm $\left\Vert \mybullet\right\Vert _{X_{i}}$
is continuous.

Define
\[
\ell_{w}^{q}\left(\left[X_{i}\right]_{i\in I}\right):=\left\{ x=\left(x_{i}\right)_{i\in I}\in\prod_{i\in I}X_{i}\with\left\Vert x\right\Vert :=\left\Vert \left(\left\Vert x_{i}\right\Vert _{X_{i}}\right)_{i\in I}\right\Vert _{\ell_{w}^{q}\left(I\right)}<\infty\right\} .
\]
Then $\left(\ell_{w}^{q}\left(\left[X_{i}\right]_{i\in I}\right),\left\Vert \mybullet\right\Vert \right)$
is a Quasi-Banach space.
\end{lem}
\begin{rem*}
The lemma applies in particular with the choice $X_{i}=V_{i}$. Indeed,
$\left\Vert \mybullet\right\Vert _{L_{v}^{p}}$ is an $s$-norm for
$s:=\min\left\{ 1,p\right\} $; since $M_{Q}\left(f+g\right)\leq M_{Q}f+M_{Q}g$,
we get $\left\Vert f+g\right\Vert _{W_{Q}\left(L_{v}^{p}\right)}^{s}\leq\left\Vert f\right\Vert _{W_{Q}\left(L_{v}^{p}\right)}^{s}+\left\Vert g\right\Vert _{W_{Q}\left(L_{v}^{p}\right)}^{s}$,
so that $\left\Vert \mybullet\right\Vert _{W_{Q}\left(L_{v}^{p}\right)}$
is also an $s$-norm and hence continuous, since $\left|\left\Vert x_{n}\right\Vert ^{s}-\left\Vert x\right\Vert ^{s}\right|\leq\left\Vert x_{n}-x\right\Vert ^{s}\xrightarrow[n\to\infty]{}0$
for any $s$-norm $\left\Vert \mybullet\right\Vert $ if $\left\Vert x_{n}-x\right\Vert \xrightarrow[n\to\infty]{}0$.
Furthermore, in case of $p\in\left[1,\infty\right]$, one can choose
$C_{i}=1$ for all $i\in I$. Finally, for $p\in\left(0,1\right)$,
Remark \ref{rem:MainAssumptionsRemark} shows that each $V_{i}$ is
a Quasi-Banach space and that we can choose $C_{i}=2^{\frac{1}{p}-1}$
for all $i\in I$. Hence, $V=\ell_{w}^{q}\left(\left[V_{i}\right]_{i\in I}\right)$
is a Quasi-Banach space.
\end{rem*}
\begin{proof}
For brevity, let $X:=\ell_{w}^{q}\left(\left[X_{i}\right]_{i\in I}\right)$.
It is clear that $X$ is closed under multiplication with scalars
and that $\left\Vert \alpha\cdot x\right\Vert =\left|\alpha\right|\cdot\left\Vert x\right\Vert $
for $\alpha\in\mathbb{K}$ (with $\mathbb{K}\in\left\{ \R,\Compl\right\} $)
and $x\in X$. Furthermore, if $\left\Vert x\right\Vert =0$ for $x=\left(x_{i}\right)_{i\in I}$,
then $\left\Vert x_{i}\right\Vert _{X_{i}}=0$ for all $i\in I$,
so that $x=0$. Finally, for $x,y\in X$, we have by solidity of $\ell_{w}^{q}\left(I\right)$
that 
\begin{align*}
\left\Vert x+y\right\Vert  & =\left\Vert \left(\left\Vert x_{i}+y_{i}\right\Vert _{X_{i}}\right)_{i\in I}\right\Vert _{\ell_{w}^{q}}\leq\left\Vert \left(C\cdot\left[\left\Vert x_{i}\right\Vert _{X_{i}}+\left\Vert y_{i}\right\Vert _{X_{i}}\right]\right)_{i\in I}\right\Vert _{\ell_{w}^{q}}\\
 & \leq C\cdot C_{q}\cdot\left[\left\Vert \left(\left\Vert x_{i}\right\Vert _{X_{i}}\right)_{i\in I}\right\Vert _{\ell_{w}^{q}}+\left\Vert \left(\left\Vert y_{i}\right\Vert _{X_{i}}\right)_{i\in I}\right\Vert _{\ell_{w}^{q}}\right]\\
 & =C\cdot C_{q}\cdot\left[\left\Vert x\right\Vert +\left\Vert y\right\Vert \right]<\infty,
\end{align*}
where $C_{q}$ is a triangle constant for $\ell_{w}^{q}\left(I\right)$.
Hence, $X$ is closed under addition (and thus a vector space as a
subspace of $\prod_{i\in I}X_{i}$) and $\left\Vert \mybullet\right\Vert $
is a quasi-norm on $X$.

\medskip{}

Now, let $\left(x^{\left(n\right)}\right)_{n\in\N}=\left[\left(\smash{x_{i}^{\left(n\right)}}\right)_{i\in I}\right]_{n\in\N}$
be a Cauchy sequence in $X$. It is not hard to see that each of the
projections $\pi_{i}:X\to X_{i},\left(x_{j}\right)_{j\in I}\mapsto x_{i}$
is a bounded linear map, so that each sequence $\left(\smash{x_{i}^{\left(n\right)}}\right)_{n\in\N}$
is Cauchy in $X_{i}$ and hence convergent to some $x_{i}\in X_{i}$.
Now, let $\varepsilon>0$ be arbitrary. There is some $N_{0}\in\N$
satisfying $\left\Vert x^{\left(n\right)}-x^{\left(m\right)}\right\Vert \leq\varepsilon$
for all $n,m\geq N_{0}$. By Fatou's lemma and by continuity of $\left\Vert \mybullet\right\Vert _{X_{i}}$,
this implies for $m\geq N_{0}$ that
\begin{align*}
\left\Vert \left(\left\Vert x_{i}-\smash{x_{i}^{\left(m\right)}}\right\Vert _{X_{i}}\right)_{i\in I}\right\Vert _{\ell_{w}^{q}} & =\left\Vert \left(\liminf_{n\to\infty}\left\Vert \smash{x_{i}^{\left(n\right)}}-\smash{x_{i}^{\left(m\right)}}\right\Vert _{X_{i}}\right)_{i\in I}\right\Vert _{\ell_{w}^{q}}\leq\liminf_{n\to\infty}\left\Vert \left(\left\Vert \smash{x_{i}^{\left(n\right)}}-\smash{x_{i}^{\left(m\right)}}\right\Vert _{X_{i}}\right)_{i\in I}\right\Vert _{\ell_{w}^{q}}\\
 & =\liminf_{n\to\infty}\left\Vert \smash{x^{\left(n\right)}}-\smash{x^{\left(m\right)}}\right\Vert \leq\varepsilon<\infty.
\end{align*}
Since $X$ is a vector space, this implies $x=\left(x_{i}\right)_{i\in I}=\left(x-\smash{x^{\left(m\right)}}\right)+x^{\left(m\right)}\in X$,
as well as $\left\Vert x-\smash{x^{\left(m\right)}}\right\Vert \xrightarrow[m\to\infty]{}0$.
\end{proof}
The next lemma is our final preparation for proving that the coefficient
map
\[
f\mapsto\left[\left(\gamma^{\left[j\right]}\ast f\right)\!\left(\delta\cdot T_{j}^{-T}k\right)\right]_{j\in I,\,k\in\Z^{\dimension}}
\]
indeed yields a Banach frame for $\DecompSp{\CalQ}p{\ell_{w}^{q}}v$.
This lemma essentially yields a replacement for the usual reproducing
kernel property which is used in the theory of coorbit spaces (cf.\@
\cite{FeichtingerCoorbit0,FeichtingerCoorbit1,FeichtingerCoorbit2,RauhutCoorbitQuasiBanach}
and \cite[Section 2]{VoigtlaenderPhDThesis}).
\begin{lem}
\label{lem:SpecialProjection}Assume that $\Gamma=\left(\gamma_{i}\right)_{i\in I}$
satisfies Assumptions \ref{assu:DiscreteBanachFrameAssumptions} and
\ref{assu:GammaCoversOrbit}. We clearly have a norm-decreasing embedding
$W_{j}\hookrightarrow V_{j}$ and hence also $\iota:\ell_{w}^{q}\left(\left[W_{j}\right]_{j\in I}\right)\hookrightarrow\ell_{w}^{q}\left(\left[V_{j}\right]_{j\in I}\right)$.

Let
\[
{\rm Ana}_{{\rm osc}}:\DecompSp{\CalQ}p{\ell_{w}^{q}}v\to\ell_{w}^{q}\left(\left[W_{j}\right]_{j\in I}\right),f\mapsto\left(\gamma^{\left(j\right)}\ast f\right)_{j\in I}
\]
as in Lemma \ref{lem:OscillationForFree}, let
\[
{\rm Synth}_{\CalD}:\ell_{w}^{q}\left(\left[V_{j}\right]_{j\in I}\right)\to\DecompSp{\CalQ}p{\ell_{w}^{q}}v,\left(f_{i}\right)_{i\in I}\mapsto\sum_{i\in I}\left[\Fourier^{-1}\left(\varphi_{i}\cdot\widehat{f_{i}}\right)\right]
\]
be defined as in Lemma \ref{lem:DecompositionSynthesis} and let 
\[
m_{\theta}:\ell_{w}^{q}\left(\left[V_{j}\right]_{j\in I}\right)\to\ell_{w}^{q}\left(\left[V_{j}\right]_{j\in I}\right),\left(f_{j}\right)_{j\in I}\mapsto\left[\left(\Fourier^{-1}\theta_{j}\right)\ast f_{j}\right]_{j\in I}
\]
be defined as in Lemma \ref{lem:LocalInverseConvolution}.

Then, the map
\[
F:\ell_{w}^{q}\left(\left[V_{j}\right]_{j\in I}\right)\to\ell_{w}^{q}\left(\left[V_{j}\right]_{j\in I}\right),\quad F:=\iota\circ{\rm Ana}_{{\rm osc}}\circ{\rm Synth}_{\CalD}\circ m_{\theta}
\]
is well-defined and bounded and satisfies the following additional
properties:

\begin{enumerate}
\item \label{enu:CompactAnalysisMapsIntoVernal}$F\left[\left(\gamma^{\left(j\right)}\ast f\right)_{j\in I}\right]=\left(\gamma^{\left(j\right)}\ast f\right)_{j\in I}$
for all $f\in\DecompSp{\CalQ}p{\ell_{w}^{q}}v$.
\item ${\rm Synth}_{\CalD}\circ m_{\theta}\circ\iota\circ{\rm Ana}_{{\rm osc}}=\identity_{\DecompSp{\CalQ}p{\ell_{w}^{q}}v}$.
\item $F\circ F=F$.
\item The space $\vernal:=\left\{ \left(f_{i}\right)_{i\in I}\in\ell_{w}^{q}\left(\left[V_{i}\right]_{i\in I}\right)\with F\left(f_{i}\right)_{i\in I}=\left(f_{i}\right)_{i\in I}\right\} $
is a closed subspace of $\ell_{w}^{q}\left(\left[V_{i}\right]_{i\in I}\right)$.
\item For each $f=\left(f_{i}\right)_{i\in I}\in\vernal$, we have that
each $f_{i}:\R^{\dimension}\to\Compl$ is continuous and furthermore
\begin{equation}
\left\Vert \left[\osc{\delta\cdot T_{i}^{-T}\left[-1,1\right]^{\dimension}}\left(M_{-b_{i}}f_{i}\right)\right]_{i\in I}\right\Vert _{\ell_{w}^{q}\left(\left[V_{i}\right]_{i\in I}\right)}\leq\vertiii{F_{0}}\cdot\delta\cdot\left\Vert f\right\Vert _{\ell_{w}^{q}\left(\left[V_{i}\right]_{i\in I}\right)}\qquad\forall\delta\in\left(0,1\right]\label{eq:OscillationEstimateOnVernal}
\end{equation}
for $F_{0}:={\rm Ana}_{{\rm osc}}\circ{\rm Synth}_{\CalD}\circ m_{\theta}:\ell_{w}^{q}\left(\left[V_{i}\right]_{i\in I}\right)\to\ell_{w}^{q}\left(\left[W_{i}\right]_{i\in I}\right)$.
Here, we have
\[
\vertiii{F_{0}}\leq2^{\frac{1}{q}}C_{\CalQ,\Phi,v_{0},p}^{2}\cdot\vertiii{\smash{\Gamma_{\CalQ}}}^{2}\cdot\left(\vertiii{\smash{\overrightarrow{A}}}^{\max\left\{ 1,\frac{1}{p}\right\} }+\vertiii{\smash{\overrightarrow{B}}}^{\max\left\{ 1,\frac{1}{p}\right\} }\right)\cdot C,
\]
for $N:=\left\lceil K+\frac{\dimension+1}{\min\left\{ 1,p\right\} }\right\rceil $
and
\[
C\!:=\!\!\begin{cases}
\frac{\left(2^{16}\cdot768/\dimension^{\frac{3}{2}}\right)^{\frac{\dimension}{p}}}{2^{42}\cdot12^{\dimension}\cdot\dimension^{15}}\!\cdot\!\left(2^{52}\!\cdot\!\dimension^{\frac{25}{2}}\!\cdot\!N^{3}\right)^{N+1}\!\!\!\cdot\!N_{\CalQ}^{2\left(\frac{1}{p}-1\right)}\!\left(1\!+\!R_{\CalQ}C_{\CalQ}\right)^{\dimension\left(\frac{4}{p}-1\right)}\!\!\cdot\Omega_{0}^{13K}\Omega_{1}^{13}\Omega_{2}^{\left(p,K\right)}, & \text{if }p<1,\\
\frac{1}{\sqrt{\dimension}\cdot2^{12+6\left\lceil K\right\rceil }}\cdot\left(2^{17}\cdot\dimension^{5/2}\cdot N\right)^{\left\lceil K\right\rceil +\dimension+2}\cdot\left(1+R_{\CalQ}\right)^{\dimension}\cdot\Omega_{0}^{3K}\Omega_{1}^{3}\Omega_{2}^{\left(p,K\right)}, & \text{if }p\geq1.
\end{cases}\qedhere
\]
\end{enumerate}
\end{lem}
\begin{proof}
As a consequence of Lemmas \ref{lem:DecompositionSynthesis}, \ref{lem:LocalInverseConvolution}
and \ref{lem:OscillationForFree}, we see that $F_{0}:\ell_{w}^{q}\left(\left[V_{i}\right]_{i\in I}\right)\to\ell_{w}^{q}\left(\left[W_{i}\right]_{i\in I}\right)$
is bounded with $\vertiii{F_{0}}\leq\vertiii{{\rm Ana}_{{\rm osc}}}\cdot\vertiii{{\rm Synth}_{\CalD}}\cdot\vertiii{m_{\theta}}$.
By plugging in the estimates for the norms of these operators which
were obtained in the respective lemmas and using elementary estimates,
we easily get the stated estimate for $\vertiii{F_{0}}$. With $F_{0}$,
also $F=\iota\circ F_{0}$ is bounded. We now verify the different
claims individually.

\begin{enumerate}
\item The assumptions of the current lemma include those of Theorem \ref{thm:SemiDiscreteBanachFrame},
where it was shown (cf.\@ equation (\ref{eq:SemiDiscreteBanachFrame}))
that $\identity_{\DecompSp{\CalQ}p{\ell_{w}^{q}}v}={\rm Synth}_{\CalD}\circ m_{\theta}\circ{\rm Ana}_{\Gamma}$,
where ${\rm Ana}_{\Gamma}=\iota\circ{\rm Ana}_{{\rm osc}}$. Hence,
\begin{equation}
{\rm Synth}_{\CalD}\circ m_{\theta}\circ\iota\circ{\rm Ana}_{{\rm osc}}=\identity_{\DecompSp{\CalQ}p{\ell_{w}^{q}}v},\label{eq:SpecialIdentitiy}
\end{equation}
which proves the second part of the current lemma.

Furthermore, for $f\in\DecompSp{\CalQ}p{\ell_{w}^{q}}v$, we have
$\left(\gamma^{\left(j\right)}\ast f\right)_{j\in I}={\rm Ana}_{{\rm osc}}f\in\ell_{w}^{q}\left(\left[W_{i}\right]_{i\in I}\right)\subset\ell_{w}^{q}\left(\left[V_{i}\right]_{i\in I}\right)$,
so that $F\left(\left[\gamma^{\left(j\right)}\ast f\right]_{j\in I}\right)\in\ell_{w}^{q}\left(\left[V_{i}\right]_{i\in I}\right)$
is well-defined. Finally, we get
\begin{align*}
F\left[\left(\gamma^{\left(j\right)}\ast f\right)_{j\in I}\right] & =\left(\iota\circ{\rm Ana}_{{\rm osc}}\right)\left[\left({\rm Synth}_{\CalD}\circ m_{\theta}\right)\left[\gamma^{\left(j\right)}\ast f\right]_{j\in I}\right]\\
 & =\left(\iota\circ{\rm Ana}_{{\rm osc}}\right)\left[\left({\rm Synth}_{\CalD}\circ m_{\theta}\circ\iota\circ{\rm Ana}_{{\rm osc}}\right)f\right]\\
\left({\scriptstyle \text{eq. }\eqref{eq:SpecialIdentitiy}}\right) & =\iota\left({\rm Ana}_{{\rm osc}}f\right)=\iota\left[\left(\gamma^{\left(j\right)}\ast f\right)_{j\in I}\right]=\left(\gamma^{\left(j\right)}\ast f\right)_{j\in I},
\end{align*}
as claimed in the first part.
\item This was proved just above.
\item As a consequence of equation (\ref{eq:SpecialIdentitiy}) (i.e., of
the second part of the lemma), we get
\[
F\circ F=\iota\circ{\rm Ana}_{{\rm osc}}\circ\underbrace{{\rm Synth}_{\CalD}\circ m_{\theta}\circ\iota\circ{\rm Ana}_{{\rm osc}}}_{=\identity_{\DecompSp{\CalQ}p{\ell_{w}^{q}}v}}\circ{\rm Synth}_{\CalD}\circ m_{\theta}=\iota\circ{\rm Ana}_{{\rm osc}}\circ{\rm Synth}_{\CalD}\circ m_{\theta}=F.
\]
\item This trivially follows from continuity and linearity of $F$.
\item For $\left(f_{i}\right)_{i\in I}\in\vernal$, we have $\left(f_{i}\right)_{i\in I}=F\left(f_{i}\right)_{i\in I}=\iota\circ F_{0}\left(f_{i}\right)_{i\in I}$
and hence $\left(f_{i}\right)_{i\in I}=F_{0}\left(f_{i}\right)_{i\in I}$,
where—strictly speaking—on the left-hand side, $\left(f_{i}\right)_{i\in I}$
is interpreted as an element of $\ell_{w}^{q}\left(\left[W_{j}\right]_{j\in I}\right)$
and on the right-hand side as an element of $\ell_{w}^{q}\left(\left[V_{i}\right]_{i\in I}\right)$.
In particular, since $W_{j}\leq C\left(\R^{\dimension}\right)$, we
see that each $f_{i}:\R^{\dimension}\to\Compl$ is continuous. Finally,
using boundedness of $F_{0}$, we get
\begin{align*}
\sup_{0<\delta\leq1}\frac{1}{\delta}\left\Vert \left(\osc{\delta\cdot T_{i}^{-T}\left[-1,1\right]^{\dimension}}\left[M_{-b_{i}}f_{i}\right]\right)_{i\in I}\right\Vert _{\ell_{w}^{q}\left(\left[V_{i}\right]_{i\in I}\right)} & \leq\left\Vert \left(f_{i}\right)_{i\in I}\right\Vert _{\ell_{w}^{q}\left(\left[W_{i}\right]_{i\in I}\right)}\\
 & =\left\Vert F_{0}\left(f_{i}\right)_{i\in I}\right\Vert _{\ell_{w}^{q}\left(\left[W_{i}\right]_{i\in I}\right)}\\
 & \leq\vertiii{F_{0}}\cdot\left\Vert \left(f_{i}\right)_{i\in I}\right\Vert _{\ell_{w}^{q}\left(\left[V_{i}\right]_{i\in I}\right)},
\end{align*}
which easily yields the claim.\qedhere
\end{enumerate}
\end{proof}
Given all of these preparations, we can finally show that we obtain
Banach frames for $\DecompSp{\CalQ}p{\ell_{w}^{q}}v$ in the expected
way:
\begin{thm}
\label{thm:DiscreteBanachFrameTheorem}Assume that $\Gamma=\left(\gamma_{i}\right)_{i\in I}$
satisfies Assumptions \ref{assu:DiscreteBanachFrameAssumptions} and
\ref{assu:GammaCoversOrbit}. Then there is some $\delta_{0}>0$ such
that for every $0<\delta\leq\delta_{0}$, the family $\left(L_{\delta\cdot T_{i}^{-T}k}\widetilde{\gamma^{\left[i\right]}}\right)_{i\in I,k\in\Z^{\dimension}}$
forms a Banach frame for $\DecompSp{\CalQ}p{\ell_{w}^{q}}v$, with
$\widetilde{\gamma^{\left[i\right]}}\left(x\right)=\gamma^{\left[i\right]}\left(-x\right)$
and
\[
\gamma^{\left[i\right]}=\left|\det T_{i}\right|^{1/2}\cdot M_{b_{i}}\left[\gamma_{i}\circ T_{i}^{T}\right]\qquad\forall i\in I.
\]
In fact, one can choose $\delta_{0}=\frac{1}{1+2\vertiii{F_{0}}^{2}}$,
with $F_{0}$ as in Lemma \ref{lem:SpecialProjection}.

Precisely, the Banach frame property has to be understood as follows:

\begin{itemize}
\item The \textbf{analysis operator}
\begin{alignat*}{2}
A_{\delta}:\: & \DecompSp{\CalQ}p{\ell_{w}^{q}}v\to\ell_{\left(\left|\det T_{i}\right|^{\frac{1}{2}-\frac{1}{p}}\cdot w_{i}\right)_{i\in I}}^{q}\!\!\!\!\!\left(\left[\vphantom{\sum}\smash{C_{i}^{\left(\delta\right)}}\right]_{i\in I}\right), & f\mapsto\left(\left[\gamma^{\left[i\right]}\ast f\right]\left(\delta\cdot T_{i}^{-T}k\right)\right)_{k\in\Z^{\dimension},i\in I}
\end{alignat*}
is well-defined and bounded for each $\delta\in\left(0,1\right]$.
\item As long as $0<\delta\leq\delta_{0}$, there is a bounded linear \textbf{reconstruction
operator}
\[
R_{\delta}:\ell_{\left(\left|\det T_{i}\right|^{\frac{1}{2}-\frac{1}{p}}\cdot w_{i}\right)_{i\in I}}^{q}\!\!\!\!\!\left(\left[\vphantom{\sum}\smash{C_{i}^{\left(\delta\right)}}\right]_{i\in I}\right)\to\DecompSp{\CalQ}p{\ell_{w}^{q}}v
\]
satisfying $R_{\delta}\circ A_{\delta}=\identity_{\DecompSp{\CalQ}p{\ell_{w}^{q}}v}$.
\end{itemize}
Finally, we also have the following \textbf{consistency property}:
If $p_{1},p_{2},q_{1},q_{2}\in\left(0,\infty\right]$, if $w^{\left(1\right)}=\left(\smash{w_{i}^{\left(1\right)}}\right)_{i\in I}$
and $w^{\left(2\right)}=\left(\smash{w_{i}^{\left(2\right)}}\right)_{i\in I}$
are $\CalQ$-moderate weights and if $v^{\left(1\right)},v^{\left(2\right)}:\R^{\dimension}\to\Compl$
are weights such that the assumptions of the current theorem are satisfied
for $\DecompSp{\CalQ}{p_{i}}{\ell_{w^{\left(i\right)}}^{q_{i}}}{v^{\left(i\right)}}$
for $i\in\left\{ 1,2\right\} $ and if $0<\delta\leq\min\left\{ \delta_{1},\delta_{2}\right\} $,
where the constant $\delta_{i}$ is equal to the constant $\delta_{0}$
for the choices $p=p_{i},q=q_{i}$, $w=w^{\left(i\right)}$ and $v=v^{\left(i\right)}$,
then we have
\[
\forall f\!\in\DecompSp{\CalQ}{p_{2}}{\ell_{w^{\left(2\right)}}^{q_{2}}}{v^{\left(2\right)}}:\,f\!\in\DecompSp{\CalQ}{p_{1}}{\ell_{w^{\left(1\right)}}^{q_{1}}}{v^{\left(1\right)}}\Longleftrightarrow\left[\!\left(\gamma^{\left[j\right]}\!\ast f\right)\!\!\left(\delta\cdot T_{j}^{-T}k\right)\!\right]_{k\in\Z^{\dimension},j\in I}\!\in\!\ell_{\left(\left|\det T_{j}\right|^{\frac{1}{2}-\frac{1}{p_{1}}}w_{j}^{\left(1\right)}\right)_{j\in I}}^{q_{1}}\!\!\!\!\!\!\!\!\!\!\left(\!\left[C_{j}^{\left(1,\delta\right)}\right]_{j\in I}\!\right)\!,
\]
with $C_{j}^{\left(1,\delta\right)}=\ell_{\left(v^{\left(1\right)}\right)^{\left(j,\delta\right)}}^{p}\left(\Z^{\dimension}\right)$
and $\left(v^{\left(1\right)}\right)_{k}^{\left(j,\delta\right)}=v^{\left(1\right)}\!\left(\delta\cdot T_{j}^{-T}k\right)$
for $j\in I$ and $k\in\Z^{\dimension}$. 
\end{thm}
\begin{rem*}

\begin{itemize}[leftmargin=0.4cm]
\item The statement of the theorem that the family $\left(L_{\delta\cdot T_{i}^{-T}k}\widetilde{\gamma^{\left[i\right]}}\right)_{i\in I,k\in\Z^{\dimension}}$
forms a Banach frame for the decomposition space $\DecompSp{\CalQ}p{\ell_{w}^{q}}v$
has to be taken with a grain of salt (i.e., as saying that $A_{\delta},R_{\delta}$
as in the statement of the theorem are bounded and $R_{\delta}\circ A_{\delta}=\identity_{\DecompSp{\CalQ}p{\ell_{w}^{q}}v}$).
But if we have $\CalO=\R^{\dimension}$, $\gamma_{i}\in\Schwartz\left(\R^{\dimension}\right)$
for all $i\in I$ and $\DecompSp{\CalQ}p{\ell_{w}^{q}}v\hookrightarrow\Schwartz'\left(\R^{\dimension}\right)$,
then this statement can be taken literally: As seen in the remark
after Theorem \ref{thm:ConvolvingDecompositionSpaceWithGammaJ}, the
definition of $\gamma^{\left[i\right]}\ast f$ given there coincides
with the usual interpretation for $f\in\Schwartz'\left(\R^{\dimension}\right)\supset\DecompSp{\CalQ}p{\ell_{w}^{q}}v$,
so that we indeed have
\[
A_{\delta}f=\left(\left(\gamma^{\left[i\right]}\ast f\right)\left(\delta\cdot T_{i}^{-T}k\right)\right)_{k\in\Z^{\dimension},i\in I}=\left(\left\langle f,\,L_{\delta\cdot T_{i}^{-T}k}\widetilde{\gamma^{\left[i\right]}}\right\rangle _{\Schwartz',\Schwartz}\right)_{k\in\Z^{\dimension},i\in I}.
\]
\item For the consistency statement, note that we only claim that an equivalence
of the form
\[
f\!\in\DecompSp{\CalQ}{p_{1}}{\ell_{w^{\left(1\right)}}^{q_{1}}}{v^{\left(1\right)}}\Longleftrightarrow\left[\left(\gamma^{\left[j\right]}\!\ast f\right)\!\!\left(\delta\cdot T_{j}^{-T}k\right)\right]_{k\in\Z^{\dimension},j\in I}\!\in\!\ell_{\left(\left|\det T_{j}\right|^{\frac{1}{2}-\frac{1}{p_{1}}}w_{j}^{\left(1\right)}\right)_{j\in I}}^{q_{1}}\!\!\!\!\!\!\!\!\!\!\left(\!\left[C_{j}^{\left(1,\delta\right)}\right]_{j\in I}\!\right)
\]
holds under the \emph{assumption} that we \emph{already know} $f\in\DecompSp{\CalQ}{p_{2}}{\ell_{w^{\left(2\right)}}^{q_{2}}}{v^{\left(2\right)}}$
for suitable $p_{2},q_{2},v^{\left(2\right)},w^{\left(2\right)}$.
In other words, we require that we already know that $f$ has a certain
\emph{minimal amount of regularity}. This is quite natural, since
for an arbitrary $f\in Z'\left(\CalO\right)$, there is no reason
why $\gamma^{\left[j\right]}\ast f$ should be defined at all.
\item As the proof will show, the action of $R_{\delta}$ on a given sequence
$\left(\smash{c_{k}^{\left(i\right)}}\right)_{i\in I,k\in\Z^{\dimension}}\in\ell_{\left(\left|\det T_{i}\right|^{\frac{1}{2}-\frac{1}{p}}\cdot w_{i}\right)_{i\in I}}^{q}\!\!\!\!\!\!\!\!\left(\left[\smash{C_{i}^{\left(\delta\right)}}\right]_{i\in I}\right)$
is actually \emph{independent} of $p,q,v,w$. The only thing which
depends on these quantities is $\delta_{0}$, so that $R_{\delta}\left(\smash{c_{k}^{\left(i\right)}}\right)_{i\in I,k\in\Z^{\dimension}}$
is only defined for $0<\delta\leq\delta_{0}=\delta_{0}\left(p,q,v,w,\gamma\right)$.
But once this is satisfied, the definition is independent of $p,q,v,w$.\qedhere
\end{itemize}
\end{rem*}
\begin{proof}
First of all, we remark that the $L^{2}$-normalized functions $\gamma^{\left[i\right]}$
yield a nice statement of the theorem, while the proof can be formulated
easier in terms of the $L^{1}$-normalized functions $\gamma^{\left(i\right)}$.
Hence, we introduce the isometric isomorphism
\[
J:\ell_{\left(\left|\det T_{i}\right|^{\frac{1}{2}-\frac{1}{p}}\cdot w_{i}\right)_{i\in I}}^{q}\!\!\!\!\!\!\!\!\!\left(\!\left[C_{i}^{\left(\delta\right)}\right]_{i\in I}\!\right)\to\ell_{\left(\left|\det T_{i}\right|^{-1/p}\cdot w_{i}\right)_{i\in I}}^{q}\!\!\left(\!\left[C_{i}^{\left(\delta\right)}\right]_{i\in I}\!\right),\left(\smash{c_{k}^{\left(i\right)}}\right)_{k\in\Z^{\dimension},i\in I}\mapsto\left(\left|\det T_{i}\right|^{1/2}\cdot c_{k}^{\left(i\right)}\right)_{k\in\Z^{\dimension},i\in I}.
\]
Then, we define $A_{\delta}^{\left(0\right)}:=J\circ A_{\delta}$
and note
\[
A_{\delta}^{\left(0\right)}f=\left(\left(\gamma^{\left(i\right)}\ast f\right)\left(\delta\cdot T_{i}^{-T}k\right)\right)_{k\in\Z^{\dimension},i\in I}\qquad\forall f\in\DecompSp{\CalQ}p{\ell_{w}^{q}}v,
\]
so it suffices to show that $A_{\delta}^{\left(0\right)}:\DecompSp{\CalQ}p{\ell_{w}^{q}}v\to\ell_{\left(\left|\det T_{i}\right|^{-1/p}\cdot w_{i}\right)_{i\in I}}^{q}\!\!\!\!\!\!\left(\left[\smash{C_{i}^{\left(\delta\right)}}\right]_{i\in I}\right)$
is well-defined and bounded. Further, if there is a bounded operator
$R_{\delta}^{\left(0\right)}:\ell_{\left(\left|\det T_{i}\right|^{-1/p}\cdot w_{i}\right)_{i\in I}}^{q}\!\!\!\!\!\!\left(\left[\smash{C_{i}^{\left(\delta\right)}}\right]_{i\in I}\right)\to\DecompSp{\CalQ}p{\ell_{w}^{q}}v$
satisfying $R_{\delta}^{\left(0\right)}\circ A_{\delta}^{\left(0\right)}=\identity_{\DecompSp{\CalQ}p{\ell_{w}^{q}}v}$,
then a suitable definition of the reconstruction operator $R_{\delta}$
in the statement of the theorem is given by $R_{\delta}:=R_{\delta}^{\left(0\right)}\circ J$,
because of $R_{\delta}\circ A_{\delta}=R_{\delta}^{\left(0\right)}\circ J\circ J^{-1}\circ A_{\delta}^{\left(0\right)}=R_{\delta}^{\left(0\right)}\circ A_{\delta}^{\left(0\right)}$.

These considerations also apply to the consistency statement at the
end of the theorem. All in all, we can thus replace $\gamma^{\left[i\right]}$
by $\gamma^{\left(i\right)}$ in the proof, as long as we replace
all occurrences of $\smash{\ell_{\left(\left|\det T_{i}\right|^{\frac{1}{2}-\frac{1}{p}}\cdot w_{i}\right)_{i\in I}}^{q}}\!\!\!\!\!\!\left(\left[\smash{C_{i}^{\left(\delta\right)}}\right]_{i\in I}\right)$
by $\ell_{\left(\left|\det T_{i}\right|^{-1/p}\cdot w_{i}\right)_{i\in I}}^{q}\!\!\!\!\!\!\left(\left[\smash{C_{i}^{\left(\delta\right)}}\right]_{i\in I}\right)$.

\medskip{}

In the whole proof, we will use the nomenclature introduced in Lemma
\ref{lem:SpecialProjection}. As noted in that lemma, every function
$f_{i}$ is continuous if $\left(f_{i}\right)_{i\in I}\in\vernal$.
Hence, for each $i\in I$, the operator
\[
{\rm Samp}_{\delta,i}:\vernal\to\Compl^{\Z^{\dimension}},\left(f_{j}\right)_{j\in I}\mapsto\left[f_{i}\left(\delta\cdot T_{i}^{-T}k\right)\right]_{k\in\Z^{\dimension}}
\]
is well-defined. Now, note with ${\rm Synth}_{\delta,i}$ as in Lemma
\ref{lem:DiscreteSynthesisOperatorIsAlmostIsometric} that
\begin{align*}
 & \left|f_{i}\left(x\right)-\left[\left({\rm Synth}_{\delta,i}\circ{\rm Samp}_{\delta,i}\right)\left(f_{j}\right)_{j\in I}\right]\left(x\right)\right|\\
 & =\left|f_{i}\left(x\right)-\left(M_{b_{i}}\left[\sum_{k\in\Z^{\dimension}}f_{i}\left(\delta\cdot T_{i}^{-T}k\right)\cdot e^{-2\pi i\left\langle b_{i},\delta\cdot T_{i}^{-T}k\right\rangle }\Indicator_{\delta\cdot T_{i}^{-T}\left(k+\left[0,1\right)^{\dimension}\right)}\right]\right)\left(x\right)\right|\\
 & =\left|\left(M_{-b_{i}}f_{i}\right)\left(x\right)-\sum_{k\in\Z^{\dimension}}\left(M_{-b_{i}}f_{i}\right)\left(\delta\cdot T_{i}^{-T}k\right)\cdot\Indicator_{\delta\cdot T_{i}^{-T}\left(k+\left[0,1\right)^{\dimension}\right)}\left(x\right)\right|\\
\left({\scriptstyle \R^{\dimension}=\biguplus_{k\in\Z^{\dimension}}\delta T_{i}^{-T}\left(k+\left[0,1\right)^{\dimension}\right)}\right) & =\left|\sum_{k\in\Z^{\dimension}}\Indicator_{\delta\cdot T_{i}^{-T}\left(k+\left[0,1\right)^{\dimension}\right)}\left(x\right)\cdot\left[\left(M_{-b_{i}}f_{i}\right)\left(x\right)-\left(M_{-b_{i}}f_{i}\right)\left(\delta\cdot T_{i}^{-T}k\right)\right]\right|\\
 & \leq\sum_{k\in\Z^{\dimension}}\Indicator_{\delta\cdot T_{i}^{-T}\left(k+\left[0,1\right)^{\dimension}\right)}\left(x\right)\cdot\left|\left(M_{-b_{i}}f_{i}\right)\left(x\right)-\left(M_{-b_{i}}f_{i}\right)\left(\delta\cdot T_{i}^{-T}k\right)\right|.
\end{align*}
Now, note that $\Indicator_{\delta\cdot T_{i}^{-T}\left(k+\left[0,1\right)^{\dimension}\right)}\left(x\right)\neq0$
implies $\delta\cdot T_{i}^{-T}k\in x-\delta T_{i}^{-T}\left[0,1\right)^{\dimension}\subset x+\delta T_{i}^{-T}\left[-1,1\right]^{\dimension}$.
Since we trivially have $x\in x+\delta T_{i}^{-T}\left[-1,1\right]^{\dimension}$,
we obtain
\[
\left|\left(M_{-b_{i}}f_{i}\right)\left(x\right)-\left(M_{-b_{i}}f_{i}\right)\left(\delta\cdot T_{i}^{-T}k\right)\right|\leq\left(\osc{\delta T_{i}^{-T}\left[-1,1\right]^{\dimension}}\left[M_{-b_{i}}f_{i}\right]\right)\left(x\right).
\]
Using again that $\R^{\dimension}=\biguplus_{k\in\Z^{\dimension}}\delta T_{i}^{-T}\left(k+\left[0,1\right)^{\dimension}\right)$,
we conclude
\[
\left|f_{i}\left(x\right)-\left[\left({\rm Synth}_{\delta,i}\circ{\rm Samp}_{\delta,i}\right)\left(f_{j}\right)_{j\in I}\right]\left(x\right)\right|\leq\left(\osc{\delta T_{i}^{-T}\left[-1,1\right]^{\dimension}}\left[M_{-b_{i}}f_{i}\right]\right)\left(x\right)\qquad\forall i\in I\;\forall x\in\R^{\dimension}\;\forall\left(f_{j}\right)_{j\in I}\in\vernal.
\]

Consequently, using the solidity of $V_{i}$, we get for
\begin{align*}
{\rm Samp}_{\delta} & :=\prod_{i\in I}{\rm Samp}_{\delta,i}:\vernal\to\left(\Compl^{\Z^{\dimension}}\right)^{I},\left(f_{j}\right)_{j\in I}\mapsto\left({\rm Samp}_{\delta,i}\left(f_{j}\right)_{j\in I}\right)_{i\in I},\\
{\rm Synth}_{\delta} & :=\bigotimes_{i\in I}{\rm Synth}_{\delta,i}:\left(\smash{\Compl^{\Z^{\dimension}}}\right)^{I}\!\to\left\{ \left(f_{i}\right)_{i\in I}\with f_{i}:\R^{\dimension}\to\Compl\text{ measurable }\forall i\in I\right\} ,\left(\smash{c_{k}^{\left(i\right)}}\right)_{k\in\Z^{\dimension},i\in I}\mapsto\left({\rm Synth}_{\delta,i}\left(\smash{c_{k}^{\left(i\right)}}\right)_{k\in\Z^{\dimension}}\right)_{i\in I}
\end{align*}
that
\begin{align}
\left\Vert \left(f_{i}\right)_{i\in I}-\left({\rm Synth}_{\delta}\circ{\rm Samp}_{\delta}\right)\left(f_{i}\right)_{i\in I}\right\Vert _{\ell_{w}^{q}\left(\left[V_{i}\right]_{i\in I}\right)} & =\left\Vert \left(\left\Vert f_{i}-{\rm Synth}_{\delta,i}\circ{\rm Samp}_{\delta,i}\left(f_{j}\right)_{j\in I}\right\Vert _{V_{i}}\right)_{i\in I}\right\Vert _{\ell_{w}^{q}}\nonumber \\
 & \leq\left\Vert \left(\left\Vert \osc{\delta T_{i}^{-T}\left[-1,1\right]^{\dimension}}\left[M_{-b_{i}}f_{i}\right]\right\Vert _{V_{i}}\right)_{i\in I}\right\Vert _{\ell_{w}^{q}}\nonumber \\
\left({\scriptstyle \text{eq. }\eqref{eq:OscillationEstimateOnVernal}}\right) & \leq\vertiii{F_{0}}\cdot\delta\cdot\left\Vert \left(f_{i}\right)_{i\in I}\right\Vert _{\ell_{w}^{q}\left(\left[V_{i}\right]_{i\in I}\right)}\quad\forall\left(f_{i}\right)_{i\in I}\in\vernal\:\forall\delta\in\left(0,1\right].\label{eq:UnsmoothedDiscretizationIsClose}
\end{align}
Using the (quasi)-triangle inequality for $\ell_{w}^{q}\left(\left[V_{i}\right]_{i\in I}\right)$
(where $2^{\frac{1}{p}+\frac{1}{q}}$ is a valid triangle constant,
thanks to \cite[Exercise 1.1.5(c)]{GrafakosClassicalFourierAnalysis}
and to (the proof of) Lemma \ref{lem:IteratedSequenceSpaceComplete}),
we conclude that
\[
T_{0}^{\left(\delta\right)}:={\rm Synth}_{\delta}\circ{\rm Samp}_{\delta}:\vernal\to\ell_{w}^{q}\left(\left[V_{i}\right]_{i\in I}\right)
\]
is well-defined and bounded, with $\vertiii{T_{0}^{\left(\delta\right)}}\leq2^{\frac{1}{p}+\frac{1}{q}}\cdot\left(1+\vertiii{F_{0}}\delta\right)\leq2^{\frac{1}{p}+\frac{1}{q}}\left(1+\vertiii{F_{0}}\right)$
for all $\delta\in\left(0,1\right]$.

\medskip{}

Boundedness of $T_{0}^{\left(\delta\right)}$—together with estimate
(\ref{eq:UnsmoothedDiscretizationIsClose})—is almost sufficient for
our purposes, but not quite: In general, it need not be the case that
$T_{0}^{\left(\delta\right)}$ maps $\vernal$ into $\vernal$. But
since Lemma \ref{lem:SpecialProjection} shows $F\circ F=F$, it is
easy to see $F:\ell_{w}^{q}\left(\left[V_{i}\right]_{i\in I}\right)\to\vernal$,
so that $T^{\left(\delta\right)}:=F\circ T_{0}^{\left(\delta\right)}:\vernal\to\vernal$.
Furthermore, since $F|_{\vernal}=\identity_{\vernal}$, we get
\begin{align}
\left\Vert \left(f_{i}\right)_{i\in I}-T^{\left(\delta\right)}\left(f_{i}\right)_{i\in I}\right\Vert _{\ell_{w}^{q}\left(\left[V_{i}\right]_{i\in I}\right)} & =\left\Vert F\left(f_{i}\right)_{i\in I}-FT_{0}^{\left(\delta\right)}\left(f_{i}\right)_{i\in I}\right\Vert _{\ell_{w}^{q}\left(\left[V_{i}\right]_{i\in I}\right)}\nonumber \\
 & \leq\vertiii F\cdot\left\Vert \left(f_{i}\right)_{i\in I}-T_{0}^{\left(\delta\right)}\left(f_{i}\right)_{i\in I}\right\Vert _{\ell_{w}^{q}\left(\left[V_{i}\right]_{i\in I}\right)}\nonumber \\
\left({\scriptstyle \text{eq. }\eqref{eq:UnsmoothedDiscretizationIsClose}}\right) & \leq\vertiii{F_{0}}^{2}\cdot\delta\cdot\left\Vert \left(f_{i}\right)_{i\in I}\right\Vert _{\ell_{w}^{q}\left(\left[V_{i}\right]_{i\in I}\right)}\qquad\forall\left(f_{i}\right)_{i\in I}\in\vernal\:\forall\delta\in\left(0,1\right].\label{eq:SmoothedDiscretizationIsClose}
\end{align}
But for $0<\delta\leq\delta_{0}=\frac{1}{1+2\vertiii{F_{0}}^{2}}$,
we have $\vertiii{F_{0}}^{2}\cdot\delta\leq\frac{1}{2}$ and hence
$\vertiii{\identity_{\vernal}-T^{\left(\delta\right)}}\leq\frac{1}{2}$.
Using a Neumann-series argument (which is also valid for Quasi-Banach
spaces, cf.\@ e.g.\@ \cite[Lemma 2.4.11]{VoigtlaenderPhDThesis}
and thus for the closed subspace $\vernal$ of the Quasi-Banach space
$\ell_{w}^{q}\left(\left[V_{i}\right]_{i\in I}\right)$ thanks to
Lemmas \ref{lem:IteratedSequenceSpaceComplete} and \ref{lem:SpecialProjection}),
we conclude that $T^{\left(\delta\right)}:\vernal\to\vernal$ is boundedly
invertible, as long as $0<\delta\leq\delta_{0}$.

\medskip{}

Now, for arbitrary $f\in\DecompSp{\CalQ}p{\ell_{w}^{q}}v$, Lemma
\ref{lem:SpecialProjection} shows that $\left({\rm Synth}_{\CalD}\circ m_{\theta}\circ\iota\circ{\rm Ana}_{{\rm osc}}\right)f=f$.
The same lemma also shows $F\left[\left(\gamma^{\left(j\right)}\ast f\right)_{j\in I}\right]=\left(\gamma^{\left(j\right)}\ast f\right)_{j\in I}$,
i.e., $\left(\iota\circ{\rm Ana}_{{\rm osc}}\right)f=\left(\gamma^{\left(j\right)}\ast f\right)_{j\in I}\in\vernal$.
Hence,
\begin{align*}
f & =\left({\rm Synth}_{\CalD}\circ m_{\theta}\circ\iota\circ{\rm Ana}_{{\rm osc}}\right)f\\
 & =\left[\left({\rm Synth}_{\CalD}\circ m_{\theta}\right)\circ\left(T^{\left(\delta\right)}\right)^{-1}\circ T^{\left(\delta\right)}\circ\iota\circ{\rm Ana}_{{\rm osc}}\right]f\\
\left({\scriptstyle \text{def. of }T^{\left(\delta\right)}}\right) & =\left[\left(\left[{\rm Synth}_{\CalD}\circ m_{\theta}\right]\circ\left(T^{\left(\delta\right)}\right)^{-1}\circ F\circ{\rm Synth}_{\delta}\right)\circ{\rm Samp}_{\delta}\circ\iota\circ{\rm Ana}_{{\rm osc}}\right]f.
\end{align*}
Now, note
\[
\left[\left(\left[{\rm Samp}_{\delta}\circ\iota\circ{\rm Ana}_{{\rm osc}}\right]f\right)_{i}\right]_{k}=\left(\left[{\rm Samp}_{\delta}\left(\gamma^{\left(j\right)}\ast f\right)_{j\in I}\right]_{i}\right)_{k}=\left(\gamma^{\left(i\right)}\ast f\right)\left(\delta\cdot T_{i}^{-T}k\right)
\]
and hence ${\rm Samp}_{\delta}\circ\iota\circ{\rm Ana}_{{\rm osc}}=A_{\delta}^{\left(0\right)}$
on $\DecompSp{\CalQ}p{\ell_{w}^{q}}v$. Thus, if we define $R_{\delta}^{\left(0\right)}:=\left[{\rm Synth}_{\CalD}\circ m_{\theta}\right]\circ\left(T^{\left(\delta\right)}\right)^{-1}\circ F\circ{\rm Synth}_{\delta}$,
we have shown $\identity_{\DecompSp{\CalQ}p{\ell_{w}^{q}}v}=R_{\delta}^{\left(0\right)}\circ A_{\delta}^{\left(0\right)}$,
as claimed. All that remains to show is that $R_{\delta}^{\left(0\right)},A_{\delta}^{\left(0\right)}$
are indeed well-defined and bounded with domains and codomains as
stated at the beginning of the proof.

\medskip{}

To this end, note that Lemma \ref{lem:DiscreteSynthesisOperatorIsAlmostIsometric}
easily implies that ${\rm Synth}_{\delta}:\ell_{\left(\left|\det T_{i}\right|^{-1/p}\cdot w_{i}\right)_{i\in I}}^{q}\!\!\!\!\!\left(\left[\smash{C_{i}^{\left(\delta\right)}}\right]_{i\in I}\right)\to\ell_{w}^{q}\left(\left[V_{i}\right]_{i\in I}\right)$
is well-defined and bounded. In fact, the lemma even shows that 
\[
\left(\smash{c_{k}^{\left(i\right)}}\right)_{k\in\Z^{\dimension},i\in I}\in\ell_{\left(\left|\det T_{i}\right|^{-1/p}\cdot w_{i}\right)_{i\in I}}^{q}\!\!\left(\left[\vphantom{F}\smash{C_{i}^{\left(\delta\right)}}\right]_{i\in I}\right)\Longleftrightarrow{\rm Synth}_{\delta}\left(\smash{c_{k}^{\left(i\right)}}\right)_{k\in\Z^{\dimension},i\in I}\in\ell_{w}^{q}\left(\left[V_{i}\right]_{i\in I}\right)
\]
and
\[
\left\Vert {\rm Synth}_{\delta}\left(\smash{c_{k}^{\left(i\right)}}\right)_{k\in\Z^{\dimension},i\in I}\right\Vert _{\ell_{w}^{q}\left(\left[V_{i}\right]_{i\in I}\right)}\asymp\left\Vert \left(\smash{c_{k}^{\left(i\right)}}\right)_{k\in\Z^{\dimension},i\in I}\right\Vert _{\ell_{\left(\left|\det T_{i}\right|^{-1/p}\cdot w_{i}\right)_{i\in I}}^{q}\!\!\!\!\!\left(\left[\vphantom{F}\smash{C_{i}^{\left(\delta\right)}}\right]_{i\in I}\right)},
\]
where the implied constant may depend on $\delta$. Consequently,
$R_{\delta}^{\left(0\right)}:\ell_{\left(\left|\det T_{i}\right|^{-1/p}\cdot w_{i}\right)_{i\in I}}^{q}\!\!\!\left(\left[\vphantom{F}\smash{C_{i}^{\left(\delta\right)}}\right]_{i\in I}\right)\to\DecompSp{\CalQ}p{\ell_{w}^{q}}v$
is indeed well-defined and bounded for $0<\delta\leq\delta_{0}$.
Furthermore, we see (now for arbitrary $\delta\in\left(0,1\right]$)
that
\begin{align*}
\left\Vert \smash{A_{\delta}^{\left(0\right)}}f\right\Vert _{\ell_{\left(\left|\det T_{i}\right|^{-1/p}\cdot w_{i}\right)_{i\in I}}^{q}\!\!\!\left(\left[\vphantom{F}\smash{C_{i}^{\left(\delta\right)}}\right]_{i\in I}\right)} & =\left\Vert \left({\rm Samp}_{\delta}\circ\iota\circ{\rm Ana}_{{\rm osc}}\right)f\right\Vert _{\ell_{\left(\left|\det T_{i}\right|^{-1/p}\cdot w_{i}\right)_{i\in I}}^{q}\!\!\!\left(\left[\vphantom{F}\smash{C_{i}^{\left(\delta\right)}}\right]_{i\in I}\right)}\\
 & \asymp_{\delta}\:\left\Vert \left({\rm Synth}_{\delta}\circ{\rm Samp}_{\delta}\circ\iota\circ{\rm Ana}_{{\rm osc}}\right)f\right\Vert _{\ell_{w}^{q}\left(\left[V_{i}\right]_{i\in I}\right)}\\
 & =\left\Vert \left(\smash{T_{0}^{\left(\delta\right)}}\circ\iota\circ{\rm Ana}_{{\rm osc}}\right)f\right\Vert _{\ell_{w}^{q}\left(\left[V_{i}\right]_{i\in I}\right)}\\
\left({\scriptstyle \iota\circ{\rm Ana}_{{\rm osc}}:\DecompSp{\CalQ}p{\ell_{w}^{q}}v\to\vernal,\text{ as seen above}}\right) & \lesssim\vertiii{\smash{T_{0}^{\left(\delta\right)}}}_{\vernal\to\ell_{w}^{q}\left(\left[V_{i}\right]_{i\in I}\right)}\cdot\vertiii{\iota\circ{\rm Ana}_{{\rm osc}}}_{\DecompSp{\CalQ}p{\ell_{w}^{q}}v\to\ell_{w}^{q}\left(\left[V_{i}\right]_{i\in I}\right)}\cdot\left\Vert f\right\Vert _{\DecompSp{\CalQ}p{\ell_{w}^{q}}v}<\infty
\end{align*}
for all $f\in\DecompSp{\CalQ}p{\ell_{w}^{q}}v$. This finally shows
that $A_{\delta}^{\left(0\right)}:\DecompSp{\CalQ}p{\ell_{w}^{q}}v\to\ell_{\left(\left|\det T_{i}\right|^{-1/p}\cdot w_{i}\right)_{i\in I}}^{q}\!\!\!\left(\left[\vphantom{F}\smash{C_{i}^{\left(\delta\right)}}\right]_{i\in I}\right)$
is well-defined and bounded for each $\delta\in\left(0,1\right]$
and thus completes the proof of the Banach frame property.

\medskip{}

It remains to verify the consistency property stated above. To this
end, first define 
\[
V_{j}^{\left(i\right)}:=\begin{cases}
L_{v^{\left(i\right)}}^{p_{i}}\left(\smash{\R^{\dimension}}\right), & \text{if }p_{i}\in\left[1,\infty\right],\\
W_{T_{j}^{-T}\left[-1,1\right]^{\dimension}}\left(L_{v^{\left(i\right)}}^{p_{i}}\right), & \text{if }p_{i}\in\left(0,1\right),
\end{cases}
\]
as well as $C_{j}^{\left(i,\delta\right)}=\ell_{\left(v^{\left(i\right)}\right)^{\left(j,\delta\right)}}^{p_{i}}\left(\Z^{\dimension}\right)$
with $\left(v^{\left(i\right)}\right)_{k}^{\left(j,\delta\right)}=v^{\left(i\right)}\left(\delta\cdot T_{j}^{-T}k\right)$
for $k\in\Z^{\dimension}$, $i\in\left\{ 1,2\right\} $ and $j\in I$.
Next, we observe that the domain and codomain of the reconstruction/analysis
operators 
\[
R_{\delta}^{\left(0,i\right)}:\ell_{\left(\left|\det T_{j}\right|^{-1/p_{i}}\cdot w_{j}^{\left(i\right)}\right)_{j\in I}}^{q_{i}}\!\!\!\left(\left[\vphantom{F}\smash{C_{j}^{\left(i,\delta\right)}}\right]_{j\in I}\right)\to\DecompSp{\CalQ}{p_{i}}{\ell_{w^{\left(i\right)}}^{q_{i}}}{v^{\left(i\right)}}
\]
and
\[
A_{\delta}^{\left(0,i\right)}:\DecompSp{\CalQ}{p_{i}}{\ell_{w^{\left(i\right)}}^{q_{i}}}{v^{\left(i\right)}}\to\ell_{\left(\left|\det T_{j}\right|^{-1/p_{i}}\cdot w_{j}^{\left(i\right)}\right)_{j\in I}}^{q_{i}}\left(\left[\vphantom{F}\smash{C_{j}^{\left(i,\delta\right)}}\right]_{j\in I}\right)
\]
do depend on $i\in\left\{ 1,2\right\} $, but the actual \emph{action}
of these mappings do not: We always have
\[
A_{\delta}^{\left(0,i\right)}f=\left[\left(\gamma^{\left(j\right)}\ast f\right)\left(\delta\cdot T_{j}^{-T}k\right)\right]_{k\in\Z^{\dimension},j\in I}\quad\overset{\text{eq. }\eqref{eq:SpecialConvolutionPointwiseDefinition}}{=}\quad\left[\sum_{\ell\in I}\Fourier^{-1}\left(\widehat{\gamma^{\left(j\right)}}\cdot\varphi_{\ell}\cdot\widehat{f}\right)\left(\delta\cdot T_{j}^{-T}k\right)\right]_{k\in\Z^{\dimension},j\in I}
\]
and
\begin{align*}
R_{\delta}^{\left(0,i\right)}\left(c_{k}^{\left(j\right)}\right)_{k\in\Z^{\dimension},j\in I} & =\left(\left[{\rm Synth}_{\CalD}\circ m_{\theta}\right]\circ\left(T^{\left(\delta\right)}\right)^{-1}\circ F\circ{\rm Synth}_{\delta}\right)\left(c_{k}^{\left(j\right)}\right)_{k\in\Z^{\dimension},j\in I}\\
 & =\left(\left[{\rm Synth}_{\CalD}\circ m_{\theta}\right]\circ\left(T^{\left(\delta\right)}\right)^{-1}\circ\iota\circ{\rm Ana}_{{\rm osc}}\circ{\rm Synth}_{\CalD}\circ m_{\theta}\circ{\rm Synth}_{\delta}\right)\left(c_{k}^{\left(j\right)}\right)_{k\in\Z^{\dimension},j\in I},
\end{align*}
where
\begin{align*}
{\rm Synth}_{\CalD}\left(f_{j}\right)_{j\in I} & =\sum_{j\in I}\left[\Fourier^{-1}\left(\varphi_{j}\cdot\widehat{f_{j}}\right)\right]\quad\text{ with unconditional convergence in }Z'\left(\CalO\right),\\
m_{\theta}\left(f_{j}\right)_{j\in I} & =\left[\left(\Fourier^{-1}\theta_{j}\right)\ast f_{j}\right]_{j\in I},\\
\iota\left(f_{j}\right)_{j\in I} & =\left(f_{j}\right)_{j\in I},\\
{\rm Ana}_{{\rm osc}}f & =\left(\gamma^{\left(j\right)}\ast f\right)_{j\in I}\quad\text{with }\gamma^{\left(j\right)}\ast f\text{ as in equation }\eqref{eq:SpecialConvolutionPointwiseDefinition},\\
{\rm Synth}_{\delta}\left(\smash{c_{k}^{\left(j\right)}}\right)_{k\in\Z^{\dimension},j\in I} & =\left(M_{b_{j}}\left[\sum_{k\in\Z^{\dimension}}c_{k}^{\left(j\right)}\cdot e^{-2\pi i\left\langle b_{j},\delta\cdot T_{j}^{-T}k\right\rangle }\Indicator_{\delta\cdot T_{j}^{-T}\left(k+\left[0,1\right)^{\dimension}\right)}\right]\right)_{j\in I}
\end{align*}
for all $\left(f_{j}\right)_{j\in I}\in\ell_{w^{\left(i\right)}}^{q_{i}}\!\!\left(\left[\smash{V_{j}^{\left(i\right)}}\right]_{j\in I}\right)$,
all $\left(\smash{c_{k}^{\left(j\right)}}\right)_{k\in\Z^{\dimension},j\in I}\in\ell_{\left(\left|\det T_{j}\right|^{-1/p_{i}}\cdot w_{j}^{\left(i\right)}\right)_{j\in I}}^{q_{i}}\!\!\!\left(\left[\vphantom{F}\smash{C_{j}^{\left(i,\delta\right)}}\right]_{j\in I}\right)$,
and all $f\in\DecompSp{\CalQ}{p_{i}}{\ell_{w^{\left(i\right)}}^{q_{i}}}{v^{\left(i\right)}}$.
Finally, we also have (since $\left(T^{\left(\delta\right)}\right)^{-1}$
can be computed by a Neumann series, as shown above)
\[
\left(T^{\left(\delta\right)}\right)^{-1}\left(f_{j}\right)_{j\in I}=\left(\identity-\left[\identity-T^{\left(\delta\right)}\right]\right)^{-1}\left(f_{j}\right)_{j\in I}=\sum_{n=0}^{\infty}\left(\identity-T^{\left(\delta\right)}\right)^{n}\left(f_{j}\right)_{j\in I},
\]
where
\[
T^{\left(\delta\right)}\left(f_{j}\right)_{j\in I}=\left(F\circ T_{0}^{\left(\delta\right)}\right)\left(f_{j}\right)_{j\in I}=\left(\iota\circ{\rm Ana}_{{\rm osc}}\circ{\rm Synth}_{\CalD}\circ m_{\theta}\right)\circ\left({\rm Synth}_{\delta}\circ{\rm Samp}_{\delta}\right)\left(f_{j}\right)_{j\in I}
\]
for 
\[
\left(f_{j}\right)_{j\in I}\in\vernal_{i}:=\left\{ \left(g_{j}\right)_{j\in I}\in\ell_{w^{\left(i\right)}}^{q_{i}}\left(\left[\vphantom{F}\smash{V_{j}^{\left(i\right)}}\right]_{j\in I}\right)\with F\left(g_{j}\right)_{j\in I}=\left(g_{j}\right)_{j\in I}\right\} .
\]

In summary, we have shown $R_{\delta}^{\left(0,1\right)}\left(\smash{c_{k}^{\left(j\right)}}\right)_{k\in\Z^{\dimension},j\in I}=R_{\delta}^{\left(0,2\right)}\left(\smash{c_{k}^{\left(j\right)}}\right)_{k\in\Z^{\dimension},j\in I}$
and $A_{\delta}^{\left(0,1\right)}f=A_{\delta}^{\left(0,2\right)}f$,
as long as both sides of the respective equations are defined. Now,
let $f\in\DecompSp{\CalQ}{p_{2}}{\ell_{w^{\left(2\right)}}^{q_{2}}}{v^{\left(2\right)}}$
be arbitrary. The implication ``$\Rightarrow$'' of the consistency
statement follows immediately from the main statement of the theorem,
so that we only need to show ``$\Leftarrow$''. Hence, assume 
\[
c:=A_{\delta}^{\left(0,2\right)}f=\left[\left(\gamma^{\left(j\right)}\ast f\right)\left(\delta\cdot T_{j}^{-T}k\right)\right]_{k\in\Z^{\dimension},j\in I}\in\ell_{\left(\left|\det T_{j}\right|^{-1/p_{1}}\cdot w_{j}^{\left(1\right)}\right)_{j\in I}}^{q_{1}}\!\!\!\left(\left[\vphantom{F}\smash{C_{j}^{\left(1,\delta\right)}}\right]_{j\in I}\right).
\]
We know from above that $f=R_{\delta}^{\left(0,2\right)}A_{\delta}^{\left(0,2\right)}f=R_{\delta}^{\left(0,2\right)}c$.
But we have $c\in\ell_{\left(\left|\det T_{j}\right|^{-1/p_{i}}\cdot w_{j}^{\left(i\right)}\right)_{j\in I}}^{q_{i}}\!\!\left(\left[\vphantom{F}\smash{C_{j}^{\left(i,\delta\right)}}\right]_{j\in I}\right)$
for both $i=1$ and $i=2$, so that we get $f=R_{\delta}^{\left(0,2\right)}c=R_{\delta}^{\left(0,1\right)}c$.
Since
\[
R_{\delta}^{\left(0,1\right)}:\ell_{\left(\left|\det T_{j}\right|^{-1/p_{1}}\cdot w_{j}^{\left(1\right)}\right)_{j\in I}}^{q_{1}}\!\!\!\left(\left[\vphantom{F}\smash{C_{j}^{\left(1,\delta\right)}}\right]_{j\in I}\right)\to\DecompSp{\CalQ}{p_{1}}{\ell_{w^{\left(1\right)}}^{q_{1}}}{v^{\left(1\right)}}
\]
is well-defined and bounded, we get $f\in\DecompSp{\CalQ}{p_{1}}{\ell_{w^{\left(1\right)}}^{q_{1}}}{v^{\left(1\right)}}$,
as claimed.
\end{proof}
The main limitation of Theorem \ref{thm:DiscreteBanachFrameTheorem}
is its somewhat opaque set of assumptions regarding $\Gamma=\left(\gamma_{i}\right)_{i\in I}$.
In Section \ref{sec:SimplifiedCriteria} (see in particular Corollary
\ref{cor:BanachFrameSimplifiedCriteria}), we will derive more transparent
criteria which ensure that Theorem \ref{thm:DiscreteBanachFrameTheorem}
is applicable.

But before that, we first consider the ``dual'' problem to the Banach
frame property, i.e., we show that the family $\left(L_{\delta\cdot T_{i}^{-T}k}\:\gamma^{\left[i\right]}\right)_{k\in\Z^{\dimension},i\in I}$
forms an \emph{atomic decomposition} for the decomposition space $\DecompSp{\CalQ}p{\ell_{w}^{q}}v$,
under suitable assumptions on $\Gamma=\left(\gamma_{i}\right)_{i\in I}$.
Proving this is the main goal of the next section.

\section{Atomic decompositions}

\label{sec:AtomicDecompositions}In this section, we show the dual
statement to the preceding section, i.e., we show that the (discretely
translated) $\gamma^{\left[j\right]}$ not only form a Banach frame
for $\DecompSp{\CalQ}p{\ell_{w}^{q}}v$, but also an atomic decomposition.
For this, we introduce still another set of assumptions:
\begin{assumption}
\label{assu:AtomicDecompositionAssumption}We assume that for each
$i\in I$, we are given functions $\gamma_{i},\gamma_{i,1},\gamma_{i,2}$
with the following properties:

\begin{enumerate}
\item We have $\gamma_{i,1},\gamma_{i,2}\in L_{\left(1+\left|\mybullet\right|\right)^{K}}^{1}\left(\R^{\dimension}\right)\hookrightarrow L_{v_{0}}^{1}\left(\R^{\dimension}\right)\hookrightarrow L^{1}\left(\R^{\dimension}\right)$
for all $i\in I$.
\item We have $\gamma_{i}=\gamma_{i,1}\ast\gamma_{i,2}$ for all $i\in I$.
\item We have $\widehat{\gamma_{i,1}},\widehat{\gamma_{i,2}}\in C^{\infty}\left(\R^{\dimension}\right)$
for all $i\in I$ and all partial derivatives of $\widehat{\gamma_{i,1}},\widehat{\gamma_{i,2}}$
have at most polynomial growth.
\item We have $\gamma_{i,2}\in C^{1}\left(\R^{\dimension}\right)$ with
$\nabla\gamma_{i,2}\in L_{v_{0}}^{1}\left(\R^{\dimension}\right)$
for all $i\in I$.
\item The constant
\begin{equation}
\Omega_{4}^{\left(p,K\right)}:=\sup_{i\in I}\left\Vert \gamma_{i,2}\right\Vert _{K_{0}}+\sup_{i\in I}\left\Vert \nabla\gamma_{i,2}\right\Vert _{K_{0}}\label{eq:AtomicDecompositionGamma2ConstantDefinition}
\end{equation}
is finite. Here, $K_{0}:=K+\frac{\dimension}{\min\left\{ 1,p\right\} }+1$
and
\[
\left\Vert f\right\Vert _{K_{0}}:=\sup_{x\in\R^{\dimension}}\left(1+\left|x\right|\right)^{K_{0}}\left|f\left(x\right)\right|\in\left[0,\infty\right].
\]
\item We have $\left\Vert \gamma_{i}\right\Vert _{K_{0}}<\infty$ for all
$i\in I$.
\item For $\ell\in\left\{ 1,2\right\} $ and $i\in I$, define
\begin{equation}
\gamma_{\ell}^{\left(i\right)}:=\Fourier^{-1}\left(\widehat{\gamma_{i,\ell}}\circ S_{i}^{-1}\right)=\left|\det T_{i}\right|\cdot M_{b_{i}}\left[\gamma_{i,\ell}\circ T_{i}^{T}\right],\label{eq:AtomicDecompositionFamilyDefinition}
\end{equation}
so that $\gamma_{\ell}^{\left(i\right)}$ is to $\gamma_{i,\ell}$
as $\gamma^{\left(i\right)}$ is to $\gamma_{i}$.
\item For $i,j\in I$ set
\begin{equation}
C_{i,j}:=\begin{cases}
\left\Vert \Fourier^{-1}\left(\varphi_{i}\cdot\widehat{\gamma_{1}^{\left(j\right)}}\right)\right\Vert _{L_{v_{0}}^{1}}, & \text{if }p\in\left[1,\infty\right],\\
\left(1+\left\Vert T_{j}^{-1}T_{i}\right\Vert \right)^{pK+\dimension}\cdot\left|\det T_{j}\right|^{1-p}\cdot\left\Vert \Fourier^{-1}\left(\varphi_{i}\cdot\widehat{\gamma_{1}^{\left(j\right)}}\right)\right\Vert _{L_{v_{0}}^{p}}^{p}, & \text{if }p\in\left(0,1\right).
\end{cases}\label{eq:GammaSynthesisMatrixEntries}
\end{equation}
\item With $r:=\max\left\{ q,\frac{q}{p}\right\} $, we assume that the
operator $\overrightarrow{C}$ induced by $\left(C_{i,j}\right)_{i,j\in I}$,
i.e.\@
\[
\overrightarrow{C}\left(c_{j}\right)_{j\in I}:=\left(\sum_{j\in I}C_{i,j}c_{j}\right)_{i\in I}
\]
defines a well-defined, bounded operator $\overrightarrow{C}:\ell_{w^{\min\left\{ 1,p\right\} }}^{r}\left(I\right)\to\ell_{w^{\min\left\{ 1,p\right\} }}^{r}\left(I\right)$.\qedhere
\end{enumerate}
\end{assumption}
\begin{rem*}
In the following, we will often use $\Gamma_{\ell}$ as a short notation
for the family $\Gamma_{\ell}=\left(\gamma_{i,\ell}\right)_{i\in I}$
($\ell\in\left\{ 1,2\right\} $), similar to the notation $\Gamma=\left(\gamma_{i}\right)_{i\in I}$.

The assumptions above are slightly redundant. In particular, since
$v_{0}\left(x\right)\leq\Omega_{1}\cdot\left(1+\left|x\right|\right)^{K}$,
it is an easy consequence of equation (\ref{eq:StandardDecayLpEstimate})
that $\left\Vert \gamma_{i,2}\right\Vert _{K_{0}}<\infty$ and $\left\Vert \nabla\gamma_{i,2}\right\Vert _{K_{0}}<\infty$
already imply $\gamma_{i,2}\in L_{\left(1+\left|\mybullet\right|\right)^{K}}^{1}\left(\R^{\dimension}\right)\hookrightarrow L_{v_{0}}^{1}\left(\R^{\dimension}\right)$
and $\nabla\gamma_{i,2}\in L_{\left(1+\left|\mybullet\right|\right)^{K}}^{1}\left(\R^{\dimension}\right)\hookrightarrow L_{v_{0}}^{1}\left(\R^{\dimension}\right)$,
respectively.

Exactly as in Remark \ref{rem:MainAssumptionsRemark}, we see that
$\gamma_{i,1},\gamma_{i,2}\in L_{\left(1+\left|\mybullet\right|\right)^{K}}^{1}\left(\R^{\dimension}\right)$
entails $\gamma_{1}^{\left(j\right)},\gamma_{2}^{\left(j\right)}\in L_{\left(1+\left|\mybullet\right|\right)^{K}}^{1}\left(\R^{\dimension}\right)\hookrightarrow L_{v_{0}}^{1}\left(\R^{\dimension}\right)$
for all $j\in I$.
\end{rem*}
Part of the definition of an atomic decomposition $\left(\theta_{\ell}\right)_{\ell\in L}$
is that the synthesis map $\left(c_{\ell}\right)_{\ell\in L}\mapsto\sum_{\ell\in L}c_{\ell}\theta_{\ell}$
is bounded, when defined on a suitable sequence space. Our next lemma
establishes a variant of this property for a certain \emph{continuous}
(as opposed to discrete) synthesis operator. This lemma should be
compared to Lemma \ref{lem:DecompositionSynthesis}.
\begin{lem}
\label{lem:GammaSynthesisBounded}Assume that the family $\Gamma_{1}=\left(\gamma_{i,1}\right)_{i\in I}$
satisfies Assumption \ref{assu:AtomicDecompositionAssumption}. Then,
the operator
\[
{\rm Synth}_{\Gamma_{1}}\::\ell_{w}^{q}\left(\left[V_{j}\right]_{j\in I}\right)\to\DecompSp{\CalQ}p{\ell_{w}^{q}}v\leq Z'\left(\CalO\right),\left(g_{j}\right)_{j\in I}\mapsto\sum_{j\in I}\gamma_{1}^{\left(j\right)}\ast g_{j}\:\overset{\text{Lem. }\ref{lem:SpecialConvolutionConsistent}}{=}\:\sum_{j\in I}\Fourier^{-1}\left(\widehat{\gamma_{1}^{\left(j\right)}}\cdot\widehat{g_{j}}\right)
\]
is well-defined and bounded with
\[
\vertiii{{\rm Synth}_{\Gamma_{1}}}\leq C\cdot\vertiii{\smash{\overrightarrow{C}}}^{\max\left\{ 1,\frac{1}{p}\right\} }\;\text{with}\;C=\begin{cases}
1, & \text{if }p\geq1\\
\frac{\left(2^{6}/\sqrt{\dimension}\right)^{\frac{\dimension}{p}}}{2^{21}\cdot\dimension^{7}}\!\cdot\!\left(2^{21}\!\cdot\!\dimension^{5}\!\cdot\!\left\lceil K\!+\!\frac{\dimension+1}{p}\right\rceil \right)^{\!\left\lceil K+\frac{\dimension+1}{p}\right\rceil +1}\!\!\!\cdot\!\left(1\!+\!R_{\CalQ}\right)^{\frac{\dimension}{p}}\!\cdot\!\Omega_{0}^{5K}\Omega_{1}^{5}, & \text{if }p<1.
\end{cases}
\]

Here, ${\rm Synth}_{\Gamma_{1}}\left(g_{j}\right)_{j\in I}$ is the
linear functional
\begin{equation}
Z\left(\CalO\right)\to\Compl,f\mapsto\sum_{j\in I}\left\langle \widehat{\gamma_{1}^{\left(j\right)}}\cdot\widehat{g_{j}},\,\Fourier^{-1}f\right\rangle _{\Schwartz',\Schwartz}=\sum_{j\in I}\left\langle \Fourier^{-1}\left(\widehat{\gamma_{1}^{\left(j\right)}}\cdot\widehat{g_{j}}\right),\,f\right\rangle _{\Schwartz',\Schwartz}\:\overset{\text{Lem. }\ref{lem:SpecialConvolutionConsistent}}{=}\:\sum_{j\in I}\left\langle \gamma_{1}^{\left(j\right)}\ast g_{j},\,f\right\rangle _{\Schwartz',\Schwartz},\label{eq:GammaSynthesisFunctionalExplicitDefinition}
\end{equation}
where each of the series converges absolutely for each $f\in Z\left(\CalO\right)$.
\end{lem}
\begin{proof}
First of all, recall from Lemma \ref{lem:SpecialConvolutionConsistent}
that $V_{j}\hookrightarrow\Schwartz'\left(\R^{\dimension}\right)$
for all $j\in I$. Thus, for $\left(g_{j}\right)_{j\in I}\in\ell_{w}^{q}\left(\left[V_{j}\right]_{j\in I}\right)$,
we see that $\widehat{g_{j}}\in\Schwartz'\left(\R^{\dimension}\right)$
is a well-defined tempered distribution for all $j\in I$. In view
of the inclusion $Z\left(\CalO\right)\hookrightarrow\Schwartz\left(\R^{\dimension}\right)$,
we thus see that every \emph{individual} term of each of the series
in equation (\ref{eq:GammaSynthesisFunctionalExplicitDefinition})
is well-defined. Here, we use that $\widehat{\gamma_{j,1}}\in C^{\infty}\left(\R^{\dimension}\right)$
with all derivatives of at most polynomial growth, so that the same
holds for $\widehat{\gamma_{1}^{\left(j\right)}}=\widehat{\gamma_{j,1}}\circ S_{j}^{-1}$.
We still have to show, however, that (each of) the series in equation
(\ref{eq:GammaSynthesisFunctionalExplicitDefinition}) converges (absolutely)
for every $f\in Z\left(\CalO\right)$ and defines a continuous linear
functional.

Since $Z\left(\CalO\right)=\Fourier\left(\TestFunctionSpace{\CalO}\right)$,
this is equivalent to showing for arbitrary $\left(g_{j}\right)_{j\in I}\in\ell_{w}^{q}\left(\left[V_{j}\right]_{j\in I}\right)$
that the series defining the functional
\[
\phi:\TestFunctionSpace{\CalO}\to\Compl,f\mapsto\sum_{j\in I}\left\langle \widehat{\gamma_{1}^{\left(j\right)}}\cdot\widehat{g_{j}},\,f\right\rangle _{\Schwartz',\Schwartz}
\]
converges absolutely for each $f\in\TestFunctionSpace{\CalO}$ and
that $\phi\in\DistributionSpace{\CalO}$. With the same reasoning
as above, we see at least that each term in the series is well-defined.
In the following, we will show that $\phi\in\DistributionSpace{\CalO}$
is indeed well-defined, with absolute convergence of the series.

But first, let us \emph{assume} that this is the case. Then we have,
for fixed $i\in I$ (by the usual formula for the (inverse) Fourier
transform of a compactly supported distribution, see e.g.\@ \cite[Theorem 7.23]{RudinFA})
\begin{align}
\left|\left[\Fourier^{-1}\left(\varphi_{i}\phi\right)\right]\left(x\right)\right| & =\left|\left\langle \phi,\,\varphi_{i}\cdot e^{2\pi i\left\langle x,\cdot\right\rangle }\right\rangle _{\DistributionSpace{\CalO},\TestFunctionSpace{\CalO}}\right|\nonumber \\
 & =\left|\sum_{j\in I}\left\langle \widehat{\gamma_{1}^{\left(j\right)}}\cdot\widehat{g_{j}},\,\varphi_{i}\cdot e^{2\pi i\left\langle x,\cdot\right\rangle }\right\rangle _{\Schwartz',\Schwartz}\right|\nonumber \\
 & \leq\sum_{j\in I}\left|\left\langle \widehat{\gamma_{1}^{\left(j\right)}}\cdot\widehat{g_{j}},\,\varphi_{i}\cdot e^{2\pi i\left\langle x,\cdot\right\rangle }\right\rangle _{\Schwartz',\Schwartz}\right|\nonumber \\
 & =\sum_{j\in I}\left|\left[\Fourier^{-1}\left(\varphi_{i}\cdot\widehat{\gamma_{1}^{\left(j\right)}}\cdot\widehat{g_{j}}\right)\right]\left(x\right)\right|\nonumber \\
\left({\scriptstyle \text{Lemma }\ref{lem:SpecialConvolutionConsistent}}\right) & =\sum_{j\in I}\left|\left[\Fourier^{-1}\left(\varphi_{i}\cdot\widehat{\gamma_{1}^{\left(j\right)}}\right)\ast g_{j}\right]\left(x\right)\right|\qquad\forall x\in\R^{\dimension},\label{eq:GammaSynthesisLocalizedFourier}
\end{align}
where all but the first three terms always make sense (as elements
of $\left[0,\infty\right]$), even without assuming that $\phi$ is
a well-defined distribution.

Now, we invoke Theorem \ref{thm:BandlimitedWienerAmalgamSelfImproving}
to obtain for all $f\in\Schwartz'\left(\R^{\dimension}\right)$ with
$\supp\widehat{f}\subset\overline{Q_{i}}\subset T_{i}\left[-R_{\CalQ},R_{\CalQ}\right]^{\dimension}+b_{i}$
that $\left\Vert f\right\Vert _{W_{T_{i}^{-T}\left[-1,1\right]^{\dimension}}\left(L_{v}^{p}\right)}\lesssim\,\left\Vert f\right\Vert _{L_{v}^{p}}$,
where the implied constant depends on $p,\dimension$ and on $K,R_{\CalQ}$,
which are fixed throughout. In combination with the embedding $W_{T_{i}^{-T}\left[-1,1\right]^{\dimension}}\left(L_{v}^{p}\right)\hookrightarrow L_{v}^{\infty}\left(\R^{\dimension}\right)\hookrightarrow L_{\left(1+\left|\mybullet\right|\right)^{-K}}^{\infty}\left(\R^{\dimension}\right)$
from equation (\ref{eq:WeightedWienerAmalgamTemperedDistribution})
(where now the norm of the embedding depends on $i$ (and on $p,\dimension,K,v$)),
we thus get for every $i\in I$ some constant $C^{\left(i\right)}=C^{\left(i\right)}\left(p,\dimension,K,v,R_{\CalQ}\right)>0$
such that $\left\Vert f\right\Vert _{\ast}\leq C^{\left(i\right)}\cdot\left\Vert f\right\Vert _{L_{v}^{p}}$
for all $f\in\Schwartz'\left(\R^{\dimension}\right)$ satisfying $\supp\widehat{f}\subset\overline{Q_{i}}$,
with $\left\Vert f\right\Vert _{\ast}:=\sup_{x\in\R^{\dimension}}\left(1+\left|x\right|\right)^{-K}\left|f\left(x\right)\right|$.
Since $\supp\Fourier\left[\Fourier^{-1}\left(\varphi_{i}\cdot\widehat{\gamma_{1}^{\left(j\right)}}\right)\ast g_{j}\right]=\supp\left(\varphi_{i}\cdot\widehat{\gamma_{1}^{\left(j\right)}}\cdot\widehat{g_{j}}\right)\subset\overline{Q_{i}}$,
this yields
\[
\left\Vert \Fourier^{-1}\left(\varphi_{i}\cdot\widehat{\gamma_{1}^{\left(j\right)}}\right)\ast g_{j}\right\Vert _{\ast}\leq C^{\left(i\right)}\cdot\left\Vert \Fourier^{-1}\left(\varphi_{i}\cdot\widehat{\gamma_{1}^{\left(j\right)}}\right)\ast g_{j}\right\Vert _{L_{v}^{p}}\qquad\forall j\in I.
\]
Now, we distinguish two cases:

\textbf{Case 1}: We have $p\in\left[1,\infty\right]$. In this case,
we can simply use the weighted Young inequality (equation (\ref{eq:WeightedYoungInequality}))
to derive
\[
\left\Vert \Fourier^{-1}\left(\varphi_{i}\cdot\widehat{\gamma_{1}^{\left(j\right)}}\right)\ast g_{j}\right\Vert _{L_{v}^{p}}\leq\left\Vert \Fourier^{-1}\left(\varphi_{i}\cdot\widehat{\gamma_{1}^{\left(j\right)}}\right)\right\Vert _{L_{v_{0}}^{1}}\cdot\left\Vert g_{j}\right\Vert _{L_{v}^{p}}=C_{i,j}\cdot\left\Vert g_{j}\right\Vert _{V_{j}}.
\]
But since we have $c=\left(c_{j}\right)_{j\in I}\in\ell_{w}^{q}\left(I\right)=\ell_{w^{\min\left\{ 1,p\right\} }}^{r}\left(I\right)$
for $c_{j}:=\left\Vert g_{j}\right\Vert _{V_{j}}$, we get by boundedness
of $\overrightarrow{C}$ that $\overrightarrow{C}c\in\ell_{w}^{q}\left(I\right)$.
In particular,
\begin{equation}
\frac{1}{C^{\left(i\right)}}\cdot\sum_{j\in I}\left\Vert \Fourier^{-1}\left(\varphi_{i}\cdot\widehat{\gamma_{1}^{\left(j\right)}}\right)\ast g_{j}\right\Vert _{\ast}\leq\sum_{j\in I}\left\Vert \Fourier^{-1}\left(\varphi_{i}\cdot\widehat{\gamma_{1}^{\left(j\right)}}\right)\ast g_{j}\right\Vert _{L_{v}^{p}}\leq\sum_{j\in I}C_{i,j}\left\Vert g_{j}\right\Vert _{V_{j}}=\left(\smash{\overrightarrow{C}}c\right)_{i}<\infty,\label{eq:GammaSynthesisBanachCaseEstimate}
\end{equation}
from which it follows that the series $\sum_{j\in I}\left[\Fourier^{-1}\left(\varphi_{i}\cdot\widehat{\gamma_{1}^{\left(j\right)}}\right)\ast g_{j}\right]\left(x\right)$
converges absolutely for all $x\in\R^{\dimension}$ (even locally
uniformly in $x$).

Furthermore, for arbitrary $\theta\in\TestFunctionSpace{\CalO}$,
we have
\begin{align*}
\sum_{j\in I}\left|\left\langle \widehat{\gamma_{1}^{\left(j\right)}}\cdot\widehat{g_{j}},\,\varphi_{i}\theta\right\rangle _{\Schwartz',\Schwartz}\right| & =\sum_{j\in I}\left|\left\langle \Fourier^{-1}\!\!\left(\varphi_{i}\cdot\widehat{\gamma_{1}^{\left(j\right)}}\cdot\widehat{g_{j}}\right),\,\widehat{\theta}\right\rangle _{\Schwartz',\Schwartz}\right|\\
 & \leq\sum_{j\in I}\left\Vert \Fourier^{-1}\!\!\left(\varphi_{i}\cdot\widehat{\gamma_{1}^{\left(j\right)}}\right)\!\ast\!g_{j}\right\Vert _{\ast}\left\Vert \smash{\widehat{\theta}}\right\Vert _{L_{\left(1+\left|\mybullet\right|\right)^{K}}^{1}}\leq C^{\left(i\right)}\cdot\left(\smash{\overrightarrow{C}}c\right)_{i}\cdot\left\Vert \smash{\widehat{\theta}}\right\Vert _{L_{\left(1+\left|\mybullet\right|\right)^{K}}^{1}}\!<\infty,
\end{align*}
so that the series $\sum_{j\in I}\left\langle \widehat{\gamma_{1}^{\left(j\right)}}\cdot\widehat{g_{j}},\,\varphi_{i}\theta\right\rangle $
defining $\phi\left(\varphi_{i}\theta\right)$ converges absolutely.
The same estimate also shows that $\theta\mapsto\phi\left(\varphi_{i}\theta\right)$
is a distribution on $\CalO$, since $\theta\mapsto\left\Vert \smash{\widehat{\theta}}\right\Vert _{L_{\left(1+\left|\mybullet\right|\right)^{K}}^{1}}$
is a continuous seminorm on $\TestFunctionSpace{\CalO}\hookrightarrow\Schwartz'\left(\R^{\dimension}\right)$.

But since $\left(\varphi_{i}\right)_{i\in I}$ is a locally finite
partition of unity on $\CalO$, we have $\theta=\sum_{i\in I_{\Upsilon}}\varphi_{i}\theta$
for every $\theta\in\TestFunctionSpace{\CalO}$ with $\supp\theta\subset\Upsilon$,
where $\Upsilon\subset\CalO$ is an arbitrary compact set and where
$I_{\Upsilon}\subset I$ is \emph{finite}. Hence, $\theta\mapsto\phi\left(\theta\right)=\sum_{i\in I_{\Upsilon}}\phi\left(\varphi_{i}\theta\right)$
is a continuous linear functional on $\left\{ \theta\in\TestFunctionSpace{\CalO}\with\supp\theta\subset\Upsilon\right\} $
for arbitrary compact $\Upsilon\subset\CalO$ and the defining series
converges absolutely (as a \emph{finite} sum of absolutely convergent
series). This shows that $\phi\in\DistributionSpace{\CalO}$ is well-defined
(with absolute convergence of the defining series), so that equation
(\ref{eq:GammaSynthesisLocalizedFourier}) is valid.

\medskip{}

As a consequence of equations (\ref{eq:GammaSynthesisLocalizedFourier})
and (\ref{eq:GammaSynthesisBanachCaseEstimate}) and of the triangle
inequality for $L_{v}^{p}\left(\R^{\dimension}\right)$, we finally
get
\[
\left\Vert \Fourier^{-1}\left(\varphi_{i}\phi\right)\right\Vert _{L_{v}^{p}}\leq\sum_{j\in I}\left\Vert \Fourier^{-1}\left(\varphi_{i}\cdot\widehat{\gamma_{1}^{\left(j\right)}}\right)\ast g_{j}\right\Vert _{L_{v}^{p}}\leq\left(\smash{\overrightarrow{C}}c\right)_{i}\qquad\forall i\in I,
\]
so that solidity of $\ell_{w}^{q}\left(I\right)$ yields
\[
\left\Vert \phi\right\Vert _{\FourierDecompSp{\CalQ}p{\ell_{w}^{q}}v}=\left\Vert \left(\left\Vert \Fourier^{-1}\left(\varphi_{i}\phi\right)\right\Vert _{L_{v}^{p}}\right)_{i\in I}\right\Vert _{\ell_{w}^{q}}\leq\left\Vert \smash{\overrightarrow{C}}c\right\Vert _{\ell_{w}^{q}}\leq\vertiii{\smash{\overrightarrow{C}}}\cdot\left\Vert c\right\Vert _{\ell_{w}^{q}}=\vertiii{\smash{\overrightarrow{C}}}\cdot\left\Vert \left(g_{j}\right)_{j\in I}\right\Vert _{\ell_{w}^{q}\left(\left[V_{i}\right]_{i\in I}\right)}<\infty.
\]
All in all, we see that $\phi\in\FourierDecompSp{\CalQ}p{\ell_{w}^{q}}v\leq\DistributionSpace{\CalO}$
is well-defined. But by definition of ${\rm Synth}_{\Gamma_{1}}$,
we have ${\rm Synth}_{\Gamma_{1}}\left(g_{j}\right)_{j\in I}=\Fourier^{-1}\phi$
for the isometric isomorphism $\Fourier^{-1}:\FourierDecompSp{\CalQ}p{\ell_{w}^{q}}v\to\DecompSp{\CalQ}p{\ell_{w}^{q}}v,\psi\mapsto\psi\circ\Fourier^{-1}$.
As a consequence, ${\rm Synth}_{\Gamma_{1}}\left(g_{j}\right)_{j\in I}\in\DecompSp{\CalQ}p{\ell_{w}^{q}}v$
is well-defined and 
\[
\left\Vert {\rm Synth}_{\Gamma_{1}}\left(g_{j}\right)_{j\in I}\right\Vert _{\DecompSp{\CalQ}p{\ell_{w}^{q}}v}=\left\Vert \phi\right\Vert _{\FourierDecompSp{\CalQ}p{\ell_{w}^{q}}v}\leq\vertiii{\smash{\overrightarrow{C}}}\cdot\left\Vert \left(g_{j}\right)_{j\in I}\right\Vert _{\ell_{w}^{q}\left(\left[V_{i}\right]_{i\in I}\right)},
\]
as desired.

\medskip{}

\textbf{Case 2}: We have $p\in\left(0,1\right)$. In this case, we
replace the application of the weighted Young inequality (equation
(\ref{eq:WeightedYoungInequality})) by an application of Corollary
\ref{cor:WienerAmalgamConvolutionSimplified} to get for $C_{1}:=\dimension^{-\frac{\dimension}{2p}}\cdot\left(972\cdot\dimension^{5/2}\right)^{K+\frac{\dimension}{p}}\cdot\Omega_{0}^{3K}\Omega_{1}^{3}$
that
\begin{align}
 & \left\Vert \Fourier^{-1}\left(\varphi_{i}\cdot\widehat{\gamma_{1}^{\left(j\right)}}\right)\ast g_{j}\right\Vert _{L_{v}^{p}}\nonumber \\
\left({\scriptstyle \text{Lemma }\ref{lem:MaximalFunctionDominatesF}}\right) & \leq\left\Vert \Fourier^{-1}\left(\varphi_{i}\cdot\widehat{\gamma_{1}^{\left(j\right)}}\right)\ast g_{j}\right\Vert _{W_{T_{j}^{-T}\left[-1,1\right]^{\dimension}}\left(L_{v}^{p}\right)}\nonumber \\
\left({\scriptstyle \text{Cor. }\ref{cor:WienerAmalgamConvolutionSimplified}}\right) & \leq C_{1}\cdot\left|\det T_{j}\right|^{\frac{1}{p}-1}\cdot\left\Vert \Fourier^{-1}\left(\varphi_{i}\cdot\widehat{\gamma_{1}^{\left(j\right)}}\right)\right\Vert _{W_{T_{j}^{-T}\left[-1,1\right]^{\dimension}}\left(L_{v_{0}}^{p}\right)}\cdot\left\Vert g_{j}\right\Vert _{W_{T_{j}^{-T}\left[-1,1\right]^{\dimension}}\left(L_{v}^{p}\right)}\nonumber \\
\left({\scriptstyle \text{eq. }\eqref{eq:WienerLinearCubeTransformationChange}}\right) & \leq C_{1}\left(6\dimension\right)^{K+\frac{\dimension}{p}}\cdot\Omega_{0}^{K}\Omega_{1}\cdot\left(1+\left\Vert T_{j}^{-1}T_{i}\right\Vert \right)^{K+\frac{\dimension}{p}}\cdot\left|\det T_{j}\right|^{\frac{1}{p}-1}\cdot\left\Vert \Fourier^{-1}\left(\varphi_{i}\cdot\widehat{\gamma_{1}^{\left(j\right)}}\right)\right\Vert _{W_{T_{i}^{-T}\left[-1,1\right]^{\dimension}}\left(L_{v_{0}}^{p}\right)}\cdot\left\Vert g_{j}\right\Vert _{V_{j}}\nonumber \\
\left({\scriptstyle \text{Thm. }\ref{thm:BandlimitedWienerAmalgamSelfImproving}}\right) & \overset{\left(\dagger\right)}{\leq}C_{1}C_{2}\left(6\dimension\right)^{K+\frac{\dimension}{p}}\cdot\Omega_{0}^{K}\Omega_{1}\cdot\left(1+\left\Vert T_{j}^{-1}T_{i}\right\Vert \right)^{K+\frac{\dimension}{p}}\cdot\left|\det T_{j}\right|^{\frac{1}{p}-1}\cdot\left\Vert \Fourier^{-1}\left(\varphi_{i}\cdot\widehat{\gamma_{1}^{\left(j\right)}}\right)\right\Vert _{L_{v_{0}}^{p}}\cdot\left\Vert g_{j}\right\Vert _{V_{j}}\nonumber \\
 & =C_{1}C_{2}\left(6\dimension\right)^{K+\frac{\dimension}{p}}\cdot\Omega_{0}^{K}\Omega_{1}\cdot C_{i,j}^{1/p}\cdot\left\Vert g_{j}\right\Vert _{V_{j}}=:C_{3}\cdot C_{i,j}^{1/p}\cdot\left\Vert g_{j}\right\Vert _{V_{j}},\label{eq:GammaSynthesisQuasiBanachCaseMainEstimate}
\end{align}
where  the step marked with $\left(\dagger\right)$ used that
\[
\supp\left(\varphi_{i}\cdot\widehat{\gamma_{1}^{\left(j\right)}}\right)\subset\overline{Q_{i}}\subset T_{i}\left[\overline{B_{R_{\CalQ}}}\left(0\right)\right]+b_{i}\subset T_{i}\left[-R_{\CalQ},R_{\CalQ}\right]^{\dimension}+b_{i},
\]
so that Theorem \ref{thm:BandlimitedWienerAmalgamSelfImproving} (with
$v_{0}$ instead of $v$) yields $\left\Vert \Fourier^{-1}\left(\varphi_{i}\cdot\widehat{\gamma_{1}^{\left(j\right)}}\right)\right\Vert _{W_{T_{i}^{-T}\left[-1,1\right]^{\dimension}}\left(L_{v_{0}}^{p}\right)}\leq C_{2}\cdot\left\Vert \Fourier^{-1}\left(\varphi_{i}\cdot\widehat{\gamma_{1}^{\left(j\right)}}\right)\right\Vert _{L_{v_{0}}^{p}}$
for
\[
C_{2}:=2^{4\left(1+\frac{\dimension}{p}\right)}s_{\dimension}^{\frac{1}{p}}\left(192\cdot\dimension^{3/2}\cdot\left\lceil K+\frac{\dimension+1}{p}\right\rceil \right)^{\left\lceil K+\frac{\dimension+1}{p}\right\rceil +1}\cdot\Omega_{0}^{K}\Omega_{1}\cdot\left(1+R_{\CalQ}\right)^{\frac{\dimension}{p}}.
\]

Next, set $c_{j}:=\left\Vert g_{j}\right\Vert _{V_{j}}^{p}$ for $j\in I$
and note that $\left(g_{j}\right)_{j\in I}\in\ell_{w}^{q}\left(\left[V_{j}\right]_{j\in I}\right)$
yields $c=\left(c_{j}\right)_{j\in I}\in\ell_{w^{p}}^{q/p}\left(I\right)=\ell_{w^{\min\left\{ 1,p\right\} }}^{r}\left(I\right)$.
Hence, we get because of $\ell^{p}\hookrightarrow\ell^{1}$ that
\begin{align*}
\frac{1}{C^{\left(i\right)}}\cdot\sum_{j\in I}\left\Vert \Fourier^{-1}\left(\varphi_{i}\cdot\widehat{\gamma_{1}^{\left(j\right)}}\right)\ast g_{j}\right\Vert _{\ast} & \leq\left(\sum_{j\in I}\left\Vert \Fourier^{-1}\left(\varphi_{i}\cdot\widehat{\gamma_{1}^{\left(j\right)}}\right)\ast g_{j}\right\Vert _{L_{v}^{p}}^{p}\right)^{1/p}\\
 & \leq C_{3}\cdot\left(\sum_{j\in I}C_{i,j}\cdot c_{j}\right)^{1/p}\\
 & =C_{3}\cdot\left(\smash{\overrightarrow{C}}\cdot c\right)_{i}^{1/p}<\infty,
\end{align*}
which is a slight variation of equation (\ref{eq:GammaSynthesisBanachCaseEstimate}).
Now, we see exactly as in case of $p\in\left[1,\infty\right]$ that
$\phi$ is a well-defined distribution $\phi\in\DistributionSpace{\CalO}$,
so that equation (\ref{eq:GammaSynthesisLocalizedFourier}) is valid.

Using this equation and the $p$-triangle inequality for $L_{v}^{p}\left(\R^{\dimension}\right)$,
we derive
\begin{align*}
\left\Vert \Fourier^{-1}\left(\varphi_{i}\phi\right)\right\Vert _{L_{v}^{p}} & \leq\left(\sum_{j\in I}\left\Vert \Fourier^{-1}\left(\varphi_{i}\cdot\widehat{\gamma_{1}^{\left(j\right)}}\right)\ast g_{j}\right\Vert _{L_{v}^{p}}^{p}\right)^{1/p}\\
\left({\scriptstyle \text{eq. }\eqref{eq:GammaSynthesisQuasiBanachCaseMainEstimate}}\right) & \leq C_{3}\cdot\left(\sum_{j\in I}C_{i,j}\cdot\left\Vert g_{j}\right\Vert _{V_{j}}^{p}\right)^{1/p}\\
 & =C_{3}\cdot\left(\smash{\overrightarrow{C}}\cdot c\right)_{i}^{1/p}<\infty
\end{align*}
and hence
\begin{align*}
\left\Vert \phi\right\Vert _{\FourierDecompSp{\CalQ}p{\ell_{w}^{q}}v} & =\left\Vert \left(\left\Vert \Fourier^{-1}\left(\varphi_{i}\phi\right)\right\Vert _{L_{v}^{p}}\right)_{i\in I}\right\Vert _{\ell_{w}^{q}}\\
 & \leq C_{3}\cdot\left\Vert \left(\smash{\overrightarrow{C}}\cdot c\right)^{1/p}\right\Vert _{\ell_{w}^{q}}\\
 & =C_{3}\cdot\left\Vert w^{p}\cdot\left[\smash{\overrightarrow{C}}\cdot c\right]\right\Vert _{\ell^{q/p}}^{1/p}\\
 & =C_{3}\cdot\left\Vert \smash{\overrightarrow{C}}\cdot c\right\Vert _{\ell_{w^{\min\left\{ 1,p\right\} }}^{r}}^{1/p}\\
 & \leq C_{3}\cdot\vertiii{\smash{\overrightarrow{C}}}^{1/p}\cdot\left\Vert c\right\Vert _{\ell_{w^{\min\left\{ 1,p\right\} }}^{r}}^{1/p}\\
 & =C_{3}\cdot\vertiii{\smash{\overrightarrow{C}}}^{1/p}\cdot\left\Vert \left(w_{j}\cdot\left\Vert g_{j}\right\Vert _{V_{j}}\right)_{j\in I}\right\Vert _{\ell^{q}}\\
 & =C_{3}\cdot\vertiii{\smash{\overrightarrow{C}}}^{1/p}\cdot\left\Vert \left(g_{j}\right)_{j\in I}\right\Vert _{\ell_{w}^{q}\left(\left[V_{j}\right]_{j\in I}\right)}<\infty.
\end{align*}
Now, we see as for $p\in\left[1,\infty\right]$ that ${\rm Synth}_{\Gamma_{1}}$
is a bounded linear operator with $\vertiii{{\rm Synth}_{\Gamma_{1}}}\leq C_{3}\cdot\vertiii{\smash{\overrightarrow{C}}}^{1/p}$.

Finally, we observe, using $s_{\dimension}\leq2^{2\dimension}$, that
\begin{align*}
C_{3} & =C_{1}C_{2}\left(6\dimension\right)^{K+\frac{\dimension}{p}}\cdot\Omega_{0}^{K}\Omega_{1}\\
 & =\dimension^{-\frac{\dimension}{2p}}\cdot\left(972\cdot\dimension^{5/2}\right)^{K+\frac{\dimension}{p}}\left(6\dimension\right)^{K+\frac{\dimension}{p}}\cdot2^{4\left(1+\frac{\dimension}{p}\right)}s_{\dimension}^{\frac{1}{p}}\left(192\cdot\dimension^{3/2}\cdot\left\lceil K+\frac{\dimension+1}{p}\right\rceil \right)^{\left\lceil K+\frac{\dimension+1}{p}\right\rceil +1}\cdot\left(1+R_{\CalQ}\right)^{\frac{\dimension}{p}}\cdot\Omega_{0}^{5K}\Omega_{1}^{5}\\
 & \leq\dimension^{-\frac{\dimension}{2p}}\cdot\left(5832\cdot\dimension^{7/2}\right)^{K+\frac{\dimension}{p}}\cdot2^{4+6\frac{\dimension}{p}}\left(192\cdot\dimension^{3/2}\cdot\left\lceil K+\frac{\dimension+1}{p}\right\rceil \right)^{\left\lceil K+\frac{\dimension+1}{p}\right\rceil +1}\cdot\left(1+R_{\CalQ}\right)^{\frac{\dimension}{p}}\cdot\Omega_{0}^{5K}\Omega_{1}^{5}\\
 & \leq\dimension^{-\frac{\dimension}{2p}}\cdot\left(5832\cdot\dimension^{7/2}\right)^{-2}\cdot2^{4+6\frac{\dimension}{p}}\left(2^{21}\cdot\dimension^{5}\cdot\left\lceil K+\frac{\dimension+1}{p}\right\rceil \right)^{\left\lceil K+\frac{\dimension+1}{p}\right\rceil +1}\cdot\left(1+R_{\CalQ}\right)^{\frac{\dimension}{p}}\cdot\Omega_{0}^{5K}\Omega_{1}^{5}\\
 & \leq\frac{\left(2^{6}/\sqrt{\dimension}\right)^{\dimension/p}}{2^{21}\cdot\dimension^{7}}\cdot\left(2^{21}\cdot\dimension^{5}\cdot\left\lceil K+\frac{\dimension+1}{p}\right\rceil \right)^{\left\lceil K+\frac{\dimension+1}{p}\right\rceil +1}\cdot\left(1+R_{\CalQ}\right)^{\frac{\dimension}{p}}\cdot\Omega_{0}^{5K}\Omega_{1}^{5},
\end{align*}
which completes the proof.
\end{proof}
In order to switch from the continuous synthesis operator from the
preceding lemma to a discrete one, our next technical result is helpful.
\begin{lem}
\label{lem:SchwartzTranslationSynthesis}Let $p\in\left(0,\infty\right]$
and assume that $\varrho:\R^{\dimension}\to\Compl$ is measurable
and satisfies $\left\Vert \varrho\right\Vert _{K_{0}}<\infty$ with
$K_{0}$ and $\left\Vert \mybullet\right\Vert _{K_{0}}$ as in Assumption
\ref{assu:AtomicDecompositionAssumption}. Let $i\in I$ and $\delta\in\left(0,1\right]$
and let $V_{i}$ be defined as in Assumption \ref{assu:MainAssumptions}.
Furthermore, let the \textbf{coefficient space} $C_{i}^{\left(\delta\right)}$
be defined as in equation (\ref{eq:CoefficientSpaceDefinition}).

Then, the maps
\[
\Psi_{\left|\varrho\right|}^{\left(i,\delta\right)}:C_{i}^{\left(\delta\right)}\to V_{i},\left(c_{k}\right)_{k\in\Z^{\dimension}}\mapsto\left(\sum_{k\in\Z^{\dimension}}c_{k}\cdot L_{\delta\cdot k}\left|\varrho\right|\right)\circ T_{i}^{T}=\sum_{k\in\Z^{\dimension}}c_{k}\cdot L_{\delta\cdot T_{i}^{-T}k}\left|\varrho\circ T_{i}^{T}\right|
\]
and
\[
\Psi_{\varrho}^{\left(i,\delta\right)}:C_{i}^{\left(\delta\right)}\to V_{i},\left(c_{k}\right)_{k\in\Z^{\dimension}}\mapsto\left(\sum_{k\in\Z^{\dimension}}c_{k}\cdot L_{\delta\cdot k}\varrho\right)\circ T_{i}^{T}=\sum_{k\in\Z^{\dimension}}c_{k}\cdot L_{\delta\cdot T_{i}^{-T}k}\left[\varrho\circ T_{i}^{T}\right]
\]
are well-defined and bounded, with pointwise absolute convergence
of the series and with
\[
\vertiii{\Psi_{\varrho}^{\left(i,\delta\right)}}\leq\vertiii{\Psi_{\left|\varrho\right|}^{\left(i,\delta\right)}}\leq\begin{cases}
\Omega_{0}^{K}\Omega_{1}\cdot\left(1+2\sqrt{\dimension}\right)^{K_{0}}\left(\frac{s_{\dimension}}{p}\right)^{1/p}\cdot\left\Vert \varrho\right\Vert _{K_{0}}\cdot\left|\det T_{i}\right|^{-\frac{1}{p}}, & \text{if }p<1,\\
\Omega_{0}^{K}\Omega_{1}\cdot12^{\dimension+1}\cdot\delta^{-\left(1-\frac{1}{p}\right)\left(\dimension+1\right)}\cdot\left\Vert \varrho\right\Vert _{K_{0}}\cdot\left|\det T_{i}\right|^{-\frac{1}{p}}, & \text{if }p\geq1.
\end{cases}
\]

In particular, if $g\in L_{v_{0}}^{1}\left(\R^{\dimension}\right)$,
then
\begin{equation}
g\ast\left[\Psi_{\varrho}^{\left(i,\delta\right)}\left(c_{k}\right)_{k\in\Z^{\dimension}}\right]=\sum_{k\in\Z^{\dimension}}\left(c_{k}\cdot\left[g\ast L_{\delta\cdot T_{i}^{-T}k}\left[\varrho\circ T_{i}^{T}\right]\right]\right).\qedhere\label{eq:SchwartzTranslationSynthesisConvolution}
\end{equation}
\end{lem}
\begin{proof}
Clearly, since $V_{i}$ and $C_{i}^{\left(\delta\right)}$ are solid,
boundedness of $\Psi_{\left|\varrho\right|}^{\left(i,\delta\right)}$
implies that of $\Psi_{\varrho}^{\left(i,\delta\right)}$, with $\vertiii{\Psi_{\varrho}^{\left(i,\delta\right)}}\leq\vertiii{\Psi_{\left|\varrho\right|}^{\left(i,\delta\right)}}$.
Furthermore, again by solidity and since $\left|\varrho\left(x\right)\right|\leq\left\Vert \varrho\right\Vert _{K_{0}}\cdot\left(1+\left|x\right|\right)^{-K_{0}}$
for all $x\in\R^{\dimension}$, it suffices to prove the claim (except
for equation (\ref{eq:SchwartzTranslationSynthesisConvolution}))
for the special case $\varrho\left(x\right)=\left(1+\left|x\right|\right)^{-K_{0}}$,
so that $\left\Vert \varrho\right\Vert _{K_{0}}=1$.

Recall from equation (\ref{eq:CoefficientSpaceDefinition}) that $v_{k}^{\left(j,\delta\right)}=v\left(\delta\cdot T_{j}^{-T}k\right)$.
Now, we first observe 
\begin{equation}
v\left(x\right)=v\left(\delta\cdot T_{i}^{-T}k+x-\delta\cdot T_{i}^{-T}k\right)\leq v\left(\delta\cdot T_{i}^{-T}k\right)\cdot v_{0}\left(x-\delta\cdot T_{i}^{-T}k\right)=v_{k}^{\left(i,\delta\right)}\cdot v_{0}\left(x-\delta\cdot T_{i}^{-T}k\right),\label{eq:SchwartzTranslationSynthesisWeightTrick}
\end{equation}
so that
\[
0\leq\frac{v\left(x\right)}{v_{k}^{\left(i,\delta\right)}}\cdot\left(L_{\delta\cdot T_{i}^{-T}k}\left[\varrho\circ T_{i}^{T}\right]\right)\left(x\right)\leq v_{0}\left(x-\delta\cdot T_{i}^{-T}k\right)\cdot\left(\varrho\circ T_{i}^{T}\right)\left(x-\delta\cdot T_{i}^{-T}k\right).
\]
By translation invariance of $\left\Vert \mybullet\right\Vert _{L^{1}}$,
this implies for $p\in\left[1,\infty\right]$ (which entails $K_{0}=K+\dimension+1$)
that
\begin{align*}
\left\Vert \frac{L_{\delta\cdot T_{i}^{-T}k}\left[\varrho\circ T_{i}^{T}\right]}{v_{k}^{\left(i,\delta\right)}}\right\Vert _{L_{v}^{1}} & \leq\left\Vert x\mapsto v_{0}\left(x-\delta\cdot T_{i}^{-T}k\right)\cdot\left(\varrho\circ T_{i}^{T}\right)\left(x-\delta\cdot T_{i}^{-T}k\right)\right\Vert _{L^{1}}\\
 & =\left\Vert v_{0}\cdot\left(\varrho\circ T_{i}^{T}\right)\right\Vert _{L^{1}}\\
\left({\scriptstyle \text{standard change of variables}}\right) & =\left|\det T_{i}^{T}\right|^{-1}\cdot\left\Vert \left(v_{0}\circ T_{i}^{-T}\right)\cdot\varrho\right\Vert _{L^{1}}\\
\left({\scriptstyle \text{assumptions on }v_{0}}\right) & \leq\Omega_{1}\cdot\left|\det T_{i}\right|^{-1}\cdot\left\Vert x\mapsto\left(1+\left|T_{i}^{-T}x\right|\right)^{K}\cdot\varrho\left(x\right)\right\Vert _{L^{1}}\\
\left({\scriptstyle \text{eq. }\eqref{eq:WeightLinearTransformationsConnection}}\right) & \leq\Omega_{0}^{K}\Omega_{1}\cdot\left|\det T_{i}\right|^{-1}\cdot\left\Vert x\mapsto\left(1+\left|x\right|\right)^{K}\cdot\varrho\left(x\right)\right\Vert _{L^{1}}\\
\left({\scriptstyle K_{0}=K+\dimension+1\text{ since }p\in\left[1,\infty\right]}\right) & =\Omega_{0}^{K}\Omega_{1}\cdot\left|\det T_{i}\right|^{-1}\cdot\left\Vert x\mapsto\left(1+\left|x\right|\right)^{-\left(\dimension+1\right)}\right\Vert _{L^{1}}\\
\left({\scriptstyle \text{eq. }\eqref{eq:StandardDecayLpEstimate}}\right) & \leq\Omega_{0}^{K}\Omega_{1}s_{\dimension}\cdot\left|\det T_{i}\right|^{-1}.
\end{align*}
Hence, we get in case of $p=1$ that
\begin{align*}
\left\Vert \Psi_{\left|\varrho\right|}^{\left(i,\delta\right)}\left(c_{k}\right)_{k\in\Z^{\dimension}}\right\Vert _{L_{v}^{1}} & \leq\sum_{k\in\Z^{\dimension}}v_{k}^{\left(i,\delta\right)}\left|c_{k}\right|\cdot\left\Vert \frac{L_{\delta\cdot T_{i}^{-T}k}\left[\varrho\circ T_{i}^{T}\right]}{v_{k}^{\left(i,\delta\right)}}\right\Vert _{L_{v}^{1}}\\
 & \leq\Omega_{0}^{K}\Omega_{1}s_{\dimension}\cdot\left|\det T_{i}\right|^{-1}\cdot\sum_{k\in\Z^{\dimension}}v_{k}^{\left(i,\delta\right)}\left|c_{k}\right|\\
 & =\Omega_{0}^{K}\Omega_{1}s_{\dimension}\cdot\left|\det T_{i}\right|^{-1}\cdot\left\Vert \left(c_{k}\right)_{k\in\Z^{\dimension}}\right\Vert _{C_{i}^{\left(\delta\right)}}<\infty.
\end{align*}
This establishes boundedness of $\Psi_{\left|\varrho\right|}^{\left(i,\delta\right)}$
for $p=1$.

\medskip{}

As our next step, we first note
\[
\left[M_{Q}\left(L_{x}f\right)\right]\left(y\right)=\left\Vert \Indicator_{y+Q}\cdot L_{x}f\right\Vert _{L^{\infty}}=\left\Vert \left(L_{-x}\Indicator_{y+Q}\right)\cdot f\right\Vert _{L^{\infty}}=\left\Vert \Indicator_{y-x+Q}\cdot f\right\Vert _{L^{\infty}}=\left(M_{Q}f\right)\left(y-x\right)=\left(L_{x}\left[M_{Q}f\right]\right)\left(y\right)
\]
for arbitrary measurable $f:\R^{\dimension}\to\Compl$ and $Q\subset\R^{\dimension}$.
Hence,
\begin{align*}
g_{i}^{\left(\delta\right)}\left(x\right) & :=v\left(x\right)\cdot\left[M_{T_{i}^{-T}\left[-1,1\right]^{\dimension}}\left(\frac{L_{\delta\cdot T_{i}^{-T}k}\left[\varrho\circ T_{i}^{T}\right]}{v_{k}^{\left(i,\delta\right)}}\right)\right]\left(x\right)\\
 & =\frac{v\left(x\right)}{v_{k}^{\left(i,\delta\right)}}\cdot\left[M_{T_{i}^{-T}\left[-1,1\right]^{\dimension}}\left(\varrho\circ T_{i}^{T}\right)\right]\left(x-\delta\cdot T_{i}^{-T}k\right)\\
\left({\scriptstyle \text{eq. }\eqref{eq:SchwartzTranslationSynthesisWeightTrick}}\right) & \leq\left[v_{0}\cdot M_{T_{i}^{-T}\left[-1,1\right]^{\dimension}}\left(\varrho\circ T_{i}^{T}\right)\right]\left(x-\delta\cdot T_{i}^{-T}k\right)\\
\left({\scriptstyle \text{Lemma }\ref{lem:WienerTransformationFormula}}\right) & =\left(v_{0}\cdot\left[\left(M_{\left[-1,1\right]^{\dimension}}\varrho\right)\circ T_{i}^{T}\right]\right)\left(x-\delta\cdot T_{i}^{-T}k\right)\\
\left({\scriptstyle \text{assumptions on }v_{0}\text{ and eq. }\eqref{eq:WeightLinearTransformationsConnection}}\right) & \leq\Omega_{0}^{K}\Omega_{1}\cdot\left[\left(1+\left|\mybullet\right|\right)^{K}\cdot\left(M_{\left[-1,1\right]^{\dimension}}\varrho\right)\right]\left(T_{i}^{T}x-\delta\cdot k\right)\\
\left({\scriptstyle \text{Lemma }\ref{lem:SchwartzFunctionsAreWiener}}\right) & \leq\Omega_{0}^{K}\Omega_{1}\cdot\left(1+2\sqrt{\dimension}\right)^{K_{0}}\cdot\left[\left(1+\left|\mybullet\right|\right)^{K-K_{0}}\right]\left(T_{i}^{T}x-\delta k\right).
\end{align*}

But this implies for $p\in\left(0,1\right)$ that
\begin{align*}
\left\Vert \frac{L_{\delta\cdot T_{i}^{-T}k}\left[\varrho\circ T_{i}^{T}\right]}{v_{k}^{\left(i,\delta\right)}}\right\Vert _{V_{i}} & =\left\Vert g_{i}^{\left(\delta\right)}\right\Vert _{L^{p}}\\
 & \leq\Omega_{0}^{K}\Omega_{1}\cdot\left(1+2\sqrt{\dimension}\right)^{K_{0}}\cdot\left\Vert \left[L_{\delta k}\left(1+\left|\mybullet\right|\right)^{K-K_{0}}\right]\circ T_{i}^{T}\right\Vert _{L^{p}}\\
 & =\Omega_{0}^{K}\Omega_{1}\cdot\left(1+2\sqrt{\dimension}\right)^{K_{0}}\cdot\left|\det T_{i}^{T}\right|^{-1/p}\cdot\left\Vert \left(1+\left|\mybullet\right|\right)^{K-K_{0}}\right\Vert _{L^{p}}\\
\left({\scriptstyle \text{eq. }\eqref{eq:StandardDecayLpEstimate}\text{ and }K-K_{0}=-\left(\frac{\dimension}{p}+1\right)\text{ since }p\in\left(0,1\right)}\right) & \leq\Omega_{0}^{K}\Omega_{1}\cdot\left(1+2\sqrt{\dimension}\right)^{K_{0}}\left(\frac{s_{\dimension}}{p}\right)^{1/p}\cdot\left|\det T_{i}\right|^{-1/p}.
\end{align*}
Now, we recall that for $p\in\left(0,1\right)$, we have the $p$-triangle
inequality $\left\Vert f+g\right\Vert _{L^{p}}^{p}\leq\left\Vert f\right\Vert _{L^{p}}^{p}+\left\Vert g\right\Vert _{L^{p}}^{p}$.
By solidity and because of $M_{Q}\left(f+g\right)\leq M_{Q}f+M_{Q}g$,
this also yields $\left\Vert f+g\right\Vert _{V_{i}}^{p}\leq\left\Vert f\right\Vert _{V_{i}}^{p}+\left\Vert g\right\Vert _{V_{i}}^{p}$,
so that we get
\begin{align*}
\left\Vert \Psi_{\left|\varrho\right|}^{\left(i,\delta\right)}\left(c_{k}\right)_{k\in\Z^{\dimension}}\right\Vert _{V_{i}}^{p} & \leq\sum_{k\in\Z^{\dimension}}\left[v_{k}^{\left(i,\delta\right)}\left|c_{k}\right|\right]^{p}\cdot\left\Vert \frac{L_{\delta\cdot T_{i}^{-T}k}\left[\varrho\circ T_{i}^{T}\right]}{v_{k}^{\left(i,\delta\right)}}\right\Vert _{V_{i}}^{p}\\
 & \leq\left[\Omega_{0}^{K}\Omega_{1}\cdot\left(1+2\sqrt{\dimension}\right)^{K_{0}}\left(\frac{s_{\dimension}}{p}\right)^{1/p}\cdot\left|\det T_{i}\right|^{-1/p}\right]^{p}\cdot\left\Vert \left(c_{k}\right)_{k\in\Z^{\dimension}}\right\Vert _{C_{i}^{\left(\delta\right)}}^{p},
\end{align*}
which yields the desired boundedness for $p\in\left(0,1\right)$.

\medskip{}

Next, we consider the case $p=\infty$. Here, we note because of
$\left|\varrho\left(x\right)\right|=\varrho\left(x\right)=\left(1+\left|x\right|\right)^{-K_{0}}$
that
\begin{align*}
v\left(x\right)\cdot\left|\left[\Psi_{\left|\varrho\right|}^{\left(i,\delta\right)}\left(c_{k}\right)_{k\in\Z^{\dimension}}\right]\left(x\right)\right| & =\left|\sum_{k\in\Z^{\dimension}}c_{k}v_{k}^{\left(i,\delta\right)}\cdot\frac{v\left(x\right)}{v_{k}^{\left(i,\delta\right)}}\cdot\left(\varrho\circ T_{i}^{T}\right)\left(x-\delta\cdot T_{i}^{-T}k\right)\right|\\
\left({\scriptstyle \text{for }c=\left(c_{k}\right)_{k\in\Z^{\dimension}}\text{ and since }p=\infty}\right) & \leq\left\Vert c\right\Vert _{C_{i}^{\left(\delta\right)}}\cdot\sum_{k\in\Z^{\dimension}}\left[\frac{v\left(x\right)}{v_{k}^{\left(i,\delta\right)}}\cdot\left(\varrho\circ T_{i}^{T}\right)\left(x-\delta\cdot T_{i}^{-T}k\right)\right]\\
\left({\scriptstyle \text{eq. }\eqref{eq:SchwartzTranslationSynthesisWeightTrick}\text{ and assumption on }v_{0}}\right) & \leq\Omega_{1}\cdot\left\Vert c\right\Vert _{C_{i}^{\left(\delta\right)}}\cdot\sum_{k\in\Z^{\dimension}}\left[\left(1+\left|x-\delta\cdot T_{i}^{-T}k\right|\right)^{K}\cdot\left(\varrho\circ T_{i}^{T}\right)\left(x-\delta\cdot T_{i}^{-T}k\right)\right]\\
\left({\scriptstyle \text{eq. }\eqref{eq:WeightLinearTransformationsConnection}\text{ and }\varrho\left(x\right)=\left(1+\left|x\right|\right)^{-K_{0}}}\right) & \leq\Omega_{0}^{K}\Omega_{1}\cdot\left\Vert c\right\Vert _{C_{i}^{\left(\delta\right)}}\cdot\sum_{k\in\Z^{\dimension}}\left[\left(1+\left|T_{i}^{T}\left(x-\delta\cdot T_{i}^{-T}k\right)\right|\right)^{K-K_{0}}\right]\\
 & =\Omega_{0}^{K}\Omega_{1}\cdot\left\Vert c\right\Vert _{C_{i}^{\left(\delta\right)}}\cdot h\left(T_{i}^{T}x\right),
\end{align*}
where we introduced $h\left(y\right):=\sum_{k\in\Z^{\dimension}}\left[\left(1+\left|y-\delta\cdot k\right|\right)^{K-K_{0}}\right]$
for $y\in\R^{\dimension}$ in the last step. Now, we recall $0<\delta\leq1$
and $K-K_{0}\leq-\left(\dimension+1\right)<0$, so that
\begin{align*}
h\left(y\right) & \leq\sum_{k\in\Z^{\dimension}}\left(1+\left|\delta\cdot\left(\frac{y}{\delta}-k\right)\right|\right)^{-\left(\dimension+1\right)}\\
 & \leq\delta^{-\left(\dimension+1\right)}\cdot\sum_{k\in\Z^{\dimension}}\left(1+\left|\frac{y}{\delta}-k\right|\right)^{-\left(\dimension+1\right)}\\
\left({\scriptstyle \text{since }\left|x\right|\geq\left\Vert x\right\Vert _{\infty}}\right) & \leq\delta^{-\left(\dimension+1\right)}\cdot\sum_{k\in\Z^{\dimension}}\left(1+\left\Vert \frac{y}{\delta}-k\right\Vert _{\infty}\right)^{-\left(\dimension+1\right)}=:\delta^{-\left(\dimension+1\right)}\cdot\widetilde{h}\left(\frac{y}{\delta}\right).
\end{align*}
Now, we note that $\widetilde{h}$ is $\Z^{\dimension}$-periodic,
so that $\left\Vert \smash{\widetilde{h}}\right\Vert _{\sup}=\left\Vert \smash{\widetilde{h}}\right\Vert _{\sup,\left[0,1\right)}$.
But for arbitrary $x\in\left[0,1\right)^{\dimension}$, we have
\[
1+\left\Vert k\right\Vert _{\infty}\leq1+\left\Vert k-x\right\Vert _{\infty}+\left\Vert x\right\Vert _{\infty}\leq2+\left\Vert x-k\right\Vert _{\infty}\leq2\left(1+\left\Vert x-k\right\Vert _{\infty}\right)
\]
and thus
\[
\widetilde{h}\left(x\right)=\sum_{k\in\Z^{\dimension}}\left(1+\left\Vert x-k\right\Vert _{\infty}\right)^{-\left(\dimension+1\right)}\leq2^{\dimension+1}\cdot\sum_{k\in\Z^{\dimension}}\left(1+\left\Vert k\right\Vert _{\infty}\right)^{-\left(\dimension+1\right)}.
\]
Next, we use the layer-cake formula (cf.\@ \cite[Proposition (6.24)]{FollandRA})
with the counting measure $\mu$ on $\Z^{\dimension}$ to estimate
for $\theta\left(k\right):=\left(1+\left\Vert k\right\Vert _{\infty}\right)^{-\left(\dimension+1\right)}$
the series
\begin{equation}
\begin{split}\sum_{k\in\Z^{\dimension}}\left(1+\left\Vert k\right\Vert _{\infty}\right)^{-\left(\dimension+1\right)}=\int_{\Z^{\dimension}}\theta\left(k\right)\d\mu\left(k\right) & =\int_{0}^{\infty}\mu\left(\left\{ k\in\Z^{\dimension}\with\theta\left(k\right)>\lambda\right\} \right)\d\lambda\\
\left({\scriptstyle \text{since }\theta\left(k\right)>\lambda\Longleftrightarrow\left\Vert k\right\Vert _{\infty}<\lambda^{\frac{-1}{\dimension+1}}-1\text{ can only hold for }\lambda<1}\right) & \leq\int_{0}^{1}\mu\left(\left[-\lambda^{-\frac{1}{\dimension+1}},\lambda^{-\frac{1}{\dimension+1}}\right]^{\dimension}\cap\Z^{\dimension}\right)\d\lambda\\
 & \leq\int_{0}^{1}\mu\left(\left\{ -\left\lfloor \lambda^{-\frac{1}{\dimension+1}}\right\rfloor ,\dots,\left\lfloor \lambda^{-\frac{1}{\dimension+1}}\right\rfloor \right\} ^{\dimension}\right)\d\lambda\\
 & =\int_{0}^{1}\left(1+2\left\lfloor \lambda^{-\frac{1}{\dimension+1}}\right\rfloor \right)^{\dimension}\d\lambda\\
 & \leq3^{\dimension}\cdot\int_{0}^{1}\lambda^{-\frac{\dimension}{\dimension+1}}\d\lambda=\left(\dimension+1\right)\cdot3^{\dimension}\cdot\lambda^{\frac{1}{\dimension+1}}\bigg|_{0}^{1}\\
\left({\scriptstyle \text{since }1+\dimension\leq2^{\dimension}}\right) & \leq6^{\dimension}.
\end{split}
\label{eq:StandardDecayLatticeSeries}
\end{equation}
Hence, we get $\widetilde{h}\left(x\right)\leq2^{\dimension+1}\cdot6^{\dimension}\leq2\cdot12^{\dimension}$
for all $x\in\left[0,1\right)^{\dimension}$ and thus for all $x\in\R^{\dimension}$,
by $\Z^{\dimension}$-periodicity. In view of the estimates from above,
this entails for $c=\left(c_{k}\right)_{k\in\Z^{\dimension}}$ that
\[
v\left(x\right)\cdot\left|\left[\Psi_{\left|\varrho\right|}^{\left(i,\delta\right)}\left(c_{k}\right)_{k\in\Z^{\dimension}}\right]\left(x\right)\right|\leq\Omega_{0}^{K}\Omega_{1}\cdot\left\Vert c\right\Vert _{C_{i}^{\left(\delta\right)}}\cdot h\left(T_{i}^{T}x\right)\leq\delta^{-\left(\dimension+1\right)}\cdot2\cdot12^{\dimension}\cdot\Omega_{0}^{K}\Omega_{1}\cdot\left\Vert c\right\Vert _{C_{i}^{\left(\delta\right)}}<\infty
\]
for all $x\in\R^{\dimension}$. In particular, the series defining
$\Psi_{\left|\varrho\right|}^{\left(i,\delta\right)}$ converges pointwise
absolutely in case of $p=\infty$. But since we have $\ell_{v^{\left(i,\delta\right)}}^{p}\left(\Z^{\dimension}\right)\hookrightarrow\ell_{v^{\left(i,\delta\right)}}^{\infty}\left(\Z^{\dimension}\right)$
for all $p\in\left(0,\infty\right]$, this implies absolute pointwise
convergence for arbitrary $p\in\left(0,\infty\right]$.

\medskip{}

Next, for $p\in\left[1,\infty\right]$, it is easy to see that the
operator
\[
\Lambda_{\left|\varrho\right|}^{\left(i,\delta\right)}:\ell^{p}\left(\smash{\Z^{\dimension}}\right)\to L^{p}\left(\smash{\R^{\dimension}}\right),\left(\zeta_{k}\right)_{k\in\Z^{\dimension}}\mapsto v\cdot\sum_{k\in\Z^{\dimension}}\frac{\zeta_{k}}{v_{k}^{\left(i,\delta\right)}}\cdot L_{\delta\cdot T_{i}^{-T}k}\left[\left|\varrho\right|\circ T_{i}^{T}\right]
\]
is well-defined and bounded if and only if $\Psi_{\left|\varrho\right|}^{\left(i,\delta\right)}$
is, with $\vertiii{\Lambda_{\left|\varrho\right|}^{\left(i,\delta\right)}}=\vertiii{\Psi_{\left|\varrho\right|}^{\left(i,\delta\right)}}$.
Hence, since we have already shown boundedness for $p\in\left\{ 1,\infty\right\} $,
we can use complex interpolation (the Riesz-Thorin theorem, \cite[Theorem (6.27)]{FollandRA})
to derive for $p\in\left[1,\infty\right]$ that $\Psi_{\left|\varrho\right|}^{\left(i,\delta\right)}$
is bounded, with
\begin{align*}
\vertiii{\Psi_{\left|\varrho\right|}^{\left(i,\delta\right)}} & \leq\left[\Omega_{0}^{K}\Omega_{1}s_{\dimension}\cdot\left|\det T_{i}\right|^{-1}\right]^{\frac{1}{p}}\cdot\left[\delta^{-\left(\dimension+1\right)}\cdot2\cdot12^{\dimension}\cdot\Omega_{0}^{K}\Omega_{1}\right]^{1-\frac{1}{p}}\\
\left({\scriptstyle \text{since }s_{\dimension}\leq4^{\dimension}}\right) & \leq\Omega_{0}^{K}\Omega_{1}\cdot12^{\dimension+1}\cdot\delta^{-\left(1-\frac{1}{p}\right)\left(\dimension+1\right)}\cdot\left|\det T_{i}\right|^{-\frac{1}{p}}.
\end{align*}

\medskip{}

Finally, for $g\in L_{v_{0}}^{1}\left(\R^{\dimension}\right)$ and
$c=\left(c_{k}\right)_{k\in\Z^{\dimension}}\in C_{i}^{\left(\delta\right)}=\ell_{v^{\left(i,\delta\right)}}^{p}\left(\Z^{\dimension}\right)\hookrightarrow\ell_{v^{\left(i,\delta\right)}}^{\infty}\left(\Z^{\dimension}\right)$,
our previous considerations, in combination with $v\left(x\right)=v\left(x-y+y\right)\leq v\left(x-y\right)v_{0}\left(y\right)$,
show for arbitrary measurable $\varrho:\R^{\dimension}\to\Compl$
with $\left\Vert \varrho\right\Vert _{K_{0}}<\infty$ that 
\begin{align*}
\int_{\R^{\dimension}}\left|g\left(y\right)\right|\!\cdot\!\sum_{k\in\Z^{\dimension}}\left|c_{k}\right|\!\cdot\!\left|\left(\!L_{\delta\cdot T_{i}^{-T}k}\left[\varrho\circ T_{i}^{T}\right]\!\right)\!\!\left(x-y\right)\right|\d y & \leq\frac{1}{v\left(x\right)}\!\cdot\!\int_{\R^{\dimension}}\left|v_{0}\left(y\right)\!\cdot\!g\left(y\right)\right|\cdot v\left(x\!-\!y\right)\cdot\left[\Psi_{\left|\varrho\right|}^{\left(i,\delta\right)}\left(\left|c_{k}\right|\right)_{k\in\Z^{\dimension}}\right]\!\!\left(x-y\right)\d y\\
 & \leq\frac{1}{v\left(x\right)}\cdot\left\Vert \Psi_{\left|\varrho\right|}^{\left(i,\delta\right)}\left(\left|c_{k}\right|\right)_{k\in\Z^{\dimension}}\right\Vert _{L_{v}^{\infty}}\cdot\left\Vert g\right\Vert _{L_{v_{0}}^{1}}<\infty
\end{align*}
for all $x\in\R^{\dimension}$, so that the interchange of summation
and integration in
\begin{align*}
\left[g\ast\Psi_{\varrho}^{\left(i,\delta\right)}\left(c_{k}\right)_{k\in\Z^{\dimension}}\right]\left(x\right) & =\int_{\R^{\dimension}}g\left(y\right)\cdot\left(\sum_{k\in\Z^{\dimension}}c_{k}\cdot L_{\delta\cdot T_{i}^{-T}k}\left[\varrho\circ T_{i}^{T}\right]\right)\left(x-y\right)\d y\\
 & =\sum_{k\in\Z^{\dimension}}c_{k}\cdot\int_{\R^{\dimension}}g\left(y\right)\cdot\left(L_{\delta\cdot T_{i}^{-T}k}\left[\varrho\circ T_{i}^{T}\right]\right)\left(x-y\right)\d y\\
 & =\sum_{k\in\Z^{\dimension}}c_{k}\cdot\left(g\ast L_{\delta\cdot T_{i}^{-T}k}\left[\varrho\circ T_{i}^{T}\right]\right)\left(x\right)
\end{align*}
is justified by the dominated convergence theorem.
\end{proof}
As a further ingredient, we will need the following ``sampling theorem''
for bandlimited functions. A very similar statement already appears
in \cite[Proposition in §1.3.3]{TriebelTheoryOfFunctionSpaces},
so no originality at all is claimed. Note, however, that in \cite{TriebelTheoryOfFunctionSpaces},
it is assumed that $\varphi\in\Schwartz\left(\R^{\dimension}\right)$
instead of $\varphi\in\Schwartz'\left(\R^{\dimension}\right)$ and
furthermore, the statement in \cite{TriebelTheoryOfFunctionSpaces}
is restricted to the \emph{unweighted} case.
\begin{lem}
\label{lem:BandlimitedSampling}For each $i\in I$, $R>0$ and $p\in\left(0,\infty\right]$,
as well as
\[
C:=2^{\max\left\{ 1,\frac{1}{p}\right\} }\cdot\Omega_{0}^{3K}\Omega_{1}^{3}\cdot\left(1+\sqrt{\dimension}\right)^{K}\cdot\left(23040\cdot\dimension^{3/2}\cdot\left(K+1+\frac{\dimension+1}{\min\left\{ 1,p\right\} }\right)\right)^{K+2+\frac{\dimension+1}{\min\left\{ 1,p\right\} }}\cdot\left(1+R\right)^{1+\frac{\dimension}{\min\left\{ 1,p\right\} }}
\]
we have
\[
\left\Vert \left[\varphi\left(\delta\cdot T_{i}^{-T}k\right)\right]_{k\in\Z^{\dimension}}\right\Vert _{C_{i}^{\left(\delta\right)}}\leq C\cdot\delta^{-\dimension/p}\cdot\left|\det T_{i}\right|^{1/p}\cdot\left\Vert \varphi\right\Vert _{L_{v}^{p}}
\]
for all $\delta\in\left(0,1\right]$ and all $\varphi\in\Schwartz'\left(\R^{\dimension}\right)$
with $\supp\widehat{\varphi}\subset T_{i}\left[-R,R\right]^{\dimension}+\xi_{0}$,
for arbitrary $\xi_{0}\in\R^{\dimension}$.
\end{lem}
\begin{proof}
Clearly, we can assume without loss of generality that $\left\Vert \varphi\right\Vert _{L_{v}^{p}}<\infty$.

Let us first consider the case $p=\infty$. Since we have $v\left(x\right)\cdot\left|\varphi\left(x\right)\right|\leq\left\Vert \varphi\right\Vert _{L_{v}^{\infty}}$
for almost all $x\in\R^{\dimension}$, and thus for a dense subset
of $\R^{\dimension}$, there is for arbitrary $k\in\Z^{\dimension}$
a sequence $\left(x_{n}\right)_{n\in\N}$ satisfying $x_{n}\xrightarrow[n\to\infty]{}\delta\cdot T_{i}^{-T}k$
as well as $v\left(x_{n}\right)\cdot\left|\varphi\left(x_{n}\right)\right|\leq\left\Vert \varphi\right\Vert _{L_{v}^{\infty}}$.
But since $\varphi$ is given by (integration against) a continuous
function by the Paley-Wiener theorem, this implies
\begin{align*}
v\left(\delta\cdot T_{i}^{-T}k\right)\cdot\left|\varphi\left(\delta\cdot T_{i}^{-T}k\right)\right| & =\lim_{n\to\infty}v\left(\delta\cdot T_{i}^{-T}k\right)\cdot\left|\varphi\left(x_{n}\right)\right|\\
 & \leq\liminf_{n\to\infty}v\left(x_{n}+\delta\cdot T_{i}^{-T}k-x_{n}\right)\cdot\left|\varphi\left(x_{n}\right)\right|\\
 & \leq\liminf_{n\to\infty}v\left(x_{n}\right)\cdot\left|\varphi\left(x_{n}\right)\right|\cdot v_{0}\left(\delta\cdot T_{i}^{-T}k-x_{n}\right)\\
 & \leq\Omega_{1}\cdot\left\Vert \varphi\right\Vert _{L_{v}^{\infty}}\cdot\liminf_{n\to\infty}\left(1+\left|\delta\cdot T_{i}^{-T}k-x_{n}\right|\right)^{K}\\
 & =\Omega_{1}\cdot\left\Vert \varphi\right\Vert _{L_{v}^{\infty}}<\infty.
\end{align*}
Since $C\geq\Omega_{1}$, since $\delta^{-\dimension/p}\cdot\left|\det T_{i}\right|^{1/p}=1$
for $p=\infty$ and since $k\in\Z^{\dimension}$ was arbitrary, this
establishes the claim for $p=\infty$. Hence, we can assume $p\in\left(0,\infty\right)$
in what follows.

\medskip{}

Let $C>0$ as in the statement of the theorem and $C_{1}:=2^{-\max\left\{ 1,\frac{1}{p}\right\} }\cdot\left[\Omega_{0}^{K}\Omega_{1}\cdot\left(1+\sqrt{\dimension}\right)^{K}\right]^{-1}\cdot C$.
Let $\varrho:=M_{-\xi_{0}}\varphi$ and note $\left|\varrho\right|=\left|\varphi\right|$.
Now, Theorem \ref{thm:BandlimitedOscillationSelfImproving} shows
\[
\left\Vert \osc{\delta\cdot T_{i}^{-T}\left[-1,1\right]^{\dimension}}\,\varrho\right\Vert _{W_{T_{i}^{-T}\left[-1,1\right]^{\dimension}}\left(L_{v}^{p}\right)}=\left\Vert \osc{\delta\cdot T_{i}^{-T}\left[-1,1\right]^{\dimension}}\left[M_{-\xi_{0}}\varphi\right]\right\Vert _{W_{T_{i}^{-T}\left[-1,1\right]^{\dimension}}\left(L_{v}^{p}\right)}\leq C_{1}\cdot\delta\cdot\left\Vert \varphi\right\Vert _{L_{v}^{p}}
\]
for all $\delta\in\left(0,1\right]$ and $\varphi$ as in the statement
of the lemma.

Now, notice for arbitrary $k\in\Z^{\dimension}$ and $x\in\delta T_{i}^{-T}\left(k+\left[0,1\right)^{\dimension}\right)$
that $\delta T_{i}^{-T}k\in x-\delta T_{i}^{-T}\left[0,1\right)^{\dimension}\subset x+\delta T_{i}^{-T}\left[-1,1\right]^{\dimension}$
and hence
\begin{align*}
\left|\varphi\left(\delta\cdot T_{i}^{-T}k\right)\right|=\left|\varrho\left(\delta\cdot T_{i}^{-T}k\right)\right| & \leq\left|\varrho\left(x\right)\right|+\left|\varrho\left(x\right)-\varrho\left(\delta\cdot T_{i}^{-T}k\right)\right|\\
 & \leq\left|\varphi\left(x\right)\right|+\left(\osc{\delta\cdot T_{i}^{-T}\left[-1,1\right]^{\dimension}}\varrho\right)\left(x\right).
\end{align*}
By multiplying this estimate with $\Indicator_{\delta T_{i}^{-T}\left(k+\left[0,1\right)^{\dimension}\right)}\left(x\right)$,
summing over $k\in\Z^{\dimension}$ and using $\R^{\dimension}=\biguplus_{k\in\Z^{\dimension}}\delta T_{i}^{-T}\left(k+\smash{\left[0,1\right)^{\dimension}}\right)$,
we obtain
\[
\sum_{k\in\Z^{\dimension}}\left(\Indicator_{\delta T_{i}^{-T}\left(k+\left[0,1\right)^{\dimension}\right)}\left(x\right)\cdot\left|\varphi\left(\delta T_{i}^{-T}k\right)\right|\right)\leq\left|\varphi\left(x\right)\right|+\left(\osc{\delta\cdot T_{i}^{-T}\left[-1,1\right]^{\dimension}}\,\varrho\right)\left(x\right)\qquad\forall x\in\R^{\dimension}.
\]
By solidity of $L_{v}^{p}\left(\R^{\dimension}\right)$, we conclude
\begin{align*}
 & \left\Vert \sum_{k\in\Z^{\dimension}}\left(\Indicator_{\delta T_{i}^{-T}\left(k+\left[0,1\right)^{\dimension}\right)}\cdot\left|\varphi\left(\delta T_{i}^{-T}k\right)\right|\right)\right\Vert _{L_{v}^{p}}\\
 & \leq\left\Vert \left|\varphi\right|+\osc{\delta\cdot T_{i}^{-T}\left[-1,1\right]^{\dimension}}\,\varrho\right\Vert _{L_{v}^{p}}\\
\left({\scriptstyle C_{2}:=2^{\max\left\{ 0,\frac{1}{p}-1\right\} }\text{ is triangle const. for }L_{v}^{p}}\right) & \leq C_{2}\cdot\left(\left\Vert \varphi\right\Vert _{L_{v}^{p}}+\left\Vert \osc{\delta\cdot T_{i}^{-T}\left[-1,1\right]^{\dimension}}\,\varrho\right\Vert _{L_{v}^{p}}\right)\\
 & \leq C_{2}\cdot\left(\left\Vert \varphi\right\Vert _{L_{v}^{p}}+C_{1}\cdot\delta\cdot\left\Vert \varphi\right\Vert _{L_{v}^{p}}\right)\\
\left({\scriptstyle \text{since }\delta\leq1}\right) & \leq C_{2}\left(1+C_{1}\right)\cdot\left\Vert \varphi\right\Vert _{L_{v}^{p}}\\
\left({\scriptstyle \text{since }C_{1}\geq1}\right) & \leq2^{\max\left\{ 1,\frac{1}{p}\right\} }\cdot C_{1}\cdot\left\Vert \varphi\right\Vert _{L_{v}^{p}}.
\end{align*}
Finally, we note for $x\in\delta\cdot T_{i}^{-T}\left(k+\left[0,1\right)^{\dimension}\right)$,
i.e., for $x=\delta T_{i}^{-T}k+\delta T_{i}^{-T}q$ with $q\in\left[0,1\right)^{\dimension}$
that
\begin{align*}
v_{k}^{\left(i,\delta\right)} & =v\left(\delta\cdot T_{i}^{-T}k\right)=v\left(x-\delta T_{i}^{-T}q\right)\leq v\left(x\right)\cdot v_{0}\left(\delta\cdot T_{i}^{-T}q\right)\\
\left({\scriptstyle \text{assump. on }v_{0}\text{ and eq. }\eqref{eq:WeightLinearTransformationsConnection}}\right) & \leq\Omega_{1}\cdot v\left(x\right)\cdot\left(1+\left|\delta\cdot T_{i}^{-T}q\right|\right)^{K}\leq\Omega_{0}^{K}\Omega_{1}\cdot v\left(x\right)\cdot\left(1+\left|\delta\cdot q\right|\right)^{K}\\
\left({\scriptstyle \text{since }\delta\leq1}\right) & \leq\Omega_{0}^{K}\Omega_{1}\cdot\left(1+\sqrt{\dimension}\right)^{K}\cdot v\left(x\right).
\end{align*}
Using this estimate and the pairwise disjointness of $\left(\delta\cdot T_{i}^{-T}\left(k+\left[0,1\right)^{\dimension}\right)\right)_{k\in\Z^{\dimension}}$,
we conclude
\begin{align*}
\left\Vert \sum_{k\in\Z^{\dimension}}\!\!\left(\Indicator_{\delta T_{i}^{-T}\!\left(k+\left[0,1\right)^{\dimension}\right)}\cdot\left|\varphi\!\left(\delta T_{i}^{-T}k\right)\right|\right)\right\Vert _{L_{v}^{p}}\!\!\! & =\left(\sum_{k\in\Z^{\dimension}}\left|\varphi\left(\delta\cdot T_{i}^{-T}k\right)\right|^{p}\!\int_{\delta T_{i}^{-T}\left(k+\left[0,1\right)\right)^{\dimension}}\left[v\left(x\right)\right]^{p}\d x\right)^{1/p}\\
 & \geq\left[\Omega_{0}^{K}\Omega_{1}\!\cdot\!\left(1\!+\!\sqrt{\dimension}\right)^{\!K}\right]^{-1}\!\!\!\cdot\!\left(\sum_{k\in\Z^{\dimension}}\!\left[v_{k}^{\left(i,\delta\right)}\!\cdot\left|\varphi\!\left(\delta T_{i}^{-T}k\right)\right|\right]^{p}\!\!\cdot\!\lambda_{\dimension}\left(\!\delta T_{i}^{-T}\!\left(k\!+\!\left[0,1\right)^{\dimension}\right)\!\right)\!\right)^{\!\!\frac{1}{p}}\\
 & =\left[\Omega_{0}^{K}\Omega_{1}\cdot\left(1+\sqrt{\dimension}\right)^{K}\right]^{-1}\cdot\delta^{\dimension/p}\cdot\left|\det T_{i}\right|^{-1/p}\cdot\left\Vert \left[\varphi\left(\delta\cdot T_{i}^{-T}k\right)\right]_{k\in\Z^{\dimension}}\right\Vert _{C_{i}^{\left(\delta\right)}}.
\end{align*}
Putting everything together, we conclude
\begin{align*}
\left\Vert \left[\varphi\left(\delta\cdot T_{i}^{-T}k\right)\right]_{k\in\Z^{\dimension}}\right\Vert _{C_{i}^{\left(\delta\right)}} & \leq\delta^{-\frac{\dimension}{p}}\cdot\left|\det T_{i}\right|^{\frac{1}{p}}\cdot\Omega_{0}^{K}\Omega_{1}\cdot\left(1+\sqrt{\dimension}\right)^{K}\cdot\left\Vert \sum_{k\in\Z^{\dimension}}\left(\Indicator_{\delta T_{i}^{-T}\left(k+\left[0,1\right)^{\dimension}\right)}\cdot\left|\varphi\left(\delta T_{i}^{-T}k\right)\right|\right)\right\Vert _{L_{v}^{p}}\\
 & \leq\delta^{-\frac{\dimension}{p}}\cdot\left|\det T_{i}\right|^{\frac{1}{p}}\cdot\Omega_{0}^{K}\Omega_{1}\cdot\left(1+\sqrt{\dimension}\right)^{K}\cdot2^{\max\left\{ 1,\frac{1}{p}\right\} }\cdot C_{1}\cdot\left\Vert \varphi\right\Vert _{L_{v}^{p}}\\
 & =\delta^{-\frac{\dimension}{p}}\cdot\left|\det T_{i}\right|^{\frac{1}{p}}\cdot C\cdot\left\Vert \varphi\right\Vert _{L_{v}^{p}}<\infty,
\end{align*}
as desired.
\end{proof}
In the proof of Theorem \ref{thm:AtomicDecomposition}, we will employ
a Neumann series argument for an operator defined on $\DecompSp{\CalQ}p{\ell_{w}^{q}}v$.
For this to be justified, we need to know that this space is a Quasi-Banach
space, i.e., complete. For $v\equiv1$, this was already shown in
\cite[Theorem 3.21]{DecompositionEmbedding}, but for $v\not\equiv1$
and general $p,q\in\left(0,\infty\right]$, it seems that the following
lemma is a novel (though not too surprising) result:
\begin{lem}
\label{lem:WeightedDecompositionSpaceComplete}The decomposition space
$\DecompSp{\CalQ}p{\ell_{w}^{q}}v$ is a Quasi-Banach space.
\end{lem}
\begin{proof}
Verifying the quasi-norm properties of $\left\Vert \mybullet\right\Vert _{\DecompSp{\CalQ}p{\ell_{w}^{q}}v}$
is relatively straightforward (and essentially identical to the verification
in \cite[Theorem 3.21]{DecompositionEmbedding}), so we only prove
completeness.

Instead of verifying completeness directly, we use a slightly more
abstract approach, employing other results from the paper. The main
new ingredient that we need to provide is boundedness of
\[
{\rm Ana}_{\ast}:\DecompSp{\CalQ}p{\ell_{w}^{q}}v\to V=\ell_{w}^{q}\left(\left[V_{i}\right]_{i\in I}\right),f\mapsto\left[\Fourier^{-1}\left(\varphi_{i}^{\ast}\widehat{f}\right)\right]_{i\in I}.
\]
To this end, let $f\in\DecompSp{\CalQ}p{\ell_{w}^{q}}v$ be arbitrary
and define $c_{i}:=\left\Vert \Fourier^{-1}\left(\varphi_{i}\widehat{f}\right)\right\Vert _{L_{v}^{p}}$
for $i\in I$. Note $c=\left(c_{i}\right)_{i\in I}\in\ell_{w}^{q}\left(I\right)$
and $\left\Vert c\right\Vert _{\ell_{w}^{q}}=\left\Vert f\right\Vert _{\DecompSp{\CalQ}p{\ell_{w}^{q}}v}$.
Recall that the clustering map $\Gamma_{\CalQ}:\ell_{w}^{q}\left(I\right)\to\ell_{w}^{q}\left(I\right)$
with $\Gamma_{\CalQ}\left(e_{i}\right)_{i\in I}=\left(e_{i}^{\ast}\right)_{i\in I}$
and $e_{i}^{\ast}=\sum_{\ell\in i^{\ast}}e_{\ell}$ is bounded.

Now, we distinguish the cases $p\in\left[1,\infty\right]$ and $p\in\left(0,1\right)$.
In case of $p\in\left[1,\infty\right]$, the triangle inequality for
$L_{v}^{p}\left(\R^{\dimension}\right)$ yields because of $V_{i}=L_{v}^{p}\left(\R^{\dimension}\right)$
for all $i\in I$ that
\begin{align*}
\left\Vert {\rm Ana}_{\ast}f\right\Vert _{V} & =\left\Vert \left(\left\Vert \Fourier^{-1}\left(\varphi_{i}^{\ast}\cdot\smash{\widehat{f}}\:\right)\right\Vert _{V_{i}}\right)_{i\in I}\right\Vert _{\ell_{w}^{q}}\leq\left\Vert \left(\sum_{\ell\in i^{\ast}}\left\Vert \Fourier^{-1}\left(\varphi_{\ell}\cdot\smash{\widehat{f}}\:\right)\right\Vert _{L_{v}^{p}}\right)_{i\in I}\right\Vert _{\ell_{w}^{q}}\\
 & =\left\Vert \Gamma_{\CalQ}\,c\right\Vert _{\ell_{w}^{q}}\leq\vertiii{\smash{\Gamma_{\CalQ}}}\cdot\left\Vert c\right\Vert _{\ell_{w}^{q}}=\vertiii{\smash{\Gamma_{\CalQ}}}\cdot\left\Vert f\right\Vert _{\DecompSp{\CalQ}p{\ell_{w}^{q}}v}<\infty.
\end{align*}
Now, we consider the case $p\in\left(0,1\right)$. We observe that
\cite[Lemma 2.7]{DecompositionEmbedding} yields some $R=R\left(R_{\CalQ},C_{\CalQ}\right)>0$
satisfying $\overline{Q_{i}^{\ast}}\subset T_{i}\overline{B_{R}\left(0\right)}+b_{i}\subset T_{i}\left[-R,R\right]^{\dimension}+b_{i}$
for all $i\in I$. Because of $\supp\left(\varphi_{i}^{\ast}\widehat{f}\right)\subset\overline{Q_{i}^{\ast}}$,
Theorem \ref{thm:BandlimitedWienerAmalgamSelfImproving} thus yields
a constant $C_{1}=C_{1}\left(\dimension,p,R,\Omega_{0},\Omega_{1},K\right)>0$
such that 
\[
\left\Vert \Fourier^{-1}\left(\varphi_{i}^{\ast}\cdot\smash{\widehat{f}}\:\right)\right\Vert _{V_{i}}=\left\Vert \Fourier^{-1}\left(\varphi_{i}^{\ast}\cdot\smash{\widehat{f}}\:\right)\right\Vert _{W_{T_{i}^{-T}\left[-1,1\right]^{\dimension}}\left(L_{v}^{p}\right)}\leq C_{1}\cdot\left\Vert \Fourier^{-1}\left(\varphi_{i}^{\ast}\cdot\smash{\widehat{f}}\:\right)\right\Vert _{L_{v}^{p}}\qquad\forall i\in I.
\]

Next, since $L_{v}^{p}\left(\R^{\dimension}\right)$ is a Quasi-Banach
space and since we have the uniform estimate $\left|i^{\ast}\right|\leq N_{\CalQ}$
for all $i\in I$, there is a constant $C_{2}=C_{2}\left(N_{\CalQ},p\right)>0$
satisfying
\[
\left\Vert \Fourier^{-1}\left(\varphi_{i}^{\ast}\cdot\smash{\widehat{f}}\:\right)\right\Vert _{L_{v}^{p}}\leq C_{2}\cdot\sum_{\ell\in i^{\ast}}\left\Vert \Fourier^{-1}\left(\varphi_{\ell}\cdot\smash{\widehat{f}}\:\right)\right\Vert _{L_{v}^{p}}=C_{2}\cdot\left(\Gamma_{\CalQ}\,c\right)_{i}\qquad\forall i\in I.
\]
All in all, this entails by solidity of $\ell_{w}^{q}\left(I\right)$
that
\[
\left\Vert {\rm Ana}_{\ast}f\right\Vert _{V}=\left\Vert \left(\left\Vert \Fourier^{-1}\left(\varphi_{i}^{\ast}\cdot\smash{\widehat{f}}\:\right)\right\Vert _{V_{i}}\right)_{i\in I}\right\Vert _{\ell_{w}^{q}}\leq C_{1}C_{2}\cdot\left\Vert \Gamma_{\CalQ}\,c\right\Vert _{\ell_{w}^{q}}\leq C_{1}C_{2}\vertiii{\smash{\Gamma_{\CalQ}}}\cdot\left\Vert f\right\Vert _{\DecompSp{\CalQ}p{\ell_{w}^{q}}v}<\infty,
\]
as above. In summary, ${\rm Ana}_{\ast}$ is well-defined and bounded
for all $p\in\left(0,\infty\right]$.

Now, using the map ${\rm Synth}_{\CalD}$ from Lemma \ref{lem:DecompositionSynthesis},
we have because of $\varphi_{i}^{\ast}\varphi_{i}=\varphi_{i}$ that
\[
\left({\rm Synth}_{\CalD}\circ{\rm Ana}_{\ast}\right)f=\sum_{i\in I}\Fourier^{-1}\left(\varphi_{i}\cdot\Fourier\left[{\rm Ana}_{\ast}f\right]_{i}\right)=\sum_{i\in I}\Fourier^{-1}\left(\varphi_{i}\varphi_{i}^{\ast}\cdot\smash{\widehat{f}}\:\right)=\sum_{i\in I}\Fourier^{-1}\left(\varphi_{i}\cdot\smash{\widehat{f}}\:\right)=f
\]
for all $f\in\DecompSp{\CalQ}p{\ell_{w}^{q}}v$. Finally, recall from
(the remark after) Lemma \ref{lem:IteratedSequenceSpaceComplete}
that $V=\ell_{w}^{q}\left(\left[V_{i}\right]_{i\in I}\right)$ is
complete.

Now, let $\left(f_{n}\right)_{n\in\N}$ be Cauchy in $\DecompSp{\CalQ}p{\ell_{w}^{q}}v$.
Since ${\rm Ana}_{\ast}$ is bounded, the sequence $\left(g_{n}\right)_{n\in\N}$
with $g_{n}:={\rm Ana}_{\ast}f_{n}$ is Cauchy in $V$. Hence, $g_{n}\to g$
for some $g\in V$. Define $f:={\rm Synth}_{\CalD}g\in\DecompSp{\CalQ}p{\ell_{w}^{q}}v$
and observe
\[
\left\Vert f_{n}-f\right\Vert _{\DecompSp{\CalQ}p{\ell_{w}^{q}}v}=\left\Vert {\rm Synth}_{\CalD}{\rm Ana}_{\ast}f_{n}-{\rm Synth}_{\CalD}g\right\Vert _{\DecompSp{\CalQ}p{\ell_{w}^{q}}v}\leq\vertiii{{\rm Synth}_{\CalD}}\cdot\left\Vert g_{n}-g\right\Vert _{V}\xrightarrow[n\to\infty]{}0.\qedhere
\]
\end{proof}
Using all of the technical lemmata collected in this section, we can
finally prove that the family $\left(L_{\delta\cdot T_{j}^{-T}k}\gamma^{\left[j\right]}\right)_{k\in\Z^{\dimension},j\in I}$
yields an atomic decomposition of $\DecompSp{\CalQ}p{\ell_{w}^{q}}v$,
if $\delta>0$ is chosen small enough.
\begin{thm}
\label{thm:AtomicDecomposition}Assume that the families $\Gamma=\left(\gamma_{i}\right)_{i\in I}$
and $\Gamma_{\ell}=\left(\gamma_{i,\ell}\right)_{i\in I}$ with $\ell\in\left\{ 1,2\right\} $
satisfy Assumption \ref{assu:AtomicDecompositionAssumption} and that
$\Gamma=\left(\gamma_{i}\right)_{i\in I}$ satisfies Assumption \ref{assu:GammaCoversOrbit}.
Define $\delta_{0}>0$ by
\[
\delta_{0}^{-1}\!:=\!\begin{cases}
\!\frac{2s_{\dimension}}{\sqrt{\dimension}}\cdot\left(2^{17}\!\cdot\!\dimension^{2}\!\cdot\!\left(K\!+\!2\!+\!\dimension\right)\right)^{\!K+\dimension+3}\!\!\!\cdot\left(1\!+\!R_{\CalQ}\right)^{\dimension+1}\cdot\Omega_{0}^{4K}\Omega_{1}^{4}\Omega_{2}^{\left(p,K\right)}\Omega_{4}^{\left(p,K\right)}\cdot\vertiii{\smash{\overrightarrow{C}}}\,, & \text{if }p\geq1,\\
\frac{\left(2^{14}/\dimension^{\frac{3}{2}}\right)^{\!\frac{\dimension}{p}}}{2^{45}\cdot\dimension^{17}}\!\cdot\!\left(\frac{s_{\dimension}}{p}\right)^{\!\frac{1}{p}}\left(2^{68}\!\cdot\!\dimension^{14}\!\cdot\!\left(K\!+\!1\!+\!\frac{\dimension+1}{p}\right)^{3}\right)^{\!K+2+\frac{\dimension+1}{p}}\!\!\!\cdot\!\left(1\!+\!R_{\CalQ}\right)^{1+\frac{3\dimension}{p}}\!\cdot\!\Omega_{0}^{16K}\Omega_{1}^{16}\Omega_{2}^{\left(p,K\right)}\Omega_{4}^{\left(p,K\right)}\cdot\vertiii{\smash{\overrightarrow{C}}}^{\frac{1}{p}}, & \text{if }p<1.
\end{cases}
\]
Then, for each $0<\delta\leq\min\left\{ 1,\delta_{0}\right\} $, the
family $\left(L_{\delta\cdot T_{j}^{-T}k}\:\gamma^{\left[j\right]}\right)_{j\in I,k\in\Z^{\dimension}}$
forms an \textbf{atomic decomposition} of $\DecompSp{\CalQ}p{\ell_{w}^{q}}v$.
Precisely, this means the following:

\begin{enumerate}
\item The \textbf{synthesis map}
\begin{align*}
\qquad\qquad S^{\left(\delta\right)}:\ell_{\left(\left|\det T_{j}\right|^{\frac{1}{2}-\frac{1}{p}}w_{j}\right)_{\!\!j\in I}}^{q}\!\!\!\!\!\!\!\left(\left[\vphantom{F}\smash{C_{j}^{\left(\delta\right)}}\right]_{j\in I}\right) & \to\DecompSp{\CalQ}p{\ell_{w}^{q}}v,\\
\left(\smash{c_{k}^{\left(j\right)}}\right)_{j\in I,k\in\Z^{\dimension}} & \mapsto\sum_{j\in I}\sum_{k\in\Z^{\dimension}}\left(\left|\det T_{j}\right|^{-\frac{1}{2}}c_{k}^{\left(j\right)}\cdot L_{\delta\cdot T_{j}^{-T}k}\gamma^{\left(j\right)}\right)=\sum_{j\in I}\sum_{k\in\Z^{\dimension}}\left(c_{k}^{\left(j\right)}\cdot L_{\delta\cdot T_{j}^{-T}k}\gamma^{\left[j\right]}\right)
\end{align*}
is well-defined and bounded for each $\delta\in\left(0,1\right]$.
\item For $0<\delta\leq\min\left\{ 1,\delta_{0}\right\} $, there is a bounded
linear \textbf{coefficient map}
\[
C^{\left(\delta\right)}:\DecompSp{\CalQ}p{\ell_{w}^{q}}v\to\ell_{\left(\left|\det T_{j}\right|^{\frac{1}{2}-\frac{1}{p}}w_{j}\right)_{\!\!j\in I}}^{q}\!\!\!\!\!\!\left(\left[\vphantom{F}\smash{C_{j}^{\left(\delta\right)}}\right]_{j\in I}\right)
\]
satisfying $S^{\left(\delta\right)}\circ C^{\left(\delta\right)}=\identity_{\DecompSp{\CalQ}p{\ell_{w}^{q}}v}$.\qedhere
\end{enumerate}
\end{thm}
\begin{rem*}
As the proof shows, convergence of the series $\sum_{j\in I}\sum_{k\in\Z^{\dimension}}\left(\left|\det T_{j}\right|^{-\frac{1}{2}}c_{k}^{\left(j\right)}\cdot L_{\delta\cdot T_{j}^{-T}k}\gamma^{\left(j\right)}\right)$
has to be understood as follows: Each of the series
\[
\sum_{k\in\Z^{\dimension}}\left(\left|\det T_{j}\right|^{-\frac{1}{2}}c_{k}^{\left(j\right)}\cdot L_{\delta\cdot T_{j}^{-T}k}\gamma^{\left(j\right)}\right)
\]
converges pointwise absolutely to a function $g_{j}\in V_{j}\hookrightarrow\Schwartz'\left(\R^{\dimension}\right)$
and the series $\sum_{j\in I}g_{j}$ converges in the weak-$\ast$-sense
in $Z'\left(\CalO\right)$, i.e., for every $\phi\in Z\left(\CalO\right)$,
the series $\sum_{j\in I}\left\langle g_{j},\,\phi\right\rangle _{\Schwartz',\Schwartz}$
converges absolutely and the functional $\phi\mapsto\sum_{j\in I}\left\langle g_{j},\,\phi\right\rangle _{\Schwartz',\Schwartz}$
is continuous on $Z\left(\CalO\right)$.

Furthermore, the proof shows that the definition of $C^{\left(\delta\right)}$
is independent of the precise choice of $p,q,v,w$, as long as $\delta>0$
is chosen small enough that $C^{\left(\delta\right)}$ is defined
at all. In fact, the proof shows that $C^{\left(\delta\right)}=D^{\left(\delta\right)}\cdot\left(T^{\left(\delta\right)}\right)^{-1}$
where $T^{\left(\delta\right)}=S^{\left(\delta\right)}\circ D^{\left(\delta\right)}:\DecompSp{\CalQ}p{\ell_{w}^{q}}v\to\DecompSp{\CalQ}p{\ell_{w}^{q}}v$
is invertible (using a Neumann series) for $0<\delta\leq\min\left\{ 1,\delta_{0}\right\} $,
with
\begin{align*}
D^{\left(\delta\right)}:\DecompSp{\CalQ}p{\ell_{w}^{q}}v & \to\ell_{\left(\left|\det T_{j}\right|^{\frac{1}{2}-\frac{1}{p}}w_{j}\right)_{\!\!j\in I}}^{q}\!\!\!\!\!\!\!\!\left(\left[\vphantom{F}\smash{C_{j}^{\left(\delta\right)}}\right]_{j\in I}\right),\\
f & \mapsto\left[\left(\delta^{\dimension}\cdot\left|\det T_{j}\right|^{-1/2}\cdot\left[\Fourier^{-1}\!\left(\theta_{j}\varphi_{j}\cdot\widehat{f}\right)\right]\!\!\left(\delta\cdot T_{j}^{-T}k\right)\right)_{k\in\Z^{\dimension}}\right]_{j\in I},
\end{align*}
where $\theta_{j}$ for $j\in I$ is defined as in Assumption \ref{assu:GammaCoversOrbit}.
\end{rem*}
\begin{proof}
We first study for arbitrary $j\in I$ and $\delta\in\left(0,1\right]$
boundedness (and well-definedness) of the map
\[
S_{\Gamma_{2}}^{\left(\delta,j\right)}:C_{j}^{\left(\delta\right)}\to V_{j},\left(c_{k}\right)_{k\in\Z^{\dimension}}\mapsto\left|\det T_{j}\right|^{-1/2}\cdot\sum_{k\in\Z^{\dimension}}c_{k}\cdot L_{\delta\cdot T_{j}^{-T}k}\:\gamma_{2}^{\left(j\right)}.
\]
Recall $\gamma_{2}^{\left(j\right)}=\left|\det T_{j}\right|\cdot M_{b_{j}}\left[\gamma_{j,2}\circ T_{j}^{T}\right]$,
so that
\begin{align}
\left[S_{\Gamma_{2}}^{\left(\delta,j\right)}\left(c_{k}\right)_{k\in\Z^{\dimension}}\right]\left(x\right) & =\left|\det T_{j}\right|^{-1/2}\cdot\sum_{k\in\Z^{\dimension}}c_{k}\cdot\gamma_{2}^{\left(j\right)}\left(x-\delta\cdot T_{j}^{-T}k\right)\nonumber \\
 & =\left|\det T_{j}\right|^{1/2}\cdot\sum_{k\in\Z^{\dimension}}e^{2\pi i\left\langle b_{j},x-\delta\cdot T_{j}^{-T}k\right\rangle }\cdot c_{k}\cdot\gamma_{j,2}\left(T_{j}^{T}x-\delta\cdot k\right)\nonumber \\
 & =\left|\det T_{j}\right|^{1/2}\cdot e^{2\pi i\left\langle b_{j},x\right\rangle }\cdot\left(\sum_{k\in\Z^{\dimension}}e^{-2\pi i\left\langle b_{j},\delta\cdot T_{j}^{-T}k\right\rangle }c_{k}\cdot L_{\delta\cdot k}\gamma_{j,2}\right)\!\!\left(T_{j}^{T}x\right)\nonumber \\
\left({\scriptstyle \text{with }\Psi_{\gamma_{j,2}}^{\left(j,\delta\right)}\text{ as in Lemma }\ref{lem:SchwartzTranslationSynthesis}}\right) & =\left|\det T_{j}\right|^{1/2}\cdot e^{2\pi i\left\langle b_{j},x\right\rangle }\cdot\Psi_{\gamma_{j,2}}^{\left(j,\delta\right)}\left[\left(e^{-2\pi i\left\langle b_{j},\delta\cdot T_{j}^{-T}k\right\rangle }c_{k}\right)_{k\in\Z^{\dimension}}\right].\label{eq:AtomicDecompositionTranslationSynthesisNormalization}
\end{align}
Thus, in terms of the isometric isomorphism $m_{j}^{\left(\delta\right)}:C_{j}^{\left(\delta\right)}\to C_{j}^{\left(\delta\right)},\left(c_{k}\right)_{k\in\Z^{\dimension}}\mapsto\left(e^{-2\pi i\left\langle b_{j},\delta\cdot T_{j}^{-T}k\right\rangle }\cdot c_{k}\right)_{k\in\Z^{\dimension}}$
and of the map $\Psi_{\gamma_{j,2}}^{\left(j,\delta\right)}$ defined
in Lemma \ref{lem:SchwartzTranslationSynthesis}, the preceding calculations
show 
\[
S_{\Gamma_{2}}^{\left(\delta,j\right)}c=\left|\det T_{j}\right|^{1/2}\cdot M_{b_{j}}\left[\Psi_{\gamma_{j,2}}^{\left(j,\delta\right)}\left(m_{j}^{\left(\delta\right)}c\right)\right]\qquad\forall c\in C_{j}^{\left(\delta\right)}.
\]
As a consequence of the solidity of $V_{j}$ and of Lemma \ref{lem:SchwartzTranslationSynthesis}
(which is applicable, since $\left\Vert \gamma_{j,2}\right\Vert _{K_{0}}\leq\Omega_{4}^{\left(p,K\right)}<\infty$
for all $j\in I$, cf.\@ Assumption \ref{assu:AtomicDecompositionAssumption}),
we thus get 
\begin{equation}
\vertiii{S_{\Gamma_{2}}^{\left(\delta,j\right)}}\leq C_{K,\delta,\dimension,p,\Gamma_{2}}\cdot\left|\det T_{j}\right|^{\frac{1}{2}-\frac{1}{p}}<\infty\qquad\forall j\in I\label{eq:DiscreteSynthesisUniformBound}
\end{equation}
for a suitable constant $C_{K,\delta,\dimension,p,\Gamma_{2}}>0$
which is independent of $j\in I$. In particular, each map $S_{\Gamma_{2}}^{\left(\delta,j\right)}$
is well-defined with pointwise absolute convergence of the defining
series.

\medskip{}

Now, we can establish boundedness of the synthesis map $S^{\left(\delta\right)}$
as follows: In view of equation (\ref{eq:DiscreteSynthesisUniformBound}),
it follows that
\[
\bigotimes_{j\in I}S_{\Gamma_{2}}^{\left(\delta,j\right)}:\ell_{\left(\left|\det T_{j}\right|^{\frac{1}{2}-\frac{1}{p}}\cdot w_{j}\right)_{\!\!j\in I}}^{q}\!\!\!\!\!\!\!\!\left(\left[\vphantom{F}\smash{C_{j}^{\left(\delta\right)}}\right]_{j\in I}\right)\to\ell_{w}^{q}\left(\left[V_{j}\right]_{j\in I}\right),\left(\smash{c_{k}^{\left(j\right)}}\right)_{j\in I,k\in\Z^{\dimension}}\mapsto\left[S_{\Gamma_{2}}^{\left(\delta,j\right)}\left(\smash{c_{k}^{\left(j\right)}}\right)_{k\in\Z^{\dimension}}\right]_{j\in I}
\]
is well-defined and bounded, with $\vertiii{\bigotimes_{j\in I}S_{\Gamma_{2}}^{\left(\delta,j\right)}}\leq C_{K,\delta,\dimension,p,\Gamma_{2}}$.
Furthermore, using $\gamma_{j}=\gamma_{j,1}\ast\gamma_{j,2}$, it
follows easily that $L_{x}\gamma^{\left(j\right)}=\gamma_{1}^{\left(j\right)}\ast L_{x}\gamma_{2}^{\left(j\right)}$
for arbitrary $x\in\R^{\dimension}$ and $j\in I$, from which it
follows (with ${\rm Synth}_{\Gamma_{1}}$ as in Lemma \ref{lem:GammaSynthesisBounded})
that
\begin{align}
\left[{\rm Synth}_{\Gamma_{1}}\circ\bigotimes_{j\in I}S_{\Gamma_{2}}^{\left(\delta,j\right)}\right]\left(\smash{c_{k}^{\left(j\right)}}\right)_{j\in I,k\in\Z^{\dimension}} & =\sum_{j\in I}\left[\gamma_{1}^{\left(j\right)}\ast S_{\Gamma_{2}}^{\left(\delta,j\right)}\left(\smash{c_{k}^{\left(j\right)}}\right)_{k\in\Z^{\dimension}}\right]\nonumber \\
 & =\sum_{j\in I}\left[\left|\det T_{j}\right|^{-1/2}\cdot\gamma_{1}^{\left(j\right)}\ast\sum_{k\in\Z^{\dimension}}c_{k}^{\left(j\right)}\cdot L_{\delta\cdot T_{j}^{-T}k}\gamma_{2}^{\left(j\right)}\right]\nonumber \\
\left({\scriptstyle \gamma_{1}^{\left(j\right)}\in L_{v_{0}}^{1}\left(\R^{\dimension}\right)\text{ and (proof of) Lemma }\ref{lem:SchwartzTranslationSynthesis}}\right) & =\sum_{j\in I}\left[\left|\det T_{j}\right|^{-1/2}\cdot\sum_{k\in\Z^{\dimension}}c_{k}^{\left(j\right)}\cdot L_{\delta\cdot T_{j}^{-T}k}\left(\gamma_{1}^{\left(j\right)}\ast\gamma_{2}^{\left(j\right)}\right)\right]\nonumber \\
 & =\sum_{j\in I}\left[\sum_{k\in\Z^{\dimension}}\left|\det T_{j}\right|^{-1/2}\cdot c_{k}^{\left(j\right)}\cdot L_{\delta\cdot T_{j}^{-T}k}\gamma^{\left(j\right)}\right]\nonumber \\
 & =S^{\left(\delta\right)}\left(\smash{c_{k}^{\left(j\right)}}\right)_{j\in I,k\in\Z^{\dimension}}.\label{eq:AtomicDecompositionSynthesisOperatorAsComposition}
\end{align}
This shows for arbitrary $\delta\in\left(0,1\right]$ that $S^{\left(\delta\right)}={\rm Synth}_{\Gamma_{1}}\circ\bigotimes_{j\in I}S_{\Gamma_{2}}^{\left(\delta,j\right)}$
is bounded, as a composition of bounded maps.

\medskip{}

Finally, we prove that convergence of the series defining $S^{\left(\delta\right)}\left(\smash{c_{k}^{\left(j\right)}}\right)_{j\in I,k\in\Z^{\dimension}}$
occurs in the sense described in the remark following the theorem:
First of all, we get exactly as in equation (\ref{eq:AtomicDecompositionTranslationSynthesisNormalization})
(but with $\gamma_{j}$ instead of $\gamma_{j,2}$) that
\[
g_{j}:=\sum_{k\in\Z^{\dimension}}\left(\left|\det T_{j}\right|^{-\frac{1}{2}}c_{k}^{\left(j\right)}\cdot L_{\delta\cdot T_{j}^{-T}k}\gamma^{\left(j\right)}\right)=\left|\det T_{j}\right|^{1/2}\cdot M_{b_{j}}\left[\Psi_{\gamma_{j}}^{\left(j,\delta\right)}\left(m_{j}^{\left(\delta\right)}\left(\smash{c_{k}^{\left(j\right)}}\right)_{k\in\Z^{\dimension}}\right)\right].
\]
Here, it is worth mentioning that the conditions in Assumption \ref{assu:AtomicDecompositionAssumption}
imply that Lemma \ref{lem:SchwartzTranslationSynthesis} is applicable
with $\varrho=\gamma_{j,2}$, as well as with $\varrho=\gamma_{j}$.
Hence, that lemma shows that the series defining $g_{j}$ converges
pointwise absolutely and that $g_{j}\in V_{j}\hookrightarrow\Schwartz'\left(\R^{\dimension}\right)$,
where the last embedding is justified by Lemma \ref{lem:SpecialConvolutionConsistent}.

Finally, Lemma \ref{lem:GammaSynthesisBounded} shows that for $\left(h_{j}\right)_{j\in I}:=\left(\bigotimes_{j\in I}S_{\Gamma_{2}}^{\left(\delta,j\right)}\right)\left(\smash{c_{k}^{\left(j\right)}}\right)_{j\in I,k\in\Z^{\dimension}}\in\ell_{w}^{q}\left(\left[V_{j}\right]_{j\in I}\right)$
and each $\phi\in Z\left(\CalO\right)$, the series
\begin{align*}
\left\langle {\rm Synth}_{\Gamma_{1}}\left(h_{j}\right)_{j\in I},\,\phi\right\rangle _{Z'\left(\CalO\right),Z\left(\CalO\right)} & =\sum_{j\in I}\left\langle \gamma_{1}^{\left(j\right)}\ast h_{j},\,\phi\right\rangle _{\Schwartz',\Schwartz}\\
 & =\sum_{j\in I}\left\langle \gamma_{1}^{\left(j\right)}\ast\left[\left|\det T_{j}\right|^{-1/2}\cdot\sum_{k\in\Z^{\dimension}}c_{k}^{\left(j\right)}L_{\delta\cdot T_{j}^{-T}k}\gamma_{2}^{\left(j\right)}\right],\,\phi\right\rangle _{\Schwartz',\Schwartz}\\
\left({\scriptstyle \text{cf. eq. }\eqref{eq:AtomicDecompositionSynthesisOperatorAsComposition}}\right) & =\sum_{j\in I}\left\langle \left|\det T_{j}\right|^{-1/2}\cdot\sum_{k\in\Z^{\dimension}}c_{k}^{\left(j\right)}L_{\delta\cdot T_{j}^{-T}k}\gamma^{\left(j\right)},\,\phi\right\rangle _{\Schwartz',\Schwartz}\\
 & =\sum_{j\in I}\left\langle g_{j},\,\phi\right\rangle _{\Schwartz',\Schwartz}
\end{align*}
converges absolutely and defines a continuous functional on $Z\left(\CalO\right)$.

\medskip{}

Next, we want to show existence of the coefficient map $C^{\left(\delta\right)}$,
for $0<\delta\leq\min\left\{ 1,\delta_{0}\right\} $. To this end,
first note that Theorem \ref{thm:BandlimitedWienerAmalgamSelfImproving}
shows for
\[
C_{1}:=\begin{cases}
1, & \text{if }p\geq1,\\
2^{4\left(1+\frac{\dimension}{p}\right)}s_{\dimension}^{\frac{1}{p}}\left(192\cdot\dimension^{3/2}\cdot\left\lceil K+\frac{\dimension+1}{p}\right\rceil \right)^{\left\lceil K+\frac{\dimension+1}{p}\right\rceil +1}\cdot\Omega_{0}^{K}\Omega_{1}\cdot\left(1+R_{\CalQ}\right)^{\frac{\dimension}{p}}, & \text{if }p<1
\end{cases}
\]
that
\[
\left\Vert \Fourier^{-1}\left(\varphi_{j}\cdot\widehat{f}\right)\right\Vert _{V_{j}}\leq C_{1}\cdot\left\Vert \Fourier^{-1}\left(\varphi_{j}\cdot\widehat{f}\right)\right\Vert _{L_{v}^{p}}\qquad\forall j\in I\qquad\forall f\in Z'\left(\CalO\right),
\]
since we have $\supp\left(\varphi_{j}\cdot\widehat{f}\right)\subset\overline{Q_{j}}\subset T_{j}\left[\overline{B_{R_{\CalQ}}}\left(0\right)\right]+b_{j}\subset T_{j}\left[-R_{\CalQ},R_{\CalQ}\right]^{\dimension}+b_{j}$.
This easily shows that the map
\[
{\rm Ana}_{\varphi}:\DecompSp{\CalQ}p{\ell_{w}^{q}}v\to\ell_{w}^{q}\left(\left[V_{j}\right]_{j\in I}\right),f\mapsto\left(\Fourier^{-1}\left[\varphi_{j}\cdot\widehat{f}\right]\right)_{j\in I}
\]
is well-defined and bounded, with $\vertiii{{\rm Ana}_{\varphi}}\leq C_{1}$.

Furthermore, Lemma \ref{lem:LocalInverseConvolution} shows that the
map
\[
m_{\theta}:\ell_{w}^{q}\left(\left[V_{j}\right]_{j\in I}\right)\to\ell_{w}^{q}\left(\left[V_{j}\right]_{j\in I}\right),\left(f_{j}\right)_{j\in I}\mapsto\left[\left(\Fourier^{-1}\theta_{j}\right)\ast f_{j}\right]_{j\in I}\;\overset{\text{Lem. }\ref{lem:SpecialConvolutionConsistent}}{=}\;\left[\Fourier^{-1}\left(\theta_{j}\cdot\widehat{f_{j}}\right)\right]_{j\in I}
\]
is well-defined and bounded, with $\vertiii{m_{\theta}}\leq C_{2}<\infty$
for
\[
C_{2}:=\begin{cases}
\Omega_{0}^{4K}\Omega_{1}^{4}\Omega_{2}^{\left(p,K\right)}\cdot\dimension^{-\frac{\dimension}{2p}}\cdot\left(972\cdot\dimension^{5/2}\right)^{K+\frac{\dimension}{p}}, & \text{if }p\in\left(0,1\right),\\
\Omega_{0}^{K}\Omega_{1}\Omega_{2}^{\left(p,K\right)}, & \text{if }p\in\left[1,\infty\right].
\end{cases}
\]

Next, it follows from Lemma \ref{lem:GammaSynthesisBounded} that
the map
\[
{\rm Synth}_{\Gamma_{1}}:\ell_{w}^{q}\left(\left[V_{j}\right]_{j\in I}\right)\to\DecompSp{\CalQ}p{\ell_{w}^{q}}v,\left(g_{j}\right)_{j\in I}\mapsto\sum_{j\in I}\gamma_{1}^{\left(j\right)}\ast g_{j}=\sum_{j\in I}\Fourier^{-1}\left(\widehat{\gamma_{1}^{\left(j\right)}}\cdot\widehat{g_{j}}\right)
\]
is well-defined and bounded with $\vertiii{{\rm Synth}_{\Gamma_{1}}}\leq C_{3}\cdot\vertiii{\smash{\overrightarrow{C}}}^{\max\left\{ 1,1/p\right\} }$,
with
\[
C_{3}:=\begin{cases}
1, & \text{if }p\geq1\\
\frac{\left(2^{6}/\sqrt{\dimension}\right)^{\frac{\dimension}{p}}}{2^{21}\cdot\dimension^{7}}\!\cdot\!\left(2^{21}\!\cdot\!\dimension^{5}\!\cdot\!\left\lceil K\!+\!\frac{\dimension+1}{p}\right\rceil \right)^{\!\left\lceil K+\frac{\dimension+1}{p}\right\rceil +1}\!\!\!\cdot\!\left(1\!+\!R_{\CalQ}\right)^{\frac{\dimension}{p}}\!\cdot\!\Omega_{0}^{5K}\Omega_{1}^{5}, & \text{if }p<1.
\end{cases}
\]

Finally, we will show below that the map
\[
m_{\Gamma_{2}}:\ell_{w}^{q}\left(\left[V_{j}\right]_{j\in I}\right)\to\ell_{w}^{q}\left(\left[V_{j}\right]_{j\in I}\right),\left(f_{j}\right)_{j\in I}\mapsto\left(\gamma_{2}^{\left(j\right)}\ast f_{j}\right)_{j\in I}
\]
is also well-defined and bounded. Once this is shown, note that we
have
\begin{align}
\left({\rm Synth}_{\Gamma_{1}}\circ m_{\Gamma_{2}}\circ m_{\theta}\circ{\rm Ana}_{\varphi}\right)f & =\sum_{j\in I}\left(\gamma_{1}^{\left(j\right)}\ast\gamma_{2}^{\left(j\right)}\ast\Fourier^{-1}\theta_{j}\ast\Fourier^{-1}\left(\varphi_{j}\cdot\widehat{f}\right)\right)\nonumber \\
\left({\scriptstyle \text{Lemma }\ref{lem:SpecialConvolutionConsistent}}\right) & =\sum_{j\in I}\Fourier^{-1}\left(\widehat{\gamma_{1}^{\left(j\right)}}\cdot\widehat{\gamma_{2}^{\left(j\right)}}\cdot\theta_{j}\cdot\varphi_{j}\cdot\widehat{f}\right)\nonumber \\
\left({\scriptstyle \text{easy consequence of }\gamma_{j}=\gamma_{j,1}\ast\gamma_{j,2}}\right) & =\sum_{j\in I}\Fourier^{-1}\left(\widehat{\gamma^{\left(j\right)}}\cdot\theta_{j}\cdot\varphi_{j}\cdot\widehat{f}\right)\nonumber \\
\left({\scriptstyle \text{since }\widehat{\gamma^{\left(j\right)}}\cdot\theta_{j}\equiv1\text{ on }\overline{Q_{j}}\supset\supp\varphi_{j}}\right) & =\sum_{j\in I}\Fourier^{-1}\left(\varphi_{j}\cdot\widehat{f}\right)\nonumber \\
\left({\scriptstyle \text{since }\widehat{f}\in\DistributionSpace{\CalO}\text{ and }\left(\varphi_{j}\right)_{j\in I}\text{ is locally finite part. of unity on }\CalO}\right) & =f\label{eq:AtomicDecompositionReproducingFormula}
\end{align}
for all $f\in\DecompSp{\CalQ}p{\ell_{w}^{q}}v$. Thus, our goal in
the remainder of the proof—once we have shown boundedness of $m_{\Gamma_{2}}$—will
be to discretize this \textbf{reproducing formula}.

But first of all, let us verify boundedness of $m_{\Gamma_{2}}$.
To this end, it suffices to show that each map
\[
J_{j}:V_{j}\to V_{j},f\mapsto\gamma_{2}^{\left(j\right)}\ast f
\]
is bounded, with $\sup_{j\in I}\vertiii{J_{j}}<\infty$. But for $p\in\left[1,\infty\right]$,
this simply follows from the weighted Young inequality (equation (\ref{eq:WeightedYoungInequality})),
since in this case, we have $K_{0}=K+\dimension+1$ and thus (cf.\@
equation (\ref{eq:AtomicDecompositionFamilyDefinition}))
\begin{align}
\left\Vert \gamma_{2}^{\left(j\right)}\right\Vert _{L_{v_{0}}^{1}} & =\left\Vert v_{0}\cdot\left|\det T_{j}\right|\cdot M_{b_{j}}\left[\gamma_{j,2}\circ T_{j}^{T}\right]\right\Vert _{L^{1}}\nonumber \\
 & =\left\Vert \left(v_{0}\circ T_{j}^{-T}\right)\cdot\gamma_{j,2}\right\Vert _{L^{1}}\nonumber \\
\left({\scriptstyle \text{assump. on }v_{0}\text{ and eq. }\eqref{eq:WeightLinearTransformationsConnection}}\right) & \leq\Omega_{0}^{K}\Omega_{1}\cdot\left\Vert \left(1+\left|\mybullet\right|\right)^{K}\cdot\gamma_{j,2}\right\Vert _{L^{1}}\nonumber \\
 & \leq\Omega_{0}^{K}\Omega_{1}\cdot\left\Vert \gamma_{j,2}\right\Vert _{K_{0}}\cdot\left\Vert \left(1+\left|\mybullet\right|\right)^{K-K_{0}}\right\Vert _{L^{1}}\nonumber \\
\left({\scriptstyle \text{eq. }\eqref{eq:StandardDecayLpEstimate}}\right) & \leq\Omega_{0}^{K}\Omega_{1}\Omega_{4}^{\left(p,K,1\right)}\cdot s_{\dimension}<\infty.\label{eq:AtomicDecompositionGamma2WeightedL1Norm}
\end{align}
Here, we defined $\Omega_{4}^{\left(p,K,1\right)}:=\sup_{j\in I}\left\Vert \gamma_{j,2}\right\Vert _{K_{0}}$
in the last step. Note that $\Omega_{4}^{\left(p,K,1\right)}\leq\Omega_{4}^{\left(p,K\right)}$,
cf.\@ Assumption \ref{assu:AtomicDecompositionAssumption}, equation
(\ref{eq:AtomicDecompositionGamma2ConstantDefinition}).

Likewise, for $p\in\left(0,1\right)$, we can simply use Corollary
\ref{cor:WienerAmalgamConvolutionSimplified} to derive for $C_{4}:=\Omega_{0}^{3K}\Omega_{1}^{3}\dimension^{-\frac{\dimension}{2p}}\cdot\left(972\cdot\dimension^{5/2}\right)^{K+\frac{\dimension}{p}}$
that
\begin{align*}
\left\Vert \gamma_{2}^{\left(j\right)}\ast f\right\Vert _{V_{j}} & =\left\Vert \gamma_{2}^{\left(j\right)}\ast f\right\Vert _{W_{T_{j}^{-T}\left[-1,1\right]^{\dimension}}\left(L_{v}^{p}\right)}\\
\left({\scriptstyle \text{Cor. }\ref{cor:WienerAmalgamConvolutionSimplified}}\right) & \leq C_{4}\cdot\left|\det T_{j}\right|^{\frac{1}{p}-1}\cdot\left\Vert \gamma_{2}^{\left(j\right)}\right\Vert _{W_{T_{j}^{-T}\left[-1,1\right]^{\dimension}}\left(L_{v_{0}}^{p}\right)}\cdot\left\Vert f\right\Vert _{W_{T_{j}^{-T}\left[-1,1\right]^{\dimension}}\left(L_{v}^{p}\right)}\\
\left({\scriptstyle \text{eq. }\eqref{eq:AtomicDecompositionFamilyDefinition}}\right) & =C_{4}\cdot\left|\det T_{j}\right|^{\frac{1}{p}}\cdot\left\Vert M_{b_{j}}\left[\gamma_{j,2}\circ T_{j}^{T}\right]\right\Vert _{W_{T_{j}^{-T}\left[-1,1\right]^{\dimension}}\left(L_{v_{0}}^{p}\right)}\cdot\left\Vert f\right\Vert _{V_{j}}\\
 & =C_{4}\cdot\left|\det T_{j}\right|^{\frac{1}{p}}\cdot\left\Vert v_{0}\cdot M_{T_{j}^{-T}\left[-1,1\right]^{\dimension}}\left[\gamma_{j,2}\circ T_{j}^{T}\right]\right\Vert _{L^{p}}\cdot\left\Vert f\right\Vert _{V_{j}}\\
\left({\scriptstyle \text{Lemma }\ref{lem:WienerTransformationFormula}}\right) & =C_{4}\cdot\left|\det T_{j}\right|^{\frac{1}{p}}\cdot\left\Vert v_{0}\cdot\left(\left[M_{\left[-1,1\right]^{\dimension}}\gamma_{j,2}\right]\circ T_{j}^{T}\right)\right\Vert _{L^{p}}\cdot\left\Vert f\right\Vert _{V_{j}}\\
 & =C_{4}\cdot\left\Vert \left(v_{0}\circ T_{j}^{-T}\right)\cdot M_{\left[-1,1\right]^{\dimension}}\gamma_{j,2}\right\Vert _{L^{p}}\cdot\left\Vert f\right\Vert _{V_{j}}\\
\left({\scriptstyle \text{assump. on }v_{0}\text{ and eq. }\eqref{eq:WeightLinearTransformationsConnection}}\right) & \leq\Omega_{0}^{K}\Omega_{1}\cdot C_{4}\cdot\left\Vert \left(1+\left|\mybullet\right|\right)^{K}\cdot M_{\left[-1,1\right]^{\dimension}}\gamma_{j,2}\right\Vert _{L^{p}}\cdot\left\Vert f\right\Vert _{V_{j}}\\
\left({\scriptstyle \text{Lemma }\ref{lem:SchwartzFunctionsAreWiener}}\right) & \leq\Omega_{0}^{K}\Omega_{1}\cdot C_{4}\cdot\left(1+2\sqrt{\dimension}\right)^{K_{0}}\cdot\left\Vert \gamma_{j,2}\right\Vert _{K_{0}}\cdot\left\Vert \left(1+\left|\mybullet\right|\right)^{K-K_{0}}\right\Vert _{L^{p}}\cdot\left\Vert f\right\Vert _{V_{j}}\\
\left({\scriptstyle \text{eq. }\eqref{eq:StandardDecayLpEstimate}\text{ and }K-K_{0}=-\left(\frac{\dimension}{p}+1\right)}\right) & \leq\Omega_{0}^{K}\Omega_{1}\Omega_{4}^{\left(p,K,1\right)}\cdot C_{4}\cdot\left(1+2\sqrt{\dimension}\right)^{K_{0}}\cdot\left(\frac{s_{\dimension}}{p}\right)^{1/p}\cdot\left\Vert f\right\Vert _{V_{j}}.
\end{align*}
Here, we used the same definition of $\Omega_{4}^{\left(p,K,1\right)}$
as above. We have thus established boundedness of $m_{\Gamma_{2}}$
in all cases.

\medskip{}

In order to discretize the reproducing formula from equation (\ref{eq:AtomicDecompositionReproducingFormula}),
we define for $\delta\in\left(0,1\right]$ the map
\begin{align*}
D^{\left(\delta\right)}:\DecompSp{\CalQ}p{\ell_{w}^{q}}v & \to\ell_{\left(\left|\det T_{j}\right|^{\frac{1}{2}-\frac{1}{p}}w_{j}\right)_{\!\!j\in I}}^{q}\!\!\!\!\!\!\!\!\left(\left[\vphantom{F}\smash{C_{j}^{\left(\delta\right)}}\right]_{j\in I}\right),\\
f & \mapsto\left[\left(\delta^{\dimension}\cdot\left|\det T_{j}\right|^{-1/2}\cdot\left[\Fourier^{-1}\!\left(\theta_{j}\varphi_{j}\cdot\widehat{f}\right)\right]\!\!\left(\delta\cdot T_{j}^{-T}k\right)\right)_{k\in\Z^{\dimension}}\right]_{j\in I}.
\end{align*}
This map is indeed well-defined and bounded, since Lemma \ref{lem:BandlimitedSampling}
yields for
\[
C_{5}:=2^{\max\left\{ 1,\frac{1}{p}\right\} }\cdot\Omega_{0}^{3K}\Omega_{1}^{3}\cdot\left(1+\sqrt{\dimension}\right)^{K}\cdot\left(23040\cdot\dimension^{3/2}\cdot\left(K+1+\frac{\dimension+1}{\min\left\{ 1,p\right\} }\right)\right)^{K+2+\frac{\dimension+1}{\min\left\{ 1,p\right\} }}\cdot\left(1+R_{\CalQ}\right)^{1+\frac{\dimension}{\min\left\{ 1,p\right\} }}
\]
that
\begin{align*}
\left\Vert D^{\left(\delta\right)}f\right\Vert _{\ell_{\left(\left|\det T_{j}\right|^{\frac{1}{2}-\frac{1}{p}}w_{j}\right)_{j\in I}}^{q}\!\!\!\!\!\!\!\!\left(\vphantom{F}\smash{C_{j}^{\left(\delta\right)}}\right)} & =\delta^{\dimension}\cdot\left\Vert \left(\left|\det T_{j}\right|^{-1/p}\cdot\left\Vert \left(\left[\Fourier^{-1}\!\left(\theta_{j}\varphi_{j}\cdot\widehat{f}\right)\right]\!\left(\delta\cdot T_{j}^{-T}k\right)\right)_{k\in\Z^{\dimension}}\right\Vert _{C_{j}^{\left(\delta\right)}}\right)_{j\in I}\right\Vert _{\ell_{w}^{q}\left(I\right)}\\
\left({\scriptstyle \text{since }\supp\left(\theta_{j}\varphi_{j}\cdot\widehat{f}\right)\subset\overline{Q_{j}}\subset T_{j}\left[-R_{\CalQ},R_{\CalQ}\right]^{\dimension}+b_{j}}\right) & \leq C_{5}\cdot\delta^{\dimension\left(1-\frac{1}{p}\right)}\cdot\left\Vert \left(\left\Vert \Fourier^{-1}\left(\theta_{j}\varphi_{j}\cdot\widehat{f}\right)\right\Vert _{L_{v}^{p}}\right)_{j\in I}\right\Vert _{\ell_{w}^{q}\left(I\right)}\\
 & \leq C_{5}\cdot\delta^{\dimension\left(1-\frac{1}{p}\right)}\cdot\left\Vert \left(\left\Vert \Fourier^{-1}\left(\theta_{j}\varphi_{j}\cdot\widehat{f}\right)\right\Vert _{V_{j}}\right)_{j\in I}\right\Vert _{\ell_{w}^{q}\left(I\right)}\\
 & =C_{5}\cdot\delta^{\dimension\left(1-\frac{1}{p}\right)}\cdot\left\Vert \left(m_{\theta}\circ{\rm Ana}_{\varphi}\right)f\right\Vert _{\ell_{w}^{q}\left(\left[V_{j}\right]_{j\in I}\right)}\\
 & \leq C_{5}\cdot\vertiii{m_{\theta}}\cdot\vertiii{{\rm Ana}_{\varphi}}\cdot\delta^{\dimension\left(1-\frac{1}{p}\right)}\cdot\left\Vert f\right\Vert _{\DecompSp{\CalQ}p{\ell_{w}^{q}}v}<\infty.
\end{align*}
Now, our goal is to show for
\[
E^{\left(\delta\right)}:=\left(m_{\Gamma_{2}}\circ m_{\theta}\circ{\rm Ana}_{\varphi}\right)-\left(\left[\smash{\bigotimes_{j\in I}}\vphantom{\sum_{i}}S_{\Gamma_{2}}^{\left(\delta,j\right)}\right]\circ D^{\left(\delta\right)}\right):\DecompSp{\CalQ}p{\ell_{w}^{q}}v\to\ell_{w}^{q}\left(\left[V_{j}\right]_{j\in I}\right)
\]
that we have $\vertiii{{\rm Synth}_{\Gamma_{1}}}\cdot\vertiii{E^{\left(\delta\right)}}\leq\frac{1}{2}$
for all $0<\delta\leq\min\left\{ 1,\delta_{0}\right\} $.

To this end, let $f\in\DecompSp{\CalQ}p{\ell_{w}^{q}}v$ be arbitrary
and for brevity, let 
\[
f_{j}:=\Fourier^{-1}\left(\theta_{j}\varphi_{j}\widehat{f}\right)=\left[\left(m_{\theta}\circ{\rm Ana}_{\varphi}\right)f\right]_{j}\in V_{j},
\]
as well as $f_{j}^{\left(2\right)}:=M_{-b_{j}}f_{j}$ and $\gamma_{2}^{\left(j,2\right)}:=M_{-b_{j}}\gamma_{2}^{\left(j\right)}=\left|\det T_{j}\right|\cdot\gamma_{j,2}\circ T_{j}^{T}$.
Note that since $f_{j}\in V_{j}$ is bandlimited with $\supp\widehat{f_{j}}\subset\overline{Q_{j}}\subset T_{j}\left[-R_{\CalQ},R_{\CalQ}\right]^{\dimension}+b_{j}$,
Theorem \ref{thm:BandlimitedWienerAmalgamSelfImproving} and equation
(\ref{eq:WeightedWienerAmalgamTemperedDistribution}) yield $f_{j}\in W_{T_{j}^{-T}\left[-1,1\right]^{\dimension}}\left(L_{v}^{p}\right)\hookrightarrow L_{v}^{\infty}\left(\R^{\dimension}\right)$.
Since $\gamma_{2}^{\left(j\right)}\in L_{v_{0}}^{1}\left(\R^{\dimension}\right)$,
this implies that the integral defining $\left(\gamma_{2}^{\left(j\right)}\ast f_{j}\right)\left(x\right)$
exists for every $x\in\R^{\dimension}$, cf.\@ equation (\ref{eq:WeightedLInftyConvolution}).
Hence, using our newly introduced notation, we have
\begin{align*}
 & \left|\left[E^{\left(\delta\right)}f\right]_{j}\left(x\right)\right|\\
 & =\left|\left[\gamma_{2}^{\left(j\right)}\ast\Fourier^{-1}\left(\theta_{j}\varphi_{j}\cdot\widehat{f}\right)\right]\left(x\right)-\left|\det T_{j}\right|^{-\frac{1}{2}}\sum_{k\in\Z^{\dimension}}\delta^{\dimension}\left|\det T_{j}\right|^{-\frac{1}{2}}\cdot\left[\Fourier^{-1}\!\left(\theta_{j}\varphi_{j}\cdot\widehat{f}\right)\right]\!\left(\delta\cdot T_{j}^{-T}k\right)\cdot\left(L_{\delta\cdot T_{j}^{-T}k}\gamma_{2}^{\left(j\right)}\right)\left(x\right)\right|\\
 & =\left|\sum_{k\in\Z^{\dimension}}\left[\int_{\delta T_{j}^{-T}\!\left(k+\left[0,1\right)^{\dimension}\right)}\gamma_{2}^{\left(j\right)}\left(x-y\right)\cdot f_{j}\left(y\right)\d y-\delta^{\dimension}\left|\det T_{j}^{-T}\right|\cdot f_{j}\left(\delta\cdot T_{j}^{-T}k\right)\cdot\gamma_{2}^{\left(j\right)}\left(x-\delta\cdot T_{j}^{-T}k\right)\right]\right|\\
 & \leq\sum_{k\in\Z^{\dimension}}\int_{\delta T_{j}^{-T}\!\left(k+\left[0,1\right)^{\dimension}\right)}\left|\gamma_{2}^{\left(j\right)}\left(x-y\right)\cdot f_{j}\left(y\right)-f_{j}\left(\delta\cdot T_{j}^{-T}k\right)\cdot\gamma_{2}^{\left(j\right)}\left(x-\delta\cdot T_{j}^{-T}k\right)\right|\d y\\
 & \overset{\left(\ast\right)}{=}\sum_{k\in\Z^{\dimension}}\int_{\delta T_{j}^{-T}\!\left(k+\left[0,1\right)^{\dimension}\right)}\left|\gamma_{2}^{\left(j,2\right)}\left(x-y\right)\cdot f_{j}^{\left(2\right)}\left(y\right)-\gamma_{2}^{\left(j,2\right)}\left(x-\delta\cdot T_{j}^{-T}k\right)f_{j}^{\left(2\right)}\left(\delta\cdot T_{j}^{-T}k\right)\right|\d y\\
 & \leq\!\sum_{k\in\Z^{\dimension}}\int_{\delta T_{j}^{-T}\!\left(k+\left[0,1\right)^{\dimension}\right)}\!\left|\gamma_{2}^{\left(j,2\right)}\!\left(x\!-\!y\right)\left[f_{j}^{\left(2\right)}\!\left(y\right)-f_{j}^{\left(2\right)}\!\left(\delta T_{j}^{-T}\!k\right)\right]\right|+\left|f_{j}^{\left(2\right)}\!\left(\delta T_{j}^{-T}\!k\right)\left[\gamma_{2}^{\left(j,2\right)}\!\left(x\!-\!y\right)-\gamma_{2}^{\left(j,2\right)}\!\left(x-\delta T_{j}^{-T}\!k\right)\right]\right|\d y.
\end{align*}
In this calculation, we used at $\left(\ast\right)$ the easily verifiable
identity $\left(M_{b}f\right)\left(x-y\right)\cdot\left(M_{b}g\right)\left(y\right)=e^{2\pi i\left\langle b,x\right\rangle }\cdot f\left(x-y\right)g\left(y\right)$.

Next, note for arbitrary $y\in\delta T_{j}^{-T}\left(k+\left[0,1\right)^{\dimension}\right)$
that $y=\delta T_{j}^{-T}k+\delta T_{j}^{-T}u$ for some $u\in\left[-1,1\right]^{\dimension}$.
This implies $\delta T_{j}^{-T}k=y-\delta T_{j}^{-T}u\in y+\delta T_{j}^{-T}\left[-1,1\right]^{\dimension}$
and hence
\[
\left|f_{j}^{\left(2\right)}\left(y\right)-f_{j}^{\left(2\right)}\left(\delta\cdot T_{j}^{-T}k\right)\right|\leq\left(\osc{\delta T_{j}^{-T}\left[-1,1\right]^{\dimension}}f_{j}^{\left(2\right)}\right)\left(y\right).
\]
Likewise, we have $x-\delta T_{j}^{-T}k=x-\left(y-\delta T_{j}^{-T}u\right)\in x-y+\delta T_{j}^{-T}\left[-1,1\right]^{\dimension}$,
which yields
\[
\left|\gamma_{2}^{\left(j,2\right)}\left(x-y\right)-\gamma_{2}^{\left(j,2\right)}\left(x-\delta\cdot T_{j}^{-T}k\right)\right|\leq\left(\osc{\delta T_{j}^{-T}\left[-1,1\right]^{\dimension}}\gamma_{2}^{\left(j,2\right)}\right)\left(x-y\right).
\]
Finally, we also have
\begin{align*}
\left|f_{j}^{\left(2\right)}\left(\delta\cdot T_{j}^{-T}k\right)\right| & \leq\left|f_{j}^{\left(2\right)}\left(\delta\cdot T_{j}^{-T}k\right)-f_{j}^{\left(2\right)}\left(y\right)\right|+\left|f_{j}^{\left(2\right)}\left(y\right)\right|\\
 & \leq\left|f_{j}^{\left(2\right)}\left(y\right)\right|+\left(\osc{\delta T_{j}^{-T}\left[-1,1\right]^{\dimension}}f_{j}^{\left(2\right)}\right)\left(y\right)\\
 & =:e_{j}\left(y\right),
\end{align*}
so that we see
\begin{align}
 & \left|\left[E^{\left(\delta\right)}f\right]_{j}\left(x\right)\right|\nonumber \\
 & \leq\sum_{k\in\Z^{\dimension}}\int_{\delta T_{j}^{-T}\left(k+\left[0,1\right)^{\dimension}\right)}\left|\gamma_{2}^{\left(j,2\right)}\left(x-y\right)\right|\cdot\left[\osc{\delta T_{j}^{-T}\left[-1,1\right]^{\dimension}}f_{j}^{\left(2\right)}\right]\left(y\right)+e_{j}\left(y\right)\cdot\left[\osc{\delta T_{j}^{-T}\left[-1,1\right]^{\dimension}}\gamma_{2}^{\left(j,2\right)}\right]\left(x-y\right)\d y\nonumber \\
 & =\left(\left|\gamma_{2}^{\left(j,2\right)}\right|\ast\left[\osc{\delta T_{j}^{-T}\left[-1,1\right]^{\dimension}}f_{j}^{\left(2\right)}\right]\right)\left(x\right)+\left(e_{j}\ast\left[\osc{\delta T_{j}^{-T}\left[-1,1\right]^{\dimension}}\gamma_{2}^{\left(j,2\right)}\right]\right)\left(x\right).\label{eq:AtomicDecompositionEDeltaPointwiseEstimate}
\end{align}

We now distinguish two cases: For $p\in\left[1,\infty\right]$, first
note $\left|\smash{\gamma_{2}^{\left(j,2\right)}}\right|=\left|\smash{\gamma_{2}^{\left(j\right)}}\right|$
and hence, thanks to equation (\ref{eq:AtomicDecompositionGamma2WeightedL1Norm}),
$\left\Vert \smash{\gamma_{2}^{\left(j,2\right)}}\right\Vert _{L_{v_{0}}^{1}}=\left\Vert \smash{\gamma_{2}^{\left(j\right)}}\right\Vert _{L_{v_{0}}^{1}}\leq\Omega_{0}^{K}\Omega_{1}\Omega_{4}^{\left(p,K,1\right)}\cdot s_{\dimension}=:C_{6}$.
Hence, we get using the triangle inequality, the weighted Young inequality
(equation (\ref{eq:WeightedYoungInequality})), the definition of
$e_{j}$ and since $K_{0}=K+\dimension+1$ that
\begin{align*}
 & \left\Vert \left[E^{\left(\delta\right)}f\right]_{j}\right\Vert _{V_{j}}=\left\Vert \left[E^{\left(\delta\right)}f\right]_{j}\right\Vert _{L_{v}^{p}}\\
\left({\scriptstyle \text{def. of }e_{j}}\right) & \leq\left\Vert \gamma_{2}^{\left(j,2\right)}\right\Vert _{L_{v_{0}}^{1}}\!\cdot\left\Vert \osc{\delta T_{j}^{-T}\left[-1,1\right]^{\dimension}}f_{j}^{\left(2\right)}\right\Vert _{L_{v}^{p}}+\left\Vert \osc{\delta T_{j}^{-T}\left[-1,1\right]^{\dimension}}\gamma_{2}^{\left(j,2\right)}\right\Vert _{L_{v_{0}}^{1}}\!\cdot\left(\left\Vert f_{j}^{\left(2\right)}\right\Vert _{L_{v}^{p}}\!+\left\Vert \osc{\delta T_{j}^{-T}\left[-1,1\right]^{\dimension}}f_{j}^{\left(2\right)}\right\Vert _{L_{v}^{p}}\right)\\
 & \leq C_{6}\!\left\Vert \osc{\delta T_{j}^{-T}\left[-1,1\right]^{\dimension}}\!\left[M_{-b_{j}}f_{j}\right]\right\Vert _{L_{v}^{p}}\!\!\!+\!\left\Vert \osc{\delta T_{j}^{-T}\left[-1,1\right]^{\dimension}}\!\left[M_{-b_{j}}\gamma_{2}^{\left(j\right)}\right]\right\Vert _{L_{v_{0}}^{1}}\!\!\left(\!\left\Vert f_{j}\right\Vert _{L_{v}^{p}}\!+\!\left\Vert \osc{\delta T_{j}^{-T}\left[-1,1\right]^{\dimension}}\!\left[M_{-b_{j}}f_{j}\right]\right\Vert _{L_{v}^{p}}\right)\\
\left({\scriptstyle \text{Thm. }\ref{thm:BandlimitedOscillationSelfImproving}}\right) & \leq C_{6}C_{7}\cdot\delta\cdot\left\Vert f_{j}\right\Vert _{L_{v}^{p}}+\left\Vert \osc{\delta T_{j}^{-T}\left[-1,1\right]^{\dimension}}\left[M_{-b_{j}}\gamma_{2}^{\left(j\right)}\right]\right\Vert _{L_{v_{0}}^{1}}\cdot\left(\left\Vert f_{j}\right\Vert _{L_{v}^{p}}+C_{7}\cdot\delta\cdot\left\Vert f_{j}\right\Vert _{L_{v}^{p}}\right)\\
\left({\scriptstyle \text{eq. }\eqref{eq:AtomicDecompositionFamilyDefinition}}\right) & =C_{6}C_{7}\cdot\delta\cdot\left\Vert f_{j}\right\Vert _{L_{v}^{p}}+\left|\det T_{j}\right|\cdot\left\Vert \osc{\delta T_{j}^{-T}\left[-1,1\right]^{\dimension}}\left[\gamma_{j,2}\circ T_{j}^{T}\right]\right\Vert _{L_{v_{0}}^{1}}\cdot\left(\left\Vert f_{j}\right\Vert _{L_{v}^{p}}+C_{7}\cdot\delta\cdot\left\Vert f_{j}\right\Vert _{L_{v}^{p}}\right)\\
\left({\scriptstyle \text{Lem. }\ref{lem:OscillationLinearChange}}\right) & =C_{6}C_{7}\cdot\delta\cdot\left\Vert f_{j}\right\Vert _{L_{v}^{p}}+\left|\det T_{j}\right|\cdot\left\Vert v_{0}\cdot\left(\left[\osc{\delta\left[-1,1\right]^{\dimension}}\gamma_{j,2}\right]\circ T_{j}^{T}\right)\right\Vert _{L^{1}}\cdot\left(\left\Vert f_{j}\right\Vert _{L_{v}^{p}}+C_{7}\cdot\delta\cdot\left\Vert f_{j}\right\Vert _{L_{v}^{p}}\right)\\
 & \overset{\left(\dagger\right)}{\leq}C_{6}C_{7}\cdot\delta\cdot\left\Vert f_{j}\right\Vert _{L_{v}^{p}}+\Omega_{0}^{K}\Omega_{1}\cdot\left\Vert \left(1+\left|\mybullet\right|\right)^{K}\cdot\osc{\delta\left[-1,1\right]^{\dimension}}\gamma_{j,2}\right\Vert _{L^{1}}\cdot\left(\left\Vert f_{j}\right\Vert _{L_{v}^{p}}+C_{7}\cdot\delta\cdot\left\Vert f_{j}\right\Vert _{L_{v}^{p}}\right)\\
\left({\scriptstyle \text{Lem. }\ref{lem:OscillationSchwartzFunction}}\right) & \leq C_{6}C_{7}\cdot\delta\cdot\left\Vert f_{j}\right\Vert _{L_{v}^{p}}+\Omega_{0}^{K}\Omega_{1}\cdot\left(3\sqrt{\dimension}\right)^{K_{0}+1}\cdot\delta\cdot\left\Vert \nabla\gamma_{j,2}\right\Vert _{K_{0}}\cdot\left\Vert \left(1+\left|\mybullet\right|\right)^{K-K_{0}}\right\Vert _{L^{1}}\cdot\left(\left\Vert f_{j}\right\Vert _{L_{v}^{p}}+C_{7}\cdot\delta\cdot\left\Vert f_{j}\right\Vert _{L_{v}^{p}}\right)\\
\left({\scriptstyle \text{eq. }\eqref{eq:StandardDecayLpEstimate}}\right) & \leq C_{6}C_{7}\cdot\delta\cdot\left\Vert f_{j}\right\Vert _{L_{v}^{p}}+\Omega_{0}^{K}\Omega_{1}\cdot s_{\dimension}\left(3\sqrt{\dimension}\right)^{K_{0}+1}\cdot\delta\cdot\left\Vert \nabla\gamma_{j,2}\right\Vert _{K_{0}}\cdot\left(\left\Vert f_{j}\right\Vert _{L_{v}^{p}}+C_{7}\cdot\delta\cdot\left\Vert f_{j}\right\Vert _{L_{v}^{p}}\right)\\
\left({\scriptstyle \text{since }\delta\leq1}\right) & \leq\delta\cdot\left\Vert f_{j}\right\Vert _{L_{v}^{p}}\cdot\left(C_{6}C_{7}+\Omega_{0}^{K}\Omega_{1}\cdot s_{\dimension}\left(3\sqrt{\dimension}\right)^{K_{0}+1}\cdot\Omega_{4}^{\left(p,K,2\right)}\cdot\left(1+C_{7}\right)\right)\\
 & =:C_{8}\cdot\delta\cdot\left\Vert f_{j}\right\Vert _{L_{v}^{p}}=C_{8}\cdot\delta\cdot\left\Vert f_{j}\right\Vert _{V_{j}}.
\end{align*}
Here, we defined $\Omega_{4}^{\left(p,K,2\right)}:=\sup_{j\in I}\left\Vert \nabla\gamma_{j,2}\right\Vert _{K_{0}}$,
which is finite thanks to equation (\ref{eq:AtomicDecompositionGamma2ConstantDefinition}).
The step marked with $\left(\dagger\right)$ used a simple change
of variables and our assumption $v_{0}\left(x\right)\leq\Omega_{1}\cdot\left(1+\left|x\right|\right)^{K}$
in combination with estimate (\ref{eq:WeightLinearTransformationsConnection}).
Furthermore, our application of Theorem \ref{thm:BandlimitedOscillationSelfImproving}
is justified, since we have $f_{j}=\Fourier^{-1}\left(\theta_{j}\varphi_{j}\widehat{f}\right)$,
which implies $\supp\widehat{f_{j}}\subset\supp\varphi_{j}\subset\overline{Q_{j}}\subset T_{j}\left[-R_{\CalQ},R_{\CalQ}\right]^{\dimension}+b_{j}$,
so that Theorem \ref{thm:BandlimitedOscillationSelfImproving} yields
\begin{equation}
\left\Vert \osc{\delta T_{j}^{-T}\left[-1,1\right]^{\dimension}}\left[M_{-b_{j}}f_{j}\right]\right\Vert _{V_{j}}\leq C_{7}\cdot\delta\cdot\left\Vert f_{j}\right\Vert _{L_{v}^{p}}\label{eq:AtomicDecompositionOscillationEstimate}
\end{equation}
for
\[
C_{7}:=\Omega_{0}^{2K}\Omega_{1}^{2}\cdot\left(\!23040\!\cdot\!\dimension^{\frac{3}{2}}\!\cdot\!\left(\!K\!+\!1\!+\!\frac{\dimension+1}{\min\left\{ 1,p\right\} }\right)\right)^{\!\!K+2+\frac{\dimension+1}{\min\left\{ 1,p\right\} }}\!\!\!\cdot\left(1\!+\!R_{\CalQ}\right)^{1+\frac{\dimension}{\min\left\{ 1,p\right\} }}\!=\!\left[\!2^{\max\left\{ 1,\frac{1}{p}\right\} }\Omega_{0}^{K}\Omega_{1}\left(1\!+\!\sqrt{\dimension}\right)^{\!K}\right]^{-1}\!\cdot C_{5}.
\]

In case of $p\in\left(0,1\right)$, we let $C_{9}:=2^{\frac{1}{p}-1}$,
so that $C_{9}$ is a triangle constant for $L^{p}\left(\R^{\dimension}\right)$.
Furthermore, we set $V_{j}^{\natural}:=W_{T_{j}^{-T}\left[-1,1\right]^{\dimension}}\left(L_{v_{0}}^{p}\right)$
for brevity. Then, we use Corollary \ref{cor:WienerAmalgamConvolutionSimplified}
to get for $C_{10}:=C_{9}\cdot\Omega_{0}^{3K}\Omega_{1}^{3}\cdot\dimension^{-\frac{\dimension}{2p}}\cdot\left(972\cdot\dimension^{5/2}\right)^{K+\frac{\dimension}{p}}$
that
\begin{align*}
 & \left\Vert \left[E^{\left(\delta\right)}f\right]_{j}\right\Vert _{V_{j}}\\
\left({\scriptstyle \text{eq. }\eqref{eq:AtomicDecompositionEDeltaPointwiseEstimate}}\right) & \leq C_{9}\cdot\left[\left\Vert \left|\gamma_{2}^{\left(j,2\right)}\right|\ast\osc{\delta T_{j}^{-T}\left[-1,1\right]^{\dimension}}f_{j}^{\left(2\right)}\right\Vert _{V_{j}}+\left\Vert e_{j}\ast\left[\osc{\delta T_{j}^{-T}\left[-1,1\right]^{\dimension}}\gamma_{2}^{\left(j,2\right)}\right]\right\Vert _{V_{j}}\right]\\
\left({\scriptstyle \text{Cor. }\ref{cor:WienerAmalgamConvolutionSimplified}}\right) & \leq C_{10}\!\cdot\!\left|\det T_{j}\right|^{\frac{1}{p}-1}\left(\left\Vert \gamma_{2}^{\left(j,2\right)}\right\Vert _{V_{j}^{\natural}}\cdot\left\Vert \osc{\delta T_{j}^{-T}\left[-1,1\right]^{\dimension}}f_{j}^{\left(2\right)}\right\Vert _{V_{j}}+\left\Vert e_{j}\right\Vert _{V_{j}}\cdot\left\Vert \osc{\delta T_{j}^{-T}\left[-1,1\right]^{\dimension}}\gamma_{2}^{\left(j,2\right)}\right\Vert _{V_{j}^{\natural}}\right)\\
\left({\scriptstyle \text{eq. }\eqref{eq:AtomicDecompositionFamilyDefinition}\text{, def. of }e_{j}}\right) & \leq C_{10}\!\cdot\!\left|\det T_{j}\right|^{\frac{1}{p}}\cdot\left(\left\Vert \gamma_{j,2}\circ T_{j}^{T}\right\Vert _{V_{j}^{\natural}}\cdot\left\Vert \osc{\delta T_{j}^{-T}\left[-1,1\right]^{\dimension}}\left[M_{-b_{j}}f_{j}\right]\right\Vert _{V_{j}}\right.\\
 & \phantom{\leq C_{10}\cdot\left|\det T_{j}\right|^{\frac{1}{p}}\cdot\bigg(}\left.+C_{9}\!\left[\left\Vert f_{j}\right\Vert _{V_{j}}\!+\!\left\Vert \osc{\delta T_{j}^{-T}\left[-1,1\right]^{\dimension}}\left[M_{-b_{j}}f_{j}\right]\right\Vert _{V_{j}}\right]\!\cdot\left\Vert \osc{\delta T_{j}^{-T}\left[-1,1\right]^{\dimension}}\left[\gamma_{j,2}\!\circ\!T_{j}^{T}\right]\right\Vert _{V_{j}^{\natural}}\right)\\
\left({\scriptstyle \text{Lem. }\ref{lem:OscillationLinearChange}\text{, eq. }\eqref{eq:AtomicDecompositionOscillationEstimate}}\right) & \leq C_{10}\!\cdot\!\left|\det T_{j}\right|^{\frac{1}{p}}\cdot\left(\left\Vert v_{0}\cdot M_{T_{j}^{-T}\left[-1,1\right]^{\dimension}}\left[\gamma_{j,2}\circ T_{j}^{T}\right]\right\Vert _{L^{p}}\cdot C_{7}\cdot\delta\cdot\left\Vert f_{j}\right\Vert _{V_{j}}\right.\\
 & \phantom{\leq C_{10}\!\cdot\!\left|\det T_{j}\right|^{\frac{1}{p}}\cdot\bigg(}\left.+C_{9}\!\left[\left\Vert f_{j}\right\Vert _{V_{j}}\!+\!C_{7}\cdot\delta\cdot\left\Vert f_{j}\right\Vert _{V_{j}}\right]\cdot\left\Vert v_{0}\!\cdot\!M_{T_{j}^{-T}\left[-1,1\right]^{\dimension}}\left[\!\left(\osc{\delta\left[-1,1\right]^{\dimension}}\gamma_{j,2}\right)\!\circ\!T_{j}^{T}\right]\right\Vert _{L^{p}}\right)\\
\left({\scriptstyle \text{Lem. }\ref{lem:WienerTransformationFormula}}\right) & =C_{10}\cdot\left(\left\Vert \left(v_{0}\circ T_{j}^{-T}\right)\cdot M_{\left[-1,1\right]^{\dimension}}\gamma_{j,2}\right\Vert _{L^{p}}\cdot C_{7}\cdot\delta\cdot\left\Vert f_{j}\right\Vert _{V_{j}}\right.\\
 & \phantom{=C_{10}\cdot\bigg(}\left.+C_{9}\left[\left\Vert f_{j}\right\Vert _{V_{j}}+C_{7}\cdot\delta\cdot\left\Vert f_{j}\right\Vert _{V_{j}}\right]\cdot\left\Vert \left(v_{0}\circ T_{j}^{-T}\right)\cdot M_{\left[-1,1\right]^{\dimension}}\left[\osc{\delta\left[-1,1\right]^{\dimension}}\gamma_{j,2}\right]\right\Vert _{L^{p}}\right)\\
\left({\scriptstyle \text{since }\delta\leq1}\right) & \leq C_{10}\Omega_{0}^{K}\Omega_{1}\cdot\left(\left\Vert \left(1+\left|\mybullet\right|\right)^{K}\cdot M_{\left[-1,1\right]^{\dimension}}\gamma_{j,2}\right\Vert _{L^{p}}\cdot C_{7}\cdot\delta\cdot\left\Vert f_{j}\right\Vert _{V_{j}}\right.\\
 & \phantom{=C_{10}\Omega_{0}^{K}\Omega_{1}\cdot\bigg(}\left.+C_{9}\left\Vert f_{j}\right\Vert _{V_{j}}\left(1+C_{7}\right)\cdot\left\Vert \left(1+\left|\mybullet\right|\right)^{K}\cdot M_{\left[-1,1\right]^{\dimension}}\left[\osc{\delta\left[-1,1\right]^{\dimension}}\gamma_{j,2}\right]\right\Vert _{L^{p}}\right).
\end{align*}
Here, the last step used as usual our assumption $v_{0}\left(x\right)\leq\Omega_{1}\cdot\left(1+\left|x\right|\right)^{K}$,
in combination with equation (\ref{eq:WeightLinearTransformationsConnection}).
We now combine Lemmas \ref{lem:OscillationSchwartzFunction} and
\ref{lem:SchwartzFunctionsAreWiener}  to obtain
\begin{align*}
\left\Vert \left(1+\left|\mybullet\right|\right)^{K}\cdot M_{\left[-1,1\right]^{\dimension}}\left[\osc{\delta\left[-1,1\right]^{\dimension}}\gamma_{j,2}\right]\right\Vert _{L^{p}} & \leq\left(3\sqrt{\dimension}\right)^{K_{0}+1}\cdot\delta\cdot\left\Vert \nabla\gamma_{j,2}\right\Vert _{K_{0}}\cdot\left\Vert \left(1+\left|\mybullet\right|\right)^{K}\cdot M_{\left[-1,1\right]^{\dimension}}\left(1+\left|\mybullet\right|\right)^{-K_{0}}\right\Vert _{L^{p}}\\
 & \leq\left(1+2\sqrt{\dimension}\right)^{K_{0}}\left(3\sqrt{\dimension}\right)^{K_{0}+1}\cdot\delta\cdot\left\Vert \nabla\gamma_{j,2}\right\Vert _{K_{0}}\cdot\left\Vert \left(1+\left|\mybullet\right|\right)^{K-K_{0}}\right\Vert _{L^{p}}\\
\left({\scriptstyle \text{eq. }\eqref{eq:StandardDecayLpEstimate}}\right) & \leq\left(3\sqrt{\dimension}\right)^{2K_{0}+1}\Omega_{4}^{\left(p,K,2\right)}\cdot\delta\cdot\left(\frac{s_{\dimension}}{p}\right)^{1/p}.
\end{align*}
Likewise, Lemma \ref{lem:SchwartzFunctionsAreWiener} also yields
\begin{align*}
\left\Vert \left(1+\left|\mybullet\right|\right)^{K}\cdot M_{\left[-1,1\right]^{\dimension}}\gamma_{j,2}\right\Vert _{L^{p}} & \leq\left\Vert \gamma_{j,2}\right\Vert _{K_{0}}\cdot\left\Vert \left(1+\left|\mybullet\right|\right)^{K}\cdot M_{\left[-1,1\right]^{\dimension}}\left(1+\left|\mybullet\right|\right)^{-K_{0}}\right\Vert _{L^{p}}\\
 & \leq\left\Vert \gamma_{j,2}\right\Vert _{K_{0}}\cdot\left(1+2\sqrt{\dimension}\right)^{K_{0}}\cdot\left\Vert \left(1+\left|\mybullet\right|\right)^{K-K_{0}}\right\Vert _{L^{p}}\\
 & \leq\left(3\sqrt{\dimension}\right)^{K_{0}}\Omega_{4}^{\left(p,K,1\right)}\cdot\left(\frac{s_{\dimension}}{p}\right)^{1/p}.
\end{align*}
Combining these estimates with our estimate for $\left\Vert \left[E^{\left(\delta\right)}f\right]_{j}\right\Vert _{V_{j}}$,
we arrive at
\begin{align*}
\left\Vert \left[E^{\left(\delta\right)}f\right]_{j}\right\Vert _{V_{j}} & \leq C_{10}\left(\frac{s_{\dimension}}{p}\right)^{1/p}\left(3\sqrt{\dimension}\right)^{2K_{0}+1}\cdot\Omega_{0}^{K}\Omega_{1}\cdot\left(C_{7}\Omega_{4}^{\left(p,K,1\right)}+C_{9}\left(1+C_{7}\right)\Omega_{4}^{\left(p,K,2\right)}\right)\cdot\delta\cdot\left\Vert f_{j}\right\Vert _{V_{j}}\\
 & =:C_{11}\cdot\delta\cdot\left\Vert f_{j}\right\Vert _{V_{j}},
\end{align*}
where $C_{11}$ is independent of $\delta$ and $j$.

\medskip{}

All in all, if we set $C_{12}:=C_{8}$ for $p\in\left[1,\infty\right]$
and $C_{12}:=C_{11}$ for $p\in\left(0,1\right)$, we have $\left\Vert \left[E^{\left(\delta\right)}f\right]_{j}\right\Vert _{V_{j}}\leq C_{12}\cdot\delta\cdot\left\Vert f_{j}\right\Vert _{V_{j}}$
for all $j\in I$ and $\delta\in\left(0,1\right]$. But this entails
\begin{align*}
\left\Vert E^{\left(\delta\right)}f\right\Vert _{\ell_{w}^{q}\left(\left[V_{j}\right]_{j\in I}\right)} & \leq C_{12}\cdot\delta\cdot\left\Vert \left(f_{j}\right)_{j\in I}\right\Vert _{\ell_{w}^{q}\left(\left[V_{j}\right]_{j\in I}\right)}\\
\left({\scriptstyle \text{since }f_{j}=\left[\left(m_{\theta}\circ{\rm Ana}_{\varphi}\right)f\right]_{j}}\right) & =C_{12}\cdot\delta\cdot\left\Vert \left(m_{\theta}\circ{\rm Ana}_{\varphi}\right)f\right\Vert _{\ell_{w}^{q}\left(\left[V_{j}\right]_{j\in I}\right)}\\
 & \leq C_{12}\cdot\vertiii{m_{\theta}}\cdot\vertiii{{\rm Ana}_{\varphi}}\cdot\delta\cdot\left\Vert f\right\Vert _{\DecompSp{\CalQ}p{\ell_{w}^{q}}v}\\
 & \leq C_{1}C_{2}C_{12}\cdot\delta\cdot\left\Vert f\right\Vert _{\DecompSp{\CalQ}p{\ell_{w}^{q}}v}.
\end{align*}
But since $\vertiii{{\rm Synth}_{\Gamma_{1}}}\leq C_{3}\cdot\vertiii{\smash{\overrightarrow{C}}}^{\max\left\{ 1,\frac{1}{p}\right\} }$,
this means in view of equation (\ref{eq:AtomicDecompositionReproducingFormula})
that
\begin{align*}
\vertiii{\identity_{\DecompSp{\CalQ}p{\ell_{w}^{q}}v}-{\rm Synth}_{\Gamma_{1}}\circ\bigotimes_{j\in I}S_{\Gamma_{2}}^{\left(\delta,j\right)}\circ D^{\left(\delta\right)}} & =\vertiii{{\rm Synth}_{\Gamma_{1}}\circ m_{\Gamma_{2}}\circ m_{\theta}\circ{\rm Ana}_{\varphi}-{\rm Synth}_{\Gamma_{1}}\circ\bigotimes_{j\in I}S_{\Gamma_{2}}^{\left(\delta,j\right)}\circ D^{\left(\delta\right)}}\\
 & \leq\vertiii{{\rm Synth}_{\Gamma_{1}}}\cdot\vertiii{\smash{E^{\left(\delta\right)}}}\\
 & \leq C_{1}C_{2}C_{3}C_{12}\cdot\vertiii{\smash{\overrightarrow{C}}}^{\max\left\{ 1,\frac{1}{p}\right\} }\cdot\delta.
\end{align*}

\medskip{}

Now, we estimate the constant $C_{1}C_{2}C_{3}C_{12}$ to see that
$\delta\leq\delta_{0}$ implies $C_{1}C_{2}C_{3}C_{12}\cdot\vertiii{\smash{\overrightarrow{C}}}^{\max\left\{ 1,\frac{1}{p}\right\} }\cdot\delta\leq\frac{1}{2}$.
First, in case of $p\in\left[1,\infty\right]$, we have because of
$C_{7}\geq1$ and $\max\left\{ 1,\frac{1}{p}\right\} =1$, as well
as $K_{0}=K+\dimension+1$ that
\begin{align*}
C_{1}C_{2}C_{3}C_{12} & =\Omega_{0}^{K}\Omega_{1}\Omega_{2}^{\left(p,K\right)}\cdot C_{8}\\
 & =\left(C_{6}C_{7}+\Omega_{0}^{K}\Omega_{1}\cdot s_{\dimension}\left(3\sqrt{\dimension}\right)^{K_{0}+1}\cdot\Omega_{4}^{\left(p,K,2\right)}\cdot\left(1+C_{7}\right)\right)\cdot\Omega_{0}^{K}\Omega_{1}\Omega_{2}^{\left(p,K\right)}\\
\left({\scriptstyle \text{since }\Omega_{0},\Omega_{1}\geq1}\right) & \leq C_{7}s_{\dimension}\left(\Omega_{4}^{\left(p,K,1\right)}+2\cdot\left(3\sqrt{\dimension}\right)^{K_{0}+1}\Omega_{4}^{\left(p,K,2\right)}\right)\cdot\Omega_{0}^{2K}\Omega_{1}^{2}\Omega_{2}^{\left(p,K\right)}\\
 & \leq\left(2^{17}\!\cdot\!\dimension^{2}\!\cdot\!\left(K\!+\!2\!+\!\dimension\right)\right)^{\!K+\dimension+3}\!\!\!\cdot2s_{\dimension}\cdot\left(3\!\cdot\!\dimension^{\frac{1}{2}}\right)^{-1}\!\!\cdot\!\left(1\!+\!R_{\CalQ}\right)^{\dimension+1}\!\cdot\!\left(\Omega_{4}^{\left(p,K,1\right)}\!+\!\Omega_{4}^{\left(p,K,2\right)}\right)\!\cdot\!\Omega_{0}^{4K}\Omega_{1}^{4}\Omega_{2}^{\left(p,K\right)}\\
 & \leq\frac{s_{\dimension}}{\sqrt{\dimension}}\cdot\left(2^{17}\!\cdot\!\dimension^{2}\!\cdot\!\left(K\!+\!2\!+\!\dimension\right)\right)^{\!K+\dimension+3}\!\!\!\cdot\left(1\!+\!R_{\CalQ}\right)^{\dimension+1}\cdot\Omega_{0}^{4K}\Omega_{1}^{4}\Omega_{2}^{\left(p,K\right)}\Omega_{4}^{\left(p,K\right)}.
\end{align*}
Next, for $p\in\left(0,1\right)$, we get because of $\max\left\{ 1,\frac{1}{p}\right\} =\frac{1}{p}$
and $s_{\dimension}\leq2^{2\dimension}$ that
\begin{align*}
C_{1}C_{3} & =\frac{\left(2^{10}/\dimension^{\frac{1}{2}}\right)^{\frac{\dimension}{p}}}{2^{21}\cdot\dimension^{7}}\cdot2^{4}s_{\dimension}^{\frac{1}{p}}\left(192\!\cdot\dimension^{\frac{3}{2}}\!\cdot\!\left\lceil K\!+\!\frac{\dimension+1}{p}\right\rceil \right)^{\!\left\lceil K+\frac{\dimension+1}{p}\right\rceil +1}\!\cdot\!\left(2^{21}\!\cdot\!\dimension^{5}\!\cdot\!\left\lceil K\!+\!\frac{\dimension+1}{p}\right\rceil \right)^{\!\left\lceil K+\frac{\dimension+1}{p}\right\rceil +1}\!\!\!\cdot\!\left(1\!+\!R_{\CalQ}\right)^{\frac{2\dimension}{p}}\!\cdot\!\Omega_{0}^{6K}\Omega_{1}^{6}\\
 & \leq2^{4}\cdot\frac{\left(2^{12}/\sqrt{\dimension}\right)^{\frac{\dimension}{p}}}{2^{21}\cdot\dimension^{7}}\cdot\left(2^{29}\cdot\dimension^{\frac{13}{2}}\cdot\left\lceil K+\frac{\dimension+1}{p}\right\rceil ^{2}\right)^{\left\lceil K+\frac{\dimension+1}{p}\right\rceil +1}\!\!\!\cdot\!\left(1\!+\!R_{\CalQ}\right)^{\frac{2\dimension}{p}}\!\cdot\!\Omega_{0}^{6K}\Omega_{1}^{6}
\end{align*}
and thus, since $\left\lceil K+\frac{\dimension+1}{p}\right\rceil +1\geq K+\frac{\dimension}{p}+2$,
\begin{align*}
C_{1}C_{2}C_{3} & \leq\dimension^{-\frac{\dimension}{2p}}\cdot\left(972\cdot\dimension^{\frac{5}{2}}\right)^{-2}\!\cdot\!2^{4}\!\cdot\!\frac{\left(2^{12}/\sqrt{\dimension}\right)^{\frac{\dimension}{p}}}{2^{21}\cdot\dimension^{7}}\!\cdot\!\left(2^{39}\!\cdot\!\dimension^{9}\!\cdot\!\left\lceil K\!+\!\frac{\dimension+1}{p}\right\rceil ^{2}\right)^{\!\left\lceil K+\frac{\dimension+1}{p}\right\rceil +1}\!\!\!\cdot\!\left(1\!+\!R_{\CalQ}\right)^{\frac{2\dimension}{p}}\!\cdot\!\Omega_{0}^{10K}\Omega_{1}^{10}\Omega_{2}^{\left(p,K\right)}\\
 & \leq\frac{\left(2^{12}/\dimension\right)^{\frac{\dimension}{p}}}{2^{36}\cdot\dimension^{12}}\cdot\left(2^{39}\cdot\dimension^{9}\cdot\left\lceil K+\frac{\dimension+1}{p}\right\rceil ^{2}\right)^{\left\lceil K+\frac{\dimension+1}{p}\right\rceil +1}\!\!\!\cdot\!\left(1\!+\!R_{\CalQ}\right)^{\frac{2\dimension}{p}}\!\cdot\!\Omega_{0}^{10K}\Omega_{1}^{10}\Omega_{2}^{\left(p,K\right)}.
\end{align*}
Now, recall that $C_{7}\geq1$ and $K_{0}=K+\frac{\dimension}{p}+1$,
so that
\begin{align*}
C_{12} & =C_{11}=C_{10}\left(\frac{s_{\dimension}}{p}\right)^{1/p}\left(3\sqrt{\dimension}\right)^{2K_{0}+1}\cdot\Omega_{0}^{K}\Omega_{1}\cdot\left(C_{7}\Omega_{4}^{\left(p,K,1\right)}+C_{9}\Omega_{4}^{\left(p,K,2\right)}\left(1+C_{7}\right)\right)\\
 & \leq C_{7}\cdot2^{\frac{1}{p}-1}\cdot\dimension^{-\frac{\dimension}{2p}}\cdot\left(972\cdot\dimension^{5/2}\right)^{K+\frac{\dimension}{p}}\cdot\left(\frac{s_{\dimension}}{p}\right)^{1/p}\left(9\dimension\right)^{K_{0}+1}\cdot\Omega_{0}^{4K}\Omega_{1}^{4}\cdot\left(\Omega_{4}^{\left(p,K,1\right)}+2^{\frac{1}{p}}\Omega_{4}^{\left(p,K,2\right)}\right)\\
 & \leq\frac{1}{2}C_{7}\cdot4^{\frac{1}{p}}\cdot\dimension^{-\frac{\dimension}{2p}}\cdot\left(972\cdot\dimension^{5/2}\right)^{-2}\cdot\left(8748\cdot\dimension^{7/2}\right)^{K_{0}+1}\cdot\left(\frac{s_{\dimension}}{p}\right)^{1/p}\cdot\Omega_{0}^{4K}\Omega_{1}^{4}\Omega_{4}^{\left(p,K\right)}
\end{align*}
and hence because of $K_{0}+1=K+\frac{\dimension}{p}+2\leq\left\lceil K+\frac{\dimension+1}{p}\right\rceil +1$,
\begin{align*}
 & 2C_{1}C_{2}C_{3}C_{12}\\
 & \leq C_{7}\!\cdot\!\frac{\left(2^{14}/\dimension^{\frac{3}{2}}\right)^{\frac{\dimension}{p}}}{2^{45}\cdot\dimension^{17}}\!\cdot\!\left(\frac{s_{\dimension}}{p}\right)^{\frac{1}{p}}\!\cdot\!\left(2^{53}\!\cdot\!\dimension^{\frac{25}{2}}\!\cdot\!\left\lceil K\!+\!\frac{\dimension+1}{p}\right\rceil ^{2}\right)^{\!\left\lceil K+\frac{\dimension+1}{p}\right\rceil +1}\!\!\!\cdot\!\left(1\!+\!R_{\CalQ}\right)^{\frac{2\dimension}{p}}\!\cdot\!\Omega_{0}^{14K}\Omega_{1}^{14}\Omega_{2}^{\left(p,K\right)}\Omega_{4}^{\left(p,K\right)}\\
 & \leq\frac{\left(2^{14}/\dimension^{\frac{3}{2}}\right)^{\frac{\dimension}{p}}}{2^{45}\cdot\dimension^{17}}\!\cdot\!\left(\frac{s_{\dimension}}{p}\right)^{\frac{1}{p}}\!\cdot\!\left(2^{68}\!\cdot\!\dimension^{14}\!\cdot\!\left[K\!+\!1\!+\!\frac{\dimension+1}{p}\right]^{3}\right)^{\!K+\frac{\dimension+1}{p}+2}\!\!\!\cdot\!\left(1\!+\!R_{\CalQ}\right)^{1+\frac{3\dimension}{p}}\!\cdot\!\Omega_{0}^{16K}\Omega_{1}^{16}\Omega_{2}^{\left(p,K\right)}\Omega_{4}^{\left(p,K\right)}.
\end{align*}
These considerations easily show that $\delta\leq\delta_{0}$ indeed
implies $C_{1}C_{2}C_{3}C_{12}\cdot\vertiii{\smash{\overrightarrow{C}}}^{\max\left\{ 1,\frac{1}{p}\right\} }\cdot\delta\leq\frac{1}{2}$.

\medskip{}

All in all, our considerations show for
\[
T^{\left(\delta\right)}:={\rm Synth}_{\Gamma_{1}}\circ\bigotimes_{j\in I}S_{\Gamma_{2}}^{\left(\delta,j\right)}\circ D^{\left(\delta\right)}\:\overset{\text{eq. }\eqref{eq:AtomicDecompositionSynthesisOperatorAsComposition}}{=}\:S^{\left(\delta\right)}\circ D^{\left(\delta\right)}:\DecompSp{\CalQ}p{\ell_{w}^{q}}v\to\DecompSp{\CalQ}p{\ell_{w}^{q}}v
\]
that $\vertiii{\identity_{\DecompSp{\CalQ}p{\ell_{w}^{q}}v}-T^{\left(\delta\right)}}\leq\frac{1}{2}$
for all $0<\delta\leq\min\left\{ 1,\delta_{0}\right\} $. Hence, since
$\DecompSp{\CalQ}p{\ell_{w}^{q}}v$ is a Quasi-Banach space by Lemma
\ref{lem:WeightedDecompositionSpaceComplete}, a Neumann series argument
(which is also valid for Quasi-Banach spaces, cf.\@ e.g.\@ \cite[Lemma 2.4.11]{VoigtlaenderPhDThesis}),
shows that $T^{\left(\delta\right)}$ is invertible for all $0<\delta\leq\min\left\{ 1,\delta_{0}\right\} $.

But then, $C^{\left(\delta\right)}:=D^{\left(\delta\right)}\circ\left(T^{\left(\delta\right)}\right)^{-1}:\DecompSp{\CalQ}p{\ell_{w}^{q}}v\to\smash{\ell_{\left(\left|\det T_{j}\right|^{\frac{1}{2}-\frac{1}{p}}w_{j}\right)_{\!\!j\in I}}^{q}\!\!\!}\!\!\!\left(\left[\vphantom{F}\smash{C_{j}^{\left(\delta\right)}}\right]_{j\in I}\right)$
is well-defined and bounded and we have for arbitrary $f\in\DecompSp{\CalQ}p{\ell_{w}^{q}}v$
that
\[
f=\left[T^{\left(\delta\right)}\circ\left(T^{\left(\delta\right)}\right)^{-1}\right]f=\left(\left[S^{\left(\delta\right)}\circ D^{\left(\delta\right)}\right]\circ\left(T^{\left(\delta\right)}\right)^{-1}\right)f=\left[S^{\left(\delta\right)}\circ C^{\left(\delta\right)}\right]f,
\]
as desired.
\end{proof}

\section{Simplified Criteria}

\label{sec:SimplifiedCriteria}In this section, we will derive simplified
conditions which ensure boundedness of the operators $\overrightarrow{A},\overrightarrow{B}$
and $\overrightarrow{C}$, mentioned in Assumptions \ref{assu:MainAssumptions},
\ref{assu:DiscreteBanachFrameAssumptions} and \ref{assu:AtomicDecompositionAssumption},
respectively.

One such general criterion is given by \textbf{Schur's test}, which
we state below. Afterwards, we will provide a convenient standard
estimate for the main term $\left\Vert \Fourier^{-1}\left(\varphi_{i}\cdot\widehat{\gamma^{\left(j\right)}}\right)\right\Vert _{L_{v_{0}}^{p}}$
occurring in the entries of $\overrightarrow{A},\overrightarrow{B}$
and $\overrightarrow{C}$. Then we use these results to formulate
simplified criteria which allow to apply Theorems \ref{thm:DiscreteBanachFrameTheorem}
(leading to Banach frames) and \ref{thm:AtomicDecomposition} (leading
to atomic decompositions).

But first of all, we introduce certain additional assumptions regarding
the partition of unity $\Phi=\left(\varphi_{i}\right)_{i\in I}$.
Recall that in the preceding sections, we only assumed $\Phi$ to
be a $\CalQ$-$v_{0}$-BAPU, but in this section we will make the
following stronger assumption:
\begin{assumption}
\label{assu:RegularPartitionOfUnity}We assume that $\Phi=\left(\varphi_{i}\right)_{i\in I}$
is a \textbf{regular partition of unity} for $\CalQ$. This means

\begin{enumerate}
\item $\varphi_{i}\in\TestFunctionSpace{\CalO}$ with $\supp\varphi_{i}\subset Q_{i}$
for all $i\in I$,
\item $\sum_{i\in I}\varphi_{i}\equiv1$ on $\CalO$,
\item the \textbf{normalized family} $\Phi^{\natural}:=\left(\smash{\varphi_{i}^{\natural}}\right)_{i\in I}$—given
by $\varphi_{i}^{\natural}:=\varphi_{i}\circ S_{i}$ for $S_{i}\xi:=T_{i}\xi+b_{i}$—satisfies
\begin{equation}
C^{\left(\alpha\right)}:=\sup_{i\in I}\left\Vert \partial^{\alpha}\smash{\varphi_{i}^{\natural}}\right\Vert _{\sup}<\infty\qquad\text{ for all }\alpha\in\N_{0}^{\dimension}.\qedhere\label{eq:RegularBAPUCondition}
\end{equation}
\end{enumerate}
\end{assumption}
\begin{rem*}
As seen in \cite[Lemma 2.5]{DecompositionEmbedding}, every regular
partition of unity is also a $\CalQ$-$v_{0}$-BAPU, as long as $v_{0}\lesssim1$.
As we will see in Corollary \ref{cor:RegularBAPUsAreWeightedBAPUs},
the same also remains true for general $v_{0}$.

Furthermore, it was shown in \cite[Theorem 2.8]{DecompositionIntoSobolev}
that every \emph{structured} admissible covering $\CalQ$ admits a
regular partition of unity. Here, the semi-structured covering $\CalQ=\left(T_{i}Q_{i}'+b_{i}\right)_{i\in I}$
is called \textbf{structured} if $Q_{i}'=Q$ for all $i\in I$ and
some fixed open set $Q\subset\R^{\dimension}$ and if additionally,
there is an open set $P\subset\R^{\dimension}$, compactly contained
in $Q$, such that the family $\left(T_{i}P+b_{i}\right)_{i\in I}$
covers all of $\CalO$.
\end{rem*}
Now that we have clarified our assumptions for this section, we state
a version of Schur's test which is suitable for our setting. We remark
that this lemma is in no way new; for example, it already appears
in \cite[Lemma 4.4]{ParabolicMolecules}.
\begin{lem}
\label{lem:SchursLemma}Let $I,J\neq\emptyset$ be two nonempty sets
and let $A=\left(A_{i,j}\right)_{\left(i,j\right)\in I\times J}\in\Compl^{I\times J}$.
Let $p\in\left(1,\infty\right)$ and assume that
\[
C_{1}:=\sup_{i\in I}\sum_{j\in J}\left|A_{i,j}\right|\qquad\text{ and }\qquad C_{2}:=\sup_{j\in J}\sum_{i\in I}\left|A_{i,j}\right|
\]
are finite. Then the operator
\[
\overrightarrow{A}:\ell^{p}\left(J\right)\to\ell^{p}\left(I\right),\left(c_{j}\right)_{j\in J}\mapsto\left(\sum_{j\in J}A_{i,j}c_{j}\right)_{i\in I}
\]
is well-defined and bounded with $\vertiii{\smash{\overrightarrow{A}}}\leq\max\left\{ C_{1},C_{2}\right\} $.

In case of $p\in\left(0,1\right]$, it suffices if
\[
C_{3}^{\left(p\right)}:=\sup_{j\in J}\sum_{i\in I}\left|A_{i,j}\right|^{p}
\]
is finite. In this case, $\vertiii{\smash{\overrightarrow{A}}}\leq\left(C_{3}^{\left(p\right)}\right)^{1/p}$.

Finally, in case of $p=\infty$, it suffices if
\[
C_{4}:=\sup_{i\in I}\sum_{j\in J}\left|A_{i,j}\right|
\]
is finite. In this case, $\vertiii{\smash{\overrightarrow{A}}}\leq C_{4}$.
\end{lem}
\begin{proof}
The statement for $p\in\left(1,\infty\right)$ follows from the more
general form of Schur's test as given e.g.\@ in \cite[Theorem 6.18]{FollandRA},
by considering $I$ and $J$ as measure spaces by equipping them with
the counting measure. Strictly speaking, that lemma assumes the underlying
measure spaces to be $\sigma$-finite (i.e., $I,J$ have to be countable),
but since Tonelli's theorem is applicable to uncountable sets equipped
with the counting measure, the proof given in \cite{FollandRA} still
works even for uncountable $I,J$.

Now, let us assume $p\in\left(0,1\right]$. In this case, we have
\begin{align*}
\left\Vert \overrightarrow{A}\left(c_{j}\right)_{j\in J}\right\Vert _{\ell^{p}}^{p}=\sum_{i\in I}\left|\left(\overrightarrow{A}\left(c_{j}\right)_{j\in J}\right)_{i}\right|^{p} & =\sum_{i\in I}\left|\sum_{j\in J}A_{i,j}\cdot c_{j}\right|^{p}\\
\left({\scriptstyle \text{since }\left(\sum a_{j}\right)^{p}\leq\sum a_{j}^{p}\text{ for }p\in\left(0,1\right]\text{ and }a_{j}\geq0}\right) & \leq\sum_{i\in I}\:\sum_{j\in J}\left|A_{i,j}\right|^{p}\left|c_{j}\right|^{p}\\
 & =\sum_{j\in J}\left(\left|c_{j}\right|^{p}\,\sum_{i\in I}\left|A_{i,j}\right|^{p}\right)\\
 & \leq C_{3}^{\left(p\right)}\cdot\sum_{j\in J}\left|c_{j}\right|^{p}\\
 & =C_{3}^{\left(p\right)}\cdot\left\Vert \left(c_{j}\right)_{j\in J}\right\Vert _{\ell^{p}}^{p}<\infty,
\end{align*}
so that $\overrightarrow{A}:\ell^{p}\left(J\right)\to\ell^{p}\left(I\right)$
is bounded with $\vertiii{\smash{\overrightarrow{A}}}_{\ell^{p}\to\ell^{p}}\leq\left(C_{3}^{\left(p\right)}\right)^{1/p}$.

Finally, let $p=\infty$. For arbitrary $i\in I$, we have
\begin{align*}
\left|\left(\overrightarrow{A}\left(c_{j}\right)_{j\in J}\right)_{i}\right| & =\left|\sum_{j\in J}A_{i,j}\cdot c_{j}\right|\\
 & \leq\sum_{j\in J}\left(\left|A_{i,j}\right|\cdot\left|c_{j}\right|\right)\\
 & \leq\left\Vert \left(c_{j}\right)_{j\in J}\right\Vert _{\ell^{\infty}}\cdot\sum_{j\in J}\left|A_{i,j}\right|\\
 & \leq C_{4}\cdot\left\Vert \left(c_{j}\right)_{j\in J}\right\Vert _{\ell^{p}}.
\end{align*}

As a further remark, we observe that the case $p\in\left(1,\infty\right)$
can be obtained by complex interpolation (i.e., by the Riesz-Thorin
Theorem \cite[Theorem 6.27]{FollandRA}) from the cases $p=1$ and
$p=\infty$, since $C_{1}=C_{4}$ and $C_{2}=C_{3}^{\left(1\right)}$.
\end{proof}
In Lemma \ref{lem:GramMatrixEstimates} below, we will derive a convenient
estimate for the main term of the ``infinite matrices'' $A,B,C$
from Assumptions \ref{assu:MainAssumptions}, \ref{assu:DiscreteBanachFrameAssumptions}
and \ref{assu:AtomicDecompositionAssumption}, namely for the term
$\left\Vert \Fourier^{-1}\left(\varphi_{i}\widehat{\gamma^{\left(j\right)}}\right)\right\Vert _{L_{v_{0}}^{p}}$.
To derive this estimate, the following lemma will be useful. It makes
precise the notion that \emph{smoothness of $f$ yields decay of $\widehat{f}$}.
The statement itself is probably folklore, so no originality is claimed.
\begin{lem}
\label{lem:PointwiseFourierDecayEstimate}Let $N\in\N_{0}$ and $g\in W^{N,1}\left(\R^{\dimension}\right)$.
Then
\begin{equation}
\left(1+\left|x\right|\right)^{N}\cdot\left|\Fourier^{-1}g\left(x\right)\right|\leq\left(1+\dimension\right)^{N}\cdot\left(\left|\Fourier^{-1}g\left(x\right)\right|+\sum_{m=1}^{\dimension}\left|\left[\Fourier^{-1}\left(\partial_{m}^{N}g\right)\right]\left(x\right)\right|\right)\qquad\forall x\in\R^{\dimension}.\qedhere\label{eq:PointwiseFourierDecayEstimate}
\end{equation}
\end{lem}
\begin{rem*}
Here, $W^{N,1}\left(\R^{\dimension}\right)$ is the \textbf{Sobolev
space} of all functions $g\in L^{1}\left(\R^{\dimension}\right)$
for which all weak derivatives $\partial^{\alpha}g$ with $\left|\alpha\right|\leq N$
satisfy $\partial^{\alpha}g\in L^{1}\left(\R^{\dimension}\right)$.
It is a Banach space when equipped with the norm $\left\Vert g\right\Vert _{W^{N,1}}:=\sum_{\left|\alpha\right|\leq N}\left\Vert \partial^{\alpha}g\right\Vert _{L^{1}}$.
For $N=0$, we use the convention $W^{N,1}\left(\R^{\dimension}\right)=L^{1}\left(\R^{\dimension}\right)$.
\end{rem*}
\begin{proof}
It is well-known (see e.g. \cite[Corollary 3.23]{AdamsSobolevSpaces})
that $\TestFunctionSpace{\R^{\dimension}}\subset W^{N,1}\left(\R^{\dimension}\right)$
is dense. Furthermore, since $\Fourier^{-1}:L^{1}\left(\R^{\dimension}\right)\to C_{0}\left(\R^{\dimension}\right)$
is well-defined and bounded, where the space $C_{0}\left(\R^{\dimension}\right)$
of continuous functions vanishing at infinity is equipped with the
norm $\left\Vert h\right\Vert _{\sup}:=\sup_{x\in\R^{\dimension}}\left|h\left(x\right)\right|$,
it is not hard to see that $\Fourier^{-1}g\left(x\right)$ and $\left[\Fourier^{-1}\left(\partial_{m}^{N}g\right)\right]\left(x\right)$
all depend continuously on $g\in W^{N,1}\left(\R^{\dimension}\right)$,
for arbitrary $x\in\R^{\dimension}$ and $m\in\underline{\dimension}$.
Hence, we can without loss of generality assume $g\in\TestFunctionSpace{\R^{\dimension}}$.

But under this assumption, we have (see e.g.\@ \cite[Theorem 8.22]{FollandRA})
the standard identity
\[
\left[\Fourier^{-1}\left(\partial_{m}^{N}g\right)\right]\left(x\right)=\left(-2\pi ix_{m}\right)^{N}\cdot\left(\Fourier^{-1}g\right)\left(x\right)\qquad\forall x\in\R^{\dimension}.
\]

In particular, since
\[
\left(\sum_{i=1}^{k}a_{i}\right)^{N}\leq\left(k\cdot\max\left\{ a_{i}\with i\in\underline{k}\right\} \right)^{N}\leq k^{N}\cdot\sum_{i=1}^{k}a_{i}^{N}
\]
holds for arbitrary $a_{1},\dots,a_{k}\geq0$ and because of $\left|x\right|\leq\left\Vert x\right\Vert _{1}$,
we get
\begin{align*}
\left(1+\left|x\right|\right)^{N}\cdot\left|\left(\Fourier^{-1}g\right)\left(x\right)\right| & \leq\left(1+\sum_{m=1}^{\dimension}\left|x_{m}\right|\right)^{N}\cdot\left|\left(\Fourier^{-1}g\right)\left(x\right)\right|\\
 & \leq\left(\dimension+1\right)^{N}\cdot\left|\left(\Fourier^{-1}g\right)\left(x\right)\right|\cdot\left(1+\sum_{m=1}^{\dimension}\left|x_{m}^{N}\right|\right)\\
 & =\left(\dimension+1\right)^{N}\cdot\left(\left|\left(\Fourier^{-1}g\right)\left(x\right)\right|+\sum_{m=1}^{\dimension}\left|\frac{\left[\Fourier^{-1}\left(\partial_{m}^{N}g\right)\right]\left(x\right)}{\left(2\pi\right)^{N}}\right|\right)\\
 & \leq\left(1+\dimension\right)^{N}\cdot\left(\left|\left(\Fourier^{-1}g\right)\left(x\right)\right|+\sum_{m=1}^{\dimension}\left|\left[\Fourier^{-1}\left(\partial_{m}^{N}g\right)\right]\left(x\right)\right|\right)
\end{align*}
for arbitrary $x\in\R^{\dimension}$, as desired.
\end{proof}
Now, we are finally in a position to derive the promised estimate
for $\left\Vert \Fourier^{-1}\left(\varphi_{i}\widehat{\gamma^{\left(j\right)}}\right)\right\Vert _{L_{v_{0}}^{p}}$.
\begin{lem}
\label{lem:GramMatrixEstimates}Suppose that $\left(\varphi_{i}\right)_{i\in I}$
satisfies Assumption \ref{assu:RegularPartitionOfUnity}. Let $\gamma\in L^{1}\left(\R^{\dimension}\right)$
and assume that $\widehat{\gamma}\in C^{\infty}\left(\R^{\dimension}\right)$.
For $j\in I$, define
\[
\gamma^{\left\llbracket j\right\rrbracket }:=\Fourier^{-1}\left(\widehat{\gamma}\circ S_{j}^{-1}\right)=\left|\det T_{j}\right|\cdot M_{b_{j}}\left[\gamma\circ T_{j}^{T}\right].
\]
Then we have for arbitrary $\varepsilon>0$, $i,j\in I$ and $p\in\left(0,\infty\right)$
the estimate
\begin{align*}
\left\Vert \Fourier^{-1}\left(\widehat{\gamma^{\left\llbracket j\right\rrbracket }}\cdot\varphi_{i}\right)\right\Vert _{L_{v_{0}}^{p}} & \leq C_{0}\cdot\left(1\!+\!\left\Vert T_{j}^{-1}T_{i}\right\Vert \right)^{\!\left\lceil K+\frac{\dimension+\varepsilon}{p}\right\rceil }\!\cdot\left|\det T_{i}\right|^{-\frac{1}{p}}\!\cdot\!\int_{Q_{i}}\:\max_{\left|\alpha\right|\leq\left\lceil K+\frac{\dimension+\varepsilon}{p}\right\rceil }\left|\left(\partial^{\alpha}\widehat{\gamma}\right)\!\left(S_{j}^{-1}\eta\right)\right|\d\eta\\
\left({\scriptstyle \xi=S_{i}^{-1}\eta}\right) & \leq C_{0}\cdot\left(1\!+\!\left\Vert T_{j}^{-1}T_{i}\right\Vert \right)^{\!\left\lceil K+\frac{\dimension+\varepsilon}{p}\right\rceil }\!\cdot\left|\det T_{i}\right|^{1-\frac{1}{p}}\!\cdot\!\int_{Q_{i}'}\:\max_{\left|\alpha\right|\leq\left\lceil K+\frac{\dimension+\varepsilon}{p}\right\rceil }\left|\left(\partial^{\alpha}\widehat{\gamma}\right)\!\left(S_{j}^{-1}S_{i}\xi\right)\right|\d\xi,
\end{align*}
with
\[
C_{0}:=\Omega_{0}^{K}\Omega_{1}\cdot\left(4\cdot\dimension\right)^{1+2\left\lceil K+\frac{\dimension+\varepsilon}{p}\right\rceil }\cdot\left(\frac{s_{\dimension}}{\varepsilon}\right)^{1/p}\cdot\max_{\left|\alpha\right|\leq\left\lceil K+\frac{\dimension+\varepsilon}{p}\right\rceil }C^{\left(\alpha\right)},
\]
where the constants $C^{\left(\alpha\right)}$ are defined in Assumption
\ref{assu:RegularPartitionOfUnity}, equation (\ref{eq:RegularBAPUCondition}).
\end{lem}
\begin{rem*}
With the notation $\gamma^{\left\llbracket j\right\rrbracket }$,
the usual notation $\gamma^{\left(j\right)}$ for a family $\Gamma=\left(\gamma_{i}\right)_{i\in I}$
takes the form $\gamma^{\left(j\right)}=\gamma_{j}^{\left\llbracket j\right\rrbracket }$.
\end{rem*}
\begin{proof}
Set $N:=\left\lceil K+\frac{\dimension+\varepsilon}{p}\right\rceil $.
 Now, recall from \cite[Lemma 2.6]{DecompositionIntoSobolev} the
identity
\[
\left(\partial^{\alpha}\left[f\circ A\right]\right)\left(x\right)=\sum_{\ell_{1},\dots,\ell_{k}\in\underline{d}}\left[A_{\ell_{1},i_{1}}\cdots A_{\ell_{k},i_{k}}\cdot\left(\partial_{\ell_{1}}\cdots\partial_{\ell_{k}}f\right)\left(Ax\right)\right]\qquad\forall x\in\R^{\dimension}
\]
for arbitrary $A\in\GL\left(\R^{\dimension}\right)$, $k\in\N$, $f\in C^{k}\left(\R^{\dimension}\right)$
and $\alpha=\sum_{m=1}^{k}e_{i_{m}}\in\N_{0}^{\dimension}$, where
$\left(e_{1},\dots,e_{\dimension}\right)$ is the standard basis of
$\R^{\dimension}$. In particular, this implies for arbitrary $k\in\N$
that
\[
\left|\left(\partial^{\alpha}\left[f\circ A\right]\right)\left(x\right)\right|\leq\dimension^{k}\cdot\left\Vert A\right\Vert ^{k}\cdot\max_{\left|\beta\right|=k}\left|\left(\partial^{\beta}f\right)\left(Ax\right)\right|\qquad\forall\,x\in\R^{\dimension}\,\forall\,\alpha\in\N_{0}^{\dimension}\text{ with }\left|\alpha\right|=k\text{ and }f\in C^{k}\left(\smash{\R^{\dimension}}\right)
\]
and this estimate obviously also holds for $k=0$.

Thus, using the notation $h^{\heartsuit}\left(\xi\right):=\max_{\left|\alpha\right|\leq N}\left|\left(\partial^{\alpha}h\right)\left(\xi\right)\right|$
for $\xi\in\R^{\dimension}$ and $h\in C^{\infty}\left(\R^{\dimension}\right)$,
we get for arbitrary $i\in I$, $m\in\underline{\dimension}$ and
$\ell\in\left\{ 0\right\} \cup\underline{N}=\left\{ 0,\dots,\left\lceil K+\frac{\dimension+\varepsilon}{p}\right\rceil \right\} $
as well as $T\in\GL\left(\R^{\dimension}\right)$ that
\begin{align*}
\left|\left[\partial_{m}^{\ell}\left(h\circ S_{i}^{-1}\circ T\right)\right]\left(\xi\right)\right| & =\left|\left[\partial_{m}^{\ell}\left(\eta\mapsto h\left(T_{i}^{-1}T\eta-T_{i}^{-1}b_{i}\right)\right)\right]\left(\xi\right)\right|\\
 & \leq\dimension^{\ell}\cdot\left\Vert T_{i}^{-1}T\right\Vert ^{\ell}\cdot\max_{\left|\alpha\right|\leq\ell}\left|\left(\partial^{\alpha}h\right)\left(T_{i}^{-1}T\xi-T_{i}^{-1}b_{i}\right)\right|\\
\left({\scriptstyle \text{since }\ell\leq N}\right) & \leq\dimension^{\ell}\cdot\left\Vert T_{i}^{-1}T\right\Vert ^{\ell}\cdot\left(h^{\heartsuit}\circ S_{i}^{-1}\circ T\right)\left(\xi\right).
\end{align*}

Now, set $g:=\varphi_{i}\cdot\widehat{\gamma^{\left\llbracket j\right\rrbracket }}\in\TestFunctionSpace{\R^{\dimension}}$
and apply Leibniz's rule to the product
\[
g\circ T=\left(\widehat{\gamma^{\left\llbracket j\right\rrbracket }}\circ T\right)\cdot\left(\varphi_{i}\circ T\right)=\left(\widehat{\gamma}\circ S_{j}^{-1}\circ T\right)\cdot\left(\varphi_{i}^{\natural}\circ S_{i}^{-1}\circ T\right),
\]
to see using the binomial theorem that
\begin{align*}
\left|\partial_{m}^{N}\left(g\circ T\right)\right| & =\left|\sum_{\ell=0}^{N}\binom{N}{\ell}\cdot\partial_{m}^{\ell}\left(\widehat{\gamma}\circ S_{j}^{-1}\circ T\right)\cdot\partial_{m}^{N-\ell}\left(\varphi_{i}^{\natural}\circ S_{i}^{-1}\circ T\right)\right|\\
 & \leq\left[\widehat{\gamma}^{\heartsuit}\circ S_{j}^{-1}\circ T\right]\cdot\left[\left(\smash{\varphi_{i}^{\natural}}\right)^{\heartsuit}\circ S_{i}^{-1}\circ T\right]\cdot\sum_{\ell=0}^{N}\left[\binom{N}{\ell}\cdot\dimension^{N}\cdot\left\Vert T_{j}^{-1}T\right\Vert ^{\ell}\left\Vert T_{i}^{-1}T\right\Vert ^{N-\ell}\right]\\
 & =\dimension^{N}\cdot\left(\left\Vert T_{j}^{-1}T\right\Vert +\left\Vert T_{i}^{-1}T\right\Vert \right)^{N}\cdot\left[\widehat{\gamma}^{\heartsuit}\circ S_{j}^{-1}\circ T\right]\cdot\left[\left(\smash{\varphi_{i}^{\natural}}\right)^{\heartsuit}\circ S_{i}^{-1}\circ T\right].
\end{align*}

Now, set $C_{2}:=\max_{\left|\alpha\right|\leq N}C^{\left(\alpha\right)}$,
with $C^{\left(\alpha\right)}$ as in Assumption \ref{assu:RegularPartitionOfUnity},
equation (\ref{eq:RegularBAPUCondition}). Because of $\supp\varphi_{i}^{\natural}\subset Q_{i}'$,
this yields $\left(\smash{\varphi_{i}^{\natural}}\right)^{\heartsuit}\leq C_{2}\cdot\Indicator_{Q_{i}'}$
and thus
\[
\left(\smash{\varphi_{i}^{\natural}}\right)^{\heartsuit}\circ S_{i}^{-1}\circ T\leq C_{2}\cdot\Indicator_{T^{-1}\left(Q_{i}\right)}=C_{2}\cdot\Indicator_{Q_{i}}\circ T.
\]
Hence,
\begin{align*}
\left|\partial_{m}^{N}\left(g\circ T\right)\right| & \leq\dimension^{N}C_{2}\cdot\left(\left\Vert T_{j}^{-1}T\right\Vert +\left\Vert T_{i}^{-1}T\right\Vert \right)^{N}\cdot\left[\widehat{\gamma}^{\heartsuit}\circ S_{j}^{-1}\circ T\right]\cdot\left(\Indicator_{Q_{i}}\circ T\right)\\
 & =\dimension^{N}C_{2}\cdot\left(\left\Vert T_{j}^{-1}T\right\Vert +\left\Vert T_{i}^{-1}T\right\Vert \right)^{N}\cdot\left[\left(\widehat{\gamma}^{\heartsuit}\circ S_{j}^{-1}\right)\cdot\Indicator_{Q_{i}}\right]\circ T,
\end{align*}
as well as
\[
\left|g\circ T\right|=\left|\left(\widehat{\gamma}\circ S_{j}^{-1}\circ T\right)\cdot\left(\varphi_{i}^{\natural}\circ S_{i}^{-1}\circ T\right)\right|\leq C_{2}\cdot\left[\left(\widehat{\gamma}^{\heartsuit}\circ S_{j}^{-1}\right)\cdot\Indicator_{Q_{i}}\right]\circ T.
\]

By combining Lemma \ref{lem:PointwiseFourierDecayEstimate}, equation
(\ref{eq:PointwiseFourierDecayEstimate}) (with $g\circ T$ instead
of $g$) with the preceding estimates, we arrive at
\begin{align*}
\left(1+\left|x\right|\right)^{N}\cdot\left|\left[\Fourier^{-1}\left(g\circ T\right)\right]\left(x\right)\right| & \leq\left(1+\dimension\right)^{N}\cdot\left(\left|\left[\Fourier^{-1}\left(g\circ T\right)\right]\left(x\right)\right|+\sum_{m=1}^{\dimension}\left|\left[\Fourier^{-1}\left(\partial_{m}^{N}\left(g\circ T\right)\right)\right]\left(x\right)\right|\right)\\
 & \leq\left(1+\dimension\right)^{N}\cdot\left(\left\Vert g\circ T\right\Vert _{L^{1}}+\sum_{m=1}^{\dimension}\left\Vert \partial_{m}^{N}\left(g\circ T\right)\right\Vert _{L^{1}}\right)\\
 & \leq\left(1\!+\!\dimension\right)^{N}C_{2}\cdot\left\Vert \left[\left(\widehat{\gamma}^{\heartsuit}\circ S_{j}^{-1}\right)\cdot\Indicator_{Q_{i}}\right]\circ T\right\Vert _{L^{1}}\cdot\left(1\!+\!\sum_{m=1}^{\dimension}\dimension^{N}\cdot\left(\left\Vert T_{j}^{-1}T\right\Vert \!+\!\left\Vert T_{i}^{-1}T\right\Vert \right)^{N}\right)\\
 & \overset{\left(\dagger\right)}{\leq}\dimension^{N+1}\left(1+\dimension\right)^{N}C_{2}\cdot\left(1+\left\Vert T_{j}^{-1}T\right\Vert +\left\Vert T_{i}^{-1}T\right\Vert \right)^{N}\cdot\left|\det T\right|^{-1}\cdot\left\Vert \left(\widehat{\gamma}^{\heartsuit}\circ S_{j}^{-1}\right)\cdot\Indicator_{Q_{i}}\right\Vert _{L^{1}}\\
 & \leq\left(1+\dimension\right)^{1+2N}\cdot C_{2}\cdot\left(1+\left\Vert T_{j}^{-1}T\right\Vert +\left\Vert T_{i}^{-1}T\right\Vert \right)^{N}\cdot\left|\det T\right|^{-1}\cdot\left\Vert \left(\widehat{\gamma}^{\heartsuit}\circ S_{j}^{-1}\right)\cdot\Indicator_{Q_{i}}\right\Vert _{L^{1}}\\
 & =:\left|\det T\right|^{-1}\cdot C_{3}^{\left(i,j,T\right)},
\end{align*}
where the step marked with $\left(\dagger\right)$ used that $1+a^{N}\leq\left(1+a\right)^{N}$
for $a\geq0$, as can be seen by expanding the right-hand side using
the binomial theorem.

\medskip{}

Now, we choose $T=T_{i}$ and note $\left[\Fourier^{-1}\left(g\circ T_{i}\right)\right]\left(x\right)=\left|\det T_{i}\right|^{-1}\cdot\left(\Fourier^{-1}g\right)\left(T_{i}^{-T}x\right)$,
so that we have shown $\left|\left(\Fourier^{-1}g\right)\left(T_{i}^{-T}x\right)\right|\leq C_{3}^{\left(i,j,T_{i}\right)}\cdot\left(1+\left|x\right|\right)^{-N}$
and thus $\left|\left(\Fourier^{-1}g\right)\left(y\right)\right|\leq C_{3}^{\left(i,j,T_{i}\right)}\cdot\left(1+\left|T_{i}^{T}y\right|\right)^{-N}$
for all $y\in\R^{\dimension}$. In conjunction with equation (\ref{eq:WeightLinearTransformationsConnection})
and because of $v_{0}\left(y\right)\leq\Omega_{1}\cdot\left(1+\left|y\right|\right)^{K}$,
we arrive at
\begin{align*}
v_{0}\left(y\right)\cdot\left|\left(\Fourier^{-1}g\right)\left(y\right)\right| & \leq\Omega_{1}\cdot\left(1+\left|y\right|\right)^{K}\cdot\left|\left(\Fourier^{-1}g\right)\left(y\right)\right|\\
\left({\scriptstyle \text{eq. }\eqref{eq:WeightLinearTransformationsConnection}}\right) & \leq\Omega_{0}^{K}\Omega_{1}\cdot C_{3}^{\left(i,j,T_{i}\right)}\cdot\left(1+\left|T_{i}^{T}y\right|\right)^{K-N}.
\end{align*}
By taking the $L^{p}$-quasi-norm of this estimate, we arrive at
\begin{align*}
\left\Vert \Fourier^{-1}g\right\Vert _{L_{v_{0}}^{p}} & \leq\Omega_{0}^{K}\Omega_{1}\cdot C_{3}^{\left(i,j,T_{i}\right)}\cdot\left\Vert \left(1+\left|T_{i}^{T}\mybullet\right|\right)^{K-N}\right\Vert _{L^{p}}\\
 & =\Omega_{0}^{K}\Omega_{1}\cdot C_{3}^{\left(i,j,T_{i}\right)}\cdot\left|\det T_{i}\right|^{-1/p}\cdot\left\Vert \left(1+\left|\mybullet\right|\right)^{-\left(N-K\right)}\right\Vert _{L^{p}}\\
\left({\scriptstyle \text{eq. }\eqref{eq:StandardDecayLpEstimate}}\right) & \leq\Omega_{0}^{K}\Omega_{1}\cdot C_{3}^{\left(i,j,T_{i}\right)}\cdot\left|\det T_{i}\right|^{-1/p}\cdot\left(\frac{s_{\dimension}}{\varepsilon}\right)^{1/p},
\end{align*}
where the last step used our choice of $N=\left\lceil K+\frac{\dimension+\varepsilon}{p}\right\rceil $.

This proves the claim, since
\begin{align*}
C_{3}^{\left(i,j,T_{i}\right)} & =\left(1+\dimension\right)^{1+2N}\cdot C_{2}\cdot\left(2+\left\Vert T_{j}^{-1}T_{i}\right\Vert \right)^{N}\cdot\left\Vert \left(\widehat{\gamma}^{\heartsuit}\circ S_{j}^{-1}\right)\cdot\Indicator_{Q_{i}}\right\Vert _{L^{1}}\\
 & \leq\left(4\dimension\right)^{1+2N}\cdot C_{2}\cdot\left(1+\left\Vert T_{j}^{-1}T_{i}\right\Vert \right)^{N}\cdot\int_{Q_{i}}\:\max_{\left|\alpha\right|\leq\left\lceil K+\frac{\dimension+\varepsilon}{p}\right\rceil }\left|\left(\partial^{\alpha}\widehat{\gamma}\right)\left(S_{j}^{-1}\xi\right)\right|\d\xi.\qedhere
\end{align*}
\end{proof}
As a consequence of the preceding estimate, we see in particular that
every regular $\CalQ$-BAPU is also a $\CalQ$-$v_{0}$-BAPU, even
for $v_{0}\not\equiv1$.
\begin{cor}
\label{cor:RegularBAPUsAreWeightedBAPUs}Every regular $\CalQ$-BAPU
$\Phi=\left(\varphi_{i}\right)_{i\in I}$ is a $\CalQ$-$v_{0}$-BAPU.

In fact, there is some $\varrho\in\TestFunctionSpace{\R^{\dimension}}$,
depending only on $Q:=\overline{\bigcup_{i\in I}Q_{i}'}$ (and thus
only on $\CalQ$), such that
\[
C_{\CalQ,\Phi,v_{0},p}\leq\Omega_{0}^{K}\Omega_{1}\cdot\left(4\cdot\dimension\right)^{1+2\left\lceil K+\frac{\dimension+\varepsilon}{p}\right\rceil }\cdot\left(\frac{s_{\dimension}}{\varepsilon}\right)^{1/p}\cdot2^{\!\left\lceil K+\frac{\dimension+\varepsilon}{p}\right\rceil }\cdot\lambda_{\dimension}\left(Q\right)\cdot\max_{\left|\alpha\right|\leq\left\lceil K+\frac{\dimension+\varepsilon}{p}\right\rceil }\left\Vert \partial^{\alpha}\varrho\right\Vert _{\sup}\cdot\max_{\left|\alpha\right|\leq\left\lceil K+\frac{\dimension+\varepsilon}{p}\right\rceil }C^{\left(\alpha\right)},
\]
where $\varepsilon>0$ can be chosen arbitrarily.
\end{cor}
\begin{proof}
The set $Q\subset\R^{\dimension}$ is compact, so that there is some
$\gamma\in\Schwartz\left(\R^{\dimension}\right)$ satisfying $\widehat{\gamma}\in\TestFunctionSpace{\R^{\dimension}}$
and $\gamma\equiv1$ on $Q$. In the notation of Lemma \ref{lem:GramMatrixEstimates},
this entails $\widehat{\gamma^{\left\llbracket j\right\rrbracket }}=\widehat{\gamma}\circ S_{j}^{-1}\equiv1$
on $S_{j}Q\supset S_{j}\overline{Q_{j}'}=\overline{Q_{j}}$. But because
of $\varphi_{j}\equiv0$ outside of $\overline{Q_{j}}$, this implies
$\widehat{\gamma^{\left\llbracket j\right\rrbracket }}\cdot\varphi_{j}=\varphi_{j}$,
so that Lemma \ref{lem:GramMatrixEstimates} yields (with $C_{0}$
as in that lemma) that 
\begin{align*}
\left\Vert \Fourier^{-1}\varphi_{j}\right\Vert _{L_{v_{0}}^{p}} & =\left\Vert \Fourier^{-1}\left(\widehat{\gamma^{\left\llbracket j\right\rrbracket }}\cdot\varphi_{j}\right)\right\Vert _{L_{v_{0}}^{p}}\\
 & \leq C_{0}\cdot\left(1\!+\!\left\Vert T_{j}^{-1}T_{j}\right\Vert \right)^{\!\left\lceil K+\frac{\dimension+\varepsilon}{p}\right\rceil }\!\cdot\left|\det T_{j}\right|^{1-\frac{1}{p}}\!\cdot\!\int_{Q_{j}'}\:\max_{\left|\alpha\right|\leq\left\lceil K+\frac{\dimension+\varepsilon}{p}\right\rceil }\left|\left(\partial^{\alpha}\widehat{\gamma}\right)\!\left(S_{j}^{-1}S_{j}\xi\right)\right|\d\xi\\
 & \leq C_{0}\cdot2^{\!\left\lceil K+\frac{\dimension+\varepsilon}{p}\right\rceil }\!\cdot\left|\det T_{j}\right|^{1-\frac{1}{p}}\!\cdot\!\lambda_{\dimension}\left(Q_{j}'\right)\cdot\max_{\left|\alpha\right|\leq\left\lceil K+\frac{\dimension+\varepsilon}{p}\right\rceil }\left\Vert \partial^{\alpha}\widehat{\gamma}\right\Vert _{\sup}\\
 & \leq C_{0}\cdot2^{\!\left\lceil K+\frac{\dimension+\varepsilon}{p}\right\rceil }\cdot\lambda_{\dimension}\left(Q\right)\cdot\max_{\left|\alpha\right|\leq\left\lceil K+\frac{\dimension+\varepsilon}{p}\right\rceil }\left\Vert \partial^{\alpha}\widehat{\gamma}\right\Vert _{\sup}\cdot\left|\det T_{j}\right|^{1-\frac{1}{p}}\\
 & =:C\cdot\left|\det T_{j}\right|^{1-\frac{1}{p}},
\end{align*}
where $C>0$ is independent of $j\in I$. Recalling the definition
of a $\CalQ$-$v_{0}$-BAPU from Subsection \ref{subsec:DecompSpaceDefinitionStandingAssumptions},
this yields the claim, with $\varrho:=\widehat{\gamma}$.
\end{proof}
Using Schur's test as well as the estimates given in Lemma \ref{lem:GramMatrixEstimates},
we can now derive simplified sufficient criteria which ensure that
a given family $\Gamma=\left(\gamma_{i}\right)_{i\in I}$ of prototypes
indeed generates a Banach frame (as in Theorem \ref{thm:DiscreteBanachFrameTheorem})
or an atomic decomposition (as in Theorem \ref{thm:AtomicDecomposition}).
We start with a simplified criterion for Banach frames.
\begin{cor}
\label{cor:BanachFrameSimplifiedCriteria}Assume that $\left(\varphi_{i}\right)_{i\in I}$
satisfies Assumption \ref{assu:RegularPartitionOfUnity}. Then, for
each $p,q\in\left(0,\infty\right]$, there are
\[
N\in\N,\qquad\sigma>0,\qquad\text{ and }\qquad\tau>0
\]
with the following property: If the family $\Gamma=\left(\gamma_{i}\right)_{i\in I}$
satisfies the following:

\begin{enumerate}
\item We have $\gamma_{i}\in L_{\left(1+\left|\mybullet\right|\right)^{K}}^{1}\left(\R^{\dimension}\right)$
and $\widehat{\gamma_{i}}\in C^{\infty}\left(\R^{\dimension}\right)$
for all $i\in I$, where all partial derivatives of $\widehat{\gamma_{i}}$
are polynomially bounded.
\item We have $\gamma_{i}\in C^{1}\left(\R^{\dimension}\right)$ and $\partial_{\ell}\gamma_{i}\in L_{v_{0}}^{1}\left(\R^{\dimension}\right)\cap L^{\infty}\left(\R^{\dimension}\right)$
for all $\ell\in\underline{\dimension}$ and $i\in I$.
\item The family $\Gamma=\left(\gamma_{i}\right)_{i\in I}$ satisfies Assumption
\ref{assu:GammaCoversOrbit}.
\item We have
\[
C_{1}:=\sup_{i\in I}\,\sum_{j\in I}M_{j,i}<\infty\qquad\text{ and }\qquad C_{2}:=\sup_{j\in I}\sum_{i\in I}M_{j,i}<\infty
\]
with
\[
M_{j,i}:=\left(\frac{w_{j}}{w_{i}}\right)^{\tau}\cdot\left(1+\left\Vert T_{j}^{-1}T_{i}\right\Vert \right)^{\sigma}\cdot\max_{\left|\beta\right|\leq1}\left(\left|\det T_{i}\right|^{-1}\cdot\int_{Q_{i}}\max_{\left|\alpha\right|\leq N}\left|\left(\partial^{\alpha}\widehat{\partial^{\beta}\gamma_{j}}\right)\!\!\left(S_{j}^{-1}\xi\right)\right|\d\xi\right)^{\tau}.
\]
\end{enumerate}
Then $\Gamma$ fulfills Assumptions \ref{assu:DiscreteBanachFrameAssumptions}
and \ref{assu:GammaCoversOrbit} and thus all assumptions of Theorem
\ref{thm:DiscreteBanachFrameTheorem}.

In fact, the following choices are possible, for an arbitrary $\varepsilon>0$:
\begin{align*}
N & =\left\lceil K+\frac{\dimension+\varepsilon}{\min\left\{ 1,p\right\} }\right\rceil \,,\\
\tau & =\min\left\{ 1,p,q\right\} =\begin{cases}
\min\left\{ 1,q\right\} , & \text{if }p\in\left[1,\infty\right],\\
\min\left\{ q,p\right\} , & \text{if }p\in\left(0,1\right),
\end{cases}\\
\sigma & =\tau\cdot\left(\frac{\dimension}{\min\left\{ 1,p\right\} }+K+\left\lceil K+\frac{\dimension+\varepsilon}{\min\left\{ 1,p\right\} }\right\rceil \right)=\begin{cases}
\min\left\{ 1,q\right\} \cdot\left(\dimension+K+\left\lceil K+\dimension+\varepsilon\right\rceil \right), & \text{if }p\in\left[1,\infty\right],\\
\min\left\{ p,q\right\} \cdot\left(\frac{\dimension}{p}+K+\left\lceil K+\frac{\dimension+\varepsilon}{p}\right\rceil \right), & \text{if }p\in\left(0,1\right).
\end{cases}
\end{align*}
With these choices, we even have $\vertiii{\smash{\overrightarrow{A}}}^{\max\left\{ 1,\frac{1}{p}\right\} }\leq C\cdot\left(C_{1}^{1/\tau}+C_{2}^{1/\tau}\right)$
and $\vertiii{\smash{\overrightarrow{B}}}^{\max\left\{ 1,\frac{1}{p}\right\} }\leq C\cdot\left(C_{1}^{1/\tau}+C_{2}^{1/\tau}\right)$
for
\[
C:=\Omega_{0}^{K}\Omega_{1}\cdot\dimension^{1/\min\left\{ 1,p\right\} }\cdot\left(4\cdot\dimension\right)^{1+2\left\lceil K+\frac{\dimension+\varepsilon}{\min\left\{ 1,p\right\} }\right\rceil }\cdot\left(\frac{s_{\dimension}}{\varepsilon}\right)^{1/\min\left\{ 1,p\right\} }\cdot\max_{\left|\alpha\right|\leq\left\lceil K+\frac{\dimension+\varepsilon}{\min\left\{ 1,p\right\} }\right\rceil }C^{\left(\alpha\right)}.\qedhere
\]
\end{cor}
\begin{rem*}
As usual, the most important special case is when $\gamma_{i}=\gamma$
is independent of $i\in I$. In this case, validity of Assumption
\ref{assu:GammaCoversOrbit} can be verified easily using Lemma \ref{lem:GammaCoversOrbitAssumptionSimplified}.
The same lemma is also highly helpful if $\left\{ \gamma_{i}\with i\in I\right\} $
is finite, i.e., if only a finite number of different prototypes is
used.
\end{rem*}
\begin{proof}
Since our assumptions clearly include those of Assumption \ref{assu:GammaCoversOrbit},
we only need to verify Assumption \ref{assu:DiscreteBanachFrameAssumptions}.
This means the following:

\begin{itemize}
\item We have $\gamma_{i}\in C^{1}\left(\R^{\dimension}\right)$ and the
gradient $\phi_{i}:=\nabla\gamma_{i}$ is bounded and satisfies $\phi_{i}\in L_{v_{0}}^{1}\left(\R^{\dimension};\Compl^{\dimension}\right)$,
as well as $\widehat{\phi_{i}}\in C^{\infty}\left(\R^{\dimension};\Compl^{\dimension}\right)$.
All of these properties except the last are included in our assumptions.
But standard properties of the Fourier transform show $\widehat{\partial_{\ell}\gamma_{i}}\left(\xi\right)=2\pi i\xi_{\ell}\cdot\widehat{\gamma_{i}}\left(\xi\right)$
for $\xi\in\R^{\dimension}$, so that $\widehat{\partial_{\ell}\gamma_{i}}\in C^{\infty}\left(\R^{\dimension}\right)$,
since $\widehat{\gamma_{i}}\in C^{\infty}\left(\R^{\dimension}\right)$.
\item Assumption \ref{assu:MainAssumptions} is satisfied. For this, it
remains—in view of our present assumptions—to check that the operator
$\overrightarrow{A}:\ell_{w^{\min\left\{ 1,p\right\} }}^{r}\left(I\right)\to\ell_{w^{\min\left\{ 1,p\right\} }}^{r}\left(I\right)$
is bounded, where $r:=\max\left\{ q,\frac{q}{p}\right\} $ and $A=\left(A_{j,i}\right)_{j,i\in I}$
is given by
\[
A_{j,i}:=\begin{cases}
\left\Vert \Fourier^{-1}\left(\varphi_{i}\cdot\widehat{\gamma^{\left(j\right)}}\right)\right\Vert _{L_{v_{0}}^{1}}, & \text{if }p\in\left[1,\infty\right],\\
\left(1+\left\Vert T_{j}^{-1}T_{i}\right\Vert \right)^{\dimension}\cdot\left|\det T_{i}\right|^{1-p}\cdot\left\Vert \Fourier^{-1}\left(\varphi_{i}\cdot\widehat{\gamma^{\left(j\right)}}\right)\right\Vert _{L_{v_{0}}^{p}}^{p}, & \text{if }p\in\left(0,1\right).
\end{cases}
\]
\item The infinite matrix $B=\left(B_{j,i}\right)_{j,i\in I}$ defines a
bounded linear operator $\overrightarrow{B}:\ell_{w^{\min\left\{ 1,p\right\} }}^{r}\left(I\right)\to\ell_{w^{\min\left\{ 1,p\right\} }}^{r}\left(I\right)$,
where $r=\max\left\{ q,\frac{q}{p}\right\} $ as above, $\phi_{i}=\nabla\gamma_{i}$
for $i\in I$ and
\[
B_{j,i}:=\begin{cases}
\left(1+\left\Vert T_{j}^{-1}T_{i}\right\Vert \right)^{K+\dimension}\cdot\left\Vert \Fourier^{-1}\left(\varphi_{i}\cdot\widehat{\phi^{\left(j\right)}}\right)\right\Vert _{L_{v_{0}}^{1}}, & \text{if }p\in\left[1,\infty\right],\\
\left(1+\left\Vert T_{j}^{-1}T_{i}\right\Vert \right)^{pK+\dimension}\cdot\left|\det T_{i}\right|^{1-p}\cdot\left\Vert \Fourier^{-1}\left(\varphi_{i}\cdot\widehat{\phi^{\left(j\right)}}\right)\right\Vert _{L_{v_{0}}^{p}}^{p}, & \text{if }p\in\left(0,1\right).
\end{cases}
\]
Here, $\phi^{\left(j\right)}$ is defined as in equation (\ref{eq:NonCompactFilterDefinition}),
with $\gamma_{j}$ replaced by $\phi_{j}$.
\end{itemize}
Hence, in the following, we verify boundedness of $\overrightarrow{A}$
and $\overrightarrow{B}$.

\medskip{}

We first make the auxiliary observation that a matrix operator $\overrightarrow{C}:\ell_{v}^{q}\left(I\right)\to\ell_{v}^{q}\left(I\right)$
is bounded if and only if the operator $\overrightarrow{C_{v}}:\ell^{q}\left(I\right)\to\ell^{q}\left(I\right)$
is bounded, where
\[
\left(C_{v}\right)_{j,i}=\frac{v_{j}}{v_{i}}\cdot C_{j,i}.
\]
This simply comes from the fact that $m_{v}:\ell_{v}^{q}\left(I\right)\to\ell^{q}\left(I\right),\left(x_{i}\right)_{i\in I}\mapsto\left(v_{i}x_{i}\right)_{i\in I}$
is an isometric isomorphism and that 
\[
\left[\overrightarrow{C}:\ell_{v}^{q}\left(I\right)\to\ell_{v}^{q}\left(I\right)\right]=m_{v}^{-1}\circ\left(m_{v}\circ\overrightarrow{C}\circ m_{v}^{-1}\right)\circ m_{v},
\]
where a direct calculation shows $m_{v}\circ\overrightarrow{C}\circ m_{v}^{-1}=\overrightarrow{C_{v}}$.
Since $m_{v}$ is isometric, we also get $\vertiii{\smash{\overrightarrow{C_{v}}}}=\vertiii{\smash{\overrightarrow{C}}}$.

\medskip{}

Now, let us first consider the case $p\in\left[1,\infty\right]$.
Here, we want to have $\overrightarrow{A}:\ell_{w}^{q}\left(I\right)\to\ell_{w}^{q}\left(I\right)$
and likewise for $\overrightarrow{B}$. Recall $\gamma^{\left(j\right)}=\gamma_{j}^{\left\llbracket j\right\rrbracket }$
in the notation of Lemma \ref{lem:GramMatrixEstimates}. Hence, an
application of that lemma (with $p=1$) yields, with $C$ as in the
statement of the present corollary,
\begin{align*}
\left(A_{w}\right)_{j,i}=\frac{w_{j}}{w_{i}}\cdot A_{j,i} & \leq\frac{C}{\dimension}\cdot\frac{w_{j}}{w_{i}}\cdot\left[1+\left\Vert T_{j}^{-1}T_{i}\right\Vert \right]^{\left\lceil K+\dimension+\varepsilon\right\rceil }\cdot\left|\det T_{i}\right|^{-1}\cdot\int_{Q_{i}}\max_{\left|\alpha\right|\leq\left\lceil K+\dimension+\varepsilon\right\rceil }\left|\left(\partial^{\alpha}\widehat{\gamma_{j}}\right)\left(S_{j}^{-1}\xi\right)\right|\d\xi\\
 & \leq C\cdot\frac{w_{j}}{w_{i}}\cdot\left[1+\left\Vert T_{j}^{-1}T_{i}\right\Vert \right]^{K+\dimension+\left\lceil K+\dimension+\varepsilon\right\rceil }\cdot\left|\det T_{i}\right|^{-1}\cdot\max_{\left|\beta\right|\leq1}\int_{Q_{i}}\max_{\left|\alpha\right|\leq N}\left|\left(\partial^{\alpha}\widehat{\partial^{\beta}\gamma_{j}}\right)\left(S_{j}^{-1}\xi\right)\right|\d\xi\\
 & \leq C\cdot M_{j,i}^{1/\min\left\{ 1,q\right\} }.
\end{align*}
Likewise, using $\left|\phi_{j}\right|=\left|\widehat{\nabla\gamma_{j}}\right|\leq\sum_{\ell=1}^{\dimension}\left|\widehat{\partial_{\ell}\gamma_{j}}\right|$,
we get
\begin{align*}
\left(B_{w}\right)_{j,i}=\frac{w_{j}}{w_{i}}\cdot B_{j,i} & \leq\dimension\cdot\frac{w_{j}}{w_{i}}\cdot\left(1+\left\Vert T_{j}^{-1}T_{i}\right\Vert \right)^{K+\dimension}\cdot\max_{\left|\beta\right|=1}\left\Vert \Fourier^{-1}\left(\varphi_{i}\cdot\widehat{\left(\partial^{\beta}\gamma_{j}\right)^{\left\llbracket j\right\rrbracket }}\right)\right\Vert _{L_{v_{0}}^{1}}\\
\left({\scriptstyle \text{Lemma }\ref{lem:GramMatrixEstimates}\text{ with }\partial^{\beta}\gamma_{j}\text{ instead of }\gamma}\right) & \leq C\!\cdot\!\frac{w_{j}}{w_{i}}\!\cdot\!\left(1\!+\!\left\Vert T_{j}^{-1}T_{i}\right\Vert \right)^{K+\dimension+\left\lceil K+\dimension+\varepsilon\right\rceil }\!\max_{\left|\beta\right|=1}\left[\left|\det T_{i}\right|^{-1}\!\!\int_{Q_{i}}\,\max_{\left|\alpha\right|\leq N}\left|\left(\partial^{\alpha}\widehat{\partial^{\beta}\gamma_{j}}\right)\!\!\left(S_{j}^{-1}\xi\right)\right|\d\xi\right]\\
 & \leq C\cdot M_{j,i}^{1/\min\left\{ 1,q\right\} }.
\end{align*}
But Lemma \ref{lem:SchursLemma} shows that $\overrightarrow{A_{w}}:\ell^{q}\left(I\right)\to\ell^{q}\left(I\right)$
is bounded as soon as we have $K_{1}:=\sup_{i\in I}\sum_{j\in I}\left(A_{w}\right)_{j,i}^{\min\left\{ 1,q\right\} }<\infty$
and $K_{2}:=\sup_{j\in I}\sum_{i\in I}\left(A_{w}\right)_{j,i}^{\min\left\{ 1,q\right\} }<\infty$.
Further, that lemma also shows $\vertiii{\smash{\overrightarrow{A_{w}}}}\leq\max\left\{ K_{1}^{1/\min\left\{ 1,q\right\} },K_{2}^{1/\min\left\{ 1,q\right\} }\right\} $.
Since we have by assumption that $C_{1}=\sup_{i\in I}\sum_{j\in I}M_{j,i}<\infty$
and $C_{2}=\sup_{j\in I}\sum_{i\in I}M_{j,i}<\infty$, we get $K_{1}^{1/\min\left\{ 1,q\right\} }\leq C\cdot C_{1}^{1/\min\left\{ 1,q\right\} }=C\cdot C_{1}^{1/\tau}$
and likewise $K_{2}^{1/\min\left\{ 1,q\right\} }\leq C\cdot C_{2}^{1/\min\left\{ 1,q\right\} }=C\cdot C_{2}^{1/\tau}$,
so that $\vertiii{\smash{\overrightarrow{A_{w}}}}\leq C\cdot\max\left\{ C_{1}^{1/\tau},C_{2}^{1/\tau}\right\} $.
The same arguments show that $\overrightarrow{B_{w}}:\ell^{q}\left(I\right)\to\ell^{q}\left(I\right)$
is bounded and satisfies $\vertiii{\smash{\overrightarrow{B_{w}}}}\leq C\cdot\max\left\{ C_{1}^{1/\tau},C_{2}^{1/\tau}\right\} $.
In view of the auxiliary observation from above, this completes the
proof in case of $p\in\left[1,\infty\right]$.

\medskip{}

Now, let $p\in\left(0,1\right)$. In this case, we want to have $\overrightarrow{A}:\ell_{w^{p}}^{q/p}\left(I\right)\to\ell_{w^{p}}^{q/p}\left(I\right)$
and likewise for $\overrightarrow{B}$. But Lemma \ref{lem:GramMatrixEstimates}
yields, because of $\gamma^{\left(j\right)}=\gamma_{j}^{\left\llbracket j\right\rrbracket }$,
\begin{align*}
\left(A_{w^{p}}\right)_{j,i} & =\left(\frac{w_{j}}{w_{i}}\right)^{p}\cdot A_{j,i}\\
 & =\left(\frac{w_{j}}{w_{i}}\right)^{p}\cdot\left(1+\left\Vert T_{j}^{-1}T_{i}\right\Vert \right)^{\dimension}\cdot\left|\det T_{i}\right|^{1-p}\cdot\left\Vert \Fourier^{-1}\left(\varphi_{i}\cdot\widehat{\gamma^{\left(j\right)}}\right)\right\Vert _{L_{v_{0}}^{p}}^{p}\\
 & \leq\left(C/\dimension\right)^{p}\cdot\left(\frac{w_{j}}{w_{i}}\right)^{p}\left(1+\left\Vert T_{j}^{-1}T_{i}\right\Vert \right)^{\dimension+p\left\lceil K+\frac{\dimension+\varepsilon}{p}\right\rceil }\left|\det T_{i}\right|^{1-p}\left|\det T_{i}\right|^{-1}\left(\int_{Q_{i}}\max_{\left|\alpha\right|\leq\left\lceil K+\frac{\dimension+\varepsilon}{p}\right\rceil }\left|\left(\partial^{\alpha}\widehat{\gamma_{j}}\right)\!\left(S_{j}^{-1}\xi\right)\right|\d\xi\right)^{p}\\
 & =\left(C/\dimension\right)^{p}\cdot\left(\frac{w_{j}}{w_{i}}\right)^{p}\left(1+\left\Vert T_{j}^{-1}T_{i}\right\Vert \right)^{\dimension+p\left\lceil K+\frac{\dimension+\varepsilon}{p}\right\rceil }\left(\left|\det T_{i}\right|^{-1}\int_{Q_{i}}\max_{\left|\alpha\right|\leq N}\left|\left(\partial^{\alpha}\widehat{\gamma_{j}}\right)\left(S_{j}^{-1}\xi\right)\right|\d\xi\right)^{p}\\
 & \leq\left(C/\dimension\right)^{p}\cdot M_{j,i}^{1/\min\left\{ 1,\,\frac{q}{p}\right\} },
\end{align*}
where the last step used
\[
\min\left\{ 1,\frac{q}{p}\right\} \cdot\left(\dimension+p\left\lceil K+\frac{\dimension+\varepsilon}{p}\right\rceil \right)=\min\left\{ p,q\right\} \cdot\left(\frac{\dimension}{p}+\left\lceil K+\frac{\dimension+\varepsilon}{p}\right\rceil \right)\leq\sigma.
\]
Furthermore, since $p\in\left(0,1\right)$, we have $\left(\sum_{\ell=1}^{\dimension}a_{\ell}\right)^{p}\leq\sum_{\ell=1}^{\dimension}a_{\ell}^{p}$
for $a_{1},\dots,a_{\dimension}\geq0$, so that the $L^{p}$-norm
of a vector valued (measurable) function $f:\R^{\dimension}\to\Compl^{\dimension}$
can be estimated as follows:
\begin{align*}
\left\Vert \left(f_{1},\dots,f_{\dimension}\right)\right\Vert _{L^{p}}^{p} & =\int_{\R^{\dimension}}\left|\left(f_{1},\dots,f_{\dimension}\right)\left(x\right)\right|^{p}\d x\leq\int_{\R^{\dimension}}\left(\sum_{\ell=1}^{\dimension}\left|f_{\ell}\left(x\right)\right|\right)^{p}\d x\\
 & \leq\int_{\R^{\dimension}}\sum_{\ell=1}^{\dimension}\left|f_{\ell}\left(x\right)\right|^{p}\d x=\sum_{\ell=1}^{\dimension}\left\Vert f_{\ell}\right\Vert _{L^{p}}^{p}\leq\dimension\cdot\max_{\ell\in\underline{\dimension}}\left\Vert f_{\ell}\right\Vert _{L^{p}}^{p}.
\end{align*}
Consequently, 
\begin{align*}
\left(B_{w^{p}}\right)_{j,i} & =\left(\frac{w_{j}}{w_{i}}\right)^{p}\cdot B_{j,i}\\
 & =\left(\frac{w_{j}}{w_{i}}\right)^{p}\cdot\left(1+\left\Vert T_{j}^{-1}T_{i}\right\Vert \right)^{pK+\dimension}\cdot\left|\det T_{i}\right|^{1-p}\cdot\left\Vert \Fourier^{-1}\left(\varphi_{i}\cdot\widehat{\phi^{\left(j\right)}}\right)\right\Vert _{L_{v_{0}}^{p}}^{p}\\
\left({\scriptstyle \text{since }\phi_{j}=\nabla\gamma_{j}}\right) & \leq\dimension\cdot\left(\frac{w_{j}}{w_{i}}\right)^{p}\cdot\left(1+\left\Vert T_{j}^{-1}T_{i}\right\Vert \right)^{pK+\dimension}\cdot\left|\det T_{i}\right|^{1-p}\cdot\max_{\left|\beta\right|=1}\left\Vert \Fourier^{-1}\left(\varphi_{i}\cdot\widehat{\left(\partial^{\beta}\gamma_{j}\right)^{\left\llbracket j\right\rrbracket }}\right)\right\Vert _{L_{v_{0}}^{p}}^{p}\\
\left({\scriptstyle \text{Lem. }\ref{lem:GramMatrixEstimates}\text{ /w }\partial^{\beta}\gamma_{j}\text{ inst. of }\gamma_{j}}\right) & \leq\dimension\cdot\left(C/\dimension^{\frac{1}{p}}\right)^{p}\cdot\left(\frac{w_{j}}{w_{i}}\right)^{p}\left(1+\left\Vert T_{j}^{-1}T_{i}\right\Vert \right)^{pK+\dimension}\cdot\left|\det T_{i}\right|^{1-p}\\
 & \phantom{\lesssim}\cdot\max_{\left|\beta\right|\leq1}\left[\left(1+\left\Vert T_{j}^{-1}T_{i}\right\Vert \right)^{\left\lceil K+\frac{\dimension+\varepsilon}{p}\right\rceil }\left|\det T_{i}\right|^{-\frac{1}{p}}\int_{Q_{i}}\max_{\left|\alpha\right|\leq\left\lceil K+\frac{\dimension+\varepsilon}{p}\right\rceil }\left|\left(\partial^{\alpha}\widehat{\partial^{\beta}\gamma_{j}}\right)\!\left(S_{j}^{-1}\xi\right)\right|\d\xi\right]^{p}\\
 & =C^{p}\cdot\left(\frac{w_{j}}{w_{i}}\right)^{p}\!\left(1\!+\!\left\Vert T_{j}^{-1}T_{i}\right\Vert \right)^{p\left(\frac{\dimension}{p}+K+\left\lceil K+\frac{\dimension+\varepsilon}{p}\right\rceil \right)}\\
 & \phantom{=}\cdot\max_{\left|\beta\right|\leq1}\left[\left|\det T_{i}\right|^{-1}\!\int_{Q_{i}}\max_{\left|\alpha\right|\leq N}\left|\left(\partial^{\alpha}\widehat{\partial^{\beta}\gamma_{j}}\right)\!\left(S_{j}^{-1}\xi\right)\right|\d\xi\right]^{p}\\
 & =C^{p}\cdot M_{j,i}^{1/\min\left\{ 1,\,\frac{q}{p}\right\} }.
\end{align*}

Now, the remainder of the proof is similar to the case $p\in\left[1,\infty\right]$:
Lemma \ref{lem:SchursLemma} shows that $\overrightarrow{A_{w^{p}}}:\ell^{q/p}\left(I\right)\to\ell^{q/p}\left(I\right)$
is bounded as soon as we have $K_{3}:=\sup_{i\in I}\sum_{j\in I}\left(A_{w^{p}}\right)_{j,i}^{\min\left\{ 1,\frac{q}{p}\right\} }<\infty$
and $K_{4}:=\sup_{j\in I}\sum_{i\in I}\left(A_{w^{p}}\right)_{j,i}^{\min\left\{ 1,\frac{q}{p}\right\} }<\infty$.
Further, that lemma also shows $\vertiii{\smash{\overrightarrow{A_{w^{p}}}}}\leq\max\left\{ K_{3}^{1/\min\left\{ 1,q/p\right\} },K_{4}^{1/\min\left\{ 1,q/p\right\} }\right\} $.
Since we have by assumption that $C_{1}=\sup_{i\in I}\sum_{j\in I}M_{j,i}<\infty$
and $C_{2}=\sup_{j\in I}\sum_{i\in I}M_{j,i}<\infty$, we get 
\begin{align*}
\vertiii{\smash{\overrightarrow{A_{w^{p}}}}}^{\max\left\{ 1,\frac{1}{p}\right\} } & =\vertiii{\smash{\overrightarrow{A_{w^{p}}}}}^{1/p}\\
 & \leq\max\left\{ K_{3}^{\frac{1}{p}\cdot\frac{1}{\min\left\{ 1,q/p\right\} }},\,K_{4}^{\frac{1}{p}\cdot\frac{1}{\min\left\{ 1,q/p\right\} }}\right\} \\
 & =\max\left\{ K_{3}^{1/\min\left\{ p,q\right\} },K_{4}^{1/\min\left\{ p,q\right\} }\right\} \\
 & \leq\max\left\{ \left(\sup_{i\in I}\sum_{j\in I}\left[C^{p\cdot\min\left\{ 1,\frac{q}{p}\right\} }\cdot M_{j,i}\right]\right)^{1/\min\left\{ p,q\right\} },\,\left(\sup_{j\in I}\sum_{i\in I}\left[C^{p\cdot\min\left\{ 1,\frac{q}{p}\right\} }\cdot M_{j,i}\right]\right)^{1/\min\left\{ p,q\right\} }\right\} \\
 & \leq\max\left\{ C\cdot C_{1}^{1/\tau},\,C\cdot C_{2}^{1/\tau}\right\} =C\cdot\max\left\{ C_{1}^{1/\tau},C_{2}^{1/\tau}\right\} .
\end{align*}
Exactly the same arguments also yield $\vertiii{\smash{\overrightarrow{B_{w^{p}}}}}^{\max\left\{ 1,\frac{1}{p}\right\} }\leq C\cdot\max\left\{ C_{1}^{1/\tau},C_{2}^{1/\tau}\right\} $.
\end{proof}
Our next result yields simplified criteria for the application of
Theorem \ref{thm:AtomicDecomposition}, which yields atomic decompositions
for $\DecompSp{\CalQ}p{\ell_{w}^{q}}v$.
\begin{cor}
\label{cor:AtomicDecompositionSimplifiedCriteria}Assume that $\left(\varphi_{i}\right)_{i\in I}$
satisfies Assumption \ref{assu:RegularPartitionOfUnity}. Then, for
each $p,q\in\left(0,\infty\right]$, there are
\[
N\in\N,\qquad\sigma>0,\qquad\vartheta\geq0\qquad\text{ and }\qquad\tau>0
\]
with the following property: If the families $\Gamma=\left(\gamma_{i}\right)_{i\in I}$
and $\Gamma_{\ell}=\left(\gamma_{i,\ell}\right)_{i\in I}$ (with $\ell\in\left\{ 1,2\right\} $)
satisfy the following properties:

\begin{enumerate}
\item All $\gamma_{i},\gamma_{i,1},\gamma_{i,2}$ are measurable functions
$\R^{\dimension}\to\Compl$,
\item We have $\gamma_{i,1}\in L_{\left(1+\left|\mybullet\right|\right)^{K}}^{1}\left(\R^{\dimension}\right)$
for all $i\in I$.
\item We have $\gamma_{i,2}\in C^{1}\left(\R^{\dimension}\right)$ for all
$i\in I$.
\item With $K_{0}:=K+\frac{\dimension}{\min\left\{ 1,p\right\} }+1$, we
have
\[
\Omega_{4}^{\left(p,K\right)}:=\sup_{i\in I}\left\Vert \gamma_{i,2}\right\Vert _{K_{0}}+\sup_{i\in I}\left\Vert \nabla\gamma_{i,2}\right\Vert _{K_{0}}<\infty,
\]
where $\left\Vert f\right\Vert _{K_{0}}:=\sup_{x\in\R^{\dimension}}\left(1+\left|x\right|\right)^{K_{0}}\left|f\left(x\right)\right|$.
\item We have $\left\Vert \gamma_{i}\right\Vert _{K_{0}}<\infty$ for all
$i\in I$.
\item We have $\gamma_{i}=\gamma_{i,1}\ast\gamma_{i,2}$ for all $i\in I$.
\item We have $\widehat{\gamma_{i,1}},\widehat{\gamma_{i,2}}\in C^{\infty}\left(\R^{\dimension}\right)$
and all partial derivatives of $\widehat{\gamma_{i,1}},\widehat{\gamma_{i,2}}$
are polynomially bounded for all $i\in I$.
\item $\Gamma=\left(\gamma_{i}\right)_{i\in I}$ satisfies Assumption \ref{assu:GammaCoversOrbit}.
\item \label{enu:AtomicDecompositionSimplifiedMainCondition}We have
\[
K_{1}:=\sup_{i\in I}\,\sum_{j\in I}N_{i,j}<\infty\qquad\text{ and }\qquad K_{2}:=\sup_{j\in I}\sum_{i\in I}N_{i,j}<\infty
\]
with
\[
\qquad\qquad N_{i,j}:=\left(\frac{w_{i}}{w_{j}}\cdot\left(\left|\det T_{j}\right|\big/\left|\det T_{i}\right|\right)^{\vartheta}\right)^{\tau}\cdot\left(1+\left\Vert T_{j}^{-1}T_{i}\right\Vert \right)^{\sigma}\cdot\left(\left|\det T_{i}\right|^{-1}\cdot\int_{Q_{i}}\:\max_{\left|\alpha\right|\leq N}\left|\left(\partial^{\alpha}\widehat{\gamma_{j,1}}\right)\left(S_{j}^{-1}\xi\right)\right|\d\xi\right)^{\tau}\!.
\]
\end{enumerate}
Then the families $\Gamma,\Gamma_{1},\Gamma_{2}$ fulfill Assumption
\ref{assu:AtomicDecompositionAssumption} and the family $\Gamma$
satisfies Assumption \ref{assu:GammaCoversOrbit}, so that Theorem
\ref{thm:AtomicDecomposition} is applicable to $\Gamma$.

In fact, the following choices are possible, for an arbitrary $\varepsilon>0$:
\begin{align*}
N & =\left\lceil K+\frac{\dimension+\varepsilon}{\min\left\{ 1,p\right\} }\right\rceil ,\\
\tau & =\min\left\{ 1,p,q\right\} =\begin{cases}
\min\left\{ 1,q\right\} , & \text{if }p\in\left[1,\infty\right],\\
\min\left\{ p,q\right\} , & \text{if }p\in\left(0,1\right),
\end{cases}\\
\sigma & =\begin{cases}
\min\left\{ 1,q\right\} \cdot\left\lceil K+\dimension+\varepsilon\right\rceil , & \text{if }p\in\left[1,\infty\right],\\
\min\left\{ p,q\right\} \cdot\left(\frac{\dimension}{p}+K+\left\lceil K+\frac{\dimension+\varepsilon}{p}\right\rceil \right), & \text{if }p\in\left(0,1\right),
\end{cases}\\
\vartheta & =\begin{cases}
0, & \text{if }p\in\left[1,\infty\right],\\
\frac{1}{p}-1, & \text{if }p\in\left(0,1\right).
\end{cases}
\end{align*}
With these choices, we even have $\vertiii{\smash{\overrightarrow{C}}}^{\max\left\{ 1,\frac{1}{p}\right\} }\leq\Omega\cdot\left(K_{1}^{1/\tau}+K_{2}^{1/\tau}\right)$,
where $\overrightarrow{C}:\ell_{w^{\min\left\{ 1,p\right\} }}^{\max\left\{ q,q/p\right\} }\left(I\right)\to\ell_{w^{\min\left\{ 1,p\right\} }}^{\max\left\{ q,q/p\right\} }\left(I\right)$
is defined as in Assumption \ref{assu:AtomicDecompositionAssumption}
and where
\[
\Omega:=\Omega_{0}^{K}\Omega_{1}\cdot\left(4\cdot\dimension\right)^{1+2\left\lceil K+\frac{\dimension+\varepsilon}{\min\left\{ 1,p\right\} }\right\rceil }\cdot\left(\frac{s_{\dimension}}{\varepsilon}\right)^{1/\min\left\{ 1,p\right\} }\cdot\max_{\left|\alpha\right|\leq N}C^{\left(\alpha\right)}\qedhere
\]
\end{cor}
\begin{proof}
First, our assumptions clearly imply that Assumption \ref{assu:GammaCoversOrbit}
is satisfied. Hence, we only need to verify Assumption \ref{assu:AtomicDecompositionAssumption},
which means the following:

\begin{itemize}
\item We have $\gamma_{i,1},\gamma_{i,2}\in L_{\left(1+\left|\mybullet\right|\right)^{K}}^{1}\left(\R^{\dimension}\right)$
for all $i\in I$. For $\gamma_{i,1}$, this is part of our assumptions.
But for $\gamma_{i,2}$, we have $\left\Vert \gamma_{i,2}\right\Vert _{K_{0}}<\infty$,
which implies $\left\Vert \left(1+\left|\mybullet\right|\right)^{K}\cdot\gamma_{i,2}\right\Vert _{K_{0}-K}<\infty$.
Because of $K_{0}-K=\frac{\dimension}{\min\left\{ 1,p\right\} }+1\geq\dimension+1$,
this easily implies $\gamma_{i,2}\in L_{\left(1+\left|\mybullet\right|\right)^{K}}^{1}\left(\R^{\dimension}\right)$;
see also the remark after Assumption \ref{assu:AtomicDecompositionAssumption}.
\item We have $\gamma_{i}=\gamma_{i,1}\ast\gamma_{i,2}$ for all $i\in I$,
which is part of our assumptions.
\item We have $\widehat{\gamma_{i,1}},\widehat{\gamma_{i,2}}\in C^{\infty}\left(\R^{\dimension}\right)$,
with all partial derivatives of $\widehat{\gamma_{i,1}},\widehat{\gamma_{i,2}}$
being polynomially bounded. Again, this is part of our assumptions.
\item We have $\gamma_{i,2}\in C^{1}\left(\R^{\dimension}\right)$ with
$\nabla\gamma_{i,2}\in L_{v_{0}}^{1}\left(\R^{\dimension}\right)$.
The first of these properties is part of our assumptions and the second
property follows easily from $\left\Vert \nabla\gamma_{i,2}\right\Vert _{K_{0}}\leq\Omega_{4}^{\left(p,K\right)}<\infty$,
cf.\@ the remark after Assumption \ref{assu:AtomicDecompositionAssumption}.
\item We have $\Omega_{4}^{\left(p,K\right)}<\infty$, where the constant
$\Omega_{4}^{\left(p,K\right)}$ is defined as in the present corollary.
Hence, this prerequisite is part of our assumptions.
\item We have $\left\Vert \gamma_{i}\right\Vert _{K_{0}}<\infty$ for all
$i\in I$, which is again part of our assumptions.
\item The operator $\overrightarrow{C}:\ell_{w^{\min\left\{ 1,p\right\} }}^{r}\left(I\right)\to\ell_{w^{\min\left\{ 1,p\right\} }}^{r}\left(I\right)$
is bounded, where $r=\max\left\{ q,\frac{q}{p}\right\} $ and where
\[
C_{i,j}:=\begin{cases}
\left\Vert \Fourier^{-1}\left(\varphi_{i}\cdot\widehat{\gamma_{1}^{\left(j\right)}}\right)\right\Vert _{L_{v_{0}}^{1}}, & \text{if }p\in\left[1,\infty\right],\\
\left(1+\left\Vert T_{j}^{-1}T_{i}\right\Vert \right)^{pK+\dimension}\cdot\left|\det T_{j}\right|^{1-p}\cdot\left\Vert \Fourier^{-1}\left(\varphi_{i}\cdot\widehat{\gamma_{1}^{\left(j\right)}}\right)\right\Vert _{L_{v_{0}}^{p}}^{p}, & \text{if }p\in\left(0,1\right).
\end{cases}
\]
\end{itemize}
Thus, in the remainder of the proof, we only need to verify boundedness
of $\overrightarrow{C}$. First of all, we recall from the proof
of Corollary \ref{cor:BanachFrameSimplifiedCriteria} that we have
$\vertiii{\smash{\overrightarrow{C}}}=\vertiii{\overrightarrow{C_{w^{\min\left\{ 1,p\right\} }}}}_{\ell^{r}\to\ell^{r}}$,
where
\[
\left(C_{w^{\min\left\{ 1,p\right\} }}\right)_{i,j}:=\frac{w_{i}}{w_{j}}\cdot C_{i,j}.
\]
To prove boundedness of $\overrightarrow{C_{w^{\min\left\{ 1,p\right\} }}}$,
we distinguish the cases $p\in\left[1,\infty\right]$ and $p\in\left(0,1\right)$.

\medskip{}

For $p\in\left[1,\infty\right]$, we want to have $\overrightarrow{C_{w}}:\ell^{q}\left(I\right)\to\ell^{q}\left(I\right)$.
But Lemma \ref{lem:GramMatrixEstimates} (with $p=1$ and with $\gamma_{j,1}$
instead of $\gamma$) yields because of $\gamma_{1}^{\left(j\right)}=\gamma_{j,1}^{\left\llbracket j\right\rrbracket }$
that
\begin{align*}
\left(C_{w}\right)_{i,j} & =\frac{w_{i}}{w_{j}}\cdot C_{i,j}=\frac{w_{i}}{w_{j}}\cdot\left\Vert \Fourier^{-1}\left(\varphi_{i}\cdot\widehat{\gamma_{1}^{\left(j\right)}}\right)\right\Vert _{L_{v_{0}}^{1}}\\
 & \leq\Omega\cdot\frac{w_{i}}{w_{j}}\cdot\left(1+\left\Vert T_{j}^{-1}T_{i}\right\Vert \right)^{\left\lceil K+\dimension+\varepsilon\right\rceil }\cdot\left|\det T_{i}\right|^{-1}\cdot\int_{Q_{i}}\max_{\left|\alpha\right|\leq N}\left|\left(\partial^{\alpha}\widehat{\gamma_{j,1}}\right)\left(S_{j}^{-1}\eta\right)\right|\,\d\eta\\
 & =\Omega\cdot N_{i,j}^{1/\min\left\{ 1,q\right\} }.
\end{align*}
Finally, Lemma \ref{lem:SchursLemma} shows that
\begin{align*}
\vertiii{\smash{\overrightarrow{C_{w}}}}_{\ell^{q}\to\ell^{q}} & \leq\max\left\{ \left(\sup_{i\in I}\sum_{j\in I}\left|\left(C_{w}\right)_{i,j}\right|^{\min\left\{ 1,q\right\} }\right)^{1/\min\left\{ 1,q\right\} },\,\left(\sup_{j\in I}\sum_{i\in I}\left|\left(C_{w}\right)_{i,j}\right|^{\min\left\{ 1,q\right\} }\right)^{1/\min\left\{ 1,q\right\} }\right\} \\
 & \leq\Omega\cdot\max\left\{ \left(\sup_{i\in I}\sum_{j\in I}N_{i,j}\right)^{1/\min\left\{ 1,q\right\} },\,\left(\sup_{j\in I}\sum_{i\in I}N_{i,j}\right)^{1/\min\left\{ 1,q\right\} }\right\} \\
 & =\Omega\cdot\max\left\{ K_{1}^{1/\tau},\,K_{2}^{1/\tau}\right\} ,
\end{align*}
as desired.

\medskip{}

In case of $p\in\left(0,1\right)$, we want to have $\overrightarrow{C_{w^{p}}}:\ell^{q/p}\left(I\right)\to\ell^{q/p}\left(I\right)$.
But Lemma \ref{lem:GramMatrixEstimates} (with $\gamma_{j,1}$ instead
of $\gamma$) yields
\begin{align*}
\left(C_{w^{p}}\right)_{i,j} & =\left(\frac{w_{i}}{w_{j}}\right)^{p}\cdot C_{i,j}=\left(\frac{w_{i}}{w_{j}}\right)^{p}\cdot\left(1+\left\Vert T_{j}^{-1}T_{i}\right\Vert \right)^{pK+\dimension}\cdot\left|\det T_{j}\right|^{1-p}\cdot\left\Vert \Fourier^{-1}\left(\varphi_{i}\cdot\widehat{\gamma_{1}^{\left(j\right)}}\right)\right\Vert _{L_{v_{0}}^{p}}^{p}\\
 & \leq\Omega^{p}\cdot\left(\frac{w_{i}}{w_{j}}\right)^{p}\left(1\!+\!\left\Vert T_{j}^{-1}T_{i}\right\Vert \right)^{p\left(\frac{\dimension}{p}+K+\left\lceil K+\frac{\dimension+\varepsilon}{p}\right\rceil \right)}\left|\det T_{j}\right|^{1-p}\left|\det T_{i}\right|^{-1}\left(\int_{Q_{i}}\,\max_{\left|\alpha\right|\leq N}\left|\left(\partial^{\alpha}\widehat{\gamma_{j,1}}\right)\!\left(S_{j}^{-1}\eta\right)\right|\,\d\eta\right)^{p}\\
 & =\Omega^{p}\cdot\left(\frac{w_{i}}{w_{j}}\cdot\left(\frac{\left|\det T_{j}\right|}{\left|\det T_{i}\right|}\right)^{\frac{1}{p}-1}\left(1+\left\Vert T_{j}^{-1}T_{i}\right\Vert \right)^{\frac{\dimension}{p}+K+\left\lceil K+\frac{\dimension+\varepsilon}{p}\right\rceil }\left|\det T_{i}\right|^{-1}\int_{Q_{i}}\,\max_{\left|\alpha\right|\leq N}\left|\left(\partial^{\alpha}\widehat{\gamma_{j,1}}\right)\left(S_{j}^{-1}\eta\right)\right|\,\d\eta\right)^{p}\\
 & \leq\Omega^{p}\cdot N_{i,j}^{1/\min\left\{ 1,\frac{q}{p}\right\} }.
\end{align*}
Finally, Lemma \ref{lem:SchursLemma} shows that
\begin{align*}
\vertiii{\smash{\overrightarrow{C_{w^{p}}}}}_{\ell^{q/p}\to\ell^{q/p}}^{1/p} & \leq\left[\max\left\{ \left(\sup_{i\in I}\sum_{j\in I}\left(C_{w^{p}}\right)_{i,j}^{\min\left\{ 1,\frac{q}{p}\right\} }\right)^{1/\min\left\{ 1,\frac{q}{p}\right\} },\left(\sup_{j\in I}\sum_{i\in I}\left(C_{w^{p}}\right)_{i,j}^{\min\left\{ 1,\frac{q}{p}\right\} }\right)^{1/\min\left\{ 1,\frac{q}{p}\right\} }\right\} \right]^{1/p}\\
 & \leq\Omega\cdot\max\left\{ \left(\sup_{i\in I}\sum_{j\in I}N_{i,j}\right)^{1/\min\left\{ p,q\right\} },\left(\sup_{j\in I}\sum_{i\in I}N_{i,j}\right)^{1/\min\left\{ p,q\right\} }\right\} \\
 & =\Omega\cdot\max\left\{ K_{1}^{1/\tau},K_{2}^{1/\tau}\right\} ,
\end{align*}
as desired.
\end{proof}
One remaining limitation of Corollary \ref{cor:AtomicDecompositionSimplifiedCriteria}
is the assumption $\gamma_{i}=\gamma_{i,1}\ast\gamma_{i,2}$ with
certain assumptions on $\gamma_{i,1}$ and $\gamma_{i,2}$. For a
given function $\gamma$ (or $\gamma_{i}$), it can be cumbersome
to verify that it can be factorized as the convolution product of
two such functions.

Hence, we close this section by providing more readily verifiable
criteria which ensure that such a factorization is possible. For reasons
that will become clear later, we begin with the following technical
result:
\begin{lem}
\label{lem:ChineseBracketDerivative}For $\xi\in\R^{\dimension}$,
set $\left\{ \xi\right\} :=1+\left|\xi\right|^{2}$. Then, for each
$\theta\in\R$ and each $\alpha\in\N_{0}^{\dimension}$, there is
a polynomial $P_{\theta,\alpha}\in\R\left[\xi_{1},\dots,\xi_{\dimension}\right]$
of degree $\deg P_{\theta,\alpha}\leq\left|\alpha\right|$ satisfying
\[
\partial^{\alpha}\left\{ \xi\right\} ^{\theta}=\left\{ \xi\right\} ^{\theta-\left|\alpha\right|}\cdot P_{\theta,\alpha}\left(\xi\right)\qquad\forall\xi\in\R^{\dimension},
\]
as well as $\left|P_{\theta,\alpha}\left(\xi\right)\right|\leq C\cdot\left(1+\left|\xi\right|\right)^{\left|\alpha\right|}$
for all $\xi\in\R^{\dimension}$, where
\[
C=\left|\alpha\right|!\cdot\left[2\cdot\left(1+\dimension+\left|\theta\right|\right)\right]^{\left|\alpha\right|}.
\]

In particular, we have
\[
\left|\partial^{\alpha}\left\{ \xi\right\} ^{\theta}\right|\leq2^{\left|\theta\right|+\left|\alpha\right|}C\cdot\left(1+\left|\xi\right|\right)^{2\theta-\left|\alpha\right|}\leq2^{\left|\theta\right|+\left|\alpha\right|}C\cdot\left(1+\left|\xi\right|\right)^{2\theta}\qquad\forall\xi\in\R^{\dimension}.\qedhere
\]
\end{lem}
\begin{proof}
We prove existence of the polynomial $P_{\theta,\alpha}$ by induction
on $N=\left|\alpha\right|\in\N_{0}$. But to do this, we need a slightly
different formulation of the claim: For $P=\sum_{\sigma\in\N_{0}^{\dimension}}c_{\sigma}\xi^{\sigma}\in\R\left[\xi_{1},\dots,\xi_{\dimension}\right]$,
we define $\left\Vert P\right\Vert _{\ast}:=\sum_{\sigma\in\N_{0}^{\dimension}}\left|c_{\sigma}\right|$.
Below, we will prove by induction on $N=\left|\alpha\right|\in\N_{0}$
that the polynomial $P_{\theta,\alpha}$ satisfying $\partial^{\alpha}\left\{ \xi\right\} ^{\theta}$=$\left\{ \xi\right\} ^{\theta-\left|\alpha\right|}\cdot P_{\theta,\alpha}\left(\xi\right)$
can be chosen to satisfy $\left\Vert P_{\theta,\alpha}\right\Vert _{\ast}\leq C$
with $C$ as in the statement of the lemma. This will imply the claim,
since, for suitable coefficients $c_{\sigma}=c_{\sigma}\left(P_{\theta,\alpha}\right)$,
we have
\begin{align*}
\left|P_{\theta,\alpha}\left(\xi\right)\right| & \leq\sum_{\left|\sigma\right|\leq\left|\alpha\right|}\left|c_{\sigma}\right|\cdot\left|\xi^{\sigma}\right|\leq\sum_{\left|\sigma\right|\leq\left|\alpha\right|}\left|c_{\sigma}\right|\cdot\left(1+\left|\xi\right|\right)^{\left|\sigma\right|}\\
 & \leq\left(1+\left|\xi\right|\right)^{\left|\alpha\right|}\cdot\left\Vert P_{\theta,\alpha}\right\Vert _{\ast}\leq C\cdot\left(1+\left|\xi\right|\right)^{\left|\alpha\right|}
\end{align*}
for all $\xi\in\R^{\dimension}$.

It remains to prove the modified claim by induction on $N$. For $N=0$,
all properties are trivially satisfied for $P_{\theta,\alpha}\equiv1$.

For the induction step, we observe that each $\alpha\in\N_{0}^{\dimension}$
with $\left|\alpha\right|=N+1$ can be written as $\alpha=\beta+e_{j}$
for some $j\in\underline{\dimension}$, where $e_{j}$ is the $j$-th
standard basis vector and where $\beta\in\N_{0}^{\dimension}$ with
$\left|\beta\right|=N$. Now, a direct calculation yields
\begin{align*}
\partial^{\alpha}\left\{ \xi\right\} ^{\theta} & =\partial_{j}\partial^{\beta}\left\{ \xi\right\} ^{\theta}=\partial_{j}\left[\left\{ \xi\right\} ^{\theta-\left|\beta\right|}\cdot P_{\theta,\beta}\left(\xi\right)\right]\\
 & =\left\{ \xi\right\} ^{\theta-\left|\beta\right|}\cdot\partial_{j}P_{\theta,\beta}\left(\xi\right)+P_{\theta,\beta}\left(\xi\right)\cdot\left(\theta-\left|\beta\right|\right)\cdot\left\{ \xi\right\} ^{\theta-\left|\beta\right|-1}\cdot\partial_{j}\left\{ \xi\right\} \\
 & =\left\{ \xi\right\} ^{\theta-\left|\alpha\right|}\left[\left\{ \xi\right\} \cdot\partial_{j}P_{\theta,\beta}\left(\xi\right)+2\left(\theta-\left|\beta\right|\right)\cdot\xi_{j}\cdot P_{\theta,\beta}\left(\xi\right)\right]\\
 & =:\left\{ \xi\right\} ^{\theta-\left|\alpha\right|}\cdot P_{\theta,\alpha}\left(\xi\right).
\end{align*}
Since $\deg\left[\left\{ \xi\right\} \cdot\partial_{j}P_{\theta,\beta}\right]\leq2+\deg P_{\theta,\beta}-1\leq\left|\beta\right|+1=\left|\alpha\right|$,
it is not hard to see that $\deg P_{\theta,\alpha}\leq\left|\alpha\right|$.

Next, observe for $\sigma\in\N_{0}^{\dimension}$ with $\sigma_{j}\geq1$
that 
\[
\partial_{j}\xi^{\sigma}=\prod_{\ell\neq j}\xi_{\ell}^{\sigma_{\ell}}\cdot\partial_{j}\xi_{j}^{\sigma_{j}}=\sigma_{j}\cdot\xi^{\sigma-e_{j}}\qquad\forall\xi\in\R^{\dimension}.
\]
Furthermore, $\partial_{j}\xi^{\sigma}\equiv0$ in case of $\sigma_{j}=0$.
For $P_{\theta,\beta}\left(\xi\right)=\sum_{\left|\sigma\right|\leq\left|\beta\right|}c_{\sigma}\xi^{\sigma}$,
this implies
\[
\left\Vert \partial_{j}P_{\theta,\beta}\right\Vert _{\ast}\leq\sum_{\left|\sigma\right|\leq\left|\beta\right|}\left|c_{\sigma}\right|\cdot\left\Vert \partial_{j}\xi^{\sigma}\right\Vert _{\ast}\leq\sum_{\substack{\left|\sigma\right|\leq\left|\beta\right|\\
\sigma_{j}\geq1
}
}\sigma_{j}\cdot\left|c_{\sigma}\right|\leq\left|\beta\right|\cdot\left\Vert P_{\theta,\beta}\right\Vert _{\ast}.
\]
Furthermore, since $\left\Vert \xi^{\sigma}\cdot P\right\Vert _{\ast}=\left\Vert P\right\Vert _{\ast}$
for each polynomial $P\in\R\left[\xi_{1},\dots,\xi_{\dimension}\right]$
and each $\sigma\in\N_{0}^{\dimension}$ and since we have $\left\{ \xi\right\} =1+\sum_{\ell=1}^{\dimension}\xi^{2e_{\ell}}$,
we get
\begin{align*}
\left\Vert P_{\theta,\alpha}\right\Vert _{\ast} & \leq\left\Vert \partial_{j}P_{\theta,\beta}\right\Vert _{\ast}+\sum_{\ell=1}^{\dimension}\left\Vert \xi^{2e_{\ell}}\cdot\partial_{j}P_{\theta,\beta}\right\Vert _{\ast}+2\left|\theta-\left|\beta\right|\right|\cdot\left\Vert \xi_{j}\cdot P_{\theta,\beta}\right\Vert _{\ast}\\
 & \leq\left(1+\dimension\right)\cdot\left\Vert \partial_{j}P_{\theta,\beta}\right\Vert _{\ast}+2\left(\left|\theta\right|+\left|\beta\right|\right)\cdot\left\Vert P_{\theta,\beta}\right\Vert _{\ast}\\
 & \leq\left\Vert P_{\theta,\beta}\right\Vert _{\ast}\left[\left(1+\dimension\right)\cdot\left|\beta\right|+2\left(\left|\theta\right|+\left|\beta\right|\right)\right]\\
 & \leq\left\Vert P_{\theta,\beta}\right\Vert _{\ast}\left[\left(1+\dimension\right)\cdot\left(1+\left|\beta\right|\right)+2\left(1+\left|\theta\right|\right)\left(1+\left|\beta\right|\right)\right]\\
 & \leq\left|\alpha\right|\cdot\left[\left(1+\dimension\right)+2\left(1+\left|\theta\right|\right)\right]\cdot\left\Vert P_{\theta,\beta}\right\Vert _{\ast}\\
\left({\scriptstyle \text{since }\dimension\geq1}\right) & \leq\left|\alpha\right|\cdot2\left(1+\dimension+\left|\theta\right|\right)\cdot\left\Vert P_{\theta,\beta}\right\Vert _{\ast}.
\end{align*}
Since $\left\Vert P_{\theta,\beta}\right\Vert _{\ast}\leq\left|\beta\right|!\cdot\left[2\left(1+\dimension+\left|\theta\right|\right)\right]^{\left|\beta\right|}$
by induction and since $\left|\alpha\right|=\left|\beta\right|+1$,
the induction is complete.

It remains to verify the final statement of the lemma. To this end,
note
\[
\frac{1}{2}\left(1+\left|\xi\right|\right)^{2}\leq\left\{ \xi\right\} =1+\left|\xi\right|^{2}\leq\left(1+\left|\xi\right|\right)^{2}\leq2\cdot\left(1+\left|\xi\right|\right)^{2},
\]
so that $\left\{ \xi\right\} ^{\varrho}\leq2^{\left|\varrho\right|}\cdot\left(1+\left|\xi\right|\right)^{2\varrho}$
for all $\varrho\in\R$ and $\xi\in\R^{\dimension}$. Hence,
\begin{align*}
\left|\partial^{\alpha}\left\{ \xi\right\} ^{\theta}\right|=\left|\left\{ \xi\right\} ^{\theta-\left|\alpha\right|}\cdot P_{\theta,\alpha}\left(\xi\right)\right| & \leq2^{\left|\theta-\left|\alpha\right|\right|}\cdot\left(1+\left|\xi\right|\right)^{2\theta-2\left|\alpha\right|}\cdot\left|P_{\theta,\alpha}\left(\xi\right)\right|\\
 & \leq2^{\left|\theta\right|+\left|\alpha\right|}C\cdot\left(1+\left|\xi\right|\right)^{2\theta-\left|\alpha\right|},
\end{align*}
as claimed.
\end{proof}
\begin{lem}
\label{lem:ConvolutionFactorization}Let $\varrho\in L^{1}\left(\R^{\dimension}\right)$
with $\varrho\geq0$. Let $N\geq\dimension+1$ and assume that $\gamma\in L^{1}\left(\R^{\dimension}\right)$
satisfies $\widehat{\gamma}\in C^{N}\left(\R^{\dimension}\right)$
with
\[
\left|\partial^{\alpha}\widehat{\gamma}\left(\xi\right)\right|\leq\varrho\left(\xi\right)\cdot\left(1+\left|\xi\right|\right)^{-\left(\dimension+1+\varepsilon\right)}\qquad\forall\xi\in\R^{\dimension}\quad\forall\alpha\in\N_{0}^{\dimension}\text{ with }\left|\alpha\right|\leq N
\]
for some $\varepsilon\in\left(0,1\right]$.

Then there are functions $\gamma_{1}\in C_{0}\left(\R^{\dimension}\right)\cap L^{1}\left(\R^{\dimension}\right)$
and $\gamma_{2}\in C^{1}\left(\R^{\dimension}\right)\cap W^{1,1}\left(\R^{\dimension}\right)$
with $\gamma=\gamma_{1}\ast\gamma_{2}$ and with the following additional
properties:

\begin{enumerate}
\item We have $\left\Vert \gamma_{2}\right\Vert _{K}\leq s_{\dimension}\cdot2^{1+\dimension+3K}\cdot K!\cdot\left(1+\dimension\right)^{1+2K}$
and $\left\Vert \nabla\gamma_{2}\right\Vert _{K}\leq\frac{s_{\dimension}}{\varepsilon}\cdot2^{4+\dimension+3K}\cdot\left(1+\dimension\right)^{2\left(1+K\right)}\cdot\left(K+1\right)!$
for all $K\in\N_{0}$, where as usual $\left\Vert g\right\Vert _{K}:=\sup_{x\in\R^{\dimension}}\left(1+\left|x\right|\right)^{K}\left|g\left(x\right)\right|$.
\item We have $\widehat{\gamma_{2}}\in C^{\infty}\left(\R^{\dimension}\right)$
with all partial derivatives of $\widehat{\gamma_{2}}$ being polynomially
bounded (even bounded).
\item If $\widehat{\gamma}\in C^{\infty}\left(\R^{\dimension}\right)$ with
all partial derivatives being polynomially bounded, the same also
holds for $\widehat{\gamma_{1}}$.
\item We have $\left\Vert \gamma_{1}\right\Vert _{N}\leq\left(1+\dimension\right)^{1+2N}\cdot2^{1+\dimension+4N}\cdot N!\cdot\left\Vert \varrho\right\Vert _{L^{1}}$
and $\left\Vert \gamma\right\Vert _{N}\leq\left(1+\dimension\right)^{N+1}\cdot\left\Vert \varrho\right\Vert _{L^{1}}$.
\item We have $\left|\partial^{\alpha}\widehat{\gamma_{1}}\left(\xi\right)\right|\leq2^{1+\dimension+4N}\cdot N!\cdot\left(1+\dimension\right)^{N}\cdot\varrho\left(\xi\right)$
for all $\xi\in\R^{\dimension}$ and $\alpha\in\N_{0}^{\dimension}$
with $\left|\alpha\right|\leq N$.\qedhere
\end{enumerate}
\end{lem}
\begin{rem*}
For concrete special cases of $\CalQ$, this lemma will be applied
as follows: In most cases, one can find a suitable function $\varrho$
as above such that property (\ref{enu:AtomicDecompositionSimplifiedMainCondition})
of Corollary \ref{cor:AtomicDecompositionSimplifiedCriteria} is satisfied
as soon as all $\gamma_{j,1}\in L^{1}\left(\R^{\dimension}\right)$
satisfy $\left|\partial^{\alpha}\widehat{\gamma_{j,1}}\left(\xi\right)\right|\lesssim\varrho\left(\xi\right)$
uniformly in $\left|\alpha\right|\leq N$, $j\in I$ and $\xi\in\R^{\dimension}$.

In this case, the preceding lemma shows that if we instead assume
$\left|\partial^{\alpha}\widehat{\gamma_{j}}\left(\xi\right)\right|\leq\varrho\left(\xi\right)\cdot\left(1+\left|\xi\right|\right)^{-\left(\dimension+1+\varepsilon\right)}$
for all $\alpha,j,\xi$ as above, then we can write $\gamma_{j}=\gamma_{j,1}\ast\gamma_{j,2}$
such that $\left|\partial^{\alpha}\widehat{\gamma_{j,1}}\left(\xi\right)\right|\lesssim\varrho\left(\xi\right)$
uniformly in $\alpha,j,\xi$ as above, so that the family $\left(\gamma_{j,1}\right)_{j\in I}$
satisfies property (\ref{enu:AtomicDecompositionSimplifiedMainCondition}).
Furthermore, the family $\left(\gamma_{j,2}\right)_{j\in I}$ satisfies
all assumptions of Corollary \ref{cor:AtomicDecompositionSimplifiedCriteria}.

By possibly enlarging $N$ for the application of the lemma, it is
then not hard to ensure that all prerequisites of Corollary \ref{cor:AtomicDecompositionSimplifiedCriteria}
are satisfied. Hence, Lemma \ref{lem:ConvolutionFactorization} essentially
solves the \emph{factorization problem} mentioned before Lemma \ref{lem:ChineseBracketDerivative}.
\end{rem*}
\begin{proof}
We will use the notation $\left\{ \xi\right\} =1+\left|\xi\right|^{2}$
from Lemma \ref{lem:ChineseBracketDerivative}, as well as $\left\langle \xi\right\rangle :=\left\{ \xi\right\} ^{1/2}$.
With this notation, define $g\in C^{\infty}\left(\R^{\dimension}\right)$
by $g:\R^{\dimension}\to\left(0,\infty\right),\xi\mapsto\left\{ \xi\right\} ^{-\frac{\dimension+1+\varepsilon}{2}}=\left\langle \xi\right\rangle ^{-\left(\dimension+1+\varepsilon\right)}$.
In view of equation (\ref{eq:StandardDecayLpEstimate}) and since
$\left\langle \xi\right\rangle \geq\frac{1}{2}\left(1+\left|\xi\right|\right)$,
it is not hard to see $g\in L^{1}\left(\R^{\dimension}\right)$, so
that $\gamma_{2}:=\Fourier^{-1}g\in C_{0}\left(\R^{\dimension}\right)$
is well-defined.

Next, let $K\in\N_{0}$ be arbitrary. For $\alpha\in\N_{0}^{\dimension}$
with $\left|\alpha\right|\leq K$, Lemma \ref{lem:ChineseBracketDerivative}
(with $\theta=-\frac{\dimension+1+\varepsilon}{2}$) shows 
\begin{equation}
\left|\partial^{\alpha}g\left(\xi\right)\right|\leq2^{1+\dimension+K}\cdot K!\cdot\left[4\cdot\left(1+\dimension\right)\right]^{K}\cdot\left(1+\left|\xi\right|\right)^{-\left(\dimension+1+\varepsilon\right)}=:C_{\dimension,K}\cdot\left(1+\left|\xi\right|\right)^{-\left(\dimension+1+\varepsilon\right)}\qquad\forall\xi\in\R^{\dimension}.\label{eq:ChineseBracketDerivativeSelfBound}
\end{equation}
In particular, this implies $g\in W^{K,1}\left(\R^{\dimension}\right)$.
In view of Lemma \ref{lem:PointwiseFourierDecayEstimate}, we thus
get
\begin{align*}
\left(1+\left|x\right|\right)^{K}\cdot\left|\gamma_{2}\left(x\right)\right|=\left(1+\left|x\right|\right)^{K}\cdot\left|\Fourier^{-1}g\left(x\right)\right| & \leq\left(1+\dimension\right)^{K}\cdot\left(\left|\left(\Fourier^{-1}g\right)\left(x\right)\right|+\sum_{m=1}^{\dimension}\left|\left[\Fourier^{-1}\left(\partial_{m}^{K}g\right)\right]\left(x\right)\right|\right)\\
\left({\scriptstyle \text{since }\left|\Fourier^{-1}h\right|\leq\left\Vert h\right\Vert _{L^{1}}}\right) & \leq\left(1+\dimension\right)^{K}\cdot C_{\dimension,K}\cdot\left\Vert \left(1+\left|\mybullet\right|\right)^{-\left(\dimension+1\right)}\right\Vert _{L^{1}}\cdot\left(1+\dimension\right)\\
\left({\scriptstyle \text{eq. }\eqref{eq:StandardDecayLpEstimate}}\right) & \leq s_{\dimension}\cdot2^{1+\dimension+3K}\cdot K!\cdot\left(1+\dimension\right)^{1+2K}
\end{align*}
for all $x\in\R^{\dimension}$, and thus $\left\Vert \gamma_{2}\right\Vert _{K}\leq s_{\dimension}\cdot2^{1+\dimension+3K}\cdot K!\cdot\left(1+\dimension\right)^{1+2K}<\infty$
for arbitrary $K\in\N_{0}$, as desired. In particular, $\gamma_{2}\in L^{1}\left(\R^{\dimension}\right)$,
so that $\widehat{\gamma_{2}}=\Fourier\Fourier^{-1}g=g$ by Fourier
inversion. Hence, equation (\ref{eq:ChineseBracketDerivativeSelfBound})
shows that all partial derivatives of $\widehat{\gamma_{2}}=g$ are
bounded.

Next, we want to estimate $\left\Vert \nabla\gamma_{2}\right\Vert _{K}$.
To this end, we observe for arbitrary $j\in\underline{\dimension}$
that 
\[
\left|\xi_{j}\cdot g\left(\xi\right)\right|\leq\left\{ \xi\right\} ^{1/2}\cdot\left|g\left(\xi\right)\right|=\left\{ \xi\right\} ^{-\frac{\dimension+\varepsilon}{2}}=\left\langle \xi\right\rangle ^{-\left(\dimension+\varepsilon\right)}\in L^{1}\left(\R^{\dimension}\right),
\]
so that we can differentiate under the integral in the definition
of $\gamma_{2}\left(x\right)=\left(\Fourier^{-1}g\right)\left(x\right)$
to conclude for $g_{j}\left(\xi\right):=\xi_{j}\cdot g\left(\xi\right)$
that $\gamma_{2}\in C^{1}\left(\R^{\dimension}\right)$ with derivative
\[
\partial_{j}\gamma_{2}\left(x\right)=2\pi i\cdot\int_{\R^{\dimension}}\xi_{j}\cdot g\left(\xi\right)\cdot e^{2\pi i\left\langle x,\xi\right\rangle }\d\xi=2\pi i\cdot\left(\Fourier^{-1}g_{j}\right)\left(x\right).
\]

Now, we want to apply Lemma \ref{lem:PointwiseFourierDecayEstimate}
again to derive a bound for $\partial_{j}\gamma_{2}\left(x\right)$,
which requires us to bound the derivatives $\partial^{\alpha}g_{j}$.
To this end, we observe $\partial^{\alpha}\xi_{j}\equiv0$ unless
$\alpha=0$, in which case we have $\partial^{\alpha}\xi_{j}=\xi_{j}$,
or unless $\alpha=e_{j}$, in which case we have $\partial^{\alpha}\xi_{j}=1$.
In combination with Leibniz's rule, this yields for $\left|\alpha\right|\leq K$
the estimate
\begin{align*}
\left|\partial^{\alpha}g_{j}\left(\xi\right)\right| & \leq\sum_{\beta\leq\alpha}\binom{\alpha}{\beta}\left|\partial^{\beta}\xi_{j}\right|\left|\partial^{\alpha-\beta}g\left(\xi\right)\right|\\
 & =\begin{cases}
\left|\xi_{j}\right|\cdot\left|\partial^{\alpha}g\left(\xi\right)\right|, & \text{if }\alpha_{j}=0,\\
\left|\xi_{j}\right|\cdot\left|\partial^{\alpha}g\left(\xi\right)\right|+\binom{\alpha}{e_{j}}\cdot\left|\partial^{\alpha-e_{j}}g\left(\xi\right)\right|, & \text{if }\alpha_{j}\geq1
\end{cases}\\
\left({\scriptstyle \text{eq. }\eqref{eq:ChineseBracketDerivativeSelfBound}}\right) & \leq\begin{cases}
C_{\dimension,K}\cdot\left(1+\left|\xi\right|\right)^{-\left(\dimension+1+\varepsilon\right)}\cdot\left|\xi_{j}\right|, & \text{if }\alpha_{j}=0,\\
C_{\dimension,K}\cdot\left(\left|\xi_{j}\right|+\binom{\alpha}{e_{j}}\right)\cdot\left(1+\left|\xi\right|\right)^{-\left(\dimension+1+\varepsilon\right)}, & \text{if }\alpha_{j}\geq1
\end{cases}\\
\left({\scriptstyle \text{since }\binom{\alpha}{e_{j}}=\binom{\alpha_{j}}{1}=\alpha_{j}\leq\left|\alpha\right|\leq K}\right) & \leq C_{\dimension,K}\cdot\left(1+K\right)\cdot\left(1+\left|\xi\right|\right)^{-\left(\dimension+\varepsilon\right)}.
\end{align*}
In particular, this implies $g_{j}\in W^{K,1}\left(\R^{\dimension}\right)$.
Hence, another application of Lemma \ref{lem:PointwiseFourierDecayEstimate}
and equation (\ref{eq:StandardDecayLpEstimate}) yields
\begin{align*}
\left(1+\left|x\right|\right)^{K}\cdot\left|\partial_{j}\gamma_{2}\left(x\right)\right| & =2\pi\cdot\left(1+\left|x\right|\right)^{K}\cdot\left|\Fourier^{-1}g_{j}\left(x\right)\right|\\
\left({\scriptstyle \text{eq. }\eqref{eq:PointwiseFourierDecayEstimate}}\right) & \leq8\cdot\left(1+\dimension\right)^{K}\cdot\left(\left|\left(\Fourier^{-1}g_{j}\right)\left(x\right)\right|+\sum_{m=1}^{\dimension}\left|\left[\Fourier^{-1}\left(\partial_{m}^{K}g_{j}\right)\right]\left(x\right)\right|\right)\\
 & \leq8\cdot\left(1+\dimension\right)^{K}\cdot C_{\dimension,K}\cdot\left(1+K\right)\cdot\left(1+\dimension\right)\cdot\left\Vert \left(1+\left|\mybullet\right|\right)^{-\left(\dimension+\varepsilon\right)}\right\Vert _{L^{1}}\\
\left({\scriptstyle \text{eq. }\eqref{eq:StandardDecayLpEstimate}}\right) & \leq\frac{s_{\dimension}}{\varepsilon}\cdot2^{4+\dimension+3K}\cdot\left(1+\dimension\right)^{1+2K}\cdot\left(K+1\right)!
\end{align*}
and hence $\left\Vert \nabla\gamma_{2}\right\Vert _{K}\leq\frac{s_{\dimension}}{\varepsilon}\cdot2^{4+\dimension+3K}\cdot\left(1+\dimension\right)^{2\left(1+K\right)}\cdot\left(K+1\right)!$
for arbitrary $K\in\N_{0}$, as claimed. In particular, $\nabla\gamma_{2}\in L^{1}\left(\R^{\dimension}\right)$,
so that $\gamma_{2}\in C^{1}\left(\R^{\dimension}\right)\cap W^{1,1}\left(\R^{\dimension}\right)$,
as claimed.

It remains to construct $\gamma_{1}$ with the desired properties.
To this end, define $h:\R^{\dimension}\to\Compl,\xi\mapsto\widehat{\gamma}\left(\xi\right)\cdot\left\{ \xi\right\} ^{\frac{\dimension+1+\varepsilon}{2}}$,
note $h\in C^{N}\left(\R^{\dimension}\right)$  and observe that Lemma
\ref{lem:ChineseBracketDerivative} shows for arbitrary $\beta\in\N_{0}^{\dimension}$
with $\left|\beta\right|\leq N$ that
\begin{equation}
\left|\partial^{\beta}\left\{ \xi\right\} ^{\frac{\dimension+1+\varepsilon}{2}}\right|\leq C_{\dimension,N}\cdot\left(1+\left|\xi\right|\right)^{\dimension+1+\varepsilon}\label{eq:FactorizationLemmaChineseBracketPositiveExponentDerivative}
\end{equation}
with the same constant $C_{\dimension,N}$ (with $N=K$) as in equation
(\ref{eq:ChineseBracketDerivativeSelfBound}). In combination with
Leibniz's rule and the $\dimension$-dimensional binomial theorem
(cf.\@ \cite[Section 8.1, Exercise 2.b]{FollandRA}), this yields
for arbitrary $\alpha\in\N_{0}^{\dimension}$ with $\left|\alpha\right|\leq N$
that
\begin{align}
\left|\partial^{\alpha}h\left(\xi\right)\right| & \leq\sum_{\beta\leq\alpha}\binom{\alpha}{\beta}\cdot\left|\partial^{\beta}\left\{ \xi\right\} ^{\frac{\dimension+1+\varepsilon}{2}}\right|\cdot\left|\partial^{\alpha-\beta}\widehat{\gamma}\left(\xi\right)\right|\nonumber \\
\left({\scriptstyle \text{assump. on }\widehat{\gamma}\text{ and eq. \eqref{eq:FactorizationLemmaChineseBracketPositiveExponentDerivative}}}\right) & \leq\varrho\left(\xi\right)\cdot\left(1+\left|\xi\right|\right)^{-\left(\dimension+1+\varepsilon\right)}C_{\dimension,N}\left(1+\left|\xi\right|\right)^{\dimension+1+\varepsilon}\cdot\sum_{\beta\leq\alpha}\binom{\alpha}{\beta}\nonumber \\
 & \leq2^{N}C_{\dimension,N}\cdot\varrho\left(\xi\right)\in L^{1}\left(\smash{\R^{\dimension}}\right).\label{eq:FactorizationLemmaFirstFactorFourierEstimate}
\end{align}
In particular, $h\in L^{1}\left(\R^{\dimension}\right)$, so that
$\gamma_{1}:=\Fourier^{-1}h\in C_{0}\left(\R^{\dimension}\right)$
is well-defined. Furthermore, we get $h\in W^{N,1}\left(\R^{\dimension}\right)$.

Hence, we can invoke Lemma \ref{lem:PointwiseFourierDecayEstimate}
once again to derive
\begin{align*}
\left(1+\left|x\right|\right)^{N}\cdot\left|\gamma_{1}\left(x\right)\right| & =\left(1+\left|x\right|\right)^{N}\cdot\left|\left(\Fourier^{-1}h\right)\left(x\right)\right|\\
\left({\scriptstyle \text{eq. }\eqref{eq:PointwiseFourierDecayEstimate}}\right) & \leq\left(1+\dimension\right)^{N}\cdot\left(\left|\Fourier^{-1}h\left(x\right)\right|+\sum_{m=1}^{\dimension}\left|\left[\Fourier^{-1}\left(\partial_{m}^{N}h\right)\right]\left(x\right)\right|\right)\\
 & \leq\left(1+\dimension\right)^{N+1}\cdot2^{N}C_{\dimension,N}\cdot\left\Vert \varrho\right\Vert _{L^{1}}\\
 & \leq\left(1+\dimension\right)^{1+2N}\cdot2^{1+\dimension+4N}\cdot N!\cdot\left\Vert \varrho\right\Vert _{L^{1}}<\infty
\end{align*}
for all $x\in\R^{\dimension}$, so that $\left\Vert \gamma_{1}\right\Vert _{N}\leq\left(1+\dimension\right)^{1+2N}\cdot2^{1+\dimension+4N}\cdot N!\cdot\left\Vert \varrho\right\Vert _{L^{1}}$,
as claimed. Since $N\geq\dimension+1$ by assumption, equation (\ref{eq:StandardDecayLpEstimate})
implies in particular that $\gamma_{1}\in L^{1}\left(\R^{\dimension}\right)$.
Hence, Fourier inversion yields $\widehat{\gamma_{1}}=\Fourier\Fourier^{-1}h=h$,
so that equation (\ref{eq:FactorizationLemmaFirstFactorFourierEstimate})
yields the claimed estimate for $\left|\partial^{\alpha}\widehat{\gamma_{1}}\right|$.

Next, the convolution theorem yields
\[
\Fourier\left[\gamma_{1}\ast\gamma_{2}\right]\left(\xi\right)=\widehat{\gamma_{1}}\left(\xi\right)\cdot\widehat{\gamma_{2}}\left(\xi\right)=h\left(\xi\right)\cdot g\left(\xi\right)=\widehat{\gamma}\left(\xi\right)\qquad\forall\xi\in\R^{\dimension}
\]
and thus $\gamma=\gamma_{1}\ast\gamma_{2}$, by injectivity of the
Fourier transform.

Next, note that if $\widehat{\gamma}\in C^{\infty}\left(\R^{\dimension}\right)$
with all derivatives being polynomially bounded, we clearly get $\widehat{\gamma_{1}}=h\in C^{\infty}\left(\R^{\dimension}\right)$,
again with all derivatives being polynomially bounded, thanks to Lemma
\ref{lem:ChineseBracketDerivative} and the Leibniz rule.

It remains to establish the estimate for $\left\Vert \gamma\right\Vert _{N}$.
But since $\left|\widehat{\gamma}\right|\leq\varrho\in L^{1}\left(\R^{\dimension}\right)$,
we get $\gamma=\Fourier^{-1}\widehat{\gamma}$ by Fourier inversion.
Furthermore, our assumptions easily yield $\partial^{\alpha}\widehat{\gamma}\in L^{1}\left(\R^{\dimension}\right)$
for all $\alpha\in\N_{0}^{\dimension}$ with $\left|\alpha\right|\leq N$.
Hence, a final application of Lemma \ref{lem:PointwiseFourierDecayEstimate},
together with our assumptions on $\widehat{\gamma}$, yields
\begin{align*}
\left(1+\left|x\right|\right)^{N}\cdot\left|\gamma\left(x\right)\right| & =\left(1+\left|x\right|\right)^{N}\cdot\left|\left(\Fourier^{-1}\widehat{\gamma}\right)\left(x\right)\right|\\
 & \leq\left(1+\dimension\right)^{N}\cdot\left(\left|\left(\Fourier^{-1}\widehat{\gamma}\right)\left(x\right)\right|+\sum_{m=1}^{\dimension}\left|\left[\Fourier^{-1}\left(\partial_{m}^{N}\widehat{\gamma}\right)\right]\left(x\right)\right|\right)\\
 & \leq\left(1+\dimension\right)^{N}\cdot\left(\left\Vert \widehat{\gamma}\right\Vert _{L^{1}}+\sum_{m=1}^{\dimension}\left\Vert \partial_{m}^{N}\widehat{\gamma}\right\Vert _{L^{1}}\right)\\
 & \leq\left(1+\dimension\right)^{N+1}\cdot\left\Vert \varrho\right\Vert _{L^{1}}<\infty,
\end{align*}
which easily yields the claim.
\end{proof}

\section{Existence of compactly supported Banach frames and atomic decompositions
for \texorpdfstring{$\alpha$}{α}-modulation spaces}

\label{sec:CompactlySupportedAlphaModulationFrames}In this section,
we show that the general theory developed in this paper can be used
to prove existence of compactly supported, structured Banach frames
for $\alpha$-modulation spaces. A brief discussion of the history
and the applications of $\alpha$-modulation spaces, as well as a
comparison of our results with the established literature will be
given at the end of the section.

We begin our considerations by recalling the definition of $\alpha$-modulation
spaces, as given by Borup and Nielsen\cite{BorupNielsenAlphaModulationSpaces}.
First of all, we have to define the associated covering. Its admissibility
was established in \cite[Theorem 2.6]{BorupNielsenAlphaModulationSpaces};
precisely, the following was shown:
\begin{thm}
\label{thm:AlphaModulationCoveringDefinition}(cf. \cite[Theorem 2.6]{BorupNielsenAlphaModulationSpaces})
Let $\dimension\in\N$ and $\alpha\in\left[0,1\right)$ be arbitrary.
Define $\alpha_{0}:=\frac{\alpha}{1-\alpha}$. Then there is a constant
$r_{1}=r_{1}\left(\dimension,\alpha\right)$ such that the family
\begin{equation}
\CalQ^{\left(\alpha\right)}:=\CalQ_{r}^{\left(\alpha\right)}:=\left(Q_{r,k}^{\left(\alpha\right)}\right)_{k\in\Z^{\dimension}\setminus\left\{ 0\right\} }:=\left(B_{r\cdot\left|k\right|^{\alpha_{0}}}\left(\left|k\right|^{\alpha_{0}}k\right)\right)_{k\in\Z^{\dimension}\setminus\left\{ 0\right\} }\label{eq:AlphaModulationCovering}
\end{equation}
is an admissible covering of $\R^{\dimension}$ for every $r>r_{1}$.
The covering $\CalQ_{r}^{\left(\alpha\right)}$ is called the \textbf{$\alpha$-modulation
covering} of $\R^{\dimension}$. If the values of $r$ and $\alpha$
are clear from the context, we also write $Q_{k}:=Q_{r,k}^{\left(\alpha\right)}$.
\end{thm}
The associated weight is defined in our next lemma. There, and in
the remainder of this section, we use the notation $\left\langle \xi\right\rangle :=\left(1+\smash{\left|\xi\right|^{2}}\right)^{1/2}$
for $\xi\in\R^{\dimension}$.
\begin{lem}
\label{lem:AlphaModulationFrequencyWeight}(cf. \cite[Lemma 9.2]{DecompositionEmbedding})
Let $\dimension\in\N$ and $\alpha\in\left[0,1\right)$ and let $r>0$
be chosen such that the $\alpha$-modulation covering $\CalQ_{r}^{\left(\alpha\right)}=\left(Q_{r,k}^{\left(\alpha\right)}\right)_{k\in\Z^{\dimension}\setminus\left\{ 0\right\} }$
is an admissible covering of $\R^{\dimension}$. We then have
\[
\left\langle \xi\right\rangle \asymp\left\langle k\right\rangle ^{\frac{1}{1-\alpha}}\qquad\text{ for all }\qquad k\in\Z^{\dimension}\setminus\left\{ 0\right\} \text{ and }\xi\in Q_{r,k}^{\left(\alpha\right)},
\]
where the implied constant only depends on $r,\alpha$.

Now, for $s\in\R$, we define the weight $w^{\left(s\right)}$ on
$\Z^{\dimension}\setminus\left\{ 0\right\} $ by
\[
w_{k}^{\left(s\right)}:=\left\langle k\right\rangle ^{s}\qquad\text{ for }k\in\Z^{\dimension}\setminus\left\{ 0\right\} .
\]
Then $w^{\left(s\right)}$ is $\CalQ_{r}^{\left(\alpha\right)}$-moderate
(cf.\@ equation (\ref{eq:IntroductionModerateWeightDefinition})).
\end{lem}
Note that Theorem \ref{thm:AlphaModulationCoveringDefinition} only
claims that $\CalQ_{r}^{\left(\alpha\right)}$ is an \emph{admissible}
covering of $\R^{\dimension}$. The next result shows that it is actually
a \emph{structured} admissible covering of $\R^{\dimension}$ (cf.\@
the remark after Assumption \ref{assu:RegularPartitionOfUnity}) and
thus in particular a semi-structured covering.
\begin{lem}
\label{lem:AlphaModulationStructuredAndBAPU}Let $\dimension\in\N$
and $\alpha\in\left[0,1\right)$ and let $r_{1}=r_{1}\left(\dimension,\alpha\right)$
be as in Theorem \ref{thm:AlphaModulationCoveringDefinition} and
let $r>r_{1}$. For $k\in\Z^{\dimension}\setminus\left\{ 0\right\} $,
set $T_{k}:=\left|k\right|^{\alpha_{0}}\cdot\identity$ and $b_{k}:=\left|k\right|^{\alpha_{0}}k$
and let $Q:=B_{r}\left(0\right)$. Then we have
\[
\CalQ_{r}^{\left(\alpha\right)}=\left(T_{k}Q+b_{k}\right)_{k\in\Z^{\dimension}\setminus\left\{ 0\right\} }
\]
and with these choices, $\CalQ_{r}^{\left(\alpha\right)}$ is a semi-structured
admissible covering of $\R^{\dimension}$.

Finally, $\CalQ_{r}^{\left(\alpha\right)}$ admits a regular partition
of unity $\left(\varphi_{k}\right)_{k\in\Z^{\dimension}\setminus\left\{ 0\right\} }$
(which thus fulfills Assumption \ref{assu:RegularPartitionOfUnity})
and $\CalQ_{r}^{\left(\alpha\right)}$ fulfills the standing assumptions
from Section \ref{subsec:DecompSpaceDefinitionStandingAssumptions};
in particular, $\left\Vert T_{k}^{-1}\right\Vert \leq1=:\Omega_{0}$
for all $k\in\Z^{\dimension}\setminus\left\{ 0\right\} $.
\end{lem}
\begin{proof}
The fact that $\CalQ_{r}^{\left(\alpha\right)}$ is a structured admissible
covering of $\R^{\dimension}$ for $r>r_{1}$ was shown in \cite[Lemma 9.3]{DecompositionEmbedding}.
Since this is the case, \cite[Theorem 2.8]{DecompositionIntoSobolev}
shows that there is a regular partition of unity $\Phi=\left(\varphi_{k}\right)_{k\in\Z^{\dimension}\setminus\left\{ 0\right\} }$
for $\CalQ_{r}^{\left(\alpha\right)}$. In view of Corollary \ref{cor:RegularBAPUsAreWeightedBAPUs},
$\Phi$ is thus also a $\CalQ_{r}^{\left(\alpha\right)}$-$v_{0}$-BAPU
for every weight $v_{0}$ satisfying the general assumptions of Section
\ref{subsec:DecompSpaceDefinitionStandingAssumptions}.

Finally, we clearly have $\left\Vert T_{k}^{-1}\right\Vert =\left\Vert \left|k\right|^{-\alpha_{0}}\identity\right\Vert =\left|k\right|^{-\alpha_{0}}\leq1$
for all $k\in\Z^{\dimension}\setminus\left\{ 0\right\} $, since $\left|k\right|\geq1$
and $\alpha_{0}=\frac{\alpha}{1-\alpha}\geq0$.
\end{proof}
Finally, we introduce the weights $v=v^{\left(\mu\right)}$ that we
will use for the weighted $L^{p}$ spaces $L_{v}^{p}\left(\R^{\dimension}\right)$.
\begin{lem}
\label{lem:AlphaModulationSpaceWeight}For $\mu\in\R$ let
\begin{align*}
v^{\left(\mu\right)}: & \R^{\dimension}\to\left(0,\infty\right),x\mapsto\left\langle x\right\rangle ^{\mu}=\left(1+\left|x\right|^{2}\right)^{\mu/2},\\
v_{0}: & \R^{\dimension}\to\left(0,\infty\right),x\mapsto\left[2\cdot\left(1+\left|x\right|\right)\right]^{\left|\mu\right|}
\end{align*}
and set $K:=\left|\mu\right|$ and $\Omega_{1}:=2^{\left|\mu\right|}$.
With these choices, $v=v^{\left(\mu\right)}$ satisfies the standing
assumptions of Section \ref{subsec:DecompSpaceDefinitionStandingAssumptions}.
\end{lem}
\begin{proof}
First of all, assume $\mu=1$. In this case, we get
\begin{align*}
v^{\left(\mu\right)}\left(x+y\right) & =\left|\left(\begin{matrix}1\\
x+y
\end{matrix}\right)\right|\leq1+\left|x+y\right|\\
 & \leq1+\left|x\right|+\left|y\right|\leq\left(1+\left|x\right|\right)\left(1+\left|y\right|\right)\\
 & \leq2\cdot\left|\left(\begin{matrix}1\\
x
\end{matrix}\right)\right|\cdot\left(1+\left|y\right|\right)=v^{\left(\mu\right)}\left(x\right)\cdot v_{0}\left(y\right),
\end{align*}
where the last step used $\mu=1$. Now, for arbitrary $\mu\geq0$,
we likewise get $v^{\left(\mu\right)}\left(x+y\right)\leq v^{\left(\mu\right)}\left(x\right)\cdot v_{0}\left(y\right)$
by taking the $\mu$-th power of the preceding estimate.

Finally, if $\mu<0$, we have
\[
v^{\left(-\mu\right)}\left(x\right)=v^{\left(-\mu\right)}\left(x+y+\left(-y\right)\right)\leq v^{\left(-\mu\right)}\left(x+y\right)\cdot\left[2\cdot\left(1+\left|-y\right|\right)\right]^{\left|-\mu\right|}.
\]
Rearranging yields
\[
v^{\left(\mu\right)}\left(x+y\right)=\left[v^{\left(-\mu\right)}\left(x+y\right)\right]^{-1}\leq\left[v^{\left(-\mu\right)}\left(x\right)\right]^{-1}\cdot\left[2\cdot\left(1+\left|y\right|\right)\right]^{\left|\mu\right|}=v^{\left(\mu\right)}\left(x\right)\cdot v_{0}\left(y\right).
\]
Hence, we have shown for all $\mu\in\R$ that $v^{\left(\mu\right)}$
is $v_{0}$-moderate.

It is clear that $v_{0}\geq1$ and that $v_{0}$ is symmetric. Furthermore,
$v_{0}\left(x\right)=2^{\left|\mu\right|}\cdot\left(1+\left|x\right|\right)^{\left|\mu\right|}=\Omega_{1}\cdot\left(1+\left|x\right|\right)^{K}$
for all $x\in\R^{\dimension}$. We do not necessarily have $K=0$,
but in Lemma \ref{lem:AlphaModulationStructuredAndBAPU}, we already
saw $\left\Vert T_{k}^{-1}\right\Vert \leq1=\Omega_{0}$ for all $k\in\Z^{\dimension}\setminus\left\{ 0\right\} $.

The only thing which remains to be verified is that $v_{0}$ is submultiplicative.
But we have
\[
2\cdot\left(1+\left|x+y\right|\right)\leq2\cdot\left(1+\left|x\right|+\left|y\right|\right)\leq2\cdot\left(1+\left|x\right|\right)\left(1+\left|y\right|\right)\leq2\cdot\left(1+\left|x\right|\right)\cdot2\cdot\left(1+\left|y\right|\right).
\]
Taking the $\left|\mu\right|$-th power of this estimate yields $v_{0}\left(x+y\right)\leq v_{0}\left(x\right)\cdot v_{0}\left(y\right)$,
as desired.
\end{proof}
Having verified all these assumptions, we conclude from Proposition
\ref{prop:WeightedDecompositionSpaceWellDefined} and Lemma \ref{lem:WeightedDecompositionSpaceComplete}
that the $\alpha$-modulation spaces defined below are indeed well-defined
Quasi-Banach spaces.
\begin{defn}
\label{def:AlphaModulationSpaces}For $\dimension\in\N$ and $\alpha\in\left[0,1\right)$,
choose some $r>r_{1}\left(\dimension,\alpha\right)$ with $r_{1}\left(\dimension,\alpha\right)$
as in Theorem \ref{thm:AlphaModulationCoveringDefinition}. Then,
for $p,q\in\left(0,\infty\right]$ and $s,\mu\in\R$, we define the
associated \textbf{(weighted) $\alpha$-modulation space} as
\[
M_{\left(s,\mu\right),\alpha}^{p,q}\left(\smash{\R^{\dimension}}\right):=\DecompSp{\CalQ_{r}^{\left(\alpha\right)}}p{\ell_{w^{\left(s/\left(1-\alpha\right)\right)}}^{q}}{v^{\left(\mu\right)}}
\]
with $w^{\left(s/\left(1-\alpha\right)\right)}$ and $v^{\left(\mu\right)}$
as in Lemmas \ref{lem:AlphaModulationFrequencyWeight} and \ref{lem:AlphaModulationSpaceWeight},
respectively.

Furthermore, we define the \textbf{classical $\alpha$-modulation
space} as $M_{s,\alpha}^{p,q}\left(\R^{\dimension}\right):=M_{\left(s,0\right),\alpha}^{p,q}\left(\R^{\dimension}\right)$.
\end{defn}
\begin{rem*}

\begin{itemize}[leftmargin=0.4cm]
\item The classical $\alpha$-modulation spaces $M_{s,\alpha}^{p,q}\left(\R^{\dimension}\right)$
defined above coincide with the $\alpha$-modulation spaces defined
in \cite[Definition 2.4]{BorupNielsenAlphaModulationSpaces}, up to
trivial identifications: The quasi-norms used in the two definitions
are precisely the same; the only difference between the two definitions
is that in \cite[Definition 2.4]{BorupNielsenAlphaModulationSpaces},
the $\alpha$-modulation spaces are defined as subspaces of $\Schwartz'\left(\R^{\dimension}\right)$.
In contrast, with our definition as a decomposition space, $M_{s,\alpha}^{p,q}\left(\R^{\dimension}\right)$
is a subspace of $Z'\left(\R^{\dimension}\right)=\left[\Fourier\left(\TestFunctionSpace{\R^{\dimension}}\right)\right]'$,
cf.\@ Section \ref{subsec:DecompSpaceDefinitionStandingAssumptions}.
But \cite[Lemma 9.15]{DecompositionEmbedding} and \cite[Theorem 9.13]{DecompositionEmbedding}
show that each $f\in M_{s,\alpha}^{p,q}\left(\R^{\dimension}\right)$
extends to a (uniquely determined) tempered distribution, which implies
that the two different definitions of $\alpha$-modulation spaces
indeed yield the same spaces, up to trivial identifications.
\item Observe that the parameter $r>r_{1}\left(\dimension,\alpha\right)$
is suppressed on the left-hand side of the definition above. This
is justified, as we show now: Since any two coverings $\CalQ_{r}^{\left(\alpha\right)},\CalQ_{t}^{\left(\alpha\right)}$
(with $r,t>r_{1}\left(\dimension,\alpha\right)$) use the same families
$\left(T_{k}\right)_{k\in\Z^{\dimension}\setminus\left\{ 0\right\} }$
and $\left(b_{k}\right)_{k\in\Z^{\dimension}\setminus\left\{ 0\right\} }$,
it follows that every regular partition of unity $\Phi=\left(\varphi_{k}\right)_{k\in\Z^{\dimension}\setminus\left\{ 0\right\} }$
(cf.\@ Assumption \ref{assu:RegularPartitionOfUnity}) for $\CalQ_{r}^{\left(\alpha\right)}$
is also a regular partition of unity for $\CalQ_{t}^{\left(\alpha\right)}$,
at least for $t\geq r$, which we can always assume by symmetry. Thus,
by choosing the \emph{same} BAPU $\Phi$ for both coverings, we see
$\DecompSp{\CalQ_{r}^{\left(\alpha\right)}}p{\ell_{w^{\left(s^{\ast}\right)}}^{q}}{v^{\left(\mu\right)}}=\DecompSp{\CalQ_{t}^{\left(\alpha\right)}}p{\ell_{w^{\left(s^{\ast}\right)}}^{q}}{v^{\left(\mu\right)}}$,
with equivalent quasi-norms. Here, $s^{\ast}:=s/\left(1-\alpha\right)$.

We finally note that this argument implicitly uses that different
choices of the BAPU yield the same space (with equivalent quasi-norms),
cf.\@ Proposition \ref{prop:WeightedDecompositionSpaceWellDefined}.\qedhere
\end{itemize}
\end{rem*}
In the remainder of this section, we will determine conditions on
the prototype $\gamma$ which ensure that Corollary \ref{cor:BanachFrameSimplifiedCriteria}
(leading to Banach frames) or Corollary \ref{cor:AtomicDecompositionSimplifiedCriteria}
(leading to atomic decompositions) is applicable to $\gamma$. We
will see that this is the case for arbitrary Schwartz functions $\gamma$,
as long as $\widehat{\gamma}$ fulfills a certain non-vanishing condition.
To be precise, recall that in Corollaries \ref{cor:BanachFrameSimplifiedCriteria}
and \ref{cor:AtomicDecompositionSimplifiedCriteria}, we allowed the
prototype to depend on $i\in I$, i.e., we used a family $\left(\gamma_{i}\right)_{i\in I}$
of prototypes. But in this section, we will only consider the case
where $\gamma_{i}=\gamma$ is independent of $i\in I$.

To begin with, we recall that in Corollary \ref{cor:BanachFrameSimplifiedCriteria},
we are imposing certain summability conditions on
\[
M_{j,i}:=\left(\frac{w_{j}^{\left(s/\left(1-\alpha\right)\right)}}{w_{i}^{\left(s/\left(1-\alpha\right)\right)}}\right)^{\tau}\cdot\left(1+\left\Vert T_{j}^{-1}T_{i}\right\Vert \right)^{\sigma}\cdot\max_{\left|\beta\right|\leq1}\left(\left|\det T_{i}\right|^{-1}\cdot\int_{Q_{i}}\max_{\left|\alpha\right|\leq N}\left|\left(\partial^{\alpha}\widehat{\partial^{\beta}\gamma}\right)\!\!\left(S_{j}^{-1}\xi\right)\right|\d\xi\right)^{\tau}
\]
for suitable values of $\tau,\sigma>0$ and $N\in\N$. To slightly
simplify this expression, we will use the following notation for the
remainder of the section:
\begin{equation}
s^{\ast}:=s/\left(1-\alpha\right).\label{eq:AlphaModulationSmoothnessExponent}
\end{equation}

For our application of Corollary \ref{cor:BanachFrameSimplifiedCriteria},
we will assume $\gamma\in L^{1}\left(\R^{\dimension}\right)\cap C^{1}\left(\R^{\dimension}\right)$
with $\partial_{\ell}\gamma\in L^{1}\left(\R^{\dimension}\right)$
for all $\ell\in\underline{\dimension}$ and with $\widehat{\gamma}\in C^{\infty}\left(\R^{\dimension}\right)$.
Under these assumptions, elementary properties of the Fourier transform
yield for $\beta=e_{\ell}$ (the $\ell$-th unit vector) that
\[
\widehat{\partial^{\beta}\gamma}\left(\xi\right)=2\pi i\xi_{\ell}\cdot\widehat{\gamma}\left(\xi\right)\qquad\forall\xi\in\R^{\dimension}.
\]
Since we clearly have $\left|\frac{\partial^{\theta}}{\partial\eta^{\theta}}\eta_{\ell}\right|\leq1+\left|\eta\right|$
for all $\eta\in\R^{\dimension}$ and $\theta\in\N_{0}^{\dimension}$,
the Leibniz rule and the $\dimension$-dimensional binomial theorem
(cf.\@ \cite[Section 8.1, Exercise 2.b]{FollandRA}) yield
\begin{align}
\left|\left(\partial^{\alpha}\widehat{\partial^{\beta}\gamma}\right)\left(\eta\right)\right| & =\left|2\pi i\cdot\sum_{\theta\leq\alpha}\binom{\alpha}{\theta}\cdot\partial^{\theta}\eta_{\ell}\cdot\left(\partial^{\alpha-\theta}\widehat{\gamma}\right)\left(\eta\right)\right|\nonumber \\
 & \leq2\pi\cdot\left(1+\left|\eta\right|\right)\cdot\sum_{\theta\leq\alpha}\binom{\alpha}{\theta}\cdot\left|\left(\partial^{\alpha-\theta}\widehat{\gamma}\right)\left(\eta\right)\right|\nonumber \\
\left({\scriptstyle \text{eq. }\eqref{eq:AlphaModulationBanachFrameDecayAssumption}}\right) & \leq\left(1+\left|\eta\right|\right)^{1-N_{0}}\cdot2\pi C\cdot\sum_{\theta\leq\alpha}\binom{\alpha}{\theta}\nonumber \\
 & \leq2^{N+1}\pi\cdot C\cdot\left(1+\left|\eta\right|\right)^{1-N_{0}},\label{eq:AlphaModulationIteratedDerivativeEstimate}
\end{align}
where we used $\left|\alpha-\theta\right|\leq\left|\alpha\right|\leq N$
and assumed that there is some $N_{0}\in\R$ satisfying
\begin{equation}
\max_{\left|\alpha\right|\leq N}\left|\left(\partial^{\alpha}\widehat{\gamma}\right)\left(\eta\right)\right|\leq C\cdot\left(1+\left|\eta\right|\right)^{-N_{0}}\qquad\forall\eta\in\R^{\dimension}.\label{eq:AlphaModulationBanachFrameDecayAssumption}
\end{equation}
Recall that equation (\ref{eq:AlphaModulationIteratedDerivativeEstimate})
holds for $\beta=e_{\ell}$, with arbitrary $\ell\in\underline{\dimension}$.
But for $\beta=0$, we simply have 
\[
\left|\left(\partial^{\alpha}\widehat{\partial^{\beta}\gamma}\right)\left(\eta\right)\right|=\left|\left(\partial^{\alpha}\widehat{\gamma}\right)\left(\eta\right)\right|\leq C\cdot\left(1+\left|\eta\right|\right)^{-N_{0}}\leq2^{N+1}\pi\cdot C\cdot\left(1+\left|\eta\right|\right)^{1-N_{0}},
\]
so that we have verified equation (\ref{eq:AlphaModulationIteratedDerivativeEstimate})
for arbitrary $\alpha,\beta\in\N_{0}^{\dimension}$ with $\left|\beta\right|\leq1$
and $\left|\alpha\right|\leq N$.

Hence, we have shown for $C':=2^{N+1}\pi\cdot C$ that
\begin{equation}
M_{j,i}\leq\left(C'\right)^{\tau}\cdot\left(\frac{w_{j}^{\left(s^{\ast}\right)}}{w_{i}^{\left(s^{\ast}\right)}}\right)^{\tau}\cdot\left(1+\left\Vert T_{j}^{-1}T_{i}\right\Vert \right)^{\sigma}\cdot\left(\left|\det T_{i}\right|^{-1}\cdot\int_{Q_{i}}\left(1+\left|S_{j}^{-1}\xi\right|\right)^{1-N_{0}}\d\xi\right)^{\tau}=:\left(C'\right)^{\tau}\cdot M_{j,i}^{\left(0\right)}\label{eq:AlphaModulationStandardEstimate}
\end{equation}
for all $i,j\in\Z^{\dimension}\setminus\left\{ 0\right\} $. In view
of this estimate, the following lemma is crucial:
\begin{lem}
\label{lem:AlphaModulationMainLemma}Let $\dimension\in\N$ and $\alpha\in\left[0,1\right)$
and set $\alpha_{0}:=\frac{\alpha}{1-\alpha}$. With $\CalQ=\CalQ_{r}^{\left(\alpha\right)}$
and with $M_{j,i}^{\left(0\right)}$ as in equation (\ref{eq:AlphaModulationStandardEstimate}),
assume
\[
N_{0}\geq\dimension+2+\frac{\dimension+1}{\tau}+\max\left\{ \left|s^{\ast}+\dimension\alpha_{0}\right|,\left|s^{\ast}+\left(\dimension-\frac{\sigma}{\tau}\right)\alpha_{0}\right|\right\} .
\]
Then we have
\[
\sup_{i\in\Z^{\dimension}\setminus\left\{ 0\right\} }\:\sum_{j\in\Z^{\dimension}\setminus\left\{ 0\right\} }M_{j,i}^{\left(0\right)}\leq\Omega\qquad\text{ and }\qquad\sup_{j\in\Z^{\dimension}\setminus\left\{ 0\right\} }\:\sum_{i\in\Z^{\dimension}\setminus\left\{ 0\right\} }M_{j,i}^{\left(0\right)}\leq\Omega
\]
for
\[
\Omega:=6^{\dimension}2^{1+\sigma+\tau\left|s^{\ast}\right|}\!\cdot\max\left\{ 4^{\alpha_{0}\left(\sigma+\dimension\tau\right)+\tau\left|s^{\ast}\right|}\cdot\left(12^{N_{0}}s_{\dimension}\right)^{\tau},\:\left(2\!+\!4r\right)^{\tau\left|s^{\ast}\right|+\alpha_{0}\left[\sigma+\tau\left(\dimension+N_{0}\right)\right]}\!\cdot\!2^{\tau\dimension}\!\cdot\!\left(1\!+\!\left(2\!+\!4r\right)^{\alpha_{0}}\!\cdot r\right)^{\tau\left(N_{0}+\dimension\right)}\!\right\} .\qedhere
\]
\end{lem}
\begin{proof}
For brevity, set $M:=N_{0}-1$. Recall $T_{j}=\left|j\right|^{\alpha_{0}}\cdot\identity$
and $b_{j}=\left|j\right|^{\alpha_{0}}j$, so that
\[
S_{j}^{-1}\xi=T_{j}^{-1}\left(\xi-b_{j}\right)=\left|j\right|^{-\alpha_{0}}\left(\xi-\left|j\right|^{\alpha_{0}}j\right)=\left|j\right|^{-\alpha_{0}}\xi-j.
\]
Hence,
\begin{align*}
\int_{Q_{i}}\left(1+\left|S_{j}^{-1}\xi\right|\right)^{1-N_{0}}\d\xi & =\int_{B_{\left|i\right|^{\alpha_{0}}r}\left(\left|i\right|^{\alpha_{0}}i\right)}\:\left(1+\left|\left|j\right|^{-\alpha_{0}}\xi-j\right|\right)^{-M}\d\xi\\
\left({\scriptstyle \text{with }\eta=\left|j\right|^{-\alpha_{0}}\xi}\right) & =\left|j\right|^{\dimension\alpha_{0}}\cdot\int_{B_{\left(\left|i\right|/\left|j\right|\right)^{\alpha_{0}}\cdot r}\left(\left(\left|i\right|/\left|j\right|\right)^{\alpha_{0}}i\right)}\:\left(1+\left|\eta-j\right|\right)^{-M}\d\eta\\
\left({\scriptstyle \text{with }\xi=\eta-j}\right) & =\left|j\right|^{\dimension\alpha_{0}}\cdot\int_{B_{\left(\left|i\right|/\left|j\right|\right)^{\alpha_{0}}\cdot r}\left(\left(\left|i\right|/\left|j\right|\right)^{\alpha_{0}}i-j\right)}\:\left(1+\left|\xi\right|\right)^{-M}\d\xi\\
 & =\left|j\right|^{\dimension\alpha_{0}}\cdot\int_{B_{R_{i,j}}\left(\xi_{i,j}\right)}\:\left(1+\left|\xi\right|\right)^{-M}\d\xi,
\end{align*}
where we defined
\[
\xi_{i,j}:=\left(\frac{\left|i\right|}{\left|j\right|}\right)^{\alpha_{0}}i-j\qquad\text{ and }\qquad R_{i,j}:=\left(\frac{\left|i\right|}{\left|j\right|}\right)^{\alpha_{0}}\cdot r
\]
for $i,j\in I=\Z^{\dimension}\setminus\left\{ 0\right\} $. Here,
$r>r_{1}\left(\dimension,\alpha\right)$ comes from the covering $\CalQ_{r}^{\left(\alpha\right)}$.

All in all, since $\left\Vert T_{j}^{-1}T_{i}\right\Vert =\left\Vert \left|i\right|^{\alpha_{0}}/\left|j\right|^{\alpha_{0}}\cdot\identity\right\Vert =\left(\left|i\right|/\left|j\right|\right)^{\alpha_{0}}$
and $\left|\det T_{i}\right|=\left|\det\left|i\right|^{\alpha_{0}}\identity\right|=\left|i\right|^{\dimension\alpha_{0}}$,
and because of $\left|i\right|\leq\left\langle i\right\rangle \leq1+\left|i\right|\leq2\left|i\right|$
for $i\in\Z^{\dimension}$, so that $\frac{1}{2}\frac{\left|j\right|}{\left|i\right|}\leq\frac{\left\langle j\right\rangle }{\left\langle i\right\rangle }\leq2\frac{\left|j\right|}{\left|i\right|}$,
we get the following estimate for $M_{j,i}^{\left(0\right)}$:
\begin{align}
M_{j,i}^{\left(0\right)} & =\left(\frac{w_{j}^{\left(s^{\ast}\right)}}{w_{i}^{\left(s^{\ast}\right)}}\right)^{\tau}\cdot\left(1+\left\Vert T_{j}^{-1}T_{i}\right\Vert \right)^{\sigma}\cdot\left(\left|\det T_{i}\right|^{-1}\cdot\int_{Q_{i}}\left(1+\left|S_{j}^{-1}\xi\right|\right)^{1-N_{0}}\d\xi\right)^{\tau}\nonumber \\
 & =\left(\frac{\left\langle j\right\rangle }{\left\langle i\right\rangle }\right)^{s^{\ast}\cdot\tau}\cdot\left(1+\left(\frac{\left|i\right|}{\left|j\right|}\right)^{\alpha_{0}}\right)^{\sigma}\cdot\left[\left(\frac{\left|j\right|}{\left|i\right|}\right)^{\dimension\alpha_{0}}\cdot\int_{B_{R_{i,j}}\left(\xi_{i,j}\right)}\left(1+\left|\xi\right|\right)^{-M}\d\xi\right]^{\tau}\nonumber \\
\left({\scriptstyle \text{since }\left(1+a\right)^{\sigma}\leq2^{\sigma}\cdot\left(1+a^{\sigma}\right)}\right) & \leq2^{\sigma+\tau\left|s^{\ast}\right|}\cdot\left(\frac{\left|j\right|}{\left|i\right|}\right)^{\tau\left(s^{\ast}+\dimension\alpha_{0}\right)}\cdot\left(1+\left(\frac{\left|j\right|}{\left|i\right|}\right)^{-\sigma\alpha_{0}}\right)\cdot\left(\int_{B_{R_{i,j}}\left(\xi_{i,j}\right)}\left(1+\left|\xi\right|\right)^{-M}\d\xi\right)^{\tau}\nonumber \\
 & =2^{\sigma+\tau\left|s^{\ast}\right|}\cdot\sum_{\lambda\in\left\{ 0,1\right\} }\left[\left(\frac{\left|j\right|}{\left|i\right|}\right)^{k_{\lambda}}\cdot\left(\int_{B_{R_{i,j}}\left(\xi_{i,j}\right)}\left(1+\left|\xi\right|\right)^{-M}\d\xi\right)^{\!\tau}\,\right],\label{eq:AlphaModulationMXiConnection}
\end{align}
where we defined $k_{\lambda}:=\tau\left(s^{\ast}+\dimension\alpha_{0}\right)-\lambda\sigma\alpha_{0}$
for $\lambda\in\left\{ 0,1\right\} $.

Thus, our main goal is to estimate the term
\[
\Xi_{i,j}^{\left(k\right)}:=\left(\frac{\left|j\right|}{\left|i\right|}\right)^{k}\cdot\left(\int_{B_{R_{i,j}}\left(\xi_{i,j}\right)}\left(1+\left|\xi\right|\right)^{-M}\d\xi\right)^{\tau}
\]
for arbitrary $i,j\in\Z^{\dimension}\setminus\left\{ 0\right\} $,
$k\in\R$ and $\tau>0$. To this end, we will distinguish three cases
concerning $i,j$ below. But before that, we introduce a useful notation
and some related estimates that will be used in several of the cases:
For $x\in\R^{\dimension}$, we set $\left[x\right]:=1+\left|x\right|$.
We then have
\begin{equation}
\left[x\right]^{z}\leq\left[y\right]^{z}\cdot\left[x-y\right]^{\left|z\right|}\qquad\forall z\in\R\quad\forall x,y\in\R^{\dimension}.\label{eq:AlphaModulationBracketModerate}
\end{equation}
Indeed, since $\left[x\right]=1+\left|x\right|\leq1+\left|y\right|+\left|x-y\right|\leq\left(1+\left|y\right|\right)\left(1+\left|x-y\right|\right)=\left[y\right]\cdot\left[x-y\right]$,
we get the claim for $z\geq0$. Finally, for $z<0$, we have
\begin{align*}
\left[x\right]^{z}=\left[x\right]^{-\left|z\right|} & =\left(\left[x\right]^{-\left|z\right|}\left[x-y\right]^{-\left|z\right|}\right)\cdot\left[x-y\right]^{\left|z\right|}\\
\left({\scriptstyle \text{eq. }\eqref{eq:AlphaModulationBracketModerate}\text{ rearranged, with }x,y\text{ interchanged and }\left|z\right|\text{ instead of }z}\right) & \leq\left[y\right]^{-\left|z\right|}\cdot\left[x-y\right]^{\left|z\right|}=\left[y\right]^{z}\cdot\left[x-y\right]^{\left|z\right|},
\end{align*}
as desired. Now, note for $i,j\in\Z^{\dimension}\setminus\left\{ 0\right\} $
that $\left|i\right|\leq\left[i\right]\leq2\left|i\right|$ and likewise
for $j$, so that
\begin{equation}
\left|i\right|^{z}\leq2^{\left|z\right|}\cdot\left[i\right]^{z}\leq2^{\left|z\right|}\cdot\left[j\right]^{z}\cdot\left[i-j\right]^{\left|z\right|}\leq4^{\left|z\right|}\cdot\left|j\right|^{z}\cdot\left[i-j\right]^{\left|z\right|}.\label{eq:AlphaModulationAbsoluteValueBracketModerate}
\end{equation}

We will also need the following estimate, which I learned from \cite{EmbeddingsOfAlphaModulationIntoSobolev}:
\begin{equation}
\left|\beta\cdot x-y\right|\geq\left|x-y\right|\qquad\text{ if }\beta\in\R_{\geq1}\text{ and }x,y\in\R^{\dimension}\text{ with }\left|x\right|\geq\left|y\right|.\label{eq:AlphaModulationKatoEstimate}
\end{equation}
For $\beta=1$, this estimate is trivial, so that we can assume $\beta>1$.
Next, note that both sides are nonnegative, so that the estimate is
equivalent to $\left|\beta\cdot x-y\right|^{2}\geq\left|x-y\right|^{2}$
and thus to
\begin{align*}
 & \beta^{2}\left|x\right|^{2}-2\beta\cdot\left\langle x,y\right\rangle +\left|y\right|^{2}\overset{!}{\geq}\left|x\right|^{2}-2\left\langle x,y\right\rangle +\left|y\right|^{2}\\
\Longleftrightarrow & \left|x\right|^{2}\cdot\left(\beta^{2}-1\right)\overset{!}{\geq}2\cdot\left\langle x,y\right\rangle \cdot\left(\beta-1\right)\\
\left({\scriptstyle \text{since }\beta-1>0}\right)\Longleftrightarrow & \left|x\right|^{2}\cdot\left(\beta+1\right)\overset{!}{\geq}2\cdot\left\langle x,y\right\rangle .
\end{align*}
But the Cauchy-Schwarz inequality yields $2\cdot\left\langle x,y\right\rangle \leq2\cdot\left|\left\langle x,y\right\rangle \right|\leq2\cdot\left|x\right|\left|y\right|\leq2\cdot\left|x\right|^{2}\leq\left(1+\beta\right)\cdot\left|x\right|^{2}$,
since $\left|y\right|\leq\left|x\right|$ and since $\beta\geq1$.
Hence, we have established equation (\ref{eq:AlphaModulationKatoEstimate}).
We remark that for this estimate, it is crucial to use a norm which
is induced by a scalar product. For other norms, equation (\ref{eq:AlphaModulationKatoEstimate})
can fail.

\medskip{}

Now, we distinguish three cases depending on $i,j$:

\textbf{Case 1}: We have $\left|i\right|\geq2\left|j\right|+4r$.
In this case, we get
\begin{align*}
\left|\left(\frac{\left|i\right|}{\left|j\right|}\right)^{\alpha_{0}}i\right| & \leq\left|\left(\frac{\left|i\right|}{\left|j\right|}\right)^{\alpha_{0}}i-j\right|+\left|j\right|\\
 & \leq\left|\xi_{i,j}\right|+\frac{\left|i\right|}{2}\\
\left({\scriptstyle \text{since }\frac{\left|i\right|}{\left|j\right|}\geq2\geq1\text{ and }\alpha_{0}\geq0}\right) & \leq\left|\xi_{i,j}\right|+\frac{1}{2}\left|\left(\frac{\left|i\right|}{\left|j\right|}\right)^{\alpha_{0}}i\right|
\end{align*}
and thus, since $\left|i\right|\geq4r$,
\begin{equation}
\left|\xi_{i,j}\right|\geq\frac{1}{2}\left|\left(\frac{\left|i\right|}{\left|j\right|}\right)^{\alpha_{0}}i\right|\geq2r\cdot\left(\frac{\left|i\right|}{\left|j\right|}\right)^{\alpha_{0}}=2\cdot R_{i,j}.\label{eq:AlphaModulationILargeCenterLarge}
\end{equation}
Hence, for arbitrary $\xi\in B_{R_{i,j}}\left(\xi_{i,j}\right)$,
we have $\left|\xi\right|\geq\left|\xi_{i,j}\right|-\left|\xi-\xi_{i,j}\right|\geq\left|\xi_{i,j}\right|-R_{i,j}\geq\frac{1}{2}\left|\xi_{i,j}\right|$
and thus
\begin{equation}
\begin{split}\int_{B_{R_{i,j}}\left(\xi_{i,j}\right)}\left(1+\left|\xi\right|\right)^{-M}\d\xi & \leq\left[\sup_{\xi\in B_{R_{i,j}}\left(\xi_{i,j}\right)}\left(1+\left|\xi\right|\right)^{\dimension+1-M}\right]\cdot\int_{B_{R_{i,j}}\left(\xi_{i,j}\right)}\left(1+\left|\xi\right|\right)^{-\left(\dimension+1\right)}\d\xi\\
\left({\scriptstyle \text{since }\dimension+1-M\leq0}\right) & \leq\left(1+\frac{1}{2}\left|\xi_{i,j}\right|\right)^{\dimension+1-M}\cdot\int_{\R^{\dimension}}\left(1+\left|\xi\right|\right)^{-\left(\dimension+1\right)}\d\xi\\
\left({\scriptstyle \text{eq. }\eqref{eq:StandardDecayLpEstimate}}\right) & \leq2^{M}\cdot s_{\dimension}\cdot\left|\xi_{i,j}\right|^{\dimension+1-M}\\
\left({\scriptstyle \text{eq. }\eqref{eq:AlphaModulationILargeCenterLarge}\text{ and }\dimension+1-M\leq0}\right) & \leq4^{M}\cdot s_{\dimension}\cdot\left|\left(\frac{\left|i\right|}{\left|j\right|}\right)^{\alpha_{0}}i\right|^{\dimension+1-M}\\
\left({\scriptstyle \text{since }\dimension+1-M\leq0\text{ and }\frac{\left|i\right|}{\left|j\right|}\geq1}\right) & \leq4^{M}\cdot s_{\dimension}\cdot\left|i\right|^{\dimension+1-M}.
\end{split}
\label{eq:AlphaModulationIntegralSupremumEstimate}
\end{equation}
Next, we observe $\left|i\right|\geq2\left|j\right|+4r\geq\left|j\right|$
and $i\in\Z^{\dimension}\setminus\left\{ 0\right\} $, so that $\left|i\right|\geq1$.
This implies
\[
\left[i-j\right]=1+\left|i-j\right|\leq1+\left|i\right|+\left|j\right|\leq1+2\left|i\right|\leq3\left|i\right|,
\]
so that we finally arrive, again using $\dimension+1-M\leq0$, at
\[
\int_{B_{R_{i,j}}\left(\xi_{i,j}\right)}\left(1+\left|\xi\right|\right)^{-M}\d\xi\leq12^{M}\cdot s_{\dimension}\cdot\left[i-j\right]^{\dimension+1-M}.
\]
Thus, using equation (\ref{eq:AlphaModulationAbsoluteValueBracketModerate}),
we conclude
\begin{equation}
\Xi_{i,j}^{\left(k\right)}=\left(\frac{\left|j\right|}{\left|i\right|}\right)^{k}\cdot\left(\int_{B_{R_{i,j}}\left(\xi_{i,j}\right)}\left(1+\left|\xi\right|\right)^{-M}\d\xi\right)^{\tau}\leq4^{\left|k\right|}\cdot\left(12^{M}\cdot s_{\dimension}\right)^{\tau}\cdot\left[j-i\right]^{\left|k\right|+\tau\left(\dimension+1-M\right)}.\label{eq:AlphaModulationXiEstimateFirstCase}
\end{equation}

\textbf{Case 2}: We have $\left|j\right|\geq2\left|i\right|+4r$.
Here, we first observe $\left|j-i\right|\geq\left|j\right|-\left|i\right|\geq\left|i\right|+4r\geq4r$.
Hence, we get
\begin{align*}
\left|\xi_{i,j}\right| & =\left(\frac{\left|i\right|}{\left|j\right|}\right)^{\alpha_{0}}\cdot\left|\left(\frac{\left|j\right|}{\left|i\right|}\right)^{\alpha_{0}}j-i\right|\\
\left({\scriptstyle \text{eq. }\eqref{eq:AlphaModulationKatoEstimate}\text{ and }\left(\left|j\right|/\left|i\right|\right)^{\alpha_{0}}\geq1,\text{ as well as }\left|j\right|\geq\left|i\right|}\right) & \geq\left(\frac{\left|i\right|}{\left|j\right|}\right)^{\alpha_{0}}\cdot\left|j-i\right|\\
 & \geq4\cdot\left(\frac{\left|i\right|}{\left|j\right|}\right)^{\alpha_{0}}r=4\cdot R_{i,j}.
\end{align*}
Further, $\left|j\right|\geq2\left|i\right|+4r\geq2\left|i\right|$
implies $\left|i\right|\leq\frac{\left|j\right|}{2}$ and thus $\left|j\right|-\left|i\right|\geq\frac{1}{2}\left|j\right|$.
Consequently, 
\[
\left|j-i\right|\leq\left|j\right|+\left|i\right|\leq2\left|j\right|\leq4\cdot\left(\left|j\right|-\left|i\right|\right),
\]
so that we get
\begin{align*}
\left|\xi_{i,j}\right| & =\left|\left(\frac{\left|i\right|}{\left|j\right|}\right)^{\alpha_{0}}i-j\right|\\
 & \geq\left|j\right|-\left(\frac{\left|i\right|}{\left|j\right|}\right)^{\alpha_{0}}\left|i\right|\\
\left({\scriptstyle \text{since }\left|i\right|/\left|j\right|\leq1}\right) & \geq\left|j\right|-\left|i\right|\geq\frac{\left|j-i\right|}{4}.
\end{align*}
Now, for arbitrary $\xi\in B_{R_{i,j}}\left(\xi_{i,j}\right)$, the
two preceding displayed estimates yield $\left|\xi\right|\geq\left|\xi_{i,j}\right|-R_{i,j}\geq\frac{3}{4}\left|\xi_{i,j}\right|\geq\frac{3}{16}\cdot\left|j-i\right|$
and hence $1+\left|\xi\right|\geq\frac{3}{16}\cdot\left[j-i\right]$.
With an estimate entirely analogous to that in equation (\ref{eq:AlphaModulationIntegralSupremumEstimate}),
this implies 
\begin{align*}
\int_{B_{R_{i,j}}\left(\xi_{i,j}\right)}\left(1+\left|\xi\right|\right)^{-M}\d\xi & \leq\left[\sup_{\xi\in B_{R_{i,j}}\left(\xi_{i,j}\right)}\left(1+\left|\xi\right|\right)^{\dimension+1-M}\right]\cdot\int_{B_{R_{i,j}}\left(\xi_{i,j}\right)}\left(1+\left|\xi\right|\right)^{-\left(\dimension+1\right)}\d\xi\\
\left({\scriptstyle \text{eq. }\eqref{eq:StandardDecayLpEstimate}\text{ and }\dimension+1-M\leq0}\right) & \leq s_{\dimension}\cdot\left(\frac{16}{3}\right)^{M}\cdot\left[j-i\right]^{\dimension+1-M}.
\end{align*}
In view of equation (\ref{eq:AlphaModulationAbsoluteValueBracketModerate})
and since $\frac{16}{3}\leq12$, we conclude 
\[
\Xi_{i,j}^{\left(k\right)}=\left(\frac{\left|j\right|}{\left|i\right|}\right)^{k}\cdot\left(\int_{B_{R_{i,j}}\left(\xi_{i,j}\right)}\left(1+\left|\xi\right|\right)^{-M}\d\xi\right)^{\tau}\leq4^{\left|k\right|}\cdot\left(12^{M}\cdot s_{\dimension}\right)^{\tau}\cdot\left[j-i\right]^{\left|k\right|+\tau\left(\dimension+1-M\right)},
\]
as in the previous case.

\textbf{Case 3}: The remaining case, i.e., $\left|i\right|<2\left|j\right|+4r$
\textbf{and} $\left|j\right|<2\left|i\right|+4r$. Since $i,j\in\Z^{\dimension}\setminus\left\{ 0\right\} $,
we have $\left|i\right|,\left|j\right|\geq1$ and thus $\left|i\right|\leq\left|j\right|\cdot\left(2+4r\right)$
and $\left|j\right|\leq\left|i\right|\cdot\left(2+4r\right)$. In
particular, we have
\[
R_{i,j}=\left(\frac{\left|i\right|}{\left|j\right|}\right)^{\alpha_{0}}\cdot r\leq r\cdot\left(2+4r\right)^{\alpha_{0}}=:C_{r,\alpha_{0}}\:,
\]
so that every $\xi\in B_{R_{i,j}}\left(\xi_{i,j}\right)$ satisfies
\begin{align*}
1+\left|\xi_{i,j}\right| & \leq1+\left|\xi_{i,j}-\xi\right|+\left|\xi\right|\leq1+R_{i,j}+\left|\xi\right|\\
 & \leq\left(1+R_{i,j}\right)\left(1+\left|\xi\right|\right)\leq\left(1+C_{r,\alpha_{0}}\right)\left(1+\left|\xi\right|\right).
\end{align*}
Consequently,
\begin{align*}
\int_{B_{R_{i,j}}\left(\xi_{i,j}\right)}\,\left(1+\left|\xi\right|\right)^{-M}\d\xi & \leq\left(1+C_{r,\alpha_{0}}\right)^{M}\cdot\lambda_{\dimension}\left(B_{R_{i,j}}\left(\xi_{i,j}\right)\right)\cdot\left(1+\left|\xi_{i,j}\right|\right)^{-M}\\
\left({\scriptstyle \text{since }\lambda_{\dimension}\left(B_{1}\left(0\right)\right)\leq\lambda_{\dimension}\left(\left[-1,1\right]^{\dimension}\right)=2^{\dimension}}\right) & \leq2^{\dimension}\cdot\left(1+C_{r,\alpha_{0}}\right)^{M}R_{i,j}^{\dimension}\cdot\left(1+\left|\xi_{i,j}\right|\right)^{-M}\\
\left({\scriptstyle \text{since }R_{i,j}\leq C_{r,\alpha_{0}}}\right) & \leq2^{\dimension}\cdot\left(1+C_{r,\alpha_{0}}\right)^{M+\dimension}\cdot\left(1+\left|\xi_{i,j}\right|\right)^{-M}=:C_{\dimension,M,r,\alpha_{0}}\cdot\left(1+\left|\xi_{i,j}\right|\right)^{-M}.
\end{align*}
Furthermore, since $\frac{1}{2+4r}\leq\frac{\left|j\right|}{\left|i\right|}\leq2+4r$,
we get $\left(\frac{\left|j\right|}{\left|i\right|}\right)^{k}\leq\left(2+4r\right)^{\left|k\right|}$
and thus
\[
\Xi_{i,j}^{\left(k\right)}\leq\left(2+4r\right)^{\left|k\right|}\cdot C_{\dimension,M,r,\alpha_{0}}^{\tau}\cdot\left(1+\left|\xi_{i,j}\right|\right)^{-\tau M}.
\]
To further estimate the right-hand side, we distinguish two sub-cases:

\begin{enumerate}
\item We have $\left|i\right|\geq\left|j\right|$. In this case, equation
(\ref{eq:AlphaModulationKatoEstimate}) yields
\[
\left|\xi_{i,j}\right|=\left|\left(\frac{\left|i\right|}{\left|j\right|}\right)^{\alpha_{0}}i-j\right|\geq\left|i-j\right|
\]
and hence
\begin{align*}
\Xi_{i,j}^{\left(k\right)} & \leq\left(2+4r\right)^{\left|k\right|}\cdot C_{\dimension,M,r,\alpha_{0}}^{\tau}\cdot\left[i-j\right]^{-\tau M}\\
 & \leq\left(2+4r\right)^{\left|k\right|}\cdot C_{\dimension,M,r,\alpha_{0}}^{\tau}\cdot\left[i-j\right]^{\left|k\right|+\tau\left(\dimension+1-M\right)}.
\end{align*}
\item We have $\left|j\right|\geq\left|i\right|$. In this case, we can
again—after some rearranging—use equation (\ref{eq:AlphaModulationKatoEstimate})
to obtain
\begin{align*}
\left|\xi_{i,j}\right| & =\left|\left(\frac{\left|i\right|}{\left|j\right|}\right)^{\alpha_{0}}i-j\right|=\left(\frac{\left|i\right|}{\left|j\right|}\right)^{\alpha_{0}}\left|\left(\frac{\left|j\right|}{\left|i\right|}\right)^{\alpha_{0}}j-i\right|\\
\left({\scriptstyle \text{eq. }\eqref{eq:AlphaModulationKatoEstimate}}\right) & \geq\left(2+4r\right)^{-\alpha_{0}}\cdot\left|j-i\right|
\end{align*}
and hence
\begin{align*}
\Xi_{i,j}^{\left(k\right)} & \leq\left(2+4r\right)^{\left|k\right|}\cdot C_{\dimension,M,r,\alpha_{0}}^{\tau}\cdot\left(1+\left(2+4r\right)^{-\alpha_{0}}\cdot\left|j-i\right|\right)^{-\tau M}\\
 & \leq\left(2+4r\right)^{\left|k\right|+\alpha_{0}\tau M}\cdot C_{\dimension,M,r,\alpha_{0}}^{\tau}\cdot\left[j-i\right]^{-\tau M}\\
 & \leq\left(2+4r\right)^{\left|k\right|+\alpha_{0}\tau M}\cdot C_{\dimension,M,r,\alpha_{0}}^{\tau}\cdot\left[j-i\right]^{\left|k\right|+\tau\left(\dimension+1-M\right)}.
\end{align*}
\end{enumerate}
All in all, the preceding case distinction has established the bound
\begin{align}
\Xi_{i,j}^{\left(k\right)} & \leq\max\left\{ 4^{\left|k\right|}\cdot\left(12^{M}\cdot s_{\dimension}\right)^{\tau},\:\left(2+4r\right)^{\left|k\right|+\alpha_{0}\tau M}\cdot C_{\dimension,M,r,\alpha_{0}}^{\tau}\right\} \cdot\left[j-i\right]^{\left|k\right|+\tau\left(\dimension+1-M\right)}\nonumber \\
 & =:C_{\dimension,M,r,\alpha_{0},k,\tau}\cdot\left[j-i\right]^{\left|k\right|+\tau\left(\dimension+1-M\right)}\label{eq:AlphaModulationFinalXiEstimate}
\end{align}
for all $i,j\in\Z^{\dimension}\setminus\left\{ 0\right\} $ and all
$k\in\R$, with $C_{\dimension,M,r,\alpha_{0}}=2^{\dimension}\cdot\left(1+C_{r,\alpha_{0}}\right)^{M+\dimension}$
and $C_{r,\alpha_{0}}=r\cdot\left(2+4r\right)^{\alpha_{0}}$.

Now, we want to utilize this estimate for $k=k_{\lambda}=\tau\left(s^{\ast}+\dimension\alpha_{0}\right)-\lambda\sigma\alpha_{0}$
for $\lambda\in\left\{ 0,1\right\} $, cf.\@ equation (\ref{eq:AlphaModulationMXiConnection}).
Note the equivalence
\begin{align*}
 & \left|k_{\lambda}\right|+\tau\left(\dimension+1-M\right)=\left|k_{\lambda}\right|+\tau\left(\dimension+2-N_{0}\right)\overset{!}{\leq}-\left(\dimension+1\right)\\
\Longleftrightarrow & N_{0}\overset{!}{\geq}\dimension+2+\frac{\left|k_{\lambda}\right|+\dimension+1}{\tau},
\end{align*}
where the last condition is satisfied for $\lambda\in\left\{ 0,1\right\} $
by definition of $k_{\lambda}$ and our assumptions regarding $N_{0}$.
Hence, we get—in view of equations (\ref{eq:StandardDecayLatticeSeries})
and (\ref{eq:AlphaModulationFinalXiEstimate}) and because of $\left|j-i\right|\geq\left\Vert j-i\right\Vert _{\infty}$—that
\begin{align*}
\sum_{j\in\Z^{\dimension}\setminus\left\{ 0\right\} }\Xi_{i,j}^{\left(k_{\lambda}\right)} & \leq\left[\max_{\lambda\in\left\{ 0,1\right\} }C_{\dimension,M,r,\alpha_{0},k_{\lambda},\tau}\right]\cdot\sum_{j\in\Z^{\dimension}}\left[j-i\right]^{-\left(\dimension+1\right)}\\
\left({\scriptstyle \text{with }\ell=j-i}\right) & \leq\left[\max_{\lambda\in\left\{ 0,1\right\} }C_{\dimension,M,r,\alpha_{0},k_{\lambda},\tau}\right]\cdot\sum_{\ell\in\Z^{\dimension}}\left(1+\left\Vert \ell\right\Vert _{\infty}\right)^{-\left(\dimension+1\right)}\leq6^{\dimension}\cdot\max_{\lambda\in\left\{ 0,1\right\} }C_{\dimension,M,r,\alpha_{0},k_{\lambda},\tau}.
\end{align*}
The same estimate also holds when taking the sum over $i\in\Z^{\dimension}\setminus\left\{ 0\right\} $
instead of over $j\in\Z^{\dimension}\setminus\left\{ 0\right\} $.
In view of equation (\ref{eq:AlphaModulationMXiConnection}), we thus
get
\[
\sup_{i\in\Z^{\dimension}\setminus\left\{ 0\right\} }\:\sum_{j\in\Z^{\dimension}\setminus\left\{ 0\right\} }M_{j,i}^{\left(0\right)}\leq2^{1+\sigma+\tau\left|s^{\ast}\right|}\cdot6^{\dimension}\cdot\max_{\lambda\in\left\{ 0,1\right\} }C_{\dimension,M,r,\alpha_{0},k_{\lambda},\tau}
\]
and the same estimate also holds for $\sup_{j\in\Z^{\dimension}\setminus\left\{ 0\right\} }\:\sum_{i\in\Z^{\dimension}\setminus\left\{ 0\right\} }M_{j,i}^{\left(0\right)}$.
This easily yields the claim.
\end{proof}
Now, we can derive readily verifiable conditions which ensure that
the structured family generated by $\gamma$ yields a Banach frame
for a given $\alpha$-modulation space.
\begin{thm}
\label{thm:AlphaModulationBanachFrame}Let $\dimension\in\N$, $\alpha\in\left[0,1\right)$
and choose $r>r_{1}\left(\dimension,\alpha\right)$ with $r_{1}\left(\dimension,\alpha\right)$
as in Theorem \ref{thm:AlphaModulationCoveringDefinition}.

Let $s_{0},\mu_{0}\geq0$ and $p_{0},q_{0}\in\left(0,1\right]$, as
well as $\varepsilon\in\left(0,1\right)$. Assume that $\gamma:\R^{\dimension}\to\Compl$
satisfies the following:

\begin{enumerate}
\item We have $\gamma\in L_{\left(1+\left|\mybullet\right|\right)^{\mu_{0}}}^{1}\left(\R^{\dimension}\right)$
and $\widehat{\gamma}\in C^{\infty}\left(\R^{\dimension}\right)$,
where all partial derivatives of $\widehat{\gamma}$ are polynomially
bounded.
\item We have $\gamma\in C^{1}\left(\R^{\dimension}\right)$ and $\partial_{\ell}\gamma\in L_{\left(1+\left|\mybullet\right|\right)^{\mu_{0}}}^{1}\left(\R^{\dimension}\right)\cap L^{\infty}\left(\R^{\dimension}\right)$
for all $\ell\in\underline{\dimension}$.
\item \label{enu:AlphaModulationFrameFourierNonVanishing}We have $\left|\widehat{\gamma}\left(\xi\right)\right|\geq c>0$
for all $\xi\in\overline{B_{r}}\left(0\right)$.
\item \label{enu:AlphaModulationFrameFourierDecay}We have
\[
\left|\left(\partial^{\beta}\widehat{\gamma}\right)\left(\xi\right)\right|\lesssim\left(1+\left|\xi\right|\right)^{-N_{0}}
\]
for all $\xi\in\R^{\dimension}$ and all $\beta\in\N_{0}^{\dimension}$
with $\left|\beta\right|\leq\left\lceil \mu_{0}+\frac{\dimension+\varepsilon}{p_{0}}\right\rceil $,
where
\[
\qquad\qquad N_{0}=\dimension+2+\frac{\dimension+1}{\min\left\{ p_{0},q_{0}\right\} }+\frac{1}{1-\alpha}\cdot\max\left\{ s_{0}+\alpha\dimension,\:s_{0}+\alpha\left(\frac{\dimension}{p_{0}}-\dimension+\mu_{0}+\left\lceil \mu_{0}+\frac{\dimension+\varepsilon}{p_{0}}\right\rceil \right)\right\} .
\]
\end{enumerate}
Then there is some $\delta_{0}>0$ such that for all $0<\delta\leq\delta_{0}$,
the family 
\[
\Gamma^{\left(\delta\right)}:=\left(L_{\delta\cdot k/\left|i\right|^{\alpha_{0}}}\:\widetilde{\gamma^{\left[i\right]}}\right)_{i\in\Z^{\dimension}\setminus\left\{ 0\right\} ,k\in\Z^{\dimension}},\quad\text{ with }\quad\gamma^{\left[i\right]}=\left|i\right|^{\frac{\dimension\cdot\alpha_{0}}{2}}\cdot M_{\left|i\right|^{\alpha_{0}}\cdot i}\left[\gamma\circ\left|i\right|^{\alpha_{0}}\identity\right]\quad\text{ and }\quad\tilde{g}\left(x\right):=g\left(-x\right)
\]
forms a Banach frame for $M_{\left(s,\mu\right),\alpha}^{p,q}\left(\R^{\dimension}\right)$
for all $\left|s\right|\leq s_{0}$, $\left|\mu\right|\leq\mu_{0}$
and all $p,q\in\left(0,\infty\right]$ with $p\geq p_{0}$ and $q\geq q_{0}$.

Precisely, this means the following: Define the coefficient space
\[
\mathscr{C}_{p,q,s,\mu}^{\left(\alpha\right)}:=\ell_{\left[\left|i\right|^{\frac{1}{1-\alpha}\left(s+\alpha\dimension\left(\frac{1}{2}-\frac{1}{p}\right)\right)}\right]_{i\in\Z^{\dimension}\setminus\left\{ 0\right\} }}^{q}\!\!\!\!\!\!\!\!\left(\left[\ell_{\left[\left(1+\left|k\right|/\left|i\right|^{\alpha_{0}}\right)^{\mu}\right]_{k\in\Z^{\dimension}}}^{p}\left(\Z^{\dimension}\right)\right]_{i\in\Z^{\dimension}\setminus\left\{ 0\right\} }\right).
\]
Then the following hold:

\begin{enumerate}
\item The \textbf{analysis map}
\[
A^{\left(\delta\right)}:M_{\left(s,\mu\right),\alpha}^{p,q}\left(\smash{\R^{\dimension}}\right)\to\mathscr{C}_{p,q,s,\mu}^{\left(\alpha\right)},f\mapsto\left[\left(\gamma^{\left[i\right]}\ast f\right)\left(\delta\cdot k/\left|i\right|^{\alpha_{0}}\right)\right]_{i\in\Z^{\dimension}\setminus\left\{ 0\right\} ,k\in\Z^{\dimension}}
\]
is well-defined and bounded for all $0<\delta\leq1$. Here, the convolution
$\left(\gamma^{\left[i\right]}\ast f\right)\left(x\right)$ has to
be understood similar to equation (\ref{eq:SpecialConvolutionPointwiseDefinition}).
\item For $0<\delta\leq\delta_{0}$, there is a bounded linear map \textbf{reconstruction
map} $R^{\left(\delta\right)}:\mathscr{C}_{p,q,s,\mu}^{\left(\alpha\right)}\to M_{\left(s,\mu\right),\alpha}^{p,q}\left(\R^{\dimension}\right)$
satisfying $R^{\left(\delta\right)}\circ A^{\left(\delta\right)}=\identity_{M_{\left(s,\mu\right),\alpha}^{p,q}\left(\R^{\dimension}\right)}$.
Furthermore, the action of $R^{\left(\delta\right)}$ on a given sequence
is independent of the precise choice of $p,q,s,\mu$.
\item We have the following \textbf{consistency statement}: If $f\in M_{\left(s,\mu\right),\alpha}^{p,q}\left(\R^{\dimension}\right)$
and if $p_{0}\leq\tilde{p}\leq\infty$ and $q_{0}\leq\tilde{q}\leq\infty$
and if furthermore $\left|\tilde{s}\right|\leq s_{0}$ and $\left|\tilde{\mu}\right|\leq\mu_{0}$,
then the following equivalence holds:
\[
f\in M_{\left(\tilde{s},\tilde{\mu}\right),\alpha}^{\tilde{p},\tilde{q}}\left(\smash{\R^{\dimension}}\right)\qquad\Longleftrightarrow\qquad A^{\left(\delta\right)}f\in\mathscr{C}_{\tilde{p},\tilde{q},\tilde{s},\tilde{\mu}}^{\left(\alpha\right)}.\qedhere
\]
\end{enumerate}
\end{thm}
\begin{proof}
First of all, we remark that it is comparatively easy to show that
the family $\Gamma^{\left(\delta\right)}$ forms a Banach frame for
$M_{\left(s,\mu\right),\alpha}^{p,q}\left(\R^{\dimension}\right)$
if $0<\delta\leq\delta_{0}$, where $\delta_{0}$ might depend on
$p,q,s,\mu$. About half of the proof will be spent on showing that
$\delta_{0}$ can actually be chosen \emph{independently} of $p,q,s,\mu$,
as long as these satisfy the restrictions mentioned in the statement
of the theorem.

Recall from Lemma \ref{lem:AlphaModulationStructuredAndBAPU} that
there is a family $\Phi=\left(\varphi_{i}\right)_{i\in\Z^{\dimension}\setminus\left\{ 0\right\} }$
associated to $\CalQ=\CalQ_{r}^{\left(\alpha\right)}$ satisfying
Assumption \ref{assu:RegularPartitionOfUnity}. Furthermore, Corollary
\ref{cor:RegularBAPUsAreWeightedBAPUs} yields a function $\varrho\in\TestFunctionSpace{\R^{\dimension}}$
such that, for $v_{0}\left(x\right)=\left[2\cdot\left(1+\left|x\right|\right)\right]^{\left|\mu\right|}$
as in Lemma \ref{lem:AlphaModulationSpaceWeight}, with $K=\left|\mu\right|\leq\mu_{0}$
and $Q=B_{r}\left(0\right)$, as well as $p\geq p_{0}$, we have
\begin{align}
C_{\CalQ_{r}^{\left(\alpha\right)},\Phi,v_{0},p} & \leq\Omega_{0}^{K}\Omega_{1}\cdot\left(4\cdot\dimension\right)^{1+2\left\lceil K+\frac{\dimension+\varepsilon}{p}\right\rceil }\cdot\left(\frac{s_{\dimension}}{\varepsilon}\right)^{1/p}\cdot2^{\!\left\lceil K+\frac{\dimension+\varepsilon}{p}\right\rceil }\cdot\lambda_{\dimension}\left(Q\right)\cdot\max_{\left|\beta\right|\leq\left\lceil K+\frac{\dimension+\varepsilon}{p}\right\rceil }\left\Vert \partial^{\beta}\varrho\right\Vert _{\sup}\cdot\max_{\left|\beta\right|\leq\left\lceil K+\frac{\dimension+\varepsilon}{p}\right\rceil }C^{\left(\beta\right)}\nonumber \\
 & \leq2^{\mu_{0}}\lambda_{\dimension}\left(Q\right)\!\cdot\!\left(8\cdot\dimension\right)^{1+2\left\lceil \mu_{0}+\frac{\dimension+\varepsilon}{p_{0}}\right\rceil }\!\cdot\!\left(1\!+\!\frac{s_{\dimension}}{\varepsilon}\right)^{\frac{1}{p_{0}}}\!\cdot\!\max_{\left|\beta\right|\leq\left\lceil \mu_{0}+\frac{\dimension+\varepsilon}{p_{0}}\right\rceil }\left\Vert \partial^{\beta}\varrho\right\Vert _{\sup}\cdot\!\max_{\left|\beta\right|\leq\left\lceil \mu_{0}+\frac{\dimension+\varepsilon}{p_{0}}\right\rceil }C^{\left(\beta\right)}=:L_{0}.\label{eq:AlphaModulationFrameBAPUConstantEstimate}
\end{align}

Now, assume that $\gamma$ satisfies all the stated properties and
let $\left|s\right|\leq s_{0}$, $\left|\mu\right|\leq\mu_{0}$ and
$p,q\in\left(0,\infty\right]$ with $p\geq p_{0}$ and $q\geq q_{0}$.
We want to verify the assumptions of Corollary \ref{cor:BanachFrameSimplifiedCriteria}
for the family $\left(\gamma_{i}\right)_{i\in\Z^{\dimension}\setminus\left\{ 0\right\} }$,
with $\gamma_{i}:=\gamma$ for all $i\in\Z^{\dimension}\setminus\left\{ 0\right\} $
and with $\CalQ=\CalQ_{r}^{\left(\alpha\right)}$. To this end, let
$\gamma_{1}^{\left(0\right)}:=\gamma$ and set $n_{i}:=1$, so that
$\gamma_{i}=\gamma=\gamma_{n_{i}}^{\left(0\right)}$ for all $i\in\Z^{\dimension}\setminus\left\{ 0\right\} $.
In the notation of Lemma \ref{lem:GammaCoversOrbitAssumptionSimplified},
we then have $Q^{\left(1\right)}=\bigcup\left\{ Q_{i}'\with i\in\Z^{\dimension}\setminus\left\{ 0\right\} \text{ and }n_{i}=1\right\} =B_{r}\left(0\right)$,
since $Q_{i}=T_{i}\left[B_{r}\left(0\right)\right]+b_{i}$ and thus
$Q_{i}'=B_{r}\left(0\right)$ for all $i\in\Z^{\dimension}\setminus\left\{ 0\right\} $.
In view of part (\ref{enu:AlphaModulationFrameFourierNonVanishing})
of our assumptions, we thus see that Lemma \ref{lem:GammaCoversOrbitAssumptionSimplified}
is applicable, so that the family $\left(\gamma_{i}\right)_{i\in\Z^{\dimension}\setminus\left\{ 0\right\} }$
satisfies Assumption \ref{assu:GammaCoversOrbit}, with $\Omega_{2}^{\left(p,K\right)}\leq L_{1}$
for some constant $L_{1}=L_{1}\left(\gamma,\mu_{0},p_{0},r,\alpha,\dimension\right)>0$,
all $p\geq p_{0}$ and all $K\leq\mu_{0}$. Recall from Lemma \ref{lem:AlphaModulationSpaceWeight}
that we can choose $K=\left|\mu\right|\leq\mu_{0}$ in our present
setting.

Finally, recall from part (\ref{enu:AlphaModulationFrameFourierDecay})
of our assumptions that there is some $L_{2}>0$ (independent of $p,q,s,\mu$)
satisfying $\left|\partial^{\beta}\widehat{\gamma}\left(\xi\right)\right|\leq L_{2}\cdot\left(1+\left|\xi\right|\right)^{-N_{0}}$
for all $\xi\in\R^{\dimension}$ and all $\beta\in\N_{0}^{\dimension}$
with $\left|\beta\right|\leq\left\lceil \mu_{0}+\frac{\dimension+\varepsilon}{p_{0}}\right\rceil $.

With these preparations, we can now verify the prerequisites of Corollary
\ref{cor:BanachFrameSimplifiedCriteria}:

\begin{enumerate}
\item As we have seen at the beginning of this section, the covering $\CalQ=\CalQ_{r}^{\left(\alpha\right)}$,
the weight $v=v^{\left(\mu\right)}$ (with $v_{0}\left(x\right)=\left[2\cdot\left(1+\left|x\right|\right)\right]^{\left|\mu\right|}$,
$\Omega_{0}=1$ and $\Omega_{1}=2^{\left|\mu\right|}\leq2^{\mu_{0}}$,
as well as $K=\left|\mu\right|\leq\mu_{0}$) satisfy all standing
assumptions of Section \ref{subsec:DecompSpaceDefinitionStandingAssumptions}.
Furthermore, the family $\Phi$ satisfies Assumption \ref{assu:RegularPartitionOfUnity}.
\item By our assumptions, we have $\gamma_{i}=\gamma\in L_{\left(1+\left|\mybullet\right|\right)^{\mu_{0}}}^{1}\left(\R^{\dimension}\right)\hookrightarrow L_{\left(1+\left|\mybullet\right|\right)^{K}}^{1}\left(\R^{\dimension}\right)$
and $\widehat{\gamma_{i}}=\widehat{\gamma}\in C^{\infty}\left(\R^{\dimension}\right)$
for all $i\in\Z^{\dimension}\setminus\left\{ 0\right\} $, where all
partial derivatives of $\widehat{\gamma_{i}}=\widehat{\gamma}$ are
polynomially bounded.
\item By our assumptions, we have $\gamma_{i}=\gamma\in C^{1}\left(\R^{\dimension}\right)$
and $\partial_{\ell}\gamma\in L_{\left(1+\left|\mybullet\right|\right)^{\mu_{0}}}^{1}\left(\R^{\dimension}\right)\cap L^{\infty}\left(\R^{\dimension}\right)$
and consequently also $\partial_{\ell}\gamma_{i}=\partial_{\ell}\gamma\in L_{\left(1+\left|\mybullet\right|\right)^{K}}^{1}\left(\R^{\dimension}\right)\cap L^{\infty}\left(\R^{\dimension}\right)$
for all $\ell\in\underline{\dimension}$ and $i\in\Z^{\dimension}\setminus\left\{ 0\right\} $.
\item As seen above, the family $\left(\gamma_{i}\right)_{i\in\Z^{\dimension}\setminus\left\{ 0\right\} }=\left(\gamma\right)_{i\in\Z^{\dimension}\setminus\left\{ 0\right\} }$
satisfies Assumption \ref{assu:GammaCoversOrbit}.
\item Let 
\[
C_{1}:=\sup_{i\in\Z^{\dimension}\setminus\left\{ 0\right\} }\:\sum_{j\in\Z^{\dimension}\setminus\left\{ 0\right\} }\:M_{j,i}\in\left[0,\infty\right]\qquad\text{ and }\qquad C_{2}:=\sup_{j\in\Z^{\dimension}\setminus\left\{ 0\right\} }\:\sum_{i\in\Z^{\dimension}\setminus\left\{ 0\right\} }\:M_{j,i}\in\left[0,\infty\right]
\]
as in Corollary \ref{cor:BanachFrameSimplifiedCriteria}, i.e., with
\[
M_{j,i}:=\left(\frac{w_{j}^{\left(s^{\ast}\right)}}{w_{i}^{\left(s^{\ast}\right)}}\right)^{\tau}\cdot\left(1+\left\Vert T_{j}^{-1}T_{i}\right\Vert \right)^{\sigma}\cdot\max_{\left|\beta\right|\leq1}\left(\left|\det T_{i}\right|^{-1}\cdot\int_{Q_{i}}\max_{\left|\theta\right|\leq N}\left|\left(\partial^{\theta}\widehat{\partial^{\beta}\gamma}\right)\!\!\left(S_{j}^{-1}\xi\right)\right|\d\xi\right)^{\tau},
\]
where, $s^{\ast}=\frac{s}{1-\alpha}$ and, since $K=\left|\mu\right|$,
\begin{align*}
N & =\left\lceil \left|\mu\right|+\frac{\dimension+\varepsilon}{\min\left\{ 1,p\right\} }\right\rceil \,,\\
\tau & =\min\left\{ 1,p,q\right\} ,\\
\sigma & =\tau\cdot\left(\frac{\dimension}{\min\left\{ 1,p\right\} }+\left|\mu\right|+\left\lceil \left|\mu\right|+\frac{\dimension+\varepsilon}{\min\left\{ 1,p\right\} }\right\rceil \right).
\end{align*}

Note $\min\left\{ 1,p\right\} \geq\min\left\{ 1,p_{0}\right\} =p_{0}$
and $\left|\mu\right|\leq\mu_{0}$, so that $N\leq\left\lceil \mu_{0}+\frac{\dimension+\varepsilon}{p_{0}}\right\rceil $.
Hence, we see that equation (\ref{eq:AlphaModulationBanachFrameDecayAssumption})
(with $C=L_{2}$) and hence also equation (\ref{eq:AlphaModulationStandardEstimate})
is satisfied, i.e., we have
\[
M_{j,i}\leq\left(2^{N+1}\pi\cdot L_{2}\right)^{\tau}\cdot M_{j,i}^{\left(0\right)}\qquad\forall i,j\in\Z^{\dimension}\setminus\left\{ 0\right\} ,
\]
with $M_{j,i}^{\left(0\right)}$ as in equation (\ref{eq:AlphaModulationStandardEstimate}).
We now want to apply Lemma \ref{lem:AlphaModulationMainLemma}. To
this end, we have to verify that
\[
N_{0}\geq\dimension+2+\frac{\dimension+1}{\tau}+\max\left\{ \left|s^{\ast}+\dimension\alpha_{0}\right|,\left|s^{\ast}+\left(\dimension-\frac{\sigma}{\tau}\right)\alpha_{0}\right|\right\} .
\]
But we have $\tau=\min\left\{ 1,p,q\right\} \geq\min\left\{ 1,p_{0},q_{0}\right\} =\min\left\{ p_{0},q_{0}\right\} =:\tau_{0}$.
Furthermore, with $s_{0}^{\ast}:=s_{0}/\left(1-\alpha\right)$, we
have $\left|s^{\ast}+\dimension\alpha_{0}\right|\leq s_{0}^{\ast}+\dimension\alpha_{0}=\frac{1}{1-\alpha}\left(s_{0}+\alpha\dimension\right)$
and
\[
\frac{\sigma}{\tau}=\frac{\dimension}{\min\left\{ 1,p\right\} }+\left|\mu\right|+\left\lceil \left|\mu\right|+\frac{\dimension+\varepsilon}{\min\left\{ 1,p\right\} }\right\rceil \geq\frac{\dimension}{\min\left\{ 1,p\right\} }\geq\dimension,
\]
so that 
\begin{align*}
\left|s^{\ast}+\left(\dimension-\frac{\sigma}{\tau}\right)\alpha_{0}\right| & \leq s_{0}^{\ast}+\alpha_{0}\left|\dimension-\frac{\sigma}{\tau}\right|=s_{0}^{\ast}+\alpha_{0}\left(\frac{\sigma}{\tau}-\dimension\right)\\
 & =\frac{1}{1-\alpha}\left[s_{0}+\alpha\left(\frac{\dimension}{\min\left\{ 1,p\right\} }+\left|\mu\right|+\left\lceil \left|\mu\right|+\frac{\dimension+\varepsilon}{\min\left\{ 1,p\right\} }\right\rceil -\dimension\right)\right]\\
 & \leq\frac{1}{1-\alpha}\left[s_{0}+\alpha\left(\frac{\dimension}{p_{0}}-\dimension+\mu_{0}+\left\lceil \mu_{0}+\frac{\dimension+\varepsilon}{p_{0}}\right\rceil \right)\right].
\end{align*}
Hence, our assumptions easily yield
\begin{align*}
N_{0} & =\dimension+2+\frac{\dimension+1}{\tau_{0}}+\frac{1}{1-\alpha}\cdot\max\left\{ s_{0}+\alpha\dimension,\:s_{0}+\alpha\left(\frac{\dimension}{p_{0}}-\dimension+\mu_{0}+\left\lceil \mu_{0}+\frac{\dimension+\varepsilon}{p_{0}}\right\rceil \right)\right\} \\
 & \geq\dimension+2+\frac{\dimension+1}{\tau}+\max\left\{ \left|s^{\ast}+\dimension\alpha_{0}\right|,\left|s^{\ast}+\left(\dimension-\frac{\sigma}{\tau}\right)\alpha_{0}\right|\right\} ,
\end{align*}
as desired.

For brevity, set $L_{3}:=2^{\dimension}\cdot\left(1+r\cdot\left(2+4r\right)^{\alpha_{0}}\right)^{N_{0}+\dimension}$.
Since Lemma \ref{lem:AlphaModulationMainLemma} is applicable, and
since $\tau_{0}\leq\tau\leq1$, we get
\begin{align*}
\quad\quad C_{1}^{1/\tau} & \leq2^{N+1}\pi\cdot L_{2}\cdot\left[\sup_{i\in\Z^{\dimension}\setminus\left\{ 0\right\} }\:\sum_{j\in\Z^{\dimension}\setminus\left\{ 0\right\} }M_{j,i}^{\left(0\right)}\right]^{1/\tau}\\
\quad\quad\quad\quad & \leq2^{N+1}\pi L_{2}\cdot2^{\left|s^{\ast}\right|+\frac{1+\sigma}{\tau}}6^{\frac{\dimension}{\tau}}\cdot\max\left\{ 4^{\alpha_{0}\left(\frac{\sigma}{\tau}+\dimension\right)+\left|s^{\ast}\right|}\cdot12^{N_{0}}\cdot s_{\dimension},\:\left(2+4r\right)^{\left|s^{\ast}\right|+\alpha_{0}\left[\frac{\sigma}{\tau}+\dimension+N_{0}\right]}L_{3}\right\} \\
\quad\quad\quad\quad & \overset{\left(\ast\right)}{\leq}6^{\frac{\dimension}{\tau_{0}}}2^{1+\left\lceil \mu_{0}+\frac{\dimension+\varepsilon}{p_{0}}\right\rceil }\pi L_{2}\cdot2^{\frac{s_{0}}{1-\alpha}+\frac{1}{\tau_{0}}+\frac{\sigma_{0}}{\tau_{0}}}\cdot\max\left\{ \!4^{\frac{1}{1-\alpha}\left[s_{0}+\alpha\left(\frac{\sigma_{0}}{\tau_{0}}+\dimension\right)\right]}\!\cdot12^{N_{0}}s_{\dimension},\:\left(2\!+\!4r\right)^{\frac{1}{1-\alpha}\left[s_{0}+\alpha\left(\frac{\sigma_{0}}{\tau_{0}}+\dimension+N_{0}\right)\right]}\!L_{3}\!\right\} \\
\quad\quad\quad\quad & =:L_{4}.
\end{align*}
Here, we recall $\tau_{0}=\min\left\{ p_{0},q_{0}\right\} $ and observe
that the step marked with $\left(\ast\right)$ used the estimate
\[
\frac{\sigma}{\tau}=\frac{\dimension}{\min\left\{ 1,p\right\} }+\left|\mu\right|+\left\lceil \left|\mu\right|+\frac{\dimension+\varepsilon}{\min\left\{ 1,p\right\} }\right\rceil \leq\frac{\dimension}{p_{0}}+\mu_{0}+\left\lceil \mu_{0}+\frac{\dimension+\varepsilon}{p_{0}}\right\rceil =:\frac{\sigma_{0}}{\tau_{0}}
\]
and $\sigma=\frac{\sigma}{\tau}\cdot\tau\leq\frac{\sigma}{\tau}\leq\frac{\sigma_{0}}{\tau_{0}}$.
Hence, we have shown $C_{1}^{1/\tau}\leq L_{4}$, where $L_{4}$ is
\emph{independent} of $p,q,s,\mu$. Exactly the same argument also
shows $C_{2}^{1/\tau}\leq L_{4}$, with the same constant $L_{4}$.
In particular, $C_{1}<\infty$ and $C_{2}<\infty$, which was the
last part of Corollary \ref{cor:BanachFrameSimplifiedCriteria} that
we needed to verify.

\end{enumerate}
All in all, Corollary \ref{cor:BanachFrameSimplifiedCriteria} shows
that the family $\left(\gamma_{i}\right)_{i\in\Z^{\dimension}\setminus\left\{ 0\right\} }=\left(\gamma\right)_{i\in\Z^{\dimension}\setminus\left\{ 0\right\} }$
satisfies all assumptions of Theorem \ref{thm:DiscreteBanachFrameTheorem}
and also that $\vertiii{\smash{\overrightarrow{A}}}^{\max\left\{ 1,\frac{1}{p}\right\} }\leq2L_{5}^{\left(0\right)}\cdot L_{4}$,
as well as $\vertiii{\smash{\overrightarrow{B}}}^{\max\left\{ 1,\frac{1}{p}\right\} }\leq2L_{5}^{\left(0\right)}\cdot L_{4}$,
where $\overrightarrow{A}$ and $\overrightarrow{B}$ are defined
as in Assumptions \ref{assu:MainAssumptions} and \ref{assu:DiscreteBanachFrameAssumptions}
and where
\begin{align*}
L_{5}^{\left(0\right)} & =\Omega_{0}^{K}\Omega_{1}\cdot\dimension^{1/\min\left\{ 1,p\right\} }\cdot\left(4\cdot\dimension\right)^{1+2\left\lceil K+\frac{\dimension+\varepsilon}{\min\left\{ 1,p\right\} }\right\rceil }\cdot\left(\frac{s_{\dimension}}{\varepsilon}\right)^{1/\min\left\{ 1,p\right\} }\cdot\max_{\left|\beta\right|\leq\left\lceil K+\frac{\dimension+\varepsilon}{\min\left\{ 1,p\right\} }\right\rceil }C^{\left(\beta\right)}\\
 & \leq2^{\mu_{0}}\cdot\dimension^{1/p_{0}}\cdot\left(4\cdot\dimension\right)^{1+2\left\lceil \mu_{0}+\frac{\dimension+\varepsilon}{p_{0}}\right\rceil }\cdot\left[1+\frac{s_{\dimension}}{\varepsilon}\right]^{\frac{1}{p_{0}}}\cdot\max_{\left|\beta\right|\leq\left\lceil \mu_{0}+\frac{\dimension+\varepsilon}{p_{0}}\right\rceil }C^{\left(\beta\right)}=:L_{5},
\end{align*}
where the constants $C^{\left(\beta\right)}=C^{\left(\beta\right)}\left(\Phi\right)$
are defined as in Assumption \ref{assu:RegularPartitionOfUnity}.

Now, Theorem \ref{thm:DiscreteBanachFrameTheorem} shows that the
family $\Gamma^{\left(\delta\right)}$ is a Banach frame for $M_{\left(s,\mu\right),\alpha}^{p,q}\left(\R^{\dimension}\right)=\DecompSp{\CalQ_{r}^{\left(\alpha\right)}}p{\ell_{w^{\left(s^{\ast}\right)}}^{q}}{v^{\left(\mu\right)}}$,
as soon as $0<\delta\leq\frac{1}{1+2\vertiii{F_{0}}^{2}}$ with $F_{0}$
as in Lemma \ref{lem:SpecialProjection}. But that lemma yields the
estimate
\begin{align*}
\vertiii{F_{0}} & \leq2^{\frac{1}{q}}C_{\CalQ_{r}^{\left(\alpha\right)},\Phi,v_{0},p}^{2}\cdot\vertiii{\smash{\Gamma_{\CalQ_{r}^{\left(\alpha\right)}}}}^{2}\cdot\left(\vertiii{\smash{\overrightarrow{A}}}^{\max\left\{ 1,\frac{1}{p}\right\} }+\vertiii{\smash{\overrightarrow{B}}}^{\max\left\{ 1,\frac{1}{p}\right\} }\right)\cdot L_{6}^{\left(0\right)}\\
\left({\scriptstyle \text{eqs. }\eqref{eq:WeightedSequenceSpaceClusteringMapNormEstimate}\text{ and }\eqref{eq:AlphaModulationFrameBAPUConstantEstimate}}\right) & \leq2^{\frac{1}{q_{0}}}\cdot C_{\CalQ_{r}^{\left(\alpha\right)},w^{\left(s^{\ast}\right)}}^{2}\cdot N_{\CalQ_{r}^{\left(\alpha\right)}}^{2\left(1+\frac{1}{q}\right)}\cdot4L_{0}^{2}L_{4}L_{5}\cdot L_{6}^{\left(0\right)}
\end{align*}
for
\[
L_{6}^{\left(0\right)}=\!\!\begin{cases}
\frac{\left(2^{16}\cdot768/\dimension^{\frac{3}{2}}\right)^{\frac{\dimension}{p}}}{2^{42}\cdot12^{\dimension}\cdot\dimension^{15}}\!\cdot\!\left(2^{52}\!\cdot\!\dimension^{\frac{25}{2}}\!\cdot\!N^{3}\right)^{N+1}\!\!\!\cdot\!N_{\CalQ_{r}^{\left(\alpha\right)}}^{2\left(\frac{1}{p}-1\right)}\!\left(1\!+\!R_{\CalQ_{r}^{\left(\alpha\right)}}C_{\CalQ_{r}^{\left(\alpha\right)}}\right)^{\dimension\left(\frac{4}{p}-1\right)}\!\!\cdot\Omega_{0}^{13K}\Omega_{1}^{13}\Omega_{2}^{\left(p,K\right)}, & \text{if }p<1,\\
\frac{1}{\sqrt{\dimension}\cdot2^{12+6\left\lceil K\right\rceil }}\cdot\left(2^{17}\cdot\dimension^{5/2}\cdot N\right)^{\left\lceil K\right\rceil +\dimension+2}\cdot\left(1+R_{\CalQ_{r}^{\left(\alpha\right)}}\right)^{\dimension}\cdot\Omega_{0}^{3K}\Omega_{1}^{3}\Omega_{2}^{\left(p,K\right)}, & \text{if }p\geq1.
\end{cases}
\]
But above we saw $\Omega_{2}^{\left(p,K\right)}\leq L_{1}$ since
$K=\left|\mu\right|\leq\mu_{0}$ and $p\geq p_{0}$. Using this estimate
and the bounds $N\leq\left\lceil \mu_{0}+\frac{\dimension+\varepsilon}{p_{0}}\right\rceil $
and $0\leq K=\left|\mu\right|\leq\mu_{0}$, as well as $\Omega_{0}=1$
and $\Omega_{1}=2^{\left|\mu\right|}\leq2^{\mu_{0}}$, we see $L_{6}^{\left(0\right)}\leq L_{6}$,
where $L_{6}$ is independent of $p,q,s,\mu$.

Finally, it is not hard to see—because of $w^{\left(s^{\ast}\right)}=\left(w^{\left(1/\left(1-\alpha\right)\right)}\right)^{s}$—that
\[
N_{\CalQ_{r}^{\left(\alpha\right)},w^{\left(s^{\ast}\right)}}\leq N_{\CalQ_{r}^{\left(\alpha\right)},w^{\left(1/\left(1-\alpha\right)\right)}}^{\left|s\right|}\leq N_{\CalQ_{r}^{\left(\alpha\right)},w^{\left(1/\left(1-\alpha\right)\right)}}^{s_{0}}=:L_{7},
\]
where $L_{7}>0$ is independent of $p,q,s,\mu$. By putting everything
together, we get $\vertiii{F_{0}}\leq L_{8}$, with $L_{8}$ independent
of $p,q,s,\mu$. Hence, $\Gamma^{\left(\delta\right)}$ is a Banach
frame for $M_{\left(s,\mu\right),\alpha}^{p,q}\left(\R^{\dimension}\right)$
in the sense of Theorem \ref{thm:DiscreteBanachFrameTheorem} as soon
as $0<\delta\leq\delta_{0}:=\frac{1}{1+2L_{8}^{2}}$, with $\delta_{0}$
independent of $p,q,s,\mu$.

\medskip{}

All that remains to verify is that the space $\ell_{\left(\left|\det T_{i}\right|^{\frac{1}{2}-\frac{1}{p}}\cdot w_{i}\right)_{i\in I}}^{q}\!\!\!\!\!\left(\left[\vphantom{\sum}\smash{C_{i}^{\left(\delta\right)}}\right]_{i\in I}\right)$
(with $w_{i}=w_{i}^{\left(s^{\ast}\right)}=\left\langle i\right\rangle ^{s^{\ast}}$)
appearing in Theorem \ref{thm:DiscreteBanachFrameTheorem} coincides
with the space $\mathscr{C}_{p,q,s,\mu}^{\left(\alpha\right)}$ from
the statement of the current theorem. To this end, first note that
$\left|i\right|\asymp\left\langle i\right\rangle $ for all $i\in\Z^{\dimension}\setminus\left\{ 0\right\} $,
so that
\[
\left|\det T_{i}\right|^{\frac{1}{2}-\frac{1}{p}}\cdot w_{i}=\left|i\right|^{\dimension\alpha_{0}\left(\frac{1}{2}-\frac{1}{p}\right)}\cdot\left\langle i\right\rangle ^{s^{\ast}}\asymp_{s}\;\left|i\right|^{\frac{1}{1-\alpha}\left[s+\alpha\dimension\left(\frac{1}{2}-\frac{1}{p}\right)\right]}\qquad\forall i\in\Z^{\dimension}\setminus\left\{ 0\right\} .
\]
Finally, recall from equation (\ref{eq:CoefficientSpaceDefinition})
that
\[
C_{i}^{\left(\delta\right)}=\ell_{v^{\left(j,\delta\right)}}^{p}\left(\smash{\Z^{\dimension}}\right)\quad\text{ with }\quad v_{k}^{\left(j,\delta\right)}=v\left(\delta\cdot T_{j}^{-T}k\right)
\]
where $v=v^{\left(\mu\right)}$ with $v^{\left(\mu\right)}\left(x\right)=\left\langle x\right\rangle ^{\mu}\asymp_{\mu}\;\left(1+\left|x\right|\right)^{\mu}$
for all $x\in\R^{\dimension}$. Furthermore, since $0<\delta\leq1$,
we have 
\[
\delta\cdot\left(1+\left|x\right|\right)\leq1+\left|\delta\cdot x\right|\leq1+\left|x\right|\qquad\forall x\in\R^{\dimension},
\]
which yields
\[
v_{k}^{\left(j,\delta\right)}=v\left(\delta\cdot T_{j}^{-T}k\right)\:\asymp_{\mu}\:\left(1+\left|\delta\cdot T_{j}^{-T}k\right|\right)^{\mu}\:\asymp_{\mu,\delta}\:\left(1+\left|T_{j}^{-T}k\right|\right)^{\mu}=\left(1+\left|k\right|/\left|j\right|^{\alpha_{0}}\right)^{\mu}
\]
for all $k\in\Z^{\dimension}$ and $j\in\Z^{\dimension}\setminus\left\{ 0\right\} $.
Here, the implied constant might depend on $\mu,\delta$, but not
on $j,k$. Combining these facts, we conclude $\ell_{\left(\left|\det T_{i}\right|^{\frac{1}{2}-\frac{1}{p}}\cdot w_{i}\right)_{i\in I}}^{q}\!\!\!\!\!\left(\left[\vphantom{\sum}\smash{C_{i}^{\left(\delta\right)}}\right]_{i\in I}\right)=\mathscr{C}_{p,q,s,\mu}^{\left(\alpha\right)}$,
with equivalent quasi-norms, as desired. Now all claims follow from
Theorem \ref{thm:DiscreteBanachFrameTheorem}.
\end{proof}
Having established convenient criteria for the existence of Banach
frames, we finally consider nice criteria which ensure that a prototype
$\gamma$ generates an atomic decomposition for a given $\alpha$-modulation
space.
\begin{thm}
\label{thm:AlphaModulationAtomicDecomposition}Let $\dimension\in\N$,
$\alpha\in\left[0,1\right)$ and choose $r>r_{1}\left(\dimension,\alpha\right)$
with $r_{1}\left(\dimension,\alpha\right)$ as in Theorem \ref{thm:AlphaModulationCoveringDefinition}.

Let $s_{0},\mu_{0}\geq0$ and $p_{0},q_{0}\in\left(0,1\right]$, as
well as $\varepsilon\in\left(0,1\right)$. Assume that $\gamma:\R^{\dimension}\to\Compl$
is measurable and satisfies the following conditions:

\begin{enumerate}
\item We have $\left\Vert \gamma\right\Vert _{\mu_{0}+\frac{\dimension}{p_{0}}+1}<\infty$,
where as usual $\left\Vert g\right\Vert _{M}=\sup_{x\in\R^{\dimension}}\left(1+\left|x\right|\right)^{M}\left|g\left(x\right)\right|$.
In particular, $\gamma\in L^{1}\left(\R^{\dimension}\right)$.
\item We have $\widehat{\gamma}\in C^{\infty}\left(\R^{\dimension}\right)$
and all partial derivatives of $\widehat{\gamma}$ are polynomially
bounded.
\item \label{enu:AlphaModulationAtomicFourierNonVanishing}We have $\left|\widehat{\gamma}\left(\xi\right)\right|\geq c>0$
for all $\xi\in\overline{B_{r}}\left(0\right)$.
\item \label{enu:AlphaModulationAtomicFourierDecay}We have
\[
\left|\left(\partial^{\beta}\widehat{\gamma}\right)\left(\xi\right)\right|\lesssim\left(1+\left|\xi\right|\right)^{-M_{0}}
\]
for all $\xi\in\R^{\dimension}$ and all $\beta\in\N_{0}^{\dimension}$
with $\left|\beta\right|\leq\left\lceil \mu_{0}+\frac{\dimension+\varepsilon}{p_{0}}\right\rceil $,
where
\[
M_{0}=\left(\dimension+1\right)\cdot\left(2+\varepsilon+\frac{1}{\min\left\{ p_{0},q_{0}\right\} }\right)+\Lambda,
\]
with
\[
\Lambda=\begin{cases}
\frac{1}{1-\alpha}\max\left\{ s_{0}+\dimension\alpha,\,s_{0}+\alpha\left(\left\lceil \mu_{0}+\dimension+\varepsilon\right\rceil -\dimension\right)\right\} , & \text{if }p_{0}=1,\\
\frac{1}{1-\alpha}\left[s_{0}+\alpha\left(\mu_{0}+\left\lceil \mu_{0}+p_{0}^{-1}\cdot\left(\dimension+\varepsilon\right)\right\rceil \right)\right], & \text{if }p_{0}\in\left(0,1\right).
\end{cases}
\]
\end{enumerate}
Then there is some $\delta_{0}>0$ such that for all $0<\delta\leq\delta_{0}$,
the family 
\[
\Gamma^{\left(\delta\right)}:=\left(L_{\delta\cdot k/\left|i\right|^{\alpha_{0}}}\:\gamma^{\left[i\right]}\right)_{i\in\Z^{\dimension}\setminus\left\{ 0\right\} ,k\in\Z^{\dimension}},\qquad\text{ with }\qquad\gamma^{\left[i\right]}=\left|i\right|^{\frac{\dimension\cdot\alpha_{0}}{2}}\cdot M_{\left|i\right|^{\alpha_{0}}\cdot i}\left[\gamma\circ\left|i\right|^{\alpha_{0}}\identity\right]
\]
forms an atomic decomposition for $M_{\left(s,\mu\right),\alpha}^{p,q}\left(\R^{\dimension}\right)$
for all $\left|s\right|\leq s_{0}$, $\left|\mu\right|\leq\mu_{0}$
and all $p,q\in\left(0,\infty\right]$ with $p\geq p_{0}$ and $q\geq q_{0}$.

Precisely, this means the following: Define the coefficient space
\[
\mathscr{C}_{p,q,s,\mu}^{\left(\alpha\right)}:=\ell_{\left[\left|i\right|^{\frac{1}{1-\alpha}\left(s+\alpha\dimension\left(\frac{1}{2}-\frac{1}{p}\right)\right)}\right]_{i\in\Z^{\dimension}\setminus\left\{ 0\right\} }}^{q}\!\!\!\!\!\!\!\!\left(\left[\ell_{\left[\left(1+\left|k\right|/\left|i\right|^{\alpha_{0}}\right)^{\mu}\right]_{k\in\Z^{\dimension}}}^{p}\left(\Z^{\dimension}\right)\right]_{i\in\Z^{\dimension}\setminus\left\{ 0\right\} }\right).
\]
Then the following hold:

\begin{enumerate}
\item For arbitrary $\delta\in\left(0,1\right]$, the \textbf{synthesis
map}
\[
S^{\left(\delta\right)}:\mathscr{C}_{p,q,s,\mu}^{\left(\alpha\right)}\to M_{\left(s,\mu\right),\alpha}^{p,q}\left(\smash{\R^{\dimension}}\right),\left(\smash{c_{k}^{\left(i\right)}}\right)_{i\in\Z^{\dimension}\setminus\left\{ 0\right\} ,k\in\Z^{\dimension}}\mapsto\sum_{i\in\Z^{\dimension}\setminus\left\{ 0\right\} }\:\sum_{k\in\Z^{\dimension}}\left[c_{k}^{\left(i\right)}\cdot L_{\delta\cdot k/\left|i\right|^{\alpha_{0}}}\:\gamma^{\left[i\right]}\right]
\]
is well-defined and bounded. Convergence of the series has to be understood
as described in the remark following Theorem \ref{thm:AtomicDecomposition}.
\item For $0<\delta\leq\delta_{0}$, there is a bounded linear \textbf{coefficient
map} $C^{\left(\delta\right)}:M_{\left(s,\mu\right),\alpha}^{p,q}\left(\smash{\R^{\dimension}}\right)\to\mathscr{C}_{p,q,s,\mu}^{\left(\alpha\right)}$
satisfying $S^{\left(\delta\right)}\circ C^{\left(\delta\right)}=\identity_{M_{\left(s,\mu\right),\alpha}^{p,q}\left(\smash{\R^{\dimension}}\right)}$.
Furthermore, the action of $C^{\left(\delta\right)}$ on a given $f\in M_{\left(s,\mu\right),\alpha}^{p,q}\left(\R^{\dimension}\right)$
is independent of the precise choice of $p,q,s,\mu$.\qedhere
\end{enumerate}
\end{thm}
\begin{rem*}
Choose $M_{0}$ as in the theorem above. If $\gamma\in C_{c}^{\left\lceil M_{0}\right\rceil }\left(\R^{\dimension}\right)$,
then $\widehat{\gamma}\in C^{\infty}\left(\R^{\dimension}\right)$
with all partial derivatives being polynomially bounded. Furthermore,
for arbitrary $\alpha\in\N_{0}^{\dimension}$, we have $\gamma_{\alpha}\in C_{c}^{\left\lceil M_{0}\right\rceil }\left(\R^{\dimension}\right)\hookrightarrow W^{\left\lceil M_{0}\right\rceil ,1}\left(\R^{\dimension}\right)$
for
\[
\gamma_{\alpha}:\R^{\dimension}\to\Compl,x\mapsto\left(-2\pi ix\right)^{\alpha}\cdot\gamma\left(x\right).
\]
But by differentiation under the integral, it is not hard to see $\partial^{\alpha}\widehat{\gamma}\left(\xi\right)=\widehat{\gamma_{\alpha}}\left(\xi\right)$
for all $\xi\in\R^{\dimension}$. Hence, Lemma \ref{lem:PointwiseFourierDecayEstimate}
yields $\left|\partial^{\alpha}\widehat{\gamma}\left(\xi\right)\right|=\left|\left(\Fourier^{-1}\gamma_{\alpha}\right)\left(-\xi\right)\right|\lesssim\left(1+\left|\xi\right|\right)^{-\left\lceil M_{0}\right\rceil }\leq\left(1+\left|\xi\right|\right)^{-M_{0}}$.
Finally, we clearly have $\left\Vert \gamma\right\Vert _{\mu_{0}+\frac{\dimension}{p_{0}}+1}<\infty$,
since $\gamma$ has compact support.

All in all, these considerations show that every prototype $\gamma\in C_{c}^{\left\lceil M_{0}\right\rceil }\left(\R^{\dimension}\right)$
with $\widehat{\gamma}\left(\xi\right)\neq0$ for all $\xi\in\overline{B_{r}}\left(0\right)$
generates an atomic decomposition for $M_{\left(s,\mu\right),\alpha}^{p,q}\left(\R^{\dimension}\right)$,
where $M_{0}=M_{0}\left(\dimension,p,q,s,\mu,\alpha,\varepsilon\right)$
has to be chosen suitably. Very similar considerations apply for the
case of Banach frames: Here, it suffices to have $\gamma\in C_{c}^{\left\lceil N_{0}\right\rceil }\left(\R^{\dimension}\right)$
with $\widehat{\gamma}\left(\xi\right)\neq0$ for all $\xi\in\overline{B_{r}}\left(0\right)$,
where $N_{0}=N_{0}\left(\dimension,p,q,s,\mu,\alpha,\varepsilon\right)$
is chosen as in Theorem \ref{thm:AlphaModulationBanachFrame}.
\end{rem*}
\begin{proof}
First of all, we remark as in the proof of Theorem \ref{thm:AlphaModulationBanachFrame}
that it is comparatively easy to show that the family $\Gamma^{\left(\delta\right)}$
forms an atomic decomposition for $M_{\left(s,\mu\right),\alpha}^{p,q}\left(\R^{\dimension}\right)$
if $0<\delta\leq\delta_{0}$, where $\delta_{0}$ might depend on
$p,q,s,\mu$. About half of the proof will be spent on showing that
$\delta_{0}$ can actually be chosen \emph{independently} of $p,q,s,\mu$,
as long as these satisfy the restrictions mentioned in the statement
of the theorem.

Our assumptions yield $L_{0}>0$ satisfying $\left|\partial^{\beta}\widehat{\gamma}\left(\xi\right)\right|\leq L_{0}\cdot\left(1+\left|\xi\right|\right)^{-M_{0}}$
for all $\xi\in\R^{\dimension}$ and all multiindices $\beta\in\N_{0}^{\dimension}$
with $\left|\beta\right|\leq\left\lceil \mu_{0}+\frac{\dimension+\varepsilon}{p_{0}}\right\rceil =:N_{0}$.
As our first step, we invoke Lemma \ref{lem:ConvolutionFactorization}
with $N=N_{0}$ and 
\[
\varrho:\R^{\dimension}\to\left(0,\infty\right),\xi\mapsto L_{0}\cdot\left(1+\left|\xi\right|\right)^{-\left[M_{0}-\left(\dimension+1\right)\left(1+\varepsilon\right)\right]}.
\]
To this end, we observe $N_{0}\geq\left\lceil \left(\dimension+\varepsilon\right)/p_{0}\right\rceil \geq\left\lceil \dimension+\varepsilon\right\rceil =\dimension+1$
and furthermore 
\[
M_{0}-\left(\dimension+1\right)\left(1+\varepsilon\right)\geq\left(\dimension+1\right)\cdot\left(2+\varepsilon\right)-\left(\dimension+1\right)\left(1+\varepsilon\right)=\dimension+1>\dimension,
\]
so that $\varrho\in L^{1}\left(\R^{\dimension}\right)$. Finally,
as we just saw, we indeed have 
\[
\left|\partial^{\beta}\widehat{\gamma}\left(\xi\right)\right|\leq L_{0}\cdot\left(1+\left|\xi\right|\right)^{-M_{0}}=\varrho\left(\xi\right)\cdot\left(1+\left|\xi\right|\right)^{-\left(\dimension+1\right)\left(1+\varepsilon\right)}\leq\varrho\left(\xi\right)\cdot\left(1+\left|\xi\right|\right)^{-\left(\dimension+1+\varepsilon\right)}
\]
for all $\xi\in\R^{\dimension}$ and all $\beta\in\N_{0}^{\dimension}$
with $\left|\beta\right|\leq N_{0}$, so that all assumptions of Lemma
\ref{lem:ConvolutionFactorization} are satisfied. Hence, there are
functions $\gamma_{1},\gamma_{2}\in L^{1}\left(\R^{\dimension}\right)$
with the following properties:

\begin{enumerate}
\item We have $\gamma=\gamma_{1}\ast\gamma_{2}$.
\item We have $\gamma_{2}\in C^{1}\left(\R^{\dimension}\right)$ with $L_{1}^{\left(M\right)}:=\left\Vert \gamma_{2}\right\Vert _{M}+\left\Vert \nabla\gamma_{2}\right\Vert _{M}<\infty$
for arbitrary $M\in\N_{0}$.
\item We have $\widehat{\gamma_{1}},\widehat{\gamma_{2}}\in C^{\infty}\left(\R^{\dimension}\right)$,
where all partial derivatives of these functions are polynomially
bounded.
\item We have $\left\Vert \gamma_{1}\right\Vert _{N_{0}}<\infty$ and $\left\Vert \gamma\right\Vert _{N_{0}}<\infty$.
\item We have $\left|\partial^{\beta}\widehat{\gamma_{1}}\left(\xi\right)\right|\leq L_{2}\cdot\varrho\left(\xi\right)=L_{0}L_{2}\cdot\left(1+\left|\xi\right|\right)^{-M_{00}}$
for all $\xi\in\R^{\dimension}$ and all $\beta\in\N_{0}^{\dimension}$
with $\left|\beta\right|\leq N_{0}$. Here, $L_{2}:=2^{1+\dimension+4N_{0}}\cdot N_{0}!\cdot\left(1+\dimension\right)^{N_{0}}$
and 
\[
M_{00}:=M_{0}-\left(\dimension+1\right)\left(1+\varepsilon\right)=\left(\dimension+1\right)\cdot\left(1+\frac{1}{\min\left\{ p_{0},q_{0}\right\} }\right)+\Lambda.
\]
\end{enumerate}
Next, recall from Lemma \ref{lem:AlphaModulationStructuredAndBAPU}
that there is a family $\Phi=\left(\varphi_{i}\right)_{i\in\Z^{\dimension}\setminus\left\{ 0\right\} }$
associated to $\CalQ=\CalQ_{r}^{\left(\alpha\right)}$ satisfying
Assumption \ref{assu:RegularPartitionOfUnity}. As in the proof of
Theorem \ref{thm:AlphaModulationBanachFrame} (cf.\@ equation (\ref{eq:AlphaModulationFrameBAPUConstantEstimate})),
we get as a consequence of Corollary \ref{cor:RegularBAPUsAreWeightedBAPUs}
a constant $L_{3}>0$ satisfying 
\begin{equation}
C_{\CalQ_{r}^{\left(\alpha\right)},\Phi,v_{0},p}\leq L_{3},\quad\forall p\geq p_{0}\,\forall K=\left|\mu\right|\leq\mu_{0},\label{eq:AlphaModulationAtomicDecompositionBAPUConstantEstimate}
\end{equation}
where $v_{0}\left(x\right)=\left[2\cdot\left(1+\left|x\right|\right)\right]^{\left|\mu\right|}$
is as in Lemma \ref{lem:AlphaModulationSpaceWeight}.

Now, let $p,q,s,\mu$ as in the statement of the Theorem. We want
to show that the family $\left(\gamma_{i}\right)_{i\in\Z^{\dimension}\setminus\left\{ 0\right\} }$
with $\gamma_{i}:=\gamma$ for all $i\in\Z^{\dimension}\setminus\left\{ 0\right\} $
satisfies all assumptions of Corollary \ref{cor:AtomicDecompositionSimplifiedCriteria},
for $\CalQ=\CalQ_{r}^{\left(\alpha\right)}$. To this end, let $\gamma_{1}^{\left(0\right)}:=\gamma$
and $n_{i}:=1$, so that $\gamma_{i}=\gamma=\gamma_{n_{i}}^{\left(0\right)}$
for all $i\in\Z^{\dimension}\setminus\left\{ 0\right\} $. As in the
proof of Theorem \ref{thm:AlphaModulationBanachFrame}, we then see
that all assumptions of Lemma \ref{lem:GammaCoversOrbitAssumptionSimplified}
are satisfied. Hence, the family $\left(\gamma_{i}\right)_{i\in\Z^{\dimension}\setminus\left\{ 0\right\} }$
satisfies Assumption \ref{assu:GammaCoversOrbit} and we also get
$\Omega_{2}^{\left(p,K\right)}\leq L_{4}$ for some constant $L_{4}=L_{4}\left(\gamma,\mu_{0},p_{0},r,\alpha,\dimension\right)>0$,
all $p\geq p_{0}$ and all $K\leq\mu_{0}$.

Set $\gamma_{i,1}:=\gamma_{1}$ and $\gamma_{i,2}:=\gamma_{2}$ for
$i\in\Z^{\dimension}\setminus\left\{ 0\right\} $. Let us verify the
assumptions of Corollary \ref{cor:AtomicDecompositionSimplifiedCriteria}
for these choices:

\begin{enumerate}
\item As we have seen at the beginning of this section, the covering $\CalQ=\CalQ_{r}^{\left(\alpha\right)}$,
the weight $v=v^{\left(\mu\right)}$ (with $v_{0}\left(x\right)=\left[2\cdot\left(1+\left|x\right|\right)\right]^{\left|\mu\right|}$,
$\Omega_{0}=1$ and $\Omega_{1}=2^{\left|\mu\right|}\leq2^{\mu_{0}}$,
as well as $K=\left|\mu\right|\leq\mu_{0}$) satisfy all standing
assumptions of Section \ref{subsec:DecompSpaceDefinitionStandingAssumptions}.
Furthermore, the family $\Phi=\left(\varphi_{i}\right)_{i\in I}$
from above satisfies Assumption \ref{assu:RegularPartitionOfUnity}.
\item Our choice of $\gamma_{1},\gamma_{2}$ from above ensures that all
$\gamma_{i}=\gamma$, $\gamma_{i,1}=\gamma_{1}$ and $\gamma_{i,2}=\gamma_{2}$
are measurable functions.
\item As seen above, we have $\left\Vert \gamma_{1}\right\Vert _{N_{0}}<\infty$.
Since $N_{0}=\left\lceil \mu_{0}+\frac{\dimension+\varepsilon}{p_{0}}\right\rceil \geq\mu_{0}+\frac{\dimension+\varepsilon}{p_{0}}\geq\mu_{0}+\dimension+\varepsilon$,
equation (\ref{eq:StandardDecayLpEstimate}) easily yields $\gamma_{i,1}=\gamma_{1}\in L_{\left(1+\left|\mybullet\right|\right)^{\mu_{0}}}^{1}\left(\R^{\dimension}\right)\hookrightarrow L_{\left(1+\left|\mybullet\right|\right)^{K}}^{1}\left(\R^{\dimension}\right)$
for all $i\in\Z^{\dimension}\setminus\left\{ 0\right\} $.
\item As seen above, we have $\gamma_{i,2}=\gamma_{2}\in C^{1}\left(\R^{\dimension}\right)$
for all $i\in\Z^{\dimension}\setminus\left\{ 0\right\} $.
\item With $K_{0}:=K+\frac{\dimension}{\min\left\{ 1,p\right\} }+1$, we
have
\begin{align*}
\Omega_{4}^{\left(p,K\right)} & =\sup_{i\in\Z^{\dimension}\setminus\left\{ 0\right\} }\left\Vert \gamma_{i,2}\right\Vert _{K_{0}}+\sup_{i\in\Z^{\dimension}\setminus\left\{ 0\right\} }\left\Vert \nabla\gamma_{i,2}\right\Vert _{K_{0}}\\
\left({\scriptstyle \text{since }\gamma_{i,2}=\gamma_{2}\text{ and }K_{0}\leq\mu_{0}+\frac{\dimension}{p_{0}}+1}\right) & \leq2\cdot\sup_{i\in\Z^{\dimension}\setminus\left\{ 0\right\} }\left(\left\Vert \gamma_{2}\right\Vert _{\left\lceil \mu_{0}+\frac{\dimension}{p_{0}}+1\right\rceil }+\left\Vert \nabla\gamma_{2}\right\Vert _{\left\lceil \mu_{0}+\frac{\dimension}{p_{0}}+1\right\rceil }\right)\\
 & =2\cdot L_{1}^{\left(\left\lceil \mu_{0}+\frac{\dimension}{p_{0}}+1\right\rceil \right)}=:L_{5}<\infty.
\end{align*}
\item By our assumptions, we have $\left\Vert \gamma_{i}\right\Vert _{K_{0}}=\left\Vert \gamma\right\Vert _{K_{0}}\leq\left\Vert \gamma\right\Vert _{\mu_{0}+\frac{\dimension}{p_{0}}+1}<\infty$
for all $i\in\Z^{\dimension}\setminus\left\{ 0\right\} $.
\item By choice of $\gamma_{1},\gamma_{2}$, we have $\gamma_{i}=\gamma=\gamma_{1}\ast\gamma_{2}=\gamma_{i,1}\ast\gamma_{i,2}$
for all $i\in\Z^{\dimension}\setminus\left\{ 0\right\} $.
\item By the properties of $\gamma_{1},\gamma_{2}$ from above, we have
$\widehat{\gamma_{i,\ell}}=\widehat{\gamma_{\ell}}$ and all partial
derivatives of this function are polynomially bounded, for arbitrary
$i\in\Z^{\dimension}\setminus\left\{ 0\right\} $ and $\ell\in\left\{ 1,2\right\} $.
\item As seen above, the family $\left(\gamma_{i}\right)_{i\in\Z^{\dimension}\setminus\left\{ 0\right\} }$
satisfies Assumption \ref{assu:GammaCoversOrbit}.
\item Let 
\[
K_{1}:=\sup_{i\in\Z^{\dimension}\setminus\left\{ 0\right\} }\:\sum_{j\in\Z^{\dimension}\setminus\left\{ 0\right\} }\:N_{i,j}\in\left[0,\infty\right]\qquad\text{ and }\qquad K_{2}:=\sup_{j\in\Z^{\dimension}\setminus\left\{ 0\right\} }\:\sum_{i\in\Z^{\dimension}\setminus\left\{ 0\right\} }\:N_{i,j}\in\left[0,\infty\right]
\]
as in Corollary \ref{cor:AtomicDecompositionSimplifiedCriteria},
i.e., with
\[
\qquad N_{i,j}=\left[\frac{w_{i}^{\left(s^{\ast}\right)}}{w_{j}^{\left(s^{\ast}\right)}}\cdot\left(\left|\det T_{j}\right|/\left|\det T_{i}\right|\right)^{\vartheta}\right]^{\tau}\cdot\left(1+\left\Vert T_{j}^{-1}T_{i}\right\Vert \right)^{\sigma}\cdot\left(\left|\det T_{i}\right|^{-1}\cdot\int_{Q_{i}}\max_{\left|\beta\right|\leq N}\left|\left(\partial^{\beta}\widehat{\gamma_{1}}\right)\left(S_{j}^{-1}\xi\right)\right|\d\xi\right)^{\tau},
\]
where $s^{\ast}=\frac{s}{1-\alpha}$ and, since $K=\left|\mu\right|$,
\begin{align*}
N & =\left\lceil \left|\mu\right|+\frac{\dimension+\varepsilon}{\min\left\{ 1,p\right\} }\right\rceil \,,\\
\tau & =\min\left\{ 1,p,q\right\} ,\\
\vartheta & =\begin{cases}
0, & \text{if }p\in\left[1,\infty\right],\\
\frac{1}{p}-1, & \text{if }p\in\left(0,1\right),
\end{cases}\\
\sigma & =\begin{cases}
\tau\cdot\left\lceil \left|\mu\right|+\dimension+\varepsilon\right\rceil , & \text{if }p\in\left[1,\infty\right],\\
\tau\cdot\left(\frac{\dimension}{p}+\left|\mu\right|+\left\lceil \left|\mu\right|+\frac{\dimension+\varepsilon}{p}\right\rceil \right), & \text{if }p\in\left(0,1\right).
\end{cases}
\end{align*}
Note that $\vartheta\geq0$ and $\tau>0$. Furthermore, since $\alpha_{0}\geq0$
and since $\left|j\right|\geq1$ for all $j\in\Z^{\dimension}\setminus\left\{ 0\right\} $,
we have 
\[
\left(2\left\langle j\right\rangle \right)^{\dimension\alpha_{0}}\geq\left\langle j\right\rangle ^{\dimension\alpha_{0}}\geq\left|\det T_{j}\right|=\left|j\right|^{\dimension\alpha_{0}}\geq\left(\frac{1}{2}\left\langle j\right\rangle \right)^{\dimension\alpha_{0}}
\]
and thus
\[
\qquad\quad\left[\frac{w_{i}^{\left(s^{\ast}\right)}}{w_{j}^{\left(s^{\ast}\right)}}\cdot\left(\frac{\left|\det T_{j}\right|}{\left|\det T_{i}\right|}\right)^{\!\vartheta}\,\right]^{\tau}\leq\left[\frac{w_{i}^{\left(s^{\ast}\right)}}{w_{j}^{\left(s^{\ast}\right)}}\cdot\left(\frac{2\left\langle j\right\rangle }{\frac{1}{2}\left\langle i\right\rangle }\right)^{\!\vartheta\dimension\alpha_{0}}\right]^{\tau}\leq4^{\tau\vartheta\dimension\alpha_{0}}\!\cdot\!\left[\frac{w_{i}^{\left(s^{\ast}-\vartheta\dimension\alpha_{0}\right)}}{w_{j}^{\left(s^{\ast}-\vartheta\dimension\alpha_{0}\right)}}\right]^{\tau}=4^{\tau\vartheta\dimension\alpha_{0}}\!\cdot\!\left[\frac{w_{j}^{\left(\vartheta\dimension\alpha_{0}-s^{\ast}\right)}}{w_{i}^{\left(\vartheta\dimension\alpha_{0}-s^{\ast}\right)}}\right]^{\tau}
\]
for all $i,j\in\Z^{\dimension}\setminus\left\{ 0\right\} $. Moreover,
$N\leq\left\lceil \mu_{0}+\frac{\dimension+\varepsilon}{p_{0}}\right\rceil =N_{0}$
and hence $\left|\left(\partial^{\beta}\widehat{\gamma_{1}}\right)\left(\xi\right)\right|\leq L_{0}L_{2}\cdot\left(1+\left|\xi\right|\right)^{-M_{00}}$
for all $\left|\beta\right|\leq N$. Combining these estimates and
noting $\vartheta\leq\frac{1}{p_{0}}$, we arrive at
\[
\qquad N_{i,j}\leq\left(L_{0}L_{2}\cdot4^{\frac{\dimension\alpha_{0}}{p_{0}}}\right)^{\tau}\cdot\left[\frac{w_{j}^{\left(\vartheta\dimension\alpha_{0}-s^{\ast}\right)}}{w_{i}^{\left(\vartheta\dimension\alpha_{0}-s^{\ast}\right)}}\right]^{\tau}\cdot\left(1+\left\Vert T_{j}^{-1}T_{i}\right\Vert \right)^{\sigma}\cdot\left(\left|\det T_{i}\right|^{-1}\cdot\int_{Q_{i}}\left(1+\left|S_{j}^{-1}\xi\right|\right)^{-M_{00}}\d\xi\right)^{\tau}
\]
for all $i,j\in\Z^{\dimension}\setminus\left\{ 0\right\} $. For brevity,
set $L_{6}:=\left(L_{0}L_{2}\cdot4^{\frac{\dimension\alpha_{0}}{p_{0}}}\right)^{\tau}$.
Note that with this notation, the preceding estimate shows $N_{i,j}\leq N_{6}\cdot M_{j,i}^{\left(0\right)}$,
where $M_{j,i}^{\left(0\right)}$ is defined as in equation (\ref{eq:AlphaModulationStandardEstimate}),
but with $s^{\ast}$ replaced by $s^{\natural}:=\vartheta\dimension\alpha_{0}-s^{\ast}$
and with $N_{0}$ replaced by $M_{00}+1$.

Now, we want to apply Lemma \ref{lem:AlphaModulationMainLemma} to
estimate $M_{j,i}^{\left(0\right)}$. To this end, we have to verify
\begin{align*}
M_{00}+1 & \overset{!}{\geq}\dimension+2+\frac{\dimension+1}{\tau}+\max\left\{ \left|s^{\natural}+\dimension\alpha_{0}\right|,\left|s^{\natural}+\left(\dimension-\frac{\sigma}{\tau}\right)\alpha_{0}\right|\right\} \\
 & =\dimension+2+\frac{\dimension+1}{\tau}+\max\left\{ \left|s^{\ast}-\dimension\alpha_{0}\left(1+\vartheta\right)\right|,\left|s^{\ast}+\left(\frac{\sigma}{\tau}-\dimension\left(1+\vartheta\right)\right)\alpha_{0}\right|\right\} \\
 & =\dimension+2+\frac{\dimension+1}{\tau}+\max\left\{ \left|s^{\ast}-\frac{\dimension\alpha_{0}}{\min\left\{ 1,p\right\} }\right|,\left|s^{\ast}+\alpha_{0}\left(\frac{\sigma}{\tau}-\frac{\dimension}{\min\left\{ 1,p\right\} }\right)\right|\right\} ,
\end{align*}
where the last line used that
\[
1+\vartheta=\begin{cases}
1=\frac{1}{\min\left\{ 1,p\right\} }, & \text{if }p\in\left[1,\infty\right],\\
1+\frac{1}{p}-1=\frac{1}{p}=\frac{1}{\min\left\{ 1,p\right\} }, & \text{if }p\in\left(0,1\right).
\end{cases}
\]
But we have $\tau=\min\left\{ 1,p,q\right\} \geq\min\left\{ 1,p_{0},q_{0}\right\} =\min\left\{ p_{0},q_{0}\right\} =:\tau_{0}$
and furthermore
\begin{align*}
\frac{\sigma}{\tau} & =\begin{cases}
\left\lceil \left|\mu\right|+\dimension+\varepsilon\right\rceil , & \text{if }p\in\left[1,\infty\right],\\
\frac{\dimension}{p}+\left|\mu\right|+\left\lceil \left|\mu\right|+\frac{\dimension+\varepsilon}{p}\right\rceil , & \text{if }p\in\left(0,1\right)
\end{cases}\\
 & \leq\begin{cases}
\left\lceil \mu_{0}+\dimension+\varepsilon\right\rceil , & \text{if }p_{0}=1,\\
\frac{\dimension}{p_{0}}+\mu_{0}+\left\lceil \mu_{0}+\frac{\dimension+\varepsilon}{p_{0}}\right\rceil , & \text{if }p_{0}\in\left(0,1\right)
\end{cases}\\
 & =:\frac{\sigma_{0}}{\tau_{0}}.
\end{align*}
Hence, in case of $p\in\left[1,\infty\right]$, we have
\[
\left|s^{\ast}-\dimension\alpha_{0}\left(1+\vartheta\right)\right|=\left|s^{\ast}-\frac{\dimension\alpha_{0}}{\min\left\{ 1,p\right\} }\right|\leq\frac{1}{1-\alpha}\left(s_{0}+\alpha\dimension\right)\leq\Lambda,
\]
as well as
\begin{align*}
\left|s^{\ast}+\left(\frac{\sigma}{\tau}-\dimension\left(1+\vartheta\right)\right)\alpha_{0}\right|=\left|s^{\ast}+\alpha_{0}\left(\frac{\sigma}{\tau}-\frac{\dimension}{\min\left\{ 1,p\right\} }\right)\right| & \leq\frac{1}{1-\alpha}\left(s_{0}+\alpha\cdot\left|\left\lceil \left|\mu\right|+\dimension+\varepsilon\right\rceil -\dimension\right|\right)\\
\left({\scriptstyle \text{since }\left\lceil \left|\mu\right|+\dimension+\varepsilon\right\rceil \geq\dimension}\right) & =\frac{1}{1-\alpha}\left[s_{0}+\alpha\cdot\left(\left\lceil \left|\mu\right|+\dimension+\varepsilon\right\rceil -\dimension\right)\right]\leq\Lambda,
\end{align*}
as can be easily verified by distinguishing the cases $p_{0}=1$ and
$p_{0}\in\left(0,1\right)$.

Similarly, in case of $p\in\left(0,1\right)$, we necessarily have
$p_{0}\in\left(0,1\right)$ and thus
\[
\left|\frac{\sigma}{\tau}-\frac{\dimension}{\min\left\{ 1,p\right\} }\right|=\left|\frac{\dimension}{p}+\left|\mu\right|+\left\lceil \left|\mu\right|+\frac{\dimension+\varepsilon}{p}\right\rceil -\frac{\dimension}{p}\right|\leq\mu_{0}+\left\lceil \mu_{0}+\frac{\dimension+\varepsilon}{p_{0}}\right\rceil ,
\]
which easily implies
\[
\left|s^{\ast}-\dimension\alpha_{0}\left(1+\vartheta\right)\right|=\left|s^{\ast}-\frac{\dimension\alpha_{0}}{\min\left\{ 1,p\right\} }\right|\leq\frac{1}{1-\alpha}\left(s_{0}+\frac{\alpha\dimension}{p_{0}}\right)\leq\Lambda,
\]
as well as
\begin{align*}
\left|s^{\ast}+\left(\frac{\sigma}{\tau}-\dimension\left(1+\vartheta\right)\right)\alpha_{0}\right| & =\left|s^{\ast}+\alpha_{0}\left(\frac{\sigma}{\tau}-\frac{\dimension}{\min\left\{ 1,p\right\} }\right)\right|\\
 & \leq\frac{1}{1-\alpha}\left(s_{0}+\alpha\left|\frac{\sigma}{\tau}-\frac{\dimension}{\min\left\{ 1,p\right\} }\right|\right)\\
 & \leq\frac{1}{1-\alpha}\left[s_{0}+\alpha\left(\mu_{0}+\left\lceil \mu_{0}+\frac{\dimension+\varepsilon}{p_{0}}\right\rceil \right)\right]\leq\Lambda.
\end{align*}

All in all, our assumptions on $M_{0}$ thus yield
\begin{align*}
M_{00}+1 & =\dimension+2+\frac{\dimension+1}{\min\left\{ p_{0},q_{0}\right\} }+\Lambda\\
 & \geq\dimension+2+\frac{\dimension+1}{\tau}+\max\left\{ \left|s^{\ast}-\dimension\alpha_{0}\left(1+\vartheta\right)\right|,\left|s^{\ast}+\left(\frac{\sigma}{\tau}-\dimension\left(1+\vartheta\right)\right)\alpha_{0}\right|\right\} \\
 & =\dimension+2+\frac{\dimension+1}{\tau}+\max\left\{ \left|s^{\natural}+\dimension\alpha_{0}\right|,\left|s^{\natural}+\left(\dimension-\frac{\sigma}{\tau}\right)\alpha_{0}\right|\right\} .
\end{align*}
Hence, Lemma \ref{lem:AlphaModulationMainLemma} is applicable. For
brevity, set $L_{7}:=2^{\dimension}\cdot\left(1+\left(2+4r\right)^{\alpha_{0}}\cdot r\right)^{M_{00}+1+\dimension}$.
Using Lemma \ref{lem:AlphaModulationMainLemma} and the estimate $\left|s^{\natural}\right|\leq\left|s^{\ast}\right|+\dimension\alpha_{0}\left|\vartheta\right|\leq\frac{1}{1-\alpha}\left(s_{0}+\alpha\frac{\dimension}{p_{0}}\right)$
and setting $L_{8}:=\frac{\sigma_{0}}{\tau_{0}}+\dimension\left(1+\frac{1}{p_{0}}\right)$,
we now get
\begin{align*}
\quad\qquad K_{1}^{\frac{1}{\tau}} & \leq L_{0}L_{2}\cdot4^{\frac{\dimension\alpha_{0}}{p_{0}}}\cdot\left(\sup_{i\in\Z^{\dimension}\setminus\left\{ 0\right\} }\sum_{j\in\Z^{\dimension}\setminus\left\{ 0\right\} }M_{j,i}^{\left(0\right)}\right)^{1/\tau}\\
 & \leq L_{0}L_{2}\cdot4^{\frac{\dimension\alpha_{0}}{p_{0}}}\cdot2^{\frac{1}{\tau}+\frac{\sigma}{\tau}+\left|s^{\natural}\right|}6^{\dimension/\tau}\cdot\max\left\{ 4^{\alpha_{0}\left(\frac{\sigma}{\tau}+\dimension\right)+\left|s^{\natural}\right|}\cdot12^{M_{00}+1}s_{\dimension},\:\left(2+4r\right)^{\left|s^{\natural}\right|+\alpha_{0}\left[\frac{\sigma}{\tau}+\dimension+1+M_{00}\right]}L_{7}\right\} \\
 & \leq L_{0}L_{2}\cdot6^{\frac{\dimension}{\tau_{0}}}2^{\frac{1}{\tau_{0}}+\frac{\sigma_{0}}{\tau_{0}}+\frac{1}{1-\alpha}\left(s_{0}+3\alpha\frac{\dimension}{p_{0}}\right)}\!\!\cdot\max\left\{ 4^{\alpha_{0}L_{8}+\frac{s_{0}}{1-\alpha}}\cdot12^{M_{00}+1}s_{\dimension},\:\left(2\!+\!4r\right)^{\frac{s_{0}}{1-\alpha}+\alpha_{0}\left[L_{8}+M_{00}+1\right]}L_{7}\right\} \!.
\end{align*}
Hence, $K_{1}^{1/\tau}\leq L_{9}$, for a constant $L_{9}$ which
is independent of $p,q,s,\mu$. Completely similar, we also get $K_{2}^{1/\tau}\leq L_{9}$,
with the same constant. In particular, $K_{1},K_{2}<\infty$.

\end{enumerate}
Having verified all prerequisites of Corollary \ref{cor:AtomicDecompositionSimplifiedCriteria},
we get $\vertiii{\smash{\overrightarrow{C}}}^{\max\left\{ 1,\frac{1}{p}\right\} }\leq\Omega\cdot\left(K_{1}^{1/\tau}+K_{2}^{1/\tau}\right)\leq2\Omega\cdot L_{9}$,
where $\overrightarrow{C}:\ell_{\left[w^{\left(s^{\ast}\right)}\right]^{\min\left\{ 1,p\right\} }}^{r}\left(\Z^{\dimension}\setminus\left\{ 0\right\} \right)\to\ell_{\left[w^{\left(s^{\ast}\right)}\right]^{\min\left\{ 1,p\right\} }}^{r}\left(\Z^{\dimension}\setminus\left\{ 0\right\} \right)$
is defined as in Assumption \ref{assu:AtomicDecompositionAssumption}
(with $r:=\max\left\{ q,\frac{q}{p}\right\} $) and where
\begin{align*}
\Omega & =\Omega_{0}^{K}\Omega_{1}\cdot\left(4\cdot\dimension\right)^{1+2\left\lceil K+\frac{\dimension+\varepsilon}{\min\left\{ 1,p\right\} }\right\rceil }\cdot\left(\frac{s_{\dimension}}{\varepsilon}\right)^{1/\min\left\{ 1,p\right\} }\cdot\max_{\left|\beta\right|\leq\left\lceil K+\frac{\dimension+\varepsilon}{\min\left\{ 1,p\right\} }\right\rceil }C^{\left(\beta\right)}\\
 & \leq2^{\mu_{0}}\cdot\left(4\cdot\dimension\right)^{1+2\left\lceil \mu_{0}+\frac{\dimension+\varepsilon}{p_{0}}\right\rceil }\cdot\left(1+\frac{s_{\dimension}}{\varepsilon}\right)^{1/p_{0}}\cdot\max_{\left|\beta\right|\leq\left\lceil \mu_{0}+\frac{\dimension+\varepsilon}{p_{0}}\right\rceil }C^{\left(\beta\right)}=:L_{10},
\end{align*}
where the constants $C^{\left(\beta\right)}=C^{\left(\beta\right)}\left(\Phi\right)$
are defined as in Assumption \ref{assu:RegularPartitionOfUnity}.
Note again that $L_{10}$ is independent of $p,q,s,\mu$. Finally,
Corollary \ref{cor:AtomicDecompositionSimplifiedCriteria} shows that
the families $\left(\gamma_{i}\right)_{i\in\Z^{\dimension}\setminus\left\{ 0\right\} },\left(\gamma_{i,1}\right)_{i\in\Z^{\dimension}\setminus\left\{ 0\right\} }$
and $\left(\gamma_{i,2}\right)_{i\in\Z^{\dimension}\setminus\left\{ 0\right\} }$
fulfill Assumption \ref{assu:AtomicDecompositionAssumption} and that
$\left(\gamma_{i}\right)_{i\in\Z^{\dimension}\setminus\left\{ 0\right\} }$
satisfies Assumption \ref{assu:GammaCoversOrbit}, so that Theorem
\ref{thm:AtomicDecomposition} is applicable.

This theorem shows that $\Gamma^{\left(\delta\right)}$ indeed forms
an atomic decomposition for $M_{\left(s,\mu\right),\alpha}^{p,q}\left(\R^{\dimension}\right)=\DecompSp{\CalQ_{r}^{\left(\alpha\right)}}p{\ell_{w^{\left(s^{\ast}\right)}}^{q}}{v^{\left(\mu\right)}}$,
as soon as $0<\delta\leq\min\left\{ 1,\delta_{00}\right\} $, with
\[
\delta_{00}^{-1}\!:=\!\begin{cases}
\!\frac{2s_{\dimension}}{\sqrt{\dimension}}\!\cdot\!\left(2^{17}\!\cdot\!\dimension^{2}\!\cdot\!\left(K\!+\!2\!+\!\dimension\right)\right)^{\!K+\dimension+3}\cdot\left(1\!+\!R_{\CalQ_{r}^{\left(\alpha\right)}}\right)^{\dimension+1}\!\cdot\!\Omega_{0}^{4K}\Omega_{1}^{4}\Omega_{2}^{\left(p,K\right)}\Omega_{4}^{\left(p,K\right)}\cdot\vertiii{\smash{\overrightarrow{C}}}\,, & \text{if }p\geq1,\\
\!\frac{\left(2^{14}/\dimension^{\frac{3}{2}}\right)^{\!\frac{\dimension}{p}}}{2^{45}\cdot\dimension^{17}}\!\cdot\!\left(\frac{s_{\dimension}}{p}\right)^{\!\frac{1}{p}}\left(2^{68}\!\cdot\!\dimension^{14}\!\cdot\!\left(K\!+\!1\!+\!\frac{\dimension+1}{p}\right)^{\!\!3}\,\right)^{\!K+2+\frac{\dimension+1}{p}}\!\!\!\cdot\!\left(1\!+\!R_{\CalQ_{r}^{\left(\alpha\right)}}\right)^{1+\frac{3\dimension}{p}}\!\cdot\!\Omega_{0}^{16K}\Omega_{1}^{16}\Omega_{2}^{\left(p,K\right)}\Omega_{4}^{\left(p,K\right)}\!\cdot\!\vertiii{\smash{\overrightarrow{C}}}^{\frac{1}{p}}, & \text{if }p<1.
\end{cases}
\]

Thus, to establish the claim of the current theorem, it suffices to
show that $\delta_{00}^{-1}$ can be bounded independently of $p,q,s,\mu$.
Strictly speaking, we then still have to verify that the coefficient
space used in Theorem \ref{thm:AtomicDecomposition} is just $\mathscr{C}_{p,q,s,\mu}^{\left(\alpha\right)}$,
but this can be done precisely as in the proof of Theorem \ref{thm:AlphaModulationBanachFrame}.

Now, for $p\geq1$, we have because of $K=\left|\mu\right|\leq\mu_{0}$
and $\Omega_{0}=1$, as well as $\Omega_{1}=2^{\mu_{0}}$ and $\Omega_{4}^{\left(p,K\right)}\leq L_{5}$,
as well as $\Omega_{2}^{\left(p,K\right)}\leq L_{4}$ that
\[
\delta_{00}^{-1}\leq\frac{2s_{\dimension}}{\sqrt{\dimension}}\cdot\left(2^{17}\!\cdot\!\dimension^{2}\!\cdot\!\left(\mu_{0}\!+\!2\!+\!\dimension\right)\right)^{\!\mu_{0}+\dimension+3}\!\!\!\cdot\left(1\!+\!R_{\CalQ_{r}^{\left(\alpha\right)}}\right)^{\dimension+1}\cdot L_{4}L_{5}\cdot16^{\mu_{0}}\cdot2L_{9}L_{10},
\]
which is independent of $p,q,s,\mu$.

Finally, in case of $p\in\left(0,1\right)$, we get similarly that
\[
\delta_{00}^{-1}\leq\frac{\left(2^{14}\right)^{\!\frac{\dimension}{p_{0}}}}{2^{45}\cdot\dimension^{17}}\!\cdot\!\left(1+\frac{s_{\dimension}}{p_{0}}\right)^{\!\frac{1}{p_{0}}}\left(2^{68}\!\cdot\!\dimension^{14}\!\cdot\!\left(\mu_{0}+1+\frac{\dimension+1}{p_{0}}\right)^{3}\right)^{\!\mu_{0}+2+\frac{\dimension+1}{p_{0}}}\!\!\!\cdot\!\left(1\!+\!R_{\CalQ_{r}^{\left(\alpha\right)}}\right)^{1+\frac{3\dimension}{p_{0}}}\!\cdot\!2^{16\mu_{0}}\!\cdot\!L_{4}L_{5}\!\cdot\!2L_{9}L_{10},
\]
which is independent of $p,q,s,\mu$, as desired.
\end{proof}
We close this section with an overview over the history and applications
of $\alpha$-modulation spaces and with a comparison of our results
to the established literature.
\begin{rem}
\label{rem:AlphamModulationLiteratureDiscussion}$\alpha$-modulation
spaces were originally introduced in Gröbner's PhD.\@ thesis \cite{GroebnerAlphaModulationSpaces},
see also \cite{DecompositionSpaces1}. The definition of these spaces
was motivated by the realization that both modulation spaces and (inhomogeneous)
Besov spaces fit into the common framework of decomposition spaces,
which Feichtinger and Gröbner developed at the time\cite{DecompositionSpaces1,DecompositionSpaces2}.
The underlying frequency coverings are the uniform covering for the
modulation spaces and the dyadic covering for the Besov spaces. Given
these two ``end-point'' types of spaces, the $\alpha$-modulation
spaces provide a continuous family of smoothness spaces which ``lie
between'' modulation and Besov spaces in the sense that the associated
$\alpha$-modulation coverings $\CalQ_{r}^{\left(\alpha\right)}$
form a continuously indexed family of coverings which yields the uniform
covering for $\alpha=0$ and the (inhomogeneous) dyadic covering in
the limit $\alpha\uparrow1$. Note though that the $\beta$-modulation
space $M_{s,\beta}^{p,q}\left(\R^{\dimension}\right)$ \emph{can not}
be obtained by complex interpolation between two $\alpha$-modulation
spaces $M_{s_{i},\alpha_{i}}^{p_{i},q_{i}}\left(\R^{\dimension}\right)$
with $\alpha_{1}\neq\alpha_{2}$, except in a few trivial special
cases\cite{AlphaModulationNotInterpolation}. We also remark that
the frequency covering associated to the $\alpha$-modulation spaces
was independently introduced by Päivärinta and Somersalo\cite[Lemma 2.1]{SomersaloAlphaModulation},
in order to generalize the Calderón-Vaillancourt boundedness result
for pseudodifferential operators to the local Hardy spaces $h_{p}$.

\medskip{}

As for the classical modulation spaces, one application of $\alpha$-modulation
spaces is that they are suitable domains for the study of pseudodifferential
operators: For the one-dimensional case, it was shown in \cite{BorupPsiDOsOnAlphaModulation}
that if $T\in{\rm OPS}_{\alpha,\delta}^{s_{0}}$, then $T:M_{s,\alpha}^{p,q}\left(\R\right)\to M_{s+\alpha-s_{0}-1,\alpha}^{p,q}\left(\R\right)$
is continuous, if $p,q\in\left(1,\infty\right)$, $0<\alpha\leq1$,
$0\leq\delta\leq\alpha$ and $\delta<1$. The multivariate case was
considered in \cite{BorupNielsenPsiDOsOnMultivariateAlphaModulationSpaces},
where it was shown for symbols $\sigma\in S_{\varrho,0}^{m}$ and
for $0\leq\alpha\leq\varrho\leq1$ that $\sigma\left(x,D\right):M_{s,\alpha}^{p,q}\left(\R^{\dimension}\right)\to M_{s-m,\alpha}^{p,q}\left(\R^{\dimension}\right)$
is continuous for $p,q\in\left(1,\infty\right)$. For related results,
see also \cite{NazaretBrunoPsiDOsOnAlphaModulation}.

\medskip{}

The embeddings between $\alpha$-modulation spaces and other function
spaces have been considered by a number of authors: Already in Gröbner's
PhD.\@ thesis \cite{GroebnerAlphaModulationSpaces}, embeddings between
$M_{s_{1},\alpha_{1}}^{p,q}\left(\R^{\dimension}\right)$ and $M_{s_{2},\alpha_{2}}^{p,q}\left(\R^{\dimension}\right)$
for $\alpha_{1},\alpha_{2}\in\left[0,1\right]$ are considered, but
the resulting criteria are not sharp. These non-sharp conditions were
improved by Toft and Wahlberg\cite{ToftWahlbergAlphaModulationEmbeddings},
shortly before the question was completely solved by Han and Wang\cite{HanWangAlphaModulationEmbeddings}.
Note, however, that the preceding results only considered the embedding
$M_{s_{1},\alpha_{1}}^{p,q}\left(\R^{\dimension}\right)\hookrightarrow M_{s_{2},\alpha_{2}}^{p,q}\left(\R^{\dimension}\right)$,
where the exponents $p,q$ are \emph{the same} on both sides. The
existence of the completely general embedding $M_{s_{1},\alpha_{1}}^{p_{1},q_{1}}\left(\R^{\dimension}\right)\hookrightarrow M_{s_{2},\alpha_{2}}^{p_{2},q_{2}}\left(\R^{\dimension}\right)$
was characterized completely in my PhD.\@ thesis \cite[Theorems 6.1.7 and 6.2.8]{VoigtlaenderPhDThesis}.
The same results also appear in my recent paper \cite[Theorems 9.7, 9.13, 9.14 and Corollary 9.16]{DecompositionEmbedding}.
We finally remark that the complete characterization of the embedding
$M_{s_{1},\alpha_{1}}^{p_{1},q_{1}}\left(\R^{\dimension}\right)\hookrightarrow M_{s_{2},\alpha_{2}}^{p_{2},q_{2}}\left(\R^{\dimension}\right)$
was also obtained independently in \cite{AlphaModulationEmbeddingFullCharacterization}.

The existence of embeddings $M_{s_{1},\alpha}^{p,q}\left(\R^{\dimension}\right)\hookrightarrow L_{s_{2}}^{p}\left(\R^{\dimension}\right)$
and $L_{s_{1}}^{p}\left(\R^{\dimension}\right)\hookrightarrow M_{s_{2},\alpha}^{p,q}\left(\R^{\dimension}\right)$
between $\alpha$-modulation spaces and Sobolev spaces (or Bessel
potential spaces) has been fully characterized in \cite{EmbeddingsOfAlphaModulationIntoSobolev};
in the same paper, the author also considers these embeddings when
the Sobolev spaces $L_{s}^{p}\left(\R^{\dimension}\right)$ are replaced
by the local hardy spaces $h_{p}\left(\R^{\dimension}\right)$, for
$p\in\left(0,1\right)$. Embeddings between $\alpha$-modulation spaces
$M_{s,\alpha}^{p_{1},q}\left(\R^{\dimension}\right)$ and the classical
Sobolev spaces $W^{k,p_{2}}\left(\R^{\dimension}\right)$ are also
considered in \cite[Example 7.3]{DecompositionIntoSobolev}, as an
application of a more general theory: For $p_{2}\in\left[1,2\right]\cup\left\{ \infty\right\} $,
a complete characterization of the existence of this embedding is
obtained, but for $p_{2}\in\left(2,\infty\right)$, the given sufficient
criteria are strictly stronger than the given necessary criteria.

\medskip{}

Finally, we discuss the existing results concerning Banach frames
and atomic decompositions for $\alpha$-modulation spaces. A large
number of results in this direction were obtained by Borup and Nielsen:
In \cite{AlphaModulationNonlinearApproximation}, they showed that
certain \textbf{brushlet orthonormal bases} of $L^{2}\left(\R\right)$
form unconditional bases for the $\alpha$-modulation spaces $M_{s,\alpha}^{p,q}\left(\R\right)$,
for $p,q\in\left(1,\infty\right)$. Furthermore, a characterization
of the $\alpha$-modulation (quasi)-norm in terms of the brushlet
coefficients is obtained for arbitrary $p,q\in\left(0,\infty\right]$.
In fact, Borup and Nielsen showed that brushlet bases even yield \textbf{greedy
bases} (i.e., they are unconditional bases which satisfy the so-called
\textbf{democracy condition}), which was then used to characterize
the associated approximation spaces. As a further application of these
brushlet bases for $\alpha$-modulation spaces, Borup and Nielsen
derived boundedness results for certain pseudodifferential operators,
as briefly discussed above. In \cite{NielsenOrthonormalBasesForAlphaModulation},
Nielsen generalized the results concerning brushlet bases from the
one-dimensional case to the case $\dimension=2$. Despite their great
utility, we remark that brushlet bases are \emph{not} generated by
a single prototype function; furthermore, brushlets are bandlimited
and can thus \emph{not} be compactly supported.

In addition to these brushlet bases for $\alpha$-modulation spaces
in dimensions $\dimension=1$ and $\dimension=2$, Borup and Nielsen
also proved existence of Banach frames for $\alpha$-modulation spaces
for the general case $\dimension\in\N$. A first construction of a
\emph{non-tight} frame with explicitly given dual was obtained in
\cite{BorupNielsenAlphaModulationSpaces} and then generalized in
\cite[Section 6.1]{BorupNielsenDecomposition} to obtain \emph{tight}
Banach frames, even for the case of general decomposition spaces,
not just for $\alpha$-modulation spaces. But again, the frames constructed
in these two papers are \emph{not} generated by a single prototype
function and are bandlimited.

These limitations were partly overcome in later papers: In \cite[Theorem 1.1]{CompactlySupportedFramesForDecompositionSpaces},
Nielsen and Rasmussen obtained compactly supported Banach frames for
$\alpha$-modulation spaces. These frames $\left(\psi_{k,n}\right)_{k\in\Z^{\dimension}\setminus\left\{ 0\right\} ,n\in\Z^{\dimension}}$,
however, are \emph{not} of the same structured form as in Theorems
\ref{thm:AlphaModulationBanachFrame} and \ref{thm:AlphaModulationAtomicDecomposition}.
Instead, $\psi_{k,n}\left(x\right)=e^{i\left\langle x,d_{k}\right\rangle }\sum_{\ell=1}^{K}a_{k,\ell}\cdot g\left(c_{k}x+b_{k,n,\ell}\right)$
for suitable (but unknown) $K\in\N$, $a_{k,\ell}\in\Compl$, $d_{k},b_{k,n,\ell}\in\R^{\dimension}$
and $c_{k}\in\R^{\ast}$, which makes it hard to say anything about
e.g.\@ the time-frequency concentration of the family $\left(\psi_{k,n}\right)_{k,n}$.
This limitation was finally removed by Nielsen in \cite{NielsenSinglyGeneratedFrames},
where it was shown—for a special class of coverings $\CalQ$, which
includes the $\alpha$ coverings $\CalQ_{r}^{\left(\alpha\right)}$
for $0\leq\alpha<1$, cf.\@ \cite[Lemma 2.9]{NielsenSinglyGeneratedFrames}—that
one can choose a single, compactly supported generator (or prototype)
$\gamma$ such that the resulting family $\Gamma^{\left(\delta\right)}$
defined as in Theorem \ref{thm:AlphaModulationBanachFrame} yields
a Banach frame for the $\alpha$-modulation space $M_{s,\alpha}^{p,q}\left(\R^{\dimension}\right)$.
The main difference in comparison to Theorems \ref{thm:AlphaModulationBanachFrame}
and \ref{thm:AlphaModulationAtomicDecomposition} is that \cite{NielsenSinglyGeneratedFrames}
merely establishes \emph{existence} of a suitable generator $\gamma$;
it does \emph{not} provide readily verifiable conditions on $\gamma$
which allow to decide if $\gamma$ is suitable. This is due to the
employed proof technique: Nielsen first shows that a certain bandlimited
generator $\gamma_{0}$ generators a Banach frame and then uses a
perturbation argument to show that $\gamma$ generates a Banach frame,
provided that $\gamma$ is close enough to $\gamma_{0}$ in a certain
sense. The most concrete criterion in \cite{NielsenSinglyGeneratedFrames}
concerning $\gamma$ is that under certain \emph{readily verifiable}
conditions on a ``prototypical prototype'' $g$ (cf.\@ \cite[equations (3.13), (3.14)]{NielsenSinglyGeneratedFrames}),
one can always obtain a \emph{suitable} $\gamma$ by taking linear
combinations of translations of $g$. In stark contrast, Theorems
\ref{thm:AlphaModulationBanachFrame} and \ref{thm:AlphaModulationAtomicDecomposition}
(plus the associated remark) show that \emph{any} compactly supported
prototype $\gamma$ generates Banach frames and atomic decompositions
for the $\alpha$-modulation spaces, assuming that $\gamma$ is sufficiently
smooth and has nonvanishing Fourier transform on a certain neighborhood
of the origin.

\medskip{}

In addition to the constructions by Borup, Nielsen and Rasmussen,
Banach frames for $\alpha$-modulation spaces have also been considered
by Fornasier\cite{FornasierFramesForAlphaModulation}, by Dahlke et
al\@.\cite{QuotientCoorbitTheoryAndAlphaModulationSpaces} and finally
by Speckbacher et al.\cite{SpeckbacherAlphaModulation}, all for the
case $\dimension=1$ and $p,q\in\left[1,\infty\right]$. The idea
in \cite{FornasierFramesForAlphaModulation} is to show that a family
$\Gamma_{\alpha}^{\left(\delta\right)}$ similar to $\Gamma^{\left(\delta\right)}$
from Theorem \ref{thm:AlphaModulationBanachFrame} is \textbf{intrinsicly
self-localized}, under suitable readily verifiable assumptions on
the generator $\gamma$, so that, for a sufficiently small sampling
density, the family $\Gamma_{\alpha}^{\left(\delta\right)}$ forms
a Banach frame and an atomic decomposition for the $\alpha$-modulation
space $M_{s,\alpha}^{p,q}\left(\R\right)$. In particular, $\gamma$
can be taken to have compact support, since any Schwartz function
with $\widehat{\gamma}\left(\xi\right)\neq0$ on $\left[-1,1\right]$
is suitable, cf.\@ \cite[Theorem 3.4]{FornasierFramesForAlphaModulation}.
Hence, Fornasier's results are very similar to Theorems \ref{thm:AlphaModulationBanachFrame}
and \ref{thm:AlphaModulationAtomicDecomposition}; the main difference
is that the results in this paper apply for the full range $p,q\in\left(0,\infty\right]$
and also for $\dimension>1$. In this context, we remark that Fornasier
notes that ``We expect that the approach illustrated in this paper
{[}i.e., in \cite{FornasierFramesForAlphaModulation}{]} can be useful
also for a frame characterization of $M_{p,q}^{s,\alpha}\left(\R^{\dimension}\right)$
for $\dimension>1$, \emph{with major technical difficulties}.''
A further difference is that Fornasier only requires decay of $\widehat{\gamma}$,
not of $\partial^{\alpha}\widehat{\gamma}$, cf.\@ \cite[equation (29)]{FornasierFramesForAlphaModulation}.

The approach by Dahlke et al\@.\cite{QuotientCoorbitTheoryAndAlphaModulationSpaces}
is still different: They consider \textbf{coorbit spaces} of certain
quotients of the \textbf{affine Weyl-Heisenberg group}, based on the
general theory of coorbit spaces of homogeneous spaces developed in
\cite{CoorbitOnHomogenousSpaces,CoorbitOnHomogenousSpaces2}. It is
shown in \cite{QuotientCoorbitTheoryAndAlphaModulationSpaces} that
the general theory applies in this setting and (for suitable quotients)
that the associated coorbit spaces coincide with certain $\alpha$-modulation
spaces. Precisely, \cite[Theorem 6.1]{QuotientCoorbitTheoryAndAlphaModulationSpaces}
shows that $\mathcal{H}_{p,v_{s-\alpha\left(1/p-1/2\right),\alpha}}=M_{s,\alpha}^{p,p}\left(\R\right)$,
up to trivial identifications. Note that this only includes $\alpha$-modulation
spaces with $p=q$. The same theorem also shows—as a consequence of
the general discretization theory for coorbit spaces—that one obtains
Banach frames and atomic decompositions for the spaces $\mathcal{H}_{p,v_{s-\alpha\left(1/p-1/2\right),\alpha}}=M_{s,\alpha}^{p,p}\left(\R\right)$
which are of a similar form as the family $\Gamma^{\left(\delta\right)}$
considered in Theorems \ref{thm:AlphaModulationBanachFrame} and \ref{thm:AlphaModulationAtomicDecomposition}.
The main difference of these two theorems to the results in \cite{QuotientCoorbitTheoryAndAlphaModulationSpaces}
is that \cite{QuotientCoorbitTheoryAndAlphaModulationSpaces} is only
applicable to \emph{bandlimited generators}, only for $p=q\in\left[1,\infty\right]$
and only for $\dimension=1$. Furthermore, while in Theorems \ref{thm:AlphaModulationBanachFrame}
and \ref{thm:AlphaModulationAtomicDecomposition} only the sampling
density \emph{in the space domain} has to be sufficiently dense, in
\cite{QuotientCoorbitTheoryAndAlphaModulationSpaces} one needs to
adjust the sampling density in space \emph{and in frequency}. Note
though that the assumption $\widehat{\gamma}\left(\xi\right)\neq0$
on $\overline{B_{r}}\left(0\right)$ is not present in \cite{QuotientCoorbitTheoryAndAlphaModulationSpaces}.

Finally, in \cite{SpeckbacherAlphaModulation}, Speckbacher et al.\@
extended the results of \cite{QuotientCoorbitTheoryAndAlphaModulationSpaces}
by showing that the theory developed in \cite{CoorbitOnHomogenousSpaces,CoorbitOnHomogenousSpaces2,QuotientCoorbitTheoryAndAlphaModulationSpaces}
is also applicable for suitable \emph{compactly supported} generators,
only subject to certain decay and smoothness conditions, cf.\@ \cite[Theorem 5.9]{SpeckbacherAlphaModulation}.
Note though that this does not remove the assumptions $p=q\in\left[1,\infty\right]$
and $\dimension=1$. Precisely, the condition on the generator $\gamma\in L^{2}\left(\R\right)$
imposed in \cite{SpeckbacherAlphaModulation} to generate a Banach
frame and an atomic decomposition for $\mathcal{H}_{p,v_{s-\alpha\left(1/p-1/2\right),\alpha}}=M_{s,\alpha}^{p,p}\left(\R\right)$
is that $\widehat{\gamma}\in C^{3}\left(\R\right)$ with $\left|\partial^{j}\widehat{\gamma}\left(\xi\right)\right|\lesssim\left(1+\left|\xi\right|\right)^{-r}$
for $j\in\left\{ 0,1,2,3\right\} $ and
\[
r>1+\frac{2+2\left[s-\alpha\left(1/p-1/2\right)\right]+7\alpha-4\alpha^{2}}{2\left(1-\alpha\right)^{2}}=1+\frac{1+s+\alpha\left(4-\frac{1}{p}\right)-2\alpha^{2}}{\left(1-\alpha\right)^{2}}\geq\frac{2+s+\alpha\left(1-\alpha\right)}{\left(1-\alpha\right)^{2}}.
\]
In comparison, at least for $\gamma\in C_{c}^{1}\left(\R\right)$,
Theorem \ref{thm:AlphaModulationBanachFrame} (with $\mu_{0}=0$,
$p_{0}=q_{0}=1$, $s_{0}=s_{1}=s$ and $\varepsilon=\frac{1}{2}$)
only requires $\left|\partial^{j}\widehat{\gamma}\left(\xi\right)\right|\lesssim\left(1+\left|\xi\right|\right)^{-N_{0}}$
for $j\in\left\{ 0,1,2\right\} $ and
\[
N_{0}=5+\frac{1}{1-\alpha}\cdot\max\left\{ s+\alpha,\:s+2\alpha\right\} =5+\frac{s+2\alpha}{1-\alpha}=\frac{5+s-3\alpha}{1-\alpha}
\]
and Theorem \ref{thm:AlphaModulationAtomicDecomposition} requires
the same, but with
\[
N_{0}=2\cdot\left(3+\varepsilon\right)+\frac{s+\alpha}{1-\alpha}=6+2\varepsilon+\frac{s+\alpha}{1-\alpha},
\]
where $\varepsilon\in\left(0,1\right)$ can be chosen arbitrarily
small. Hence, in particular for $\alpha\approx1$ or for $\alpha>0$
and large $s$, the prerequisites of Theorems \ref{thm:AlphaModulationBanachFrame}
and \ref{thm:AlphaModulationAtomicDecomposition} are easier to fulfill
than those in \cite{SpeckbacherAlphaModulation}. We finally remark
that Fornasier\cite{FornasierFramesForAlphaModulation} requires $\left|\widehat{\gamma}\left(\xi\right)\right|\lesssim\left(1+\left|\xi\right|\right)^{-M_{0}}$
where $M_{0}$ satisfies $M_{0}\geq5+\frac{1}{2}+\frac{2s+2\alpha}{1-\alpha}$,
which is very close to the conditions in Theorems \ref{thm:AlphaModulationBanachFrame}
and \ref{thm:AlphaModulationAtomicDecomposition}.

Hence, our very general theory yields results comparable to specialized
treatments like \cite{FornasierFramesForAlphaModulation} and (in
many cases) better results than those obtained by using general coorbit
theory. Further, it is applicable for general $\dimension\in\N$ and
also for general $p,q\in\left(0,\infty\right]$ instead of only for
$p=q\in\left[1,\infty\right]$.
\end{rem}

\section{Existence of compactly supported Banach frames and atomic decompositions
for inhomogeneous Besov spaces}

\label{sec:BesovFrames}In this section, we investigate assumptions
on the \textbf{scaling function} $\varphi:\R^{\dimension}\to\Compl$
and the \textbf{mother wavelet} $\psi:\R^{\dimension}\to\Compl$ which
ensure that the associated inhomogeneous \textbf{wavelet system} with
sampling density $c>0$,
\[
W\left(\varphi,\psi;c\right):=\left(\varphi\left(\bullet-c\cdot k\right)\right)_{k\in\Z^{\dimension}}\cup\left(2^{j\frac{\dimension}{2}}\cdot\psi\left(2^{j}\bullet-c\cdot k\right)\right)_{j\in\N,k\in\Z^{\dimension}},
\]
generates Banach frames or atomic decompositions for a subclass of
the class of inhomogeneous \textbf{Besov spaces}.

The inhomogeneous Besov spaces are decomposition spaces which are
defined using a certain dyadic covering of $\R^{\dimension}$, which
we introduce now.
\begin{lem}
\label{lem:InhomogeneousBesovCovering}For $j\in\N$, let $T_{j}:=2^{j}\cdot\identity$
and $b_{j}:=0$, as well as $Q_{j}':=B_{4}\left(0\right)\setminus\overline{B_{1/4}\left(0\right)}$.
Furthermore, set $T_{0}:=\identity$ and $b_{0}:=0$, as well as $Q_{0}':=B_{2}\left(0\right)$.
The \textbf{(inhomogeneous) Besov covering} of $\R^{\dimension}$
is given by
\[
\mathscr{B}:=\left(Q_{j}\right)_{j\in\N_{0}}:=\left(T_{j}Q_{j}'+b_{j}\right)_{j\in\N_{0}}.
\]
This covering is a semi-structured admissible covering of $\R^{\dimension}$.
Furthermore, $\mathscr{B}$ admits a regular partition of unity $\Phi=\left(\varphi_{j}\right)_{j\in\N_{0}}$
(which thus fulfills Assumption \ref{assu:RegularPartitionOfUnity}),
which we fix for the remainder of the section.

Finally, $\mathscr{B}$ fulfills the standing assumptions from Section
\ref{subsec:DecompSpaceDefinitionStandingAssumptions}; in particular,
$\left\Vert T_{j}^{-1}\right\Vert \leq1=:\Omega_{0}$ for all $j\in\N_{0}$.
\end{lem}
\begin{proof}
It was shown in \cite[Example 7.2]{DecompositionIntoSobolev} that
$\mathscr{B}$ is a semi-structured covering of $\R^{\dimension}$.
In the same example, it was also shown that $\mathscr{B}$ is in fact
a regular covering of $\R^{\dimension}$, i.e., $\mathscr{B}$ admits
a regular partition of unity $\Phi$, as claimed. Thanks to Corollary
\ref{cor:RegularBAPUsAreWeightedBAPUs}, $\Phi$ is also a $\mathscr{B}$-$v_{0}$-BAPU
for each weight $v_{0}$ satisfying the assumptions from Section \ref{subsec:DecompSpaceDefinitionStandingAssumptions}.

To verify the standing assumptions from Section \ref{subsec:DecompSpaceDefinitionStandingAssumptions}
pertaining to the covering $\CalQ=\mathscr{B}$, we thus only have
to verify $\left\Vert T_{j}^{-1}\right\Vert \leq1$ for all $j\in\N_{0}$.
But since $T_{j}=2^{j}\cdot\identity$ for all $j\in\N_{0}$, we simply
have $\left\Vert T_{j}^{-1}\right\Vert =2^{-j}\leq1$, as claimed.
\end{proof}
Now, we can define the inhomogeneous Besov spaces:
\begin{defn}
\label{def:InhomogeneousBesovSpaces}For $p,q\in\left(0,\infty\right]$
and $s,\mu\in\R$, we define the associated \textbf{(weighted) inhomogeneous
Besov space} as
\[
\mathcal{B}_{s,\mu}^{p,q}\left(\smash{\R^{\dimension}}\right):=\DecompSp{\mathscr{B}}p{\ell_{\left(2^{js}\right)_{j\in\N_{0}}}^{q}}{v^{\left(\mu\right)}},
\]
with $v^{\left(\mu\right)}$ as in Lemma \ref{lem:AlphaModulationSpaceWeight}.
The \textbf{classical inhomogeneous Besov spaces} are given by $\mathcal{B}_{s}^{p,q}\left(\R^{\dimension}\right):=\mathcal{B}_{s,0}^{p,q}\left(\R^{\dimension}\right)$.
\end{defn}
\begin{rem*}

\begin{itemize}[leftmargin=0.4cm]
\item It is not hard to see that the weight $\left(2^{js}\right)_{j\in\N_{0}}$
is indeed $\mathscr{B}$-moderate (cf.\@ equation (\ref{eq:IntroductionModerateWeightDefinition})).
Furthermore, we saw in Lemma \ref{lem:AlphaModulationSpaceWeight}
that the weight $v^{\left(\mu\right)}$ satisfies all assumptions
from Section \ref{subsec:DecompSpaceDefinitionStandingAssumptions},
with $K:=\left|\mu\right|$ and $\Omega_{1}:=2^{\left|\mu\right|}$,
as well as $v_{0}:\R^{\dimension}\to\left(0,\infty\right),x\mapsto\left[2\cdot\left(1+\left|x\right|\right)\right]^{\left|\mu\right|}$.
All in all, we thus see that all standing assumptions from Section
\ref{subsec:DecompSpaceDefinitionStandingAssumptions} are satisfied.
In particular, Lemma \ref{lem:WeightedDecompositionSpaceComplete}
and Proposition \ref{prop:WeightedDecompositionSpaceWellDefined}
show that the spaces $\mathcal{B}_{s,\mu}^{p,q}\left(\R^{\dimension}\right)$
are well-defined Quasi-Banach spaces.
\item It is not hard to see that the usual inhomogeneous Besov spaces (e.g.\@
as defined in \cite[Definition 6.5.1]{GrafakosModernFourierAnalysis})
coincide with the spaces $\mathcal{B}_{s}^{p,q}\left(\R^{\dimension}\right)$
defined above, up to trivial identifications. The main difference
is that the usual Besov spaces are defined as subspaces of the space
of tempered distributions, while $\mathcal{B}_{s}^{p,q}\left(\R^{\dimension}\right)$
is a subspace of $Z'\left(\R^{\dimension}\right)=\left[\Fourier\left(\TestFunctionSpace{\R^{\dimension}}\right)\right]'$,
cf.\@ Subsection \ref{subsec:DecompSpaceDefinitionStandingAssumptions}.
Claiming that the two spaces coincide amounts to claiming that each
$f\in\mathcal{B}_{s}^{p,q}\left(\R^{\dimension}\right)$ extends to
a (uniquely determined) tempered distribution. As shown in \cite[Lemma 9.15]{DecompositionEmbedding},
this is indeed fulfilled.\qedhere
\end{itemize}
\end{rem*}
Note that if we define the family $\Gamma=\left(\gamma_{j}\right)_{j\in\N_{0}}$
by $\gamma_{0}:=\varphi$ and $\gamma_{j}:=\psi$ for $j\in\N$, where
$\varphi,\psi:\R^{\dimension}\to\Compl$ are given, then the family
$\Gamma^{\left(\delta\right)}=\left(L_{\delta\cdot T_{j}^{-T}k}\,\gamma^{\left[j\right]}\right)_{j\in\N_{0},k\in\Z^{\dimension}}$
considered in Theorem \ref{thm:AtomicDecomposition} (and in a slightly
modified form also in Theorem \ref{thm:DiscreteBanachFrameTheorem})
satisfies
\[
\Gamma^{\left(\delta\right)}=\left(\varphi\left(\bullet-\delta k\right)\right)_{k\in\Z^{\dimension}}\cup\left(2^{j\frac{\dimension}{2}}\cdot\psi\left(2^{j}\bullet-\delta k\right)\right)=W\left(\varphi,\psi;\delta\right),
\]
at least up to an obvious re-indexing. Consequently, we can use Corollaries
\ref{cor:BanachFrameSimplifiedCriteria} and \ref{cor:AtomicDecompositionSimplifiedCriteria}
to derive conditions on $\varphi,\psi$ which ensure that the wavelet
system $W\left(\varphi,\psi;\delta\right)$ yields a Banach frame,
or an atomic decomposition for the (weighted) inhomogeneous Besov
spaces $\mathcal{B}_{s,\mu}^{p,q}\left(\R^{\dimension}\right)$. We
begin with the case of Banach frames.
\begin{prop}
\label{prop:BesovBanachFrames}Let $p_{0},q_{0}\in\left(0,1\right]$,
$\varepsilon>0$, $\mu_{0}\geq0$ and $-\infty<s_{0}\leq s_{1}<\infty$.
Assume that $\varphi,\psi:\R^{\dimension}\to\Compl$ satisfy the following
conditions:

\begin{enumerate}
\item We have $\varphi,\psi\in L_{\left(1+\left|\mybullet\right|\right)^{\mu_{0}}}^{1}\left(\R^{\dimension}\right)$
and $\widehat{\varphi},\widehat{\psi}\in C^{\infty}\left(\R^{\dimension}\right)$,
where all partial derivatives of $\widehat{\varphi}$ and $\widehat{\psi}$
are polynomially bounded.
\item We have $\varphi,\psi\in C^{1}\left(\R^{\dimension}\right)$ and $\nabla\varphi,\nabla\psi\in L_{\left(1+\left|\mybullet\right|\right)^{\mu_{0}}}^{1}\left(\R^{\dimension}\right)\cap L^{\infty}\left(\R^{\dimension}\right)$.
\item We have $\widehat{\varphi}\left(\xi\right)\neq0$ for all $\xi\in\overline{B_{2}}\left(0\right)$
and $\widehat{\psi}\left(\xi\right)\neq0$ for all $\xi\in\overline{B_{4}\left(0\right)}\setminus B_{1/4}\left(0\right)$.
\item We have
\begin{align*}
\left|\partial^{\alpha}\widehat{\varphi}\left(\xi\right)\right| & \leq G_{1}\cdot\left(1+\left|\xi\right|\right)^{-L},\\
\left|\partial^{\alpha}\widehat{\psi}\left(\xi\right)\right| & \leq G_{2}\cdot\left(1+\left|\xi\right|\right)^{-L_{1}}\cdot\min\left\{ 1,\left|\xi\right|^{L_{2}}\right\} 
\end{align*}
for all $\xi\in\R^{\dimension}$ and all $\alpha\in\N_{0}^{\dimension}$
with $\left|\alpha\right|\leq\left\lceil \mu_{0}+\frac{\dimension+\varepsilon}{p_{0}}\right\rceil $
for suitable $G_{1},G_{2}>0$ and certain $L_{2}\geq0$ and $L,L_{1}\geq1$
which satisfy
\[
L>1-s_{0}+\vartheta,\qquad L_{1}>1-s_{0}+\vartheta,\qquad\text{ and }\qquad L_{2}>s_{1},
\]
where $\vartheta:=\frac{\dimension}{p_{0}}+\mu_{0}+\left\lceil \mu_{0}+\frac{\dimension+\varepsilon}{p_{0}}\right\rceil $.
\end{enumerate}
Then there is some $\delta_{0}=\delta_{0}\left(p_{0},q_{0},s_{0},s_{1},\mu_{0},\varepsilon,\dimension,\varphi,\psi\right)>0$
such that for each $0<\delta\leq\delta_{0}$, the family 
\[
\Gamma^{\left(\delta\right)}=\left(2^{j\frac{\dimension}{2}}\cdot\widetilde{\gamma_{j}}\left(2^{j}\bullet-\delta k\right)\right)_{j\in\N_{0},k\in\Z^{\dimension}},\qquad\text{ with }\qquad\widetilde{\gamma_{0}}:=\varphi\left(-\bullet\right)\text{ and }\widetilde{\gamma_{j}}:=\psi\left(-\bullet\right)\text{ for }j\in\N,
\]
forms a Banach frame for $\mathcal{B}_{s,\mu}^{p,q}\left(\R^{\dimension}\right)$,
for arbitrary $p,q\in\left(0,\infty\right]$ and $s,\mu\in\R$ satisfying
$p\geq p_{0}$, $q\geq q_{0}$, $s_{0}\leq s\leq s_{1}$ and $\left|\mu\right|\leq\mu_{0}$.

Precisely, this means the following: Define the coefficient space
\[
\mathscr{C}_{p,q,s,\mu}:=\ell_{\left(2^{j\left(s+\dimension\left(\frac{1}{2}-\frac{1}{p}\right)\right)}\right)_{j\in\N_{0}}}^{q}\!\!\!\!\!\!\!\!\!\!\!\!\left(\left[\ell_{\left[\left(1+\left|k\right|/2^{j}\right)^{\mu}\right]_{k\in\Z^{\dimension}}}^{p}\left(\Z^{\dimension}\right)\right]_{j\in\N_{0}}\right)\leq\Compl^{\N_{0}\times\Z^{\dimension}}.
\]
Then the following hold:

\begin{itemize}
\item The \textbf{analysis operator}
\[
A^{\left(\delta\right)}:\mathcal{B}_{s,\mu}^{p,q}\left(\smash{\R^{\dimension}}\right)\to\mathscr{C}_{p,q,s,\mu},f\mapsto\left[\left(\left[2^{j\frac{\dimension}{2}}\cdot\gamma_{j}\left(2^{j}\bullet\right)\right]\ast f\right)\left(\delta\cdot\frac{k}{2^{j}}\right)\right]_{j\in\N_{0},k\in\Z^{\dimension}}
\]
is well-defined and bounded for all $0<\delta\leq1$. Here, $\gamma_{0}:=\varphi$
and $\gamma_{j}:=\psi$ for $j\in\N$. The convolution considered
here is defined as in equation (\ref{eq:SpecialConvolutionPointwiseDefinition}).
\item For $0<\delta\leq1$, there is a bounded linear \textbf{reconstruction
operator} $R^{\left(\delta\right)}:\mathscr{C}_{p,q,s,\mu}\to\mathcal{B}_{s,\mu}^{p,q}\left(\R^{\dimension}\right)$
satisfying $R^{\left(\delta\right)}\circ A^{\left(\delta\right)}=\identity_{\mathcal{B}_{s,\mu}^{p,q}}$.
Furthermore, the action of $R^{\left(\delta\right)}$ on a given sequence
is independent of the precise choice of $p,q,s,\mu$.
\item We have the following \textbf{consistency statement}: If $f\in\mathcal{B}_{s,\mu}^{p,q}\left(\R^{\dimension}\right)$
and if $p_{0}\leq\tilde{p}\leq\infty$ and $q_{0}\leq\tilde{q}\leq\infty$
and if furthermore $s_{0}\leq\tilde{s}\leq s_{1}$ and $\left|\tilde{\mu}\right|\leq\mu_{0}$,
then the following equivalence holds:
\[
f\in\mathcal{B}_{\tilde{s},\tilde{\mu}}^{\tilde{p},\tilde{q}}\left(\smash{\R^{\dimension}}\right)\qquad\Longleftrightarrow\qquad A^{\left(\delta\right)}f\in\mathscr{C}_{\tilde{p},\tilde{q},\tilde{s},\tilde{\mu}}.\qedhere
\]
\end{itemize}
\end{prop}
\begin{proof}
Let $p,q,s,\mu$ as in the statement of the proposition. Our first
goal is to provide suitable estimates for the quantity $M_{j,i}$
which appears in Corollary \ref{cor:BanachFrameSimplifiedCriteria},
i.e.,
\[
M_{j,i}=\left(\frac{w_{j}}{w_{i}}\right)^{\tau}\cdot\left(1+\left\Vert T_{j}^{-1}T_{i}\right\Vert \right)^{\sigma}\cdot\max_{\left|\beta\right|\leq1}\left(\left|\det T_{i}\right|^{-1}\cdot\int_{Q_{i}}\max_{\left|\alpha\right|\leq N}\left|\left(\partial^{\alpha}\widehat{\partial^{\beta}\gamma_{j}}\right)\!\!\left(S_{j}^{-1}\xi\right)\right|\d\xi\right)^{\tau},
\]
where $K=\left|\mu\right|$ (cf.\@ the remark after Definition \ref{def:InhomogeneousBesovSpaces})
and
\begin{align*}
N & =\left\lceil K+\frac{\dimension+\varepsilon}{\min\left\{ 1,p\right\} }\right\rceil \,,\\
\tau & =\min\left\{ 1,p,q\right\} ,\\
\sigma & =\tau\cdot\left(\frac{\dimension}{\min\left\{ 1,p\right\} }+K+\left\lceil K+\frac{\dimension+\varepsilon}{\min\left\{ 1,p\right\} }\right\rceil \right).
\end{align*}
We immediately observe $N\leq N_{0}:=\left\lceil \mu_{0}+\frac{\dimension+\varepsilon}{p_{0}}\right\rceil $,
so that our estimates regarding $\left|\partial^{\alpha}\widehat{\varphi}\left(\xi\right)\right|$
and $\left|\partial^{\alpha}\widehat{\psi}\left(\xi\right)\right|$
can be applied.

Indeed, since $\varphi,\psi\in C^{1}\left(\R^{\dimension}\right)$
with $\nabla\varphi,\nabla\psi\in L_{\left(1+\left|\mybullet\right|\right)^{\mu_{0}}}^{1}\left(\R^{\dimension}\right)\hookrightarrow L^{1}\left(\R^{\dimension}\right)$
and since $\gamma_{j}=\varphi$ for $j=0$ and $\gamma_{j}=\psi$
otherwise, standard properties of the Fourier transform show
\[
\widehat{\partial^{\beta}\gamma_{j}}\left(\xi\right)=\left(2\pi i\xi\right)^{\beta}\cdot\widehat{\gamma_{j}}\left(\xi\right)\qquad\forall\xi\in\R^{\dimension}\qquad\forall\beta\in\N_{0}^{\dimension}\text{ with }\left|\beta\right|\leq1.
\]
Consequently, we get for $i\in\N$ that
\begin{equation}
M_{j,i}\leq2^{\tau s\left(j-i\right)}\cdot\left(1+2^{i-j}\right)^{\sigma}\cdot\left(2\pi\cdot\max_{\left|\beta\right|\leq1}2^{-i\cdot\dimension}\int_{2^{i-2}<\left|\eta\right|<2^{i+2}}\:\max_{\left|\alpha\right|\leq N}\left|\left(\partial^{\alpha}\left[\xi\mapsto\xi^{\beta}\cdot\widehat{\gamma_{j}}\left(\xi\right)\right]\right)\left(\eta/2^{j}\right)\right|\d\eta\right)^{\tau}\!\!.\label{eq:BesovBanachFrameBasicEstimate}
\end{equation}

Next, we observe for $\beta\in\N_{0}^{\dimension}$ with $\left|\beta\right|=1$,
i.e., $\beta=e_{j}$ for some $j\in\underline{\dimension}$, that
\[
\left|\partial^{\nu}\xi^{\beta}\right|=\begin{cases}
\left|\xi^{\beta}\right|\leq\left|\xi\right|\leq1+\left|\xi\right|, & \text{if }\nu=0,\\
1\leq1+\left|\xi\right|, & \text{if }\nu=e_{j},\\
0\leq1+\left|\xi\right| & \text{otherwise}.
\end{cases}
\]
Likewise, in case of $\beta=0$, we have
\[
\left|\partial^{\nu}\xi^{\beta}\right|=\begin{cases}
1\leq1+\left|\xi\right|, & \text{if }\nu=0,\\
0\leq1+\left|\xi\right|, & \text{otherwise}.
\end{cases}
\]
In connection with Leibniz's rule and the $\dimension$-dimensional
binomial theorem (cf.\@ \cite[Section 8.1, Exercise 2.b]{FollandRA}),
this yields for $j=0$ and $\alpha\in\N_{0}^{\dimension}$ with $\left|\alpha\right|\leq N\leq N_{0}$
that
\begin{equation}
\begin{split}\left|\partial^{\alpha}\left[\xi\mapsto\xi^{\beta}\cdot\widehat{\gamma_{0}}\left(\xi\right)\right]\left(\eta\right)\right| & \leq\sum_{\nu\leq\alpha}\binom{\alpha}{\nu}\cdot\left|\partial^{\nu}\eta^{\beta}\right|\cdot\left|\partial^{\alpha-\nu}\widehat{\gamma_{0}}\left(\eta\right)\right|\\
\left({\scriptstyle \text{assumption for }\widehat{\gamma_{0}}=\widehat{\varphi}}\right) & \leq G_{1}\cdot\left(1+\left|\eta\right|\right)^{1-L}\cdot\sum_{\nu\leq\alpha}\binom{\alpha}{\nu}\\
 & \leq2^{N_{0}}G_{1}\cdot\left(1+\left|\eta\right|\right)^{1-L}.
\end{split}
\label{eq:BesovBanachFrameLowPassEstimate}
\end{equation}
Likewise, for $j\in\N$ and $\alpha\in\N_{0}^{\dimension}$ with $\left|\alpha\right|\leq N\leq N_{0}$,
we get
\begin{equation}
\begin{split}\left|\partial^{\alpha}\left[\xi\mapsto\xi^{\beta}\cdot\widehat{\gamma_{j}}\left(\xi\right)\right]\left(\eta\right)\right| & \leq\sum_{\nu\leq\alpha}\binom{\alpha}{\nu}\cdot\left|\partial^{\nu}\eta^{\beta}\right|\cdot\left|\partial^{\alpha-\nu}\widehat{\gamma_{j}}\left(\eta\right)\right|\\
\left(\text{\ensuremath{{\scriptstyle \text{assumption for }\widehat{\gamma_{j}}=\widehat{\psi}}}}\right) & \leq G_{2}\cdot\left(1+\left|\eta\right|\right)^{1-L_{1}}\cdot\min\left\{ 1,\left|\eta\right|^{L_{2}}\right\} \cdot\sum_{\nu\leq\alpha}\binom{\alpha}{\nu}\\
 & \leq2^{N_{0}}G_{2}\cdot\left(1+\left|\eta\right|\right)^{1-L_{1}}\cdot\min\left\{ 1,\left|\eta\right|^{L_{2}}\right\} .
\end{split}
\label{eq:BesovBanachFrameMotherWaveletEstimate}
\end{equation}

Now, we first consider the case $i\in\N$ and note
\begin{equation}
\begin{split}\lambda_{\dimension}\left(\left\{ \eta\in\R^{\dimension}\with2^{i-2}<\left|\eta\right|<2^{i+2}\right\} \right) & =v_{\dimension}\cdot\left(2^{\dimension\left(i+2\right)}-2^{\dimension\left(i-2\right)}\right)\\
 & =2^{i\cdot\dimension}\cdot v_{\dimension}\left(4^{\dimension}-4^{-\dimension}\right)\\
 & \leq2^{i\cdot\dimension}\cdot4^{\dimension}v_{\dimension},
\end{split}
\label{eq:BesovAnnulusMeasureEstimate}
\end{equation}
so that
\begin{equation}
2^{-i\cdot\dimension}\cdot\int_{2^{i-2}<\left|\eta\right|<2^{i+2}}h\left(\eta\right)\d\eta\leq4^{\dimension}v_{\dimension}\cdot\sup_{2^{i-2}<\left|\eta\right|<2^{i+2}}\;h\left(\eta\right)\label{eq:BesovBanachFrameMeanEstimate}
\end{equation}
for each nonnegative (measurable) function $h$. Now, we distinguish
two subcases for estimating $M_{j,i}$.

\textbf{Case 1}: We have $j\in\N$. In this case, we distinguish two
additional subcases:

\begin{enumerate}
\item We have $j\leq i$. For $2^{i-2}<\left|\eta\right|<2^{i+2}$, this
implies $\left|\eta/2^{j}\right|\geq2^{i-2-j}=\frac{2^{i-j}}{4}=\frac{2^{\left|i-j\right|}}{4}$.
Since we have $L_{1}\geq1$, a combination of this estimate with equations
(\ref{eq:BesovBanachFrameBasicEstimate}), (\ref{eq:BesovBanachFrameMotherWaveletEstimate})
and (\ref{eq:BesovBanachFrameMeanEstimate}) yields
\begin{align*}
M_{j,i} & \leq2^{\sigma}\cdot2^{-\tau s\left|j-i\right|}\cdot2^{\sigma\left|i-j\right|}\cdot\left(2^{N_{0}}G_{2}\cdot2\pi\cdot4^{\dimension}v_{\dimension}\cdot4^{L_{1}-1}\cdot2^{\left|i-j\right|\left(1-L_{1}\right)}\right)^{\tau}\\
 & \leq2^{\sigma}\cdot\left(2^{N_{0}}G_{2}\cdot2\pi\cdot4^{\dimension}v_{\dimension}\cdot4^{L_{1}}\right)^{\tau}\cdot2^{\left|i-j\right|\left[\sigma-\tau\left(L_{1}-1+s\right)\right]}\\
 & =:2^{\sigma}\cdot H_{1}^{\tau}\cdot2^{\left|i-j\right|\left[\sigma-\tau\left(L_{1}-1+s\right)\right]}.
\end{align*}
\item We have $i\leq j$. For $2^{i-2}<\left|\eta\right|<2^{i+2}$, this
implies because of $L_{2}\geq0$ that
\[
\min\left\{ 1,\left|\eta/2^{j}\right|^{L_{2}}\right\} \leq\left(2^{i+2-j}\right)^{L_{2}}\leq4^{L_{2}}\cdot2^{-L_{2}\left|i-j\right|}.
\]
Furthermore, $\left(1+\left|\eta/2^{j}\right|\right)^{1-L_{1}}\leq1$,
since $L_{1}\geq1$. Hence, as in the previous case, we can combine
equations (\ref{eq:BesovBanachFrameBasicEstimate}), (\ref{eq:BesovBanachFrameMotherWaveletEstimate})
and (\ref{eq:BesovBanachFrameMeanEstimate}) to derive
\begin{align*}
M_{j,i} & \leq2^{\tau s\left|j-i\right|}\cdot2^{\sigma}\cdot\left(2^{N_{0}}G_{2}\cdot2\pi\cdot4^{\dimension}v_{\dimension}\cdot4^{L_{2}}\cdot2^{-L_{2}\left|i-j\right|}\right)^{\tau}\\
 & =:2^{\sigma}\cdot H_{2}^{\tau}\cdot2^{\tau\left|j-i\right|\left(s-L_{2}\right)}.
\end{align*}
\end{enumerate}
\textbf{Case 2}: We have $j=0$. In this case, we have for $2^{i-2}<\left|\eta\right|<2^{i+2}$
that $\left|\eta/2^{j}\right|=\left|\eta\right|\geq2^{i-2}$ and hence
$\left(1+\left|\eta/2^{j}\right|\right)^{1-L}\leq4^{L-1}2^{-i\left(L-1\right)}\leq4^{L}\cdot2^{-\left(L-1\right)\left|i-j\right|}$.
Here, we used $L\geq1$.

Now, a combination of equations (\ref{eq:BesovBanachFrameBasicEstimate}),
(\ref{eq:BesovBanachFrameLowPassEstimate}) and (\ref{eq:BesovBanachFrameMeanEstimate})
yields
\begin{align*}
M_{0,i} & \leq2^{-\tau s\left|i-j\right|}\cdot2^{\sigma}2^{\sigma\left|i-j\right|}\cdot\left(2^{N_{0}}G_{1}\cdot2\pi\cdot4^{\dimension}v_{\dimension}\cdot4^{L}\cdot2^{-\left(L-1\right)\left|i-j\right|}\right)^{\tau}\\
 & =:2^{\sigma}\cdot H_{3}^{\tau}\cdot2^{\left|i-j\right|\left[\sigma-\tau\left(L-1+s\right)\right]}.
\end{align*}

These are the desired estimates in case of $i\in\N$. It remains to
consider the case $i=0$. Here, equation (\ref{eq:BesovBanachFrameBasicEstimate})
takes on the slightly modified form
\begin{equation}
\begin{split}M_{j,0} & \leq2^{\tau sj}\cdot\left(1+2^{-j}\right)^{\sigma}\cdot\left(2\pi\cdot\max_{\left|\beta\right|\leq1}\int_{B_{2}\left(0\right)}\max_{\left|\alpha\right|\leq N}\left|\left(\partial^{\alpha}\left[\xi\mapsto\xi^{\beta}\cdot\widehat{\gamma_{j}}\left(\xi\right)\right]\right)\left(\eta/2^{j}\right)\right|\d\eta\right)^{\tau}\\
 & \leq2^{\sigma}\cdot2^{\tau sj}\cdot\left(2\pi\cdot\lambda_{\dimension}\left(B_{2}\left(0\right)\right)\cdot\max_{\left|\beta\right|\leq1}\sup_{\left|\eta\right|<2}\left|\left(\partial^{\alpha}\left[\xi\mapsto\xi^{\beta}\cdot\widehat{\gamma_{j}}\left(\xi\right)\right]\right)\left(\eta/2^{j}\right)\right|\right)^{\tau}.
\end{split}
\label{eq:BesovBanachFrameBasicEstimateAtOrigin}
\end{equation}
Now, we again distinguish two cases:

\textbf{Case 1}: We have $j\in\N$. Here, we observe for $\left|\eta\right|<2$
that
\[
\min\left\{ 1,\left|\eta/2^{j}\right|^{L_{2}}\right\} \leq2^{\left(1-j\right)L_{2}}=2^{L_{2}}\cdot2^{-L_{2}\left|i-j\right|}.
\]
Since we also have $L_{1}\geq1$, a combination of equations (\ref{eq:BesovBanachFrameBasicEstimateAtOrigin})
and (\ref{eq:BesovBanachFrameMotherWaveletEstimate}) yields
\begin{align*}
M_{j,0} & \leq2^{\sigma}\cdot2^{\tau s\left|i-j\right|}\cdot\left(2^{N_{0}}G_{2}\cdot2\pi\cdot\lambda_{\dimension}\left(B_{2}\left(0\right)\right)\cdot2^{L_{2}}\cdot2^{-L_{2}\left|i-j\right|}\right)^{\tau}\\
 & =:2^{\sigma}\cdot H_{4}^{\tau}\cdot2^{\tau\left|i-j\right|\left(s-L_{2}\right)}.
\end{align*}

\textbf{Case 2}: We have $j=0$. Because of $L\geq1$, we have $\left(1+\left|\eta/2^{j}\right|\right)^{1-L}\leq1$
for arbitrary $\eta\in\R^{\dimension}$, so that a combination of
equations (\ref{eq:BesovBanachFrameBasicEstimateAtOrigin}) and (\ref{eq:BesovBanachFrameLowPassEstimate})
yields
\begin{align*}
M_{0,0} & \leq2^{\sigma}\cdot\left(2^{N_{0}}G_{1}\cdot2\pi\cdot\lambda_{\dimension}\left(B_{2}\left(0\right)\right)\right)^{\tau}\\
 & =:2^{\sigma}\cdot H_{5}^{\tau}\\
 & =2^{\sigma}\cdot H_{5}^{\tau}\cdot2^{-\left|i-j\right|\zeta},
\end{align*}
where $\zeta\in\R$ can be chosen arbitrarily, since $\left|i-j\right|=0$.

\medskip{}

All in all, if we set $H_{6}:=\max\left\{ H_{1},\dots,H_{5}\right\} $,
a combination of the preceding cases shows
\begin{equation}
M_{j,i}\leq2^{\sigma}\cdot H_{6}^{\tau}\cdot2^{-\tau\left|i-j\right|\min\left\{ 1,L_{2}-s,\,L-1+s-\frac{\sigma}{\tau},\,L_{1}-1+s-\frac{\sigma}{\tau}\right\} }.\label{eq:BesovBanachFrameSummary}
\end{equation}
Note that $H_{1},\dots,H_{5}$, and hence also $H_{6}$, are all independent
of $p,q,\mu,s$.

Furthermore, we have
\[
\frac{\sigma}{\tau}=\frac{\dimension}{\min\left\{ 1,p\right\} }+K+\left\lceil K+\frac{\dimension+\varepsilon}{\min\left\{ 1,p\right\} }\right\rceil \leq\frac{\dimension}{p_{0}}+\mu_{0}+\left\lceil \mu_{0}+\frac{\dimension+\varepsilon}{p_{0}}\right\rceil =\vartheta
\]
and $\sigma=\tau\cdot\frac{\sigma}{\tau}\leq\tau\vartheta$, as well
as $s_{0}\leq s\leq s_{1}$. Hence,
\[
M_{j,i}\leq\left(2^{\vartheta}\cdot H_{6}\right)^{\tau}\cdot2^{-\tau\left|i-j\right|\min\left\{ 1,\,L_{2}-s_{1},\,L-1+s_{0}-\vartheta,\,L_{1}-1+s_{0}-\vartheta\right\} }.
\]
But our assumptions on $L,L_{1},L_{2}$ imply that the exponent $\lambda:=\min\left\{ 1,\,L_{2}-s_{1},\,L-1+s_{0}-\vartheta,\,L_{1}-1+s_{0}-\vartheta\right\} $
is positive. Hence, we get, for the constants $C_{1},C_{2}$ defined
in Corollary \ref{cor:BanachFrameSimplifiedCriteria},
\begin{align*}
C_{1}^{1/\tau}=\left(\sup_{i\in\N_{0}}\sum_{j\in\N_{0}}M_{j,i}\right)^{1/\tau} & \leq2^{\vartheta}\cdot H_{6}\cdot\sup_{i\in\N_{0}}\left(\sum_{j\in\N_{0}}2^{-\tau\lambda\left|i-j\right|}\right)^{1/\tau}\\
\left({\scriptstyle \text{for }\ell=i-j}\right) & \leq2^{\vartheta}\cdot H_{6}\cdot\left(\sum_{\ell\in\Z}2^{-\tau\lambda\left|\ell\right|}\right)^{1/\tau}\\
 & \leq2^{\vartheta}\cdot H_{6}\cdot\left(2\cdot\sum_{\ell=0}^{\infty}2^{-\tau\lambda\ell}\right)^{1/\tau}\\
\left({\scriptstyle \text{since }\ell^{\tau_{0}}\hookrightarrow\ell^{\tau}\text{ is norm-decreasing for }\tau_{0}:=\min\left\{ p_{0},q_{0}\right\} }\right) & \leq2^{\vartheta}\cdot H_{6}\cdot2^{1/\tau}\cdot\left(\sum_{\ell=0}^{\infty}2^{-\tau_{0}\lambda\ell}\right)^{1/\tau_{0}}\\
 & \leq2^{\vartheta+\frac{1}{\tau_{0}}}\cdot H_{6}\cdot\left(\frac{1}{1-2^{-\tau_{0}\lambda}}\right)^{1/\tau_{0}}=:H_{7}<\infty.
\end{align*}
Observe again that $H_{7}$ is independent of $p,q,\mu,s$. Exactly
the same estimate also yields $C_{2}^{1/\tau}\leq H_{7}$.

\medskip{}

Next, we set $\gamma_{1}^{\left(0\right)}:=\varphi$ and $\gamma_{2}^{\left(0\right)}:=\psi$.
Furthermore, we set $n_{j}:=2$ for $j\in\N$ and $n_{0}:=1$, so
that $\gamma_{j}=\gamma_{n_{j}}^{\left(0\right)}$ for all $j\in\N_{0}$.
In the notation of Lemma \ref{lem:GammaCoversOrbitAssumptionSimplified},
these definitions entail
\begin{align*}
Q^{\left(1\right)} & =\bigcup\left\{ Q_{i}'\with i\in\N_{0}\text{ and }n_{i}=1\right\} =Q_{0}'=B_{2}\left(0\right),\\
Q^{\left(2\right)} & =\bigcup\left\{ Q_{i}'\with i\in\N_{0}\text{ and }n_{i}=2\right\} =\bigcup_{j\in\N}Q_{j}'=B_{4}\left(0\right)\setminus\overline{B_{1/4}}\left(0\right).
\end{align*}
But the prerequisites of the current proposition include the assumptions
$\smash{\widehat{\gamma_{1}^{\left(0\right)}}}\left(\xi\right)=\widehat{\varphi}\left(\xi\right)\neq0$
for all $\xi\in\overline{Q^{\left(1\right)}}$ and $\widehat{\gamma_{2}^{\left(0\right)}}\left(\xi\right)=\widehat{\psi}\left(\xi\right)\neq0$
for all $\xi\in\overline{Q^{\left(2\right)}}$. By continuity of $\widehat{\varphi},\widehat{\psi}$
and by compactness of $\overline{Q^{\left(1\right)}},\overline{Q^{\left(2\right)}}$,
we thus see that all assumptions of Lemma \ref{lem:GammaCoversOrbitAssumptionSimplified}
are satisfied. Consequently, the family $\Gamma=\left(\gamma_{i}\right)_{i\in I}$
satisfies Assumption \ref{assu:GammaCoversOrbit} and there is a constant
$\Omega_{3}=\Omega_{3}\left(\mathscr{B},\varphi,\psi,p_{0},\mu_{0},\dimension\right)>0$
satisfying $\Omega_{2}^{\left(p,K\right)}\leq\Omega_{3}$ for all
$K\leq\mu_{0}$ and $p\geq p_{0}$, with $\Omega_{2}^{\left(p,K\right)}$
as in Assumption \ref{assu:GammaCoversOrbit}. Recall that in our
case, we indeed have $K=\left|\mu\right|\leq\mu_{0}$.

In view of the assumptions of the proposition, it is now not hard
to see that all prerequisites for Corollary \ref{cor:BanachFrameSimplifiedCriteria}
are satisfied. Hence, that corollary implies that $\Gamma^{\left(\delta\right)}$
forms a Banach frame (in the sense of Theorem \ref{thm:DiscreteBanachFrameTheorem})
for $\mathcal{B}_{s,\mu}^{p,q}\left(\R^{\dimension}\right)=\DecompSp{\CalQ}p{\ell_{\left(2^{js}\right)_{j}}^{q}}{v^{\left(\mu\right)}}$,
as soon as $0<\delta\leq\delta_{00}$, where (cf.\@ Lemma \ref{lem:SpecialProjection}
for the definition of $F_{0}$ and the estimate for $\vertiii{F_{0}}$
used here)
\[
\delta_{00}=\frac{1}{1+2\vertiii{F_{0}}^{2}}
\]
and, with $w=\left(2^{js}\right)_{j\in\N_{0}}$,
\begin{align*}
\vertiii{F_{0}} & \leq2^{\frac{1}{q}}C_{\mathscr{B},\Phi,v_{0},p}^{2}\cdot\vertiii{\smash{\Gamma_{\mathscr{B}}}}_{\ell_{w}^{q}\to\ell_{w}^{q}}^{2}\cdot\left(\vertiii{\smash{\overrightarrow{A}}}^{\max\left\{ 1,\frac{1}{p}\right\} }+\vertiii{\smash{\overrightarrow{B}}}^{\max\left\{ 1,\frac{1}{p}\right\} }\right)\cdot C_{3}\\
\left({\scriptstyle \text{Corollary }\ref{cor:BanachFrameSimplifiedCriteria}}\right) & \leq2^{\frac{1}{q_{0}}}C_{\mathscr{B},\Phi,v_{0},p}^{2}\cdot\vertiii{\smash{\Gamma_{\mathscr{B}}}}_{\ell_{w}^{q}\to\ell_{w}^{q}}^{2}\cdot4H_{7}\cdot C_{3}C_{4}\\
\left({\scriptstyle \text{eq. }\eqref{eq:WeightedSequenceSpaceClusteringMapNormEstimate}}\right) & \leq4\cdot2^{\frac{1}{q_{0}}}C_{\mathscr{B},\Phi,v_{0},p}^{2}\cdot\left[C_{\mathscr{B},\left(2^{js}\right)_{j\in\N_{0}}}\cdot N_{\mathscr{B}}^{1+\frac{1}{q}}\right]^{2}\cdot H_{7}\cdot C_{3}C_{4}\\
 & \leq4\cdot2^{\frac{1}{q_{0}}}C_{\mathscr{B},\Phi,v_{0},p}^{2}\cdot\left[C_{\mathscr{B},\left(2^{j}\right)_{j\in\N_{0}}}^{\left|s\right|}\cdot N_{\mathscr{B}}^{1+\frac{1}{q_{0}}}\right]^{2}\cdot H_{7}\cdot C_{3}C_{4}\\
 & \leq4\cdot2^{\frac{1}{q_{0}}}C_{\mathscr{B},\Phi,v_{0},p}^{2}\cdot\left[C_{\mathscr{B},\left(2^{j}\right)_{j\in\N_{0}}}^{\max\left\{ s_{1},-s_{0}\right\} }\cdot N_{\mathscr{B}}^{1+\frac{1}{q_{0}}}\right]^{2}\cdot H_{7}\cdot C_{3}C_{4},
\end{align*}
where
\[
C_{3}=\begin{cases}
\frac{\left(2^{16}\cdot768/\dimension^{\frac{3}{2}}\right)^{\frac{\dimension}{p}}}{2^{42}\cdot12^{\dimension}\cdot\dimension^{15}}\!\cdot\!\left(2^{52}\!\cdot\!\dimension^{\frac{25}{2}}\!\cdot\!\tilde{N}^{3}\right)^{\tilde{N}+1}\!\!\!\cdot\!N_{\mathscr{B}}^{2\left(\frac{1}{p}-1\right)}\!\left(1\!+\!R_{\mathscr{B}}C_{\mathscr{B}}\right)^{\dimension\left(\frac{4}{p}-1\right)}\!\!\cdot\Omega_{0}^{13K}\Omega_{1}^{13}\Omega_{2}^{\left(p,K\right)}, & \text{if }p<1,\\
\frac{1}{\sqrt{\dimension}\cdot2^{12+6\left\lceil K\right\rceil }}\cdot\left(2^{17}\cdot\dimension^{5/2}\cdot\tilde{N}\right)^{\left\lceil K\right\rceil +\dimension+2}\cdot\left(1+R_{\mathscr{B}}\right)^{\dimension}\cdot\Omega_{0}^{3K}\Omega_{1}^{3}\Omega_{2}^{\left(p,K\right)}, & \text{if }p\geq1,
\end{cases}
\]
with $\tilde{N}=\left\lceil K+\frac{\dimension+1}{\min\left\{ 1,p\right\} }\right\rceil \leq\left\lceil \mu_{0}+\frac{\dimension+1}{p_{0}}\right\rceil $
and 
\[
C_{4}=\Omega_{0}^{K}\Omega_{1}\cdot\dimension^{1/\min\left\{ 1,p\right\} }\cdot\left(4\cdot\dimension\right)^{1+2\left\lceil K+\frac{\dimension+\varepsilon}{\min\left\{ 1,p\right\} }\right\rceil }\cdot\left(\frac{s_{\dimension}}{\varepsilon}\right)^{1/\min\left\{ 1,p\right\} }\cdot\max_{\left|\alpha\right|\leq\left\lceil K+\frac{\dimension+\varepsilon}{\min\left\{ 1,p\right\} }\right\rceil }C^{\left(\alpha\right)},
\]
where the constants $C^{\left(\alpha\right)}=C^{\left(\alpha\right)}\left(\Phi\right)$
are as in Assumption \ref{assu:RegularPartitionOfUnity}.

To establish that $\delta_{0}$ can be chosen independently of $p,q,s,\mu$,
it thus suffices to estimate $C_{3}C_{4}$ and $C_{\mathscr{B},\Phi,v_{0},p}$
independently of these quantities. But above, we estimated $\Omega_{2}^{\left(p,K\right)}\leq\Omega_{3}$
with $\Omega_{3}$ independent of $p,q,s,\mu$. Since we also have
$K=\left|\mu\right|\leq\mu_{0}$ and $0\leq\frac{1}{p}\leq\frac{1}{p_{0}}$,
as well as $\Omega_{0}=1$ and $\Omega_{1}=2^{\left|\mu\right|}\leq2^{\mu_{0}}$,
it is straightforward to see that $C_{3}$ can be estimated independently
of $p,q,s,\mu$. The same arguments also allow us to estimate $C_{4}$
independently of these quantities.

Finally, Corollary \ref{cor:RegularBAPUsAreWeightedBAPUs} shows that
there is a suitable $\varrho\in\TestFunctionSpace{\R^{\dimension}}$
(depending only on $\mathscr{B}$) satisfying
\[
C_{\mathscr{B},\Phi,v_{0},p}\leq\Omega_{0}^{K}\Omega_{1}\cdot\left(4\cdot\dimension\right)^{1+2\left\lceil K+\frac{\dimension+\varepsilon}{p}\right\rceil }\cdot\left(\frac{s_{\dimension}}{\varepsilon}\right)^{1/p}\cdot2^{\!\left\lceil K+\frac{\dimension+\varepsilon}{p}\right\rceil }\cdot\lambda_{\dimension}\left(Q\right)\cdot\max_{\left|\alpha\right|\leq\left\lceil K+\frac{\dimension+\varepsilon}{p}\right\rceil }\left\Vert \partial^{\alpha}\varrho\right\Vert _{\sup}\cdot\max_{\left|\alpha\right|\leq\left\lceil K+\frac{\dimension+\varepsilon}{p}\right\rceil }C^{\left(\alpha\right)},
\]
where $Q:=\overline{\bigcup_{i\in\N_{0}}Q_{i}'}\subset\overline{B_{4}}\left(0\right)$.
As above, since $0\leq\frac{1}{p}\leq\frac{1}{p_{0}}$ and $K=\left|\mu\right|\leq\mu_{0}$,
it is then not hard to see that $C_{\mathscr{B},\Phi,v_{0},p}$ can
be estimated independently of $p,q,s,\mu$.

\medskip{}

It remains to show that the sequence space $\mathscr{C}_{p,q,s,\mu}$
is identical to the coefficient space $\ell_{\left(\left|\det T_{i}\right|^{\frac{1}{2}-\frac{1}{p}}\cdot w_{i}\right)_{i\in I}}^{q}\!\!\!\!\!\!\!\!\left(\left[\vphantom{\sum}\smash{C_{i}^{\left(\delta\right)}}\right]_{i\in I}\right)$
mentioned in Theorem \ref{thm:DiscreteBanachFrameTheorem}. To this
end, recall from equation (\ref{eq:CoefficientSpaceDefinition}) that
$C_{j}^{\left(\delta\right)}=\ell_{v^{\left(j,\delta\right)}}^{p}\left(\smash{\Z^{\dimension}}\right)$
with $v=v^{\left(\mu\right)}$ and
\begin{align*}
v_{k}^{\left(j,\delta\right)} & =v^{\left(\mu\right)}\left(\delta\cdot T_{j}^{-T}k\right)\\
 & =\left(1+\left|\delta\cdot T_{j}^{-T}k\right|\right)^{\mu}\\
 & =\left(1+\left|\delta\cdot k/2^{j}\right|\right)^{\mu}.
\end{align*}
But since $0<\delta\leq1$, we have $\delta\cdot\left(1+\left|k/2^{j}\right|\right)\leq1+\left|\delta\cdot\frac{k}{2^{j}}\right|\leq1+\left|k/2^{j}\right|$,
which implies
\[
\delta^{\mu_{0}}\cdot\left(1+\left|k/2^{j}\right|\right)^{\mu}\leq\delta^{\left|\mu\right|}\cdot\left(1+\left|k/2^{j}\right|\right)^{\mu}\leq v_{k}^{\left(j,\delta\right)}\leq\delta^{-\left|\mu\right|}\cdot\left(1+\left|k/2^{j}\right|\right)^{\mu}\leq\delta^{-\mu_{0}}\cdot\left(1+\left|k/2^{j}\right|\right)^{\mu}
\]
for all $k\in\Z^{\dimension}$, $j\in\N_{0}$ and $0<\delta\leq1$.
Finally, since $w_{i}=2^{si}$ and $\left|\det T_{i}\right|=2^{i\cdot\dimension}$
for $i\in\N_{0}$, we see
\[
\left|\det T_{i}\right|^{\frac{1}{2}-\frac{1}{p}}\cdot w_{i}=2^{i\left(s+\dimension\left(\frac{1}{2}-\frac{1}{p}\right)\right)}.
\]
Taken together, these considerations easily show $\ell_{\left(\left|\det T_{i}\right|^{\frac{1}{2}-\frac{1}{p}}\cdot w_{i}\right)_{i\in I}}^{q}\!\!\!\!\!\left(\left[\vphantom{\sum}\smash{C_{i}^{\left(\delta\right)}}\right]_{i\in I}\right)=\mathscr{C}_{p,q,s,\mu}$
with equivalent quasi-norms. Here, the implicit constant is allowed
to depend on $\delta$.
\end{proof}
Next, we derive concrete conditions on $\varphi,\psi$ which ensure
that the generated wavelet system yields atomic decompositions for
the (weighted) Besov spaces $\mathcal{B}_{s,\mu}^{p,q}\left(\R^{\dimension}\right)$.
\begin{prop}
\label{prop:BesovAtomicDecomposition}Let $p_{0},q_{0}\in\left(0,1\right]$,
$\varepsilon>0$, $\mu_{0}\geq0$ and $-\infty<s_{0}\leq s_{1}<\infty$.
Assume that $\varphi,\psi\in L^{1}\left(\R^{\dimension}\right)$ satisfy
the following conditions:

\begin{enumerate}
\item We have $\left\Vert \varphi\right\Vert _{K_{00}}<\infty$ and $\left\Vert \psi\right\Vert _{K_{00}}<\infty$
for $K_{00}:=\mu_{0}+\frac{\dimension}{p_{0}}+1$, where $\left\Vert g\right\Vert _{M}=\sup_{x\in\R^{\dimension}}\left(1+\left|x\right|\right)^{M}\cdot\left|g\left(x\right)\right|$.
\item We have $\widehat{\varphi},\widehat{\psi}\in C^{\infty}\left(\R^{\dimension}\right)$,
with all partial derivatives of $\widehat{\varphi},\widehat{\psi}$
being polynomially bounded.
\item We have $\widehat{\varphi}\left(\xi\right)\neq0$ for all $\xi\in\overline{B_{2}}\left(0\right)$
and $\widehat{\psi}\left(\xi\right)\neq0$ for all $\xi\in\overline{B_{4}}\left(0\right)\setminus B_{1/4}\left(0\right)$.
\item We have
\begin{align*}
\left|\partial^{\alpha}\widehat{\varphi}\left(\xi\right)\right| & \leq G_{1}\cdot\left(1+\left|\xi\right|\right)^{-L},\\
\left|\partial^{\alpha}\widehat{\psi}\left(\xi\right)\right| & \leq G_{2}\cdot\left(1+\left|\xi\right|\right)^{-L_{1}}\cdot\min\left\{ 1,\left|\xi\right|^{L_{2}}\right\} 
\end{align*}
for all $\alpha\in\N_{0}^{\dimension}$ with $\left|\alpha\right|\leq N_{0}:=\left\lceil \mu_{0}+\frac{\dimension+\varepsilon}{p_{0}}\right\rceil $,
all $\xi\in\R^{\dimension}$, suitable $G_{1},G_{2}>0$ and certain
$L,L_{1}\geq2\dimension+1+2\varepsilon$ and $L_{2}\geq0$ which furthermore
satisfy
\[
L>s_{1}+\kappa+\dimension+1+\varepsilon,\qquad L_{1}>s_{1}+\kappa+\dimension+1+\varepsilon,\qquad\text{ and }\qquad L_{2}>\vartheta_{0}\dimension-s_{0}
\]
for
\[
\vartheta_{0}:=\begin{cases}
0, & \text{if }p_{0}=1\\
\frac{1}{p_{0}}-1, & \text{if }p_{0}\in\left(0,1\right),
\end{cases}\qquad\text{ and }\qquad\kappa:=\begin{cases}
\left\lceil \mu_{0}+\dimension+\varepsilon\right\rceil , & \text{if }p_{0}=1,\\
\dimension+\mu_{0}+\left\lceil \mu_{0}+\frac{\dimension+\varepsilon}{p_{0}}\right\rceil , & \text{if }p_{0}\in\left(0,1\right).
\end{cases}
\]
\end{enumerate}
Then there is some $\delta_{0}=\delta_{0}\left(\dimension,p_{0},q_{0},\varepsilon,\mu_{0},s_{0},s_{1},\varphi,\psi\right)>0$
such that for each $0<\delta\leq\delta_{0}$, the family
\[
\Gamma^{\left(\delta\right)}=\left(2^{j\frac{\dimension}{2}}\cdot\gamma_{j}\left(2^{j}\bullet-\delta k\right)\right)_{j\in\N_{0},k\in\Z^{\dimension}},\qquad\text{ with }\qquad\gamma_{0}:=\varphi\quad\text{ and }\quad\gamma_{j}:=\psi\text{ for }j\in\N,
\]
forms an atomic decomposition for $\mathcal{B}_{s,\mu}^{p,q}\left(\R^{\dimension}\right)$,
for arbitrary $p,q\in\left(0,\infty\right]$ and $s,\mu\in\R$ satisfying
$p\geq p_{0}$, $q\geq q_{0}$ as well as $s_{0}\leq s\leq s_{1}$
and $\left|\mu\right|\leq\mu_{0}$.

Precisely, this means the following: With the space $\mathscr{C}_{p,q,s,\mu}\leq\Compl^{\N_{0}\times\Z^{\dimension}}$
as in Proposition \ref{prop:BesovBanachFrames}, the following are
true:

\begin{enumerate}
\item The \textbf{synthesis map}
\[
S^{\left(\delta\right)}:\mathscr{C}_{p,q,s,\mu}\to\mathcal{B}_{s,\mu}^{p,q}\left(\smash{\R^{\dimension}}\right),\left(\smash{c_{k}^{\left(i\right)}}\right)_{i\in\N_{0},k\in\Z^{\dimension}}\mapsto\sum_{i\in\N_{0}}\sum_{k\in\Z^{\dimension}}\left[c_{k}^{\left(i\right)}\cdot2^{i\frac{\dimension}{2}}\cdot\gamma_{i}\left(2^{i}\bullet-\delta k\right)\right]
\]
is well-defined and bounded for each $0<\delta\leq1$.

Convergence of the series has to be understood as described in the
remark after Theorem \ref{thm:AtomicDecomposition}.
\item For $0<\delta\leq\delta_{0}$, there is a bounded linear \textbf{coefficient
map}
\[
C^{\left(\delta\right)}:\mathcal{B}_{s,\mu}^{p,q}\left(\smash{\R^{\dimension}}\right)\to\mathscr{C}_{p,q,s,\mu}
\]
satisfying $S^{\left(\delta\right)}\circ C^{\left(\delta\right)}=\identity_{\mathcal{B}_{s,\mu}^{p,q}}$.
Furthermore, the action of $C^{\left(\delta\right)}$ on a given $f\in\mathcal{B}_{s,\mu}^{p,q}\left(\R^{\dimension}\right)$
is independent of the precise choice of $p,q,s,\mu$.\qedhere
\end{enumerate}
\end{prop}
\begin{proof}
Define $\widetilde{L}:=L-\left(\dimension+1+\varepsilon\right)$ and
$\widetilde{L_{1}}:=L_{1}-\left(\dimension+1+\varepsilon\right)$.
Now, an application of Lemma \ref{lem:ConvolutionFactorization} (with
$\gamma=\psi$, $N=N_{0}\geq\left\lceil \frac{\dimension+\varepsilon}{p_{0}}\right\rceil \geq\left\lceil \dimension+\varepsilon\right\rceil \geq\dimension+1$
and $\varrho\left(\xi\right):=G_{2}\cdot\left(1+\left|\xi\right|\right)^{-\widetilde{L_{1}}}\cdot\min\left\{ 1,\left|\xi\right|^{L_{2}}\right\} $,
where $\varrho\in L^{1}\left(\R^{\dimension}\right)$ since $\widetilde{L_{1}}\geq\dimension+\varepsilon$)
yields functions $\psi_{1},\psi_{2}\in L^{1}\left(\R^{\dimension}\right)$
with the following properties:

\begin{enumerate}
\item We have $\psi=\psi_{1}\ast\psi_{2}$.
\item We have $\psi_{2}\in C^{1}\left(\R^{\dimension}\right)$ with $H_{1}^{\left(M\right)}:=\left\Vert \psi_{2}\right\Vert _{M}+\left\Vert \nabla\psi_{2}\right\Vert _{M}<\infty$
for all $M\in\N_{0}$.
\item We have $\widehat{\psi_{1}},\widehat{\psi_{2}}\in C^{\infty}\left(\R^{\dimension}\right)$,
where all partial derivatives of these functions are polynomially
bounded.
\item We have $\left\Vert \psi_{1}\right\Vert _{N_{0}}<\infty$ and $\left\Vert \psi\right\Vert _{N_{0}}<\infty$.
In particular, since $N_{0}\geq\mu_{0}+\frac{\dimension+\varepsilon}{p_{0}}\geq\mu_{0}+\dimension+\varepsilon$,
we have $\psi_{1}\in L_{\left(1+\left|\mybullet\right|\right)^{\mu_{0}+\dimension+\varepsilon}}^{\infty}\left(\R^{\dimension}\right)\hookrightarrow L_{\left(1+\left|\mybullet\right|\right)^{\mu_{0}}}^{1}\left(\R^{\dimension}\right)\hookrightarrow L_{\left(1+\left|\mybullet\right|\right)^{K}}^{1}\left(\R^{\dimension}\right)$
for all $K=\left|\mu\right|\leq\mu_{0}$.
\item We have
\begin{equation}
\begin{split}\left|\partial^{\alpha}\widehat{\psi_{1}}\left(\xi\right)\right| & \leq2^{1+\dimension+4N_{0}}\cdot N_{0}!\cdot\left(1+\dimension\right)^{N_{0}}\cdot\varrho\left(\xi\right)\\
 & \leq H_{2}\cdot\left(1+\left|\xi\right|\right)^{-\widetilde{L_{1}}}\cdot\min\left\{ 1,\left|\xi\right|^{L_{2}}\right\} 
\end{split}
\label{eq:BesovAtomicWaveletConvolutionFactorEstimate}
\end{equation}
with $H_{2}:=G_{2}\cdot2^{1+\dimension+4N_{0}}\cdot N_{0}!\cdot\left(1+\dimension\right)^{N_{0}}$
for all $\xi\in\R^{\dimension}$ and all $\alpha\in\N_{0}^{\dimension}$
with $\left|\alpha\right|\leq N_{0}$.
\end{enumerate}
Likewise, another application of Lemma \ref{lem:ConvolutionFactorization}
(this time with $\gamma=\varphi$, $N=N_{0}\geq\dimension+1$ and
$\varrho\left(\xi\right):=G_{1}\cdot\left(1+\left|\xi\right|\right)^{-\widetilde{L}}$,
where $\varrho\in L^{1}\left(\R^{\dimension}\right)$, since $\widetilde{L}\geq\dimension+\varepsilon$)
yields certain functions $\varphi_{1},\varphi_{2}\in L^{1}\left(\R^{\dimension}\right)$
with the following properties:

\begin{enumerate}
\item We have $\varphi=\varphi_{1}\ast\varphi_{2}$.
\item We have $\varphi_{2}\in C^{1}\left(\R^{\dimension}\right)$ with $H_{3}^{\left(M\right)}:=\left\Vert \varphi_{2}\right\Vert _{M}+\left\Vert \nabla\varphi_{2}\right\Vert _{M}<\infty$
for all $M\in\N_{0}$.
\item We have $\widehat{\varphi_{1}},\widehat{\varphi_{2}}\in C^{\infty}\left(\R^{\dimension}\right)$,
where all partial derivatives of these functions are polynomially
bounded.
\item We have $\left\Vert \varphi_{1}\right\Vert _{N_{0}}<\infty$ and $\left\Vert \varphi\right\Vert _{N_{0}}<\infty$.
As for $\psi_{1}$, this implies $\varphi_{1}\in L_{\left(1+\left|\mybullet\right|\right)^{\mu_{0}}}^{1}\left(\R^{\dimension}\right)\hookrightarrow L_{\left(1+\left|\mybullet\right|\right)^{K}}^{1}\left(\R^{\dimension}\right)$
for all $K=\left|\mu\right|\leq\mu_{0}$.
\item We have
\begin{equation}
\begin{split}\left|\partial^{\alpha}\widehat{\varphi_{1}}\left(\xi\right)\right| & \leq2^{1+\dimension+4N_{0}}\cdot N_{0}!\cdot\left(1+\dimension\right)^{N_{0}}\cdot\varrho\left(\xi\right)\\
 & \leq H_{4}\cdot\left(1+\left|\xi\right|\right)^{-\widetilde{L}}
\end{split}
\label{eq:BesovAtomicScalingConvolutionFactorEstimate}
\end{equation}
with $H_{4}:=G_{1}\cdot2^{1+\dimension+4N_{0}}\cdot N_{0}!\cdot\left(1+\dimension\right)^{N_{0}}$
for all $\xi\in\R^{\dimension}$ and all $\alpha\in\N_{0}^{\dimension}$
with $\left|\alpha\right|\leq N_{0}$.
\end{enumerate}
Now, set $\gamma_{0}:=\varphi$ and $\gamma_{0,\ell}:=\varphi_{\ell}$,
as well as $\gamma_{j}:=\psi$ and $\gamma_{j,\ell}:=\psi_{\ell}$
for $\ell\in\left\{ 1,2\right\} $ and $j\in\N$.

As a further preparation, set $\gamma_{1}^{\left(0\right)}:=\varphi$
and $\gamma_{2}^{\left(0\right)}:=\psi$, as well as $n_{0}:=1$ and
$n_{j}:=2$ for $j\in\N$, so that $\gamma_{j}=\gamma_{n_{j}}^{\left(0\right)}$
for all $j\in\N_{0}$. Then, in the notation of Lemma \ref{lem:GammaCoversOrbitAssumptionSimplified},
we have $Q^{\left(1\right)}=B_{2}\left(0\right)$ and $Q^{\left(2\right)}=B_{4}\left(0\right)\setminus\overline{B_{1/4}}\left(0\right)$,
cf.\@ the proof of Proposition \ref{prop:BesovBanachFrames}. Exactly
as in that proof, we see that all prerequisites of Lemma \ref{lem:GammaCoversOrbitAssumptionSimplified}
are satisfied, so that the family $\Gamma=\left(\gamma_{i}\right)_{i\in\N_{0}}$
satisfies Assumption \ref{assu:GammaCoversOrbit}, where we furthermore
have $\Omega_{2}^{\left(p,K\right)}\leq\Omega_{3}$ for all $K\leq\mu_{0}$
and all $p\geq p_{0}$, for a suitable constant $\Omega_{3}=\Omega_{3}\left(\mathscr{B},\varphi,\psi,p_{0},\mu_{0},\dimension\right)>0$.
Observe that indeed $K=\left|\mu\right|\leq\mu_{0}$ in the cases
which are of interest to us.

\medskip{}

Now, let $p,q,s,\mu$ be as in the statement of the proposition. We
want to verify the prerequisites of Corollary \ref{cor:AtomicDecompositionSimplifiedCriteria}
for the choices which we just made. For most of these assumptions,
this is not hard:

\begin{enumerate}
\item All $\gamma_{i},\gamma_{i,1},\gamma_{i,2}$ are measurable functions,
as required.
\item Since $K=\left|\mu\right|\leq\mu_{0}$, we have $\gamma_{i,1}=\psi_{1}\in L_{\left(1+\left|\mybullet\right|\right)^{K}}^{1}\left(\R^{\dimension}\right)$
for any $i\in\N$ and also $\gamma_{0,1}=\varphi_{1}\in L_{\left(1+\left|\mybullet\right|\right)^{K}}^{1}\left(\R^{\dimension}\right)$,
as seen above.
\item We have $\gamma_{i,2}=\psi_{2}\in C^{1}\left(\R^{\dimension}\right)$
for any $i\in\N$ and also $\gamma_{0,2}=\varphi_{2}\in C^{1}\left(\R^{\dimension}\right)$,
as noted above.
\item With $K_{0}:=K+\frac{\dimension}{\min\left\{ 1,p\right\} }+1\leq\mu_{0}+\frac{\dimension}{p_{0}}+1=K_{00}$,
we have
\begin{equation}
\begin{split}\Omega_{4}^{\left(p,K\right)} & =\sup_{i\in I}\left\Vert \gamma_{i,2}\right\Vert _{K_{0}}+\sup_{i\in I}\left\Vert \nabla\gamma_{i,2}\right\Vert _{K_{0}}\\
 & \leq\max\left\{ \left\Vert \varphi_{2}\right\Vert _{K_{00}},\left\Vert \psi_{2}\right\Vert _{K_{00}}\right\} +\max\left\{ \left\Vert \nabla\varphi_{2}\right\Vert _{K_{00}},\left\Vert \nabla\psi_{2}\right\Vert _{K_{00}}\right\} \\
 & \leq H_{1}^{\left(\left\lceil K_{00}\right\rceil \right)}+H_{3}^{\left(\left\lceil K_{00}\right\rceil \right)}=:H_{5}<\infty.
\end{split}
\label{eq:BesovAtomicDecompositionOmega4Estimate}
\end{equation}
\item We have $\left\Vert \gamma_{i}\right\Vert _{K_{0}}\leq\left\Vert \gamma_{i}\right\Vert _{K_{00}}<\infty$
for all $i\in\N_{0}$, since $\left\Vert \varphi\right\Vert _{K_{00}}<\infty$
and $\left\Vert \psi\right\Vert _{K_{00}}<\infty$ by assumption.
\item We have $\gamma_{i}=\psi=\psi_{1}\ast\psi_{2}=\gamma_{i,1}\ast\gamma_{i,2}$
for all $i\in\N$ and likewise $\gamma_{0}=\varphi=\varphi_{1}\ast\varphi_{2}=\gamma_{0,1}\ast\gamma_{0,2}$.
\item We have $\widehat{\gamma_{i,1}},\widehat{\gamma_{i,2}}\in C^{\infty}\left(\R^{\dimension}\right)$
for all $i\in\N_{0}$, and all partial derivatives of these functions
are polynomially bounded.
\item As we showed above, the family $\Gamma=\left(\gamma_{i}\right)_{i\in\N_{0}}$
satisfies Assumption \ref{assu:GammaCoversOrbit}.
\end{enumerate}
Hence, the only prerequisite of Corollary \ref{cor:AtomicDecompositionSimplifiedCriteria}
which still needs to be verified is that

\[
K_{1}:=\sup_{i\in\N_{0}}\,\sum_{j\in\N_{0}}N_{i,j}<\infty\qquad\text{ and }\qquad K_{2}:=\sup_{j\in\N_{0}}\:\sum_{i\in\N_{0}}N_{i,j}<\infty,
\]
where
\[
N_{i,j}:=\left(\frac{w_{i}}{w_{j}}\cdot\left(\left|\det T_{j}\right|\big/\left|\det T_{i}\right|\right)^{\vartheta}\right)^{\tau}\cdot\left(1+\left\Vert T_{j}^{-1}T_{i}\right\Vert \right)^{\sigma}\cdot\left(\left|\det T_{i}\right|^{-1}\cdot\int_{Q_{i}}\max_{\left|\alpha\right|\leq N}\left|\left(\partial^{\alpha}\widehat{\gamma_{j,1}}\right)\left(S_{j}^{-1}\xi\right)\right|\d\xi\right)^{\tau},
\]
with
\begin{align*}
N & =\left\lceil K+\frac{\dimension+\varepsilon}{\min\left\{ 1,p\right\} }\right\rceil ,\\
\tau & =\min\left\{ 1,p,q\right\} ,\\
\sigma & =\begin{cases}
\min\left\{ 1,q\right\} \cdot\left\lceil K+\dimension+\varepsilon\right\rceil , & \text{if }p\in\left[1,\infty\right],\\
\min\left\{ p,q\right\} \cdot\left(\frac{\dimension}{p}+K+\left\lceil K+\frac{\dimension+\varepsilon}{p}\right\rceil \right), & \text{if }p\in\left(0,1\right),
\end{cases}\\
\vartheta & =\begin{cases}
0, & \text{if }p\in\left[1,\infty\right],\\
\frac{1}{p}-1, & \text{if }p\in\left(0,1\right).
\end{cases}
\end{align*}

To prove this, we begin with several auxiliary observations: First
of all, we observe $N\leq N_{0}$. Furthermore, we have $\vartheta\leq\vartheta_{0}$,
since $p_{0}=1$ implies $p\in\left[1,\infty\right]$. Finally, we
also have
\begin{equation}
\begin{split}\frac{\sigma}{\tau}-\vartheta\dimension & =\begin{cases}
\left\lceil K+\dimension+\varepsilon\right\rceil , & \text{if }p\in\left[1,\infty\right],\\
\frac{\dimension}{p}+K+\left\lceil K+\frac{\dimension+\varepsilon}{p}\right\rceil -\dimension\left(\frac{1}{p}-1\right), & \text{if }p\in\left(0,1\right)
\end{cases}\\
 & \leq\begin{cases}
\left\lceil \mu_{0}+\dimension+\varepsilon\right\rceil , & \text{if }p_{0}=1,\\
\dimension+\mu_{0}+\left\lceil \mu_{0}+\frac{\dimension+\varepsilon}{p_{0}}\right\rceil , & \text{if }p_{0}\in\left(0,1\right)
\end{cases}\\
 & =\kappa,
\end{split}
\label{eq:BesovAtomicDecompositionKappaSignificance}
\end{equation}
where we used that $K=\left|\mu\right|\leq\mu_{0}$ and also that
$p_{0}=1$ entails $p\in\left[1,\infty\right]$.

We divide our estimates of $N_{i,j}$ into two main cases. The first
case is $i\in\N$. Here, we recall from the proof of Proposition \ref{prop:BesovBanachFrames}
(cf.\@ equation (\ref{eq:BesovAnnulusMeasureEstimate})) for $Q_{i}=\left\{ \xi\in\R^{\dimension}\with2^{i-2}<\left|\xi\right|<2^{i+2}\right\} $
that $\lambda_{\dimension}\left(Q_{i}\right)\leq2^{i\cdot\dimension}\cdot4^{\dimension}v_{\dimension}$.
Since we also have $N\leq N_{0}$, we conclude
\begin{equation}
N_{i,j}\leq\left[2^{s\left(i-j\right)+\vartheta\dimension\left(j-i\right)}\right]^{\tau}\cdot\left(1+2^{i-j}\right)^{\sigma}\cdot\left(4^{\dimension}v_{\dimension}\cdot\sup_{2^{i-2}<\left|\xi\right|<2^{i+2}}\;\max_{\left|\alpha\right|\leq N_{0}}\left|\left(\partial^{\alpha}\widehat{\gamma_{j,1}}\right)\left(\xi/2^{j}\right)\right|\right)^{\tau}.\label{eq:BesovAtomicDecompositionMainTerm}
\end{equation}
To further estimate this term, we consider two subcases:

\textbf{Case 1}: We have $j\in\N$. Here, we again consider two subcases:

\begin{enumerate}
\item We have $i\leq j$. For $2^{i-2}<\left|\xi\right|<2^{i+2}$, equation
(\ref{eq:BesovAtomicWaveletConvolutionFactorEstimate}) implies because
of $\widetilde{L_{1}}\geq0$ and $L_{2}\geq0$ that
\begin{align*}
\max_{\left|\alpha\right|\leq N_{0}}\left|\left(\partial^{\alpha}\widehat{\gamma_{j,1}}\right)\left(\xi/2^{j}\right)\right|=\max_{\left|\alpha\right|\leq N_{0}}\left|\left(\partial^{\alpha}\widehat{\psi_{1}}\right)\left(\xi/2^{j}\right)\right| & \leq H_{2}\cdot\min\left\{ 1,\left|\xi/2^{j}\right|^{L_{2}}\right\} \\
 & \leq H_{2}\cdot4^{L_{2}}\cdot2^{L_{2}\left(i-j\right)}\\
 & =4^{L_{2}}H_{2}\cdot2^{-L_{2}\left|i-j\right|}.
\end{align*}
In combination with equation (\ref{eq:BesovAtomicDecompositionMainTerm}),
we get
\[
N_{i,j}\leq2^{\sigma}\cdot\left(4^{\dimension+L_{2}}v_{\dimension}\cdot H_{2}\right)^{\tau}\cdot2^{-\tau\left|i-j\right|\left(L_{2}+s-\vartheta\dimension\right)}.
\]
\item We have $j\leq i$. Here, equation (\ref{eq:BesovAtomicWaveletConvolutionFactorEstimate})
implies for $2^{i-2}<\left|\xi\right|<2^{i+2}$ that
\begin{align*}
\max_{\left|\alpha\right|\leq N_{0}}\left|\left(\partial^{\alpha}\widehat{\gamma_{j,1}}\right)\left(\xi/2^{j}\right)\right|=\max_{\left|\alpha\right|\leq N_{0}}\left|\left(\partial^{\alpha}\widehat{\psi_{1}}\right)\left(\xi/2^{j}\right)\right| & \leq H_{2}\cdot\left(1+\left|\xi/2^{j}\right|\right)^{-\widetilde{L_{1}}}\\
 & \leq H_{2}\cdot\left(2^{i-j}/4\right)^{-\widetilde{L_{1}}}\\
 & =4^{\widetilde{L_{1}}}H_{2}\cdot2^{-\widetilde{L_{1}}\left|i-j\right|}.
\end{align*}
In combination with equation (\ref{eq:BesovAtomicDecompositionMainTerm}),
we get
\[
N_{i,j}\leq2^{\sigma}\cdot\left(4^{\dimension+\widetilde{L_{1}}}v_{\dimension}\cdot H_{2}\right)^{\tau}\cdot2^{-\tau\left|i-j\right|\left(-\frac{\sigma}{\tau}-s+\vartheta\dimension+\widetilde{L_{1}}\right)}.
\]
\end{enumerate}
\textbf{Case 2}: We have $j=0$. Here, equation (\ref{eq:BesovAtomicScalingConvolutionFactorEstimate})
shows for $2^{i-2}<\left|\xi\right|<2^{i+2}$ that
\begin{align*}
\max_{\left|\alpha\right|\leq N_{0}}\left|\left(\partial^{\alpha}\widehat{\gamma_{j,1}}\right)\left(\xi/2^{j}\right)\right|=\max_{\left|\alpha\right|\leq N_{0}}\left|\partial^{\alpha}\widehat{\varphi_{1}}\left(\xi\right)\right| & \leq H_{4}\cdot\left(1+\left|\xi\right|\right)^{-\widetilde{L}}\\
 & \leq H_{4}\cdot\left(2^{i}/4\right)^{-\widetilde{L}}\\
 & =4^{\widetilde{L}}H_{4}\cdot2^{-i\widetilde{L}}.
\end{align*}
In combination with equation (\ref{eq:BesovAtomicDecompositionMainTerm}),
this implies
\[
N_{i,j}\leq2^{\sigma}\cdot\left(4^{\dimension+\widetilde{L}}v_{\dimension}\cdot H_{4}\right)^{\tau}\cdot2^{-\tau i\left(\widetilde{L}-s+\vartheta\dimension-\frac{\sigma}{\tau}\right)}=2^{\sigma}\cdot\left(4^{\dimension+\widetilde{L}}v_{\dimension}\cdot H_{4}\right)^{\tau}\cdot2^{-\tau\left|i-j\right|\left(\widetilde{L}-s+\vartheta\dimension-\frac{\sigma}{\tau}\right)}.
\]
These are the desired estimates for the case $i\in\N$.

If otherwise $i=0$, so that $Q_{i}=T_{0}Q_{0}'+b_{0}=Q_{0}'=B_{2}\left(0\right)$,
estimate (\ref{eq:BesovAtomicDecompositionMainTerm}) takes the modified
form
\begin{equation}
\begin{split}N_{i,j} & \leq2^{\tau j\left(\vartheta\dimension-s\right)}\cdot\left(1+2^{-j}\right)^{\sigma}\cdot\left(\int_{B_{2}\left(0\right)}\max_{\left|\alpha\right|\leq N}\left|\left(\partial^{\alpha}\widehat{\gamma_{j,1}}\right)\left(S_{j}^{-1}\xi\right)\right|\d\xi\right)^{\tau}\\
 & \leq2^{\sigma}\cdot\left[\lambda_{\dimension}\left(B_{2}\left(0\right)\right)\right]^{\tau}\cdot2^{\tau j\left(\vartheta\dimension-s\right)}\cdot\sup_{\left|\xi\right|<2}\,\left(\max_{\left|\alpha\right|\leq N_{0}}\left|\left(\partial^{\alpha}\widehat{\gamma_{j,1}}\right)\left(\xi/2^{j}\right)\right|\right)^{\tau}.
\end{split}
\label{eq:BesovAtomicDecompositionModifiedMainTerm}
\end{equation}
To further estimate this quantity, we again distinguish two cases:

\textbf{Case 1}: We have $j\in\N$. In this case, equation (\ref{eq:BesovAtomicWaveletConvolutionFactorEstimate})
shows for $\left|\xi\right|<2$ that
\begin{align*}
\max_{\left|\alpha\right|\leq N_{0}}\left|\left(\partial^{\alpha}\widehat{\gamma_{j,1}}\right)\left(\xi/2^{j}\right)\right|=\max_{\left|\alpha\right|\leq N_{0}}\left|\left(\partial^{\alpha}\widehat{\psi_{1}}\right)\left(\xi/2^{j}\right)\right| & \leq H_{2}\cdot\min\left\{ 1,\left|\xi/2^{j}\right|^{L_{2}}\right\} \\
 & \leq2^{L_{2}}H_{2}\cdot2^{-jL_{2}}.
\end{align*}
In combination with equation (\ref{eq:BesovAtomicDecompositionModifiedMainTerm}),
we get
\begin{align*}
N_{i,j} & \leq2^{\sigma}\cdot\left[2^{L_{2}+\dimension}v_{\dimension}\cdot H_{2}\right]^{\tau}\cdot2^{\tau j\left(\vartheta\dimension-s-L_{2}\right)}\\
 & =2^{\sigma}\cdot\left[2^{L_{2}+\dimension}v_{\dimension}\cdot H_{2}\right]^{\tau}\cdot2^{-\tau\left|i-j\right|\left(L_{2}+s-\vartheta\dimension\right)}.
\end{align*}

\textbf{Case 2}: We have $j=0$. In this case, equation (\ref{eq:BesovAtomicScalingConvolutionFactorEstimate})
shows for $\left|\xi\right|<2$ that
\[
\max_{\left|\alpha\right|\leq N_{0}}\left|\left(\partial^{\alpha}\widehat{\gamma_{j,1}}\right)\left(\xi/2^{j}\right)\right|=\max_{\left|\alpha\right|\leq N_{0}}\left|\left(\partial^{\alpha}\widehat{\varphi_{1}}\right)\left(\xi\right)\right|\leq H_{4},
\]
since $\widetilde{L}\geq0$. In combination with equation (\ref{eq:BesovAtomicDecompositionModifiedMainTerm}),
this yields
\[
N_{i,j}\leq2^{\sigma}\cdot\left[2^{\dimension}v_{\dimension}\cdot H_{4}\right]^{\tau}=2^{\sigma}\cdot\left[2^{\dimension}v_{\dimension}\cdot H_{4}\right]^{\tau}\cdot2^{-\tau\left|i-j\right|\zeta}
\]
for arbitrary $\zeta\in\R$, since $\left|i-j\right|=0$.

\medskip{}

All in all, our considerations show that there is some constant $H_{6}$
which is independent of $p,q,s,\mu$, such that
\[
N_{i,j}\leq2^{\sigma}\cdot H_{6}^{\tau}\cdot2^{-\tau\left|i-j\right|\lambda}\quad\text{ where }\quad\lambda:=\min\left\{ 1,\,L_{2}+s-\vartheta\dimension,\,\widetilde{L}-s+\vartheta\dimension-\frac{\sigma}{\tau},\,\widetilde{L_{1}}-s+\vartheta\dimension-\frac{\sigma}{\tau}\right\} .
\]
But in view of equation (\ref{eq:BesovAtomicDecompositionKappaSignificance})
and since $s_{0}\leq s\leq s_{1}$, as well as $\vartheta\leq\vartheta_{0}$,
we have
\[
\lambda\geq\min\left\{ 1,\,L_{2}+s_{0}-\vartheta_{0}\dimension,\,\widetilde{L}-s_{1}-\kappa,\,\widetilde{L_{1}}-s_{1}-\kappa\right\} =:\lambda_{0}>0,
\]
by our assumptions regarding $L,L_{1},L_{2}$. Note that $\lambda_{0}$
is independent of $p,q,s,\mu$.

All in all, we thus get, for $K_{1},K_{2}$ as in Corollary \ref{cor:AtomicDecompositionSimplifiedCriteria},
\begin{align*}
K_{1}^{1/\tau}=\sup_{i\in\N_{0}}\,\left(\sum_{j\in\N_{0}}N_{i,j}\right)^{1/\tau} & \leq2^{\sigma/\tau}\cdot H_{6}\cdot\sup_{i\in\N_{0}}\,\left(\sum_{j\in\N_{0}}2^{-\tau\lambda\left|i-j\right|}\right)^{1/\tau}\\
\left({\scriptstyle \text{since }\ell^{\tau_{0}}\hookrightarrow\ell^{\tau}\text{ is norm-decreasing for }\tau_{0}:=\min\left\{ 1,p_{0},q_{0}\right\} \leq\tau}\right) & \leq2^{\sigma/\tau}\cdot H_{6}\cdot\sup_{i\in\N_{0}}\,\left(\sum_{j\in\N_{0}}2^{-\tau_{0}\lambda\left|i-j\right|}\right)^{1/\tau_{0}}\\
\left({\scriptstyle \text{since }\lambda\geq\lambda_{0}\text{ and with }\ell=i-j}\right) & \leq2^{\sigma/\tau}\cdot H_{6}\cdot\left(\sum_{\ell\in\Z}2^{-\tau_{0}\lambda_{0}\left|\ell\right|}\right)^{1/\tau_{0}}\\
 & \leq2^{\frac{\sigma}{\tau}+\frac{1}{\tau_{0}}}\cdot H_{6}\cdot\left(\sum_{\ell=0}^{\infty}2^{-\tau_{0}\lambda_{0}\ell}\right)^{1/\tau_{0}}\\
\left({\scriptstyle \text{since }\frac{\sigma}{\tau}\leq\kappa+\vartheta\dimension\leq\kappa+\vartheta_{0}\dimension\text{ by eq. }\eqref{eq:BesovAtomicDecompositionKappaSignificance}}\right) & \leq2^{\kappa+\vartheta_{0}\dimension+\frac{1}{\tau_{0}}}\cdot H_{6}\cdot\left(\frac{1}{1-2^{-\tau_{0}\lambda_{0}}}\right)^{1/\tau_{0}}=:H_{7},
\end{align*}
where $H_{7}$ is independent of $p,q,s,\mu$. Precisely the same
estimate also yields $K_{2}^{1/\tau}\leq H_{7}$.

All in all, we see that Corollary \ref{cor:AtomicDecompositionSimplifiedCriteria}
is applicable, so that the operator $\overrightarrow{C}$ from Assumption
\ref{assu:AtomicDecompositionAssumption} satisfies $\vertiii{\smash{\overrightarrow{C}}}^{\max\left\{ 1,\frac{1}{p}\right\} }\leq\Omega\cdot\left(K_{1}^{1/\tau}+K_{2}^{1/\tau}\right)\leq2\Omega H_{7}$
for $\Omega=\Omega_{0}^{K}\Omega_{1}\cdot\left(4\cdot\dimension\right)^{1+2\left\lceil K+\frac{\dimension+\varepsilon}{\min\left\{ 1,p\right\} }\right\rceil }\cdot\left(\frac{s_{\dimension}}{\varepsilon}\right)^{1/\min\left\{ 1,p\right\} }\cdot\max_{\left|\alpha\right|\leq N}C^{\left(\alpha\right)}$,
where the constants $C^{\left(\alpha\right)}$ are as in Assumption
\ref{assu:RegularPartitionOfUnity}. Since $N\leq N_{0}$, $K=\left|\mu\right|\leq\mu_{0}$
and $\frac{1}{\min\left\{ 1,p\right\} }\leq\frac{1}{p_{0}}$, as well
as $\Omega_{0}=1$ and $\Omega_{1}=2^{\left|\mu\right|}\leq2^{\mu_{0}}$,
it is not hard to see $\Omega\leq H_{8}$, where $H_{8}$ is independent
of $p,q,s,\mu$.

Corollary \ref{cor:AtomicDecompositionSimplifiedCriteria} ensures
that Theorem \ref{thm:AtomicDecomposition} is applicable, i.e., the
family $\Gamma^{\left(\delta\right)}$ from the statement of the current
proposition is indeed an atomic decomposition for $\mathcal{B}_{s,\mu}^{p,q}\left(\R^{\dimension}\right)=\DecompSp{\mathscr{B}}p{\ell_{\left(2^{js}\right)_{j\in\N_{0}}}^{q}}{v^{\left(\mu\right)}}$
as soon as $0<\delta\leq\min\left\{ 1,\delta_{00}\right\} $, where
$\delta_{00}$ is defined by
\[
\delta_{00}^{-1}\!:=\!\begin{cases}
\!\frac{2s_{\dimension}}{\sqrt{\dimension}}\cdot\left(2^{17}\!\cdot\!\dimension^{2}\!\cdot\!\left(K\!+\!2\!+\!\dimension\right)\right)^{\!K+\dimension+3}\!\!\!\cdot\left(1\!+\!R_{\mathscr{B}}\right)^{\dimension+1}\cdot\Omega_{0}^{4K}\Omega_{1}^{4}\Omega_{2}^{\left(p,K\right)}\Omega_{4}^{\left(p,K\right)}\cdot\vertiii{\smash{\overrightarrow{C}}}\,, & \text{if }p\geq1,\\
\frac{\left(2^{14}/\dimension^{\frac{3}{2}}\right)^{\!\frac{\dimension}{p}}}{2^{45}\cdot\dimension^{17}}\!\cdot\!\left(\frac{s_{\dimension}}{p}\right)^{\!\frac{1}{p}}\left(2^{68}\!\cdot\!\dimension^{14}\!\cdot\!\left(K\!+\!1\!+\!\frac{\dimension+1}{p}\right)^{3}\right)^{\!K+2+\frac{\dimension+1}{p}}\!\!\!\cdot\!\left(1\!+\!R_{\mathscr{B}}\right)^{1+\frac{3\dimension}{p}}\!\cdot\!\Omega_{0}^{16K}\Omega_{1}^{16}\Omega_{2}^{\left(p,K\right)}\Omega_{4}^{\left(p,K\right)}\cdot\vertiii{\smash{\overrightarrow{C}}}^{\frac{1}{p}}, & \text{if }p<1.
\end{cases}
\]
The verification that the discrete sequence space from Theorem \ref{thm:AtomicDecomposition}
coincides with $\mathscr{C}_{p,q,s,\mu}$ is exactly as in the proof
of Proposition \ref{prop:BesovBanachFrames}. Hence, to complete the
proof, we only have to verify $\delta_{00}^{-1}\leq\delta_{0}^{-1}$,
where $\delta_{0}>0$ is independent of $p,q,s,\mu$.

But above, we showed $\Omega_{2}^{\left(p,K\right)}\leq\Omega_{3}$,
with $\Omega_{3}$ independent of $p,q,s,\mu$, since $K=\left|\mu\right|\leq\mu_{0}$
and $p\geq p_{0}$. Furthermore, equation (\ref{eq:BesovAtomicDecompositionOmega4Estimate})
shows $\Omega_{4}^{\left(p,K\right)}\leq H_{5}$ for $K=\left|\mu\right|\leq\mu_{0}$
and with $H_{5}$ independent of $p,q,s,\mu$. Using these estimates,
the estimate for $\vphantom{\overrightarrow{C}}\vertiii{\smash{\overrightarrow{C}}}$
from above, and the straightforward inequalities $0\leq\frac{1}{p}\leq\frac{1}{p_{0}}$
and $K=\left|\mu\right|\leq\mu_{0}$, as well as $\Omega_{0}=1$ and
$\Omega_{2}=2^{\left|\mu\right|}\leq2^{\mu_{0}}$, we see that indeed
$\delta_{00}^{-1}\leq\delta_{0}^{-1}$ for some $\delta_{0}>0$ independent
of $p,q,s,\mu$.
\end{proof}
\begin{rem}
\label{rem:BesovRemarks}We close this section by showing that our
results indeed imply the existence of compactly supported Banach frames
and atomic decompositions for Besov spaces. Finally, we compare our
results with the literature.

\begin{itemize}[leftmargin=0.4cm]
\item The conditions in Propositions \ref{prop:BesovBanachFrames} and
\ref{prop:BesovAtomicDecomposition} are still not completely straightforward
to verify. Thus, let $k,N\in\N$ and $L_{1},L_{2}\geq0$ be arbitrary.
We will show how one can construct a compactly supported function
$\psi\in C_{c}^{k}\left(\R^{\dimension}\right)$ satisfying $\widehat{\psi}\left(\xi\right)\neq0$
for $\xi\in\overline{B_{4}}\left(0\right)\setminus B_{1/4}\left(0\right)$
as well as
\[
\left|\partial^{\alpha}\widehat{\psi}\left(\xi\right)\right|\lesssim\left(1+\left|\xi\right|\right)^{-L_{1}}\cdot\min\left\{ 1,\left|\xi\right|^{L_{2}}\right\} \qquad\forall\alpha\in\N_{0}^{\dimension}\text{ with }\left|\alpha\right|\leq N.
\]
To this end, let $\ell:=\max\left\{ k,\left\lceil L_{1}\right\rceil \right\} $
and $L_{3}:=\max\left\{ 0,\left\lceil \left(L_{2}-N-1\right)/2\right\rceil \right\} $
and choose $\psi_{0}\in C_{c}^{\ell+2\left(L_{3}+N+1\right)}\left(\R^{\dimension}\right)$
with $\psi_{0}\geq0$ and $\psi_{0}\not\equiv0$. This implies $\widehat{\psi_{0}}\left(0\right)=\left\Vert \psi_{0}\right\Vert _{L^{1}}>0$.
By continuity of $\widehat{\psi_{0}}$, there is thus some $c_{0}>0$
and some $\varepsilon>0$ satisfying $\left|\smash{\widehat{\psi_{0}}}\left(\xi\right)\right|\geq c_{0}$
for $\left|\xi\right|\leq\varepsilon$.

Define $\psi_{1}:=\psi_{0}\circ\frac{4}{\varepsilon}\identity$ and
note $\vphantom{\widehat{\psi_{1}}}\left|\smash{\widehat{\psi_{1}}}\left(\xi\right)\right|=\left(\varepsilon/4\right)^{\dimension}\cdot\left|\smash{\widehat{\psi_{0}}}\left(\frac{\varepsilon}{4}\xi\right)\right|\geq c_{0}\cdot\left(\varepsilon/4\right)^{\dimension}=:c_{1}$
as long as $\left|\xi\right|\leq4$. Next, set 
\[
\psi:=\Delta^{L_{3}+N+1}\psi_{1}\in C_{c}^{\ell}\left(\smash{\R^{\dimension}}\right)\subset C_{c}^{k}\left(\smash{\R^{\dimension}}\right),
\]
where $\Delta=\sum_{j=1}^{\dimension}\frac{\partial^{2}}{\partial x_{j}^{2}}$
denotes the Laplace operator. An easy calculation using partial integration
shows $\Fourier\left[\Delta g\right]\left(\xi\right)=-4\pi^{2}\cdot\left|\xi\right|^{2}\cdot\widehat{g}\left(\xi\right)$
for $g\in C_{c}^{2}\left(\R^{\dimension}\right)$, so that we get
\begin{equation}
\widehat{\psi}\left(\xi\right)=\left(-4\pi^{2}\right)^{L_{3}+N+1}\cdot\left|\xi\right|^{2\left(L_{3}+N+1\right)}\cdot\widehat{\psi_{1}}\left(\xi\right)\qquad\forall\xi\in\R^{\dimension}.\label{eq:BesovExplicitConstructionDerivative}
\end{equation}
In particular, $\widehat{\psi}\left(\xi\right)\in o\left(\smash{\left|\xi\right|^{2\left(L_{3}+N\right)+1}}\right)$
as $\xi\to0$. Furthermore, since $\psi\in C_{c}\left(\R^{\dimension}\right)$,
we have $\widehat{\psi}\in C^{\infty}\left(\R^{\dimension}\right)$,
so that Lemma \ref{lem:HighOrderVanishingYieldVanishingDerivatives}
shows $\partial^{\alpha}\widehat{\psi}\left(\xi\right)\in\vphantom{\left|\xi\right|^{2L_{3}}}o\left(\smash{\left|\xi\right|^{2\left(L_{3}+N\right)+1-\left|\alpha\right|}}\right)\subset o\left(\smash{\left|\xi\right|^{2L_{3}+N+1}}\right)\subset o\left(\smash{\left|\xi\right|^{L_{2}}}\right)$
as $\xi\to0$, for $\left|\alpha\right|\leq N$. By continuity of
$\partial^{\alpha}\widehat{\psi}$, this implies that there is a constant
$C>0$ satisfying 
\[
\left|\partial^{\alpha}\widehat{\psi}\left(\xi\right)\right|\leq C\cdot\left|\xi\right|^{L_{2}}\leq2^{L_{1}}C\cdot\left(1+\left|\xi\right|\right)^{-L_{1}}\cdot\min\left\{ 1,\left|\xi\right|^{L_{2}}\right\} 
\]
for all $\left|\xi\right|\leq1$ and all $\alpha\in\N_{0}^{\dimension}$
with $\left|\alpha\right|\leq N$.

Furthermore, since $\psi\in C_{c}^{\ell}\left(\R^{\dimension}\right)$
and thus also $\psi^{\left(\alpha\right)}:=\left[x\mapsto\left(-2\pi ix\right)^{\alpha}\cdot\psi\left(x\right)\right]\in C_{c}^{\ell}\left(\R^{\dimension}\right)$,
Lemma \ref{lem:PointwiseFourierDecayEstimate} and elementary properties
of the Fourier transform imply
\[
\left|\partial^{\alpha}\widehat{\psi}\left(\xi\right)\right|=\left|\left(\Fourier^{-1}\smash{\psi^{\left(\alpha\right)}}\right)\left(-\xi\right)\right|\lesssim\left(1+\left|\xi\right|\right)^{-\ell}\leq\left(1+\left|\xi\right|\right)^{-L_{1}}=\left(1+\left|\xi\right|\right)^{-L_{1}}\cdot\min\left\{ 1,\left|\xi\right|^{L_{2}}\right\} \qquad\text{ for }\left|\xi\right|\geq1
\]
for arbitrary $\alpha\in\N_{0}^{\dimension}$. Hence, $\partial^{\alpha}\widehat{\psi}$
satisfies the desired decay properties.

Finally, note that eq.\@ (\ref{eq:BesovExplicitConstructionDerivative})
also yields $\left|\smash{\widehat{\psi}}\left(\xi\right)\right|=\left(4\pi^{2}\right)^{L_{3}+N+1}\cdot\left|\xi\right|^{2\left(L_{3}+N+1\right)}\cdot\left|\smash{\widehat{\psi_{1}}}\left(\xi\right)\right|\geq c_{1}\cdot\left(\pi^{2}/4\right)^{L_{3}+N+1}\geq c_{1}$
for all $\xi\in\overline{B_{4}}\left(0\right)\setminus B_{1/4}\left(0\right)$.
We have thus constructed $\psi$ as desired. Similarly, but easier,
one can construct $\varphi\in C_{c}^{k}\left(\R^{\dimension}\right)$
satisfying $\widehat{\varphi}\left(\xi\right)\neq0$ for $\xi\in\overline{B_{2}}\left(0\right)$
and $\left|\partial^{\alpha}\widehat{\varphi}\left(\xi\right)\right|\lesssim\left(1+\left|\xi\right|\right)^{-L}$
for all $\alpha\in\N_{0}^{\dimension}$ with $\left|\alpha\right|\leq N$.
It is then straightforward to check that $\varphi,\psi$ satisfy all
assumptions of Propositions \ref{prop:BesovBanachFrames} and \ref{prop:BesovAtomicDecomposition}
(for proper choices of $N,L,L_{1},L_{2}$). Hence, our general theory
indeed yields \emph{compactly supported} wavelet systems that form
atomic decompositions and Banach frames for inhomogeneous Besov spaces.
\item We observe that the assumptions of Propositions \ref{prop:BesovBanachFrames}
and \ref{prop:BesovAtomicDecomposition} are \emph{structurally} very
similar, but the precise values of $L,L_{1},L_{2}$ differ greatly.
Indeed, in order to get a \emph{Banach frame} for $\mathcal{B}_{s}^{p,q}\left(\R^{\dimension}\right)$
using Proposition \ref{prop:BesovBanachFrames}, the mother wavelet
$\psi$ has to have at least $L_{2}>s$ vanishing moments, which increases
with the \emph{smoothness parameter} $s\in\R$. In contrast, in order
to obtain an \emph{atomic decomposition} of $\mathcal{B}_{s}^{p,q}\left(\R^{\dimension}\right)$
using Proposition \ref{prop:BesovAtomicDecomposition}, the mother
wavelet $\psi$ only has to have $L_{2}>\vartheta_{0}\dimension-s$
vanishing moments, where $\vartheta_{0}=\left(p^{-1}-1\right)_{+}$.
In particular, once the smoothness parameter $s$ satisfies $s>\dimension\left(p^{-1}-1\right)_{+}$,
one can choose $L_{2}=0$, so that it is possible for $\psi$ to have
\emph{no vanishing moments at all}, i.e., $\widehat{\psi}\left(0\right)\neq0$
is allowed.

In this case, one can even choose $\varphi=\psi$ to show that the
system $\left(2^{j\frac{\dimension}{2}}\cdot\varphi\left(2^{j}\bullet-\delta k\right)\right)_{j\in\N_{0},k\in\Z^{\dimension}}$
yields an atomic decompositions of $\mathcal{B}_{s}^{p,q}\left(\R^{\dimension}\right)$.
A peculiar property of this system is that it \emph{does not even
form a frame for $L^{2}\left(\R^{\dimension}\right)$}, due to the
missing vanishing moments.
\item Wavelet characterizations of inhomogeneous Besov spaces have already
been considered by many other authors: In \cite[equations (10.1) and (10.2)]{MeyerWaveletsAndOperators},
as well as in \cite[Theorem 3.5(i)]{TriebelTheoryOfFunctionSpaces3},
it is shown that certain wavelet \emph{orthonormal bases} yield atomic
decompositions and Banach frames for the Besov spaces $\mathcal{B}_{s}^{p,q}\left(\R^{\dimension}\right)$.
We remark that of the two mentioned books, only Triebel's book \cite{TriebelTheoryOfFunctionSpaces3}
covers the whole range $p,q\in\left(0,\infty\right]$, while Meyer\cite{MeyerWaveletsAndOperators}
only considers the case $p,q\in\left[1,\infty\right]$.

As explained in \cite[Theorem 1.61(ii)]{TriebelTheoryOfFunctionSpaces3},
the wavelet bases considered by Triebel in \cite[Theorem 3.5]{TriebelTheoryOfFunctionSpaces3}
are \emph{compactly supported} and are $C^{k}$, with $k$ vanishing
moments, where it is assumed that
\[
k>\max\left\{ s,\frac{2\dimension}{p}+\frac{\dimension}{2}-s\right\} .
\]
Hence, Triebel needs a large amount of vanishing moments if $s$ is
large, but also if $-s$ is large. As observed in the previous point,
this is not needed for the theory developed in this paper, at least
if one only wants to have \emph{either} Banach frames or atomic decompositions.
But since Triebel uses wavelet orthonormal bases, he obtains atomic
decompositions and Banach frames \emph{simultaneously}, which explains
the dependence of $k$ on $s$ observed above.

In addition to orthonormal bases, Triebel also considers wavelet \emph{frames},
cf.\@ \cite[Sections 1.8 and 3.2]{TriebelTheoryOfFunctionSpaces3}.
But for these, Triebel restricts to the case $p=q$. Then, for $s>\sigma_{p}=\dimension\left(p^{-1}-1\right)_{+}$,
he derives atomic decomposition results using certain compactly supported
wavelet frames (cf.\@ \cite[Theorem 1.69]{TriebelTheoryOfFunctionSpaces3}).
As seen above, this is the range in which Proposition \ref{prop:BesovAtomicDecomposition}
does not need any vanishing moments. Additional atomic decomposition
results are obtained in \cite[Theorem 1.71]{TriebelTheoryOfFunctionSpaces3},
but these use \emph{bandlimited} wavelets and require $p>1$ as well
as $s<0$.

In a different approach, Rauhut and Ullrich showed \cite{GeneralizedCoorbit2}
(based upon previous work by Ullrich\cite{UllrichContinuousCharacterizationsOfBesovTriebelLizorkin})
that the inhomogeneous Besov spaces $\mathcal{B}_{s}^{p,q}\left(\R^{\dimension}\right)$
can be obtained as certain \emph{generalized coorbit spaces}. Using
the theory of these spaces (cf.\@ \cite{GeneralizedCoorbit1,GeneralizedCoorbit2}),
they then again show that suitable wavelet \emph{orthonormal bases}
yield Banach frames and atomic decompositions for the spaces $\mathcal{B}_{s}^{p,q}\left(\R^{\dimension}\right)$,
cf.\@ \cite[Theorem 5.8 and Remark 5.9]{GeneralizedCoorbit2}. Their
assumptions on the scaling function $\varphi$ and the mother wavelet
$\psi$ are very similar to the ones imposed in this paper: $\psi$
needs to have a suitable number of vanishing moments, and $\varphi,\psi$
are required to have a suitable decay in space, as well as in Fourier
domain. Furthermore, the decay in Fourier domain also needs to hold
for certain derivatives of $\widehat{\varphi},\widehat{\psi}$, cf.\@
\cite[Definition 1.1]{GeneralizedCoorbit2}. We remark, however, that
in \cite[Theorem 5.8]{GeneralizedCoorbit2}, only the range $p,q\in\left[1,\infty\right]$
is considered. 

Finally, Frazier and Jawerth\cite{FrazierJawerthDecompositionOfBesovSpaces,FrazierJawerthDiscreteTransform,FrazierJawerthThePhiTransform}
also obtained atomic decompositions for Besov spaces, cf.\@ \cite[Theorem 7.1]{FrazierJawerthDecompositionOfBesovSpaces}.
In contrast to our approach, Frazier and Jawerth use a sampling density
which is fixed \emph{a priori}. This, however, requires the mother
wavelet $\psi$ to be \emph{bandlimited} to $\left\{ \xi\in\R^{\dimension}\with\left|\xi\right|\leq\pi\right\} $
(cf.\@ \cite[between eq. (1.8) and eq. (1.9)]{FrazierJawerthDecompositionOfBesovSpaces});
in particular, $\psi$ can \emph{not} be compactly supported. We remark
that Frazier and Jawerth assume $\psi$ to have $N$ vanishing moments
with $N\geq\max\left\{ -1,\,\dimension\left(p^{-1}-1\right)_{+}-s\right\} $.
This is very similar to the vanishing moment condition which we impose
in Proposition \ref{prop:BesovAtomicDecomposition}, cf.\@ the preceding
point. Finally, we mention the so-called \textbf{\emph{generalized}}\textbf{
$\varphi$-transform} of Frazier and Jawerth (cf.\@ \cite[Section 4]{FrazierJawerthDiscreteTransform})
which yields results that are very similar to Propositions \ref{prop:BesovBanachFrames}
and \ref{prop:BesovAtomicDecomposition}, but for the case of (homogeneous)
Triebel-Lizorkin spaces instead of inhomogeneous Besov spaces, cf.\@
\cite[Theorem 4.5]{FrazierJawerthThePhiTransform} and \cite[Corollaries 4.5 and 4.3]{FrazierJawerthDiscreteTransform}.
For the case of \emph{inhomogeneous} Triebel-Lizorkin spaces, see
\cite[Section 12, page 132]{FrazierJawerthDiscreteTransform}.

In summary, we have seen that the description of (inhomogeneous) Besov
spaces through wavelet systems—in particular through wavelet orthonormal
bases—was very well developed prior to this paper. Nevertheless, it
seems that in the case of compactly supported wavelet \emph{frames}
(as opposed to orthonormal bases), our results slightly improve the
state of the art: In \cite{TriebelTheoryOfFunctionSpaces3}, comparable
results are only derived for $p=q$ and $s>\dimension\left(p^{-1}-1\right)_{+}$
and in \cite{FrazierJawerthDecompositionOfBesovSpaces}, only bandlimited
wavelet systems are considered. Finally, in \cite{FrazierJawerthDiscreteTransform},
the authors allow compactly supported wavelet frames, but consider
Triebel-Lizorkin spaces instead of Besov spaces.

We close our comparison with the literature by comparing the advantages
and disadvantages of wavelet orthonormal bases compared to more general
wavelet systems. As noted in \cite[Example 5.6(a)]{GroechenigDescribingFunctions},
``\emph{both types of description are useful {[}...{]}: The orthogonal
bases, when a concise characterization of a function without redundancy
is important, but the form of the basic wavelet $g$ is not essential;
the non-orthogonal expansions and frames, when the basic function
$g$ is given by the problem and flexibility is required.}'' Indeed,
if one is willing to sample sufficiently densely, Propositions \ref{prop:BesovBanachFrames}
and \ref{prop:BesovAtomicDecomposition} allow a \emph{very wide variety}
of scaling functions $\varphi$ and mother wavelets $\psi$ to be
used. In contrast, to obtain an orthonormal wavelet basis, $\varphi$
and $\psi$ need to be selected \emph{very carefully}. However, using
such an orthonormal basis has several advantages\cite{TriebelTheoryOfFunctionSpaces3}
that frames lack:
\begin{itemize}
\item the sampling density is known and fixed a priori,
\item the synthesis coefficients are uniquely determined and equal to the
analysis coefficients,
\item the analysis map yields an isomorphism of $\mathcal{B}_{s}^{p,q}\left(\R^{\dimension}\right)$
\textbf{onto} the associated sequence space $b_{s}^{p,q}$.
\end{itemize}
\item Finally, we remark that we discussed inhomogeneous Besov spaces in
the general framework presented here mainly to indicate that—\emph{and
how}—the framework can be applied in concrete cases. More novel and
interesting applications of the general theory, in particular to shearlets,
will be discussed in the companion paper \cite{StructuredBanachFrames2}.\qedhere
\end{itemize}
\end{rem}

\appendix

\section{Lemmas needed to get explicit constants}
\begin{lem}
\label{lem:SmoothRampEstimates}For each $N\in\N$, there is a polynomial
$p_{N}\in\R\left[X\right]$ satisfying $0\leq p_{N}\left(x\right)\leq1$
for $x\in\left[0,1\right]$ and 
\[
p_{n}\left(0\right)=0,\quad p_{N}\left(1\right)=1,\quad\text{ as well as }\quad p_{N}^{\left(\ell\right)}\left(0\right)=0=p_{N}^{\left(\ell\right)}\left(1\right)\qquad\forall\:\ell\in\underline{N}.
\]
Furthermore, $p_{N}$ satisfies $\left\Vert \smash{p_{N}^{\left(\ell\right)}}\right\Vert _{\sup,\left[0,1\right]}\leq24^{N+1}\cdot\left(N+1\right)!$
for all $\ell\in\left\{ 0,\dots,N\right\} $.
\end{lem}
\begin{proof}
First, recall the well-known identity $\frac{\d^{\ell}}{\d x^{\ell}}x^{N}=\frac{N!}{\left(N-\ell\right)!}\cdot x^{N-\ell}=\binom{N}{\ell}\cdot\ell!\cdot x^{N-\ell}$
for $\ell\in\left\{ 0,\dots,N\right\} $. Now, define
\[
q_{N}\left(x\right):=x^{N}\cdot\left(1-x\right)^{N}=x^{N}\cdot\sum_{m=0}^{N}\binom{N}{m}\left(-x\right)^{m}=\sum_{m=0}^{N}\left[\binom{N}{m}\cdot\left(-1\right)^{m}\cdot x^{N+m}\right]
\]
and note for $x\in\left(0,1\right]$ and $\ell\in\left\{ 0,\dots,N\right\} $
that
\begin{align*}
\left|q_{N}^{\left(\ell\right)}\left(x\right)\right| & \leq\sum_{m=0}^{N}\binom{N}{m}\binom{N+m}{\ell}\cdot\ell!\cdot x^{N+m-\ell}\\
\left({\scriptstyle \text{since }\binom{a}{b}\leq2^{a}}\right) & \leq\ell!\cdot x^{-\ell}\cdot\sum_{m=0}^{N}\left(\binom{N}{m}\cdot2^{N+m}\cdot x^{N+m}\right)\\
 & =\ell!\cdot x^{N-\ell}\cdot2^{N}\cdot\sum_{m=0}^{N}\left(\binom{N}{m}\cdot\left(2x\right)^{m}\right)\\
 & =\ell!\cdot x^{N-\ell}\cdot2^{N}\cdot\left(1+2x\right)^{N}\\
\left({\scriptstyle \text{since }0<x\leq1\text{ and }N-\ell\geq0}\right) & \leq\ell!\cdot2^{N}\cdot3^{N}\leq6^{N}\cdot N!.
\end{align*}
By continuity, this also holds for $x=0$.

Furthermore, since we have $\frac{\d^{\ell}}{\d x^{\ell}}\bigg|_{x=0}x^{N+m}=0$
for all $\ell\leq N-1$ and all $m\geq0$, we see $q_{N}^{\left(\ell\right)}\left(0\right)=0$
for all $\ell\in\left\{ 0,\dots,N-1\right\} $. Likewise, note that
\begin{align*}
q_{N}\left(x\right) & =\left(-1\right)^{N}\cdot\left(x-1\right)^{N}\cdot\left(1+\left(x-1\right)\right)^{N}\\
 & =\left(-1\right)^{N}\cdot\sum_{m=0}^{N}\binom{N}{m}\cdot\left(x-1\right)^{N+m},
\end{align*}
which implies $q_{N}^{\left(\ell\right)}\left(1\right)=0$ for all
$\ell\in\left\{ 0,\dots,N-1\right\} $, since $\frac{\dimension^{\ell}}{\dimension x^{\ell}}\bigg|_{x=1}\left(x-1\right)^{N+m}=0$
for all $\ell\in\left\{ 0,\dots,N-1\right\} $ and $m\geq0$.

Next, note for $x\in\left[\frac{1}{2}\left(1-\frac{1}{2N}\right),\frac{1}{2}\left(1+\frac{1}{2N}\right)\right]\subset\left[0,1\right]$
that 
\[
x^{N}\geq2^{-N}\cdot\left(1-\frac{1}{2N}\right)^{N}=2^{-N}\cdot\sqrt{\left(1-\frac{1}{2N}\right)^{2N}}\geq2^{-N}\cdot\sqrt{\left(1-\frac{1}{2}\right)^{2}}=2^{-\left(N+1\right)},
\]
where we used the well-known fact that the sequence $\left[\left(1-\frac{1}{n}\right)^{n}\right]_{n\in\N}$
is nondecreasing\footnote{One way to see this is to note $\frac{\d}{\d x}\left(1-\frac{1}{x}\right)^{x}=\left(1-\frac{1}{x}\right)^{x}\cdot\left[\ln\left(1-\frac{1}{x}\right)+\frac{1}{x-1}\right]$
as well as $\frac{\d}{\d x}\left[\ln\left(1-\frac{1}{x}\right)+\frac{1}{x-1}\right]=\frac{1}{x-1}\left(\frac{1}{x}-\frac{1}{x-1}\right)<0$
for $x\in\left(1,\infty\right)$ and $\ln\left(1-\frac{1}{x}\right)+\frac{1}{x-1}\xrightarrow[x\to\infty]{}0$.
Together, these facts show $\frac{\d}{\d x}\left(1-\frac{1}{x}\right)^{x}>0$
on $\left(1,\infty\right)$.}. Likewise, we get
\[
\left(1-x\right)^{N}\geq\left[1-\frac{1}{2}\left(1+\frac{1}{2N}\right)\right]^{N}=\left[\frac{1}{2}\left(1-\frac{1}{2N}\right)\right]^{N}\geq2^{-\left(N+1\right)}
\]
and thus $q_{N}\left(x\right)\geq4^{-\left(N+1\right)}$, which yields
\[
C_{N}:=\int_{0}^{1}q_{N}\left(t\right)\d t\geq\int_{\frac{1}{2}\left(1-\frac{1}{2N}\right)}^{\frac{1}{2}\left(1+\frac{1}{2N}\right)}4^{-\left(N+1\right)}\d t=\frac{4^{-\left(N+1\right)}}{2N}=\frac{2^{-\left(2N+3\right)}}{N}.
\]

Now, we finally define for $x\in\left[0,1\right]$
\[
p_{N}\left(x\right):=\frac{1}{C_{N}}\cdot\int_{0}^{x}q_{N}\left(t\right)\d t=\frac{1}{C_{N}}\cdot\sum_{m=0}^{N}\left[\binom{N}{m}\cdot\frac{\left(-1\right)^{m}}{N+m+1}\cdot x^{N+m+1}\right]
\]
and note $p_{N}\left(0\right)=0$, as well as $p_{N}\left(1\right)=\frac{1}{C_{N}}\cdot C_{N}=1$,
as desired. Also, the fundamental theorem of calculus shows
\[
p_{N}^{\left(\ell\right)}\left(x\right)=\frac{1}{C_{N}}\cdot q_{N}^{\left(\ell-1\right)}\left(x\right)=0\qquad\forall\:\ell\in\underline{N}\text{ and }x\in\left\{ 0,1\right\} .
\]
Furthermore, since $q_{N}\geq0$, we see that $p_{N}$ is nondecreasing
and hence $0=p_{N}\left(0\right)\leq p_{N}\left(x\right)\leq p_{N}\left(1\right)=1$
for all $x\in\left[0,1\right]$. Finally, we get
\[
\left\Vert p_{N}^{\left(\ell\right)}\right\Vert _{\sup,\left[0,1\right]}=\frac{1}{C_{N}}\cdot\left\Vert q_{N}^{\left(\ell-1\right)}\right\Vert _{\sup,\left[0,1\right]}\leq2^{2N+3}\cdot N\cdot6^{N}\cdot N!\leq24^{N}\cdot8N\cdot N!\leq24^{N+1}\cdot\left(N+1\right)!
\]
for all $\ell\in\underline{N}$. For $\ell=0$, this estimate is trivially
satisfied since $\left\Vert p_{N}\right\Vert _{\sup,\left[0,1\right]}=1$.
\end{proof}
\begin{lem}
\label{lem:SmoothCutOffFunctionConstants}For all $\dimension,N\in\N$
and $R,s>0$ there is a function $\psi\in\TestFunctionSpace{\R^{\dimension}}$
satisfying

\begin{itemize}
\item $0\le\psi\leq1$,
\item $\supp\psi\subset\left(-\left(R+s\right),R+s\right)^{\dimension}$,
\item $\psi\equiv1$ on $\left[-R,R\right]^{\dimension}$,
\item $\left\Vert \frac{\partial^{\ell}}{\partial x_{i}^{\ell}}\psi\right\Vert _{\sup}\leq\max\left\{ 1,\left(\frac{3}{s}\right)^{\ell}\right\} \cdot24^{N+1}\cdot\left(N+1\right)!$
for all $i\in\underline{d}$ and all $\ell\in\left\{ 0,\dots,N\right\} $.\qedhere
\end{itemize}
\end{lem}
\begin{proof}
Choose $p_{N}$ as in Lemma \ref{lem:SmoothRampEstimates} and define
\[
\psi^{\left(0\right)}:\R\to\left[0,1\right],x\mapsto\begin{cases}
0, & \text{if }x\leq-\left(R+\frac{2}{3}s\right),\\
p_{N}\left(\frac{3}{s}\cdot\left(x+R+\frac{2}{3}s\right)\right), & \text{if }-\left(R+\frac{2}{3}s\right)\leq x\leq-\left(R+\frac{s}{3}\right),\\
1, & \text{if }-\left(R+\frac{s}{3}\right)\leq x\leq R+\frac{s}{3},\\
p_{N}\left(\frac{3}{s}\cdot\left[R+\frac{2}{3}s-x\right]\right), & \text{if }R+\frac{s}{3}\leq x\leq R+\frac{2}{3}s,\\
0, & \text{if }x\geq R+\frac{2}{3}s.
\end{cases}
\]
Since we have $0\leq p_{N}\leq1$ and $p_{N}\left(0\right)=0$, as
well as $p_{N}\left(1\right)=1$, it follows that $\psi^{\left(0\right)}$
is well-defined and continuous. Furthermore, it is well-known that
if $f,g:\R\to\R$ are differentiable with $f'\left(a\right)=g'\left(a\right)$
and $f\left(a\right)=g\left(a\right)$, then so is
\[
x\mapsto\begin{cases}
f\left(x\right), & \text{if }x\leq a,\\
g\left(x\right), & \text{if }x\geq a.
\end{cases}
\]
By applying this inductively to higher derivatives and since $p_{N}^{\left(\ell\right)}\left(1\right)=0=p_{N}^{\left(\ell\right)}\left(0\right)$
for all $\ell\in\left\{ 1,\dots,N\right\} $, we conclude $\psi^{\left(0\right)}\in C^{N}\left(\R\right)$,
where the derivatives are obtained by differentiation of the individual
``pieces'' defining $\psi^{\left(0\right)}$. In particular, we
get $\left\Vert \frac{\d^{\ell}}{\d x^{\ell}}\psi^{\left(0\right)}\right\Vert _{\sup}\leq\max\left\{ 1,\left(\frac{3}{s}\right)^{\ell}\right\} \cdot24^{N+1}\cdot\left(N+1\right)!$
for all $\ell\in\left\{ 0,\dots,N\right\} $, cf.\@ the estimate
for the derivatives of $p_{N}$.

Now, let $\theta\in C_{c}^{\infty}\left(\left(-1,1\right)\right)$
be a standard mollifier, i.e., $\theta\geq0$ with $\int_{\R}\theta\left(t\right)\d t=1$.
As usual, for $\varepsilon>0$, let $\theta_{\varepsilon}\left(x\right):=\frac{1}{\varepsilon}\cdot\theta\left(\frac{x}{\varepsilon}\right)$
and $\psi_{1}:=\theta_{s/3}\ast\psi^{\left(0\right)}$. Using standard
properties of convolution products, we see $\psi_{1}\in\TestFunctionSpace{\R}$
with $0\leq\psi_{1}\leq1$, as well as 
\[
\supp\psi_{1}\subset\left(-\frac{s}{3},\frac{s}{3}\right)+\supp\psi^{\left(0\right)}\subset\left(-\left(R+s\right),R+s\right)
\]
and with
\[
\left\Vert \frac{\d^{\ell}}{\d x^{\ell}}\psi_{1}\right\Vert _{\sup}=\left\Vert \theta_{s/3}\ast\left[\frac{\d^{\ell}}{\d x^{\ell}}\psi^{\left(0\right)}\right]\right\Vert _{\sup}\leq\left\Vert \frac{\d^{\ell}}{\d x^{\ell}}\psi^{\left(0\right)}\right\Vert _{\sup}\leq\max\left\{ 1,\left(\frac{3}{s}\right)^{\ell}\right\} \cdot24^{N+1}\cdot\left(N+1\right)!
\]
for all $\ell\in\left\{ 0,\dots,N\right\} $.

Finally, for $x\in\left[-R,R\right]$, we have
\begin{align*}
\psi_{1}\left(x\right) & =\int_{-\infty}^{\infty}\theta_{s/3}\left(y\right)\cdot\psi^{\left(0\right)}\left(x-y\right)\d y\\
 & =\int_{-s/3}^{s/3}\theta_{s/3}\left(y\right)\cdot\psi^{\left(0\right)}\left(x-y\right)\d y\\
\left({\scriptstyle \psi^{\left(0\right)}\left(x-y\right)=1\text{ since }x-y\in\left[-\left(R+\frac{s}{3}\right),R+\frac{s}{3}\right]}\right) & =\int_{-s/3}^{s/3}\theta_{s/3}\left(y\right)\d y\\
 & =\int_{-\infty}^{\infty}\theta_{s/3}\left(y\right)\d y=1,
\end{align*}
as desired.

\medskip{}

The preceding considerations establish the claim for $\dimension=1$.
In case of $\dimension>1$, set $\psi:=\psi_{1}\otimes\cdots\otimes\psi_{1}$
and note
\[
\frac{\partial^{\ell}}{\partial x_{i}^{\ell}}\psi=\bigotimes_{j=1}^{i-1}\psi_{1}\otimes\psi_{1}^{\left(\ell\right)}\otimes\bigotimes_{j=i+1}^{\dimension}\psi_{1},
\]
which yields the desired estimate for the derivative, since $0\leq\psi_{1}\leq1$.
\end{proof}
\begin{cor}
\label{cor:CutoffInverseFourierEstimate}The function $\psi$ from
Lemma \ref{lem:SmoothCutOffFunctionConstants} satisfies
\[
\left|\left(\partial^{\alpha}\left[\Fourier^{-1}\psi\right]\right)\left(x\right)\right|\leq2\pi\cdot2^{\dimension}\cdot\max\left\{ 1,\left(3/s\right)^{N}\right\} \cdot\left(48\dimension\right)^{N+1}\left(N+2\right)!\cdot\left(1+R+s\right)^{\left|\alpha\right|}\left(R+s\right)^{\dimension}\cdot\left(1+\left|x\right|\right)^{-N},
\]
for all $x\in\R^{\dimension}$ and all $\alpha\in\N_{0}^{d}$ with
$\left|\alpha\right|\leq1$.

In particular, we have for $N>\dimension$ that
\[
\left\Vert \nabla\left[\Fourier^{-1}\psi\right]\right\Vert _{L^{1}}\leq\frac{2\pi\cdot s_{\dimension}}{N-\dimension}\cdot\sqrt{\dimension}\cdot2^{\dimension}\cdot\max\left\{ 1,\left(3/s\right)^{N}\right\} \cdot\left(48\dimension\right)^{N+1}\left(N+2\right)!\cdot\left(1+R+s\right)\left(R+s\right)^{\dimension}.\qedhere
\]
\end{cor}
\begin{proof}
We first recall the elementary identity
\begin{align*}
\left(\partial^{\alpha}\left[\Fourier^{-1}\psi\right]\right)\left(x\right) & =\int_{\R^{\dimension}}\psi\left(\xi\right)\cdot\partial_{x}^{\alpha}e^{2\pi i\left\langle x,\xi\right\rangle }\d\xi\\
 & =\int_{\R^{\dimension}}\psi\left(\xi\right)\cdot\left(2\pi i\xi\right)^{\alpha}\cdot e^{2\pi i\left\langle x,\xi\right\rangle }\d\xi\\
 & =\left(\Fourier^{-1}\left[\xi\mapsto\left(2\pi i\xi\right)^{\alpha}\cdot\psi\left(\xi\right)\right]\right)\left(x\right).
\end{align*}
Hence, we let $g:\R^{\dimension}\to\Compl,\xi\mapsto\left(2\pi i\xi\right)^{\alpha}\cdot\psi\left(\xi\right)$,
note $g\in\TestFunctionSpace{\R^{\dimension}}$ and recall from Lemma
\ref{lem:PointwiseFourierDecayEstimate}, equation (\ref{eq:PointwiseFourierDecayEstimate})
that
\begin{equation}
\left|\left(\Fourier^{-1}g\right)\left(x\right)\right|\leq\left(1+\left|x\right|\right)^{-N}\cdot\left(1+\dimension\right)^{N}\cdot\left(\left|\left(\Fourier^{-1}g\right)\left(x\right)\right|+\sum_{m=1}^{\dimension}\left|\left[\Fourier^{-1}\left(\partial_{m}^{N}g\right)\right]\left(x\right)\right|\right)\label{eq:CutoffInverseFourierProof}
\end{equation}
for all $x\in\R^{\dimension}$. Thus, it remains to estimate the right-hand
side.

But for the first term, we simply have because of $\supp g\subset\supp\psi\subset\left(-\left(R+s\right),R+s\right)^{\dimension}$
and $0\leq\psi\leq1$, which entails $\left|g\left(\xi\right)\right|\leq\left[2\pi\left(R+s\right)\right]^{\left|\alpha\right|}$
for all $\xi\in\R^{\dimension}$, that
\[
\left|\left(\Fourier^{-1}g\right)\left(x\right)\right|\leq\left\Vert g\right\Vert _{L^{1}}\leq\left[2\pi\left(R+s\right)\right]^{\left|\alpha\right|}\cdot\left[2\left(R+s\right)\right]^{\dimension}\leq2\pi\cdot2^{\dimension}\cdot\left(1+R+s\right)^{\left|\alpha\right|}\left(R+s\right)^{\dimension}.
\]
For the second term, we have to work harder: In case of $\alpha=0$,
we simply have $g=\psi$ and hence—as above—that
\begin{align*}
\left|\left[\Fourier^{-1}\left(\partial_{m}^{N}g\right)\right]\left(x\right)\right| & \leq\left\Vert \partial_{m}^{N}g\right\Vert _{L^{1}}\leq\left[2\left(R+s\right)\right]^{\dimension}\cdot\left\Vert \partial_{m}^{N}g\right\Vert _{L^{\infty}}\\
 & \leq\max\left\{ 1,\left(\frac{3}{s}\right)^{N}\right\} \cdot24^{N+1}\cdot\left(N+1\right)!\cdot\left[2\left(R+s\right)\right]^{\dimension}\\
 & \leq2\pi\cdot2^{\dimension}\cdot\max\left\{ 1,\left(\frac{3}{s}\right)^{N}\right\} \cdot24^{N+1}\cdot\left(N+2\right)!\cdot\left(R+s\right)^{\dimension}\left(1+R+s\right)^{\left|\alpha\right|},
\end{align*}
cf.\@ Lemma \ref{lem:SmoothCutOffFunctionConstants} for the estimate
regarding $\left\Vert \partial_{m}^{N}g\right\Vert _{L^{\infty}}=\left\Vert \partial_{m}^{N}\psi\right\Vert _{L^{\infty}}$.

It remains to consider the case $\left|\alpha\right|=1$, i.e., $\alpha=e_{j}$
for some $j\in\underline{d}$. In this case, we have $g\left(\xi\right)=2\pi i\cdot\xi_{j}\cdot\psi\left(\xi\right)$,
so that Leibniz's rule yields
\begin{align*}
\left|\left(\partial_{m}^{N}g\right)\left(\xi\right)\right| & =2\pi\cdot\left|\sum_{\ell=0}^{N}\binom{N}{\ell}\cdot\left[\partial_{m}^{\ell}\xi_{j}\right]\cdot\left(\partial_{m}^{N-\ell}\psi\right)\left(\xi\right)\right|\\
\left({\scriptstyle \text{since }\left|\xi_{j}\right|\leq R+s\text{ on }\supp\psi\text{ and }\partial_{m}^{\ell}\xi_{j}=\delta_{m,j}\cdot\delta_{\ell,1}\text{ for }\ell\geq1}\right) & \leq2\pi\cdot\left[\left(R+s\right)\cdot\left|\left(\partial_{m}^{N}\psi\right)\left(\xi\right)\right|+N\cdot\left|\left(\partial_{m}^{N-1}\psi\right)\left(\xi\right)\right|\right]\\
\left({\scriptstyle \text{cf. Lemma }\ref{lem:SmoothCutOffFunctionConstants}}\right) & \leq2\pi N\cdot\left(1+R+s\right)\cdot\max\left\{ 1,\left(\frac{3}{s}\right)^{N}\right\} \cdot24^{N+1}\cdot\left(N+1\right)!\\
\left({\scriptstyle \text{since }\left|\alpha\right|=1}\right) & \leq2\pi\cdot\left(1+R+s\right)^{\left|\alpha\right|}\cdot\max\left\{ 1,\left(\frac{3}{s}\right)^{N}\right\} \cdot24^{N+1}\cdot\left(N+2\right)!\,.
\end{align*}
Here, we used that $\max\left\{ 1,\left(3/s\right)^{\ell}\right\} $
is nondecreasing with respect to $\ell\in\left\{ 0,\dots,N\right\} $.
In fact, for $3/s\leq1$, we have $\max\left\{ 1,\left(3/s\right)^{\ell}\right\} =1$
for all $\ell\in\left\{ 0,\dots,N\right\} $ and for $3/s>1$, we
have $\max\left\{ 1,\left(3/s\right)^{\ell}\right\} =\left(3/s\right)^{\ell}$,
which is increasing with respect to $\ell$. All in all, we arrive
at
\begin{align*}
\left|\left[\Fourier^{-1}\left(\partial_{m}^{N}g\right)\right]\left(x\right)\right| & \leq\left\Vert \partial_{m}^{N}g\right\Vert _{L^{1}}\leq\left[2\left(R+s\right)\right]^{\dimension}\cdot\left\Vert \partial_{m}^{N}g\right\Vert _{L^{\infty}}\\
 & \leq2\pi\cdot\max\left\{ 1,\left(\frac{3}{s}\right)^{N}\right\} \cdot24^{N+1}\cdot\left(N+2\right)!\cdot\left[2\left(R+s\right)\right]^{\dimension}\cdot\left(1+R+s\right)^{\left|\alpha\right|}\\
 & \leq2\pi\cdot2^{\dimension}\cdot\max\left\{ 1,\left(\frac{3}{s}\right)^{N}\right\} \cdot24^{N+1}\cdot\left(N+2\right)!\cdot\left(R+s\right)^{\dimension}\cdot\left(1+R+s\right)^{\left|\alpha\right|}.
\end{align*}

Recalling equation (\ref{eq:CutoffInverseFourierProof}), we conclude
\begin{align*}
 & \left(1+\left|x\right|\right)^{N}\cdot\left|\left(\Fourier^{-1}g\right)\left(x\right)\right|\\
 & \leq\left(1+\dimension\right)^{N}\!\left(2\pi\cdot2^{\dimension}\left(1+R+s\right)^{\left|\alpha\right|}\left(R+s\right)^{\dimension}+\dimension\cdot2\pi\cdot2^{\dimension}\max\left\{ 1,\left(\frac{3}{s}\right)^{N}\right\} \cdot24^{N+1}\left(N+2\right)!\cdot\left(R+s\right)^{\dimension}\left(1+R+s\right)^{\left|\alpha\right|}\!\right)\\
 & \leq2\pi\cdot2^{\dimension}\left(1+\dimension\right)^{N+1}\cdot\max\left\{ 1,\left(\frac{3}{s}\right)^{N}\right\} \cdot24^{N+1}\left(N+2\right)!\cdot\left(1+R+s\right)^{\left|\alpha\right|}\left(R+s\right)^{\dimension}\\
 & \leq2\pi\cdot2^{\dimension}\cdot\max\left\{ 1,\left(\frac{3}{s}\right)^{N}\right\} \cdot\left(48\dimension\right)^{N+1}\left(N+2\right)!\cdot\left(1+R+s\right)^{\left|\alpha\right|}\left(R+s\right)^{\dimension},
\end{align*}
as claimed.

\medskip{}

For the additional claim, recall from equation (\ref{eq:StandardDecayLpEstimate})
for $N>\dimension$ that
\begin{align*}
\left\Vert \nabla\left[\Fourier^{-1}\psi\right]\right\Vert _{L^{1}} & \leq2\pi\sqrt{\dimension}\cdot2^{\dimension}\cdot\max\left\{ 1,\left(\frac{3}{s}\right)^{N}\right\} \cdot\left(48\dimension\right)^{N+1}\left(N+2\right)!\cdot\left(1+R+s\right)\left(R+s\right)^{\dimension}\cdot\left\Vert \left(1+\left|\cdot\right|\right)^{-N}\right\Vert _{L^{1}}\\
 & \leq\frac{2\pi\cdot s_{\dimension}}{N-\dimension}\cdot\sqrt{\dimension}\cdot2^{\dimension}\cdot\max\left\{ 1,\left(\frac{3}{s}\right)^{N}\right\} \cdot\left(48\dimension\right)^{N+1}\left(N+2\right)!\cdot\left(1+R+s\right)\left(R+s\right)^{\dimension},
\end{align*}
where the first step is justified by a combination of our previous
estimates with the Cauchy-Schwarz inequality.
\end{proof}

\section{Vanishing of a function implies vanishing of derivatives}

In this section, we show that if a sufficiently smooth function $f:U\subset\R^{\dimension}\to\R$
satisfies $f\left(x\right)\in o\left(\smash{\left|x-a\right|^{N}}\right)$
as $x\to a$, then the partial derivatives of $f$ also vanish to
a suitable order at $a$, i.e., $\partial^{\alpha}f\left(x\right)\in\vphantom{\left|x-a\right|^{N-\left|\alpha\right|}}o\left(\smash{\left|x-a\right|^{N-\left|\alpha\right|}}\right)$
as $x\to a$, for $\left|\alpha\right|\leq N$. Our starting point
is the following consequence of Taylor's theorem:
\begin{lem}
\label{lem:HighOrderVanishingYieldZeroDerivatives}Let $a\in\R^{\dimension}$,
$r>0$ and $N\in\N_{0}$. Assume that $f\in C^{N}\left(B_{r}\left(a\right);\R\right)$
satisfies $f\left(x\right)\in o\left(\smash{\left|x-a\right|^{N}}\right)$
as $x\to a$, i.e.,
\[
\frac{f\left(x\right)}{\left|x-a\right|^{N}}\xrightarrow[x\to a]{}0.
\]
Then $\partial^{\alpha}f\left(a\right)=0$ for all $\alpha\in\N_{0}^{\dimension}$
with $\left|\alpha\right|\leq N$.
\end{lem}
\begin{proof}
Let $f$ as in the statement of the lemma. We will show by induction
on $\ell=\left|\alpha\right|\in\left\{ 0,\dots,N\right\} $ that $\partial^{\alpha}f\left(a\right)=0$.
By translating everything, we can clearly assume $a=0$.

For $\left|\alpha\right|=0$, we simply note by continuity of $f$
at $a=0$ that 
\[
f\left(0\right)=\lim_{x\to0}f\left(x\right)=\lim_{x\to0}\frac{f\left(x\right)}{\left|x\right|^{N}}\cdot\lim_{x\to0}\left|x\right|^{N}=0,
\]
since $\lim_{x\to0}\left|x\right|^{N}\in\left\{ 0,1\right\} $ because
of $N\in\N_{0}$.

Now, assume $\partial^{\alpha}f\left(0\right)=0$ for all $\left|\alpha\right|<\ell$
for some $\ell\in\left\{ 1,\dots,N\right\} $. In view of Taylor's
theorem (cf.\@ \cite[Theorem 5.11]{AmannEscherAnalysis2} for the
precise version used here), we get
\[
f\left(x\right)=\sum_{\left|\alpha\right|\leq N}\frac{\partial^{\alpha}f\left(0\right)}{\alpha!}x^{\alpha}+R_{N}\left(x\right)\qquad\forall x\in B_{r}\left(0\right),
\]
where $R_{N}:B_{r}\left(0\right)\to\R$ satisfies $R_{N}\left(x\right)\in o\left(\smash{\left|x\right|^{N}}\right)$,
i.e., $R_{N}\left(x\right)/\left|x\right|^{N}\xrightarrow[x\to0]{}0$.
By rearranging, and since $\partial^{\alpha}f\left(0\right)=0$ for
all $\left|\alpha\right|<\ell$, we get 
\[
p\left(x\right):=\sum_{\left|\alpha\right|=\ell}\frac{\partial^{\alpha}f\left(0\right)}{\alpha!}x^{\alpha}=f\left(x\right)-\sum_{\ell+1\leq\left|\alpha\right|\leq N}\frac{\partial^{\alpha}f\left(0\right)}{\alpha!}x^{\alpha}-R_{N}\left(x\right)=:g\left(x\right)\qquad\forall x\in B_{r}\left(0\right).
\]
But we have $R_{N}\left(x\right)\in o\left(\smash{\left|x\right|^{N}}\right)\subset o\left(\smash{\left|x\right|^{\ell}}\right)$
and likewise $f\left(x\right)\in o\left(\smash{\left|x\right|^{N}}\right)\subset o\left(\smash{\left|x\right|^{\ell}}\right)$
as $x\to0$. Also, $\left|x^{\alpha}\right|\leq\left|x\right|^{\left|\alpha\right|}\leq\left|x\right|^{\ell+1}$
for $\left|x\right|<1$ and $\ell+1\leq\left|\alpha\right|\leq N$,
so that $\frac{\left|x^{\alpha}\right|}{\left|x\right|^{\ell}}\leq\left|x\right|\xrightarrow[x\to0]{}0$.
All in all, we thus see $g\left(x\right)\in o\left(\smash{\left|x\right|^{\ell}}\right)$
as $x\to0$ and hence also $p\left(x\right)\in o\left(\smash{\left|x\right|^{\ell}}\right)$.

Next, if we can show $p\equiv0$, it follows from standard properties
of polynomials that $\frac{\partial^{\alpha}f\left(0\right)}{\alpha!}=0$
for all $\left|\alpha\right|=\ell$ and hence $\partial^{\alpha}f\left(0\right)=0$
for all $\left|\alpha\right|=\ell$. One possibility of proving this
elementary fact is to note that
\[
\partial^{\beta}x^{\alpha}=\begin{cases}
0, & \text{if }\alpha\neq\beta,\\
c_{\alpha}>0, & \text{if }\alpha=\beta,
\end{cases}
\]
for $\alpha,\beta\in\N_{0}^{\dimension}$ with $\left|\alpha\right|=\left|\beta\right|$.

Thus, all we need to show is that if $p\left(x\right)=\sum_{\left|\alpha\right|=\ell}c_{\alpha}x^{\alpha}$
satisfies $p\in o\left(\smash{\left|x\right|^{\ell}}\right)$ as $x\to0$
then $p\equiv0$. It is clear that $p\left(0\right)=0$, since $\ell\geq1$.
Now let $x\in\R^{\dimension}\setminus\left\{ 0\right\} $ be arbitrary.
Since $p$ is homogeneous of degree $\ell$, we have $p\left(rx\right)=r^{\ell}\cdot p\left(x\right)$
for all $r>0$. Using this and $p\in\vphantom{\left|x\right|^{\ell}}o\left(\smash{\left|x\right|^{\ell}}\right)$,
we get
\[
0=\lim_{r\downarrow0}\frac{p\left(rx\right)}{\left|rx\right|^{\ell}}=\lim_{r\downarrow0}\frac{p\left(x\right)}{\left|x\right|^{\ell}}=\frac{p\left(x\right)}{\left|x\right|^{\ell}}
\]
and hence $p\left(x\right)=0$ for all $x\in\R^{\dimension}$, as
desired.
\end{proof}
\begin{lem}
\label{lem:HighOrderVanishingYieldVanishingDerivatives}Let $U\subset\R^{\dimension}$
be open, let $N\in\N_{0}$ and $f\in C^{N}\left(U;\R\right)$. If
$f$ satisfies $f\left(x\right)/\left|x-a\right|^{N}\to0$ as $x\to a$
for some $a\in U$, then
\[
\frac{\partial^{\alpha}f\left(x\right)}{\left|x-a\right|^{N-\left|\alpha\right|}}\xrightarrow[x\to a]{}0\qquad\forall\alpha\in\N_{0}^{\dimension}\text{ with }\left|\alpha\right|\leq N.\qedhere
\]
\end{lem}
\begin{rem*}
The lemma remains true for complex-valued functions, since one can
simply apply it to the real- and imaginary parts separately.
\end{rem*}
\begin{proof}
Lemma \ref{lem:HighOrderVanishingYieldZeroDerivatives} shows $\partial^{\gamma}f\left(a\right)=0$
for all $\gamma\in\N_{0}^{\dimension}$ with $\left|\gamma\right|\leq N$.
Hence, for $\left|\alpha\right|=N$ we get by continuity of $\partial^{\alpha}f$
that
\[
\frac{\partial^{\alpha}f\left(x\right)}{\left|x-a\right|^{N-\left|\alpha\right|}}=\partial^{\alpha}f\left(x\right)\xrightarrow[x\to a]{}\partial^{\alpha}f\left(a\right)=0,
\]
as desired. Hence, we can assume $\left|\alpha\right|<N$ from now
on. This implies $g:=\partial^{\alpha}f\in C^{N-\left|\alpha\right|}\left(U;\R\right)$,
so that Taylor's theorem (see \cite[Theorem 5.11]{AmannEscherAnalysis2}
for the precise version used here) yields because of $\partial^{\beta}g\left(a\right)=\partial^{\alpha+\beta}f\left(a\right)=0$
for $\left|\beta\right|\leq N-\left|\alpha\right|$ that
\[
g\left(x\right)=\sum_{\left|\beta\right|\leq N-\left|\alpha\right|}\frac{\partial^{\beta}g\left(a\right)}{\beta!}\left(x-a\right)^{\beta}+o\left(\left|x-a\right|^{N-\left|\alpha\right|}\right)=o\left(\left|x-a\right|^{N-\left|\alpha\right|}\right)\qquad\text{ as }x\to a,
\]
as claimed.
\end{proof}

\section{Necessity of vanishing moment conditions for discrete cone-adapted
shearlet frames}

\label{sec:ShearletFrameVanishingMomentNecessity}
\begin{prop}
\label{prop:ShearletFramesNeedVanishingMoments}Let $\varphi,\gamma,\tilde{\gamma}\in L^{2}\left(\R^{2}\right)$
such that for some $\delta>0$, the \textbf{(discrete, cone-adapted)
shearlet system} (cf.\@ \cite[Definition 2.2]{CompactlySupportedShearletFrames})
with sampling density $\delta$,
\[
\mathcal{SH}\left(\varphi,\gamma,\tilde{\gamma};\delta\right)=\Phi\left(\varphi;\delta\right)\cup\Psi\left(\gamma;\delta\right)\cup\tilde{\Psi}\left(\tilde{\gamma};\delta\right)
\]
with
\begin{align*}
\Phi\left(\varphi;\delta\right) & =\left\{ \phi\left(\bullet-\delta m\right)\with m\in\Z^{2}\right\} ,\\
\Psi\left(\gamma;\delta\right) & =\left\{ \gamma_{j,k,m}=2^{\frac{3}{4}j}\cdot\gamma\left(S_{k}A_{2^{j}}\bullet-\delta m\right)\with\left(j,k\right)\in I\text{ and }m\in\Z^{2}\right\} \\
\tilde{\Psi}\left(\tilde{\gamma};\delta\right) & =\left\{ \tilde{\gamma}_{j,k,m}=2^{\frac{3}{4}j}\cdot\tilde{\gamma}\left(S_{k}^{T}\tilde{A_{2^{j}}}\bullet-\delta m\right)\with\left(j,k\right)\in I\text{ and }m\in\Z^{2}\right\} 
\end{align*}
and $I:=\left\{ \left(j,k\right)\in\N_{0}\times\Z\,:\,\left|k\right|\leq\left\lceil 2^{j/2}\right\rceil \right\} $,
as well as
\[
S_{k}=\left(\begin{matrix}1 & k\\
0 & 1
\end{matrix}\right),\quad A_{a}={\rm diag}\left(a,\sqrt{a}\right),\quad\text{ and }\quad\tilde{A}_{a}={\rm diag}\left(\sqrt{a},a\right)
\]
is a \textbf{Bessel system} in $L^{2}\left(\R^{2}\right)$. Then we
have
\[
\int_{\left\{ \xi\in\R^{2}\with\left|\xi_{1}\right|<1\text{ and }\left|\xi_{2}\right|<\frac{1}{8}\left|\xi_{1}\right|^{1/2}\right\} }\left|\xi_{1}\right|^{-2}\cdot\left|\widehat{\gamma}\left(\xi\right)\right|^{2}\d\xi<\infty.\qedhere
\]
\end{prop}
\begin{rem*}

\begin{itemize}[leftmargin=0.5cm]
\item Here, a system $\left(\theta_{i}\right)_{i\in I}$ in a Hilbert space
$\mathcal{H}$ is called a Bessel system if there is a constant $C>0$
satisfying $\sum_{i\in I}\left|\left\langle f,\,\theta_{i}\right\rangle _{\mathcal{H}}\right|^{2}\leq C\cdot\left\Vert f\right\Vert _{\mathcal{H}}^{2}$
for all $f\in\mathcal{H}$.
\item Likewise, one can show 
\[
\int_{\left\{ \xi\in\R^{2}\with\left|\xi_{2}\right|<1\text{ and }\left|\xi_{1}\right|<\frac{1}{8}\left|\xi_{2}\right|^{1/2}\right\} }\left|\xi_{2}\right|^{-2}\cdot\left|\left(\Fourier\smash{\tilde{\gamma}}\right)\left(\xi\right)\right|^{2}\d\xi<\infty.
\]
\item In particular, if $\widehat{\gamma}$ is continuous (e.g.\@ if $\gamma$
is compactly supported), then necessarily $\widehat{\gamma}\left(0\right)=0$,
since otherwise $\left|\widehat{\gamma}\left(\xi\right)\right|\geq c$
for $\left|\xi\right|<2\varepsilon$ with $c>0$ and $\varepsilon\in\left(0,1\right)$
suitable. But this yields
\begin{align*}
\int_{\left\{ \xi\in\R^{2}\with\left|\xi_{1}\right|<1\text{ and }\left|\xi_{2}\right|<\frac{1}{8}\left|\xi_{1}\right|^{1/2}\right\} }\left|\xi_{1}\right|^{-2}\cdot\left|\widehat{\gamma}\left(\xi\right)\right|^{2}\d\xi & \geq c^{2}\cdot\int_{\left\{ \xi\in\R^{2}\with\left|\xi_{1}\right|<\varepsilon^{2}\text{ and }\left|\xi_{2}\right|<\frac{1}{8}\left|\xi_{1}\right|^{1/2}\right\} }\left|\xi_{1}\right|^{-2}\d\xi\\
\left({\scriptstyle \text{Fubini's theorem}}\right) & \geq c^{2}\cdot\int_{0}^{\varepsilon^{2}}\left|\xi_{1}\right|^{-2}\cdot\frac{1}{4}\left|\xi_{1}\right|^{1/2}\,\d\xi_{1}\\
 & =\frac{c^{2}}{4}\cdot\int_{0}^{\varepsilon^{2}}\left|\xi_{1}\right|^{-\frac{3}{2}}\,\d\xi_{1}=\infty.\qedhere
\end{align*}
\end{itemize}
\end{rem*}
\begin{proof}
The following proof is heavily inspired by the proof of \cite[Theorem 3.3.1]{DaubechiesTenLecturesOnWavelets},
generalized from wavelets to shearlets and from homogeneous systems
to inhomogeneous systems.

In the following, we will consider the \textbf{shearlet group}
\[
H=\left\{ \varepsilon\left(\begin{matrix}a & b\\
0 & \sqrt{a}
\end{matrix}\right)\with a>0,\,b\in\R,\,\varepsilon\in\left\{ \pm1\right\} \right\} =\left\{ \left(\begin{matrix}a & b\\
0 & {\rm sgn}\left(a\right)\cdot\sqrt{\left|a\right|}
\end{matrix}\right)\with a\in\R\setminus\left\{ 0\right\} ,\,b\in\R\right\} ,
\]
which contains all of the matrices $S_{k}$ and $A_{a}$. We let
$\mu_{H}$ denote the \textbf{Haar measure} (cf.\@ \cite[Section 2.2]{FollandAHA})
on the locally compact topological group $H$. Based on $\mu_{H}$,
we define a new measure $\nu$ on the Borel $\sigma$-algebra of $\R^{2}\times H$
by
\[
\nu\left(A\right)=\int_{H}\int_{\R^{2}}\frac{\Indicator_{A}\left(x,h\right)}{\left|\det h\right|}\d x\d\mu_{H}\left(h\right)\quad\text{ for each Borel set }A\subset\R^{2}\times H.
\]

For $f,\psi\in L^{2}\left(\R^{2}\right)$, the \textbf{continuous
shearlet transform} $W_{\psi}f$ is given by
\[
W_{\psi}f:\R^{2}\times H\to\Compl,\left(x,h\right)\mapsto\left\langle f,\,\pi\left(x,h\right)\psi\right\rangle ,
\]
where $\pi\left(x,h\right)\psi:=L_{x}D_{h}\psi$, with $L_{x}f\left(y\right)=g\left(y-x\right)$
and $D_{h}f=\left|\det h\right|^{-1/2}\cdot f\circ h^{-1}$. It is
not hard to show that the inverse of the operator $\pi\left(x,h\right)=L_{x}D_{h}$
is given by $\left[\pi\left(x,h\right)\right]^{-1}=\pi\left(-h^{-1}x,h^{-1}\right)$
and furthermore that $\pi\left(x,h\right)\pi\left(y,g\right)=\pi\left(x+hy,\,hg\right)$
for arbitrary $x,y\in\R^{2}$ and $g,h\in H$.

Since we have $\left|W_{\psi}f\left(x,h\right)\right|\leq\left\Vert f\right\Vert _{L^{2}}\cdot\left\Vert \psi\right\Vert _{L^{2}}$
and since $\left\Vert \pi\left(x,h\right)\psi\right\Vert _{L^{2}}=\left\Vert \psi\right\Vert _{L^{2}}$
for all $\left(x,h\right)\in\R^{2}\times H$, it is not hard to see
for $F\in L^{1}\left(\nu;\Compl\right)$ that
\[
T_{F}f:=\int_{\R^{2}\times H}F\left(x,h\right)\cdot W_{\psi}f\left(x,h\right)\cdot\pi\left(x,h\right)\psi\d\nu\left(x,h\right)\in L^{2}\left(\smash{\R^{2}}\right)
\]
is well-defined with $\left\Vert T_{F}f\right\Vert _{L^{2}}\leq\left\Vert F\right\Vert _{L^{1}\left(\nu\right)}\cdot\left\Vert \psi\right\Vert _{L^{2}}^{2}\cdot\left\Vert f\right\Vert _{L^{2}}$.
Furthermore, in case of $F\geq0$, we have
\begin{align*}
\left\langle T_{F}f,\,f\right\rangle _{L^{2}} & =\int_{\R^{2}\times H}F\left(x,h\right)\cdot W_{\psi}f\left(x,h\right)\cdot\left\langle \pi\left(x,h\right)\psi,\,f\right\rangle _{L^{2}}\d\nu\left(x,h\right)\\
\left({\scriptstyle \text{since }\left\langle \pi\left(x,h\right)\psi,\,f\right\rangle _{L^{2}}=\overline{\left\langle f,\,\pi\left(x,h\right)\psi\right\rangle _{L^{2}}}=\overline{W_{\psi}f\left(x,h\right)}}\right) & =\int_{\R^{2}\times H}F\left(x,h\right)\cdot\left|W_{\psi}f\left(x,h\right)\right|^{2}\d\nu\left(x,h\right)\geq0,
\end{align*}
so that the operator $T_{F}:L^{2}\left(\R^{2}\right)\to L^{2}\left(\R^{2}\right)$
is bounded and nonnegative and in particular self-adjoint.

Finally, if $F\in L^{1}\left(\nu\right)$ is nonnegative and if $\left(u_{\ell}\right)_{\ell}$
is an arbitrary orthonormal basis of $L^{2}\left(\R^{2}\right)$,
then
\begin{align*}
\sum_{\ell}\left\langle T_{F}u_{\ell},\,u_{\ell}\right\rangle _{L^{2}} & =\int_{\R^{2}\times H}F\left(x,h\right)\cdot\sum_{\ell}\left|W_{\psi}u_{\ell}\left(x,h\right)\right|^{2}\d\nu\left(x,h\right)\\
\left({\scriptstyle \text{since }\sum_{\ell}\left|W_{\psi}u_{\ell}\left(x,h\right)\right|^{2}=\sum_{\ell}\left|\left\langle u_{\ell},\,\pi\left(x,h\right)\psi\right\rangle \right|^{2}=\left\Vert \pi\left(x,h\right)\psi\right\Vert _{L^{2}}^{2}=\left\Vert \psi\right\Vert _{L^{2}}^{2}}\right) & =\left\Vert \psi\right\Vert _{L^{2}}^{2}\cdot\left\Vert F\right\Vert _{L^{1}\left(\nu\right)}<\infty.
\end{align*}
Thus, $T_{F}$ is a \textbf{trace-class operator} (cf.\@ \cite[Appendix 2]{FollandAHA})
and in particular a compact operator.

Hence, if $F\in L^{1}\left(\nu\right)$ and $F\geq0$, then the spectral
theorem for compact self-adjoint operators yields an orthonormal basis
$\left(u_{\ell}\right)_{\ell\in\N}$ of $L^{2}\left(\R^{2}\right)$
satisfying $T_{F}=\sum_{\ell=1}^{\infty}c_{\ell}\cdot\left\langle \bullet,u_{\ell}\right\rangle \cdot u_{\ell}$
where $c_{\ell}\geq0$ and $\sum_{\ell=1}^{\infty}c_{\ell}=\left\Vert \psi\right\Vert _{L^{2}}^{2}\cdot\left\Vert F\right\Vert _{L^{1}\left(\nu\right)}<\infty$.

Now, if $\left(g_{i}\right)_{i\in I}$ is an arbitrary Bessel-sequence
in $L^{2}\left(\R^{2}\right)$, i.e., if $\sum_{i\in I}\left|\left\langle f,g_{i}\right\rangle \right|^{2}\leq C\cdot\left\Vert f\right\Vert _{L^{2}}^{2}$
for each $f\in L^{2}\left(\R^{2}\right)$, then
\begin{equation}
\begin{split}\sum_{i\in I}\left\langle T_{F}g_{i},g_{i}\right\rangle  & =\sum_{i\in I}\sum_{\ell=1}^{\infty}c_{\ell}\cdot\left\langle g_{i},u_{\ell}\right\rangle \left\langle u_{\ell},g_{i}\right\rangle \\
 & =\sum_{\ell=1}^{\infty}c_{\ell}\sum_{i\in I}\left|\left\langle u_{\ell},g_{i}\right\rangle \right|^{2}\\
\left({\scriptstyle \text{since }\left\Vert u_{\ell}\right\Vert _{L^{2}}=1}\right) & \leq C\cdot\sum_{\ell=1}^{\infty}c_{\ell}=C\cdot\left\Vert F\right\Vert _{L^{1}\left(\nu\right)}\cdot\left\Vert \psi\right\Vert _{L^{2}}^{2}<\infty.
\end{split}
\label{eq:ShearletFrameVanishingMomentsBesselTraceClassCombo}
\end{equation}

Next, choose an arbitrary compact set $\Lambda\subset H_{+}:=\left\{ \left(\begin{smallmatrix}a & b\\
0 & \sqrt{a}
\end{smallmatrix}\right)\with a\in\left(0,\infty\right)\text{ and }b\in\R\right\} $ with nonempty interior. By compactness, the constant $C_{0}:=\sup_{\lambda\in\Lambda}\left\Vert \lambda^{-1}\right\Vert _{\ell^{\infty}\to\ell^{\infty}}$
is finite. Now, define $\omega:=\Indicator_{\left[-3\delta C_{0},3\delta C_{0}\right]^{2}}\in L^{1}\left(\R^{2}\right)$
and let
\[
F:\R^{2}\times H\to\left[0,\infty\right),\left(x,h\right)\mapsto\omega\left(h^{-1}x\right)\cdot\Indicator_{\Lambda}\left(h\right).
\]
With this definition, we have
\begin{align*}
\left\Vert F\right\Vert _{L^{1}\left(\nu\right)} & =\int_{\Lambda}\int_{\R^{2}}\left|\det h\right|^{-1}\cdot\omega\left(h^{-1}x\right)\d x\d\mu_{H}\left(h\right)\\
\left({\scriptstyle \text{for }y=h^{-1}x}\right) & =\int_{\Lambda}\int_{\R^{2}}\omega\left(y\right)\d y\d\mu_{H}\left(h\right)=\left(6\delta C_{0}\right)^{2}\cdot\mu_{H}\left(\Lambda\right)<\infty.
\end{align*}

Next, we define
\[
Q_{0}:=\left\{ \xi\in\R^{\ast}\times\R\with\left|\xi_{1}\right|\in\left[2^{-1},1\right)\text{ and }\frac{\xi_{2}}{\xi_{1}}\in\left[-2^{-1},2^{-1}\right)\right\} 
\]
and set $Q:=\Lambda^{T}\overline{Q_{0}}$, where it is not hard to
see that $Q\subset\R^{\ast}\times\R$ is compact. Consequently, there
is some $\psi\in\Schwartz\left(\R^{2}\right)\subset L^{2}\left(\R^{2}\right)$
with $\widehat{\psi}\in\TestFunctionSpace{\R^{\ast}\times\R}$ and
$\widehat{\psi}\geq0$, as well as $\widehat{\psi}\equiv1$ on $Q$.

Now, we observe that the part $\Psi\left(\gamma;\delta\right)=\left\{ \gamma_{j,k,m}\with\left(j,k\right)\in I\text{ and }m\in\Z^{2}\right\} $
of the shearlet system $\mathcal{SH}\left(\varphi,\gamma,\tilde{\gamma};\delta\right)$
satisfies
\begin{align*}
\gamma_{j,k,m} & =2^{\frac{3}{4}j}\cdot\gamma\left(S_{k}A_{2^{j}}\bullet-\delta m\right)\\
 & =2^{\frac{3}{4}j}\cdot\gamma\left(S_{k}A_{2^{j}}\left[\bullet-\left(S_{k}A_{2^{j}}\right)^{-1}\delta m\right]\right)\\
 & =L_{x_{j,k,m}}D_{\left(S_{k}A_{2^{j}}\right)^{-1}}\gamma=\pi\left(x_{j,k,m},\,\left(S_{k}A_{2^{j}}\right)^{-1}\right)\gamma\quad\text{ with }\quad x_{j,k,m}:=\left(S_{k}A_{2^{j}}\right)^{-1}\delta m\in\R^{2}
\end{align*}
for $\left(j,k\right)\in I$ and $m\in\Z^{2}$. Consequently, since
each map $\pi\left(x,h\right)$ is unitary,
\begin{align*}
\left[W_{\psi}\gamma_{j,k,m}\right]\left(x,h\right)=\left\langle \gamma_{j,k,m},\,\pi\left(x,h\right)\psi\right\rangle  & =\left\langle \pi\left(x_{j,k,m},\,\left(S_{k}A_{2^{j}}\right)^{-1}\right)\gamma,\,\pi\left(x,h\right)\psi\right\rangle \\
 & =\left\langle \gamma,\,\pi\left(-S_{k}A_{2^{j}}x_{j,k,m},\,S_{k}A_{2^{j}}\right)\pi\left(x,h\right)\psi\right\rangle \\
 & =\left\langle \gamma,\,\pi\left(S_{k}A_{2^{j}}x-\delta m,\,S_{k}A_{2^{j}}\cdot h\right)\psi\right\rangle \qquad\forall\left(x,h\right)\in\R^{2}\times H.
\end{align*}
Since $\mathcal{SH}\left(\varphi,\gamma,\tilde{\gamma};\delta\right)$
is a Bessel system, so is $\Psi\left(\gamma;\delta\right)$. In view
of equation (\ref{eq:ShearletFrameVanishingMomentsBesselTraceClassCombo}),
this implies
\[
\sum_{\left(j,k\right)\in I,\,m\in\Z^{2}}\left\langle T_{F}\gamma_{j,k,m},\,\gamma_{j,k,m}\right\rangle <\infty.
\]
Now, let $\left(j,k\right)\in I$ and $m\in\Z^{2}$ be arbitrary and
note
\begin{align}
\left\langle T_{F}\gamma_{j,k,m},\,\gamma_{j,k,m}\right\rangle  & =\int_{\R^{2}\times H}F\left(x,h\right)\cdot\left|W_{\psi}\gamma_{j,k,m}\left(x,h\right)\right|^{2}\d\nu\left(x,h\right)\nonumber \\
 & =\int_{H}\Indicator_{\Lambda}\left(h\right)\cdot\int_{\R^{2}}\left|\det h\right|^{-1}\cdot\omega\left(h^{-1}x\right)\cdot\left|\left\langle \gamma,\,\pi\left(S_{k}A_{2^{j}}x-\delta m,\,S_{k}A_{2^{j}}\cdot h\right)\psi\right\rangle \right|^{2}\d x\,\d\mu_{H}\left(h\right)\nonumber \\
\left({\scriptstyle \text{for }z=h^{-1}x}\right) & =\int_{H}\Indicator_{\Lambda}\left(h\right)\cdot\int_{\R^{2}}\omega\left(z\right)\cdot\left|\left\langle \gamma,\,\pi\left(S_{k}A_{2^{j}}\cdot hz-\delta m,\,S_{k}A_{2^{j}}\cdot h\right)\psi\right\rangle \right|^{2}\d z\,\d\mu_{H}\left(h\right)\nonumber \\
\left({\scriptstyle \text{for }g=S_{k}A_{2^{j}}h}\right) & =\int_{H}\Indicator_{\Lambda}\left(\left(S_{k}A_{2^{j}}\right)^{-1}g\right)\cdot\int_{\R^{2}}\omega\left(z\right)\cdot\left|\left\langle \gamma,\,\pi\left(g\cdot z-\delta m,\,g\right)\psi\right\rangle \right|^{2}\d z\,\d\mu_{H}\left(g\right)\nonumber \\
\left({\scriptstyle \text{for }x=g\cdot z-\delta m}\right) & =\int_{H}\Indicator_{S_{k}A_{2^{j}}\Lambda}\left(g\right)\cdot\int_{\R^{2}}\omega\left(g^{-1}\left[x+\delta m\right]\right)\cdot\left|\left\langle \gamma,\,\pi\left(x,\,g\right)\psi\right\rangle \right|^{2}\d x\,\frac{\d\mu_{H}\left(g\right)}{\left|\det g\right|}.\label{eq:ShearletFrameVanishingMomentsTfGammaIdentity}
\end{align}

Next, let $\left(j,k\right)\in I$ and $x\in\R^{2}$ be fixed and
let $g\in S_{k}A_{2^{j}}\Lambda$, i.e., $g=S_{k}A_{2^{j}}\lambda$
for some $\lambda\in\Lambda$. Observe that $\Z^{2}\to\Z^{2},m\mapsto S_{k}m$
is a bijection, so that
\begin{align}
\sum_{m\in\Z^{2}}\omega\left(g^{-1}\left[x+\delta m\right]\right) & =\sum_{m\in\Z^{2}}\omega\left(\lambda^{-1}A_{2^{j}}^{-1}\left[S_{k}^{-1}x+\delta S_{k}^{-1}m\right]\right)\nonumber \\
\left({\scriptstyle \text{with }n=S_{k}^{-1}m}\right) & =\sum_{n\in\Z^{2}}\omega\left(\lambda^{-1}A_{2^{j}}^{-1}\left[S_{k}^{-1}x+\delta n\right]\right).\label{eq:ShearletFrameVanishingMomentsOmegaSeriesFirstStep}
\end{align}
But we have $\R^{2}=\biguplus_{n\in\Z^{2}}\left[n+\smash{\left[0,1\right)^{2}}\right]$,
so that there is some $n_{0}=n_{0}\left(x,k,\delta\right)\in\Z^{2}$
satisfying $-\frac{1}{\delta}S_{k}^{-1}x\in n_{0}+\left[0,1\right)^{2}$.
Hence, $\left\Vert S_{k}^{-1}x+\delta n_{0}\right\Vert _{\infty}=\delta\left\Vert -\frac{1}{\delta}S_{k}^{-1}x-n_{0}\right\Vert _{\infty}\leq\delta$.
Now, let
\[
Z_{j}:=\left\{ -2^{j},\dots,2^{j}\right\} \times\left\{ -\left\lceil \smash{2^{j/2}}\right\rceil ,\dots,\left\lceil \smash{2^{j/2}}\right\rceil \right\} 
\]
and observe for $n\in n_{0}+Z_{j}$ that
\begin{align*}
S_{k}^{-1}x+\delta n=S_{k}^{-1}x+\delta n_{0}+\delta\left(n-n_{0}\right) & \in\delta\left(\left[-1,1\right]^{2}+Z_{j}\right)\\
 & \subset\delta\cdot\left(\left[-\left(1+2^{j}\right),1+2^{j}\right]\times\left[-\left(1+\left\lceil \smash{2^{j/2}}\right\rceil \right),1+\left\lceil \smash{2^{j/2}}\right\rceil \right]\right)\\
 & \subset\delta\cdot\left(\left[-2^{1+j},2^{1+j}\right]\times\left[-3\cdot2^{j/2},\,3\cdot2^{j/2}\right]\right),
\end{align*}
since $1+\left\lceil 2^{j/2}\right\rceil \leq2+2^{j/2}\leq3\cdot2^{j/2}$.
Consequently, we get $A_{2^{j}}^{-1}\left[S_{k}^{-1}x+\delta n\right]\in\delta\cdot\left(\left[-2,2\right]\times\left[-3,3\right]\right)\subset\left[-3\delta,3\delta\right]^{2}$.
But by choice of $C_{0}$ from above, we have $\left\Vert \lambda^{-1}\right\Vert _{\ell^{\infty}\to\ell^{\infty}}\leq C_{0}$
and thus $\lambda^{-1}A_{2^{j}}^{-1}\left[S_{k}^{-1}x+\delta n\right]\in\left[-3\delta C_{0},3\delta C_{0}\right]^{2}$,
so that $\omega\left(\lambda^{-1}A_{2^{j}}^{-1}\left[S_{k}^{-1}x+\delta n\right]\right)=1$
for all $n\in n_{0}+Z_{j}$. In combination with equation (\ref{eq:ShearletFrameVanishingMomentsOmegaSeriesFirstStep})
and in view of $\left|Z_{j}\right|=\left(1+2\cdot2^{j}\right)\left(1+2\cdot\left\lceil 2^{j/2}\right\rceil \right)\geq2^{\frac{3}{2}j}$,
we thus get
\[
\sum_{m\in\Z^{2}}\omega\left(g^{-1}\left[x+\delta m\right]\right)\geq2^{\frac{3}{2}j}\qquad\forall x\in\R^{2},\left(j,k\right)\in I\text{ and }g\in S_{k}A_{2^{j}}\Lambda.
\]

Hence, equation (\ref{eq:ShearletFrameVanishingMomentsTfGammaIdentity})
yields the following estimate:
\begin{align*}
\sum_{m\in\Z^{2}}\left\langle T_{F}\gamma_{j,k,m},\,\gamma_{j,k,m}\right\rangle  & \geq2^{\frac{3}{2}j}\cdot\int_{H}\Indicator_{S_{k}A_{2^{j}}\Lambda}\left(g\right)\cdot\int_{\R^{2}}\left|\left\langle \gamma,\,\pi\left(x,\,g\right)\psi\right\rangle \right|^{2}\d x\,\frac{\d\mu_{H}\left(g\right)}{\left|\det g\right|}\\
 & =2^{\frac{3}{2}j}\cdot\int_{H}\Indicator_{S_{k}A_{2^{j}}\Lambda}\left(g\right)\cdot\left\Vert \left\langle \gamma,\,\pi\left(\bullet,\,g\right)\psi\right\rangle \right\Vert _{L^{2}}^{2}\,\frac{\d\mu_{H}\left(g\right)}{\left|\det g\right|}.
\end{align*}
Now, with the modulation operator $\left[M_{x}f\right]\left(\xi\right)=e^{2\pi i\left\langle x,\xi\right\rangle }\cdot f\left(\xi\right)$,
Plancherel's theorem yields
\[
\left\langle \gamma,\,\pi\left(x,\,g\right)\psi\right\rangle =\left\langle \gamma,\,L_{x}D_{g}\psi\right\rangle =\left\langle \widehat{\gamma},\,M_{-x}\Fourier\left[D_{g}\psi\right]\right\rangle =\Fourier^{-1}\left[\widehat{\gamma}\cdot\overline{\Fourier\left[D_{g}\psi\right]}\right]\left(x\right),
\]
where $\widehat{\gamma}\cdot\overline{\Fourier\left[D_{g}\psi\right]}\in L^{2}\left(\R^{2}\right)\cap L^{1}\left(\R^{2}\right)$,
since $\widehat{\gamma}\in L^{2}\left(\R^{2}\right)$ and $\overline{\Fourier\left[D_{g}\psi\right]}=\overline{D_{g^{-T}}\widehat{\psi}}\in\TestFunctionSpace{\R^{2}}$.
Consequently, another application of Plancherel's theorem shows
\begin{align*}
\sum_{m\in\Z^{2}}\left\langle T_{F}\gamma_{j,k,m},\,\gamma_{j,k,m}\right\rangle  & \geq2^{\frac{3}{2}j}\cdot\int_{H}\Indicator_{S_{k}A_{2^{j}}\Lambda}\left(g\right)\cdot\left\Vert \left\langle \gamma,\,\pi\left(\bullet,\,g\right)\psi\right\rangle \right\Vert _{L^{2}}^{2}\,\frac{\d\mu_{H}\left(g\right)}{\left|\det g\right|}\\
 & =2^{\frac{3}{2}j}\cdot\int_{H}\Indicator_{S_{k}A_{2^{j}}\Lambda}\left(g\right)\cdot\int_{\R^{2}}\left|\widehat{\gamma}\left(\xi\right)\right|^{2}\cdot\left|\left(D_{g^{-T}}\,\smash{\widehat{\psi}}\,\,\right)\left(\xi\right)\right|^{2}\d\xi\,\frac{\d\mu_{H}\left(g\right)}{\left|\det g\right|}\\
 & =2^{\frac{3}{2}j}\cdot\int_{H}\Indicator_{S_{k}A_{2^{j}}\Lambda}\left(g\right)\cdot\int_{\R^{2}}\left|\widehat{\gamma}\left(\xi\right)\right|^{2}\cdot\left|\widehat{\psi}\left(g^{T}\xi\right)\right|^{2}\d\xi\,\d\mu_{H}\left(g\right)\\
 & =\int_{\R^{2}}\left|\widehat{\gamma}\left(\xi\right)\right|^{2}\int_{H}2^{\frac{3}{2}j}\cdot\Indicator_{S_{k}A_{2^{j}}\Lambda}\left(g\right)\cdot\left|\widehat{\psi}\left(g^{T}\xi\right)\right|^{2}\d\mu_{H}\left(g\right)\,\d\xi
\end{align*}
for arbitrary $\left(j,k\right)\in I$. Now, set 
\begin{align*}
G & :=\left\{ \xi\in\R^{2}\with\left|\xi_{1}\right|\in\left(0,1\right)\text{ and }\left|\xi_{2}\right|\leq8^{-1}\cdot\left|\xi_{1}\right|^{1/2}\right\} \\
\text{and }G_{j} & :=\left\{ \xi\in G\with\left|\xi_{1}\right|\in\left[2^{-j-1},2^{-j}\right)\right\} ,\quad\text{ for }j\in\N_{0}
\end{align*}
and observe $G=\biguplus_{j=0}^{\infty}G_{j}$. In view of equation
(\ref{eq:ShearletFrameVanishingMomentsBesselTraceClassCombo}), the
preceding estimate implies
\begin{equation}
\infty>\sum_{\left(j,k\right)\in I}\int_{G}\left|\widehat{\gamma}\left(\xi\right)\right|^{2}\int_{H}2^{\frac{3}{2}j}\cdot\Indicator_{S_{k}A_{2^{j}}\Lambda}\left(g\right)\cdot\left|\widehat{\psi}\left(g^{T}\xi\right)\right|^{2}\d\mu_{H}\left(g\right)\,\d\xi.\label{eq:ShearletFrameVanishingMomentAlmostDone}
\end{equation}

Now, we need the following auxiliary claim:
\begin{equation}
\forall j\in\N_{0}\text{ and }k\in\Z\text{ with }\left|k\right|\leq\frac{1}{4}2^{j/2}:\quad G_{j}\subset\left(S_{k}A_{2^{j}}\right)^{-T}Q_{0}.\label{eq:ShearletFrameVanishingMomentSpecialSetInclusion}
\end{equation}
To see that this is true, first note for $\xi\in\R^{2}$ that $\xi\in\left(S_{k}A_{2^{j}}\right)^{-T}Q_{0}$
is equivalent to $\left(S_{k}A_{2^{j}}\right)^{T}\xi\in Q_{0}$. By
computing $\left(S_{k}A_{2^{j}}\right)^{T}\xi$ explicitly, we thus
get the following equivalence:
\begin{align*}
\xi\in\left(S_{k}A_{2^{j}}\right)^{-T}Q_{0} & \Longleftrightarrow\left(\begin{matrix}2^{j}\xi_{1}\\
2^{j/2}\left(\xi_{2}+k\xi_{1}\right)
\end{matrix}\right)\in Q_{0}\\
 & \Longleftrightarrow2^{j}\left|\xi_{1}\right|\in\left[2^{-1},1\right)\text{ and }\frac{2^{\frac{j}{2}}\left(\xi_{2}+k\xi_{1}\right)}{2^{j}\xi_{1}}\in\left[-2^{-1},2^{-1}\right)\\
 & \Longleftrightarrow\left|\xi_{1}\right|\in\left[2^{-j-1},2^{-j}\right)\text{ and }\frac{\xi_{2}}{\xi_{1}}\in\left[-2^{\frac{j}{2}-1},2^{\frac{j}{2}-1}\right)-k.
\end{align*}
Hence, to prove the auxiliary claim (\ref{eq:ShearletFrameVanishingMomentSpecialSetInclusion}),
we only need to verify that this last condition is fulfilled for $\xi\in G_{j}$
and $\left|k\right|\leq\frac{1}{4}2^{j/2}$. But $\left|\xi_{1}\right|\in\left[2^{-j-1},2^{-j}\right)$
simply holds by definition of $G_{j}$. For the second condition,
we note from the definition of $G$ that 
\[
\left|\frac{\xi_{2}}{\xi_{1}}\right|\leq\frac{1}{8}\cdot\left|\xi_{1}\right|^{-1/2}\leq\frac{1}{8}\cdot\left(2^{-j-1}\right)^{-1/2}=\frac{1}{8}\cdot2^{\frac{j}{2}+\frac{1}{2}}<2^{\frac{j}{2}-2}
\]
and hence $\frac{\xi_{2}}{\xi_{1}}\in\left(-2^{\frac{j}{2}-2},\,2^{\frac{j}{2}-2}\right)\subset\left[-2^{\frac{j}{2}-1},2^{\frac{j}{2}-1}\right)-k$.
Here, the last inclusion is indeed valid, since we have $\left|k\right|\leq\frac{1}{4}2^{j/2}$
and thus
\[
-2^{\frac{j}{2}-1}-k\leq-2^{\frac{j}{2}-1}+\left|k\right|\leq-\frac{1}{4}2^{\frac{j}{2}}=-2^{\frac{j}{2}-2},\quad\text{ as well as }\quad2^{\frac{j}{2}-1}-k\geq2^{\frac{j}{2}-1}-\left|k\right|\geq\frac{1}{4}\cdot2^{\frac{j}{2}}=2^{\frac{j}{2}-2}.
\]

Finally, note that $\left|k\right|\leq\frac{1}{4}2^{j/2}$ in particular
implies $\left|k\right|\leq\left\lceil 2^{j/2}\right\rceil $ and
thus $\left(j,k\right)\in I$. Hence, a combination of equation (\ref{eq:ShearletFrameVanishingMomentAlmostDone})
with the auxiliary claim (\ref{eq:ShearletFrameVanishingMomentSpecialSetInclusion})
yields
\begin{align*}
\infty & >\sum_{\left(j,k\right)\in I}\int_{G}\left|\widehat{\gamma}\left(\xi\right)\right|^{2}\int_{H}2^{\frac{3}{2}j}\cdot\Indicator_{S_{k}A_{2^{j}}\Lambda}\left(g\right)\cdot\left|\widehat{\psi}\left(g^{T}\xi\right)\right|^{2}\d\mu_{H}\left(g\right)\,\d\xi\\
 & \geq\sum_{j=0}^{\infty}\int_{G_{j}}2^{\frac{3}{2}j}\cdot\left|\widehat{\gamma}\left(\xi\right)\right|^{2}\sum_{k\in\Z,\,\left|k\right|\leq2^{\frac{j}{2}-2}}\int_{\Lambda}\left|\widehat{\psi}\left(\left(S_{k}A_{2^{j}}\lambda\right)^{T}\cdot\xi\right)\right|^{2}\d\mu_{H}\left(\lambda\right)\,\d\xi\\
 & \geq\mu_{H}\left(\Lambda\right)\cdot\sum_{j=0}^{\infty}\int_{G_{j}}2^{\frac{3}{2}j}\cdot\left|\widehat{\gamma}\left(\xi\right)\right|^{2}\cdot\left|\left\{ k\in\Z\,:\,\left|k\right|\leq2^{\frac{j}{2}-2}\right\} \right|\,\d\xi,
\end{align*}
where the last step used that each $\xi\in G_{j}$ satisfies $\xi\in\left(S_{k}A_{2^{j}}\right)^{-T}Q_{0}$
and thus $\left(S_{k}A_{2^{j}}\lambda\right)^{T}\xi\in\Lambda^{T}Q_{0}\subset Q$.
Since $\widehat{\psi}\equiv1$ on $Q$, this implies $\widehat{\psi}\left(\left(S_{k}A_{2^{j}}\lambda\right)^{T}\xi\right)=1$
for all $\lambda\in\Lambda$ and $\xi\in G_{j}$ and $k\in\Z$ with
$\left|k\right|\leq\frac{1}{4}2^{j/2}$.

Now, there are two cases: For $2^{\frac{j}{2}-2}\leq1$, we have $\left|\left\{ k\in\Z\,:\,\left|k\right|\leq2^{\frac{j}{2}-2}\right\} \right|\geq1\geq2^{\frac{j}{2}-2}$.
If otherwise $2^{\frac{j}{2}-2}>1$, then $2^{\frac{j}{2}-2}\leq1+\left\lfloor \smash{2^{\frac{j}{2}-2}}\right\rfloor \leq2\cdot\left\lfloor \smash{2^{\frac{j}{2}-2}}\right\rfloor $,
so that
\[
\left|\left\{ k\in\Z\,:\,\left|k\right|\leq2^{\frac{j}{2}-2}\right\} \right|\geq\left|\left\{ -\left\lfloor \smash{2^{\frac{j}{2}-2}}\right\rfloor ,\dots,\left\lfloor \smash{2^{\frac{j}{2}-2}}\right\rfloor \right\} \right|=1+2\cdot\left\lfloor \smash{2^{\frac{j}{2}-2}}\right\rfloor \geq2^{\frac{j}{2}-2}
\]
as well. Finally, for $\xi\in G_{j}$, we have $\left|\xi_{1}\right|\geq2^{-j-1}=2^{-\left(j+1\right)}$
and thus
\[
2^{\frac{3}{2}j}\cdot2^{\frac{j}{2}-2}=\frac{1}{16}\cdot2^{2\left(j+1\right)}\geq\frac{1}{16}\cdot\left|\xi_{1}\right|^{-2}.
\]
Putting everything together, we arrive at
\begin{align*}
\infty & >\mu_{H}\left(\Lambda\right)\cdot\sum_{j=0}^{\infty}\int_{G_{j}}2^{\frac{3}{2}j}\cdot\left|\widehat{\gamma}\left(\xi\right)\right|^{2}\cdot\left|\left\{ k\in\Z\,:\,\left|k\right|\leq2^{\frac{j}{2}-2}\right\} \right|\,\d\xi\\
 & \geq\frac{\mu_{H}\left(\Lambda\right)}{16}\cdot\sum_{j=0}^{\infty}\int_{G_{j}}\left|\xi_{1}\right|^{-2}\cdot\left|\widehat{\gamma}\left(\xi\right)\right|^{2}\,\d\xi\\
 & =\frac{\mu_{H}\left(\Lambda\right)}{16}\cdot\int_{G}\left|\xi_{1}\right|^{-2}\cdot\left|\widehat{\gamma}\left(\xi\right)\right|^{2}\,\d\xi,
\end{align*}
as desired.
\end{proof}

\section*{Acknowledgments}

I would like to thank Jackie Ma for raising the question whether decomposition
spaces (in particular shearlet smoothness spaces) can be characterized
using compactly supported functions. I warmly thank Hartmut Führ,
Charly Gröchenig and Anne Pein for fruitful discussions related to
the topics in this paper. Furthermore, I especially thank Hartmut
Führ for several suggestions which led to great improvements of this
paper, in particular of the abstract, the introduction, and the discussion
of related literature.

The author acknowledges support from the European Commission through
DEDALE (contract no.\@ 665044) within the H2020 Framework Program.

\bibliographystyle{plain}
\bibliography{felixbib}

\end{document}